\documentclass[reqno,english]{amsart}
\usepackage[]{amsmath,amssymb,amsfonts,latexsym,amsthm,enumerate}
\usepackage{mathtools}
\usepackage{upgreek}
\usepackage{mathrsfs}
\usepackage{enumerate,graphicx}
\usepackage{tikz}
\usetikzlibrary{arrows.meta}
\usepackage[font=footnotesize]{caption}
\usepackage[noadjust]{cite}

\numberwithin{equation}{section}
\usepackage[psdextra=true,colorlinks=true, pdfstartview=Fit, linkcolor=blue,
  citecolor=blue, urlcolor=blue,pagebackref=false]{hyperref}
\usepackage[shortlabels]{enumitem}
\usepackage[capitalise]{cleveref}
\usepackage{xspace}
\usepackage[letterpaper, width=145mm, margin=2.5cm, centering]{geometry}

\usepackage[color,final]{showkeys} %
\definecolor{refkey}{gray}{.75}
\definecolor{labelkey}{gray}{.5}

\newtheorem{maintheorem}{Theorem}
\newtheorem{theorem}{Theorem}[section]
\newtheorem{lemma}[theorem]{Lemma}

\newtheorem{proposition}[theorem]{Proposition}

\newtheorem{corollary}[theorem]{Corollary}
\newtheorem{remark}[theorem]{Remark}
\newtheorem*{remark*}{Remark}

\theoremstyle{definition}{

\newtheorem{definition}[theorem]{Definition}
\newtheorem*{definition*}{Definition}

}

\crefname{maintheorem}{Theorem}{Theorems}
\crefname{claim}{Claim}{Claims}
\crefname{observation}{Observation}{Observations}
\crefname{fact}{Fact}{Facts}
\crefname{condition}{Condition}{Conditions}
\crefname{step}{Step}{Steps}
\AddToHook{env/lemma/begin}{\crefalias{theorem}{lemma}}
\AddToHook{env/claim/begin}{\crefalias{theorem}{claim}}
\AddToHook{env/proposition/begin}{\crefalias{theorem}{proposition}}
\AddToHook{env/observation/begin}{\crefalias{theorem}{observation}}
\AddToHook{env/fact/begin}{\crefalias{theorem}{fact}}
\AddToHook{env/corollary/begin}{\crefalias{theorem}{corollary}}
\AddToHook{env/remark/begin}{\crefalias{theorem}{remark}}
\AddToHook{env/example/begin}{\crefalias{theorem}{example}}
\AddToHook{env/definition/begin}{\crefalias{theorem}{definition}}
\AddToHook{env/question/begin}{\crefalias{theorem}{question}}

\newcommand{\C}{\mathbb C}
\newcommand{\E}{\mathbb E}
\newcommand{\N}{\mathbb N}
\renewcommand{\P}{\mathbb P}
\newcommand{\Q}{\mathbb Q}
\newcommand{\R}{\mathbb R}
\newcommand{\T}{\mathbb T}
\newcommand{\Z}{\mathbb Z}
\newcommand{\cD}{\mathcal D}
\newcommand{\cK}{\mathcal K}
\newcommand{\cA}{{\mathcal A}}
\newcommand{\cB}{{\mathcal B}}
\newcommand{\cC}{{\mathcal C}}
\newcommand{\cE}{{\mathcal E}}
\newcommand{\cF}{{\mathcal F}}
\newcommand{\cG}{{\mathcal G}}
\newcommand{\cH}{{\mathcal H}}
\newcommand{\cI}{{\mathcal I}}
\newcommand{\cL}{{\mathcal L}}
\newcommand{\cM}{{\mathcal M}}
\newcommand{\cP}{{\mathcal P}}
\newcommand{\cR}{{\mathcal R}}
\newcommand{\cS}{{\mathcal S}}
\newcommand{\cT}{{\mathcal T}}
\newcommand{\cV}{{\mathcal V}}
\newcommand{\cW}{{\mathcal W}}
\newcommand{\cX}{{\mathcal X}}
\newcommand{\cZ}{{\mathcal Z}}
\newcommand{\fa}{\mathfrak a}
\newcommand{\fh}{\mathfrak h}
\newcommand{\fH}{\mathfrak H}
\newcommand{\fS}{\mathfrak S}
\newcommand{\fd}{\mathfrak d}
\newcommand{\fg}{\mathfrak g}
\newcommand{\fm}{\mathfrak m}
\newcommand{\fn}{\mathfrak n}
\renewcommand{\d}{\mathrm{d}}
\newcommand{\sB}{\mathsf B}
\newcommand{\sA}{\mathsf A}
\newcommand{\sfT}{\mathsf{T}}
\newcommand{\Tone}{\sfT_{\!1}}
\newcommand{\Ttwo}{\sfT_{\!2}}
\newcommand{\Tthree}{\sfT_{\!3}}
\newcommand{\Ti}{\sfT_{\!i}}
\newcommand{\xor}{\vartriangle}
\newcommand{\trps}{{\scriptscriptstyle\top}}
\newcommand{\fZ}{\mathfrak Z}
\newcommand{\fE}{\mathfrak E}
\newcommand{\scW}{\mathscr W}
\newcommand{\scH}{\mathscr H}
\newcommand{\scS}{\mathscr S}
\newcommand{\scD}{\mathscr D}
\newcommand{\scR}{\mathscr R}
\newcommand{\bN}{\boldsymbol N}
\newcommand{\bn}{\boldsymbol n}
\newcommand{\bx}{\boldsymbol x}

\newcommand{\be}{\boldsymbol e}

\newcommand{\scC}{\mathscr C}

\newcommand{\llb}{\mathopen{[\![}}
\newcommand{\rrb}{\mathclose{]\!]}}
\newcommand{\Dim}{\textup{\textsc{d}} }
\DeclareMathOperator{\sgn}{sgn}
\DeclareMathOperator{\cov}{Cov}
\DeclareMathOperator{\var}{Var}
\DeclareMathOperator{\Dt}{dot}
\DeclareMathOperator{\Arr}{arr}
\newcommand{\Arrm}[1]{\Arr_{#1}^{\wedge}}
\DeclareMathOperator{\diam}{diam}
\DeclareMathOperator{\dist}{dist}
\DeclareMathOperator{\supp}{supp}
\DeclareMathOperator{\adj}{adj}
\DeclareMathOperator{\pv}{pv}
\DeclareMathOperator{\conv}{conv}
\DeclarePairedDelimiter{\ceil}{\lceil}{\rceil}
\DeclarePairedDelimiter{\floor}{\lfloor}{\rfloor}
\newcommand{\given}{\mid}
\newenvironment{widequote}{\list{}{\leftmargin=1em\rightmargin=1em\topsep=0pt\partopsep=0pt\parsep=0pt\itemsep=0pt}\item\relax}{\endlist\noindent\ignorespacesafterend}
\newcommand{\one}{\mathbf{1}}
\newcommand{\GFF}{\textsc{GFF}\xspace}
\newcommand{\ZGFF}{\ensuremath{\Z}\textsc{GFF}\xspace}
\newcommand{\SOS}{\textsc{SOS}\xspace}
\newcommand{\HMSOS}{\textsc{HMSOS}\xspace}
\newcommand{\RG}{\textsc{RG}\xspace}
\newcommand{\Tetris}{\textsf{Tetris}\xspace}
\renewcommand{\restriction}{\mathord{\upharpoonright}}
\renewcommand{\epsilon}{\varepsilon}

\title{GFF for 3D Ising interfaces in three moves of Tetris}

\author{Beno\^{\i}t Laslier}
\address{B.\ Laslier\hfill\break
Universit\'e de Rouen Normandie\\
LMRS\\
76801 Saint-\'Etienne-du-Rouvray\\ France.}
\email{benoit.laslier@univ-rouen.fr}
\author{Eyal Lubetzky}
\address{E.\ Lubetzky\hfill\break
Courant Institute\\ New York University\\
251 Mercer Street\\ New York, NY 10012, USA.}
\email{eyal@courant.nyu.edu}

\begin{document}

\begin{abstract}
We show that at sufficiently large inverse-temperature $\beta>0$, the 3\Dim Ising model on the torus, under boundary conditions that have a positive slope $\theta_1,\theta_2>0$, has a \GFF scaling limit for its interface. The proof goes through three steps of a \Tetris map, which lets cubes drop according to gravity in each of the three axes, yielding Solid-On-Solid (\SOS), Half-Monotone \SOS (\HMSOS) and monotone surfaces (dimers). We establish \GFF limits sequentially, in particular resolving the analogous problem for the $(2+1)$\Dim \SOS.
A crucial ingredient is an improved \RG analysis for interacting dimers, which is of independent interest: it relaxes the assumptions on interactions from being uniformly bounded in $L^\infty$  to only being typically small.
\end{abstract}

 {\mbox{}
 \vspace{-.5cm}
\maketitle
 }

\vspace{-0.5cm}

\setcounter{tocdepth}{2}
\tableofcontents

\section{Introduction}\label{sec:intro}

\begin{figure}
\vspace{-0.2cm}
\centering
\begin{tikzpicture}[x=0.81\textwidth,y=0.81\textwidth]
\node[anchor=south west,inner sep=0pt,outer sep=0pt] at (0.0000,0.0234) {\includegraphics[width=0.2552\textwidth]{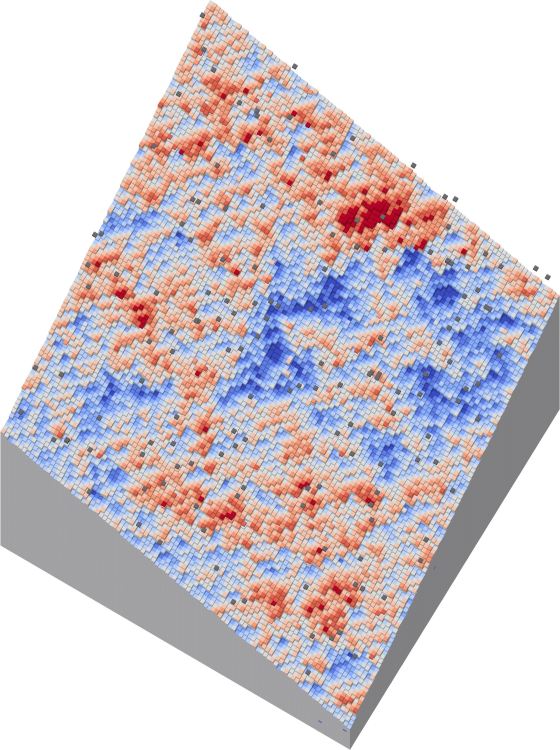}};
\node[anchor=south west,inner sep=0pt,outer sep=0pt] at (0.6882,0.0234) {\includegraphics[width=0.2525\textwidth]{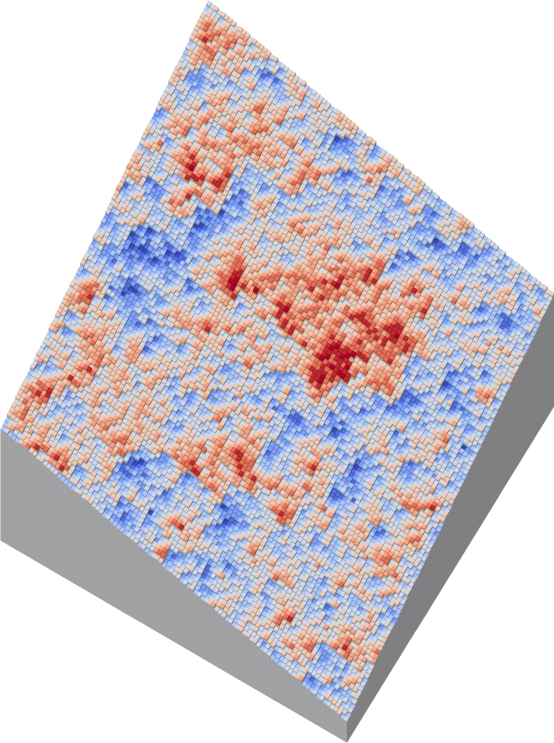}};
\node[anchor=south west,inner sep=0pt,outer sep=0pt] at (0.4056,0.3218) {\includegraphics[width=0.1529\textwidth]{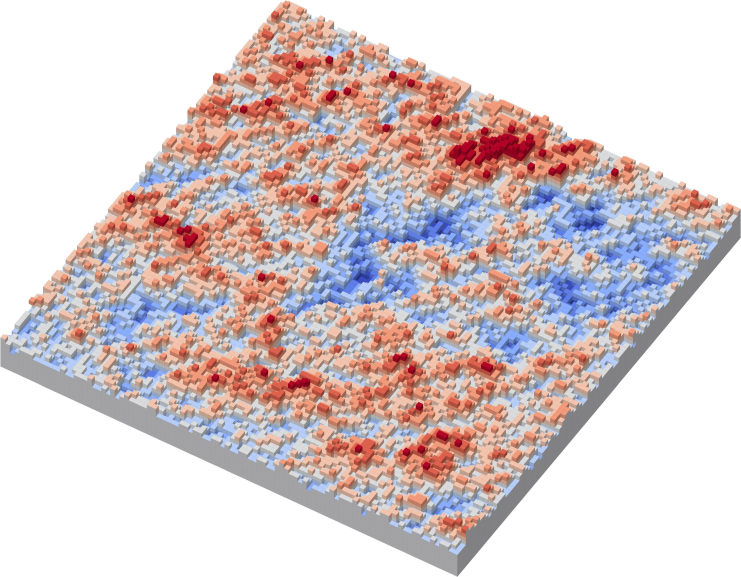}};
\node[anchor=south west,inner sep=0pt,outer sep=0pt] at (0.4047,0.1591) {\includegraphics[width=0.1544\textwidth]{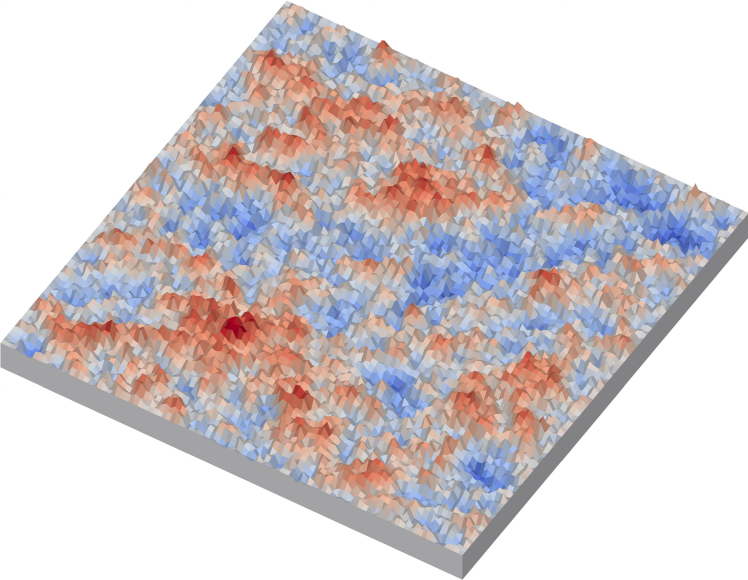}};
\node[anchor=south west,inner sep=0pt,outer sep=0pt] at (0.4066,-0.0000) {\includegraphics[width=0.1514\textwidth]{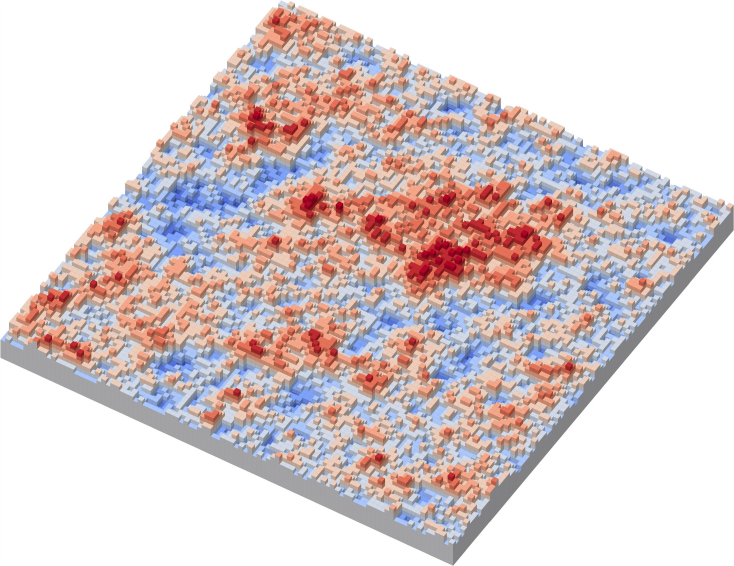}};
\draw[black!40,line width=.6pt] (0.3148,0.2781) -- (0.4057,0.3818);
\fill[black!40] (0.3148,0.2781) circle[radius=.9pt];
\fill[black!40] (0.4057,0.3818) circle[radius=.9pt];
\draw[black!40,line width=.6pt] (0.6884,0.1998) -- (0.5934,0.0910);
\fill[black!40] (0.6884,0.1998) circle[radius=.9pt];
\fill[black!40] (0.5934,0.0910) circle[radius=.9pt];
\end{tikzpicture}

\vspace{-0.3cm}
\caption{3\Dim Ising (left) and \SOS (right) at $\beta=1.2$ with $(1,1,1)$-tilted boundary conditions. Middle: their interfaces with the reference plane subtracted, and a discrete \GFF (rescaled to the \SOS variance) in between.}
\label{fig:intro-gff}
\vspace{-0.2cm}
\end{figure}
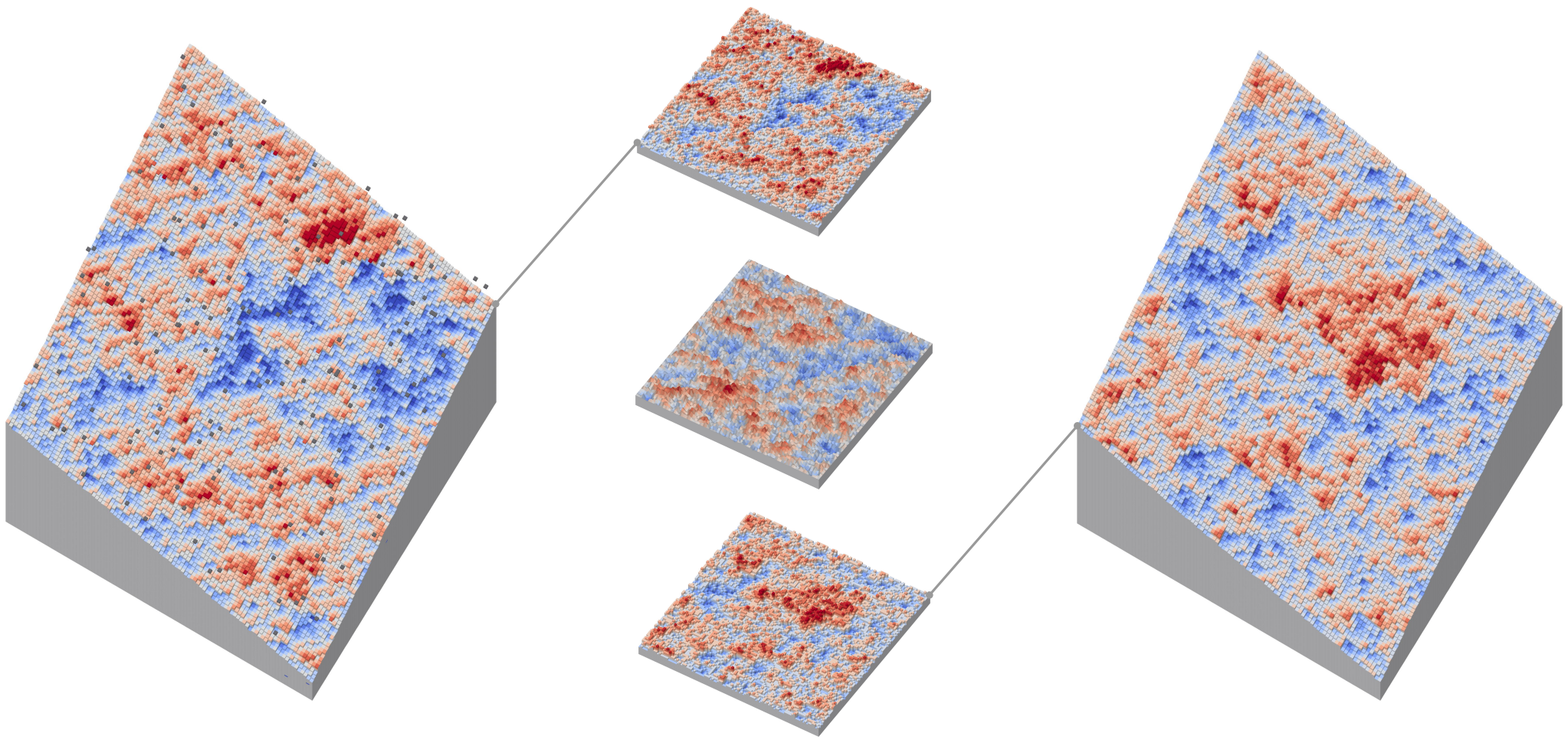

The low-temperature 3\Dim Ising model under Dobrushin boundary conditions---plus spins below a reference plane and minus spins above it---exhibits an interface separating the two pure phases. When the reference plane is horizontal, this interface is rigid at low temperature by the classical work of Dobrushin~\cite{Dobrushin72a}. When the reference plane is tilted, the interface is expected to be rough, with fluctuations of logarithmic order and a Gaussian free field (\GFF) scaling limit. This paper establishes the \GFF scaling limit of the tilted interface on the torus, for every pair of positive slopes and for all sufficiently large $\beta$, together with the analogous statements for the bubble-free Ising Interface model and for two solid-on-solid approximations of the 3\Dim Ising model: the $(2+1)$\Dim Solid-On-Solid (\SOS) model and the Half-Monotone \SOS (\HMSOS) model. In what follows, $x\sim y$ denotes a pair of nearest-neighbor sites in the lattice under consideration.

\begin{definition}[Ising model]\label{def:3dIsing}
The 3\Dim Ising model on a box $\Lambda_n=\llb -n,n\rrb^3$ with inverse-temperature $\beta>0$ and boundary conditions $\tau$ is the probability measure on configurations $\sigma$, which assign $\pm1$ spins to the cubes of  $\Lambda_n$ (i.e., to their centers in $(\frac12,\frac12,\frac12)+\Z^3$) and agree with $\tau$ on the external cube boundary $\partial\Lambda_n$, given by
\[
\P( \sigma ) \propto \exp \bigg(-\beta \sum_{x\sim y} \one_{\{\sigma(x)\neq\sigma(y)\}}\bigg)\,.
\]
Here each unoriented bond meeting the domain is counted once.

For the model on the finite cylinder $\T_n\times \llb -M,M\rrb$, where $\T_n=\llb1,n\rrb^2$ is a torus, one replaces the boundary condition $\tau$ by $+1$ on $x=(x_1,x_2,x_3)$ for $x_3 < -M$, by $-1$ for $x_3 > M$, and by the shifted copies
\[
\sigma(\cdot+(n,0,0))\,,\;\sigma(\cdot-(n,0,0))\,,\;\sigma(\cdot+(0,n,0))\,,\;\sigma(\cdot-(0,n,0))\,,
\]
accounting for the toric boundary. Equivalently, configurations are extended horizontally with period $n$; the quotient is a finite cylinder with fixed ends.

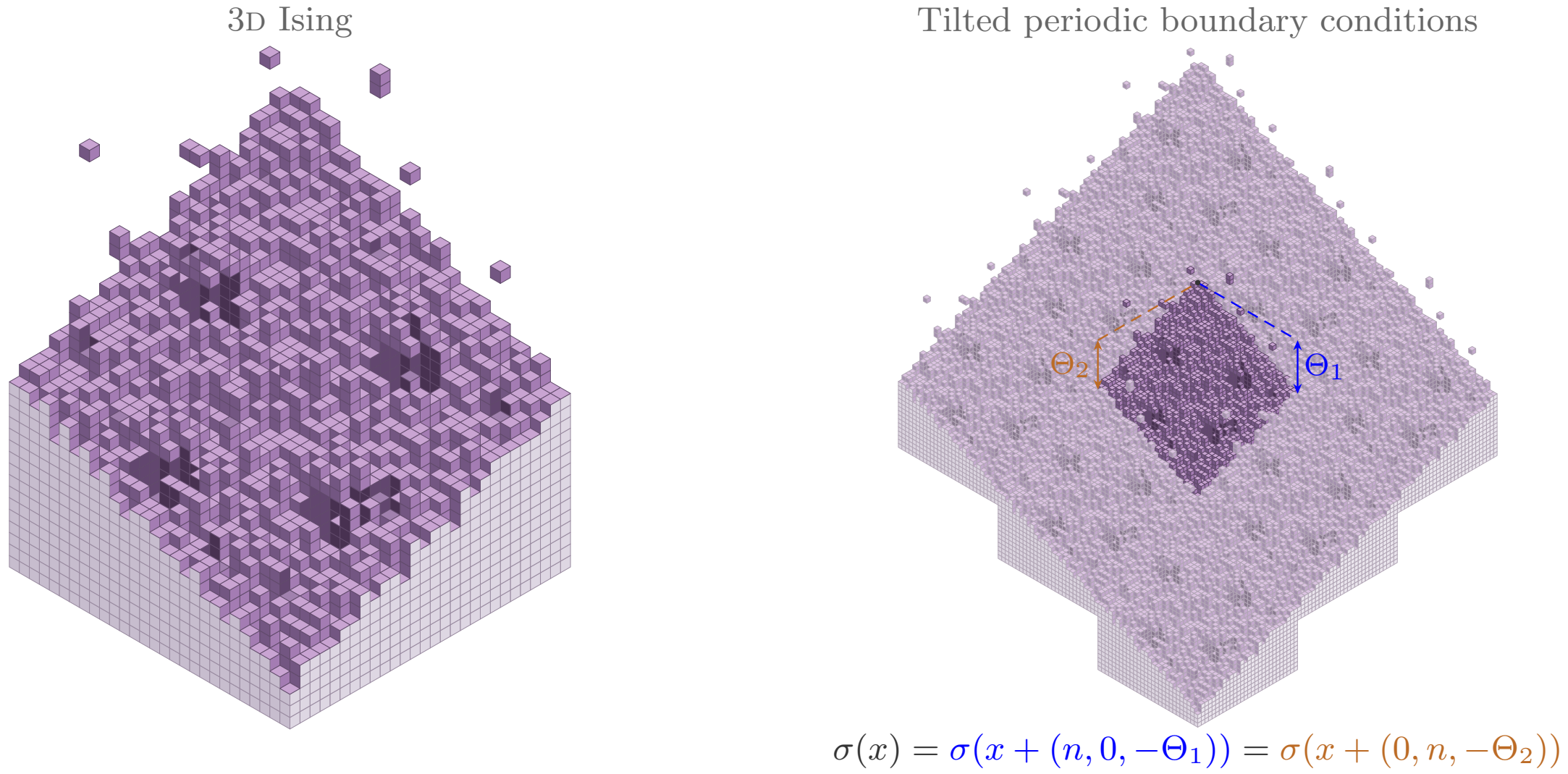
\begin{figure}
      \vspace{-0.2cm}
\centering
\definecolor{accent}{HTML}{BC6B25}
\begin{tikzpicture}[x=\textwidth,y=\textwidth]
\node[text=black!60,anchor=north,font=\small,inner sep=0pt] at (.22,-0.015) {3\Dim Ising};
\node[text=black!60,anchor=north,font=\small,inner sep=0pt] at (.72,-0.015) {Tilted periodic boundary conditions};
\node[anchor=north,inner sep=0pt,outer sep=0pt] at (.22,-.035)
    {\includegraphics[height=.38\textwidth]{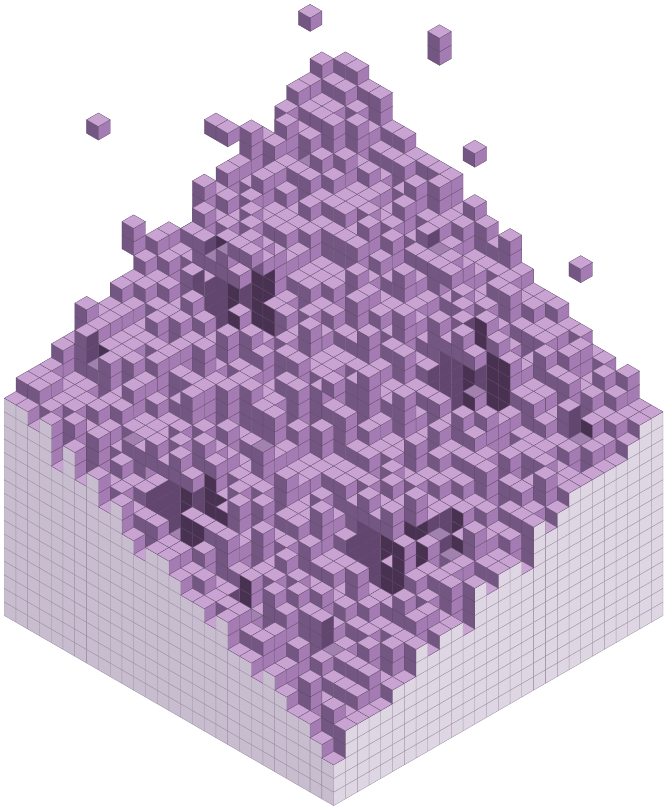}};
\node[anchor=north,inner sep=0pt,outer sep=0pt] at (.72,-.035)
    {\includegraphics[height=.38\textwidth]{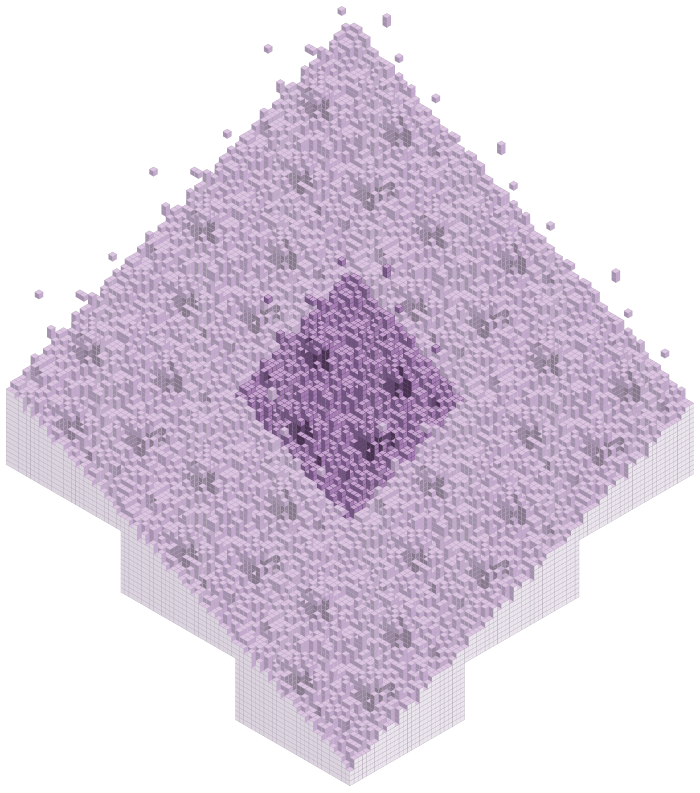}};
\begin{scope}[black!80,line width=.45pt,arr/.style={{Stealth[length=2.6pt,width=2.2pt]}-{Stealth[length=2.6pt,width=2.2pt]}}]
\draw[accent,densely dashed] (.6650,-.1989) -- (.7200,-.1672);
\draw[blue,densely dashed] (.7200,-.1672) -- (.7750,-.1989);
\fill (.7200,-.1672) circle[radius=.7pt];
\draw[arr,blue] (.7750,-.1989) -- node[text=blue,right,font=\footnotesize,inner sep=1.5pt]{$\Theta_1$} (.7750,-.2284);
\draw[arr,accent] (.6650,-.1989) -- node[text=accent,left,font=\footnotesize,inner sep=1.5pt]{$\Theta_2$} (.6650,-.2261);

\node[font=\small] at (0.72,-.425) {
$\sigma(x)={\color{blue}\sigma(x+(n,0,-\Theta_1))}={\color{accent}\sigma(x+(0,n,-\Theta_2))}$};

\end{scope}
\end{tikzpicture}
\vspace{-0.45cm}
\caption{A 3\Dim Ising configuration with overhangs and bubbles,
and its tilted periodic extension to $\Z^3$.
}
\label{fig:ising-toric}
\vspace{-0.15cm}
\end{figure}

The model on $\T_n$ with \emph{windings} $\Theta_1,\Theta_2$, which are positive integers, illustrated in~\cref{fig:ising-toric},
is defined similarly, except that the shifted copies are
\[
\sigma(\cdot+(n,0,-\Theta_1))\,,\;\sigma(\cdot-(n,0,-\Theta_1))\,,\;\sigma(\cdot+(0,n,-\Theta_2))\,,\;\sigma(\cdot-(0,n,-\Theta_2))\,,
\]
with the bottom boundary replaced by $+1$ if $x_3 < -M - (\Theta_1x_1 + \Theta_2x_2)/n$ and the top boundary by $-1$ when $x_3 > M - (\Theta_1x_1 + \Theta_2x_2)/n$.
That is, the configuration is extended to all of $\Z^3$ so that
it is invariant under the translations $(n,0,-\Theta_1)$ and
$(0,n,-\Theta_2)$: the interface descends by $\Theta_i$ over one
period in the $x_i$-direction, and we call $(\Theta_1/n,\Theta_2/n)$
its \emph{slope}. (The horizontal axes are oriented so that the
interface descends along the $+x_1$ and $+x_2$ directions.) The tilted slabs of fixed spins play
the role of the height cutoff $M$ of the horizontal case.
An \emph{admissible winding sequence} approaching a target slope
$(\theta_1,\theta_2)\in(0,\infty)^2$ assigns to each $n$ positive
integer windings $\Theta_1=\Theta_1(n)$ and $\Theta_2=\Theta_2(n)$ such that
$\Theta_1/n\to\theta_1$ and $\Theta_2/n\to\theta_2$; no integrality is
required of $\theta_1,\theta_2$ themselves.
We suppress the dependence of the windings on $n$.
For the full Ising limit, let $n\to\infty$ along admissible windings and choose a joint cutoff
\begin{equation}\label{eq:ising-height-cutoff}
                         M=M_n=\ceil{e^{n^2}}\,.
\end{equation}
With this choice and sufficiently large $\beta$, with high probability
(w.h.p.) there is a unique interface separating the two slabs, and its
law (with its bubbles
filled, and rooted as below) is asymptotically independent of the
cutoff: as $n\to\infty$ it is at total-variation distance tending
to zero from a law on bubble-free interfaces which does not involve
$M$ (this is the content of \cref{prop:ising-joint-cutoff}).
Since the energy is invariant under vertical shifts of the
configuration as long as the slabs are not reached,
we root the interface and the surrounding spins by translating it to have $\fh(0)=0$ for the height function $\fh$ associated to it (specified in \cref{def:height-function-in-models} below).
\end{definition}

\begin{definition}[Ising Interface model]\label{def:interface-model}
This is 3\Dim Ising from \cref{def:3dIsing} restricted to configurations without ``bubbles'': in either model on $\T_n$, every $+1$ spin of $\sigma$ must be connected to the bottom slab via a path of $+1$ spins, and similarly for $-1$ spins and the top slab. Spins are assigned to cube centers $(\frac12,\frac12,\frac12)+\Z^3$, and said paths can move between \emph{$*$-adjacent} cubes, defined here as cubes sharing a face or an edge (not just a corner), equivalently cubes whose centers are at Euclidean distance at most~$\sqrt{2}$. (The cylinder height is any $M=M_n$ with $M_n/n^2\to\infty$, including $M=\infty$.)
\end{definition}

\begin{figure}
\centering
\begin{minipage}[t]{.2\textwidth}
\centering
\small {\color{black!60}Ising interface}\par\smallskip
\includegraphics[width=\linewidth]{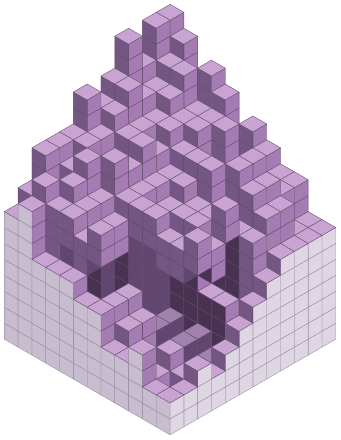}
\end{minipage}\hspace{.075\textwidth}%
\begin{minipage}[t]{.2\textwidth}
\centering
\small {\color{black!60}SOS}\par\smallskip
\includegraphics[width=\linewidth]{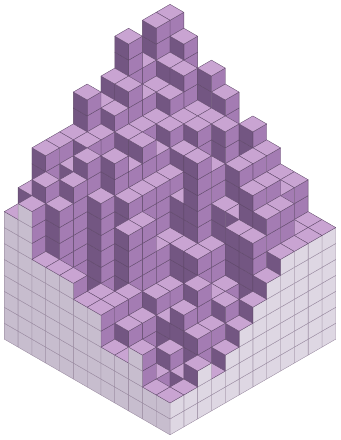}
\end{minipage}\hspace{.075\textwidth}%
\begin{minipage}[t]{.2\textwidth}
\centering
\small {\color{black!60}HMSOS}\par\smallskip
\includegraphics[width=\linewidth]{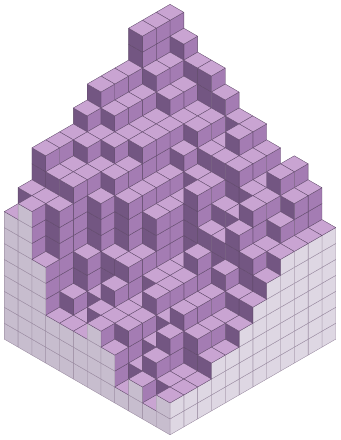}
\end{minipage}
\vspace{-0.3cm}
\caption{Tilted interfaces: the Ising interface (left) has
overhangs (darker faces mark recesses under them); \SOS (middle) is a graph of a height function;
 \HMSOS (right) heights are monotone in the $x_2$-direction~($\swarrow$).
  }
\label{fig:model-interfaces}
\end{figure}

The restriction thus forbids finite islands of either sign: the plus
spins form a single cluster attached to the bottom slab, the minus
spins a single cluster attached to the top slab, and the
configuration is determined by the surface separating them, which may have overhangs (see
\cref{fig:model-interfaces}).

\begin{remark*}
    The 3\Dim Ising Interface model can equivalently be seen as a model on a set of plaquettes (the dual faces of edges) of $\Z^3$, by describing a configuration via the set of edges $x\sim y$ of $\Z^3$ with $\sigma(x) \neq \sigma(y)$. The energy of a configuration is then the number of plaquettes in it.
\end{remark*}

\begin{definition}[Solid-On-Solid model]\label{def:sos}
The $(2+1)$\Dim \SOS model is the 3\Dim Ising model in the limit where the $x_3$-direction coupling constant (between $(x_1,x_2,x_3)$ and $(x_1,x_2,x_3+1)$) goes to~$\infty$. Equivalently, one can take the 3\Dim Ising Interface model and impose the same condition on the vertical coupling, or require that every column $(x_1, x_2) \times \Z$ contains exactly one disagreement.
\end{definition}

\begin{remark*}
    Equivalently, the \SOS model is an integer height
    $\fh:\Z^2\to\Z$ satisfying
    $\fh(x+ne_i)=\fh(x)-\Theta_i$, $i=1,2$, with energy
    $\sum_{x\in\T_n}\sum_{i=1}^2|\fh(x+e_i)-\fh(x)|$.
    Indeed, $\fh(x)$ records the height of the unique disagreement
    in the column above $x$; the vertical plaquettes of the interface
    between the columns above $x$ and $x+e_i$ number
    $|\fh(x+e_i)-\fh(x)|$, while the horizontal plaquettes
    contribute exactly $n^2$, dropped
    from the energy.
    The tilted periodic boundary conditions of \cref{def:3dIsing} translate into the stated periodicity of $\fh$.
\end{remark*}

\begin{definition}[Half-Monotone \SOS model]\label{def:HMSOS}
The $(2+1)$\Dim \HMSOS model is the 3\Dim Ising model in the limit where the $x_2$-direction and $x_3$-direction coupling constants both go to $\infty$. That is, it is the \SOS model with the additional constraint that $\fh(x+e_2)\leq\fh(x)$ on the height function $\fh$.
\end{definition}

\begin{figure}
\centering
\begin{minipage}[c]{.25\textwidth}
\centering
\includegraphics[scale=.6]{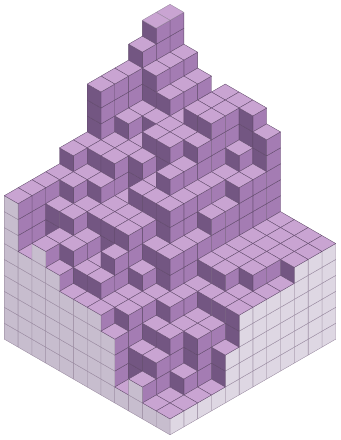}
\end{minipage}\hspace{.05\textwidth}%
\begin{minipage}[c]{.25\textwidth}
\centering
\includegraphics[scale=.65]{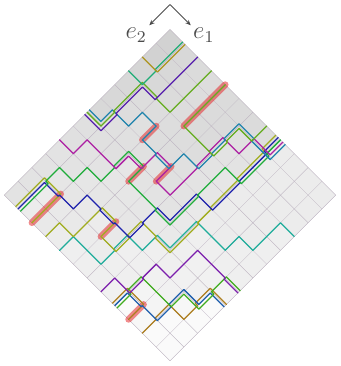}
\end{minipage}\hspace{.05\textwidth}%
\begin{minipage}[c]{.21\textwidth}
\centering
\includegraphics[scale=.65]{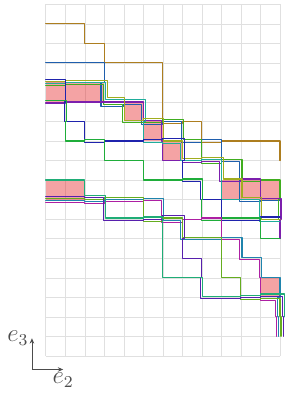}
\end{minipage}
\vspace{-0.2cm}
\caption{\HMSOS as a model of line ensembles: the surface (left);
noncrossing $\nearrow,\searrow,\swarrow$ paths with a penalty for $\swarrow$ (middle); and $\rightarrow,\downarrow$
paths with a penalty for the area above $\gamma_i$ and below
$\gamma_{i+1}$ (right).}
\label{fig:hmsos-line-ensembles}
\vspace{-0.1cm}
\end{figure}

\begin{remark*}
    The \HMSOS has two equivalent descriptions as a model of line ensembles (see \cref{fig:hmsos-line-ensembles}):
    \begin{enumerate}[(1)]
      \item \label{it:hmsos-ll-1}
In the $x_1 x_2$-plane: its lines (the boundaries of $\{x:\fh(x)>k\}$ for $k\in\Z$) are pairwise noncrossing lattice paths
(they may share edges) with steps $e_1, \pm e_2$. %
Up to an additive constant, the energy is twice the number of
$e_2$-steps per period cell (each $e_2$-step, with multiplicities, costs
a factor~$e^{-2\beta}$).
\item \label{it:hmsos-ll-2}
In the $x_2 x_3$-plane: the sections $\gamma_i$ of the interface by the planes $\{x_1=i\}$, i.e., the
boundaries of $\{(x_2,x_3):x_3<\fh(i,x_2)\}$, are lattice paths with steps
$e_2, -e_3$.
Unlike the paths in \cref{it:hmsos-ll-1}, these paths have no
backward steps but may cross, and the weight of a configuration is
$e^{-2\beta A}$, where $A$ is the number of unit squares lying above
$\gamma_i$ and below $\gamma_{i+1}$, summed over $i$.
    \end{enumerate}
    In both descriptions, the law without a height cutoff is taken
     modulo a global vertical shift.
     The model is trivial at slope $(\theta_1,0)$---reducing to a $(1+1)$\Dim \SOS (the surface is constant along every $x_2$-row)---but under boundary conditions with slope $\theta_1,\theta_2>0$ it is highly nontrivial.
\end{remark*}

In summary, the Ising Interface model has interactions $(\beta,\beta,\beta)$ in the three axes; \SOS has the interactions $(\beta,\beta,\infty)$; and \HMSOS has interactions $(\beta,\infty,\infty)$. The model that corresponds to $(\infty,\infty,\infty)$ (the limit where all interactions go to $\infty$) is zero-temperature 3\Dim Ising, i.e., monotone surfaces, which is in bijection with lozenge tilings (dimers on the honeycomb lattice).

\begin{definition}\label{def:height-function-in-models}
For \SOS and \HMSOS, the height function $\fh$ is the object at hand, as we defined above. For the Ising Interface model, define $\fh(x_1,x_2) = \sum_{x_3\geq 0} \one_{\{\sigma(x_1,x_2,x_3)=+1\}} - \sum_{x_3<0} \one_{\{\sigma(x_1,x_2,x_3)=-1\}}$, the column magnetization. For full 3\Dim Ising, define it via the interface obtained by erasing all bubbles (see \cref{def:ising-erasure}).
\end{definition}

\subsection{Previous work}

The study of low-temperature interfaces of the 3\Dim Ising model, and of the crystal surfaces they model, has a long history in both the physics and the mathematical physics literatures. In~1972, Dobrushin~\cite{Dobrushin72a} proved that under horizontal Dobrushin boundary conditions, and at low temperature, the plus/minus interface is rigid: it is localized at a bounded height, so that the corresponding non-translation-invariant Gibbs state exists in three dimensions. Predictions from the 1970's, which remain unproved, assert that the 3\Dim Ising interface undergoes a \emph{roughening transition} at an inverse-temperature $\beta_{\textsc r}$ strictly above the critical point $\beta_c$, from a rigid phase at $\beta>\beta_{\textsc r}$ to a rough phase at $\beta<\beta_{\textsc r}$ in which it is delocalized with logarithmic fluctuations (cf.~\cite{ChuiWeeks76,Abraham86,FPS81}; see~\cite{WGL73} for early numerical evidence of this transition and~\cite{ABKLLS80} for an experimental comparison of roughening in $^4\mathrm{He}$ crystals to 3\Dim Ising, as well as \cite{Weeks1980,NHB84,MWLB88,HasenbuschPinn1997} and the references therein, for a sample of such studies). The rigorous literature on these models has, apart from the celebrated work of Fr\"ohlich and Spencer~\cite{FrohlichSpencer81b}, been mostly confined to the rigid regime, proceeding from Dobrushin's work to delicate aspects such as entropic repulsion, layering and wetting; see, e.g., \cite{BEF86,KDM11,Lacoin18,Lacoin20,CLMST14,CLMST16,DinaburgMazel94,CesiMartinelli96a,CesiMartinelli96b} for a small sample, the more recent works~\cite{CKL_SOSlevelline,GheissariLubetzky21,GheissariLubetzky22,GheissariLubetzky23a,GheissariLubetzky23b,GheissariLubetzky_SOSLayering,CGL24} and the survey~\cite{IoffeVelenik18}.

When the boundary conditions are instead tilted (so the interface has a nonzero average slope), no such rigidity is expected, and the interface should be rough at every temperature, with a \GFF scaling limit. Giuliani~\cite{Giuliani23} singled this out as the last of the open problems in his ICM 2022 survey:
\begin{widequote}
\emph{Prove that the fluctuations of the interface of the 3D Ising model with tilted Dobrushin boundary conditions at low temperatures converges in the scaling limit to a GFF.}
\end{widequote}

The $(2+1)$\Dim \SOS model goes back to the early 1950's, to works of Burton, Cabrera and Frank~\cite{BCF51} and of Temperley~\cite{Temperley52} on crystal formation, and it approximates low-temperature 3\Dim Ising interfaces: viewing the height function $\fh$ as the surface having a horizontal face at height $\fh(x)$ above each $x$, completed by the minimal number of vertical faces, the \SOS weight is $e^{-\beta\sA(\fh)}$ where $\sA(\fh)$ is the number of faces of this surface, exactly the Ising weight of a bubble-free interface without overhangs. At zero slope, the roughening transition of \SOS is rigorous: Fr\"ohlich and Spencer~\cite{FrohlichSpencer81b} proved that $\var(\fh(x))\asymp\log|x|$ at small $\beta$, as in the \GFF, a Peierls argument (see~\cite{BrandenbergerWayne82,GMM73}) gives $\var(\fh(x))=O(1)$ at large $\beta$, and Lammers~\cite{Lammers22} proved the sharpness of this transition. That the rough phase should have a Gaussian scaling limit appears in many of these works; see, e.g., Fr\"ohlich, Pfister and Spencer~\cite{FPS81}, and the breakthrough of Bauerschmidt, Park and Rodriguez~\cite{BPR22a,BPR22b} proving the \GFF limit of the \ZGFF model (whose potential is $|\nabla \fh|^2$ instead of $|\nabla \fh|$) at small $\beta$.

The Peierls argument, however, only controls contractible level lines, and so it says nothing once the surface is placed on a nonzero slope $\theta=(\theta_1,\theta_2)$: the slope forces $(\theta_1+\theta_2)n$ macroscopic, non-contractible level lines, and it is these interacting lines that are expected to produce the rough phase at every $\beta$. The question of the scaling limit of such tilted interfaces was raised repeatedly (see, e.g., the work of Bodineau, Giacomin and Velenik~\cite{BGV01}). In his 2006 survey, Velenik~\cite{Velenik06} wrote that, for \SOS with any nonzero slope:
\begin{widequote}\emph{the large-scale behavior of these random interfaces is identical to that of their continuous counterpart. In particular they should have Gaussian asymptotics. This turns out to be quite delicate... I am not aware of a single rigorous proof for finite $\beta$.}\end{widequote}
Unlike 3\Dim Ising, the roughness (rather than limit) of tilted \SOS interfaces is rigorously known: Sheffield~\cite{Sheffield05} showed in 2005, via a remarkable application of the Ergodic Theorem, that no ergodic translation-invariant Gibbs measure exists for a wide class of integer-valued height function models (including all $|\nabla\phi|^p$ models) at an irrational slope, whence $\var(\fh(x))\to\infty$, without any bound on the rate. At high enough temperatures, Lammers and Ott~\cite{LammersOtt24} extended this to every slope, and Ott and Schweiger~\cite{OttSchweiger25} showed a logarithmic height divergence at the origin. However, at low temperatures, no improvement on Sheffield's result~\cite{Sheffield05} was known for the classical \SOS on a slope.

At zero temperature ($\beta=\infty$), the tilted 3\Dim Ising interface is a uniformly random monotone surface, i.e., a lozenge tiling, or equivalently a dimer configuration on the honeycomb lattice, which is exactly solvable through the determinantal formulas of Kasteleyn and of Temperley and Fisher. A landmark result of Kenyon~\cite{Kenyon00} %
is that the height function of the uniform dimer model converges to a \GFF, and Kenyon, Okounkov and Sheffield~\cite{KOS06} classified the ergodic translation-invariant Gibbs measures of the dimer model on any periodic bipartite planar graph by their slope, establishing \GFF fluctuations throughout the liquid phase. This captures the tilted 3\Dim Ising interface for every slope with $\theta_1,\theta_2>0$ at $\beta=\infty$.

Whether this \GFF behavior survives at finite $\beta$ was asked in several works, e.g., \cite{BGV01,GiulianiToninelli19} and the more recent work of Giuliani, Renzi and Toninelli~\cite{GRT23}, that pointed out: \begin{widequote}\emph{It is very likely that the \GFF behavior survives the presence of a small but positive temperature; however, the techniques underlying the proof at zero temperature, based on the exact solvability of the planar dimer problem, break down.}\end{widequote}
The main rigorous tool beyond exact solvability is the renormalization group (\RG) analysis of Giuliani, Mastropietro and Toninelli~\cite{GMT17,GMTH17,GMT20}, deriving the \GFF limit for \emph{interacting dimers}: perturbations of the uniform dimer measure by a local interaction whose strength is uniformly bounded in $L^\infty$ by a small constant.

In a prequel work~\cite{LaslierLubetzky24}, we established a \GFF scaling limit for tilted \SOS model in the presence of an additional potential $\mathsf V$, weak yet long-range and complicated. The proof there consisted of the following steps:
\begin{enumerate}[(i)]
\crefalias{enumi}{step}
    \item\label{it:prev-approx} given an \SOS height function $\fh$, approximate it by a random tiling $\varphi$, rewarded for its overlap with~$\fh$;
    \item\label{it:prev-interacting} show that the marginal law of the approximating $\varphi$ can be written as an interacting tiling model with small interactions; and as importantly, show that given $\varphi$, the conditional law of $\fh-\varphi$ is well-behaved;
    \item\label{it:prev-rg} use \RG analysis to derive the \GFF scaling limit for $\varphi$, then extend it from $\varphi$ to $\fh$.
\end{enumerate}
\Cref{it:prev-interacting} employed a delicate Markov chain temporal-to-spatial analysis to control the long range interactions of $\fh-\varphi$ given $\varphi$; \Cref{it:prev-rg} hinged on a novel extension of the \RG analysis of~\cite{GMT17,GMT20} covering interactions of unbounded range with exponentially decaying strength (see the follow-up work \cite{LaslierLubetzky26} developing one of the key ideas supporting that extension in a generic cluster expansion setup) as well as microcanonical ensembles.

The main shortcoming of \cite{LaslierLubetzky24} was its inability to treat microscopic regions wherein the \SOS height function~$\fh$ would be ``frozen,'' which are rare but nevertheless do occur. There, the approximating tiling $\varphi$ would no longer be weakly interacting, and the analysis of \cref{it:prev-interacting,it:prev-rg} would completely break down. This was the reason for needing the spurious additional potential $\mathsf V$ (which artificially suppresses such frozen regions in the surface~$\fh$), thus making the main theorem of~\cite{LaslierLubetzky24} not applicable to the actual \SOS model.

Overall, at low temperature, qualitative roughness of tilted
$(2+1)$\Dim \SOS (without a quantitative
rate) was known at irrational slopes; for tilted 3\Dim Ising, even qualitative roughness was not known.

\subsection{Main results}

Our main result is a \GFF limit for all models defined in \cref{def:3dIsing,def:interface-model,def:sos,def:HMSOS}.

\begin{maintheorem}[\GFF scaling limits]\label{thm:gff-limits}
Fix $\theta_1,\theta_2>0$ and consider the \HMSOS, \SOS, Ising Interface model or the 3\Dim Ising model on tilted tori with windings $\Theta_1,\Theta_2$ such that $\Theta_i/n\to\theta_i$ ($i=1,2$). There exists $\beta_0$ such that, for all $\beta\geq\beta_0$, the following hold for the height function $\fh$ of \cref{def:height-function-in-models}:
\begin{enumerate}[(i)]
\item \label{it:thm1-full-plane-law} The gradients of $\fh$ converge locally, as $n\to\infty$, to a full-plane law.
\item \label{it:thm1-gff} For the height observable $\fh$ in this limiting law, the centered rescaled fields $\fh_k$ of \cref{eq:test-field} converge as $k\to\infty$, modulo additive constants, to $\upsigma\,\GFF\circ L$ for some $\upsigma>0$ and some invertible linear map~$L$.
\end{enumerate}
\end{maintheorem}

The local convergence in \labelcref{it:thm1-full-plane-law} is
that of the gradient field: the height itself is only defined up
to a global vertical shift, and the limiting object is a
translation-invariant law on gradients, from which $\fh$ is
recovered modulo an additive constant. Convergence modulo additive
constants in \labelcref{it:thm1-gff} means that the pairings of $\fh_k$ with smooth
compactly supported test functions of integral zero converge (see
\cref{eq:test-field} for the precise rescaling, in which additive
constants cancel exactly). The amplitude $\upsigma$ and the linear
map $L$ depend on $\beta$ and on the slope: the map $L$ encodes the
anisotropy of the model, which is already present at $\beta=\infty$
through the slope dependence of the Kenyon--Okounkov--Sheffield
limit, and $\upsigma$ determines the strength of its fluctuations. Neither is explicit; they are
identified through the renormalized dispersion and source
coefficients in \cref{sec:rg}. %

\begin{remark}[The bulk spins in the full Ising case]\label{rem:ising-bulk-spins}
For the 3\Dim Ising model, \cref{thm:gff-limits} establishes the limit law of the height function of the interface $\cI$ obtained by erasing the bubbles from the full Ising configuration. In fact, we prove more: the erased interface $\cI$ and the original spins converge jointly in every fixed window around the root (see \cref{thm:ising-gff-complete}). Given the interface $\cI$, the spins are recovered by sampling the gas of bubbles around it from the law $\cK(\cI)$ (see \cref{prop:ising-contour-kernel}). Conversely, summing over these bubbles gives the marginal law of $\cI$, namely the bubble-free weight $e^{-\beta\sA(\cI)}$ multiplied by $e^{U_{\beta,\Gamma}(\cI)}$, where $U_{\beta,\Gamma}$ collects the contributions of clusters incompatible with the interface; this is the contact perturbation of \cref{thm:ising-rg-input}.
\end{remark}

A key ingredient is a renormalization theorem for interacting dimers
that allows large contributions on rare patterns. A \emph{pattern} $P$
is a finite connected patch of a lozenge tiling of the triangular lattice, and $|P|$ is its
number of lozenges. We consider a probability measure on tilings $\varphi$ given by a hard-core gas of patterns:
\begin{equation}\label{eq:pattern-gas}
 \P(\varphi)\propto F_\beta(\varphi)\qquad\mbox{where}\qquad F_\beta(\varphi)=
 \sum_{\substack{\{P_i\}:\ P_i\subset\varphi\\
                         P_i\text{ pairwise disjoint}}}
                         \prod_i c_\beta(P_i)>0\,,
\end{equation}
in which the sum in $F_\beta(\varphi)$ runs over all finite unordered collections, including the
empty collection with weight $1$. Thus the
interaction is specified by coefficients attached to patterns, which need not all be small.
The renormalization machinery of \cite{GMT17,GMT20} (extended to long-range interactions in \cite{LaslierLubetzky24}) established a \GFF limit for the dimer height under uniform weak-interaction assumptions, exemplified in the above notation by $|c_\beta(P)|\le e^{-A_\beta|P|}$ for every pattern $P$, with $A_\beta$ sufficiently large. The following theorem relaxes this $L^\infty$ control by allowing large coefficients on suitably restricted families of rare patterns.

\begin{maintheorem}[Informal; see \cref{thm:coefficient-rg}]
\label{thm:rg-typical}
There exist $A_\beta$ and $\epsilon_\beta$ with $A_\beta\to\infty$ and $\epsilon_\beta\to 0$ as~$\beta\to\infty$ so that the following holds. Suppose that the coefficients $c_\beta(P)$ in \cref{eq:pattern-gas} are translation invariant and admit
nonnegative \texttt{thin} and \texttt{fat} bounds
$q_\beta^{\mathrm t},q_\beta^{\mathrm f}$ such that:
\begin{enumerate}[leftmargin=1.5em, label=\arabic*.]
      \item For every pattern $P$,
      \[
 |c_\beta(P)|\leq q_\beta^{\mathrm t}(P)+q_\beta^{\mathrm f}(P),
 \qquad
 q_\beta^{\mathrm t}(P)\leq e^{-A_\beta|P|},
 \qquad
 q_\beta^{\mathrm f}(P)\leq e^{\epsilon_\beta|P|}\,.
\]
\item The $q_\beta^{\mathrm f}$ are supported on large
patterns $P$ whose area is almost entirely covered by large disks
contained in~$P$; in the infinite triangular lattice, at most
$e^{\epsilon_\beta m}$ such patterns of size $m$ contain any specified~triangle.
\end{enumerate}
Then, under appropriate bounds on pattern sizes and the
fraction of each pattern left uncovered by these disks, and suitable
consistency of the coefficients
between different volumes, for large enough
$\beta$ and every prescribed slope with all $3$ lozenge densities
positive, the measures in \cref{eq:pattern-gas} have a unique full-plane limit, and the centered height in this limit has
scaling limit $\upsigma\,\GFF\circ L$ for some $\upsigma>0$ and some
invertible linear map $L$.
\end{maintheorem}

Roughly, a \texttt{thin}
coefficient is small exponentially in the number of lozenges of
its pattern, as in the uniformly weak interactions of
\cite{GMT17,GMT20,LaslierLubetzky24}, and the expansion methods
of those works handle these. A \texttt{fat} coefficient may even grow with the area of
its support, at an exponential rate $\epsilon_\beta$ which is small
but positive. What compensates for this growth is a combination
of geometry and counting. The geometric condition in the theorem
says that a \texttt{fat} pattern is almost entirely a union of
large disks; the counting condition says that, up to translation,
there are at most $e^{\epsilon_\beta m}$ such patterns of size
$m$, far fewer than the exponentially many patterns of that size
in general. The point is that under the free dimer law, which has
positive entropy per lozenge in the liquid phase, a prescribed
pattern that fills a large disk is exponentially unlikely in its
area, at a rate which is bounded away from zero
(\cref{prop:fat_pattern}); the gain $e^{-c|P|}$ from this rarity
beats both the growth $e^{\epsilon_\beta|P|}$ of the coefficient
and the entropy $e^{\epsilon_\beta|P|}$ of the family, once
$\epsilon_\beta$ is small. The theorem thus replaces a uniform
smallness assumption by a typical one: the interaction is small
on the patterns that the free measure actually produces. The
full theorem also gives exponential tails for the
geometric witnesses used to reconstruct the original surfaces, under
an additional bound on overlapping collections; these tails are
proved under the interacting law (not only under the free one),
and this is what our application needs.

\begin{remark}\label{rem:individual-coefficient-bounds}
On first glance, \cref{thm:rg-typical} appears almost too good to be true: it allows coefficients that are exponentially \emph{growing} in the pattern size, on an exponentially growing family of patterns, and yet it still gives a \GFF limit. Morally, what makes the interaction small is that the patterns are rare under the uniform dimer law $\mu_0$, so the natural weight of a pattern is $|c_\beta(P)|\mu_0(P\subset\varphi)$. The reason for the covering criterion in the definition of a \texttt{fat} set is nonobvious but crucial in order to introduce a factor similar to $\mu_0(P\subset \varphi)$ into the \RG analysis; we give more details below on this criterion and its role in the proof.
\end{remark}

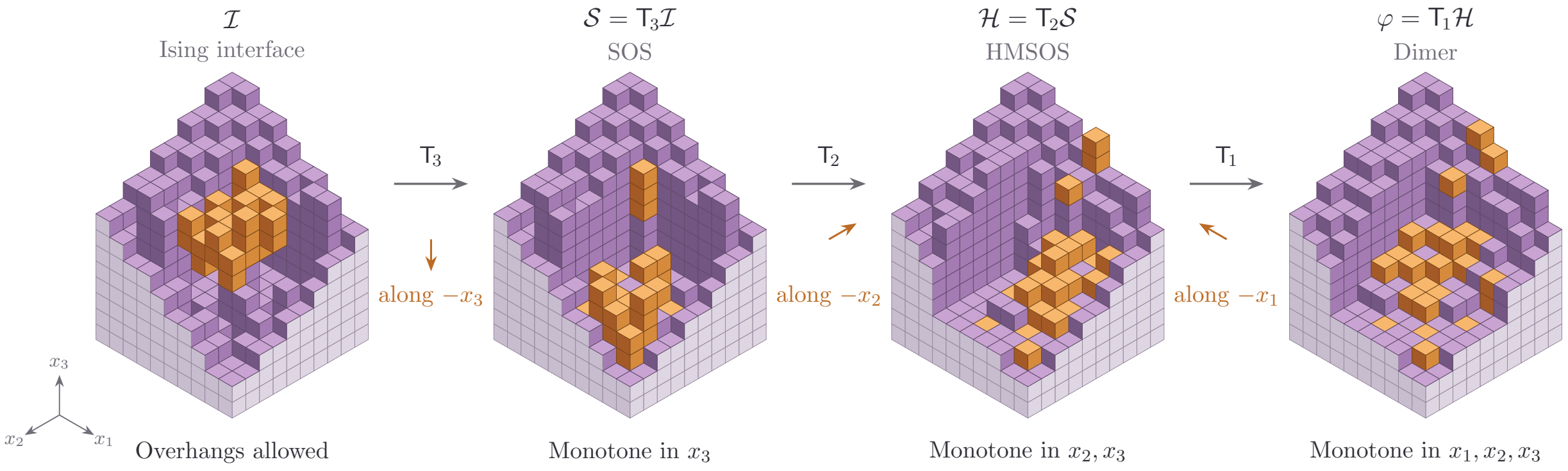
\begin{figure}
      \vspace{-0.05cm}
\centering
\begingroup
\definecolor{ink}{HTML}{302D36}
\definecolor{muted}{HTML}{716C78}
\definecolor{accent}{HTML}{BC6B25}
\resizebox{\textwidth}{!}{%
\begin{tikzpicture}[x=1cm,y=1cm,text=ink,>=Stealth]
\path[use as bounding box] (-.85,.8) rectangle (25.1,8.55);

\node[font=\Large] at (2.800,8.17) {$\mathcal I$};
\node[text=muted,font=\large] at (2.800,7.65) {Ising interface};
\node[inner sep=0pt,outer sep=0pt] at (2.80000,4.50500) {\includegraphics{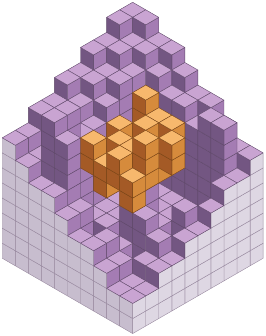}};
\node[font=\large] at (2.800,1.15) {Overhangs allowed};
\node[font=\Large] at (9.250,8.17) {$\mathcal S=\Tthree\mathcal I$};
\node[text=muted,font=\large] at (9.250,7.65) {SOS};
\node[inner sep=0pt,outer sep=0pt] at (9.25000,4.50500) {\includegraphics{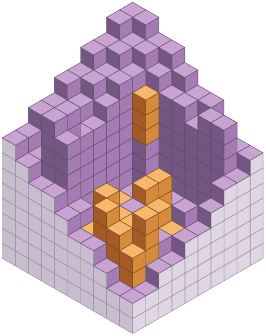}};
\node[font=\large] at (9.250,1.15) {Monotone in $x_3$};
\node[font=\Large] at (15.700,8.17) {$\mathcal H=\Ttwo\mathcal S$};
\node[text=muted,font=\large] at (15.700,7.65) {HMSOS};
\node[inner sep=0pt,outer sep=0pt] at (15.70000,4.50500) {\includegraphics{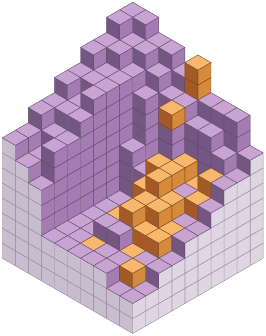}};
\node[font=\large] at (15.700,1.15) {Monotone in $x_2,x_3$};
\node[font=\Large] at (22.150,8.17) {$\varphi=\Tone\mathcal H$};
\node[text=muted,font=\large] at (22.150,7.65) {Dimer};
\node[inner sep=0pt,outer sep=0pt] at (22.15000,4.50500) {\includegraphics{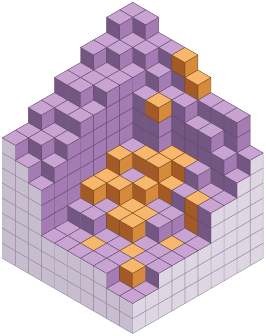}};
\node[font=\large] at (22.150,1.15) {Monotone in $x_1,x_2,x_3$};
\draw[->,line width=1pt,ink!70] (5.425,5.5) -- (6.625,5.5);
\node[font=\large,fill=white,inner sep=2pt] at (6.025,5.95) {$\Tthree$};
\draw[->,accent,line width=.9pt] (6.025,4.600) -- (6.025,4.080);
\node[text=accent,font=\large] at (6.025,3.68) {along $-x_3$};
\draw[->,line width=1pt,ink!70] (11.875,5.5) -- (13.075,5.5);
\node[font=\large,fill=white,inner sep=2pt] at (12.475,5.95) {$\Ttwo$};
\draw[->,accent,line width=.9pt] (12.475,4.600) -- (12.925,4.860);
\node[text=accent,font=\large] at (12.475,3.68) {along $-x_2$};
\draw[->,line width=1pt,ink!70] (18.325,5.5) -- (19.525,5.5);
\node[font=\large,fill=white,inner sep=2pt] at (18.925,5.95) {$\Tone$};
\draw[->,accent,line width=.9pt] (18.925,4.600) -- (18.475,4.860);
\node[text=accent,font=\large] at (18.925,3.68) {along $-x_1$};
\begin{scope}[yshift=-0.2cm]
\draw[->,muted,line width=.55pt] (0.000,1.950) -- (0.563,1.625);
\node[font=\small,text=muted] at (0.727,1.530) {$x_1$};
\draw[->,muted,line width=.55pt] (0.000,1.950) -- (-0.563,1.625);
\node[font=\small,text=muted] at (-0.727,1.530) {$x_2$};
\draw[->,muted,line width=.55pt] (0.000,1.950) -- (0.000,2.600);
\node[font=\small,text=muted] at (0.000,2.790) {$x_3$};
\end{scope}
\end{tikzpicture}
}
\endgroup
\vspace{-0.2in}
\caption{The three \Tetris moves take an Ising interface $\cI$
to an \SOS surface $\cS$, an \HMSOS surface $\cH$, and a
monotone surface $\varphi$, equivalently a dimer configuration.}
\label{fig:tetris-three-moves}
\vspace{-0.15in}
\end{figure}

\subsection{Proof ideas}

Our proof strategy originated from that of \cite{LaslierLubetzky24} recalled above, but with key differences:
\begin{enumerate}[(i)]
\item The random approximation of $\fh$ by a monotone surface $\varphi$ is replaced by a \emph{deterministic} function, which we call a \Tetris map: it lets cubes drop according to gravity along each of the three axes in turn, taking a 3\Dim Ising interface to an \SOS surface, then to an \HMSOS surface, and finally to a monotone surface (see~\cref{fig:tetris-three-moves} for an illustration of the three moves).

\item Each move of \Tetris sorts the occupied cubes along the lines parallel to one axis, so it preserves the column magnetization on each such line and only rearranges the cubes within it; this is what allows its fibers to be described exactly. The combinatorial control of the allowed differences between the base model (either \HMSOS, \SOS, or the Ising Interface model) and $\varphi$ is carried out in what we call a  \emph{dot-and-arrow} representation, in which it becomes a local enumeration problem: Dots record vertical spin disagreements, arrows record horizontal ones (each costs a factor $e^{-\beta}$), and the enumeration is a stars-and-bars counting of interleaving the two between adjacent squares (see \cref{fig:dot-arrow}).

The Markov chain analysis from the prequel is no longer needed, because the conditional law of $\fh-\varphi$ given $\varphi$ is much simpler in the deterministic setting: it is an explicit hard-core gas of independent local defects rather than a measure with long-range interactions, so its locality can be proved by a static two-copy argument rather than through the mixing time analysis of an MCMC sampler.
\item The \RG part must be improved again in order to handle interactions that are only typically small: In~\cite{LaslierLubetzky24}, the artificial added potential guaranteed that all interaction coefficients were small. Here, the new \cref{thm:rg-typical} handles rare patterns with large coefficients, that are unavoidable in the true model.
\end{enumerate}

\emph{The transition to a deterministic \Tetris map.}
The initial approach for the \SOS model in the prequel~\cite{LaslierLubetzky24} was to take the monotone surface $\varphi$ that is \emph{closest} (maximizing its overlap in faces) to the \SOS function~$\fh$, show that it is a weakly interacting dimer model, hence has a \GFF limit, and by extension, so does $\fh$ (\cref{it:prev-interacting,it:prev-rg} of that paper's program). The obvious hard core combinatorial constraints associated with such an optimizer $\varphi$ made this approach seem intractable, leading to the choice of a random $\varphi$. Evidently, this created long-range interactions that were difficult to control, and in hindsight did not remove the hard core constraints completely, because we had to work with families of defects that would make the optimization problem local. In the present work, we therefore move back to our original idea of a deterministic approximation, namely to the \Tetris map (\cref{sec:tetris}). Rather than a solution to a combinatorial optimization problem, it is an explicit operation on the interfaces, which is tractable with the help of the dot-and-arrow diagrams described next.

\emph{Dot-and-arrow diagrams.}
Consider first a single move of \Tetris, say from an Ising interface $\cI$ to an \SOS surface $\cS=\Tthree\cI$.
The surface $\cS$ is described by its local gradients, which we can view as a collection of \emph{arrows} between adjacent squares (each marking a gradient of $1$, or equivalently, a single disagreement in the associated spin configuration $\sigma$). Similarly, $\cI$ is described by such an arrow configuration, but now we must also add a set of \emph{dots} to each square, recording the number of horizontal faces above it (again, each corresponds to a single disagreement in $\sigma$). This turns out to be a complete description of the interface once we attach to it a consistent weak order of the dots-and-arrows at each corner of four incident squares.

Enumerating these weak orders reduces to a stars-and-bars-type factor (dots are stars, arrows are bars, except that the orders are weak and one accounts for where dots and arrows are placed around a vertex). Hence,
given that each extra dot and each extra arrow cost in energy, it is disadvantageous to have many dots and arrows. In particular, we can be quite precise in the enumeration of preimages.

Note that, if $\cS$ has a very large gradient between two squares, the number of preimages is large (dots can be placed in many ways around the arrows), and this is the mechanism that creates large coefficients in the pattern-gas representation of the interaction.
This is the same type of rare but problematic (``frozen'') region that foiled the analysis in the prequel paper: the \Tetris map did not remove this obstacle; however, the dot-and-arrow analysis is still valid in its presence (unlike the Markov chain analysis of the prequel), and allows us to pass it to the \RG analysis (where it is treated as a rare pattern with a large coefficient).

One technical difficulty is that we need to understand the superposition of the three \Tetris moves, which is nontrivial. The different notions of \emph{footprints} in \cref{sec:preimage-witnesses} allow us to group defects that interact through successive \Tetris moves, in a way similar to the way the prequel paper grouped defects to make the optimization problem of approximating $\fh$ by  $\varphi$ local (footprints in this paper are analogous to bubble groups in the prequel).

\begin{figure}
\vspace{-0.25cm}
\centering
\raisebox{-.5\height}{\includegraphics[height=.3\textwidth]{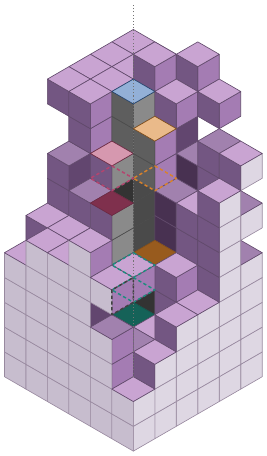}}\hspace{.09\textwidth}
\raisebox{-.5\height}{\includegraphics[height=.3\textwidth]{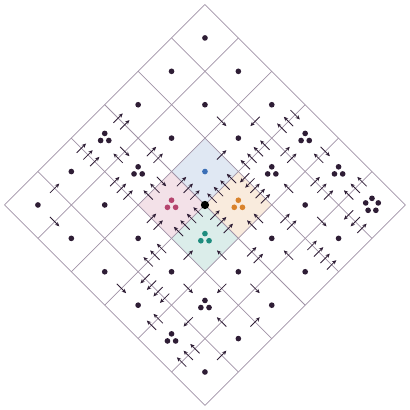}}\hspace{.09\textwidth}
\raisebox{-.5\height}{\includegraphics[height=.3\textwidth]{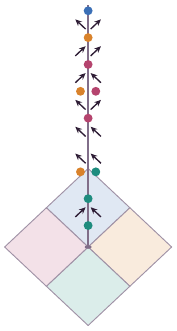}}
\vspace{-0.2cm}
\caption{An Ising interface $\cI$, where the highlighted faces are ones meeting the dotted vertical line (left); its dot-and-arrow diagram (middle); and the decoration at the marked corner (right).}
\label{fig:intro-dot-arrow}
\vspace{-0.4cm}
\end{figure}

\emph{The pattern-gas form of the interaction.}
The weight in \labelcref{eq:pattern-gas} is a hard-core gas of
patterns, whereas an interacting dimer model is more customarily
described as a tilt of the uniform measure by an exponential of a
sum of local functions,
$\exp(\sum_x\sum_r\fg_r(\varphi\restriction_{B(x,r)}))$, as in
\cite{GMT17,GMT20,LaslierLubetzky24}. For weak interactions the two
are interchangeable: writing $e^{\fg_r}=1+(e^{\fg_r}-1)$ and merging
overlapping patterns turns the exponential into a hard-core gas of
connected patterns, and the logarithm of a gas with small activities
is a convergent cluster expansion. Typically, such transformations entail
losses in the estimates. Here the gas form arises naturally by
grouping the local defects of each preimage, and our Grassmann
construction starts directly from this form (\cref{sec:model-grassmann}). Both in the prequel~\cite{LaslierLubetzky24} and in this paper, moving from the defect gas to the exponential form, and then moving back, are both nontrivial. Here, the \texttt{fat} coefficients do not justify a uniformly convergent cluster expansion for $\log F_\beta$; in the prequel, the initial defect representation had long-range
many-body interactions, which also prevented the standard cluster expansion. (Note that, in that work, the long-range interactions were localized in the exponential form by the proof technique, so that was the natural choice for the \RG input.) Here, going from hard-core to exponential and then back would incur an unsustainable loss.

\emph{The obstruction to uniformly small interactions.}
As we mentioned above, a key difficulty is that the coefficients $c_\beta(P)$ are not uniformly small in $|P|$. 
The essential
geometric observation is that these large coefficients are associated to
long straight rows of lozenges. 

Crucially, a straight row of $m$ lozenges has probability at most $e^{-C m}$ under the free dimer law $\mu_0$. Analyzing the Kasteleyn determinant shows that the probability of this event is at most $e^{-c m}$ for a small $c>0$, but that is suboptimal. Looking at the constraints imposed by such a straight row on its surroundings, one can obtain a bound as above with $C>0$ arbitrarily large: 
Thickening a row produces a tube $R$ in
which one lozenge orientation has very small density, by
\cref{lem:straight-tube}; there are at most $e^{\epsilon |R|}$ such tilings of $R$. 
Any such tiling freezes many interior degrees of
freedom, and its probability under $\mu_0$ is exponentially small in its area. Overall, this gives an upper bound of $e^{-C m}$ on the probability of the original row.

\emph{Using the rarity of a pattern in the first integration.}
The proof of
\cref{prop:fat_pattern} expresses the probability bound $e^{-Cm}$ from above through
Kasteleyn determinants. Recall that the probability of a
prescribed set of edges under the free dimer law is a minor of the
inverse Kasteleyn matrix, so the statement to prove is that this minor is exponentially small in its size when the support of the pattern is almost entirely covered by large disks contained in it. The argument proceeds by filling most of the pattern support with
interior patches having large partition functions, cancelling their
interior rows against the inverse Kasteleyn matrix, so that each
patch contributes the reciprocal of its partition function, and bounding the
small remaining minor, which involves only the boundary of the
patches and the uncovered fraction of the pattern. The same mechanism is then built into the
renormalization: the first Grassmann integration integrates out the
modes of the field away from the two singularities of the free
propagator (the ``fast'' field), and a corresponding determinant estimate in
this integration compensates both the possible
growth of a \texttt{fat} coefficient and the number of such
patterns, exactly as the probability estimate would for a single
pattern under the free law. 

After the first integration, we are essentially able to use the machinery of \cite{GMT17,GMT20} as a black box:
the resulting effective kernels, including kernels with
edge sources, are uniformly small in the norms needed by the infrared
analysis (see \cref{thm:bulk-first-step,thm:torus-first-step} for more details).

\emph{From the infrared flow to the prescribed slope.}
As in the prequel \cite{LaslierLubetzky24}, the models are naturally defined on a microcanonical ensemble with a prescribed slope, while the \RG works in a canonical ensemble. We perform Fourier inversion to extract the desired slope from the canonical law: the Kasteleyn weights are tuned to the desired slope, and we prove that the
target microcanonical slope has probability bounded below under the canonical law
(\cref{prop:interacting-central-sector}), so that conditioning on
it changes probabilities by at most a bounded factor. 
This requires  technical estimates such as verifying that different bounds are uniform in certain neighborhoods of the target slopes, which is part of the reason that \cref{sec:rg} is fairly technical.

\emph{Returning from dimers to surfaces.}
Once the tiling $\varphi$ has a \GFF limit, it remains to show that
the original surface, which is obtained from $\varphi$ by sampling
the fiber, has the same limit. The fiber is a hard-core gas of
defects, and the defects are small except on the rare witnesses
described above, so the question is how often large witnesses
occur under the interacting tiling law. The subtlety is that the
reconstruction needs to control regions with strong interactions under the interacting
tiling law, not just under the free law used in the first integration. This is an extra output required from the \RG analysis, and another part of the reason that \cref{sec:rg} is so technical.

\emph{From the Ising Interface model to full Ising.}
Finally, for full Ising, a low-temperature contour expansion
integrates out closed bubbles and produces a local contact
perturbation of the bubble-free interface law. We allow such a perturbation as a generalization of the Ising Interface model.  One notable technicality is that the resulting pattern expansion can now have signed coefficients (unlike the models from \cref{def:interface-model,def:sos,def:HMSOS}, where the coefficients are all nonnegative). Finally, the reconstruction of the original interface $\cI$ extends to the full Ising configuration, readily providing the law of the bubbles conditional on the interface $\cI$.

\subsection{Open problems}\label{sec:open-problems}

A natural open problem is to extend the 3\Dim Ising conclusion of
\cref{thm:gff-limits} to every $\beta>\beta_c$, where $\beta_c$ is
the bulk critical inverse temperature, for each fixed slope
$\theta_1,\theta_2>0$. Another is to determine the fluctuation scale
and scaling limit of \SOS and 3\Dim Ising interfaces at degenerate slopes
$(\theta_1,0)$ with $\theta_1>0$, even at low temperature; see the
discussion in \cite[\S1.3]{LaslierLubetzky24}. The ground states have
one-dimensional fluctuations of order $\sqrt n$, whereas even the
order of the positive-temperature height variance remains open.
Related results include the logarithmic variance upper bound of
Bodineau, Giacomin and Velenik~\cite{BGV01} for a real-valued
$(2+1)$\Dim \SOS/Gaussian hybrid, and the recent
rigidity result of Adhikari, Gheissari and Peled~\cite{AGP26}
for $\Z$-valued height function models in $(3+1)$\textsc{d}, including \SOS, at low-temperature and slopes $(\theta_1,0,0)$.

\subsection{Organization}

\Cref{sec:tetris} develops the sorting map, its image law and the
common one-step enumeration and bounds on the fiber partition function.
\Cref{sec:preimage-witnesses} analyzes the fibers of the one-, two-
and three-step histories and constructs their geometric witnesses.
\Cref{sec:model-grassmann} collects their pattern and Grassmann
inputs to \cref{thm:coefficient-rg}.
\Cref{sec:rg} proves the \RG theorem, the precise
version of \cref{thm:rg-typical}.
Finally, \cref{sec:transfer} uses these results to reconstruct
the original models and transfer the tiling limits to their physical
height fields.

\section{The Tetris approximation}\label{sec:tetris}

This section defines the \Tetris map and its image law and proves
the common one-step enumeration and bounds on the fiber partition function.
\Cref{sec:preimage-witnesses} applies these estimates to the three
histories and constructs their geometric witnesses, retaining the
exact constraints for successive maps.

Throughout, regard a cube $x$ as \emph{occupied} if $\sigma(x)=+1$ and
\emph{vacant} if $\sigma(x)=-1$. A \emph{fiber} is a line of cubes parallel
to a coordinate axis, an \emph{$i$-fiber} when the axis is $e_i$, and a
\emph{column} when it is vertical, i.e., of the form $\{(x_1,x_2,\cdot)\}$.
It will sometimes be convenient to refer
to the \emph{$\{0,1\}$-spin configuration} $(1+\sigma)/2$, the indicator
of the occupied cubes, and to its restriction to a fiber as the
$\{0,1\}$-spin sequence of that fiber.
Two unit cubes of $\Z^3$, or two unit squares of $\Z^2$, are
\emph{adjacent} if their centers are at distance $1$ (degree $6$,
respectively $4$) and \emph{$*$-adjacent} if their centers are at
distance at most $\sqrt2$, that is, if they share a face or an edge in
$\Z^3$ (degree $18$), or an edge or a corner in $\Z^2$ (degree $8$); the
former is the convention of \cref{def:interface-model}. We write
$\dist_*$ for the graph distance of $*$-adjacency, and call a set
\emph{$*$-connected} if it is connected under $*$-adjacency.

\subsection{The Tetris map}

By the boundary conditions of \cref{def:3dIsing}, the occupied cubes are the ones below the interface, and with its coordinate conventions they form, at zero
temperature, a decreasing set, i.e., a set which contains,
together with any cube $(x_1,x_2,x_3)$, also the three cubes $(x_1,x_2,x_3)-e_i$ ($i=1,2,3$).
This is the reason for our choice of orientation for the axes: below a graph which rises
along $+x_1$ and $+x_2$, say $\{x_3<\theta_1x_1+\theta_2x_2\}$, the
occupied set grows with $x_1$ and $x_2$, so a fiber parallel to
$e_1$ or $e_2$ is occupied sufficiently far in the positive direction
and vacant sufficiently far in the negative one, and sorting such a
fiber into a decreasing one would move infinitely many cubes. Below
the plane $x_3=-\theta_1x_1-\theta_2x_2$, every fiber is occupied
far in the negative direction and vacant far in the positive one,
which is the end behavior assumed in \cref{def:tetris-charge}
below.

Write $\cP_{111}=\{x\in\R^3:x_1+x_2+x_3=0\}$ for the
diagonal plane and let
\[
 \Upsilon:\R^3\longrightarrow\cP_{111},\qquad
 \Upsilon x=x-\frac{x_1+x_2+x_3}{3}(1,1,1)
\]
be orthogonal projection. It sends the unit plaquettes of a
monotone surface to a lozenge tiling. When using lattice
coordinates on this plane, we identify $\Upsilon x$ with
$t=(x_1-x_3,x_2-x_3)$. Thus $\Upsilon\Z^3$ is identified
with $\Z^2$, and the projections of the three coordinate
directions have coordinates $(1,0),(0,1),(-1,-1)$.

\begin{definition}[\Tetris map]\label{def:tetris-charge}
For a configuration $\omega:\Z+\frac12\to\{0,1\}$ equal to $1$ sufficiently far to
the left and $0$ sufficiently far to the right, define the \emph{magnetization}, or \emph{charge}, of $\omega$, to be
\[
 q(\omega):=\sum_{k\in\Z+\frac12}(\omega(k)-\one_{\{k<0\}})=\sum_{k\geq0}\one_{\{\omega(k)=1\}}
           -\sum_{k<0}\one_{\{\omega(k)=0\}} \in \Z\,.
\]
For a vertical column of an Ising Interface configuration $\sigma$,
if $\omega=\one_{\{\sigma(x_1,x_2,\cdot)=+1\}}$, then $q(\omega)=\fh(x_1,x_2)$ is precisely the column magnetization
of \cref{def:height-function-in-models}. With this convention,
a graph column of height $h\in\Z$ has $q(\omega)=h$: for
$\omega=\one_{\{\cdot<h\}}$ the sum defining $q(\omega)$ has exactly $|h|$
nonzero terms, each equal to $\sgn(h)$. Define the sorted configuration by
\[
 \omega^\downarrow(k)=\one_{\{k<q(\omega)\}}\,.
\]
It is the unique decreasing configuration with the same ends
and the same charge as $\omega$; that is, it is obtained by letting all
occupied cubes of the fiber drop as far as possible.
For a configuration whose coordinate fibers have these ends, $\Ti$ replaces
each fiber parallel to $e_i$ by its decreasing rearrangement. Equivalently,
for $i\in\{1,2,3\}$ and $x\in(\frac12,\frac12,\frac12)+\Z^3$, let $\omega$ be
the $i$-fiber through $x$, $\omega(k)=\one_{\{\sigma(x+(k-x_i)e_i)=+1\}}$ for
$k\in\Z+\frac12$; then
\begin{equation}\label{eq:tetris-map}
 (\Ti\sigma)(x)=2\,\omega^\downarrow(x_i)-1=\sgn\big(q(\omega)-x_i\big)\,,
\end{equation}
where the last expression is never zero since $q(\omega)\in\Z$ and
$x_i\in\Z+\frac12$.
\end{definition}

\begin{proposition}[Sorting]\label{prop:sorting}
The \Tetris map $\Ti$ preserves the charge of every $i$-fiber, does not increase
the number of disagreeing bonds, and preserves monotonicity in each
coordinate other than $i$ in which the input is already monotone.
Consequently $\Tthree$, $\Ttwo\Tthree$, and $\Tone\Ttwo\Tthree$ produce, respectively, an
\SOS, an \HMSOS, and a monotone surface. The same statements hold per
period cell whenever the configuration has the tilted periodicity and
finite energy per cell.
\end{proposition}
\begin{proof}
The charge of a fiber is preserved by definition: since
$\omega^\downarrow=\one_{\{\cdot<q(\omega)\}}$,
\[
 q(\omega^\downarrow)=\sum_{k\in\Z+\frac12}\big(\one_{\{k<q(\omega)\}}-\one_{\{k<0\}}\big)
 =q(\omega)\,.
\]

The disagreeing bonds split into two classes: bonds parallel to $e_i$,
which join two consecutive cubes of the same fiber, and transverse
bonds, which join cubes at the same position in two neighboring
parallel fibers. For the first class, the number of disagreements
along a fiber with the prescribed ends is a positive odd integer,
since the fiber starts occupied and ends vacant, and after sorting it
is exactly one. For the second class, consider two neighboring fibers
$\omega,\omega'$ with the prescribed ends. The sorted fibers are two step
sequences, so they disagree exactly on the $|q(\omega)-q(\omega')|$ positions
between their two steps, and
\[
 \sum_k|\omega^\downarrow(k)-(\omega')^\downarrow(k)|
 =|q(\omega)-q(\omega')|
 =\bigg|\sum_k(\omega(k)-\omega'(k))\bigg|
 \leq\sum_k|\omega(k)-\omega'(k)|\,,
\]
where the middle identity is the definition of the charges, whose
flat-step normalizations cancel. Summing the two classes over all
fibers and all pairs of neighboring fibers proves that the number of
disagreeing bonds does not increase.

For the monotonicity assertion, observe that if $\omega\leq\omega'$ pointwise,
then $q(\omega)\leq q(\omega')$ by the definition of the charge, and hence
$\omega^\downarrow\leq (\omega')^\downarrow$, since a step sequence with a higher
step dominates one with a lower step. Monotonicity in a coordinate
$j\neq i$ means that for every $j$-bond the spin does not
increase in the positive $e_j$ direction; equivalently, consecutive
$i$-fibers along $e_j$ are pointwise ordered. Applying the above to
each such pair of consecutive parallel fibers proves preservation of
all previously established coordinate monotonicities. The three
successive images are obtained by applying this to $\Tthree$ on an
arbitrary input (producing an \SOS surface, i.e., monotone in
$x_3$), then to $\Ttwo$ on this \SOS input (preserving
monotonicity in $x_3$ while adding it in $x_2$),
and finally to $\Tone$.

It remains to prove the statements for a configuration with tilted periodicity (as per
\cref{def:3dIsing}) and finite energy per cell. Since
$\Theta_1,\Theta_2>0$, no nonzero period is parallel to a coordinate
axis: a period $k(n,0,-\Theta_1)+l(0,n,-\Theta_2)$ is parallel to
$e_1$ only if $l=0$ and $k\Theta_1=0$, hence only if $k=l=0$; the
same holds for $e_2$; and parallelism to $e_3$ forces $kn=ln=0$.
Consequently, a nonzero period maps every $i$-fiber to a different
$i$-fiber, and the $i$-fibers fall into finitely many orbits. The
prescribed ends follow from the tilted boundary conditions: along
each coordinate line, $x_3+(\Theta_1x_1+\Theta_2x_2)/n$ tends to
$-\infty$ and $+\infty$ at its two ends, where the spins are fixed
to $+1$ and $-1$, respectively. Thus the sums over $k$ above are
finite on each fiber. Sorting commutes with every lattice
translation: a translation carrying one fiber to another acts on
the sorted configuration in the same way, while a translation of a
fiber along itself by $\ell$ changes its charge by $\ell$,
$q(\omega(\cdot-\ell))=q(\omega)+\ell$, so that the sorted translate
is the translate of the sorted fiber. Hence, $\Ti$ maps a periodic
configuration to a periodic configuration with the same periods.
Applying the inequalities above to a set of representatives of the
fiber orbits and of the orbits of pairs of neighboring fibers, and
summing, counts every bond orbit exactly once; the number of bond
orbits in the disagreement set is the energy per period cell, which
therefore does not increase. The monotonicity assertions are
pointwise and need no modification. The same argument applies to
each of the three successive maps, and in a finite box it applies
verbatim with the occupied and vacant ends fixed outside the box.
\end{proof}

\begin{definition}[Image law and excess energy]\label{def:image-law}
Let $f$, in the context of this definition, denote $\Tone$ for an \HMSOS input monotone in $x_2,x_3$,
$\Tone\Ttwo$ for an \SOS input monotone in $x_3$, and $\Tone\Ttwo\Tthree$
for an interface input. Heights are rooted, or configurations are
identified modulo common vertical translation. On a fixed finite
period cell, recall that $\sA(\cI)$ denotes the plaquette area of $\cI$, namely its
energy; let $\sA_0$ be the area of a monotone surface in the prescribed
slope sector, and set $E(\cI)=\sA(\cI)-\sA_0$. The image law is
\[
 Z_\beta(\varphi)=\sum_{\cI\in f^{-1}(\varphi)}e^{-\beta E(\cI)},\qquad
 \mu_\beta(\varphi)=\frac{Z_\beta(\varphi)}{\sum_\psi Z_\beta(\psi)}\,.
\]
The area $\sA_0$ is constant in the sector because the number of plaquettes
of each orientation is fixed for a monotone surface: a monotone
surface crosses each fiber of each coordinate direction exactly once,
and the tilted periodicity fixes the number of fibers of each
direction per period cell. By \cref{prop:sorting}, every preimage
satisfies $E(\cI)\geq0$, since
$\sA(\cI)\geq\sA(\Tthree\cI)\geq\sA(\Ttwo\Tthree\cI)\geq\sA(\varphi)=\sA_0$
(with the obvious shorter chains for the \SOS and \HMSOS inputs).
Equality in this chain forces equality in each of the sorting
energy inequalities, and equality in the inequality for $\Ti$ forces
every $i$-fiber of the input to carry exactly one disagreement, i.e.,
to be already sorted, because the number of disagreements along a
fiber is odd and only decreases to one. Hence each successive
map acts as the identity, and the unique zero-excess preimage of
$\varphi$ is $\varphi$ itself. The weight $e^{-\beta E(\cI)}$ is
thus the original Boltzmann weight normalized by the ground-state
weight of the sector, and $Z_\beta(\varphi)\geq1$ with the term $1$
accounting for $\cI=\varphi$.
\end{definition}

\begin{remark*}[Local interaction input]
We will show (see
\cref{thm:hmsos-rg-input,thm:sos-rg-input,thm:interface-rg-input,thm:ising-rg-input})
that, for large enough $\beta$, the \Tetris image weights for \HMSOS, \SOS and the Ising Interface model admit exact expansions
in connected local tiling patterns. Their absolute coefficients
split into exponentially small \texttt{thin} terms and
\texttt{fat} terms supported on large patterns of vanishing rooted
entropy. These expansions satisfy the source, torus and
marked-collection conditions needed for the \RG theorem.
For full Ising, the same assertion applies to the normalized
single-surface comparison weights obtained after integrating out
closed contours; the erased interface law is compared with this
law in total variation as the width tends to infinity.
\end{remark*}
The precise common conditions are in \cref{def:rg-admissible},
and their model-specific verification is the content of the four
theorems cited above. This is the input to
\cref{thm:coefficient-rg}: the exact pattern expansion has bounds
on absolute internal weights and on collections overlapping a witness.

The enumeration and partition function estimates in
\cref{sec:enumeration} underlie
the local interaction input just described. In \cref{sec:preimage-witnesses}, complete
\Tetris histories are grouped into local patterns, and large footprints
are controlled by geometric witnesses. The resulting coefficient
bounds, together with their source and torus compatibility in
\cref{sec:model-grassmann}, supply the \RG input.

\subsection{One-step preimages: the dot-and-arrow enumeration}\label{sec:enumeration}

Since the \Tetris map is applied up to three times in passing from an Ising interface to a monotone surface, we analyze the set of preimages in the general setting of an arbitrary admissible starting surface.
Throughout this subsection the sorting direction is $x_3$, so the
map under consideration is $\Tthree$; the other two
moves are treated by the same analysis after relabeling the
coordinates, with the additional monotonicity constraints of their
inputs imposed in \cref{sec:preimage-witnesses}.
The reference surface is thus an \SOS height function $h$, and a
preimage is an interface, possibly with overhangs, all of whose
vertical columns have the charges prescribed by $h$.

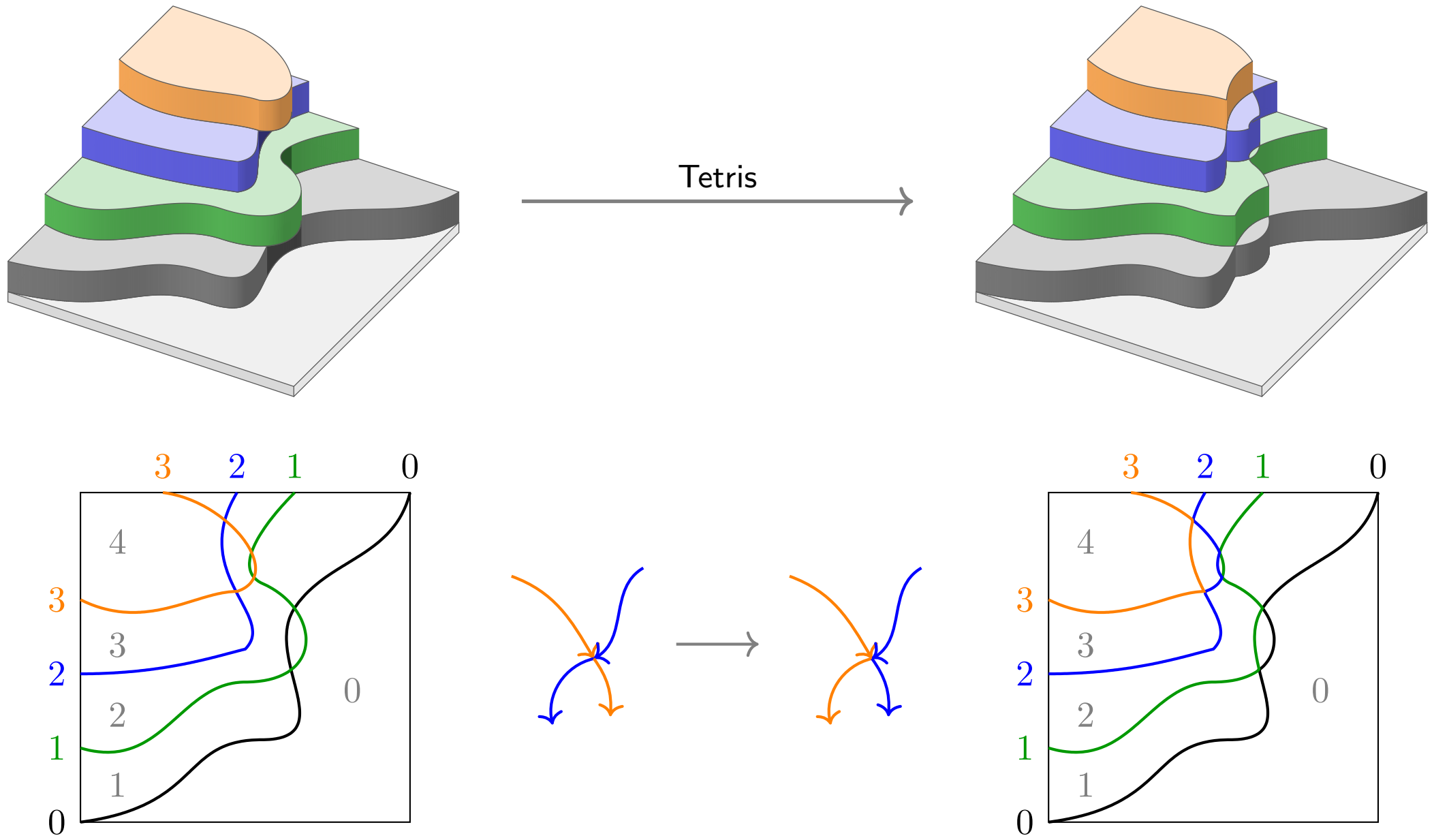
\begin{figure}
\vspace{-0.2cm}
\centering
\begin{tikzpicture}[thick]
  \node[inner sep=0pt] at (2.3,6.3) {\includegraphics[width=4.5cm]{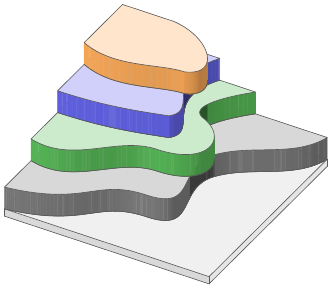}};
  \node[inner sep=0pt] at (11.7,6.3) {\includegraphics[width=4.5cm]{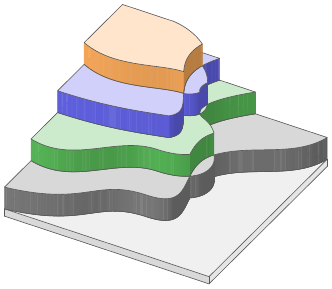}};
  \draw[->, gray, line width=1pt] (5.1,6.3) -- node[above, black] {\small\Tetris} (8.9,6.3);
  \node[inner sep=0pt] at (2.3,2) {\includegraphics{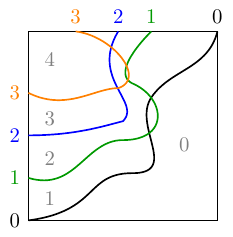}};
  \node[inner sep=0pt] at (11.7,2) {\includegraphics{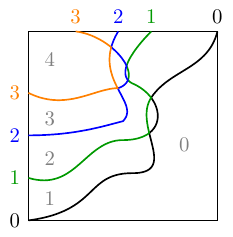}};
  \begin{scope}[xshift=5.0cm, yshift=0.9cm, scale=.8]
  \draw[->,orange] (0,2.2) to[out=-20,in=120] (1,1.2);
  \draw[->,blue] (1.6,2.3) to[out=-150,in=30] (1,1.2);
  \draw[->,blue] (1,1.2) to[bend right=45] (0.5,0.4);
  \draw[->,orange] (1,1.2) to[bend left=20] (1.2,0.5);
  \end{scope}
  \draw[->, gray, line width=.8pt] (6.6,2) -- (7.4,2);
  \begin{scope}[xshift=7.7cm, yshift=0.9cm, scale=.8]
  \draw[->,orange] (0,2.2) to[out=-20,in=120] (1,1.2);
  \draw[->,blue] (1.6,2.3) to[out=-150,in=30] (1,1.2);
  \draw[->,orange] (1,1.2) to[bend right=45] (0.5,0.4);
  \draw[->,blue] (1,1.2) to[bend left=20] (1.2,0.5);
  \end{scope}
\end{tikzpicture}
\vspace{-0.3cm}
\caption{A description of an interface $\cI$ as a set of level lines. They can cross each other if the interface has overhangs, so the colors are important (no loops are pictured, but they can occur). The effect of the \Tetris map is simply the unique rerouting of the same set of lines which makes them ordered.
At the bottom central pane: an oriented line reaching an intersection has two options to continue. The decorated diagram records, at every intersection, the local matching of the colors, which suffices to reconstruct the global picture.}
\label{fig:oct-levels}
\end{figure}

\begin{figure}
\vspace{-0.3cm}
\begin{tikzpicture}[thick]
  \draw (0,0) .. controls (0.2,1.6) and (1.8,1.9) .. (2,0.9) -- (2,0);
  \draw (3.2,1.7) .. controls (2.4,1.6) and (2.12,1.2) .. (2.12,0.9) -- (2.12,0);
  \draw[->, thin] (1.82,0.15) -- (1.82,0.7);
  \draw[->, thin] (2.3,0.7) -- (2.3,0.15);
\end{tikzpicture}
\vspace{-0.2cm}
\caption{Two level lines with opposite signs running along a common stretch cancel. This, with the fact that many level lines can meet on a single line, makes the combinatorics of counting preimages nontrivial.}
\label{fig:oct-cancel}
\vspace{-0.2cm}
\end{figure}
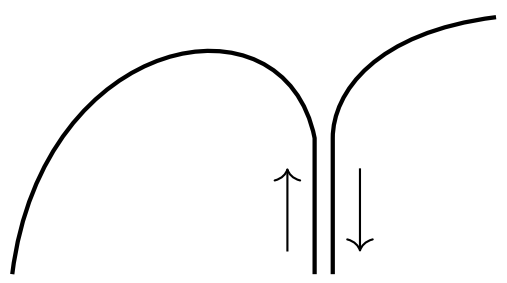

The difficulty in enumerating such preimages is visible in
\cref{fig:oct-levels,fig:oct-cancel}: described through its level
lines, an interface with overhangs is a family of oriented curves
which may cross, run along one another and cancel, and the \Tetris
map reroutes them into a nested family. We therefore replace the
level-line description by a local one, which records for each
column how many vertical disagreements it contains (dots) and for
each edge between columns how many horizontal disagreements it
contains, together with their directions (arrows), and then records
the vertical order in which these disagreements occur around each
corner of the squares. The point of this representation is twofold. First,
the energy is simply the total number of dots and arrows, so that a
bound on the excess energy bounds the number of symbols that differ
from those of the reference surface. Second, by
\cref{thm:reconstruction}, the local vertical orders determine the
interface, so that counting preimages reduces to counting
decorations, which is a problem about finite words
(\cref{lem:dot-first}).

In the one-step estimates below, $h$ denotes the fixed reference height;
$\fh$ denotes the physical height observable in the limit theorems.
We describe an \SOS reference configuration $h$ by its gradients; this determines the configuration only up to a global constant on the torus, which is the equivalence we work modulo in any case. Graphically, we identify a column with the unit square of the plane below it and address the square by its center, a point of $(\frac12,\frac12)+\Z^2$; the integer labels of $\T_n$ differ from these centers by $(\frac12,\frac12)$. The corners of the squares are the points of $\Z^2$, each incident to four squares and four edges.

We represent a general interface by placing a dot ($\bullet$) in a
square for each vertical spin disagreement in its column, and an
arrow ($\downarrow/\uparrow$) across an edge for each horizontal spin
disagreement between its two columns, directed from the $-1$ cube
to the $+1$ cube. For an \SOS
reference $h$, there is exactly one dot in each square, and there are $|\nabla_e h|$
arrows across each edge $e$, from the lower value to the
higher one.

\begin{definition}[Admissible interfaces]
Here an interface is the full disagreement surface of a spin
configuration with the prescribed ends, with an edge-connected
plaquette set, and with each spin connected through $*$-adjacent cubes to the corresponding occupied or vacant end.
For the following one-step statements the ends are vertical. A bare
plaquette set crossing lines an odd number of times is not a substitute
for this definition: it need not be the boundary of a spin set, and
the spin configuration is the object whose fibers are sorted. The
connectivity requirements are those of \cref{def:interface-model},
restated for the occupied and vacant ends of the present
coordinates; by \cref{lem:phase-tree} below, the requirement that
every spin be connected to its end already implies that the
disagreement surface is edge-connected.
\end{definition}

For a general interface, a column may contain several vertical
disagreements.
\begin{definition}[Dot-and-arrow]
We associate the following \emph{dot-and-arrow diagram} to an interface corresponding to a spin configuration $\sigma$. For a square $x\in(\frac12,\frac12)+\Z^2$ and a pair of adjacent squares $x,y$, let
\begin{equation}\label{eq:dot-arrow}
 \Dt_\sigma(x)=\sum_{x_3} \one_{\{\sigma(x,x_3) \neq \sigma(x,x_3+1)\}}\,,\qquad
 \Arr_\sigma(x,y)=\sum_{x_3} \one_{\{\sigma(x,x_3) = -1,\, \sigma(y,x_3) = +1\}}\,,
\end{equation}
and put $\Dt_\sigma(x)$ dots in the square $x$ and $\Arr_\sigma(x,y)$ arrows across the directed edge $\vec e = (x,y)$. Further let, for a set $S$ of squares,
\[
 \Arr_\sigma(S)=
 \sum_{\substack{\vec e=(x,y)\\\vec e\text{ incident to }S}}
       \Arr_\sigma(x,y)
\]
(summing over directed edges $\vec e$).
For an \SOS reference $h$, this is defined vis-\`a-vis  $\sigma_h(x,x_3)=\sgn(h(x)-x_3)$, whence
\[
 \Dt_h(x)=1,\qquad \Arr_h(x,y)=(h(y)-h(x))_+\,,\qquad
 \Arr_h(S)=\sum_{e:\ e\text{ incident to }S}|\nabla_e h|\,.
\]
\end{definition}
NB.\ Each arrow or dot again represents a disagreement between two cubes, with vertical disagreements represented by dots and horizontal disagreements represented by arrows, whose directions also encode which side of the disagreement is $+1$ and which is $-1$. So, the total energy is the number of arrows and dots, $\sA(\cI)=\sum_x\Dt_\sigma(x)+\sum_{(x,y)}\Arr_\sigma(x,y)$, the second sum running over ordered pairs of adjacent squares.

\begin{figure}
\vspace{-0.2cm}
\centering
\begin{minipage}[c]{.3\textwidth}
\centering
\includegraphics[scale=.6]{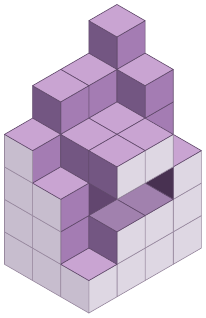}
\end{minipage}\hspace{-0.5cm}
\begin{minipage}[c]{.3\textwidth}
\centering
\includegraphics[scale=.6]{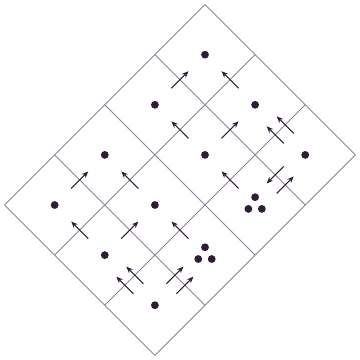}
\end{minipage}\hspace{-0.5cm}
\begin{minipage}[c]{.3\textwidth}
\centering
\includegraphics[scale=.6]{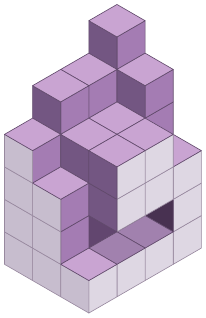}
\end{minipage}
\caption{A dot-and-arrow diagram (middle) shared by two interfaces with an overhang (left and right).}
\label{fig:dot-arrow}
\vspace{-0.1cm}
\end{figure}

\begin{figure}
\centering
\begin{minipage}[c]{.36\textwidth}
\centering
\includegraphics[scale=.6]{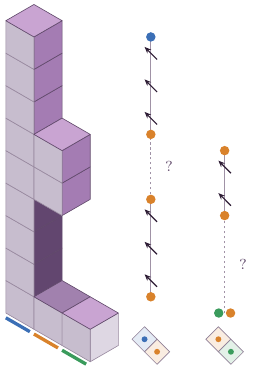}
\end{minipage}%
\begin{minipage}[c]{.26\textwidth}
\centering
\includegraphics[scale=.6]{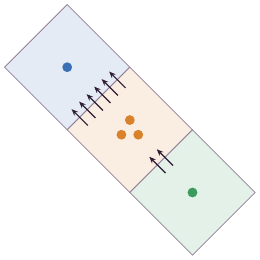}
\end{minipage}%
\begin{minipage}[c]{.36\textwidth}
\centering
\includegraphics[scale=.6]{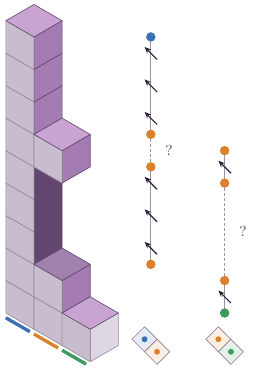}
\end{minipage}
\caption{Two pieces of interface (left and right) with the same dot-and-arrow diagram (middle), along with their individual dot-and-arrow decorations. Dashed stretches are gaps which cannot be read from the arrows of that edge alone. The two decorations agree on the first edge, where the dashed gap is the thickness of the overhang, and they differ on the second edge, which determines it.}
\label{fig:oct-ribbon}
\end{figure}

The map from general interfaces to dot-and-arrow diagrams is not injective (see \cref{fig:dot-arrow}), although if every square has exactly one dot then the arrows determine an \SOS configuration exactly; indeed, in that case each column is a step sequence whose height is
its charge, and the arrows record the height differences. In general, the dot-and-arrow diagram forgets the heights at which the
disagreements occur, which is needed, e.g., to tell whether an overhang in one column
lies above or below a disagreement in the neighboring column. The
following decoration, depicted
in \cref{fig:oct-ribbon}, restores exactly the missing information.
\begin{definition}[Decorated dot-and-arrow]\label{def:decorated}
    A \emph{decorated dot-and-arrow diagram} is a dot-and-arrow diagram in which, additionally, for every corner $v\in\Z^2$ a weak order between the arrows and dots on its four incident edges and squares is specified. Ties are only allowed between dots on different faces or between arrows on different edges. We emphasize that before the ordering, the dots and arrows of each type at one location are not distinguished, so if there are no ties the ordering is equivalently a word on a $12$-letter alphabet (one for each $\bullet$ in one of the $4$ incident squares, and one for each $\downarrow$ and $\uparrow$ on the $4$ incident edges).
\end{definition}
The map from interfaces to dot-and-arrow diagrams extends to decorated diagrams using the order induced by $x_3$.
Orders use the vertical coordinate of face centers: dots have integer
height and arrows have half-integer height, so that a dot and an
arrow are never tied, which is why ties are only allowed within the
two families. A run of arrows on a fixed
edge, with no intervening dot from either incident square, has constant
direction. Indeed, between two consecutive dot heights of the two
incident columns both columns have constant spin, so the vertical
faces between them all separate the same occupied column from the
same vacant one. This assertion concerns the order restricted to that edge
and its two incident squares, not an arbitrary interleaving at a vertex:
dots on the other two squares incident to the vertex may be interleaved
with such a run without affecting its direction.

\begin{theorem}[Interface reconstruction]\label{thm:reconstruction}
The map from admissible interfaces to decorated dot-and-arrow diagrams
is injective modulo vertical translation. Specifying one face height
removes this ambiguity.
\end{theorem}
\begin{proof}
The decoration does not record the heights of the faces, so we cannot
hope to recover them directly; the strategy is instead to recover,
from the data at each corner $v\in\Z^2$, which faces are adjacent and the
relative vertical position of adjacent faces, and then to assemble
the surface from these local incidences by propagating along it.

Consider a vertical lattice line through a vertex of the square grid.
There are four incident cube columns. List the distinct heights of dots
in these columns as $z_1<\cdots<z_k$, keeping the ties. The weak order
at the vertex tells which dots are tied and in which order the
distinct dot heights occur, although not their values. The four spins
are all occupied below $z_1$, and their values in every interval
$(z_j,z_{j+1})$ are determined by toggling a column at each of its dots,
so that the spin pattern of the four columns in each interval is read
off from the dot order alone.
Consequently the set of incident edges carrying vertical faces in this
interval is determined. If it is nonempty, each such edge carries
exactly $z_{j+1}-z_j$ arrows, one at each half-integer height, and the
decoration places these arrows between the dots bounding the
interval. Their
number determines the interval length and their ranks determine which
arrows on different edges are at the same height. If it is empty,
there are no faces in the open interval and no adjacency can cross it;
its length is not needed to determine incidences. (The question marks
in \cref{fig:oct-ribbon} indicate gaps which cannot be read from the
arrows of a single edge; when some incident edge is active in the
interval, its arrows determine the gap, and otherwise the gap is
irrelevant for the incidences.)

This reconstructs every local face incidence. Dots on neighboring
squares give adjacent horizontal faces exactly when their heights are
tied. Vertical faces on different incident edges are adjacent exactly
when their heights are tied. On the same edge consecutive unit
intervals share a horizontal edge, including intervals on opposite
sides of a dot height. A horizontal and a vertical face share an edge
exactly when the square is incident to that edge and the unit vertical
interval has the dot as an endpoint. The preceding reconstruction
specifies all these conditions and the relative vertical displacement
of every adjacent pair. Every shared lattice edge is seen at one of
these vertical lines: a horizontal lattice edge is incident to a grid
vertex of the square lattice, and a vertical lattice edge lies on a
vertical line through such a vertex. Now fix one face and assign it
height zero. Propagating the relative displacements along paths of
adjacent faces assigns a height to every face of the surface, because
the surface is edge-connected (this is part of the admissibility requirement).
Consistency around cycles is automatic for data coming from an
interface, since the relative displacements are those of the actual
surface. Finally the spins are recovered from the surface and the
occupied lower end, by toggling the spin each time a column crosses
a horizontal face. Two interfaces with the same decoration therefore
differ by the vertical translation taking the chosen face of one to
that of the other, and prescribing the height of one face removes
this ambiguity. This proves the assertion.
\end{proof}

Fix an \SOS height function $h$ and an admissible interface $\cI$
with spin configuration $\sigma$ and $\Tthree(\cI)=h$.
The number $\Dt_\sigma(x)$ of dots above $x$, see \cref{eq:dot-arrow}, is
the number of vertical disagreements in the column above $x$.
Since the column is occupied far below and vacant far above, $\Dt_\sigma(x)$ is
odd. If $\Dt_\sigma(x)=1$, the column is a step sequence, and
charge conservation forces its step to be at height $h(x)$, i.e.,
the column equals $\sigma_h$. The \emph{defect support} is therefore
exactly $\{x:\Dt_\sigma(x)>1\}$, the set of columns in which $\cI$ differs from
$\sigma_h$, and its components are taken in the square lattice with
$*$-adjacency, so that two defect squares sharing only
a corner belong to the same component. For an edge $e=xy$, let
\begin{equation}\label{eq:arrm}
 \Arrm{\sigma}(e)=\Arr_\sigma(x,y)\;\wedge\;\Arr_\sigma(y,x)\,.
\end{equation}
The signed difference of the two arrow counts
is
\[
 \Arr_\sigma(x,y)-\Arr_\sigma(y,x)
 =\sum_{x_3}\big(\omega_y(x_3)-\omega_x(x_3)\big)
 =q(\omega_y)-q(\omega_x)=h(y)-h(x)\,,
\]
where $\omega_x=\one_{\{\sigma(x,\cdot)=+1\}}$ is the fiber above $x$ as in
\cref{def:tetris-charge}, the flat-step normalizations of the two charges
cancel, and the last identity is charge conservation, so that
$\Arr_\sigma(x,y)+\Arr_\sigma(y,x)=|h(y)-h(x)|+2\Arrm{\sigma}(e)$. Since
the energy of $\cI$ is the number of dots plus the number of arrows,
whereas the energy of $h$ is the number of squares plus
$\sum_e|h(y)-h(x)|$, we obtain
\begin{equation}\label{eq:exact-excess}
 E(\cI;h):=\sA(\cI)-\sA(h)=\sum_x(\Dt_\sigma(x)-1)+2\sum_e \Arrm{\sigma}(e)\,.
\end{equation}
All the nonzero terms are supported on the defect squares and their
incident edges: an edge between two non-defect columns separates two
step sequences, whose arrows all point the same way, so that $\Arrm{\sigma}(e)=0$.
Every nonempty defect column has $\Dt_\sigma(x)\geq3$, by the parity observed
above, so that the first sum contributes at least $2$ per defect
square.

\begin{lemma}[The phase tree]\label{lem:phase-tree}
Let $\sigma$ be an Ising configuration on the tilted cylinder of
\cref{def:3dIsing} which agrees with the boundary conditions outside a
finite tilted height band. Let $G$ be the graph whose
vertices are its $*$-connected spin components and whose edges
are the edge-connected components of its disagreement plaquettes.
Each plaquette component is incident to exactly two spin
components, one of each sign, so the graph is well defined.
Then $G$ is a finite tree.
In particular, the no-bubble condition implies that the entire
disagreement surface is edge-connected.

The same statements hold for an Ising configuration $\sigma$ on $\Z^3$ which differs from $\sigma_h=\sgn(h(x_1,x_2)-x_3)$, for some \SOS height function $h:(\frac12,\frac12)+\Z^2\to\Z$, in finitely many cubes.
\end{lemma}
\begin{proof}
The idea is that each plaquette component $F$ is by itself the
disagreement surface of a binary configuration $\tau_F$,
equal to $1$ on one side of $F$ and $0$ on the other, and that $\tau_F$ is
constant on each spin component; the tree structure then follows
from the fact that removing $F$ disconnects exactly the two spin
components on its two sides.

Fix a plaquette component $F$. Mark a nearest-neighbor bond
when its dual plaquette lies in $F$. Around every elementary
primal square the total number of disagreement bonds is even,
since going around the square returns to the starting spin.
Their dual plaquettes all share the dual edge of this square,
so they are pairwise edge-adjacent and either all of them belong
to $F$ or none does; in particular the number of marked bonds
around every elementary square is even. The marked
bonds therefore form a closed $\Z/2$ cochain on the
cover $\Z^3$. Integrating along lattice paths, i.e., starting
from the value $1$ at a cube deep in the occupied end and toggling
the value each time a marked bond is crossed, defines a binary
function $\tau_F$ whose disagreement set is exactly the lift
of $F$: path independence follows by filling a closed lattice
path with elementary squares, each of which is crossed by an even
number of marked bonds. The function is constant on
each of the two tilted ends, which contain no plaquettes of $F$.
In particular it is periodic:
its difference from its translate by a period is constant, being
the integral of the same cochain along translated paths,
and this constant is zero by comparison sufficiently far in
the upper end. This deals with possible noncontractible components
of the quotient without assuming that they bound a ball.

The function $\tau_F$ is constant on each $*$-connected
component of the original spin configuration. For a face-step
of equal spins no disagreement bond is crossed. For an
edge-step, consider the two length-two face-paths around the
elementary square joining the two cubes. If an intermediate
spin agrees with their sign, that path crosses no boundary.
Otherwise both intermediate spins are opposite, all four
bonds disagree, and all four dual faces belong to one
plaquette component. Each length-two path then crosses $F$
either twice or not at all. This proves constancy also for
edge-steps.

If two disagreement plaquettes share an edge, their occupied
incident cubes lie among the four cubes around that edge and
are $*$-connected (any two of these four cubes share a face or an
edge); the same is true of their vacant incident
cubes. Following a plaquette path through $F$ proves that all
occupied cubes incident to $F$ lie in one spin component, and all
vacant ones in another, so $F$ is incident
to a fixed pair of spin components. Its function $\tau_F$
takes opposite values on this pair and equal values across
every other plaquette component, by construction. Consequently
any path in the graph from one of these components to the
other must use the edge $F$, i.e., this edge is a bridge.
The graph is connected,
since a finite face-path between any two cubes projects to
a path in it (each face-step either stays in a spin component or
crosses one disagreement plaquette). There are finitely many
vertices and edges on
the quotient, since only a finite height band is nonconstant.
A connected finite graph all of whose edges are bridges is
a tree. Under the no-bubble condition its only vertices are
the occupied and vacant end components, so it has exactly
one edge, i.e., the disagreement surface is edge-connected.
The same cochain and bridge argument applies on $\Z^3$ to an
spin configuration that differs from a graph configuration in
finitely many cubes. There are still only finitely many spin components:
changing finitely many vertices in the locally finite adjacency graph
can split either of the two connected reference phases into only
finitely many components. Thus the phase graph is again a finite tree.
Under the no-bubble condition, both spin phases are connected to their
connected reference ends, so $\tau_F$ must equal the indicator $\one_{\{\sigma=+1\}}$ of the
original plus phase or its complement $\one_{\{\sigma=-1\}}$.
\end{proof}

\begin{remark*}
By \cref{lem:phase-tree}, requiring the disagreement surface to be edge-connected, as in the definition of admissible interfaces above, is no additional restriction on the Ising Interface model: the no-bubble condition of \cref{def:interface-model} already implies it.
\end{remark*}

\begin{lemma}[Independent defect components]\label{lem:polymer-factorization}
Let $h$ be an \SOS height function on the plane or on the torus $\T_n$, and let $D_1,\ldots,D_k$ be finite $*$-connected sets of squares which are pairwise non-$*$-adjacent; a set $D_i$ may wrap around the torus, and may be the whole torus. Call an admissible preimage of $h$ whose defect support is exactly $D$ a \emph{$D$-defect}, or an admissible modification of $h$ on $D$; by definition, it agrees with $\sigma_h$ on every column outside $D$. A \emph{defect} is a $D$-defect for some $*$-connected $D$. Then the following hold.
\begin{enumerate}[(i)]
\item If $\cI_i$ is an admissible modification of $h$ on $D_i$ for each $i$, then the configuration $\cI$ which agrees with $\cI_i$ on the columns of $D_i$ for every $i$, and with $\sigma_h$ on all other columns, is an admissible preimage of $h$ with defect support $D_1\cup\cdots\cup D_k$, and its excess energy is $E(\cI;h)=\sum_{i=1}^k E(\cI_i;h)$.
\item Conversely, every admissible preimage of $h$ with defect support $D_1\cup\cdots\cup D_k$ arises in this way from a unique choice of $\cI_1,\ldots,\cI_k$, namely its restriction to the columns of $D_i$ completed by $\sigma_h$ elsewhere.
\end{enumerate}
The same holds, with the same proof, for the \SOS and \HMSOS inputs of the later steps, where admissibility includes the additional monotonicity constraints.
\end{lemma}
\begin{proof}
A preimage of $h$ is specified by its columns, and the constraints
it must satisfy are of three kinds: the charge of every column is
prescribed by $h$; any additional coordinate monotonicities (for
the \SOS and \HMSOS inputs of the later steps) are constraints on
pairs of neighboring columns; and the no-bubble condition
of \cref{def:interface-model} requires every cube to be
$*$-connected to its end. The first two kinds are local. Indeed,
the charge constraint concerns one column at a time, and a
monotonicity constraint between two neighboring columns involves
either two columns of the same support or one changed and one
unchanged column (two distinct supports are not $*$-adjacent, hence
not adjacent), so that it is checked within the isolated
modification on each $D_i$. Likewise, every changed bond joins two
cubes of one support or a cube of one support to a cube of an
unchanged column, and \cref{eq:exact-excess}, whose terms are
supported on the defect squares and their incident edges, shows that
the excess energies of the isolated modifications add.

It remains to see that the no-bubble condition also separates. Follow
a same-sign path of $*$-steps from a cube
in $D_i$ to its required end. At its first exit from $D_i$ the
path enters an unchanged column: a $*$-step between cubes projects to at
most one $*$-step between squares, so it cannot enter another support. Continue vertically
along that graph column to the same end; this is possible because
an unchanged column is the step sequence $\sigma_h$, whose spin of
the given sign forms a half-line reaching its end. If the path never
exits, it already reaches the end inside $D_i$. The resulting path uses
only the isolated patch and unchanged columns. It is thus
valid both before and after inserting any of the other patches, since
these only change columns of other supports.
Thus each spin reaches its end in the combined configuration if
and only if it does in its isolated patch. This proves independent
insertion and deletion on either transverse topology, including
when one support is the whole torus, in which case there is nothing
to separate. \Cref{lem:phase-tree}
then gives edge-connectedness of the disagreement surface of the
combined configuration, so it is an admissible interface.
\end{proof}

\begin{lemma}[Counting orders]\label{lem:dot-first}
Fix $v\in\Z^2$. Let $d_v$ be the total number of dots on the $4$
squares incident to $v$ and let $a_v$ be the total number of arrows on the $4$
edges incident to $v$. For fixed dot-and-arrow multiplicities, there are at most $8^{d_v}\binom{a_v+d_v}{d_v}$
weak orders at $v$ in \cref{def:decorated} realizable by
admissible interfaces.
\end{lemma}
\begin{proof}
First choose the order of the dots. Recording their column labels
as a word on $4$ letters, together with ties between successive
entries, gives at most $4^{d_v}2^{d_v-1}\leq8^{d_v}$ possibilities.
We may include unrealizable choices in this upper bound.

Fix this dot order. The four spins are initially $+1$ below all
dots, and each column's spin flips when one of its dots is crossed.
Thus the dot order determines the spins in every gap between
consecutive distinct dot levels, and hence which incident edges
carry arrows there. Call a gap \emph{active} if at least one edge
does.

In an active gap, each such edge carries one arrow at every
half-integer height, with direction determined by the spins.
Consequently the gap lengths sum to at most $a_v$. The active gaps
form a fixed list of at most $d_v-1$ intervals, so stars-and-bars
bounds the number of possible length vectors by
$\binom{a_v+d_v}{d_v}$. These lengths determine the order and ties
of all arrows relative to the dots. Inactive gaps contain no
arrows, so their lengths do not affect the weak order. Multiplying
the two bounds proves the claim.
\end{proof}

In what follows and throughout the paper, for a set $D$ of squares, we define its \emph{one-square
neighborhood} to be  $D^\dagger=\{x:\dist_*(x,D)\leq1\}$,
and its \emph{two-square neighborhood} to be $D^\ddagger=\{x:\dist_*(x,D)\leq2\}$.

\begin{proposition}[Enumeration with multiplicities]\label{prop:correct-count}
Let $D$ be a nonempty $*$-connected set of $m$ squares.
For $E\geq2m$ the number $N_h(D,E)$ of one-step preimages of $h$ with
single defect support $D$ and excess $E$ has
\begin{equation}\label{eq:correct-count}
 N_h(D,E)\leq2^{64E}\binom{2\Arr_h(D^\dagger)+16E}{16E}\,.
\end{equation}
It is zero for odd $E$ or $E<2m$.
\end{proposition}
\begin{proof}
We count first the dot-and-arrow diagram, then
its decoration, and then invoke \cref{thm:reconstruction}.
Let $\sigma$ be the spin configuration of a preimage being counted.

Write $\Dt_\sigma(x)=1+2k_x$, which is possible since $\Dt_\sigma(x)$ is odd. By
\cref{eq:exact-excess},
\[
 \sum_{x\in D}k_x+\sum_{e\text{ incident to }D}\Arrm{\sigma}(e)=E/2\,.
\]
There are at most $5m$ slots for these nonnegative integers, namely
the $m$ squares of $D$ and the at most $4m$ edges incident to them.
The stars-and-bars count bounds
their number by $\binom{E/2+5m}{5m}\leq2^{E/2+5m}\leq2^{3E}$, using
$m\leq E/2$. The net gradient fixes both
arrow multiplicities once $\Arrm{\sigma}(e)$ is fixed, since by \cref{eq:arrm} the two counts on
$e=xy$ are $\Arrm{\sigma}(e)$ and $\Arrm{\sigma}(e)+|h(y)-h(x)|$, with the larger one in the
direction of increasing $h$. Every column outside $D$ carries one
dot and the arrows of $h$. Thus the dot-and-arrow diagram is
determined by these integers.

The one-square neighborhood in $\Arr_h(D^\dagger)$ includes all edges at affected vertices.
Only vertices incident to a square of $D$ can have a nontrivial choice
of decoration: at any other vertex all four columns are step
sequences, and the order is the one of $h$. There are at most
$4m\leq2E$ such vertices. Let $d_v,a_v$ be the dot and arrow counts
of $\sigma$ at these vertices, as in \cref{lem:dot-first}. Then
\begin{align*}
 \sum_v d_v&\leq4(4m)+4\sum_x(\Dt_\sigma(x)-1)\leq12E,\\
 \sum_v a_v&\leq2\Arr_\sigma(D^\dagger)
                 \leq2\Arr_h(D^\dagger)+2E\,.
\end{align*}
In the first bound, each of the at most $4m$ vertices sees four
squares each carrying one baseline dot, including the dots on
neighboring good squares, and each excess dot $\Dt_\sigma(x)-1$ is seen at
the four vertices of its square; then $16m\leq8E$ and
$4\sum_x(\Dt_\sigma(x)-1)\leq4E$.
The second line counts each incident arrow of $\sigma$ at most twice, since an
edge has two endpoints: the arrows of $h$ contribute at most
$2\Arr_h(D^\dagger)$ and the $2\sum_e\Arrm{\sigma}(e)\leq E$ additional arrows at most
$2E$. Therefore
\[
 \prod_v8^{d_v}\binom{a_v+d_v}{d_v}
 \leq2^{36E}\binom{2\Arr_h(D^\dagger)+14E}{12E}
 \leq2^{50E}\binom{2\Arr_h(D^\dagger)+16E}{16E}\,.
\]
The product inequality chooses subsets in disjoint sets and then
increases the total arrow and dot counts separately. For the last
inequality, split a $12E$-element subset into its parts in sets of
sizes $2\Arr_h(D^\dagger)$ and $14E$. Record the second part arbitrarily, at
cost $2^{14E}$, and pad the first part to size $16E$ using the first
available elements of a new ordered $16E$-element set. This injects
the choices into the claimed bound. Multiplying by the multiplicity count
gives $2^{53E}$ times the displayed binomial coefficient, which is
bounded by \cref{eq:correct-count}. Not every choice of local orders
is realized by an interface, and globally inconsistent orders
are simply discarded, as we only need an upper bound. By
\cref{thm:reconstruction}, a realizable decoration determines the
interface up to a vertical shift, and
the prescribed charges fix that shift. This completes the count.
Finally, $E$ is even by \cref{eq:exact-excess}, and $E\geq2m$ since
every defect column contributes at least $\Dt_\sigma(x)-1\geq2$.
\end{proof}

Recalling the definition of $N_h(D,E)$ from \cref{prop:correct-count}, define the \emph{activity} of a defect support $D$ to be
\begin{equation}\label{eq:activity}
 z_{\beta,h}(D):=\sum_E N_h(D,E)e^{-\beta E}\,.
\end{equation}

\begin{corollary}[Defect activities]\label{cor:activities}
There is an absolute $\beta_0$ such that for $\beta\geq\beta_0$ and  every defect support~$D$,
\begin{equation}\label{eq:activity-bound}
 z_{\beta,h}(D)\leq
 2\exp(2e^{-\beta/64}\Arr_h(D^\dagger)-\beta|D|)\,.
\end{equation}
\end{corollary}
\begin{proof}
The activity is the weighted number of modifications supported
on $D$, and the point of the bound is that the dependence on the
reference gradients $\Arr_h(D^\dagger)$, which may be arbitrarily large, only
enters through a factor $e^{\delta \Arr_h(D^\dagger)}$ with $\delta\to0$ as
$\beta\to\infty$. For $a=\Arr_h(D^\dagger)$ and $0<t<1$, the binomial theorem gives
$\binom{2a+16E}{16E}\leq t^{-16E}(1+t)^{2a+16E}$, since the left-hand
side is one of the terms of the expansion of $(1+t)^{2a+16E}$
divided by $t^{16E}$.
Take $t=e^{-\beta/64}$ and use $\log(1+t)\leq t$, so that
$(1+t)^{2a}\leq e^{2e^{-\beta/64}a}$ is the first factor in
\cref{eq:activity-bound}. The logarithm of
the remaining factor $2^{64E}t^{-16E}(1+t)^{16E}e^{-\beta E}$ is at most
$(64\log2+\beta/4+16t-\beta)E=(64\log2-3\beta/4+16t)E\leq-\beta E/2$
when $\beta_0$ is large
enough. Summing $e^{-\beta E/2}$ over all integers $E\geq2|D|$ gives
at most $2e^{-\beta|D|}$ for $\beta\geq2\log 2$, which is the claimed
bound; summing also over the forbidden odd integers is harmless.
\end{proof}

\begin{proposition}[Bound on the partition function]\label{prop:positive-budget}
Let $W$ be any finite set of reference plaquettes of $h$. Sum the
weights of all compatible one-step defect families for which every
reference dot and arrow incident to $D_i^\dagger$ belongs to $W$. Then,
for all sufficiently large $\beta$,
\[
 \sum_{\{D_i\}\text{ compatible, supported in }W}
                  \prod_i z_{\beta,h}(D_i)
          \leq\exp(Ce^{-\beta/64}|W|)\,.
\]
The empty family is included, and the supports $D_i^\dagger$ may overlap.
\end{proposition}
\begin{proof}
By \cref{lem:polymer-factorization}, a preimage of $h$ is the same
as a family of pairwise non-$*$-adjacent supports together with a
$D$-defect on each support $D$, and the excess energies add. Thus the
one-step fiber partition function in a finite transverse region is
\[
 Z_{\beta,h}=\sum_{\{D_i\}\text{ pairwise non-$*$-adjacent}}
                         \prod_i z_{\beta,h}(D_i)\,.
\]
The sum in the statement restricts this gas to supports whose
specified reference dots and arrows lie in $W$.

We first control the total gradient seen by the neighborhoods. Each edge
is incident to at most eighteen distinct neighborhoods $D_i^\dagger$:
an incidence requires $D_i$ to meet one of the at most $2\cdot9=18$
squares within one $*$-step of an endpoint, and distinct compatible
supports are disjoint. Every arrow counted belongs to $W$, so
$\sum_i\Arr_h(D_i^\dagger)\leq18|W|$. Extracting the corresponding
factor from \cref{eq:activity-bound}, then dropping compatibility,
bounds the sum in the statement by
\[
 e^{36e^{-\beta/64}|W|}\prod_D(1+2e^{-\beta|D|})
 \leq \exp\bigg(36e^{-\beta/64}|W|+2\sum_De^{-\beta|D|}\bigg)\,,
\]
where the product and sum run over the allowed supports.
There are at most $|W|$ possible root squares, since each square of
a support must have its reference dot in $W$. Root a spanning tree
of a support of size $m$ at its lexicographically first square and
encode its depth-first traversal by $2(m-1)$ steps, each with at most
eight choices. Thus there are at most $|W|64^m$ such supports, and
\[
 2\sum_De^{-\beta|D|}
 \leq2|W|\sum_{m\geq1}64^me^{-\beta m}
 \leq Ce^{-\beta}|W|
\]
for all sufficiently large $\beta$. Combining the two estimates
proves the bound.
\end{proof}

This collection-level estimate controls the positive one-step fiber
sum, including compatible families of defects. Its role is the
following: in \cref{sec:preimage-witnesses} the defects of several
\Tetris steps are grouped into patterns, and the coefficient of a
pattern is a sum of absolute weights over all defect families whose
footprint is that pattern. The present bound, which is an estimate
on sums of weights rather than on individual terms, is what allows
the sum over all defect families inside a pattern to be bounded in
terms of the size of the pattern.

For a one-step fiber with reference height $h$, define its
\emph{fiber truncation} at distance $J\in\N$ by requiring, on every sorting
column $v$,
\begin{equation}\label{eq:fiber-truncation}
 \omega_v(k)=\one_{\{k<h(v)\}}
                   \qquad\text{whenever }|k-h(v)|>J.
\end{equation}
Here $k$ is the coordinate along the sorting direction, after
relabeling coordinates if needed. This is an auxiliary restriction
relative to the reference surface. We reserve $M$ for the physical
height cutoff of \cref{def:3dIsing}, in the period-invariant band
$|x_3+(\Theta_1x_1+\Theta_2x_2)/n|\leq M$.

\begin{lemma}[One-step estimates with periodic boundary conditions]
\label{lem:quotient-one-step}
The preceding one-step factorization and estimates remain valid,
with the same constants, under our tilted periodic boundary
conditions, including the auxiliary quotients arising from the
different sorting directions.
Precisely, \cref{lem:polymer-factorization,prop:correct-count,cor:activities,prop:positive-budget}
hold on every finite transverse quotient of the square lattice with
prescribed height offsets, provided no nonzero physical period is
parallel to the sorting direction. One may retain any coordinate
monotonicities of the input. The constants are independent of the
transverse periods. The quotient graph may have loops or parallel
edges, and defect supports may wrap around the quotient or occupy
it entirely.
For large enough $\beta$, the normalized fiber truncated at distance $J$ as in \cref{eq:fiber-truncation}
converges in total variation as $J\to\infty$
to the normalized unrestricted fiber.
\end{lemma}
\begin{proof}
Keep geometric edge orbits and their directional incidences: parallel
edges are distinct and a loop has two endpoint incidences. The
first-exit proof of \cref{lem:polymer-factorization} is valid on this
quotient, and each additional coordinate monotonicity is a constraint
on neighboring columns. It therefore gives the same independent
insertion and deletion of non-$*$-adjacent supports.

For a support of $m$ squares and excess $E$, the exact excess identity
has at most $5m$ slots and gives $E\geq2m$. At each affected vertex
choose a lift, distinguishing its four square and four edge incidences
by their directions, even when some columns are translates of one
another. The counts in \cref{prop:correct-count} remain
\[
 \#\{v\text{ affected}\}\leq4m,\qquad
 \sum_v d_v\leq16m+4E\leq12E,\qquad
 \sum_v a_v\leq2\Arr_h(D^\dagger)+2E.
\]
Indeed every dot has four vertex incidences and every arrow has two.
Thus \cref{lem:dot-first} and the same stars-and-bars calculation
give \cref{eq:correct-count}. Period offsets impose consistency
conditions and only discard encodings. To see injectivity, reconstruct
along the connected quotient surface with these prescribed offsets;
its face incidences fix all heights up to a common shift, and the
column charges fix that shift. This also covers a support equal to
the whole quotient. The activity estimate now follows by
\cref{cor:activities}.

The support graph has at most eight neighbors per square, so the
same directional spanning-tree encoding gives at most $64^m$ rooted
supports. A reference arrow belongs to neighborhoods of at most eighteen
compatible supports: each such support contains a distinct square
among the at most $2\cdot9$ squares neighboring its two endpoint
incidences. Identifying some of these squares cannot increase the
count. There are at most $|W|$ possible roots when all required
reference dots and arrows lie in $W$. Consequently the proof of
\cref{prop:positive-budget} gives, without an injectivity-radius
assumption,
\[
 \sum_{\{D_i\}\text{ compatible, supported in }W}\prod_i z_{t,h}(D_i)
 \leq\exp\bigg(36e^{-t/64}|W|
                 +2|W|\sum_{m\geq1}64^m e^{-tm}\bigg)
 \leq e^{Ce^{-t/64}|W|}.
\]
Taking $W$ to be the whole reference surface proves finiteness of
the fiber partition function. Let $\mu$ be the probability law on
this fiber obtained by normalizing the weights $e^{-tE}$, and let
$A_J$ be the event that every column satisfies the restriction in
\cref{eq:fiber-truncation}. A prescribed period shifts both $k$ and
$h(v)$ by the same amount, so this event is well defined on the
quotient. At fixed $J$ there are finitely many column orbits and
finitely many sites on each at which the configuration can differ
from its reference step; hence $A_J$ contains finitely many
configurations. Every configuration in the unrestricted fiber
belongs to $A_J$ for all sufficiently large $J$: each of its
columns agrees with its reference step outside a finite interval,
and we take the maximum of the required distances over the finitely
many column orbits. Thus $A_J$ increases to the whole fiber and
$\mu(A_J^c)\to0$.

The truncated law $\mu_J$ has the same Boltzmann weights, restricted
to $A_J$ and normalized, so $\mu_J={\mu(\,\cdot\mid A_J)}$.
Consequently,
\[
 \|\mu_J-\mu\|_{\mathrm{TV}}
 =\frac12\sum_{\cI}|\mu_J(\cI)-\mu(\cI)|
 =\mu(A_J^c)\longrightarrow0.
\]
Indeed, conditioning removes mass $\mu(A_J^c)$ from $A_J^c$
and adds the same mass to $A_J$ by renormalization. The two
contributions to the sum are therefore both $\mu(A_J^c)$.
\end{proof}

For successive \Tetris maps, defects at different stages are coupled
through the intermediate surfaces. In \cref{sec:preimage-witnesses}
we group them into history components with footprints on the final
monotone surface. \Cref{lem:sos-carried-components,lem:interface-components}
prove that these components can be inserted and deleted independently,
with all intermediate admissibility constraints preserved and with
additive excess energies.

\section{Fibers, histories and geometric witnesses}
\label{sec:preimage-witnesses}

We study the fiber over a monotone reference surface $\varphi$ after
one, two or three \Tetris steps, that is, the set of surfaces of the
model under study which the composition of the corresponding sorting
maps sends to $\varphi$, weighted by their Boltzmann factors. In each
case a \emph{history} consists of the surface together with its
intermediate images, and it has a total excess $E$, the sum of its
one-step area excesses as in \cref{eq:exact-excess}. Its local
components are recorded by connected reference footprints $P$ on
$\varphi$, identified with their diagonal tiling patterns. The
intermediate surfaces and, for interfaces, the phase
connectivity constraints are retained in the history. This is
unavoidable: the one-step estimates of \cref{sec:enumeration} are
stated for an arbitrary admissible starting surface, but the
intermediate surface of a two- or three-step history is itself
random, and the excess of the second step is measured relative to it
rather than relative to $\varphi$; the one-step estimate cannot
simply be iterated after discarding the intermediate surfaces.

The same dichotomy governs the three constructions. In a typical
region, a component has enough excess relative to its footprint:
its Boltzmann weight $e^{-\beta E}$ pays for the possible choices of
the component, yielding exponentially small \emph{\texttt{thin}}
coefficients. The difficulty is that a component with small excess
may nevertheless have a large footprint, since a few modified columns
can meet long stacks of wall plaquettes of the reference surface, and
the activity of such a component is not small in terms of its
footprint. A large footprint with small excess is therefore enlarged
to a connected \emph{witness} $W\supset P$. These atypical patterns
may support many low-energy preimages, but their witnesses are
\texttt{fat}---almost all of their area is covered by large disks---and
they belong to a family of small entropy, in the sense that the
number of such witnesses of a given size grows at a rate that
vanishes as $\beta\to\infty$. Both properties are needed in the
first integration of \cref{sec:first-integration}: the \texttt{fat} property makes
a prescribed witness exponentially unlikely under the free dimer law
(\cref{prop:fat_pattern}), and the small entropy ensures that this
rarity survives the sum over all possible witnesses (recall
\cref{thm:rg-typical}).

We use $B$ for the excess-to-size threshold,
$R$ for the witness radius and $L$ for the long-stack threshold. The
precise classification and constants are stated in each construction:
the one-step argument first compares $|P|$ with its support size,
whereas the two- and three-step arguments compare it with $E$.
The same witness and entropy estimates include noncontractible patterns on
the torus. We stress that this geometric dichotomy concerns the
coefficients of the pattern expansion only; the probability of a
pattern under the interacting tiling law is a different matter,
which is settled only after the \RG theorem, in
\cref{sec:interacting-witnesses}.
The common estimates in \cref{sec:enumeration}
apply to a single step over an arbitrary admissible intermediate
surface. The three constructions below carry the relevant footprints
to the final monotone surface.

\subsection{One step: HMSOS fibers}\label{sec:hmsos}

\subsubsection{Reference footprints}
We record each defect on the diagonal tiling through a footprint
large enough to determine its activity. The point of this
bookkeeping is that the pattern expansion of \cref{thm:rg-typical}
attaches coefficients to patches of the tiling $\varphi$; we must
therefore express the weight of a defect as a function of a finite
patch of $\varphi$, and the footprint is chosen large enough to
do so. Unlike the random approximation in
\cite{LaslierLubetzky24}, where the conditional law of the surface
given the tiling involved long-range interactions, the sorting map
gives these footprints directly and deterministically from the
reference surface.

\begin{definition}[Reference footprint]\label{def:footprint}
For this one-step analysis the input is \HMSOS and the reference
is its monotone image $\varphi$. Relabel the sorting direction as
the vertical direction in the dot-and-arrow description, so that
the one-step estimates of \cref{sec:enumeration} apply with the
column magnetization of $\varphi$ as the reference height; we call
the plane of the two remaining coordinates, whose unit squares index
the fibers in the sorting direction, the \emph{sorting plane}, so
that defect supports are sets of its squares; and we
indicate the reference by the subscript $\varphi$ in place of $h$,
writing $z_{\beta,\varphi}(D)$, $\Arr_\varphi(D^\dagger)$ and
$Z_{\beta,\varphi}$. For a defect support $D$, let $P_\varphi(D)$
consist of the reference dots at $D^\ddagger$ and all the reference wall plaquettes on
edges with both endpoints in $D^\ddagger$.
This is a connected patch of the monotone surface; its diagonal
projection is a tiling pattern with $|P_\varphi(D)|$ lozenges. We identify
reference plaquettes with these lozenges.
In a geometric domain $\Lambda$ we allow exactly those polymers with
$\Upsilon(P_\varphi(D))\subset\Lambda$, retaining the original
incompatibility relation on their sorting-plane supports.
\end{definition}

The radius two in the definition of $D^\ddagger$ is dictated by the
enumeration of \cref{prop:correct-count}: the quantity $\Arr_\varphi(D^\dagger)$
appearing there sums the arrows of $\varphi$ on the edges incident to the
one-square neighborhood $D^\dagger$, and both endpoints of such an edge lie
within $*$-distance two of $D$. Every arrow counted by $\Arr_\varphi(D^\dagger)$
therefore belongs to $P_\varphi(D)$, so
$\Arr_\varphi(D^\dagger)\leq |P_\varphi(D)|$. More importantly, the footprint contains all
the reference columns needed to determine the activity
$z_{\beta,\varphi}(D)$ of \cref{cor:activities}, including their relative
heights: the $D$-defects, and their excess
energies, depend on $\varphi$ only through the heights of the columns in
$D^\dagger$ (which determine the charges to be conserved and the
baseline arrows at the affected vertices), and these heights are
read off from the dots and wall plaquettes of $P_\varphi(D)$ up to a
common additive constant, which does not affect the weight.
For \HMSOS the additional monotonicity constraints are local
constraints on neighboring columns, namely coordinatewise inclusion
of the $\{0,1\}$-spin sequences of adjacent columns in the two
non-sorting directions. Thus restricting the one-step gas to these
constraints preserves its factorization
(\cref{lem:polymer-factorization}) and only decreases the
enumeration bound \cref{eq:correct-count}, since it removes
preimages without adding any. No connectedness condition on distant
bubbles is needed in this case, in contrast with the interface
histories of \cref{sec:interface}.

\begin{lemma}[Bounded footprint overlap]\label{lem:footprint-overlap}
For a compatible defect family, every reference plaquette belongs to
at most $25$ of its footprints. Given a reference triangle $x$, every
$D$ whose footprint contains $x$ intersects a fixed set of at most
$25$ sorting-plane squares.
\end{lemma}
\begin{proof}
Assign to every reference plaquette a sorting-plane square $v$: for a
dot take the square it lies in, and for a wall plaquette take one of
the two squares separated by its edge, fixed once and for all.
If the plaquette belongs to $P_\varphi(D)$ then $v\in D^\ddagger$ by
\cref{def:footprint}, so $D$ intersects the set of $25$ squares at
$*$-distance at most two from $v$ (a $5\times5$ block centered at
$v$). The supports of a compatible family are pairwise disjoint,
being non-$*$-adjacent, so at most $25$ of them can meet this block.
This proves the first assertion. For the second, note that a
reference triangle $x$ lies in exactly one reference lozenge, which
is the projection of exactly one reference plaquette, and apply the
same argument to that plaquette.
\end{proof}

\subsubsection{Fattening the \texorpdfstring{\texttt{bad}}{bad} sets}\label{sec:fattening}

Fix a threshold $B\geq1$. Call a sorting-plane support $D$
\emph{\texttt{bad}} if its reference footprint is large compared with
its support, $|P_\varphi(D)|\geq B|D|$, and \emph{\texttt{good}} otherwise. The
threshold will be fixed as $B=B_\beta$ in \cref{eq:history-scales}.
The rationale is the following. The excess of any modification on
$D$ is at least $2|D|$ (\cref{prop:correct-count}), so its
Boltzmann factor is at most $e^{-2\beta|D|}$. For a \texttt{good}
support this is at most $e^{-2\beta |P_\varphi(D)|/B}$, an exponential
decay in the size of the footprint, as a \texttt{thin}
coefficient requires. For a \texttt{bad} support the energy does not
pay for the footprint, and the geometry of the footprint must be
used instead.

We say that a planar set $W$ is \emph{$(\epsilon,r)$-\texttt{fat}} if
$W=W^{\tt f}\cup W^\partial$, where $W^{\tt f}$ is a union of radius-$r$ balls
and $|W^\partial|\leq\epsilon|W|$. For a tiling pattern we apply this
definition to the region it covers. Geometric area and the number of
lozenges differ by a fixed factor, which is absorbed in absolute
constants below. All patterns in this subsection are planar;
\cref{sec:periodic-hmsos} extends the construction to the
torus, using balls in its flat metric, including noncontractible patterns.

The reason for working with reference footprints, rather
than with arbitrary subsets of the tiling having few lozenges of one
orientation, is the following. A footprint $P_\varphi(D)$ consists of the
reference dots in $D^\ddagger$, of which there are at most
$|D^\ddagger|\leq25|D|$, and of the reference wall plaquettes on the edges
with both endpoints in $D^\ddagger$, of which there are at most $50|D|$
edges, each carrying a complete stack of wall plaquettes, as many
as the absolute difference of the column magnetizations of $\varphi$ across
$e$. Thus a \texttt{bad} footprint, having at least $B|D|$
plaquettes, owes most of its size to a bounded number of complete
stacks, some of which must be long; this is the structural
information that the straight-tube construction below exploits.
For an arbitrary subpattern no such conclusion is available: the
number of dots in it does not bound the number of edges carrying its
arrows, since a subpattern may contain a few plaquettes from each of
many stacks. This is exactly the extra information supplied
by the sorting-plane defect support.

\begin{lemma}[A tube around a straight stack]\label{lem:straight-tube}
Let $J$ be a row of $\ell$ identical lozenges, consecutive along one
of their edge directions, in a full-plane tiling. Place a radius-$R$
Euclidean ball at the center of each lozenge of $J$, and let $U_R(J)$
be the tiling closure of their union: all tiles meeting that union
in positive area are included. For $R\geq2$ there is an
orientation $b$, determined by the orientation and row direction of
$J$, such that
\[
 \#\{b\text{-lozenges in }U_R(J)\}\leq CR^2\,.
\]
The neighborhood before taking tiling closure is a union of
radius-$R$ balls. There are at most six choices of the ordered pair
$(a,b)$, where $a$ is the orientation of the row.
\end{lemma}
\begin{proof}
Let $a$ be the orientation of the lozenges of $J$. Successive
lozenges of the row share edges parallel to the edge direction
\emph{other than} the direction of translation along the row.
Let $b$ be the other lozenge orientation having edges in this shared
direction. We use the usual
representation of a lozenge tiling by nonintersecting lattice paths
(see, e.g., \cite{Kenyon09}), with the two orientations $a,b$
as its steps: across each lozenge of orientation $a$ or $b$, join
the midpoints of its pair of edges in the direction common to those
two orientations. Each such edge midpoint is shared by exactly two
lozenges of orientations in $\{a,b\}$, so the continuation at every
midpoint is unique, and the resulting paths have no endpoints and do
not cross. The lozenges of the third orientation are not visited by
any path. Since consecutive lozenges of $J$ meet across these
shared-direction edges, the row is a segment of a single
path, consisting of $\ell$ consecutive $a$-steps.

Choose affine coordinates $(u,v)$ in which the $a$ and $b$ steps
are $(1,1)$ and $(1,-1)$, respectively. We view $u$ as a time
coordinate: every step increases $u$ by one. In the sheared
coordinate $q=(v-u)/2$, an $a$ step preserves $q$ and a $b$ step
decreases $q$ by one; thus $q$ is nonincreasing along every path,
and the number of $b$-lozenges visited by a path during a time
interval is exactly the decrease of $q$ over that interval. Translate
so that the segment $J$ lies on $q=0$, for $0\leq u\leq\ell$. Since
the change of coordinates is a fixed linear map, the union of the
radius-$R$ balls centered on $J$ is contained, for a suitable
absolute $C$, in the region $\{-CR\leq q\leq CR,\ -CR\leq u\leq
\ell+CR\}$. The two end pieces $\{-CR\leq u\leq0\}$ and
$\{\ell\leq u\leq\ell+CR\}$ of this region have area $O(R^2)$, so
they contain $O(R^2)$ lozenges, and it suffices to count the
$b$-lozenges in the two bands $0\leq q\leq CR$ and $-CR\leq q\leq0$
over the time interval $0\leq u\leq\ell$.

Consider a path meeting the upper band $0<q\leq CR$ at some time
$u_0\in[0,\ell]$. It cannot subsequently reach $q=0$ before time
$\ell$, since at $q=0$ it would meet the path containing $J$, and
distinct paths do not intersect. Nor can it leave through the upper
boundary $q=CR$, since $q$ is nonincreasing. Hence it remains in the
band up to time $\ell$ and crosses the right boundary $u=\ell$ of the
band. At a fixed time, distinct paths have distinct values of $q$,
so at most $CR$ paths cross that boundary inside the band. Each of
them has at most $CR$ decreasing steps while in the band, since $q$
decreases from at most $CR$ to at least $0$ there. The upper band
therefore contains at most $(CR)^2$ lozenges of orientation $b$.
In the lower band $-CR\leq q<0$ the same argument applies with time
reversed: going backward in time, $q$ is nondecreasing, a path in
the band cannot reach $q=0$ before time $0$ (again because it cannot
meet the path containing $J$) and cannot leave through $q=-CR$, so it
crosses the left boundary $u=0$ inside the band; there are at most
$CR$ such paths, each with at most $CR$ decreasing steps in the
band. The tiling closure adds only tiles meeting the union of balls,
all of which lie within a bounded distance of it; applying the band
argument with the width $CR$ increased by a fixed constant covers
these enlarged tiles as well. This proves the estimate.

Finally, the pair $(a,b)$ is determined by the orientation $a$ of
the row and its direction---one of the $2$ edge
directions of an $a$-lozenge; the $3$ tile orientations with their $2$ row directions each give at most $6$ pairs $(a,b)$.
\end{proof}

A \texttt{bad} footprint contains many long stacks, and the witness
will be the union of their tubes. The key to combining the tubes is
to form their unions separately for the six classes $(a,b)$.
Individual tubes can overlap arbitrarily, and a direct enumeration
of a union of tubes, tube by tube, would cost an exponential in the
total area of the tubes, which is not affordable. Instead we only
record, for each of the six classes, the union of the tubes in that
class, which is the tiling closure of a union of balls, together
with the positions of its $b$-lozenges. The class unions cover every
point at most six times, and the total number of these exceptional
tiles is bounded by $CR^2$ times the number of original stacks by
\cref{lem:straight-tube}, regardless of how the tubes overlap, so
that no area is charged twice. The following lemma shows that such
a union is determined by a small amount of information.

\begin{lemma}[Sparse-orientation entropy]\label{lem:sparse-entropy}
Fix $0<\delta\leq\frac12$. On any finite subset of the honeycomb
lattice with $n$ black vertices, the number of perfect matchings
using at most $\delta n$ edges of one specified orientation is at
most $\exp(2\delta\log(1/\delta)n)$.

Consequently, a family of supports with at most $e^{\eta n}$ rooted
choices carries at most $\exp[(\eta+2\delta\log(1/\delta))n]$ patterns with at
most $\delta n$ edges of that orientation. On a honeycomb torus the
same bounds hold up to a factor $2^w$, where $w$ is the number of
noncontractible cycles formed by the two other orientations.
\end{lemma}
\begin{proof}
We show that a perfect matching is determined by the set of its edges
of the specified orientation. Every black vertex of the honeycomb
lattice is incident to exactly one edge of each of the three
orientations. Hence the vertex set has at most $n$ edges of the
specified orientation, and a matching using at most $\delta n$ of them
is obtained by first choosing which of these edges it contains, for
which there are at most $\sum_{j\leq\delta n}\binom nj$
possibilities, and then completing the matching with edges of the two
other orientations on the vertices not yet covered. We claim that
this completion is unique. Removing all edges of one orientation from
the honeycomb lattice leaves every vertex with degree two, and the
resulting graph is a disjoint union of bi-infinite zigzag paths. Its
restriction to the finitely many uncovered vertices is therefore a
disjoint union of finite paths, and a finite path has at most one
perfect matching: the edge at an endpoint is forced, and removing it
leaves a shorter path. This proves the first assertion using the
well-known bound that $\sum_{j\leq\delta n}\binom nj\leq\exp(H_e(\delta)n)$,
where $H_e(\delta)=\delta\log(1/\delta)+(1-\delta)\log(\frac1{1-\delta})$
is the binary Shannon entropy, satisfying
$H_e(\delta)\leq2\delta\log(1/\delta)$ for $0<\delta\leq\frac12$.
Multiplying by the number of supports gives the second assertion.

On a honeycomb torus the argument is the same, except that the
restriction of the two-orientation graph to the uncovered vertices
may also contain noncontractible cycles of the torus, each of which has
exactly two perfect matchings. This accounts for the factor $2^w$.
\end{proof}

\begin{lemma}[Enumeration of class unions]\label{lem:six-unions}
There exist an absolute constant $C\geq1$ and, for each $l\geq2$,
a constant $1\leq C_l\leq Cl^2$ depending only on $l$, such that
the following holds.
Suppose a connected tiling pattern $W$, with $s$ tiles, is the union
of at most six patterns $U_i$ and a remainder of at most $\delta s$
tiles. Each $U_i$ is the tiling closure of a finite union of
radius-$R$ balls. For each class fix an orientation $b_i$, and suppose
\[
 \sum_i\#\{b_i\text{-tiles in }U_i\}\leq\kappa s\,.
\]
If $2\leq l\leq R/C$, $\delta l^2\leq1$ and $0\leq\kappa\leq1$,
the number of patterns of size $s$, up to translation, is
at most
\begin{equation}\label{eq:six-union-count}
 C_l\exp\big(C\big(l^{-2}+l/R+\delta l^2
                   +\kappa\log(e/\kappa)\big)s\big)\,,
\end{equation}
where $\kappa\log(e/\kappa)$ is interpreted as zero when $\kappa=0$.
The prefactor $C_l$ accounts for the grid offsets and the choices of
the six orientations. No bound on the multiplicity of the individual balls or tubes
is assumed.
\end{lemma}
\begin{proof}
We first bound the boundary layer of a union $A$ of radius-$R$
balls. Write $A^{-t}$ for its points at distance at least $t$ from
its complement, and $A^{+t}$ for its $t$-neighborhood. Partition
$A$ into nearest-center Voronoi pieces and contract each toward its
center by $1-t/R$. By convexity, it stays in its Voronoi cell,
and it lies in the concentric radius-$(R-t)$ ball, hence in
$A^{-t}$. The images are disjoint, so
$|A^{-t}|\geq(1-t/R)^2|A|$. Similarly, contracting the Voronoi
pieces of $A^{+t}$ by $R/(R+t)$ puts them in $A$. Consequently,
\begin{equation}\label{eq:ball-union-boundary-layer}
 |A\setminus A^{-t}|\leq(2t/R)|A|,\qquad
 |A^{+t}|\leq(1+t/R)^2|A|\,.
\end{equation}
The first bound holds for $0\leq t\leq R$, the second for all
$t\geq0$. For $t\leq R/2$, both boundary layers have area
$O(t|A|/R)$; in particular, tiling closure adds $O(|A|/R)$ area,
since tiles have bounded diameter.

We now encode the supports, i.e., their elementary triangles.
Use a grid of parallelograms with side length $\floor{l}$ along
two lattice directions. Every cell meeting a class lies in the
$Cl$-neighborhood of its union of balls, so
\cref{eq:ball-union-boundary-layer}, with $t=Cl\leq R/2$, and
$\sum_i|U_i|\leq6s$
bound the number of such cells by $Cs/l^2$. The remainder adds at
most $C\delta s\leq Cs/l^2$ cells. Since $W$ is connected, these
cells form a lattice animal under corner adjacency. Record its
shape up to translation and,
for each cell, whether each class is full, empty or mixed, and
whether the remainder meets it. The bounded number of neighbors
and labels per cell gives at most $C_l\exp(Cs/l^2)$ choices,
including the $O(l^2)$ grid offsets and the fixed orientation
choices. Only mixed cells and cells meeting the remainder need
finer information. For each class, its mixed cells lie in the
$Cl$-boundary layer of its union of balls, since tiles have bounded
diameter; together these cells contain at most $C(l/R)s$
triangles. Cells meeting the remainder contain at most
$C\delta l^2s$ triangles. Record the exact membership of these
triangles in each class and in the uncovered remainder, at cost
$\exp(C(l/R+\delta l^2)s)$. This determines all the supports.

For each fixed class support, \cref{lem:sparse-entropy} says that
its matching is determined by its $b_i$-edges. Across the six
classes there are at most $6s$ candidate edges and at most
$\kappa s$ chosen ones. The remainder has at most $3^{\delta s}$
matchings, by choosing one edge at each black vertex. Thus the
binomial estimate in the proof of that lemma gives at most
\[
 3^{\delta s}\sum_{j\leq\kappa s}\binom{6s}{j}
 \leq\exp\big(C(\delta+\kappa\log(e/\kappa))s\big)
\]
matching choices. Requiring agreement on overlapping classes can
only reduce this count. These matchings determine $W$; combining
their count with the support count and using $\delta\leq\delta l^2$
proves \cref{eq:six-union-count}, with $C_l\leq Cl^2$.
\end{proof}

\begin{proposition}[Geometric witnesses from stack covers]\label{prop:stack-witness}
Fix $C_0\geq1$. Let $T\subset\varphi$ be a finite connected pattern
in a full-plane tiling, expressed as the union of at most $C_0m$
singleton tiles and straight rows as in \cref{lem:straight-tube},
for some $m\geq1$, each tile
belonging to at most $C_0$ pieces. There are $C,c>0$, depending
only on $C_0$, such that if $R\geq2$, $L\geq CR^2$, $B\geq CL$
and $|T|\geq Bm$,
this cover determines a connected enlargement $T\subset W\subset\varphi$ with
\[
 \diam W\leq\diam T+CR,\qquad |W|\leq CR|T|,\qquad
 W\text{ is }(C/R+CL/B,R/2)\text{-\texttt{fat}}\,.
\]
The patterns can be chosen in a family satisfying
\cref{eq:six-union-count}, for every admissible $l$ in that lemma, with
\[
                  \delta=CL/B,\qquad \kappa=CR^2/L\,.
\]
More precisely, $W$ is the union of at most six class unions of
radius-$R$ tubes and a remainder of at most $\delta|T|$ tiles;
the sum of their exceptional-orientation counts is at most
$\kappa|T|$.
If $\delta\leq c/R^2$, then also
\begin{equation}\label{eq:witness-area-diameter}
                     |W|\geq cR\diam W\geq cR\diam T\,.
\end{equation}
In particular one may choose $R=R_\beta=\beta^{1/100}$,
$L=\beta^{1/10}$ and $B=\beta^{1/2}$, rounded to integers.
There is then a family $\cW_\beta$ whose pattern count of
size $s$, modulo translations, is $C_\beta e^{\eta_\beta s}$,
where $\eta_\beta\to0$, and whose members are
$(\epsilon_\beta,R_\beta/2)$-\texttt{fat} with $\epsilon_\beta\to0$ and
$R_\beta\to\infty$.
\end{proposition}
\begin{proof}
Retain all rows $J_1,\ldots,J_k$ of length at least $L$, and let
$T_0$ consist of the tiles not covered by these rows. Each tile of
$T_0$ belongs to a discarded piece of fewer than $L$ tiles, whereas
the total length of the retained rows is bounded by the cover's
multiplicity. Thus
\[
 |T_0|\leq C_0Lm\leq(C_0L/B)|T|,\qquad
 kL\leq\sum_j|J_j|\leq C_0|T|\,.
\]
Around each $J_j$ take the tube $U_R(J_j)$ of
\cref{lem:straight-tube}, and let $W$ be their union with $T_0$.
This uses every long row and hence makes no auxiliary choices.
Since each tube contains its row, $W$ contains $T$, and it is
connected: each tube
is connected and meets the connected pattern $T$, and $T_0\subset T$.
Every point of $W$ is within distance $CR$ of $T$, so
$\diam W\leq\diam T+CR$. The tube around a row of $\ell$ lozenges
has at most $C(R\ell+R^2)$ tiles, and $R^2\leq R|J_j|$ because
$|J_j|\geq L\geq CR^2$, so
\[
 |W|\leq C\sum_{j=1}^k\big(R|J_j|+R^2\big)+|T_0|
     \leq CR\sum_{j=1}^k|J_j|+|T_0|\leq CR|T|\,.
\]

Group the tubes by their ordered pair $(a,b)$ of
\cref{lem:straight-tube}, and call the at most six class unions
$U_{a,b}$. Their geometric cores, the unions of the radius-$R$ balls
centered on the rows of the class, are unions of radius-$R$
balls, and the tiling closure of a union of sets is the union of
their tiling closures, so each $U_{a,b}$ is the tiling closure of
its geometric core, as \cref{lem:six-unions} requires. The total
number of exceptional-orientation tiles in the classes is at most
\[
 \sum_{a,b}\#\{b\text{-tiles of }U_{a,b}\}
 \leq\sum_{j=1}^k\#\{b_j\text{-tiles of }U_R(J_j)\}
 \leq CR^2 k\leq (CR^2/L)|T|\leq\kappa|W|\,,
\]
where $b_j$ is the exceptional orientation of the tube of $J_j$:
a $b$-tile of $U_{a,b}$ belongs to at least one tube of that class,
in which it is exceptional. This estimate remains valid with
arbitrarily many overlaps between the tubes. It is
the absolute exceptional-tile bound, together with the multiplicity
six of the class unions, that permits \cref{lem:six-unions} to
be applied with the displayed $\kappa$ and with $\delta=CL/B$, the
remainder being $T_0$. Counting the original tubes individually is
unnecessary.

For the \texttt{fat} assertion, let $A$ be the union of all the geometric cores. It
is itself a union of radius-$R$ balls. Its tiling closure, which is
the union of the tubes, differs from $A$ by an area at most $C|A|/R$,
by \cref{eq:ball-union-boundary-layer}. The remainder $T_0$
contributes at most $\delta|W|$ further tiles. Every radius-$R$ ball
is a union of the radius-$(R/2)$ balls contained in it, so $W$ is
$(C/R+\delta,R/2)$-\texttt{fat}, proving the \texttt{fat} assertion (with a
fixed adjustment of the radii if sets are recorded as unions of
elementary triangles rather than as geometric regions).

For \cref{eq:witness-area-diameter}, use a grid of squares of side
$R$. The number of its squares meeting $A$ is at most $C|A|/R^2\leq
C|W|/R^2$, by the outer-layer estimate in
\cref{eq:ball-union-boundary-layer} with $t=CR$; the remainder adds
at most $C\delta|W|$ squares. Include a fixed number of neighboring
layers of squares to cover the tiling closure of $A$; since tiles
have bounded diameter and $R\geq2$, this increases the count by only
an absolute factor. The resulting squares form a connected animal
covering $W$, so $\diam W$ is at most $CR$ times their
number, whence
\[
 \diam T\leq\diam W\leq
 CR\big(|W|/R^2+\delta|W|\big)\leq C|W|/R\,,
\]
where the last step uses $\delta\leq c/R^2$. Note that there is at
least one long row: otherwise $|T|=|T_0|\leq C_0Lm<Bm$,
contrary to the hypothesis, after choosing $C>C_0$ in $B\geq CL$.
Hence $|W|\geq L$,
so the possible additive error of order $R$ in this grid estimate is
absorbed by $C|W|/R$.
Finally take $l=\floor{\sqrt R}$ in
\cref{lem:six-unions}. With the displayed choices of $R,L,B$ as
powers of $\beta$, we have $l^{-2}+l/R\to0$, $\delta l^2\leq
C\beta^{1/10-1/2+1/100}\to0$ and $\kappa\leq C\beta^{2/100-1/10}\to0$,
so every term in the entropy exponent of \cref{eq:six-union-count}
tends to zero, while $W$ is $(C/R+\delta,R/2)$-\texttt{fat} with
$C/R+\delta\to0$ and $R/2\to\infty$. This proves the
last assertion.
\end{proof}

For the \HMSOS footprint $T=\Upsilon(P_\varphi(D))$, the cover
consists of at most $25|D|$ reference dots and $50|D|$ complete
wall stacks, as explained after \cref{def:footprint}. A stack
projects to a straight row of identical lozenges, and these pieces
are disjoint. Thus \cref{prop:stack-witness} applies with $m=|D|$
and $C_0=75$ whenever $|P_\varphi(D)|\geq B|D|$; write
$W_\varphi(D)$ for its enlargement. The same construction will
apply to the bounded-multiplicity covers of the longer histories.

\subsubsection{Periodic HMSOS fibers and noncontractible witnesses}\label{sec:periodic-hmsos}

We now allow noncontractible supports, still for a single \HMSOS sorting
map. The planar analysis above applies verbatim to supports that
are small compared with the torus, but the fiber over a periodic
reference surface also contains defects whose supports wind around
the torus, and even the whole transverse torus. Fix positive
integers $n,\Theta_1,\Theta_2$ and the physical period lattice
\[
 \Gamma=\langle(-n,0,\Theta_1),(0,-n,\Theta_2)\rangle\subset\Z^3\,.
\]
Thus the slope is $(\Theta_1/n,\Theta_2/n)$, as per \cref{def:3dIsing}.
Sorting in direction $x_1$ gives the transverse torus
\[
 \T_\mathrm{s}=\Z^2/\Gamma_\mathrm{s},\qquad
 \Gamma_\mathrm{s}=\langle(0,\Theta_1),(-n,\Theta_2)\rangle\,,
\]
obtained by projecting the two period vectors along $x_1$, and the
diagonal tiling torus is
$\T_\mathrm{d}=\cP_{111}/\Upsilon\Gamma$.
Write $h^{(1)}_\varphi(v)$ for the column magnetization of
$\varphi$ along the $x_1$-fiber over the transverse square $v$, so
that $\varphi$ is the graph of $h^{(1)}_\varphi$ in the sorting
coordinates. This function is not periodic on the transverse cover
but obeys
\[
 h^{(1)}_\varphi(v+(0,\Theta_1))=h^{(1)}_\varphi(v)-n,\qquad
 h^{(1)}_\varphi(v+(-n,\Theta_2))=h^{(1)}_\varphi(v)\,,
\]
since the first period vector has a component $-n$ along the sorting
direction, which shifts the spin sequence of the column and
hence its column magnetization by $-n$, whereas the second has
none. Assume throughout that $(\Theta_1/n,\Theta_2/n)$ ranges over a fixed compact
subset of $(0,\infty)^2$. Both period lattices then have shortest
nonzero vector at least $c_*n$ and fundamental parallelograms of
diameter at most $C_*n$. We take $n$ large enough that every
radius-five lattice neighborhood embeds in its covering plane, so
that the local definitions of $D^\dagger$, $D^\ddagger$ and $P_\varphi(D)$ make sense
on the quotient.
All set sizes in the following statements count distinct objects
in the quotient, rather than the objects in a chosen cut-open lift.

\begin{proposition}[The complete periodic fiber]\label{prop:periodic-fiber}
For a fixed monotone periodic reference surface, its \HMSOS fiber
is exactly the hard-core gas on all nonempty $*$-connected subsets
$D\subset\T_\mathrm{s}$, including noncontractible subsets and the
whole torus. The enumeration bound \cref{eq:correct-count}, the
activity bound \cref{eq:activity-bound}, and the a priori bound of
\cref{prop:positive-budget} hold with the same absolute
constants. In particular, writing $N_0=n(n+\Theta_1+\Theta_2)$,
\begin{equation}\label{eq:periodic-fiber-pressure}
     1\leq Z_{\beta,\varphi}\leq\exp(Ce^{-\beta/64}N_0)\,.
\end{equation}
The unrestricted fiber is the total-variation limit of its
normalized fiber truncations at distance $J$ as $J\to\infty$. The footprint overlap
estimate of \cref{lem:footprint-overlap} also holds on the quotient,
uniformly in the available supports.
\end{proposition}
\begin{proof}
The fiber consists of binary $x_1$-columns with the prescribed charges
and period offsets, subject to inclusion of neighboring occupied
sets in directions $x_2,x_3$. These are exactly the local
monotonicity constraints allowed in \cref{lem:quotient-one-step}.
They also ensure admissibility: monotonicity in $x_3$ gives a graph
surface, with both phases connected to their ends. Its vertical
ends follow from the sorting ends and
$\cI(x_1,x_2,x_3+k\Theta_1)=\cI(x_1+kn,x_2,x_3)$.
Positive windings exclude a period parallel to $x_1$, so the lemma
gives factorization, enumeration, activities and the limit $J\to\infty$.

The neighborhood argument of \cref{lem:footprint-overlap} uses at most the
$25$ squares within two transverse steps of a fixed square. It remains
valid on the quotient, where identifications only reduce this number.
Finally the numbers of reference plaquettes of the three normals are
$n\Theta_1,n\Theta_2,n^2$, the indices of the corresponding projected period
lattices. Apply the a priori bound with the whole reference surface,
of area $N_0=n(n+\Theta_1+\Theta_2)$, to obtain \cref{eq:periodic-fiber-pressure}.
\end{proof}

\begin{lemma}[Periodic class-union enumeration]\label{lem:periodic-unions}
There exist absolute constants $C\geq1$, $c>0$ and, for each $l\geq2$,
a constant $1\leq C_l\leq Cl^2$, independent of the torus and all
parameters other than $l$, such that the following holds.
Let $\lambda$ be the shortest period length of $\T_\mathrm{d}$.
For $2\leq l\leq R/C$ and $R\leq c\lambda$, the rooted count in
\cref{lem:six-unions}, now on this torus, is at most
\begin{equation}\label{eq:periodic-union-count}
 C_l\exp\big(
 C\big(l^{-2}+l/R+\delta l^2+\lambda^{-1}
                 +\kappa\log(e/\kappa)\big)s\big)\,.
\end{equation}
Here the root is a specified elementary triangle, and the hypotheses
$\delta l^2\leq1$, $0\leq\kappa\leq1$ are unchanged. The union
and the classes may be noncontractible.
\end{lemma}
\begin{proof}
We follow the proof of \cref{lem:six-unions} and indicate the two
places where the torus topology enters. The boundary-layer estimate
\cref{eq:ball-union-boundary-layer} holds on the torus: lift the
centers periodically and contract their nearest-center Voronoi
pieces in the plane. Each contraction stays in its Voronoi cell and
commutes with the period translations, so the contracted pieces
remain disjoint on the quotient. The condition $R+t<\lambda/2$
ensures that each relevant ball is embedded. In place of the planar grid
use a fixed maximal $l$-separated net on the torus (which need not
carry a lattice structure). Its Voronoi cells have area at least
$cl^2$, diameter at most $Cl$, and bounded adjacency degree, which
is all that the support encoding of that lemma used. That encoding
therefore applies verbatim, with cost
$C_l\exp(C(l^{-2}+l/R+\delta l^2)s)$, where the root is now a
specified elementary triangle rather than a translation class.
The root meets only a bounded number of net cells, so the prefactor
can again be bounded by $Cl^2$, uniformly in the torus.

Only the two-orientation matching count of \cref{lem:sparse-entropy}
changes. After the exceptional-orientation edges of a class are
specified and deleted together with their endpoints, the remaining
graph consists of paths and cycles. Paths have at most one perfect
matching, whereas cycles have two. Every such cycle is
noncontractible, since its lift to the plane lies on a bi-infinite
zigzag path, which contains no cycle, and a noncontractible cycle
has at least $c\lambda$ vertices. The six classes have at most
$12s$ vertices altogether, so there are at most $Cs/\lambda$ cycles
and the additional choices cost at most $2^{Cs/\lambda}$. The
exceptional edges and the remainder still cost
$\exp(C\kappa\log(e/\kappa)s)$ and $3^{\delta s}$, respectively. Combining
these bounds proves \cref{eq:periodic-union-count}.
\end{proof}

For the remaining one-step estimates, fix
\begin{equation}\label{eq:history-scales}
 B=B_\beta=\ceil{\sqrt\beta},\qquad
 R=R_\beta=\ceil{\beta^{1/100}},\qquad
 L=\ceil{\beta^{1/10}},\qquad d_\beta=\beta^{-1/10}\,.
\end{equation}
These satisfy the hypotheses $R\geq2$, $L\geq CR^2$ and $B\geq CL$
of \cref{prop:stack-witness} for all large $\beta$. The parameter $d_\beta$ will be the threshold in the
\texttt{thin}/\texttt{fat} split.

\begin{proposition}[Periodic geometric witnesses]\label{prop:periodic-witness}
Under the parameter and cover hypotheses of \cref{prop:stack-witness},
its construction and geometric estimates hold on the diagonal torus
if $R\leq c\lambda$ and every row in
the cover is the injective projection of a finite straight row
in the lifted tiling. Diameter is measured in the torus metric,
and the entropy bound is \cref{eq:periodic-union-count}.
In particular this applies to $T=\Upsilon(P_\varphi(D))$ for
every connected $D\subset\T_\mathrm{s}$ with
$|P_\varphi(D)|\geq B|D|$, defining $D^\ddagger$ and its footprint
on the quotient.
For the choices $B_\beta,R_\beta,L_\beta$ used in the mixed
expansion, if $D$ contains a noncontractible loop and
$|P_\varphi(D)|\geq B_\beta|D|$, its witness has
\begin{equation}\label{eq:wrapping-witness-area}
                        |W_\varphi(D)|\geq cR_\beta n\,.
\end{equation}
\end{proposition}
\begin{proof}
Apply the construction of \cref{prop:stack-witness}, lifting each
retained row, taking its planar tube, and projecting the tube back.
A projected tile of the exceptional orientation has at least one
exceptional preimage in
the planar tube, so the exceptional-orientation count of the
projected tube is still at most $CR^2$. This argument does not
require the tube itself to project injectively, and indeed it need
not when $R$ is comparable to $\lambda$.

Injectivity of each row and bounded multiplicity of the cover give
the same bounds $\sum_j|J_j|\leq C_0|T|$ and
$kL\leq C_0|T|$ as in the plane. Hence the remainder, the
six exceptional-orientation counts and $|W|$ obey the same estimates.
The torus boundary-layer estimate from the proof of
\cref{lem:periodic-unions} proves that $W$ is \texttt{fat} with the same parameters, with balls in
the flat metric of the torus. For the area--diameter estimate, use
the net construction of that proof at scale $R$ instead of the
planar grid. The number of cells meeting the witness is at most
$C|W|/R^2+C\delta|W|$, and the connected adjacency graph of these
cells bounds the torus diameter of $W$ by $CR$ times this number.
When $\delta\leq c/R^2$ this gives $\diam W\leq C|W|/R$, as before.
\Cref{lem:periodic-unions} gives the entropy bound.

For an \HMSOS footprint, the same local rule gives at most
$25|D|$ dots and $50|D|$ complete wall stacks on the quotient.
Each stack lifts to a finite straight row. Two of its plaquettes
could represent the same quotient face only if a nonzero physical
period were parallel to $x_1$, which is excluded by the rank-two
projection $\Gamma_\mathrm{s}$ of $\Gamma$. Thus the preceding
construction applies with $m=|D|$.
Finally, suppose that $D$ contains a noncontractible loop of
$\T_\mathrm{s}$, and let us show that its witness is large. Such a
loop can be realized on the reference footprint by joining
consecutive reference dots along the intervening complete walls,
where a diagonal step of the $*$-path in $D$ is replaced by two
coordinate steps inside $D^\ddagger$. Its lift to the cover has a nonzero
vector of $\Gamma$ as displacement, whose diagonal projection
$\Upsilon\gamma$ is nonzero since $\Upsilon$ is injective on
$\Gamma$. The resulting footprint loop is therefore noncontractible
on $\T_\mathrm{d}$. A connected set containing such a loop has
diameter at least $\lambda/4$: a set of smaller diameter lies in an
embedded ball of radius $\lambda/4$, which is contractible, so every
loop in it is contractible. Since $\lambda\geq c_*n$, the
area--diameter bound $|W|\geq cR\diam W$ proves
\cref{eq:wrapping-witness-area}; its hypothesis
$CL_\beta/B_\beta\leq c/R_\beta^2$ holds for all sufficiently large
$\beta$ by \cref{eq:history-scales}.
\end{proof}

These statements settle the periodic fiber and geometric estimates
for \HMSOS. What remains for the torus is of a different nature:
the weights $Z_{\beta,\varphi}$ are defined only on the tilings of the
prescribed sector, whereas the \RG analysis requires a weight on
every tiling of the torus. \Cref{sec:canonical-extension} proves
locality with prescribed monodromy, constructs such a canonical
extension, and gives its exact fixed-count Fourier projection. The
uniform first integration on the torus is proved in
\cref{thm:torus-first-step}. \Cref{sec:fixed-sector-limit} later
bounds the projected partition function after the remaining
integrations. The probabilities of noncontractible witnesses under the
interacting law remain a separate task, taken up in
\cref{sec:interacting-witnesses}.

\subsubsection{The mixed positive expansion and \texorpdfstring{\RG}{RG} coefficients}
\label{sec:hmsos-rg-input}
The positive expansion below is combinatorial. Its Grassmann
representation is given in \cref{sec:model-grassmann}, and the
analytic estimates begin in \cref{sec:fat-determinant}.

For \HMSOS we work directly with the positive fiber partition
function $Z_{\beta,\varphi}$ of \cref{prop:positive-budget}, which is the weight
of the tiling $\varphi$ in the image law of
\cref{thm:hmsos-rg-input}, and we include both \texttt{good} and
\texttt{bad} defects in one expansion. The goal of this subsection
is to rewrite $Z_{\beta,\varphi}$ in the form required by
\cref{thm:rg-typical}, namely as a hard-core gas of patterns of
$\varphi$ with coefficients that depend only on the pattern, and to
prove the \texttt{thin}/\texttt{fat} bounds and the entropy bounds
on these coefficients. Throughout this construction,
$|P|$ counts lozenges, whereas the
Grassmann support $S(P)$ consists of their $2|P|$ endpoints.

With the parameters fixed in \cref{eq:history-scales}, define for a support $D$
its enlarged pattern by
\[
 Q_\varphi(D)=
 \begin{cases}
   \Upsilon(P_\varphi(D)),&|P_\varphi(D)|<B|D|,\\
   W_\varphi(D),&|P_\varphi(D)|\geq B|D|.
 \end{cases}\,.
\]
Thus a \texttt{good} defect is recorded by its footprint, and a
\texttt{bad} one by the witness $W_\varphi(D)\supset\Upsilon(P_\varphi(D))$ of
\cref{prop:stack-witness}, which is the deterministic union of the
tubes around its long stacks and of the removed pieces; since the
construction in that proof makes no arbitrary choices, $W_\varphi(D)$ is
a function of $D$ and of $\varphi$ near $D$.
We join two defects exactly when their enlarged patterns share a
lozenge, and we group a compatible family of defects into the
connected components of this relation. Thus distinct components
have disjoint triangle supports, but they may be adjacent, in the
sense that their supports share edges or vertices of the triangular
lattice. This convention is dictated by the Grassmann identity of
\cref{sec:grassmann-input}: there, each pattern is represented by a
monomial in the fields at the honeycomb vertices dual to its triangles.
The product of two such monomials vanishes exactly when the patterns
share a triangle, so a hard-core gas of patterns with
disjoint supports is precisely what the exponential of the
interaction produces.

\begin{lemma}[Locality of the enlarged-pattern coefficients]\label{lem:mixed-locality}
For a finite connected pattern $P$ in a monotone tiling,
let
\begin{equation}\label{eq:mixed-coefficient}
 c_\beta(P)=
 \sum_{\substack{\cD\ne\varnothing\text{ compatible}\\
  \text{overlap graph of }\{Q_\varphi(D):D\in\cD\}\text{ is connected}\\
                       \bigcup_{D\in\cD}Q_\varphi(D)=P}}
                         \prod_{D\in\cD}z_{\beta,\varphi}(D)\,.
\end{equation}
This coefficient is determined by $P$, independently of the
tiling containing it. With any finite, locally specified restriction on available
defects, its restricted counterpart satisfies the exact identity
\begin{equation}\label{eq:mixed-gas}
 Z_{\beta,\varphi}=
 \sum_{\{P_i\}\text{ with disjoint triangle supports}}
                           \prod_i c_\beta(P_i),
 \qquad P_i\subset\varphi\,.
\end{equation}
All coefficients are nonnegative. The bounds below hold uniformly
under additional restrictions on the available supports; the
unrestricted coefficients are translation covariant.
\end{lemma}
\begin{proof}
The proof has three parts. We first show that the activity of a
defect is a function of its footprint, then that the set of
summands in \cref{eq:mixed-coefficient} can be read off from $P$
alone, and finally we derive the identity \cref{eq:mixed-gas} from
the hard-core gas of \cref{lem:polymer-factorization} by grouping
defects into components.

\begin{enumerate}[label=\textbf{Step~\arabic*.}, ref=\arabic*, wide=0pt, itemsep=1ex, listparindent=\parindent, parsep=0pt]
\item\label[step]{st:mixed-locality-activity}
(\emph{The activity depends only on the footprint.})
A connected lozenge pattern which extends to a monotone surface
specifies its lift to $\Z^3$ up to a translation in the kernel of
the diagonal projection $\Upsilon$: the lift of one lozenge fixes
the lifts of its neighbors, and so on along the connected pattern.
In monotone lattice coordinates the ambiguity is a multiple of
$(1,1,1)$, which translates the transverse origin and adds a common constant to the column magnetizations of all the columns. Now fix a candidate
support $D$. By \cref{def:footprint}, the footprint $P_\varphi(D)$
specifies the reference columns on $D^\ddagger$ and their relative
heights, that is, their charges up to the common constant. The
activity $z_{\beta,\varphi}(D)$ of \cref{cor:activities}
sums $e^{-\beta E}$ over the $D$-defects, and such a defect is
constrained only by charge
conservation on the columns of $D$ and by the two \HMSOS
monotonicity constraints with respect to the neighboring columns,
which are nearest-neighbor constraints; its excess energy
\cref{eq:exact-excess} is determined by the modified columns and
the reference columns adjacent to them. All of these are
determined by the charges on $D^\ddagger$, and a common height
translation of all the columns changes none of them. Consequently
$P_\varphi(D)$ determines $z_{\beta,\varphi}(D)$ and $|D|$, hence whether $D$ is \texttt{good} or \texttt{bad}.

\item\label[step]{st:mixed-locality-test}
(\emph{The summands are determined by $P$.})
The candidate supports $D$ with $Q_\varphi(D)\subset P$ can be enumerated
from $P$ alone, as follows. Lift $P$ to the monotone surface as in
\cref{st:mixed-locality-activity}; its dots then have well-defined
transverse coordinates, and a candidate $D$ must be a
$*$-connected set of these coordinates, since every square of $D$
carries a dot of $P_\varphi(D)\subset P$. Retain a candidate only if $P$
contains the dot of every square of $D^\ddagger$ and every wall plaquette
between the indicated adjacent column heights, for every edge with
both endpoints in $D^\ddagger$. This finite test is exactly the condition
$P_\varphi(D)\subset P$, and when it succeeds it also recovers the
footprint $P_\varphi(D)$, hence by
\cref{st:mixed-locality-activity} the activity and the
\texttt{good}/\texttt{bad} status of $D$. For a \texttt{good}
candidate $Q_\varphi(D)=\Upsilon(P_\varphi(D))$ is then known. For a
\texttt{bad} candidate the centers of all the tube balls are now
fixed by its complete wall stacks, which are part of its footprint.
The witness $W_\varphi(D)$ is the tiling closure of the union of these
balls together with the removed pieces, and the removed pieces are
already in $P_\varphi(D)\subset P$. The tiling closure is contained in
$P$ exactly when the ball cores are contained in the geometric
region covered by $P$: indeed, the interiors of different tiles of
the tiling are disjoint, and a tile is included in the
closure only if it meets the union of balls in positive area, so
the closure consists precisely of the tiles of $P$ having such an
intersection whenever the cores lie in the region covered by $P$,
and otherwise contains a tile outside $P$. Thus both the condition
$Q_\varphi(D)\subset P$ and the resulting pattern $Q_\varphi(D)$ are determined
by $P$, and so are the overlap-connectivity and union tests in
\cref{eq:mixed-coefficient}. There are finitely many candidates
because $P$ has finitely many dots. This proves that $c_\beta(P)$
is determined by $P$, including when the available defects are
subject to a finite, locally specified restriction, since such a
restriction is again a test on the candidates. Translation
covariance of the unrestricted coefficients follows since every
step of the construction commutes with lattice translations.

\item\label[step]{st:mixed-locality-gas}
(\emph{The pattern expansion.})
By \cref{lem:polymer-factorization} and the local restriction in
\cref{def:footprint},
\[
 Z_{\beta,\varphi}=\sum_{\cD\text{ compatible}}\prod_{D\in\cD}z_{\beta,\varphi}(D)\,,
\]
the sum running over compatible families of available supports. Every
compatible defect family $\cD$ has a unique decomposition
$\cD=\cD_1\cup\cdots\cup\cD_k$ into the components of its
enlarged-pattern overlap graph, and the unions
$P_i=\bigcup_{D\in\cD_i}Q_\varphi(D)$ are connected patterns with
pairwise disjoint triangle supports. This maps the summand for
$\cD$ to the collection $\{P_i\}$ together with, for each $i$, a
summand of \cref{eq:mixed-coefficient} for $P_i$. Conversely, let
$\{P_i\}$ be a collection of patterns of $\varphi$ with pairwise
disjoint triangle supports and let $\cD_i$ be, for each $i$, a
compatible family with connected overlap graph and union $P_i$. We
claim that $\cD=\bigcup_i\cD_i$ is compatible and that the $\cD_i$
are exactly its components. Indeed, if two supports $D\in\cD_i$ and
$D'\in\cD_j$ with $i\neq j$ were incompatible, their two-square
neighborhoods $D^\ddagger,(D')^\ddagger$ would share a square and hence their reference
footprints would share the dot of that square; their enlarged
patterns, which contain the footprints, would then overlap,
contradicting the disjointness of $P_i$ and $P_j$. The same
observation shows that no edge of the overlap graph joins $\cD_i$
to $\cD_j$, so the components of $\cD$ are the $\cD_i$. Hence the
two maps are inverse bijections, and summing independently inside
the components gives
\[
 Z_{\beta,\varphi}
 =\sum_{\{P_i\}}\prod_i
   \sum_{\substack{\cD_i\text{ compatible, connected}\\
                   \bigcup_{D\in\cD_i}Q_\varphi(D)=P_i}}
     \prod_{D\in\cD_i}z_{\beta,\varphi}(D)
 =\sum_{\{P_i\}}\prod_ic_\beta(P_i)\,,
\]
which is \cref{eq:mixed-gas}. Note that the components of a
compatible family can be adjacent: two enlarged patterns may share
an edge or a vertex of the triangular lattice without sharing a
lozenge. Adjacent but disjoint enlarged patterns must therefore be
allowed in that identity, as stated. Finally, the coefficients are
nonnegative as sums of products of nonnegative activities.
\qedhere
\end{enumerate}
\end{proof}

\begin{proposition}[Bounds on mixed coefficients]\label{prop:mixed-envelopes}
There are nonnegative coefficients
$c_\beta(P)=c^\mathrm{t}_\beta(P)+c^\mathrm{f}_\beta(P)$ and numbers
$\epsilon_\beta\to0$ such that, for all sufficiently large $\beta$,
\begin{align}
 c^\mathrm{t}_\beta(P)&\leq e^{-c\beta^{2/5}|P|},
                                      \label{eq:mixed-thin-bound}\\
 c^\mathrm{f}_\beta(P)&\leq e^{Ce^{-\beta/64}|P|}\,.
                                      \label{eq:mixed-fat-bound}
\end{align}
There are at most $C_0^s$ rooted patterns with
$c^\mathrm{t}_\beta(P)\ne0$ and $|P|=s$. Every pattern with
$c^\mathrm{f}_\beta(P)\ne0$ has $|P|\geq B_\beta$, is
$(\epsilon_\beta,R_\beta/2)$-\texttt{fat}, and belongs to a family with
at most $e^{\epsilon_\beta s}$ rooted patterns of size $s$.
The coefficients sum all defect decompositions of a pattern;
those decompositions are not additional pattern labels.
The same bounds hold in every fixed reference torus tiling, with
rooted \texttt{fat}-pattern entropy rate $\epsilon_\beta+C/\lambda$.
\end{proposition}
\begin{proof}
The split is defined summand by summand in
\cref{eq:mixed-coefficient}, according to whether the \texttt{good}
defects account for a definite fraction of the pattern. For a
summand with defect family $\cD$ put $s=|P|$ and
$q(\cD)=\sum_{D\text{ \texttt{good}}}|P_\varphi(D)|$, counting multiplicity
when footprints overlap. Assign the summand to $c^\mathrm{t}$ if
$q(\cD)\geq d_\beta s$, and to $c^\mathrm{f}$ otherwise. This is an
exhaustive disjoint split, so $c_\beta=c^\mathrm{t}_\beta
+c^\mathrm{f}_\beta$ with both parts nonnegative. The superscripts
indicate the estimate used, not the geometry of $P$; a pattern
contributing to the first part can itself be geometrically
\texttt{fat}, and a pattern can contribute to both parts.

\subsubsection*{The \texttt{thin} bound} The idea is that a
\texttt{good} defect has excess at least $2|D|>2|P_\varphi(D)|/B$, so its
Boltzmann factor can pay a small exponential in the size of its
footprint. Let $\lambda=\beta/(4B)$. For a \texttt{good} support and
any internal excess $E\geq2|D|$ (recall from
\cref{prop:correct-count} that smaller excesses do not occur), we
have $\lambda |P_\varphi(D)|<\lambda B|D|=\beta|D|/4\leq\beta E/8$. Hence,
multiplying each term of $z_{\beta,\varphi}(D)=\sum_EN_\varphi(D,E)e^{-\beta E}$
by $e^{\lambda |P_\varphi(D)|}\leq e^{\beta E/8}$,
\[
 e^{\lambda |P_\varphi(D)|}z_{\beta,\varphi}(D)
             \leq z_{7\beta/8,\varphi}(D)\leq z_{\beta/2,\varphi}(D)\,.
\]
For a \texttt{bad} support simply use the monotonicity
$z_{\beta,\varphi}(D)\leq z_{\beta/2,\varphi}(D)$ in $\beta$. For a summand
assigned to $c^\mathrm{t}$ we thus have
$e^{\lambda d_\beta s}\prod_Dz_{\beta,\varphi}(D)\leq
e^{\lambda q(\cD)}\prod_Dz_{\beta,\varphi}(D)\leq\prod_Dz_{\beta/2,\varphi}(D)$.
All the reference footprints in a summand lie in $P$, so every
reference dot and arrow incident to the neighborhoods $D^\dagger$ belongs to
the plaquette set $W=P$ of \cref{prop:positive-budget}. Discarding
the overlap-connectivity and union conditions, which only enlarges
the sum, and applying that proposition at inverse temperature
$\beta/2$ gives
\[
 c^\mathrm{t}_\beta(P)
 \leq e^{-\lambda d_\beta s}
       \sum_{\substack{\cD\text{ compatible}\\P_\varphi(D)\subset P}}
                      \prod_Dz_{\beta/2,\varphi}(D)
 \leq \exp([-\beta d_\beta/(4B)+Ce^{-\beta/128}]s)\,.
\]
Since $\beta d_\beta/B\asymp\beta^{1-1/10-1/2}=\beta^{2/5}$, this
proves \cref{eq:mixed-thin-bound} for large $\beta$. At inverse
temperature $\beta$, the same proposition, applied to all the summands
without the extracted factor, proves \cref{eq:mixed-fat-bound}.
The rooted count of the patterns with $c^\mathrm{t}_\beta(P)\neq0$
is crude: connected patterns of $s$ lozenges have at most $C_0^s$
rooted choices, since one may enumerate their $2s$ triangle vertices
as an animal in a fixed bounded-degree adjacency graph of the
triangular lattice, and then choose one of at most three matching
edges at each black vertex.

\subsubsection*{The \texttt{fat} bound} The size and entropy
assertions concern the patterns having at least one summand
assigned to $c^\mathrm{f}$. Consider such a summand. There must be a
\texttt{bad} support in $\cD$: otherwise $P$ is the union of the
\texttt{good} footprints $\Upsilon(P_\varphi(D))$, whence
$s\leq\sum_D|P_\varphi(D)|=q(\cD)<d_\beta s<s$, a contradiction. A
\texttt{bad} support $D$ has $s\geq |P_\varphi(D)|\geq B|D|\geq B$, proving
the size assertion. The union of the \texttt{good} footprints
contains fewer than $q(\cD)<d_\beta s$ tiles. For the \texttt{bad}
supports, the bounded footprint overlap of
\cref{lem:footprint-overlap} gives
$\sum_{D\text{ \texttt{bad}}}|P_\varphi(D)|\leq25s$. Repeat the six-class
construction in the proof of \cref{prop:stack-witness}
simultaneously for the long stacks of all the \texttt{bad}
footprints: the enlarged patterns $Q_\varphi(D)=W_\varphi(D)$ of these supports
are the unions of the tubes around their stacks of length at least
$L$ with their removed pieces, so their union is the union of the
six class unions and of the removed pieces. The geometric cores of
the class unions are unions of radius-$R$ balls, and their total
exceptional-orientation count is at most $CR^2$ times the number of
long stacks, which is at most $\sum_{D\text{ \texttt{bad}}}|P_\varphi(D)|/L
\leq25s/L$. The removed pieces, of size at most
$CL\sum_{D\text{ \texttt{bad}}}|D|\leq(CL/B)s$, together with the
\texttt{good} footprints, form a remainder of at most
$\delta_\beta s$ tiles, where
\[
        \delta_\beta=d_\beta+CL/B\,.
\]
Consequently \cref{lem:six-unions} applies to $P$ with
$\delta=\delta_\beta$, $\kappa=CR^2/L$, and
$l=\floor{\sqrt R}$. Its entropy rate tends to zero because
\begin{equation}\label{eq:entropy-rate-vanishing}
 l^{-2}+l/R+\delta_\beta l^2\longrightarrow0,
                    \qquad R^2/L\longrightarrow0\,,
\end{equation}
as computed in the proof of \cref{prop:stack-witness}, with the
additional term $d_\beta l^2\leq\beta^{-1/10+1/100}\to0$.
The same geometric-core argument shows that these patterns are
$(C/R+\delta_\beta,R/2)$-\texttt{fat}, with $C/R+\delta_\beta\to0$. The count in
\cref{eq:six-union-count} is modulo translations; rooting at a
specified triangle multiplies it by at most $2s$, and this factor
and the prefactor $C_l\leq Cl^2$ are absorbed in a vanishing entropy
rate, since $s\geq B\to\infty$ and $\log l/B\to0$ for the chosen
scales. Enlarge a single $\epsilon_\beta$ to
dominate all these rates and errors. This enumerates every pattern
having at least one summand in $c^\mathrm{f}$, while
\cref{eq:mixed-fat-bound} already bounds the sum of all its
summands. For a fixed reference torus tiling, the same proof applies using
\cref{prop:periodic-fiber,prop:periodic-witness}. The a priori bound
and the overlap bound twenty-five are unchanged, so both numerical
coefficient bounds remain valid. In the last counting step,
\cref{lem:periodic-unions} replaces \cref{lem:six-unions} and adds
only $C/\lambda$ to the entropy rate, including for noncontractible patterns.
\end{proof}

\subsection{Two steps: SOS histories}\label{sec:sos}

Let $\sA(\cI)$ denote the plaquette area, as in \cref{sec:tetris},
and for a \Tetris move put
\[
 E_i(\cI;\Ti\cI)=\sA(\cI)-\sA(\Ti\cI)\,,
\]
the excess energy of the move, which is nonnegative by
\cref{prop:sorting}.

An \SOS input $\cS$, monotone in $x_3$, has the two-step history
\[
 \cS\xrightarrow{\Ttwo}\cH\xrightarrow{\Tone}\varphi\,,
\]
where $\cH=\Ttwo\cS$ is \HMSOS by \cref{prop:sorting} and
$\varphi=\Tone\cH$ is monotone in all three coordinates. The
one-step analysis of \cref{sec:hmsos} describes the fiber of
$\Tone$ over $\varphi$ as a gas of defects on $\varphi$, and we
show in \cref{sec:sos-histories} that it applies equally to the
fiber of $\Ttwo$ over each resulting $\cH$, giving a gas of defects
on $\cH$. The difficulty in combining the two descriptions is that
the inner gas lives on a reference surface which itself depends on
the outer defects, so that the two families cannot be superposed
on $\varphi$. We therefore begin by defining two-step histories
and the footprints by which they are recorded on the final tiling
$\varphi$, so that the dependence of the inner defects on the outer
ones is kept; \cref{sec:sos-histories} then shows that the
resulting components behave like the polymers of a single
hard-core gas on $\varphi$, and the weights enter only from
\cref{sec:sos-exact-fiber} on.

\begin{definition}[Two-step history]\label{def:two-step-history}
Fix a monotone surface $\varphi$. A \emph{two-step history} over
$\varphi$ consists of an outer family $\cA$ of
\emph{$x_1$-defects} on $\varphi$, namely $D$-defects of the
reference $\varphi$ for the sorting direction $x_1$ in the sense
of \cref{lem:polymer-factorization}, with pairwise
non-$*$-adjacent supports $D_\alpha$, $\alpha\in\cA$, which
produce the \HMSOS surface $\cH$ with $\Tone\cH=\varphi$; and of an
inner family $\cB$ of \emph{$x_2$-defects} on $\cH$, defined in the
same way with $\cH$ as the reference and $x_2$ as the sorting
direction, with pairwise non-$*$-adjacent supports $D_\gamma$,
$\gamma\in\cB$, which produce the \SOS surface $\cS$ with
$\Ttwo\cS=\cH$. Its excess is $E_1(\cH;\varphi)+E_2(\cS;\cH)$.
We retain the defects themselves, not just their supports, since
the inner reference depends on the outer defects: two outer
defects with the same support produce different surfaces $\cH$,
and hence different inner gases.
\end{definition}

We use the lattice
$\Gamma=\langle(-n,0,\Theta_1),(0,-n,\Theta_2)\rangle$ of
\cref{sec:periodic-hmsos}, with both positive slopes in
the specified compact set and $n$ large enough that the fixed
neighborhoods embed in the cover. The inner \Tetris move sorts in
direction $x_2$, so its transverse coordinates are $(x_1,x_3)$,
and projecting the two generators of $\Gamma$ along $x_2$ shows
that their period lattice is
\[
 \Gamma_2=\langle(-n,\Theta_1),(0,\Theta_2)\rangle,\qquad
                 |\Z^2/\Gamma_2|=n\Theta_2\,.
\]
If $\cH$ is \HMSOS, write $j(x_1,x_3)$ for its column
magnetization in direction $x_2$, as in \cref{def:tetris-charge}.
Since $\cH$ is monotone in $x_2$, each $x_2$-fiber consists of a
single step, and $j(x_1,x_3)$ is simply the position of that
step: the surface $\cH$ is the graph of $j$ over the transverse
plane, and $j$ is an ordinary graph height in that direction.
Monotonicity of $\cH$ in $x_3$ makes $j$ nonincreasing in $x_3$,
but $\cH$ is not monotone in $x_1$, so $j$ need not be monotone
in $x_1$. We use the dot-and-arrow description \cref{eq:dot-arrow} of
\cref{sec:enumeration} relative to this reference: a dot is the
normal-$x_2$ plaquette at the step of a fiber, and the arrows on a
transverse edge are the wall plaquettes between the two
neighboring fibers. These dots and arrows count the actual
plaquettes of $\cH$, with multiplicity if their diagonal
projections overlap: because $j$ is not monotone in $x_1$, the
projection $\Upsilon$ may fold $\cH$, and distinct plaquettes of
$\cH$ may then project onto the same lozenge.

For a $*$-connected set $D$ of squares of this transverse torus
let $P_{\cH}^{(2)}(D)$ consist of the reference dots on its
two-square neighborhood and the complete reference walls between
squares of that neighborhood, exactly as in \cref{def:footprint}
with the sorting direction $x_2$ in place of $x_1$. In this
notation $P_{\cH}^{(2)}(D)$ is a set of plaquettes on $\cH$; it is
not being identified with a lozenge pattern, precisely because of
the possible folding just described. Likewise, $P_\varphi^{(1)}(D)$
denotes the reference footprint of \cref{def:footprint} on
$\varphi$ for the sorting direction $x_1$.

We now record a two-step history on the final tiling. The
footprints of the outer defects are reference footprints of
$\varphi$. The footprint of an inner defect lives on $\cH$, and it
is transported to $\varphi$ as follows: a face $q$ of $\cH$ either
is a face of $\varphi$, in which case it is its own witness on the
final tiling, or was created by a unique outer defect $\alpha$, in
which case we replace it by the whole footprint $P_\alpha$ of that
defect.

\begin{definition}[Footprints and components of a two-step history]
\label{def:history-footprint}
Let $(\cA,\cB)$ be a two-step history over $\varphi$ with
intermediate surface $\cH$. For $\alpha\in\cA$ with support
$D_\alpha$ write
\[
 P_\alpha=P_\varphi^{(1)}(D_\alpha),\qquad
 A_\alpha^-=\{\text{reference faces removed by }\alpha\},\qquad
 A_\alpha^+=\{\text{faces added by }\alpha\}\,,
\]
where faces are oriented plaquettes on the common lift, or
their quotient orbits: their occupied and vacant sides are
part of the label, so that a reversed face at the same geometric
location is removed and added, rather than declared
unchanged; cardinalities still count each surface plaquette
once. The \emph{carrier} of a face $q\in \cH$ is
\[
 \cC_{\cH}(q)=
 \begin{cases}
       \{q\},&q\in \cH\cap\varphi,\\
       P_\alpha,&q\in A_\alpha^+,
 \end{cases}
\]
which is well defined by the remark following this definition.
For $\gamma\in\cB$ with support $D_\gamma$, put
\[
 Q_\gamma=P_{\cH}^{(2)}(D_\gamma),\qquad
             R_\gamma=\bigcup_{q\in Q_\gamma}\cC_{\cH}(q)\,;
\]
thus $Q_\gamma$ is the footprint of the inner defect on the
intermediate surface, and $R_\gamma$ is its \emph{carried
footprint} on the final tiling. The \emph{assigned footprint} of
an outer defect is $P_\alpha$ and that of an inner defect is
$R_\gamma$. Join two defects when their assigned footprints share
a reference plaquette. A \emph{history component} is a connected
component of this overlap graph, and its \emph{footprint} is the
union $P$ of the assigned footprints of its defects; we also call
a history component with footprint $P$ a \emph{compound defect
with support $P$}.
\end{definition}

The sets $A_\alpha^-$ are contained in $P_\alpha$, since
a removed reference face is a dot or wall plaquette of
$\varphi$ at a changed column or at an edge incident to one.
Every face in $\cH\setminus\varphi$ belongs to exactly one
$A_\alpha^+$: a face whose status changes meets a changed
$x_1$-column or an edge incident to one, and such an incidence
cannot meet two nonadjacent defect supports. Hence the carrier
is well defined. All footprints of history components are
subsets of $\varphi$ and hence give lozenge patterns under
diagonal projection; unlike the inner footprints $Q_\gamma$ on
$\cH$, they are not affected by folding.

The carrier records a dependence on the intermediate surface.
For example, if an inner footprint $Q_\gamma$ uses a face created
by an outer defect $\alpha$, deleting $\alpha$ alone would change
the reference on which $\gamma$ is defined, and the inner defect
$\gamma$ would no longer make sense. Replacing that face
by $P_\alpha$ ensures that $R_\gamma$ contains $P_\alpha$, so
the two defects belong to the same component and are retained or
deleted together. The footprints also record the unchanged
faces needed to test neighboring-column constraints: an inner
defect whose footprint $Q_\gamma$ consists of faces of $\varphi$ only is
defined relative to a part of $\cH$ which no outer defect can
alter without overlapping it. This is why
the following independence statement concerns carried components,
rather than the two families of transverse supports separately.

\subsubsection{Independence of history components}
\label{sec:sos-histories}

\begin{lemma}[The inner fiber as a hard-core gas]\label{lem:sos-inner-gas}
For every \HMSOS surface $\cH$ with $\Tone\cH=\varphi$, the fiber
$\{\cS:\Ttwo\cS=\cH,\ \cS\text{ \SOS}\}$ is an exact hard-core gas
in the sense of \cref{lem:polymer-factorization,lem:quotient-one-step}
on the nonempty $*$-connected sets of squares of the transverse
torus $\Z^2/\Gamma_2$, including noncontractible ones and the whole
torus: a preimage is the same as a family of pairwise
non-$*$-adjacent supports $D$ together with an \SOS $D$-defect on
each, and the excess energies $E_2$ add.
\end{lemma}
\begin{proof}
In the sorting coordinates, the fiber consists of binary $x_2$-columns
with charges $j(x_1,x_3)$ and the given period offsets, constrained
by inclusion of the occupied column at $(x_1,x_3+1)$ in that at
$(x_1,x_3)$. This is precisely monotonicity in $x_3$. The relation
$\cI(x_1,x_2,x_3+k\Theta_2)=\cI(x_1,x_2+kn,x_3)$ turns the sorting ends
into vertical ends, so the configuration is an \SOS graph.

It is admissible also in the sorting direction. At fixed $x_1$,
tilted periodicity bounds its graph height below on each left
$x_2$ half-line and above on each right half-line. An occupied cube
can move downward below the first bound and then left to its end;
a vacant cube can move upward above the second bound and then right.
Thus both phases reach their sorting ends. Since $\Theta_2>0$ excludes a
period parallel to $x_2$, \cref{lem:quotient-one-step} gives the
gas structure.
\end{proof}

\begin{lemma}[Carried footprints and independent history components]
\label{lem:sos-carried-components}
Every assigned footprint and every history footprint is nonempty
and connected. Distinct history footprints have disjoint
triangle supports on the final tiling. Deleting any collection
of whole history components, in both stages, leaves an
admissible two-step history and preserves the retained defects
and their excess energies.
Conversely, independently admissible history components with
disjoint reference footprints can be combined, and their
total excess energies add. These assertions hold on the
periodic quotient, including for noncontractible components.
\end{lemma}
\begin{proof}
We prove the three assertions in turn: connectedness of the
footprints, deletion of components, and insertion of components.
\begin{enumerate}[label=\textbf{Step~\arabic*.}, ref=\arabic*, wide=0pt, itemsep=1ex, listparindent=\parindent, parsep=0pt]
\item\label[step]{st:sos-carried-connected}
(\emph{Connectedness.})
The footprint $P_\alpha$ is connected, being a reference footprint
(\cref{def:footprint}). So is $Q_\gamma$ on the
graph surface $\cH$, even though its diagonal projection may
fold: its transverse neighborhood is nearest-neighbor connected,
and the complete wall between consecutive graph columns
connects their reference dots.
To pass from $Q_\gamma$ to $R_\gamma$ we need to understand how
new faces of $\cH$ sit next to common faces.
If a new face in $A_\alpha^+$ is edge-adjacent to a common
face $q\in \cH\cap\varphi$, then $q\in P_\alpha$. Indeed their
projections in the $x_1$ sorting plane put the dot of $q$,
or both endpoints of its wall, within the two-square neighborhood
of $D_\alpha$, and the complete reference wall is included in
$P_\alpha$ by its definition.
Two new faces from distinct outer defects cannot be
edge-adjacent, since their changed transverse supports would
then be $*$-adjacent, contradicting compatibility. Thus along
any path of faces in $Q_\gamma$, the new faces come in runs,
each run created by a single outer defect $\alpha$ and flanked
by common faces lying in $P_\alpha$; replacing each run by the
connected footprint $P_\alpha$ therefore turns the path into a
connected union of footprints. This proves
connectedness and nonemptiness of $R_\gamma$. Unions along
the overlap graph preserve connectedness, which proves the
assertion for history footprints. The disjointness
assertion follows from the definition of the overlap graph:
two footprints sharing a reference plaquette would be joined.

\item\label[step]{st:sos-carried-deletion}
(\emph{Deletion.})
Retain any collection of history components and insert just
their outer defects, producing a surface $\cH'$. Outer
factorization (\cref{prop:periodic-fiber}) makes this a valid
\HMSOS preimage of $\varphi$. We claim that every face
of a retained inner footprint $Q_\gamma$ is still present in
$\cH'$. Indeed, a common face of $Q_\gamma$ was removed by none
of the original outer defects, hence by none of the retained
ones; and a new face of $Q_\gamma$ belongs to an outer defect
$\alpha$ whose footprint $P_\alpha$ occurs in $R_\gamma$, so that
$\alpha$ lies in the same component as $\gamma$ and is retained
as well. In particular, all the normal-$x_2$
dots on the footprint of $D_\gamma$ are unchanged.
Since $\cH'$ is monotone in $x_2$, these dots specify its
column magnetizations on the two-square neighborhood of $D_\gamma$ and therefore
specify exactly the same footprint $Q_\gamma$ as before.
The retained inner defects consequently have the same charges,
the same neighboring-column constraints and the same
excess-energy terms, all of which are read off the neighborhood.
Their binary columns and reference heights are unchanged, so the
fiber truncation of \cref{eq:fiber-truncation} is preserved as well. The
hard-core structure of \cref{lem:sos-inner-gas} on $\cH'$
allows deletion of the other inner defects. This proves the
deletion assertion with the actual intermediate reference
recomputed, rather than held fixed artificially.

\item\label[step]{st:sos-carried-insertion}
(\emph{Insertion.})
Consider now independently admissible history components with
disjoint reference footprints. The outer supports from distinct
components are compatible: incompatible supports would have
overlapping two-square neighborhoods and hence overlapping
$P_\alpha$, contradicting disjointness of the footprints.
Form their joint intermediate surface $\cH$ by inserting all the
outer defects, which is admissible by outer factorization.
A common face of a component's inner footprint cannot be
removed by an outer defect in another component, since that
removed face belongs to the other defect's $P_\alpha$,
which would again make the reference footprints overlap. A new
face belonging to its own outer defect is also unchanged
by the other, compatible outer defects. The same neighborhood-dot
argument as in \cref{st:sos-carried-deletion} therefore
preserves each inner footprint and the inner defects on it in the
combined $\cH$.

It remains to check that the inner supports from different
components are compatible on $\cH$. If two of them were
incompatible, their footprints would share a normal-$x_2$ dot
$q$ of this combined surface. When $q\in\varphi$, both
reference footprints contain $q$. Otherwise both contain
the carrier $P_\alpha$ of the unique outer defect that
created $q$. Both alternatives contradict disjointness of the
footprints. The inner supports are consequently compatible.
Apply their one-step factorization on $\cH$ to conclude
admissibility and additivity of $E_2$; outer factorization
gives additivity of $E_1$. All arguments use local incidences
and the fixed period offsets, so noncontractible supports create no exception.
\end{enumerate}
This completes the proof.
\end{proof}

\subsubsection{The exact two-step fiber}
\label{sec:sos-exact-fiber}

The law of $\varphi$ induced by the \SOS model is obtained by
summing the \SOS weight over the full fiber
$\{\cS:\Tone\Ttwo\cS=\varphi\}$, and this fiber is naturally
organized by the intermediate surface: every $\cS$ in it
determines $\cH$, and $\cH$ ranges over the one-step fiber of
$\varphi$ studied in \cref{sec:hmsos}. The one-step analysis alone
is however not sufficient, because the intermediate surface $\cH$
is not an \HMSOS sample: as we make explicit after
\cref{prop:sos-history-fiber}, it carries an additional weight,
namely the partition function of its own \SOS preimages. The full
two-step history must therefore be counted.

The analysis proceeds in five stages. We first write the two-step
fiber exactly as an iterated sum (\cref{prop:sos-history-fiber})
and control the inner sum, the \SOS fiber over an intermediate
surface, by an a priori bound (\cref{prop:sos-first-fiber}). We
then use the history components of \cref{sec:sos-histories} to
write the full weight as an exact hard-core gas on $\varphi$ with
nonnegative local coefficients (\cref{prop:sos-positive-expansion}).
The remaining stages split these coefficients according to the
excess energy per unit of footprint
(\cref{prop:sos-thin-coefficients}), construct \texttt{fat}
witnesses for the low-excess components
(\cref{lem:sos-stack-cover,prop:sos-witness}), and assemble the
\texttt{thin}/\texttt{fat} envelopes (\cref{prop:sos-mixed-envelopes})
that \cref{sec:model-grassmann} feeds into the \RG theorem.

We work first in finite volume, imposing the fiber truncation at distance $J$ of
\cref{eq:fiber-truncation} in each inverse step relative to its
immediate reference. For fixed $\varphi$ there are finitely many
choices of $\cH$ and then finitely many choices of $\cS$ for each
$\cH$, so all sums are finite. These restrictions increase to all
histories as $J\to\infty$ and are removed by monotone convergence.
All sums retain the chosen admissibility restrictions; truncation
subscripts are suppressed. Write
\[
 Z^{\textsc{sos}\to\textsc{hmsos}}_{\beta,2}(\cH)
   =\sum_{\substack{\cS\text{ an \SOS surface}\\\Ttwo\cS=\cH}}
                                  e^{-\beta E_2(\cS;\cH)}\,.
\]

\begin{proposition}[Two-step fiber identity]\label{prop:sos-history-fiber}
For each monotone $\varphi$, the full \SOS fiber partition function is
\begin{equation}\label{eq:sos-history-fiber}
 Z^{\textsc{sos}}_\beta(\varphi)
   =\sum_{\substack{\cH\text{ an \HMSOS surface}\\\Tone\cH=\varphi}}
        e^{-\beta E_1(\cH;\varphi)}
             Z^{\textsc{sos}\to\textsc{hmsos}}_{\beta,2}(\cH)\,.
\end{equation}
The identity extends to $J\to\infty$ as an equality of
nonnegative sums. It defines unrestricted probability laws
only when the resulting partition functions are finite.
\end{proposition}
\begin{proof}
Recall that $Z^{\textsc{sos}}_\beta(\varphi)$ is the sum of
$e^{-\beta(\sA(\cS)-\sA(\varphi))}$ over the \SOS surfaces $\cS$
with $\Tone\Ttwo\cS=\varphi$. Since $\Ttwo$ is a map, each such $\cS$
has a unique intermediate surface $\cH=\Ttwo\cS$, and
\cref{prop:sorting} makes it an \HMSOS surface with
$\Tone\cH=\varphi$. Adding and subtracting $\sA(\cH)$, the
energy telescopes exactly:
\[
 \sA(\cS)-\sA(\varphi)
   =\big(\sA(\cS)-\sA(\cH)\big)+\big(\sA(\cH)-\sA(\varphi)\big)
                   =E_2(\cS;\cH)+E_1(\cH;\varphi)\,.
\]
Partitioning the defining sum according to the value of
$\cH=\Ttwo\cS$, and noting that the first excess depends on $\cS$
only through the inner sum defining
$Z^{\textsc{sos}\to\textsc{hmsos}}_{\beta,2}(\cH)$, gives
\cref{eq:sos-history-fiber}. Both sides are sums of nonnegative
terms, so monotone convergence proves the assertion about
$J\to\infty$.
\end{proof}

Summing the \SOS weight $e^{-\beta\sA(\cS)}$ over the preimages of a
given \HMSOS surface shows that the marginal weight of $\cH$ induced
by the \SOS model is proportional to
\[
 e^{-\beta\sA(\cH)}
 Z^{\textsc{sos}\to\textsc{hmsos}}_{\beta,2}(\cH)\,.
\]
The first factor is the \HMSOS weight of $\cH$; the second factor,
which is at least one and depends on the geometry of $\cH$, is
absent from the \HMSOS model. Thus the second factor must be
included in the interaction analysis, and the unweighted \HMSOS
theorem alone does not give the \SOS theorem. The strategy below is
to control this factor by an a priori bound which is exponential
in the area of $\cH$ with a small rate (\cref{prop:sos-first-fiber}), and
then to absorb it into the coefficients of a single hard-core gas
on the final tiling.

\subsubsection{The first fiber over an intermediate surface}
\label{sec:sos-first-fiber}

Let $z_{t,\cH}^{(2)}(D)$ be the sum of $e^{-tE_2}$ over the \SOS
$D$-defects of \cref{lem:sos-inner-gas}, that is, over the
preimages of $\cH$ under $\Ttwo$ with defect support exactly $D$
which are monotone in $x_3$. We write $t$ rather than $\beta$ for
the inverse temperature because the bound below will be applied
at several values of this parameter.

\begin{proposition}[The first \SOS fiber and its a priori bound]
\label{prop:sos-first-fiber}
The bounds \cref{eq:correct-count,eq:activity-bound} hold for the
activities $z_{t,\cH}^{(2)}(D)$ of the gas of
\cref{lem:sos-inner-gas}, with the reference arrow count computed
on $\cH$. There are absolute constants $C,t_0$ such that, with
$\epsilon_t=Ce^{-t/64}$, for every finite set $W$ of reference
plaquettes of $\cH$,
\[
 \sum_{\substack{\cD\text{ compatible}\\
             P_{\cH}^{(2)}(D)\subset W\ (D\in\cD)}}
                  \prod_{D\in\cD}z_{t,\cH}^{(2)}(D)
                          \leq e^{\epsilon_t|W|},\qquad t\geq t_0\,.
\]
The empty family is included. In particular,
\begin{equation}\label{eq:sos-first-fiber-pressure}
 1\leq Z^{\textsc{sos}\to\textsc{hmsos}}_{t,2}(\cH)
                          \leq e^{\epsilon_t\sA(\cH)}\,.
\end{equation}
The normalized fiber truncated at distance $J$ converges in total
variation to the unrestricted normalized fiber as $J\to\infty$.
\end{proposition}
\begin{proof}
\Cref{lem:sos-inner-gas} identifies the fiber with a one-step gas
on the transverse quotient $\Z^2/\Gamma_2$, and $\Theta_2>0$
excludes a period parallel to $x_2$, so \cref{lem:quotient-one-step}
gives the enumeration, the activity bounds and the convergence
as $J\to\infty$. Every arrow counted in $\Arr_{\cH}(D^\dagger)$ and every required
reference dot belongs to the footprint $P_{\cH}^{(2)}(D)$. Its
inclusion in $W$ therefore permits applying that lemma's a priori
bound. Taking $W$ to be the entire quotient surface gives
\cref{eq:sos-first-fiber-pressure}. All these counts are on $\cH$
before projection, so folds of its diagonal projection cause no
change.
\end{proof}

\subsubsection{The positive expansion on the final tiling}
\label{sec:sos-positive-expansion}

We now allocate the excess of a two-step history to the triangles
of the final tiling and derive the exact expansion of the full
\SOS weight as a hard-core gas of history components. We first
locate all changes on a vertical column in one history component,
which will also be used in reconstructing the surface.

\begin{lemma}[A common carrier for changes on a vertical column]
\label{lem:sos-column-component}
In a two-step history $\cS\xrightarrow{\Ttwo}\cH\xrightarrow{\Tone}\varphi$,
let $\phi$ denote the vertical height of the monotone reference $\varphi$,
and write $q_z$ for its horizontal face over $z$.
If either stage changes a spin on this vertical column,
all changes on it belong to one base history component, whose
footprint contains $q_z$. The statement also holds on the periodic
quotient.
\end{lemma}
\begin{proof}
Consider first the outer displacement. When $\cH(z)\ne\phi(z)$,
the cubes over $z$ at which the spins of $\cH$ and $\varphi$
differ form one vertical interval, between the two heights. Each cube lies
on a changed $x_1$-sorting fiber, and the transverse
$x_1$-sorting squares of these fibers are consecutive on the line
of constant $x_2$, so they belong to a single outer defect
$\alpha$. The complete wall between the squares immediately above
and below $\phi(z)$ contains $q_z$. Indeed, put
$v_+=(x_2,\phi(z))$ and $v_-=(x_2,\phi(z)-1)$ for
$z=(x_1,x_2)$. The reference sorting charges satisfy
$h(v_+)\leq x_1<h(v_-)$: the reference column is vacant at
height $\phi(z)$ and occupied one level below. Thus the complete
wall between these adjacent sorting columns contains $q_z$, and
both columns lie in the neighborhood whenever either is changed.
Thus $q_z\in P_\alpha$. This argument concerns connectivity, not
distinctness of the squares, and so it remains valid on the
periodic quotient for intervals longer than a transverse period,
where the consecutive squares may repeat.

If $\cS(z)\ne \cH(z)$, the same argument in the $x_2$ sorting
coordinates puts all inner changes on this column in one
inner defect $\gamma$. Its footprint $Q_\gamma$
contains the horizontal face $q_z^{\cH}$ of $\cH$. When
$\cH(z)=\phi(z)$, this is the common face $q_z$ and its
carrier is $\{q_z\}$. Otherwise $q_z^{\cH}$ is a newly created
face of the outer defect $\alpha$ just identified: one of its
two adjacent cubes belongs to the differing vertical interval.
Its carrier is $P_\alpha$, so $\gamma$ and $\alpha$ lie
in the same base history component. This proves the first
assertion in all cases, including cancellation between the two
displacements. Distinct base footprints have disjoint triangle supports by
\cref{lem:sos-carried-components}, so this carrier is unique.

\end{proof}

The additivity in \cref{lem:sos-carried-components} gives an exact
allocation of the two-stage excess to the final reference triangles:
spread the excess of each history component uniformly over its carried
footprint (a footprint of $|P_C|$ lozenges covers $2|P_C|$
triangles). If $\scC$ denotes the set of history components,
define
\[
 e_x^{\textsc{sos}}(\cS,\cH;\varphi)
   =\sum_{C\in\scC}
       \frac{E_1(C)+E_2(C)}{2|P_C|}
                    \one_{\{x\in\Upsilon(P_C)\}}\,.
\]
At most one summand is nonzero, because distinct footprints have
disjoint triangle supports by \cref{lem:sos-carried-components},
and summing over $x$ recovers the total excess per period cell,
$\sum_x e_x^{\textsc{sos}}=E_\Gamma(\cS)$.
The normalization is by the carried diagonal footprint,
not by the number of transverse defect columns: a component with
few defect columns may have a large footprint, and its excess is
then spread thinly. A small uniform expectation of this
allocation is an additional estimate; it does not follow just
from disjointness, since disjointness bounds the number of
components meeting a triangle but not their weights.

\begin{proposition}[A positive local expansion for the full \SOS weight]
\label{prop:sos-positive-expansion}
For each connected reference lozenge pattern $P$, let
$c^{\textsc{sos}}_t(P)$ be the sum of $e^{-t(E_1+E_2)}$
over nonempty histories consisting of one history component
whose carried footprint is exactly $P$. For sufficiently
large $t$ these coefficients are finite, nonnegative and
determined by $P$. On the torus its inherited period
identifications are part of this data. They satisfy
\begin{equation}\label{eq:sos-positive-expansion}
 Z^{\textsc{sos}}_t(\varphi)
   =\sum_{\{P_i\}\text{ with disjoint triangle supports}}
                              \prod_i c^{\textsc{sos}}_t(P_i),
                         \qquad P_i\subset\varphi\,.
\end{equation}
There is $\widetilde\epsilon_t\leq Ce^{-t/64}$ such that,
for every set $W$ of reference plaquettes,
\begin{equation}\label{eq:sos-positive-collection-budget}
 \sum_{\substack{\{P_i\}\text{ disjoint}\\P_i\subset W}}
                      \prod_i c^{\textsc{sos}}_t(P_i)
                      \leq e^{\widetilde\epsilon_t|W|}\,.
\end{equation}
In particular
\begin{equation}\label{eq:sos-single-coefficient-budget}
          c^{\textsc{sos}}_t(P)
                           \leq e^{\widetilde\epsilon_t|P|}-1\,.
\end{equation}
These estimates sum all histories and decompositions
absolutely. They are not \texttt{thin}/\texttt{fat} witness bounds.
\end{proposition}
\begin{proof}
The proof has three steps. We first bound the total weight of
the histories whose carried footprints lie in a prescribed set
$W$, by applying the two one-step a priori bounds in succession; this
is where the inner weight of the intermediate surface is
absorbed. We then verify that the coefficient
$c^{\textsc{sos}}_t(P)$ depends only on the pattern $P$. Finally
we group the histories into components using
\cref{lem:sos-carried-components}.
\begin{enumerate}[label=\textbf{Step~\arabic*.}, ref=\arabic*, wide=0pt, itemsep=1ex, listparindent=\parindent, parsep=0pt]
\item\label[step]{st:sos-positive-budget}
(\emph{An a priori bound for histories with footprints in $W$.})
First restrict the outer defects by $P_\alpha\subset W$.
For their intermediate surface $\cH$ define the replacement patch
\[
            W^{\cH}=(W\cap \cH)\cup(\cH\setminus\varphi)\,,
\]
consisting of the faces of $W$ that survive on $\cH$ together
with all the new faces. All removed reference faces lie in
$W$, since $A_\alpha^-\subset P_\alpha\subset W$, so the exact
area identity is
\[
             |W^{\cH}|=|W|-|\varphi\setminus \cH|
                              +|\cH\setminus\varphi|
                    =|W|+E_1(\cH;\varphi)\,.
\]
This identity uses sets of plaquettes on the surfaces, not
areas of their possibly overlapping projections. The
condition that every inner carried footprint lie in $W$
is equivalent to $Q_\gamma\subset W^{\cH}$ for every inner
defect: common faces of $Q_\gamma$ are their own carriers and
must lie in $W$, while every new face has as carrier the footprint
$P_\alpha$ of an outer defect already allowed in $W$.

Apply \cref{prop:sos-first-fiber} on $\cH$, conditionally
on the specified outer defects, with the set $W^{\cH}$ in the
role of $W$ there. The total inner weight of the admissible
inner families is at most
\[
 \exp(\epsilon_t|W^{\cH}|)=\exp(\epsilon_t(|W|+E_1))\,.
\]
Let $\cF_t(W)$ denote the full nonnegative sum of
these restricted histories, including the empty history.
The factor $e^{\epsilon_tE_1}$ depends on the outer defects
only through their excess, so it can be absorbed into the inverse
temperature of the outer gas, and the one-step a priori bound
\cref{prop:positive-budget}, in its periodic form
\cref{prop:periodic-fiber}, applies at inverse temperature
$t-\epsilon_t$ to the outer defects with footprints in $W$:
\begin{align*}
 \cF_t(W)
  &\leq e^{\epsilon_t|W|}
       \sum_{\substack{\cA\text{ compatible}\\P_\alpha\subset W}}
                              e^{-(t-\epsilon_t)E_1}
       \leq\exp[(\epsilon_t+\epsilon_{t-\epsilon_t})|W|]\,.
\end{align*}
Increase the absolute constant in $\epsilon_t$ if necessary
so that it bounds both one-step a priori bounds. For large $t$,
$t-\epsilon_t\geq t_0$ and
$\epsilon_t+\epsilon_{t-\epsilon_t}\leq Ce^{-t/64}$.
This proves a uniform finite bound without treating $\cH$
as an independent unweighted \HMSOS sample; in particular the
coefficients $c^{\textsc{sos}}_t(P)$, which are sub-sums of
$\cF_t(P)$, are finite.

\item\label[step]{st:sos-positive-locality}
(\emph{Locality of the coefficients.})
We check locality before grouping the sum, including its
monodromy on the fixed quotient. The claim is that the set of
histories consisting of one component with footprint exactly
$P$, and their weights, can be read off the abstract pattern
$P$ with its period identifications, without reference to the
tiling in which $P$ sits. Align one face of two
tilings containing the connected pattern $P$.
The face labels and adjacency determine all relative lifts
along paths in $P$. A closed path in $P$ on the quotient has a
prescribed diagonal displacement, and its physical displacement
is the unique element of $\Gamma$ with that projection, since
$\Upsilon|_\Gamma$ is injective. Thus the lifted pattern,
including its noncontractible identifications, is the same in both
tilings up to a common diagonal translation. Projection to
either sorting plane retains these same prescribed period
offsets. No separate lift is chosen for different pieces of
a noncontractible component.

This lift of $P$ specifies the relative charges of every
candidate outer support whose footprint lies in $P$.
The possible defects on it and their weights are determined
by these charges and the local monotonicity constraints, as in
\cref{lem:mixed-locality}. Given such
a state family, the finite patch $P^{\cH}$ is explicitly known:
remove the specified reference faces and add the specified
new faces. Candidate inner supports have their normal-$x_2$
dots in this patch. Keep a candidate precisely when the
patch contains its neighborhood dots and every wall between their
specified charges. This is a finite test of
$P_{\cH}^{(2)}(D)\subset P^{\cH}$ and determines the full
footprint. Its inner state weights depend on the same charges
and the local $x_3$ inclusion constraints. The carrier labels of
new faces are known from the outer defects, so the overlap
graph and the requirement that its union equal $P$ are
determined as well. Consequently the sum defining
$c^{\textsc{sos}}_t(P)$ is independent of the tiling containing $P$.
A change of common lift translates the entire construction
and changes no weight. On a quotient the prescribed offsets
are retained in these tests, as in the one-step enumeration.

\item\label[step]{st:sos-positive-grouping}
(\emph{Grouping into components.})
\Cref{lem:sos-carried-components} now gives a unique
decomposition of every history into independently admissible
components with disjoint reference footprints, and its converse
permits every such collection: a collection of components with
disjoint footprints is a history, and its weight is the product
of the component weights by additivity of the excess. Summing
over the histories according to their collections of
footprints proves \cref{eq:sos-positive-expansion}, initially
as an identity between nonnegative sums and then, by
\cref{st:sos-positive-budget}, between finite sums. The same
argument restricted to footprints in $W$ gives
$\cF_t(W)=\sum_{\{P_i\}\subset W}\prod_i
c^{\textsc{sos}}_t(P_i)$, which together with the bound of
\cref{st:sos-positive-budget} proves
\cref{eq:sos-positive-collection-budget} with
$\widetilde\epsilon_t=\epsilon_t+\epsilon_{t-\epsilon_t}$.
Finally, taking $W=P$, the single nonempty component of
footprint $P$ is one term of $\cF_t(P)-1$, the subtracted one
being the empty history, which proves
\cref{eq:sos-single-coefficient-budget}.
\end{enumerate}
This completes the proof.
\end{proof}

\begin{proposition}[\texttt{Thin} coefficients and the precise low-excess class]
\label{prop:sos-thin-coefficients}
Split $c^{\textsc{sos}}_t(P)=c^{\mathrm{thin}}_t(P)
+c^{\mathrm{low}}_t(P)$ according to whether the total excess
of its history satisfies $E_1+E_2\geq |P|/\sqrt t$ or
$E_1+E_2<|P|/\sqrt t$, respectively. For sufficiently large
$t$,
\begin{equation}\label{eq:sos-thin-bound}
                 c^{\mathrm{thin}}_t(P)
                                  \leq e^{-\sqrt t\,|P|/4}\,.
\end{equation}
There are at most $C_0^s$ connected reference lozenge patterns
of size $s$ containing a prescribed triangle, in the plane
or on the fixed quotient. Moreover
\[
 c^{\mathrm{low}}_t(P)>0\quad\Longrightarrow\quad |P|>2\sqrt t,
 \qquad
 c^{\mathrm{low}}_t(P)\leq e^{\widetilde\epsilon_t|P|}-1\,.
\]
No \texttt{fat} property or improved entropy bound is asserted for this
low-excess class.
\end{proposition}
\begin{proof}
The idea for the \texttt{thin} bound is to pay for the footprint
with half of the Boltzmann factor and to control the remaining
sum by the a priori bound at half the inverse temperature.
For every history in the \texttt{thin} class, the defining
inequality $E_1+E_2\geq|P|/\sqrt t$ gives
\[
 e^{-t(E_1+E_2)}
   =e^{-(t/2)(E_1+E_2)}e^{-(t/2)(E_1+E_2)}
   \leq e^{-\sqrt t\,|P|/2}
                              e^{-(t/2)(E_1+E_2)}\,.
\]
Summing over these histories, and then enlarging the sum to all
histories of footprint $P$, gives
\[
 c^{\mathrm{thin}}_t(P)
   \leq e^{-\sqrt t\,|P|/2}c^{\textsc{sos}}_{t/2}(P)
   \leq\exp[(-\sqrt t/2+\widetilde\epsilon_{t/2})|P|]\,,
\]
by \cref{eq:sos-single-coefficient-budget} at inverse temperature
$t/2$, which proves \cref{eq:sos-thin-bound} for large $t$.

For the count, note that a connected pattern of $s$ lozenges
has a triangle support which is a connected animal of $2s$
elementary triangles containing the root, and that the pattern
is determined by this animal together with the choice, for each
triangle, of the neighbor it is paired with. The graph of
elementary triangles has bounded degree, including on the
quotient. A spanning-tree traversal bounds its rooted
animals of $2s$ triangles by $C^{2s}$. On a given such
animal the matching choices are at most $3^{2s}$. These
two bounds prove the claimed $C_0^s$ count, including
noncontractible patterns.

Finally, every nonempty history contains a nonempty one-step
defect, whose changed fiber has at least three transitions, so
its total excess is at least two by \cref{eq:exact-excess}. The
strict low-excess inequality $2\leq E_1+E_2<|P|/\sqrt t$
therefore implies $|P|>2\sqrt t$. The coefficient bound for
the low-excess part is inherited from
\cref{eq:sos-single-coefficient-budget}, since
$c^{\mathrm{low}}_t(P)\leq c^{\textsc{sos}}_t(P)$.
\end{proof}

\subsubsection{Straight stacks and \texorpdfstring{\texttt{fat}}{fat} witnesses for SOS histories}
\label{sec:sos-witnesses}

The low-excess class of \cref{prop:sos-thin-coefficients} is the
two-step analogue of the \texttt{bad} supports of
\cref{sec:fattening}: its components have footprints much larger
than their excess, and their coefficients may be as large as
$e^{\widetilde\epsilon_t|P|}$. As in the one-step case, the
strategy is to show that such a footprint is made of a bounded
number of long straight rows of lozenges, to thicken these rows
into tubes with one rare orientation
(\cref{lem:straight-tube}), and to deduce from
\cref{lem:six-unions} that there are few such patterns. The new
difficulty is that a carried footprint is built from pieces that
live on two different surfaces, $\varphi$ and $\cH$, and the
following lemma shows that the pieces coming from $\cH$ still
produce straight rows on $\varphi$.

\begin{lemma}[Interval intersections and a bounded-multiplicity cover]
\label{lem:sos-stack-cover}
Let $C$ be a carried history component, with footprint $P$
and total excess $E=E_1(C)+E_2(C)$. The pattern $P$ has a
cover by at most $CE$ pieces, each either a singleton tile
or a straight row of identical lozenges. Every tile belongs
to at most fifty of these pieces, counting repeated pieces.
The cover is determined by the history. The assertion holds
on the fixed quotient as well; its rows have injective lifts
modulo the period identifications.
\end{lemma}
\begin{proof}
By the carrier definition and the fact that every referenced
outer defect belongs to the same history component,
\[
 P=\bigcup_{\alpha\in C\cap\cA}P_\alpha
       \ \cup\!
       \bigcup_{\gamma\in C\cap\cB}(Q_\gamma\cap\varphi)\,.
\]
The first union runs over the outer defects of the component
and the second over its inner defects, whose carried footprints
$R_\gamma$ consist of the common faces of $Q_\gamma$ together
with the footprints $P_\alpha$ of the outer defects in the same
component. An outer footprint contains at most $25|D_\alpha|$
dots and $50|D_\alpha|$ complete wall stacks, each a straight
row in $\varphi$: the two-square neighborhood $D_\alpha^\ddagger$ has at
most $25|D_\alpha|$ squares, and at most twice as many edges
with both endpoints in it. An inner footprint has the
corresponding decomposition on $\cH$ in direction $x_2$. Its dots
contribute singletons when they belong to $\varphi$, and its wall
stacks are treated next.

Consider one inner wall stack, using a lift on the cover.
It is a row of plaquettes along $x_2$, of normal axis $x_1$ or
$x_3$, between two neighboring $x_2$-fibers of $\cH$. We need to
understand its intersection with $\varphi$.
A negatively oriented row has empty intersection with the
monotone reference, which has no negatively oriented faces.
For a positively oriented normal-$x_1$
row, recall from \cref{sec:periodic-hmsos} that
$h^{(1)}_\varphi(x_2,x_3)$ is the height of
$\varphi$ in direction $x_1$, a graph over the $(x_2,x_3)$-plane.
The faces of $\varphi$ on the supporting line of the row, at fixed
$x_1,x_3$, are those with $h^{(1)}_\varphi(x_2,x_3)=x_1$.
The left side is nonincreasing in $x_2$, so this level set
is an interval in $x_2$. For a normal-$x_3$ row the
corresponding condition is $\phi(x_1,x_2)=x_3$, where $\phi$ is
the height of $\varphi$ in direction $x_3$, also a level set of
a nonincreasing function of $x_2$. Intersecting with the
finite inner stack therefore gives at most one interval
in either case. Its diagonal projection is again a
straight row of identical lozenges. In particular an
inner stack does not split into arbitrarily many common
reference pieces when the outer defects are deleted; this is
the point of the lemma.

The total piece count is at most
$C(\sum_\alpha|D_\alpha|+\sum_\gamma|D_\gamma|)\leq CE$,
using that the excess of each defect is at least twice the
number of its defect squares. For the multiplicity, within
either family, any given reference plaquette lies
in at most twenty-five footprints, by the argument of
\cref{lem:footprint-overlap}: fix one
endpoint of its transverse projection and charge each
incidence to a distinct defect square within distance two.
For inner footprints this argument takes place on $\cH$;
restricting to the common faces cannot increase the count.
The pieces within each individual footprint are disjoint, being
distinct plaquettes. Adding the two families gives multiplicity
at most fifty.

On the quotient, choose any lift of each transverse edge
and its wall stack. The common-face test is unchanged by
a physical period because both surfaces have period
lattice $\Gamma$. The interval argument applies on that
lift. Two faces of one row can represent the same quotient
face only if a nonzero physical period is parallel to its
row direction. Positive windings exclude such a period
in either $x_1$ or $x_2$. All counts and the multiplicity
bound therefore survive projection.
\end{proof}

\begin{proposition}[\texttt{Fat} witnesses for low-excess \SOS histories]
\label{prop:sos-witness}
Let $R\geq2$, $L\geq CR^2$ and $B\geq CL$. If a history
component satisfies $|P|\geq BE$, there is a deterministic
connected enlargement $W_C\supset P$ such that
\[
 |W_C|\leq CR|P|,\qquad
 \diam W_C\leq\diam P+CR\,.
\]
It consists of at most six class unions of radius-$R$ tubes
and a remainder of at most $(CL/B)|P|$ tiles. The sum of
the exceptional-orientation counts in the six classes is
at most $(CR^2/L)|P|$. Consequently it is
$(C/R+CL/B,R/2)$-\texttt{fat} and satisfies the entropy bound
\cref{eq:six-union-count}, with $\delta=CL/B$ and
$\kappa=CR^2/L$. On the torus, for $R\leq c\lambda$,
the entropy bound is \cref{eq:periodic-union-count},
where $\lambda$ is the shortest diagonal period.
If $CL/B\leq c/R^2$, then
$|W_C|\geq cR\diam P$ in the relevant metric.
\end{proposition}
\begin{proof}
Apply \cref{prop:stack-witness} to the specified cover of
\cref{lem:sos-stack-cover}, with $m=E$ and a fixed $C_0$ bounding
both its piece count divided by $E$ and its multiplicity fifty.
This gives the deterministic enlargement and all the planar estimates.
On the quotient the rows project injectively by that lemma, so
\cref{prop:periodic-witness} gives the same geometric bounds and
\cref{eq:periodic-union-count}.
\end{proof}

\subsubsection{The mixed positive SOS coefficients}
\label{sec:sos-mixed-coefficients}

We now combine the previous results into the
\texttt{thin}/\texttt{fat} coefficients required by
\cref{thm:rg-typical}, following the one-step construction of
\cref{lem:mixed-locality}. Take $B=\sqrt\beta$ and the geometric
parameters $R,L,d_\beta$ as in \cref{eq:history-scales}.
The threshold $B$ is the one of \cref{prop:sos-thin-coefficients}
at $t=\beta$, and the parameters $R,L,B$ satisfy the hypotheses
$R\geq2$, $L\geq CR^2$ and $B\geq CL$ of \cref{prop:sos-witness}
for large $\beta$. For each specified base history component
$C$, use the enlarged pattern
\[
 Q(C)=\begin{cases}
       P_C,&E(C)\geq |P_C|/B,\\
       W_C,&E(C)<|P_C|/B
       \end{cases}\,,
\]
that is, high-excess components keep their footprint while
low-excess components are replaced by their \texttt{fat}
witnesses. Join base components when their enlarged patterns
share a lozenge. We call each connected component of this enlarged
overlap graph an \emph{enlarged history component}; it retains all
its base components and their defects, and its pattern is the
union of their enlarged patterns. Define $\bar c^{\textsc{sos}}_\beta(P)$ by
summing the weights of histories for which this enlarged
overlap graph is connected and its union is exactly $P$.
The histories themselves, and not only their footprints, enter
this sum, because $Q(C)$ can depend on them: the witness $W_C$ is
built from the stack cover of the component, which is determined
by its defects.

\begin{proposition}[\SOS \texttt{thin}/\texttt{fat} coefficient envelopes]
\label{prop:sos-mixed-envelopes}
For all sufficiently large $\beta$, the enlarged coefficients are
local and nonnegative, and
the exact identity \cref{eq:sos-positive-expansion} holds
with $c^{\textsc{sos}}_\beta$ replaced by
$\bar c^{\textsc{sos}}_\beta$. They split into nonnegative
parts $\bar c^{\mathrm t}_\beta+\bar c^{\mathrm f}_\beta$
satisfying
\begin{align}
 \bar c^{\mathrm t}_\beta(P)&\leq e^{-c\beta^{2/5}|P|},\notag\\
 \bar c^{\mathrm f}_\beta(P)&\leq
                       e^{\widetilde\epsilon_\beta|P|}\,.
                                      \label{eq:sos-mixed-envelopes}
\end{align}
\texttt{Thin} patterns have at most $C_0^s$ rooted choices of size
$s$. For some $\eta_\beta\to0$, every pattern with a
nonzero \texttt{fat} coefficient has $|P|>2\sqrt\beta$, is
$(\eta_\beta,R/2)$-\texttt{fat}, and belongs to a family with at
most $e^{\eta_\beta s}$ rooted patterns of size $s$.
On the fixed torus the rooted entropy rate is
$\eta_\beta+C/\lambda$, provided $R\leq c\lambda$.
All internal history decompositions have already been
summed in these numerical coefficient bounds.
\end{proposition}
\begin{proof}
We prove, in order, the exact expansion, the locality of the
enlarged coefficients, the \texttt{thin} bound, and the
\texttt{fat} geometry and entropy.

\subsubsection*{Exact expansion}
Every base history has a unique enlarged overlap
decomposition, obtained by grouping its history components
along the enlarged overlap graph. Since $Q(C)\supset P_C$,
disjoint enlarged patterns have disjoint base footprints, so
\cref{lem:sos-carried-components} permits independent insertion
and deletion across these groups, and the excess is additive.
The exact expansion with the enlarged coefficients follows
exactly as in \cref{st:sos-positive-grouping} of the proof of
\cref{prop:sos-positive-expansion}.

\subsubsection*{Locality}
For locality, enumerate the base histories with carried
footprints contained in $P$, as in
\cref{st:sos-positive-locality} of the proof of
\cref{prop:sos-positive-expansion}. Their
defects determine their interval stack covers,
their energies and which rows are retained. The tiling
closure of a prescribed ball union lies in $P$ precisely
when its geometric core lies in the region covered by
$P$; in that case its tiles are those of $P$ meeting the
core in positive area. This determines every $Q(C)$,
and hence the enlarged overlap and union tests, using
only $P$. Period offsets are retained as before.
Since every history counted in $\bar c^{\textsc{sos}}_\beta(P)$
has all its base footprints inside $P$, the entire coefficient is
bounded by the a priori bound of
\cref{eq:sos-positive-collection-budget},
$\bar c^{\textsc{sos}}_\beta(P)\leq\cF_\beta(P)\leq
e^{\widetilde\epsilon_\beta|P|}$, which is the asserted bound
on the \texttt{fat} part.

\subsubsection*{The \texttt{thin} part}
In a summand write $s=|P|$ and let
$q=\sum_{C:\,E(C)\geq |P_C|/B}|P_C|$ be the total footprint of
its high-excess components. Put the summand in the \texttt{thin}
part when $q\geq d_\beta s$. Each high-excess component has
$E(C)\geq|P_C|/B$, so the total excess of the summand is then
at least $q/B\geq d_\beta s/B$.
Taking half its Boltzmann factor out of the sum gives
\[
 \bar c^{\mathrm t}_\beta(P)
    \leq e^{-\beta d_\beta s/(2B)}\cF_{\beta/2}(P)
    \leq\exp([-\tfrac12\beta^{2/5}
                         +\widetilde\epsilon_{\beta/2}]s)\,,
\]
since $\beta d_\beta/(2B)=\tfrac12\beta^{2/5}$ for the chosen
parameters. Here the geometric classification is held fixed at
$\beta$; the last sum simply drops its restrictions and is
bounded by \cref{eq:sos-positive-collection-budget} at inverse
temperature $\beta/2$. This proves the \texttt{thin} bound. The
ordinary connected-pattern count was proved in
\cref{prop:sos-thin-coefficients}.

\subsubsection*{The \texttt{fat} part}
In the complementary summands $q<d_\beta s$. There must
be a low-excess component, since otherwise every component is
of high excess, the base footprints cover $P$, and $q=s>d_\beta
s$. Such a component has $|P_C|>2B$ by
\cref{prop:sos-thin-coefficients}, so $s>2\sqrt\beta$. The base
footprints are disjoint and contained in $P$, so their total
size is at most $s$. Apply \cref{prop:sos-witness}
simultaneously to the low-excess components and unite each of
their six classes across components. The exceptional-orientation
count is at most $(CR^2/L)\sum_C|P_C|\leq(CR^2/L)s$, and the
uncovered remainder, consisting of the remainders of the
witnesses and of all \texttt{thin} base footprints, has at most
\[
                 (d_\beta+CL/B)s
\]
tiles. \Cref{lem:six-unions} applies with
$\delta=d_\beta+CL/B$, $\kappa=CR^2/L$ and
$l=\floor{\sqrt R}$. Each term in its entropy
rate tends to zero as $\beta\to\infty$, as in
\cref{eq:entropy-rate-vanishing}, since $l^{-2}\asymp R^{-1}$, $l/R\asymp R^{-1/2}$,
$\delta l^2\leq C(\beta^{-1/10}+\beta^{-2/5})\beta^{1/100}$ and
$R^2/L\asymp\beta^{-2/25}$. The geometric cores show that the
pattern is $(C/R+\delta,R/2)$-\texttt{fat} with $C/R+\delta\to0$,
which is the asserted $(\eta_\beta,R/2)$-\texttt{fat} property. Fixed prefactors and the choice of a
triangle root are absorbed in a vanishing entropy rate,
using $s>2\sqrt\beta$. Enlarge $\eta_\beta$ to include
all these errors. The same argument on the quotient uses
\cref{lem:periodic-unions} and adds $C/\lambda$ to
the entropy rate. The coefficient bound itself follows
from the full a priori bound, as shown above, and already
includes every history giving this pattern.
\end{proof}

This completes the two-step analysis. The coefficients
$\bar c^{\textsc{sos}}_\beta$, their \texttt{thin}/\texttt{fat}
envelopes and the witness family of \cref{prop:sos-witness} are
the input to the \SOS case of the \RG theorem; their
consistency between volumes and the resulting exact sector law
are established in \cref{thm:sos-rg-input}.

\subsection{Three steps: Ising Interface histories}\label{sec:interface}

\subsubsection{The three-step interface fiber}
An interface $\cI$ of the Ising Interface model of
\cref{def:interface-model} may have overhangs, and is subject to
the no-bubble condition: every occupied cube is $18$-connected to
the occupied lower end and every vacant cube to the vacant upper
end. The first \Tetris move sorts the cubes along the vertical
lines, and produces a surface without overhangs; the two further
moves are those already analyzed in \cref{sec:sos}. Thus a
bubble-free interface has the three-step sorting history
\[
 \cI\xrightarrow{\Tthree}\cS\xrightarrow{\Ttwo}\cH\xrightarrow{\Tone}\varphi\,,
\]
in which, by \cref{prop:sorting}, the surface $\cS$ is \SOS,
$\cH$ is \HMSOS, and $\varphi$ is a monotone surface. The lower
part $\cS\to\cH\to\varphi$ of this history is exactly the two-step
\SOS history of \cref{sec:sos}; what is new is the third stage,
namely the fiber of bubble-free interfaces above a given \SOS
surface. We use the plaquette area $\sA$ and the one-step excesses
$E_i(\cI;\Ti\cI)=\sA(\cI)-\sA(\Ti\cI)$ of \cref{sec:sos}, and keep
a fiber truncation at finite distance $J$ in each inverse step initially, as there. The third-stage
fiber partition function over an \SOS surface $\cS$ is
\[
 Z^{\textsc{interface}\to\textsc{sos}}_{\beta,3}(\cS)
   =\sum_{\substack{\cI\text{ an admissible bubble-free interface}\\\Tthree\cI=\cS}}
                                      e^{-\beta E_3(\cI;\cS)}\,.
\]

Each interface $\cI$ determines its intermediate surfaces
$\cS=\Tthree\cI$ and $\cH=\Ttwo\cS$ uniquely, and the energy telescopes
exactly along the history:
\[
 \sA(\cI)-\sA(\varphi)
   =E_3(\cI;\cS)+E_2(\cS;\cH)+E_1(\cH;\varphi)\,.
\]
Grouping interfaces by their intermediate surfaces, as in
\cref{prop:sos-history-fiber}, therefore gives
\[
 Z^{\textsc{interface}}_\beta(\varphi)
 =\sum_{\Tone\cH=\varphi}e^{-\beta E_1(\cH;\varphi)}
       \sum_{\Ttwo\cS=\cH}e^{-\beta E_2(\cS;\cH)}
                Z^{\textsc{interface}\to\textsc{sos}}_{\beta,3}(\cS)\,.
\]
The sums retain all admissibility restrictions and the fiber truncation. In
particular they retain the no-bubble condition, which by
\cref{lem:phase-tree} is the same as edge-connectivity of the
disagreement surface; this restriction is carried by the third
stage alone, since the lower two surfaces are graphs, hence
connected. As in \cref{prop:sos-history-fiber},
monotone convergence extends the identity to $J\to\infty$ as
an equality of nonnegative sums, whereas normalizing it into a
probability law requires the finiteness proved below. Note
that the marginal weight of $\cS$ induced by the interface model is
proportional to
$e^{-\beta\sA(\cS)}Z^{\textsc{interface}\to\textsc{sos}}_{\beta,3}(\cS)$,
so the third-stage factor must enter the interaction analysis: the
\SOS theorem alone does not give the interface theorem.

\subsubsection{The first interface fiber}
\label{sec:interface-first-fiber}

Fix an \SOS reference $\cS$. The third \Tetris move sorts along
the vertical lines, so its transverse coordinates are the physical
columns, and its transverse torus is $\Z^2/n\Z^2$. Let $D$ denote
a nonempty $*$-connected set of physical columns in $\Z^2/n\Z^2$,
the support of a third-stage defect: the set of columns on which
a preimage $\cI$ differs from $\cS$. In analogy with the first
\SOS fiber of \cref{sec:sos-first-fiber}, write
$Q_{\cS}(D)=P_{\cS}^{(3)}(D)$ for its two-square reference
footprint on $\cS$, consisting of the reference dots of $\cS$ on
the two-square neighborhood of $D$ and of the complete reference walls
between squares of that neighborhood, with oriented faces as in
\cref{def:history-footprint}. The activity
$z^{(3)}_{t,\cS}(D)$ sums, with weight $e^{-tE_3}$, the bubble-free
configurations with $\Tthree\cI=\cS$, defect support exactly $D$, and
all other columns equal to $\cS$. Unlike the two lower fibers, a
third-stage defect is not given by a graph: a column of $\cI$ is an
arbitrary binary sequence with the prescribed ends and the
prescribed column magnetization, and the no-bubble
condition is a genuine constraint linking neighboring columns.

\begin{proposition}[The complete interface-to-SOS fiber]
\label{prop:interface-first-fiber}
The fiber over $\cS$ is an exact hard-core gas of these states
on pairwise non-$*$-adjacent supports, including noncontractible
supports and the whole transverse torus. Its excesses add.
The enumeration and activity bounds
\cref{eq:correct-count,eq:activity-bound} hold
with the actual reference arrow counts of $\cS$.
For $\epsilon_t=Ce^{-t/64}$ and $t$ sufficiently large,
\begin{equation}\label{eq:interface-first-budget}
 \sum_{\substack{\{D_i\}\text{ compatible}\\Q_{\cS}(D_i)\subset W}}
          \prod_i z^{(3)}_{t,\cS}(D_i)
                      \leq e^{\epsilon_t|W|}\,,
\end{equation}
for every finite set $W$ of reference plaquettes. In particular,
\begin{equation}\label{eq:interface-first-pressure}
 1\leq Z^{\textsc{interface}\to\textsc{sos}}_{t,3}(\cS)
                            \leq e^{\epsilon_t\sA(\cS)}\,.
\end{equation}
The normalized fiber converges in total variation as $J\to\infty$.
\end{proposition}
\begin{proof}
Apply \cref{lem:quotient-one-step} to the vertical fibers of $\cS$.
The input is constrained by its column charges and the no-bubble
condition, exactly as in \cref{lem:polymer-factorization}; no
additional coordinate monotonicity is imposed. The physical periods
have no nonzero vertical element. Hence the lemma gives the exact
gas, its enumeration and activity bounds, and its limit as $J\to\infty$,
including for a support equal to the whole transverse torus.

The footprint $Q_{\cS}(D)$ contains the reference dots and
arrows required by the a priori bound of that lemma. Its inclusion
in $W$ gives \cref{eq:interface-first-budget}. Set $W=\cS$ for the
upper bound in \cref{eq:interface-first-pressure}; the empty history
gives its lower bound.
\end{proof}

\subsubsection{Three-step components anchored in physical columns}
\label{sec:interface-components}

The \RG input requires the fiber partition function
$Z^{\textsc{interface}}_\beta(\varphi)$ to be written as a
hard-core gas of patterns on the final tiling $\varphi$, with every
intermediate state summed into the coefficients. For the two-step
\SOS history this was achieved in \cref{def:history-footprint,lem:sos-carried-components}
by carrying the footprints of both stages to $\varphi$. The
difficulty with the third stage is that its natural footprint
$Q_{\cS}(D_\delta)$ is read on the intermediate surface $\cS$,
which is itself random. We now explain how to anchor it on
$\varphi$ instead.

Fix $\varphi$, with physical height $\phi$. Decompose the
lower history $\cS\xrightarrow{\Ttwo}\cH\xrightarrow{\Tone}\varphi$
into the base \SOS components $C$ of
\cref{lem:sos-carried-components}, with footprints $P_C$.
For a specified third-stage defect $\delta$ of $\cI$ over $\cS$,
with support $D_\delta$, use the footprint
\[
                   F_\delta=P_\varphi^{(3)}(D_\delta)\,,
\]
that is, the two-square footprint of the same physical
columns, but with the reference dots and complete walls of
$\varphi$ in place of those of $\cS$. This footprint is read on
the final monotone surface, using physical columns, rather than on
the intermediate $\cS$. Join any two of these objects, lower
components or third-stage defects, when their assigned footprints
overlap; in particular two third-stage defects are joined when
$F_\delta$ and $F_{\delta'}$ share a reference face, and a
third-stage defect is joined to a lower component $C$ when
$F_\delta$ meets $P_C$. Call a connected group a base interface
history component and denote its footprint by $P$ (the union of
the assigned footprints of its members) and its total excess by
$E=E_1+E_2+E_3$.

The reason this anchoring loses no information is
\cref{lem:sos-column-component}: whenever either
$\cH(z)-\phi(z)$ or $\cS(z)-\cH(z)$ is nonzero, the whole lower
history on that physical column lies in a single base \SOS
component whose footprint contains the horizontal reference face
$q_z$ of $\varphi$ over $z$. Consequently, on every column of the
neighborhood of $D_\delta$ where $\cS$ differs from $\phi$, the lower
component responsible for the difference contains $q_z\in F_\delta$
and is thus joined to $\delta$; on all other neighborhood columns
$\cS=\phi$, and the footprint of $\varphi$ is the footprint of
$\cS$.

\begin{lemma}[Independent three-step history components]
\label{lem:interface-components}
The footprints just defined are connected; distinct components
have disjoint triangle supports. Whole components may be
inserted and deleted independently, with the intermediate
surfaces recomputed and retained states and excesses
preserved. Every component is determined by its footprint
and defects, from the reference pattern
there and the prescribed periods.
\end{lemma}
\begin{proof}
We prove connectedness, then deletion, then insertion, and finally
locality. Both kinds of assigned footprint are connected: the
lower footprints by \cref{lem:sos-carried-components}, and
$F_\delta$ because the two-square neighborhood of the $*$-connected set
$D_\delta$ is nearest-neighbor connected and consecutive neighborhood dots are joined
by their complete reference wall. A union of connected patterns
along a connected overlap graph is connected, and distinct overlap
components have disjoint triangle supports by definition.

For deletion, suppose $\delta$ is retained, and consider a column
$z$ in the two-square neighborhood of $D_\delta$. Its footprint $F_\delta$ contains
$q_z$. Any lower defect changing either of the two intermediate
heights at $z$ belongs to a lower component $C$ with $q_z\in P_C$,
by the column assertion recalled above, so that $P_C$ meets
$F_\delta$ and $C$ is in the same component as $\delta$.
Deleting any collection of other components therefore leaves $\cS$ unchanged on the whole two-square neighborhood of $D_\delta$. The data entering
the defect $\delta$, namely the reference charges of
$\cS$ on $D_\delta$, the complete walls on its neighborhood, and all
neighbor energy terms of \cref{eq:exact-excess}, are consequently
unchanged. Since its binary columns and reference heights are
unchanged, so is the fiber truncation of \cref{eq:fiber-truncation}.
Its no-bubble test is also unchanged: the proof of
\cref{lem:polymer-factorization}, used in
\cref{prop:interface-first-fiber}, follows a same-sign path
inside $D_\delta$ until its first exit, where it enters an
unchanged column of $\cS$ and continues vertically, and these
first-exit columns lie in the one-square neighborhood, where $\cS$ is
unchanged. Thus we may first delete the lower defects of the
discarded components, obtaining an admissible lower history by
\cref{lem:sos-carried-components}, and then delete their
third-stage states over the new $\cS$ by the hard-core gas
property of \cref{prop:interface-first-fiber}; the retained
states and their excesses are preserved at both steps.

For independent insertion, combine the lower histories of the
proposed groups by \cref{lem:sos-carried-components}.
A lower component from another group cannot change a column $z$ in the two-square neighborhood of $D_\delta$, since its footprint would then
contain the same $q_z$ and meet $F_\delta$, contradicting the
disjointness of the footprints. Hence every third-stage state
sees the same reference charges and walls after the combination
as in isolation. Two third-stage supports from different groups
cannot be $*$-adjacent: their reference footprints would
share a horizontal face, so their footprints would overlap.
\Cref{prop:interface-first-fiber} therefore permits inserting
all third-stage states simultaneously over the combined $\cS$,
and gives additivity of $E_3$; the lower excesses add by
\cref{lem:sos-carried-components}.

For locality, enumerate the lower histories with footprints in
$P$ and the candidate supports $D_\delta$ whose reference
footprints $F_\delta$ lie in $P$. The horizontal reference faces
of $P$ specify these physical columns and their neighborhoods. The
lower history supplies $\cS$ on each neighborhood column, and
the absence of a lower defect on a column means that $\cS=\phi$
there, by the column assertion. The third-stage states are then
binary column sequences with these fixed charges, subject to the
local no-bubble test just described and carrying their exact
excesses. All overlap, union and admissibility tests are thus
determined by $P$ together with the prescribed periods. This
argument retains the fixed identifications for noncontractible supports
and does not impose a finite vertical range on a third-stage
defect: a transition in a changed column may be
arbitrarily far from the reference height. Indeed, the monodromy
argument in \cref{prop:sos-positive-expansion} identifies the
lifted reference pattern in any two tilings with these
periods containing $P$. The lower histories then correspond, together with their
footprint dots and walls. Each third-stage column sequence
is anchored at its horizontal reference dot and is transported
with its integer transition offsets relative to that dot. The
prescribed periods give the same identifications of columns, the
same $*$-neighborhood tests and the same first-exit tests in both
tilings. This is a bijection of the defects, preserving
their energies, also when the footprint is noncontractible.
\end{proof}

In analogy with \cref{prop:sos-positive-expansion}, let
$\cF_t^{\mathrm I}(W)$ be the positive sum of $e^{-tE}$ over all
compatible collections of base interface history components whose
assigned footprints lie in a given finite set $W$ of reference plaquettes, including
the empty collection, and let $c_t^{\mathrm I}(P)$ be the sum of
$e^{-tE}$ over the histories consisting of exactly one nonempty
base component of footprint $P$. By \cref{lem:interface-components}
these coefficients are determined by $P$ and the prescribed
periods, and they sum all defects of all three stages.

\begin{proposition}[An a priori bound for the complete interface history]
\label{prop:interface-positive-budget}
For all sufficiently large $t$, with
$\epsilon_t^{\mathrm I}=\epsilon_t+
\widetilde\epsilon_{t-\epsilon_t}=O(e^{-t/64})$,
\begin{equation}\label{eq:interface-positive-budget}
       \cF_t^{\mathrm I}(W)
                       \leq e^{\epsilon_t^{\mathrm I}|W|}\,.
\end{equation}
The full fiber is exactly
\[
 Z_t^{\textsc{interface}}(\varphi)
   =\sum_{\substack{\{P_i\}\subset\varphi\\
                        S(P_i)\text{ disjoint}}}
                         \prod_i c_t^{\mathrm I}(P_i)\,.
\]
In particular $1\leq Z_t^{\textsc{interface}}(\varphi)
\leq e^{\epsilon_t^{\mathrm I}N_0}$, and the normalized fiber
converges in total variation as $J\to\infty$.
\end{proposition}
\begin{proof}
The strategy is the one of \cref{prop:sos-positive-expansion}:
we fix the lower history, bound the third-stage sum by the
first-fiber a priori bound of \cref{prop:interface-first-fiber} on a
suitable replacement patch of $\cS$, and then sum over the lower
histories using the \SOS a priori bound. The only new point is to compare
the footprints $F_\delta$, read on $\varphi$, with the footprints
$Q_{\cS}(D_\delta)$, read on $\cS$, to which that bound applies.

Fix a lower \SOS history with all its footprints in $W$, and put
$E_{12}=E_1+E_2$. Every oriented face of $\varphi$ removed in
forming $\cS$ belongs to $W$. Indeed, such a face was either
removed in the first inverse step by an outer defect, whose
footprint contains all the reference faces it removes, or it
survived into $\cH$ and was then removed by an inner defect, whose
carried footprint contains it by
\cref{def:history-footprint}. Define the replacement patch
\[
 W^{\cS}=(W\cap \cS)\cup(\cS\setminus\varphi)\,,
\]
a set of plaquettes of $\cS$. Since $\varphi\setminus\cS\subset W$,
we have $|W\cap\cS|=|W|-|\varphi\setminus\cS|$, and therefore
\[
 |W^{\cS}|=|W|-|\varphi\setminus\cS|+|\cS\setminus\varphi|
         =|W|+\sA(\cS)-\sA(\varphi)=|W|+E_{12}\,,
\]
where $\sA(\varphi)=N_0$ is the common area of all monotone
surfaces in the sector.

Now let $\delta$ be a third-stage defect with $F_\delta\subset W$.
Every common face of $Q_{\cS}(D_\delta)$ and $\varphi$ belongs to
$F_\delta$: it is either the dot at a neighborhood column $z$, which is
then the common horizontal face $q_z$, or part of a complete
reference wall of $\cS$ on a neighborhood edge, and a face of $\varphi$
on that wall belongs to the complete reference wall of $\varphi$ on
the same transverse edge. Every other face of $Q_{\cS}(D_\delta)$
belongs to $\cS\setminus\varphi$. Thus
$Q_{\cS}(D_\delta)\subset W^{\cS}$, and the third-stage a priori bound
\cref{eq:interface-first-budget}, applied with reference $\cS$ and
patch $W^{\cS}$, bounds the sum over the allowed third-stage states
by $e^{\epsilon_t|W^{\cS}|}=e^{\epsilon_t(|W|+E_{12})}$. Summing
over the lower histories, whose own footprints lie in $W$, gives
\[
 \cF_t^{\mathrm I}(W)
 \leq e^{\epsilon_t|W|}
      \sum_{\text{lower histories in }W} e^{-(t-\epsilon_t)E_{12}}
 = e^{\epsilon_t|W|}
          \cF^{\textsc{sos}}_{t-\epsilon_t}(W)
 \leq e^{(\epsilon_t+\widetilde\epsilon_{t-\epsilon_t})|W|}\,,
\]
by \cref{eq:sos-positive-collection-budget}. This is only an upper
bound: the restriction $F_\delta\subset W$ implies
$Q_{\cS}(D_\delta)\subset W^{\cS}$, and no assertion that the two
patch restrictions are equivalent is required. The exact expansion
of $Z_t^{\textsc{interface}}(\varphi)$ follows from the independent
decomposition of \cref{lem:interface-components}: every history
decomposes uniquely into its base components, and conversely any
collection of components with disjoint footprints can be
combined, with additive excess. Taking $W$ to be all of $\varphi$
proves the finiteness bound $Z_t^{\textsc{interface}}(\varphi)\leq
e^{\epsilon_t^{\mathrm I}N_0}$, the lower bound being the empty
history, and the limit $J\to\infty$ follows by monotone convergence as
in \cref{prop:interface-first-fiber}.
\end{proof}

\subsubsection{The three-stage stack cover}
\label{sec:interface-coefficients}

\begin{proposition}[Three-stage stack cover]
\label{prop:interface-mixed-envelopes}
A base interface component of footprint $P$ and excess $E$
has a cover by at most $CE$ singleton tiles and straight rows
of identical lozenges, with multiplicity at most seventy-five.
It has the low-excess witnesses of
\cref{prop:sos-witness}, with changed absolute constants,
including for noncontractible footprints.
\end{proposition}
\begin{proof}
By \cref{prop:stack-witness,prop:periodic-witness}, it suffices to
construct the stated cover and check that its rows lift injectively.
The footprint $P$ is the union of the lower footprints $P_C$
of its \SOS components and the third-stage footprints $F_\delta$.

The lower footprints are disjoint, and the union of their stack
covers from \cref{lem:sos-stack-cover} has at most $C(E_1+E_2)$
pieces and multiplicity at most fifty. Each $F_\delta$ is the
two-square footprint of $D_\delta$ on $\varphi$: the
neighborhood has at most $25|D_\delta|$ squares, contributing at most
$25|D_\delta|$ reference dots, which we take as singletons, and
at most $50|D_\delta|$ neighborhood edges, each carrying a complete
vertical wall of $\varphi$, which is a straight row of identical
lozenges. A square belongs to the two-square neighborhood of a support
only if that support contains one of the twenty-five squares
within two $*$-steps of it, and compatible supports are disjoint,
so a tile of $\varphi$ lies in the third-stage footprints of at
most twenty-five defects. Since $E_3(\delta)\geq2|D_\delta|$ for
every defect, by \cref{prop:correct-count}, the third-stage pieces
number at most $75\sum_\delta|D_\delta|\leq CE_3$ in total. The combined cover
thus has at most $CE$ pieces and multiplicity at most $50+25=75$,
as asserted. On the quotient, a vertical row projects injectively,
because $\Gamma$ has no nonzero vertical period, and the lower rows
have injective lifts by \cref{lem:sos-stack-cover}. Thus the cover
and its bounds hold for noncontractible patterns as well.

Apply \cref{prop:stack-witness} to this cover with $m=E$ and
$C_0$ large enough for its piece count and multiplicity seventy-five.
The cover is determined by the history, so its enlargement is
also deterministic. The injective row lifts just verified allow
\cref{prop:periodic-witness} to be applied on the quotient.
This gives all the asserted witness estimates.
\end{proof}

\subsubsection{Local contact perturbations}
\label{sec:interface-perturbations}

The full 3\Dim Ising model differs from the Ising Interface model
by the presence of closed contours (bubbles) away from the
interface. In \cref{sec:ising} these are integrated out by a
low-temperature contour expansion, which multiplies the weight of
a bubble-free interface $\cI$ by $e^{U(\cI)}$, where $U$ is a sum of
local terms: each term is attached to a cluster of contours, decays
exponentially in the total size of the cluster, and depends on the
interface only through whether the interface touches the cluster.
To avoid repeating the three-step analysis for the perturbed
model, we prove the interface estimates for the following class of
interactions, which includes the closed-contour correction of
\cref{sec:ising} and allows either sign for each local term.

\begin{definition}[Admissible contact perturbation]
\label{def:interface-perturbation}
For a torus matching $M$, let $\Gamma_M$ be the period lattice of
its lifted monotone surface, as constructed in
\cref{sec:interface-perturbed-positive-completion}.
On the full lattice and on every cylinder $\Z^3/\Gamma_M$
take a countable family of labels $X$.
A label has a nonempty edge-connected plaquette carrier
$N_X$, an integer mass $m_X\geq1$ with $|N_X|\leq C_0m_X$,
and a real number $b_X$. Put
\begin{equation}\label{eq:interface-contact-functional}
 a_t(X)=b_Xe^{-tm_X},\qquad
 h_X(\cI)=\one_{\{\cI\cap N_X\ne\varnothing\}},\qquad
 U_{t,\Gamma_M}(\cI)=\sum_X a_t(X)h_X(\cI)\,.
\end{equation}
This defines the energy on finite-area quotient interfaces. On the
plane we use its differences under finite modifications; these sums
converge by the anchored bound below, even for an infinite reference.
The data are independent of $t$ and the auxiliary dimer weights,
translation covariant, and satisfy, for all sufficiently large $t$,
\begin{equation}\label{eq:interface-contact-norm}
 \sup_{\Gamma_M}\sup_p
       \sum_{X:p\in N_X}|b_X|e^{-3tm_X/4}\leq Ce^{-ct}\,.
\end{equation}
The planar family satisfies the same bound. Constants are uniform
over auxiliary sectors, including sectors with a zero face density.

We also require the following precise compatibility of the labels
under lifting. For a connected partial matching $P$, choose a
connected component of its planar lift, with stabilizer $H_P$ in
the diagonal period lattice. Once the height of one face of this
component is fixed, so that the component lifts to a patch of a
monotone surface in $\Z^3$, labels whose carrier projection
$\Upsilon(N_X)$ lies in that component
are determined by the lifted pattern and the period map on $H_P$.
Changing the matching outside $P$ gives a bijection of
these anchored labels preserving $N_X,m_X,b_X$. For $H_P=\{0\}$
this is the bijection with finite planar labels. Relative integer
heights of labels are retained; no vertical-range cutoff is imposed.
\end{definition}

Let us comment on the ingredients of this definition. The
carrier $N_X$ is the set of plaquettes through which the label
``feels'' the interface, and the mass $m_X$ measures its size, so
that $a_t(X)$ decays exponentially in the size of the label. The
norm \cref{eq:interface-contact-norm} is an anchored bound: it
controls the labels touching a given plaquette $p$, which is the
form in which such sums are used below, and it is stated with the
exponent $3t/4$ rather than $t$ so as to leave a factor
$e^{-tm_X/4}$ per label unused; we call such an unused fraction of
a decaying factor, kept to pay for a later sum, a \emph{reserve}.
This reserve absorbs the combinatorial
factors, exponential in $m_X$, that arise when a label is
distributed over its contact faces or expanded by
inclusion--exclusion over the components it touches; see
\cref{prop:interface-perturbed-reference-potential,prop:interface-perturbed-absolute-history}.
The lifting condition is what makes the perturbed coefficients
functions of the pattern alone, as the \RG input of
\cref{def:rg-admissible} requires: the labels live on the cylinder
$\Z^3/\Gamma_M$ of the matching $M$, and the condition says
that those labels whose carrier projects inside a component of the planar lift
of a pattern $P$ can be transported between any two
matchings containing $P$. This last condition is part of the
hypothesis, not a consequence of an anchored norm. In the contour
application it follows by transporting the complete boundaries and
their zero-end fillings; this is checked in
\cref{prop:ising-admissible-perturbation} below. The unperturbed
model is the empty label family, $U=0$. The interface estimates
are uniform for fixed constants in \cref{def:interface-perturbation}.

For later use let $\epsilon_t=Ce^{-ct}$, increasing $C$ as needed.
Anchoring a term at a face it touches gives
\begin{equation}\label{eq:interface-contact-variation}
 |U_t(\cI)|\leq\epsilon_t\sA(\cI),\qquad
 |U_t(\cI)-U_t(\cI')|\leq\epsilon_t|\cI\xor \cI'|\,.
\end{equation}
Indeed, a label contributes to $U_t(\cI)$ only if its carrier
meets $\cI$, so
\[
 |U_t(\cI)|\leq\sum_{p\in\cI}\sum_{X:\,p\in N_X}|a_t(X)|
          \leq\sA(\cI)\sup_p\sum_{X:\,p\in N_X}|b_X|e^{-tm_X}
          \leq\epsilon_t\sA(\cI)\,,
\]
by \cref{eq:interface-contact-norm}. Similarly,
$h_X(\cI)-h_X(\cI')$ vanishes unless $N_X$ meets exactly one of
$\cI,\cI'$, in which case $N_X$ meets $\cI\xor\cI'$, and
anchoring at such a face gives the second bound. The second bound
is also the sense in which $U_t$ is defined on the plane: for an
infinite reference interface the sum defining $U_t$ need not
converge, but its differences under finite modifications do, by
this anchored estimate. Both assertions hold after any restriction
of the label family. The same proof bounds the absolute sums, so
the rearrangements below are justified whenever their stated
envelopes are finite.

\begin{proposition}[Weighted a priori bound and reference factor]
\label{prop:interface-perturbed-history-budget}
\label{prop:interface-perturbed-reference-potential}
For histories with old footprints in a finite set $W$ of reference plaquettes,
\begin{equation}\label{eq:interface-perturbed-history-budget}
 \sum_{\cX\text{ in }W}
 e^{-tE(\cX)+U_t(\cI_{\cX})-U_t(\varphi)}
 \leq \exp[(2\epsilon_t+
                     \epsilon^{\mathrm I}_{t-\epsilon_t})|W|]\,.
\end{equation}
There are real local reference coefficients with
\[
 U_{t,\Gamma}(\varphi)=\sum_{P\subset\varphi}v_{t,\Gamma}(P),
 \qquad
 \sup_o\sum_{P\ni o}|v_{t,\Gamma}(P)|e^{c_0t|P|}\leq Ce^{-ct}\,,
\]
where the norm sums over all connected patterns through $o$.
Coefficients agree under the anchored lifting in
\cref{def:interface-perturbation}.
\end{proposition}
\begin{proof}
We first show that the symmetric difference of $\cI_{\cX}$ and
$\varphi$ is controlled by the patch and the excess. Every oriented
reference face removed in forming $\cS$ belongs to $W$, as shown
in the proof of \cref{prop:interface-positive-budget}. A common
face of $\cS$ and $\varphi$ removed in the third step is removed by
a defect $\delta$ whose support contains one of its two adjacent
columns, so it belongs to the footprint of $D_\delta$ and
hence to $F_\delta\subset W$. Thus $\varphi\setminus\cI\subset W$,
whence $|\varphi\setminus\cI|\leq|W|$ and
\[
 |\cI\setminus\varphi|
   =|\varphi\setminus\cI|+\sA(\cI)-\sA(\varphi)\leq|W|+E\,,
\]
so that $|\cI\xor\varphi|\leq2|W|+E$, also for the unoriented
face sets. By \cref{eq:interface-contact-variation},
\begin{equation}\label{eq:contact-variation-history}
 U_t(\cI_{\cX})-U_t(\varphi)
   \leq\epsilon_t|\cI_{\cX}\xor\varphi|
   \leq2\epsilon_t|W|+\epsilon_tE(\cX)\,,
\end{equation}
and therefore the left-hand side of
\cref{eq:interface-perturbed-history-budget} is at most
$e^{2\epsilon_t|W|}\sum_{\cX}e^{-(t-\epsilon_t)E(\cX)}
=e^{2\epsilon_t|W|}\cF^{\mathrm I}_{t-\epsilon_t}(W)$, and the
a priori bound \cref{eq:interface-positive-budget} at
$t-\epsilon_t$ proves the first bound.

For the reference factor, let $A_X^\varphi$ be the tiling closure
of the diagonal projection of $N_X$, that is, the set of lozenges
of $\varphi$ whose interiors meet $\Upsilon(N_X)$. Since
$N_X$ is edge-connected, so is its projection, and $A_X^\varphi$
is a connected pattern with at most $Cm_X$ tiles, by
$|N_X|\leq C_0m_X$. For a label $X$ touching $\varphi$, distribute
its weight $a_t(X)$ equally among its contact faces
$p\in\varphi\cap N_X$, and assign each share
$a_t(X)/|\varphi\cap N_X|$ to the pattern $P=A_X^\varphi$. Summing
the shares assigned to $P$ defines $v_{t,\Gamma}(P)$, and the
identity $U_{t,\Gamma}(\varphi)=\sum_{P\subset\varphi}v_{t,\Gamma}(P)$
holds because the shares of each touching label sum back to
$a_t(X)$. To see that $v(P)$ is determined by the pattern, note
that a contact face $p$ fixes the height of $X$ relative to the
lifted pattern, and that no reference face outside $P$ can meet
$N_X$, since the projection of such a face would then meet
$\Upsilon(N_X)$; hence the lifted pattern determines, for each
label assigned to it, all its contact faces and their number.

For the norm, fix $o$ and a mass $m$. At mass $m$ there are at most
$C_1^m$ tiling closures and marked lifts, and when $P\ni o$ the
marked face $p$ has at most $Cm^2$ possible projected positions,
since $\diam P\leq Cm$. Anchoring at $p$, the labels of mass $m$
touching $p$ contribute at most
$\sup_p\sum_{X:p\in N_X}|b_X|e^{-3tm_X/4}\leq Ce^{-ct}$, by
\cref{eq:interface-contact-norm}, times the reserved factor
$e^{-tm/4}$. Hence
\[
 \sum_{P\ni o}|v_{t,\Gamma}(P)|e^{c_0t|P|}
   \leq Ce^{-ct}\sum_{m\geq1}Cm^2C_1^me^{c_0tCm}e^{-tm/4}
   \leq Ce^{-ct}\,,
\]
after first choosing $c_0$ small, so that $c_0C<1/8$, and then $t$
large. This proves the all-pattern norm. Finally, the stipulated
label bijection of \cref{def:interface-perturbation} preserves the
carriers and masses, hence the contact faces and the denominator
$|\varphi\cap N_X|$, and therefore also preserves these
coefficients.
\end{proof}

\subsubsection{Projection-controlled history footprints}
\label{sec:interface-perturbed-projected-history}

The next observation controls the distant transition heights
that were retained in the third-stage defects. Projections of
sets of cubes mean their unions of elementary diagonal
triangles; boundaries of triangles play no role in the counts.

\begin{lemma}[Excess bounds the projected changes outside a footprint]
\label{lem:interface-perturbed-projection-excess}
Let $L$ be a periodic interface with fixed spins outside a
finite tilted height band, and let $\varphi$ be a monotone reference
in the same homology class. Suppose a
reference pattern $P$ contains every oriented face of
$\varphi$ removed in $L$. Let $\Omega(L,\varphi)$ be the
set of cubes whose spins differ. Then
\begin{equation}\label{eq:interface-perturbed-projection-excess}
 |\Upsilon\Omega(L,\varphi)\setminus S(P)|
                          \leq \sA(L)-N_0\,.
\end{equation}
The same holds for finite modifications in the plane, with the
area difference interpreted relatively.
\end{lemma}
\begin{proof}
Orient every face of $L$ from its occupied to its vacant cube,
and call it positive when the resulting normal is one of
$e_1,e_2,e_3$, negative otherwise. All faces of the monotone
surface $\varphi$ are positive. Consider a generic line in
direction $(1,1,1)$, which projects to an interior point of an
elementary triangle $q$. Traversed in the direction $(1,1,1)$, it
starts in the occupied lower end and ends in the vacant upper
end, and it crosses a face of $L$ exactly when the spin
changes along it: from occupied to vacant at a positive face and
from vacant to occupied at a negative face. The crossings therefore
alternate in sign, beginning and ending with a positive one, so
the numbers of positive and negative crossings satisfy
$n_+(q)-n_-(q)=1$. Each plaquette projects onto two triangles,
and the quotient has $2N_0$ triangles, so summing over $q$ on the
quotient gives $\sum_q(n_++n_-)=2\sA(L)$ and
$\sum_q(n_+-n_-)=2N_0$, whence
\[
 2\#\{\text{negative faces of }L\}=\sA(L)-N_0,
 \qquad \sum_q n_-(q)=\sA(L)-N_0\,.
\]
Now let $q$ lie outside $S(P)$. The reference face of $\varphi$
over $q$ is not removed in $L$, by hypothesis, so it remains a
positive crossing of the line. If the line has no other crossing,
the line is occupied below this face and vacant
above it, exactly as for $\varphi$, so no cube on this line belongs
to $\Omega(L,\varphi)$. Hence, if spins differ somewhere on
the line, it has at least two crossings, and by the alternation at
least one negative crossing. Thus the indicator of
$\Upsilon\Omega(L,\varphi)$ at $q$ is at most $n_-(q)$, and summing over
$q\notin S(P)$ proves \cref{eq:interface-perturbed-projection-excess}.
The planar proof is the same, summing only over the triangles of
the finite changed region; all other lines have their one common
positive crossing and contribute zero to both sides.
\end{proof}

The footprints of \cref{sec:interface-components} record the
reference footprints of the changed columns, but not the positions of
the changed cubes themselves, which may lie arbitrarily far from
the reference surface. This was harmless for the unperturbed
model, whose energy is local in the dot-and-arrow diagram. The
contact perturbation, however, depends on the actual positions of
the changed cubes through the carriers $N_X$: a label touching a
distant transition of some component must be attributed to that
component. We therefore enlarge each footprint by the projections of
all the cubes the component changes, at any of the three stages, and
use \cref{lem:interface-perturbed-projection-excess} to show that
the enlargement costs only the excess.

For an isolated base interface component $C$, let
$(\cI_C,\cS_C,\cH_C,\varphi)$ be its reconstructed history and
let $P_C,E_C$ be its old footprint and excess. Define $Q_C$
to be the reference tiling closure of
\begin{equation}\label{eq:interface-perturbed-expanded-footprint}
 S(P_C)\ \cup\ \Upsilon\Omega(\cH_C,\varphi)
       \ \cup\ \Upsilon\Omega(\cS_C,\varphi)
       \ \cup\ \Upsilon\Omega(\cI_C,\varphi)\,,
\end{equation}
that is, the set of lozenges of $\varphi$ whose interiors meet this
union. These are finite sets on the quotient, since a component
changes finitely many cubes; no bound on transition heights is
imposed in their definition.

\begin{proposition}[Expanded footprints and independent components]
\label{prop:interface-perturbed-expanded-footprints}
The pattern $Q_C$ is connected and determined by the reference
pattern $P_C$ and the defects of $C$. Moreover
\[
                       |Q_C\setminus P_C|\leq3E_C\,.
\]
Regroup base components by overlap of their $Q_C$'s. The
resulting \emph{expanded components} are independently insertable
and deletable, have disjoint expanded footprints $Q$, and retain
additive excesses, and the sum of $e^{-tE}$ over the compatible
families of expanded components with all $Q\subset W$ is at most
$\exp(\epsilon_t^{\mathrm I}|W|)$.
All changed cubes in every intermediate configuration project
into the footprint of their expanded component. Thus distinct
expanded components change disjoint sets of cubes, with no
adjacent changed cubes belonging to different expanded
components, and each changed plaquette belongs to exactly one
expanded component.
\end{proposition}
\begin{proof}
We prove the size bound, then connectivity, then the independence
of the regrouped components, and finally the assertions on changed
cubes. For the size bound, all reference faces removed in
$\cH_C,\cS_C,\cI_C$ belong to $P_C$, by the carried construction
and the proof of \cref{prop:interface-perturbed-history-budget}.
Their excesses over $\varphi$, that is $\sA(\cdot)-N_0$, are
$E_1$, $E_1+E_2$ and $E_C=E_1+E_2+E_3$, respectively. Applying
\cref{lem:interface-perturbed-projection-excess} to each of the
three surfaces, with the pattern $P_C$, bounds the number of
projection triangles outside $S(P_C)$ in the union
\cref{eq:interface-perturbed-expanded-footprint} by
$E_1+(E_1+E_2)+(E_1+E_2+E_3)=3E_1+2E_2+E_3\leq3E_C$. Each
reference lozenge added by the tiling closure contains at least one
such triangle, which proves $|Q_C\setminus P_C|\leq3E_C$.

Here are the connectivity details. For $\cH_C$ and $\cS_C$, which
are graphs, the changed spins in a physical column $z$ form
a vertical interval one of whose ends is the reference horizontal
face $q_z$. That face belongs to $P_C$ by
\cref{lem:sos-column-component}, so the projected changes of the two
lower surfaces attach to $S(P_C)$. For the third step, consider an
$18$-connected component of one sign of $\cI_C-\cS_C$, say of
occupied cubes of $\cI_C$ that are vacant in $\cS_C$. By the
no-bubble condition, an admissible same-sign path in $\cI_C$ joins
it to the occupied end. This path must eventually leave the
component, since the component is finite: the changed-column
support and the transition sets of the component are finite. At the
first exit it reaches an occupied cube of $\cI_C$ outside the
component, which is then occupied in $\cS_C$ as well. If the exit
is a face-step, the crossed face is a boundary face of $\cS_C$; if
it is an edge-step, then of the two length-two face-paths around
the common cube edge, one crosses such a boundary face incident to
that edge. In either case the projection of the changed component
meets the projection of a face in the footprint
$Q_{\cS_C}(D_\delta)$ of the defect $\delta$ containing the exit
column. A common face of $\cS_C$ and $\varphi$ in this footprint
belongs to $F_\delta\subset P_C$; every other face of the footprint
is a new face of $\cS_C$, hence adjacent to a cube in
$\Omega(\cS_C,\varphi)$, whose projection is already attached.
Thus all third-stage projected changes attach to the already
connected set. Finally, a cube of $\Omega(\cI_C,\varphi)$ either
differs between $\cI_C$ and $\cS_C$ or between $\cS_C$ and
$\varphi$, so
$\Omega(\cI_C,\cS_C)\cup\Omega(\cS_C,\varphi)
 =\Omega(\cI_C,\varphi)\cup\Omega(\cS_C,\varphi)$, which proves
connectivity of \cref{eq:interface-perturbed-expanded-footprint}.
This argument uses closed projections; taking their reference
tiling closure gives the stated connected pattern.

The original components and their isolated reconstructions are
local by \cref{lem:interface-components}. Their projections
and the tiling-closure test use only the enlarged pattern.
More explicitly, combining components preserves their changed
cubes, not just their excesses. At the first inverse step the
reference dots and the binary columns of the defects are fixed. At the
second step, every dot and wall of a retained inner footprint
survives combination: a common face can be removed only by
an outer defect whose reference footprint contains it, while
a new face carries its creating outer defect. Such a defect
is in the same base component by
\cref{lem:sos-carried-components}. Hence the reference charges
and the full binary sequences on these columns are exactly
those in isolation. At the third step, a lower defect changing
$\cS(z)$ for $z$ in the two-square neighborhood of $D_\delta$ has
$q_z\in P_C\cap F_\delta$ and is again in the same component.
The reference graph columns and the third-stage column
sequences are therefore unchanged. Induction over the three
steps identifies each component's changed cubes, at each
intermediate stage, with its isolated changed cubes. In
particular its isolated expanded footprint still contains
their projections in the combined history.
Grouping overlaps is unique. Disjoint enlarged groups have
disjoint old footprints, so the insertion and deletion lemma
still applies. The extra restriction $Q\subset W$ implies
$P_C\subset W$ for every original component. Dropping the extra
restriction and using~\cref{eq:interface-positive-budget}
proves the a priori bound. Finally the projection of
either of two adjacent changed cubes contains the projected
lozenge of their common face. They cannot belong to expanded
components with disjoint expanded triangle supports. The same
observation assigns every changed plaquette to a unique expanded
component.
\end{proof}

\subsubsection{Auxiliary sorting and anchored lifting}
\label{sec:interface-perturbed-positive-completion}

The \RG input of \cref{def:rg-admissible} requires the weight to
be defined on every matching sector of the diagonal torus, not only
on the physical sector $\cM_\Gamma$ in which the three orientation
counts are those prescribed by the boundary conditions. The
interface histories, however, were defined through the physical
period lattice $\Gamma$ and its cylinder $\Z^3/\Gamma$, whereas a
matching in another sector lifts to a surface with a different
period lattice. We now define auxiliary histories on every sector,
with their own period lattices, and check that they satisfy the
same a priori bound; on the physical sector they reduce to the original
histories. The sectors in which some orientation count vanishes
require separate conventions, since a sorting step in a frozen
direction is degenerate.

Keep the diagonal torus fixed, and write $\Lambda=\Upsilon\Gamma$
for its period lattice. Every perfect matching $M$ lifts to a
monotone stepped surface. After fixing one lifted face, its
period homomorphism is
\[
 g_M:\Lambda\longrightarrow\Z^3,\qquad
 \Upsilon g_M(\lambda)=\lambda,
 \qquad \Gamma_M=g_M(\Lambda)\,.
\]
Thus $\Gamma_M$ is the period lattice of this matching, and
need not equal the prescribed physical lattice $\Gamma$.
The quotient still has two ends in the diagonal direction.
We use $N_i(M)$ for the number of faces normal to $x_i$;
$\sum_iN_i(M)=N_0$ in every matching sector.

\begin{lemma}[Sorting on an auxiliary sector]
\label{lem:interface-perturbed-auxiliary-sorting}
If $N_i(M)>0$, projection of $\Gamma_M$ onto the coordinates
other than $x_i$ is injective, has index $N_i(M)$, and defines
the usual infinite sorting fibers with their charges and period
offsets. If $N_i(M)=0$, a nonzero multiple of $e_i$ is a
period and the reference is constant along each $i$-fiber.

For $N_3(M)>0$, define the auxiliary histories by the usual
three inverse sorting steps, keeping the subsequent coordinate
monotonicities and the third-step no-bubble test. A step with
$N_i(M)=0$ is required to be the identity. For $N_3(M)=0$
retain only the empty history. These are conventions for the
auxiliary model; on $\cM_\Gamma$ they give exactly the
original interface histories. Uniformly over all $M$, they
satisfy the a priori bound
\begin{equation}\label{eq:interface-perturbed-auxiliary-budget}
 \sum_{\cX:\ Q_r\subset W\ (r\in\cX)}
                   e^{-tE(\cX)}
                   \leq e^{\epsilon_t^*|W|},
             \qquad \epsilon_t^*\leq Ce^{-ct}\,,
\end{equation}
where the expanded components and their energy are defined as in
\cref{prop:interface-perturbed-expanded-footprints}. The constants
do not require a lower bound on the nonzero $N_i(M)/N_0$.
\end{lemma}
\begin{proof}
We first prove the assertions on the period lattice, then the
counting bounds on a possibly degenerate transverse quotient, and
finally the a priori bound. Choose an oriented basis
$(\lambda_1,\lambda_2)$ of $\Lambda$. Summing oriented plaquette
areas over a fundamental surface gives
\[
 g_M(\lambda_1)\mathbin\times g_M(\lambda_2)
                         =(N_1(M),N_2(M),N_3(M))\,.
\]
This proves the area identity.
The projection of $\Gamma_M$ onto the coordinates other than $x_i$
has as its index the $i$th coordinate of this cross product, that
is, the $i$th minor of the matrix with rows
$g_M(\lambda_1),g_M(\lambda_2)$. If that minor is nonzero the
projection has rank two and no kernel, which is the injectivity
assertion. Otherwise the rational plane spanned by the periods
contains $e_i$, and hence a nonzero integer multiple of $e_i$
belongs to $\Gamma_M$. Monotonicity in $x_i$ together with
periodicity along this line make the reference constant on it.
At the beginning of inverse step $i$ the intermediate surface is
still monotone in $x_i$. Its signed normal-$i$ flux is $N_i(M)$,
since the signed flux is determined by the period lattice, which
all surfaces of the history share, so when $N_i(M)=0$ that surface
too has no normal-$i$ face and the sorting in direction $x_i$ has
nothing to do. Requiring this step to be
the identity is therefore consistent with the preceding steps.

Suppose $N_3(M)>0$. A monotone vertical fiber crosses the
reference once, so $M$, $\cH$ and $\cS$ are vertical graphs on the
cover, as in the physical sector. Finite modifications on the
cylinder have the prescribed ends on every active sorting fiber.
Their column magnetizations, local inclusion
constraints and first-exit end tests are exactly those used in
the three fiber propositions
(\cref{prop:sos-first-fiber,prop:interface-first-fiber} and the
\HMSOS fiber of \cref{sec:hmsos}). For each active direction the projected lattice has positive
index, so \cref{lem:quotient-one-step} supplies its enumeration and
a priori bound, uniformly even when the transverse quotient has loops
or parallel edges. The first-exit test still uses unchanged graph
columns, and the phase-tree proof applies with the two diagonal ends.
An inactive step has only its empty history and partition function one.

Carry lower footprints to the final reference and assign the
third-step complete physical-column footprints there, as before.
The proofs of independence and of the replacement-patch identity
$|W^{\cS}|=|W|+E_{12}$ use incidences and charges only.
Here is the successive use of the three a priori bounds. Enlarge the
constants in $\epsilon_s=Ce^{-cs}$ so that it bounds each
one-step a priori bound, and write $\cF^{\mathrm{aux}}_t(W)$ for the
full sum with old footprints in $W$. For a fixed lower history
$M\leftarrow\cH\leftarrow\cS$, put $E_{12}=E_1+E_2$.
Summing the third step costs at most
$e^{\epsilon_t(|W|+E_{12})}$. Its factor involving $E_{12}$
reduces the inverse temperature of the lower history to
$u=t-\epsilon_t$. Next fix $\cH$: the second-step replacement
patch has size $|W^{\cH}|=|W|+E_1$, so its a priori bound at $u$
similarly reduces the inverse temperature of the first step
to $u-\epsilon_u$. Thus
\[
 \begin{split}
 \cF^{\mathrm{aux}}_t(W)
 &\leq e^{\epsilon_t|W|}
       \sum_{\cH,\cS}e^{-u(E_1+E_2)}\\
 &\leq e^{(\epsilon_t+\epsilon_u)|W|}
       \sum_{\cH}e^{-(u-\epsilon_u)E_1}\\
 &\leq
   \exp[(\epsilon_t+\epsilon_u+
                       \epsilon_{u-\epsilon_u})|W|]\,.
 \end{split}
\]
The sums retain the indicated old-footprint restrictions; enlarging
the allowed one-step families to those in the replacement patches
only increases them. For large $t$ both shifted temperatures are
at least $t-1$, so the final coefficient is at most $Ce^{-ct}$.
An inactive step contributes one and obeys the same upper bound.
This proves \cref{eq:interface-perturbed-auxiliary-budget} for
old footprints.
Requiring their expanded footprints to lie in $W$ only restricts
the sum. The projection-excess proof and the connectivity proof
of \cref{prop:interface-perturbed-expanded-footprints} apply on
this cylinder and give the same expanded components. When $N_3=0$
the sum consists of its empty term, which proves the assertion
in the remaining sectors.
\end{proof}

For such a history write $\cI_{M,\cX}$ for the reconstructed
surface, with $\cI_{M,\varnothing}=M$. Define
\begin{equation}\label{eq:interface-perturbed-positive-primitive-completion}
 F^+_{t,\Lambda}(M)=
    \sum_{\cX\text{ auxiliary over }M}
       \exp(-tE(\cX)+U_{t,\Gamma_M}(\cI_{M,\cX}))\,.
\end{equation}
In a sector with $N_3=0$ this means $F^+(M)=e^{U_{t,\Gamma_M}(M)}$.
The perturbation in this definition is evaluated on the cylinder
of the matching itself. Changing it from sector to sector is
allowed in a canonical extension: its value on the prescribed
sector is fixed, and no physical meaning is required elsewhere.

\begin{lemma}[Anchored lifting of history components]
\label{lem:interface-perturbed-anchored-lifting}
Let $P$ be a connected partial matching on the diagonal torus.
The contact-marked connected components whose full footprint
is contained in $P$ have a weight-preserving correspondence
in any two matchings containing $P$, using their own
period lattices and the preceding auxiliary conventions.
The correspondence preserves subset histories, excesses,
mark masses, contact fractions, inclusion--exclusion marks,
and all witness tests. If $P$ has an embedded connected planar
lift, the same correspondence holds with its bulk components.
\end{lemma}
\begin{proof}
The proof has four parts. We first show that the lift of $P$ to
a monotone surface fixes the period map on the stabilizer $H_P$ of
a component of its planar lift, so that the two matchings
agree there. We then show that every piece of a history component,
including arbitrarily distant transitions, is anchored in one
component of the planar lift, so that it can
be copied with the period map on $H_P$. Next we check that the
copied components pass the same admissibility tests, in particular
the no-bubble test. Finally we treat the contact marks and the
bulk case. Here it is essential to use the expanded footprints,
which contain the projections of all changed cubes, rather than
just the old reference footprints: a transition far from the
reference surface is anchored through its changed cubes, and the
old footprints say nothing about where these lie.

Let $\widetilde P$ be one connected
component of the lift of $S(P)$ to the diagonal plane, and let
$H_P\subset\Lambda$ be its stabilizer. The face labels of $P$
determine its lift to a monotone surface from one face. Along any path in
$\widetilde P$ its three-dimensional increments are fixed by
these labels. In particular, for $\lambda\in H_P$ they fix
$g_M(\lambda)$. Two matchings containing $P$ have the
same homomorphism on $H_P$, after the initial lifts are aligned.
They may differ on periods outside $H_P$.

We spell out why all parts of a history component can be anchored in this
one component of the lift. Footprints connect the reference dots on neighboring sorting columns by their whole reference
walls. In a lower history a run of new faces on such a path
can be replaced by its carrier; this is the path construction
in \cref{lem:sos-carried-components}, performed on the
cover. The changed vertical intervals of $\cH$ and $\cS$ end at
their reference horizontal faces, so their cube projections
attach to these same lifted carriers.

For a third-step changed cube, follow a same-sign path until it
first reaches a cube agreeing with $\cS$. Such a finite first exit
exists by the end condition. Before this exit all its columns
belong to the same defect. The exit crosses a
face of $\cS$, or a face incident to the common edge in a $*$-step.
Within the footprint of this defect, that face
is connected back to the reference dot of the starting column by
the complete walls along a path of these columns. Common faces of
this footprint are in the assigned reference footprint; new faces have the lower carriers just
described. Thus the changed-cube component and its column anchors
belong to the same component of the lifted expanded footprint. This also covers
a path reaching its end without first leaving the transverse
support: it still leaves the finite changed set. No straight
segment across an unchanged gap between two transitions is
being assumed to lie in the footprint.

It follows that every transition in a changed column is anchored
to its dot in this component of the lift, even when the transition is
arbitrarily distant. Its integer displacement from that anchor,
the ordered transitions, and the reference footprints are
retained as data. All charge sums and local comparisons use these
anchored data. An identification within the component is an
element of $H_P$, on which the two period maps agree. Other
components of the lifted $P$ have disjoint triangle interiors.
Their changed cubes cannot be adjacent: the projection of
the common face would be contained in both footprints. This
checks energy incidences; the no-bubble test requires more.
At each inverse step, transverse supports in different lifted
components remain non-$*$-adjacent. Otherwise their complete
footprints would share a reference dot on the same intermediate
graph. A common dot has its reference face in both lifted
patterns; a new dot has its creating lower carrier in both,
by \cref{lem:sos-carried-components}. Either case joins the
two components of the lift, a contradiction. The same argument
excludes a new identification of distinct sorting columns.
These tests apply successively to the first and second inverse
steps, and then to the third, whose reference graph $\cS$ has
already been constructed. A lower defect changing $\cS$ on a
third-stage neighborhood shares its horizontal reference face with
$F_\delta$, so it too belongs to the same component of the lift.
Thus copying by the new period map preserves the
reference footprints as well as all charge and energy tests.

In particular, third-stage changed cubes in different lifted
components cannot be $*$-adjacent: their physical columns
would coincide or be $*$-adjacent. Every column in the
one-square neighborhood of a third-stage support, outside the support itself, remains
free of other third-stage defects. Transport a same-sign
first-exit path, including its possible edge-step, using the
anchored coordinates. Its finite portion within the support
is unchanged; its first exterior cube has the same spin and
lies on an unchanged column of the transported graph $\cS$.
Continue vertically there to the appropriate end. If the
path reaches its end inside the support, retain that portion
instead. This proves the no-bubble condition by the full
$*$-adjacency test of \cref{lem:polymer-factorization}.
All arguments are symmetric in the two matchings,
so the correspondence preserves the independent decomposition
and its energies.

There is no hidden dependence on a frozen direction. If $P$
contains no reference face of normal $i$, the intermediate
surface just before inverse step $i$ cannot acquire such a dot
in this component. Indeed it is monotone in $x_i$, all removed
reference faces lie in $P$, and conservation of signed normal-$i$
area equates its added and removed positive faces. The removed
number is zero. Hence this component has no active step $i$.
If $P$ contains no horizontal reference face, the two lower
graph changes are empty as well: each nonzero changed physical
column would have its horizontal reference face in its carrier.
Step three then has no reference dot. Thus a nonempty history
always contains such a horizontal face in its footprint. These
facts show that the identity-step and $N_3=0$ conventions do
not remove a component present in another matching containing $P$.

For a contact mark the projection $\Upsilon(N_X)$ of its carrier lies in $P$.
\Cref{def:interface-perturbation} transports its label,
mass and coefficient in the same component of the lift. A reference
mark is anchored at its marked contact face; an interaction mark
is anchored at a changed face of a required component. Its connected
carrier supplies the other anchors. The history correspondence
preserves every subset history, hence every contact predicate,
inclusion--exclusion integer and contact denominator.

These operations are inverse to one another. Enlargement cores,
tile closures and their overlap tests live on the common
diagonal pattern, so they too agree. When $P$ has an embedded
planar lift, $H_P=\{0\}$. The same construction is then lifting
to, or copying from, that finite planar component. All changed
cubes and all mark carriers have finite lifts there; their
heights may be large but finite. This proves the bulk assertion
without a lower bound on an auxiliary transverse period.
\end{proof}

\subsubsection{Contact marks and their a priori bound}
\label{sec:interface-perturbed-absolute-history}

The factor $e^{U_t(\cI_{\cX})}$ in
\cref{eq:interface-perturbed-positive-primitive-completion} couples
the components of a history: a label whose carrier touches the changed
cubes of several components contributes a term that depends on all
of them jointly. To obtain a hard-core gas of local objects, as the
\RG input requires, we expand this factor by inclusion--exclusion
into a reference part, depending on $\varphi$ alone, and a sum of
``marks,'' each attached to a set of components and vanishing unless
the label touches every one of them. The point of the expanded
footprints is that a label touches a component only if the projection
$\Upsilon(N_X)$ of its carrier meets the component's expanded footprint, so that a mark is
local on the tiling.

Use the expanded components of \cref{prop:interface-perturbed-expanded-footprints}
before any witness enlargement. For a compatible family $\cX$
of such components, let $\cI_{\cX}$ be its interface. Every
unoriented changed plaquette belongs to exactly one component, so
\begin{equation}\label{eq:interface-perturbed-disjoint-face-changes}
 \cI_{\cX}\xor\varphi
             =\mathop{\dot\bigcup}_{r\in\cX}
                           (\cI_{\{r\}}\xor\varphi)\,.
\end{equation}
For a contact label $X$ and a nonempty compatible set $\cV$
of components define the real, integer-valued coefficient
\[
 D_X(\cV)=\sum_{\cT\subset\cV}
               (-1)^{|\cV|-|\cT|}h_X(\cI_{\cT}),
                 \qquad \cI_\varnothing=\varphi\,.
\]
This is the M\"obius inversion of $\cV\mapsto h_X(\cI_{\cV})$
over the subsets of $\cV$; it is the part of the contact indicator
that is genuinely due to the joint presence of all components in
$\cV$. Such a mark requires all components in $\cV$ to be present.
Its weight at inverse temperature $t$ is
$b_XD_X(\cV)e^{-tm_X}$ and its reference footprint
is $A_X^\varphi$. A reference mark is the marked label
$(X,p)$ from \cref{prop:interface-perturbed-reference-potential},
of weight $b_Xe^{-tm_X}/|\varphi\cap N_X|$ and the same
footprint. It requires no history component.

\begin{proposition}[Contact marks and their a priori bound]
\label{prop:interface-perturbed-absolute-history}
The primitive sum in~\cref{eq:interface-perturbed-positive-primitive-completion}
equals the sum over compatible families of expanded components
and finite multisets of the two kinds of contact marks just
defined. A mark of multiplicity $l$ has its weight to the
power $l$, divided by $l!$. There is $\delta_t\leq Ce^{-ct}$
such that, for every finite set $W$ of reference plaquettes, the sum of absolute
weights of all these objects whose footprints lie in $W$ is
at most
\begin{equation}\label{eq:interface-perturbed-absolute-history-budget}
             \exp\big(
             (3\delta_t+\epsilon^{\mathrm I}_{t-\delta_t})|W|
                         \big)\,.
\end{equation}
This includes the reference factor and every interaction joining
distinct components. If $E$ is the total history excess and $M$
the sum of the mark masses, with multiplicities,
each weight is $B e^{-t(E+M)}$ with $t$-independent
$B$. The displayed bound is on the physical periodic quotient and also holds for
finite planar windows.
\end{proposition}
\begin{proof}
For fixed $X,\cX$, let
\[
 \cR_X=\{r\in\cX:
                (\cI_{\{r\}}\xor\varphi)\cap N_X\ne\varnothing\}\,.
\]
\Cref{eq:interface-perturbed-disjoint-face-changes} implies that
$h_X(\cI_{\cT})$ depends only on
$\cT\cap\cR_X$. Hence $D_X(\cV)=0$
unless $\cV\subset\cR_X$. Each component in
$\cR_X$ has a changed face whose projection lies in
both its expanded footprint and $\Upsilon(N_X)$. These expanded footprints have disjoint triangle
supports, so $|\cR_X|\leq Cm_X$. In particular
\begin{equation}\label{eq:interface-perturbed-mobius-budget}
 \sum_{\varnothing\ne\cV\subset\cX}
               |D_X(\cV)|
    \leq 3^{|\cR_X|}-1
    \leq 3^{Cm_X}\one_{\{N_X\cap(\cI_{\cX}\xor\varphi)
                                             \ne\varnothing\}}\,.
\end{equation}
The first bound follows from $|D_X(\cV)|\leq2^{|\cV|}$, since
$D_X(\cV)$ is a sum of $2^{|\cV|}$ terms of modulus at most one,
and from $\sum_{\varnothing\ne\cV\subset\cR_X}2^{|\cV|}
=3^{|\cR_X|}-1$. The indicator in the second bound records that
no mark is produced by a label whose carrier misses every changed
face.

The finite inclusion--exclusion identity gives
\[
 h_X(\cI_{\cX})-h_X(\varphi)
           =\sum_{\varnothing\ne\cV\subset\cX}
                                                    D_X(\cV)\,.
\]
Indeed, if we include $\cV=\varnothing$, whose coefficient is
$h_X(\varphi)$, and interchange the two finite subset sums, the
coefficient of $h_X(\cI_{\cT})$ is
\[
 \sum_{\cV:\,\cT\subset\cV\subset\cX}
                 (-1)^{|\cV|-|\cT|}
       =(1-1)^{|\cX\setminus\cT|}.
\]
Only $\cT=\cX$ survives. Likewise, if $\cV$ contains a component
$r\notin\cR_X$, then $h_X(\cI_{\cT})=h_X(\cI_{\cT\setminus\{r\}})$
for every $\cT\subset\cV$, so pairing each subset containing $r$
with the subset omitting it cancels the sum and proves the
vanishing assertion for $D_X$ above. This explains why the marks
retain all orders of the interaction between components: the
identity is exact, not a first-order approximation. Multiplying
by $a_t(X)$ and summing over $X$ therefore gives the exact
interaction difference,
\[
 U_t(\cI_{\cX})-U_t(\varphi)
   =\sum_X\sum_{\varnothing\ne\cV\subset\cX}
        b_XD_X(\cV)e^{-tm_X}\,,
\]
a sum of mark weights. Exponentiating a sum of commuting terms
$w_j$ as $\prod_j\sum_{l\geq0}w_j^l/l!$ produces precisely the
finite multisets of interaction marks, each mark of multiplicity
$l$ carrying its weight to the power $l$ divided by $l!$. The
reference factor $e^{U_t(\varphi)}$ is expanded in the same way:
by the contact-sharing identity
$U_t(\varphi)=\sum_{(X,p)}b_Xe^{-tm_X}/|\varphi\cap N_X|$ of
\cref{prop:interface-perturbed-reference-potential}, it
exponentiates to the multisets of reference marks. Multiplying by
$e^{-tE(\cX)}$ and summing over histories proves the proposed
identity, once absolute convergence is checked.

For the absolute convergence and the a priori bound, the factor
$3^{Cm_X}$ in \cref{eq:interface-perturbed-mobius-budget} is
absorbed in the reserved $e^{tm_X/4}$ of
\cref{eq:interface-contact-norm}, for $t$ large. Anchoring each
label at a face of $\cI_{\cX}\xor\varphi$ that its carrier
meets, as in \cref{eq:interface-contact-variation}, the total
absolute weight of all interaction marks at a fixed history is
consequently bounded by
\[
 \sum_X |a_t(X)|\sum_{\varnothing\ne\cV\subset\cX}
                     |D_X(\cV)|
            \leq\delta_t|\cI_{\cX}\xor\varphi|\,.
\]
If the expanded footprints of $\cX$ lie in $W$, all removed
oriented reference faces lie in $W$, and hence
$|\cI_{\cX}\xor\varphi|\leq2|W|+E$, as in the proof of
\cref{prop:interface-perturbed-history-budget}. Reference marks
whose footprints lie in $W$ have their marked face in $W$, and
summing their contact fractions costs at most $\delta_t|W|$, by
the same anchored contact bound without the inclusion--exclusion
factor. Since the sum of the absolute weights of all multisets is
the exponential of the sum of the absolute weights of the
individual marks, the marks at a fixed history contribute at most
$\exp(2\delta_t|W|+\delta_tE+\delta_t|W|)
=\exp(3\delta_t|W|+\delta_t E)$. The history sum is therefore
bounded by
\[
 e^{3\delta_t|W|}
       \sum_{\cX:\ Q_r\subset W}
                       e^{-(t-\delta_t)E(\cX)}
 \leq e^{(3\delta_t+\epsilon^{\mathrm I}_{t-\delta_t})|W|}\,,
\]
using the a priori bound of
\cref{prop:interface-perturbed-expanded-footprints} at
$t-\delta_t$. This proves
\cref{eq:interface-perturbed-absolute-history-budget} and
justifies all the series rearrangements above. Finally, the
signs, contact fractions, inclusion--exclusion integers and
multiset factorials are temperature independent, and the only
dependence on $t$ is through $e^{-tE}$ for the history and
$e^{-tm_X}$ for each mark occurrence; this proves the last
assertion.
\end{proof}

\section{Grassmann input for the three histories}
\label{sec:model-grassmann}

\begin{remark*}[Reader's guide]
\normalfont
The purpose of this section is to turn the preceding history
expansions into the hypotheses of the abstract dimer theorem,
\cref{thm:coefficient-rg}. The three resulting input statements are
\cref{thm:hmsos-rg-input,thm:sos-rg-input,thm:interface-rg-input}.
There are three tasks.
\begin{enumerate}[(i)]
\item \emph{Represent the weights and sources exactly.}
\Cref{lem:pattern-algebra} converts a disjoint-pattern polynomial
into a Grassmann integral, with the exact edge-source factors and
vertex gauge identities. This common algebra allows signed
coefficients and singular Kasteleyn matrices.
\item \emph{Define compatible coefficients on all torus sectors.}
For \HMSOS, \cref{lem:periodic-locality} proves independence of the
tiling containing the pattern within the prescribed sector and agreement of small
torus coefficients with their planar counterparts.
\Cref{prop:canonical-completion,prop:exact-sector} then give a
positive extension to all sectors and recover the physical sector
by exact coefficient extraction. For \SOS, the additional
confinement argument is \cref{lem:sos-small-locality}, leading to
\cref{prop:sos-canonical-fermions}. For the Ising Interface model,
\cref{prop:interface-perturbed-positive-completion} constructs a
positive matching weight even though its pattern coefficients may
be signed. It includes the contact perturbations needed later for
full Ising. These constructions retain the absolute
\texttt{thin}/\texttt{fat} bounds from the history enumeration.
\item \emph{Prepare the witnesses for the interacting law.}
\Cref{lem:marked-union-envelopes} controls collections of ordinary
patterns overlapping a witness. Its \SOS and interface versions
are \cref{lem:sos-marked-input,lem:interface-marked-input}.
These are coefficient and counting estimates; the interacting
probability bound is subsequently obtained from
\cref{thm:interacting-witness-tail} and used to reconstruct the
surfaces in \cref{sec:transfer}.
\end{enumerate}
\end{remark*}

Throughout this section a pattern $P$ is a nonempty partial matching
of the honeycomb graph with connected triangle support; $|P|$ counts
its lozenges (equivalently, its edges) and $S(P)$ is the set of
triangles incident to those lozenges. Since every triangle of the
triangular lattice is a vertex of the honeycomb graph, $S(P)$ is
precisely the set of vertices covered by $P$, and $|S(P)|=2|P|$.
In particular, for two patterns contained in a common matching $M$,
disjointness of their supports is the same as vertex-disjointness.
For a matching $M$ the common pattern polynomial is
\[
 F(M)=\sum_{\substack{\{P_i\}:\ S(P_i)\text{ pairwise disjoint}\\
                                  P_i\subset M}}\prod_i c(P_i),
 \qquad
 V_A(\psi)=\sum_P c(P)
       \prod_{e=(b,w)\in P}K_y(e)e^{A(e)}\psi_b\psi_w\,.
\]
Here $K_y$ is the Kasteleyn matrix including the torus twist,
$A$ is a collection of edge sources, and the sum defining $F(M)$
includes the empty collection, with weight one. Thus $F$ is the
partition function of a hard-core gas of patterns inside $M$, whereas
$V_A$ is the Grassmann potential obtained by replacing the indicator
$\one_{\{P\subset M\}}$ of each pattern by the monomial of its edges.
The precise common requirements and their Grassmann formulation are in
\cref{def:rg-admissible,prop:fermion-input}.

\begin{remark}[Why a pattern polynomial rather than an exponential]
\label{rem:pattern-polynomial-vs-exponential}
In \cite{LaslierLubetzky24} the interacting tiling law was written
as a tilt of the uniform measure by
$\exp(\sum_x\sum_r\fg_r(\varphi\restriction_{B(x,r)}))$ with
$\|\fg_r\|_\infty$ decaying exponentially in $r$, and the first task
there was to expand this exponential into a sum over patterns.
Here the hard-core sum $F$ is instead the primary object: for each
of our models it is literally the partition function of the fiber
of surfaces above the tiling, which the history expansions of
\cref{sec:preimage-witnesses} write as an exact hard-core gas of
enlarged patterns. A uniformly convergent cluster expansion for
$\log F$ would recover a sum of local functions, provided suitable
smallness and summability bounds held~\cite{KoteckyPreiss86}.
Our bounds on the \texttt{fat} coefficients do not ensure this
(compare \cref{rem:individual-coefficient-bounds}). We therefore
retain $F$: the first integration in \cref{sec:first-integration}
uses the rarity of a \texttt{fat} pattern to compensate its possibly
large coefficient. The exact identity in \cref{lem:pattern-algebra}
then represents $F$ by $e^{V_A}$ in the Grassmann algebra, where
overlapping supports vanish by nilpotence.
\end{remark}

\begin{lemma}[Common pattern algebra]\label{lem:pattern-algebra}
Let $c(P)$ be coefficients on nonempty partial matchings of a finite
bipartite graph, independent of its edge weights and sources, and let
$F$ be the disjoint-pattern polynomial above. Write
\[
 \Psi_{P,A}=\prod_{e=(b,w)\in P}K(e)e^{A(e)}\psi_b\psi_w,
 \qquad V_A=\sum_Pc(P)\Psi_{P,A},\qquad
 J_A=\sum_eK(e)(e^{A(e)}-1)\psi_b\psi_w\,.
\]
In paired edge order, with the Berezin orientation giving $\det K$
for the free integral,
\[
 \int\cD\psi\,e^{\sum_eK(e)e^{A(e)}\psi_b\psi_w+V_A}
 =\sum_M\sgn(M)F(M)\prod_{e\in M}K(e)
                                      e^{\sum_{e\in M}A(e)}\,.
\]
The identity permits signed coefficients and singular $K$. For the
free dimer law $\mu_0$ of a finite planar region all of whose inner
faces are hexagons, so that every perfect matching has the same
Kasteleyn phase, with covariance $G$, write $\cE_G$ for the
normalized Gaussian Grassmann expectation with covariance $G$. Then
\[
 \mu_0\big(F(M)e^{\sum_{e\in M}A(e)}\big)
                  =\cE_G(e^{J_A+V_A})\,.
\]
This also holds for the infinite-volume free law when patterns
and source edges lie in a finite window.

For sources on an edge set $F_0$, the absolute Taylor factor of an
interaction monomial at radius $\rho$ is exactly
$e^{\rho|P\cap F_0|}$. Shifting all sources incident to $v$ by $t$
multiplies that monomial by $e^{t\one_{\{v\in S(P)\}}}$ and
the full partition function by $e^t$. These identities survive
fixed-count extraction. The free source $J_A$ is kept separately.
\end{lemma}
\begin{proof}
We first expand both exponentials in the integrand, then identify
the terms surviving the Berezin integration, and finally deduce the
normalized, infinite-volume, source and gauge statements.

Write $q_e=K(e)e^{A(e)}\psi_b\psi_w$ for the quadratic monomial of
the edge $e=(b,w)$. Each $q_e$, and each $\Psi_{P,A}$, is a product
of an even number of Grassmann generators; hence all these monomials
commute with one another, and each squares to zero because it
contains some generator twice. For commuting nilpotent elements the
exponential series terminates after its linear term, so that
\[
 e^{\sum_eq_e}=\prod_e(1+q_e),\qquad
 e^{V_A}=\prod_P\big(1+c(P)\Psi_{P,A}\big)\,.
\]
Expanding both products, the integrand becomes a sum over pairs
consisting of an edge set $E$ and a family of patterns
$\{P_1,\ldots,P_k\}$, with term
\[
 \bigg(\prod_ic(P_i)\bigg)\prod_{e\in E}q_e\prod_i\Psi_{P_i,A}\,.
\]
A product of generators vanishes as soon as a generator is repeated,
and the Berezin integral of a monomial is nonzero only when every
generator appears exactly once. Thus the term survives integration
precisely when the edges of $E$ together with those of the $P_i$
cover every vertex exactly once: the $P_i$ are pairwise
vertex-disjoint, hence have pairwise disjoint supports $S(P_i)$;
no edge of $E$ meets $P_*=\bigcup_iP_i$; and
$M=E\cup P_*$ is a perfect matching. Conversely, for any perfect
matching $M$ and any family of patterns with disjoint supports
contained in $M$, the set $E=M\setminus P_*$ is forced. The surviving
terms are therefore indexed by pairs $(M,\{P_i\})$ with
$P_i\subset M$ and $S(P_i)$ pairwise disjoint, which is exactly the
index set of $F(M)$.

Consider one such term. Its monomial is the product, over all edges
$e\in M$, of $K(e)e^{A(e)}\psi_b\psi_w$, each edge appearing once,
regardless of whether $e$ came from $E$ or from a pattern. The pairs
$\psi_b\psi_w$ are even and commute, so this product can be brought
to the paired edge order of $M$ without any sign, and integrating it
in the fixed Berezin orientation gives the permutation sign
$\sgn(M)$ appearing in $\det K=\sum_M\sgn(M)\prod_{e\in M}K(e)$.
The contribution of the term is therefore
\[
 \sgn(M)\bigg(\prod_i c(P_i)\bigg)
                 \prod_{e\in M}K(e)e^{A(e)}\,.
\]
Summing over the allowed pattern families inside the fixed $M$
produces exactly $F(M)$ as the coefficient of
$\sgn(M)\prod_{e\in M}K(e)e^{\sum_{e\in M}A(e)}$, which is the
displayed identity. In particular each matching edge carries one
weight and one source factor, even when it belongs to an interaction
pattern: this is the reason for including $K(e)e^{A(e)}$ in
$\Psi_{P,A}$. The argument is an identity between polynomials in the
entries $K(e)$ and the coefficients $c(P)$, and it uses neither the
sign of the coefficients nor the invertibility of $K$; this is what
permits singular twists later on.

For the normalized statement, take a planar region whose inner
faces are hexagons. The matching terms have a common phase $\theta$:
\[
 \sgn(M)\prod_{e\in M}K(e)=\theta\prod_{e\in M}|K(e)|\,.
\]
The product on the right is the unnormalized weight of $M$ under
the free law $\mu_0$.
Hence dividing the identity by its value at $c\equiv0$, $A\equiv0$,
which is $\det K=\theta\sum_M\prod_{e\in M}|K(e)|$, cancels the
common phase and gives
\[
 \frac{\int\cD\psi\,e^{\sum_eq_e+V_A}}{\det K}
 =\sum_M\mu_0(M)F(M)e^{\sum_{e\in M}A(e)}
 =\mu_0\big(F(M)e^{\sum_{e\in M}A(e)}\big)\,.
\]
On the other hand $\sum_eq_e=\sum_eK(e)\psi_b\psi_w+J_A$, so the
left-hand side is the normalized Gaussian Grassmann expectation
$\cE_G(e^{J_A+V_A})$ with covariance $G=K^{-1}$. In infinite volume,
restrict the patterns and the source edges to a finite window. The
expansion of $F(M)e^{\sum A}$ is then a finite linear combination of
indicators $\one_{\{P_1\cup\cdots\cup P_k\subset M\}}$ times fixed
source factors, and the $\mu_0$-probability of each such event is
the determinant of the corresponding minor of the inverse Kasteleyn
matrix, by the local statistics formula of
Kenyon~\cite{Kenyon97}; the Gaussian Grassmann expectation of the
corresponding monomial is the same determinant by Wick's rule, so the
identity persists.

Next, the source statements. Expanding $e^{A(e)}$ for each
$e\in P\cap F_0$ as an absolutely convergent series and bounding
$|e^{A(e)}|\leq e^{\rho}$ on $|A(e)|\leq\rho$ gives the absolute
Taylor factor $e^{\rho|P\cap F_0|}$; edges of $P$ outside $F_0$
carry no source. For the gauge identities, shift $A(e)$ by $t$ on
every edge $e$ incident to a vertex $v$. A partial matching $P$
contains at most one edge at $v$, and exactly one when
$v\in S(P)$, so $\Psi_{P,A}$ is multiplied by
$e^{t\one_{\{v\in S(P)\}}}$; a perfect matching contains exactly one
edge at $v$, so every term of the matching sum, and hence the full
partition function, is multiplied by $e^t$. Since the identity holds
termwise in the matchings, extracting the coefficient of a prescribed
orientation count (which is a function of $M$ alone) commutes with
it. Finally, if the coefficients $c(P)$ are represented by convergent
sums over internal labels, truncate these sums. There are only finitely
many patterns in the graph or window, and both sides are polynomials
in their total coefficients, so the identity passes to the limit.
No positivity or absolute convergence of the internal sums is needed
for this algebraic passage.
\end{proof}

We state the common coefficient conditions before verifying them for
the three histories. The parameter $\beta$ indexes a family and need
not be an inverse temperature. The explicit rates are sufficient
rates shared by our applications. A support is
$(\epsilon,R)$-\texttt{fat} if, apart from at most an $\epsilon$
fraction of its area, it is covered from inside by a union of balls
of radius $R$. Balls, supports and boundary layers are taken in
the quotient on a torus. Rooted counts count patterns containing a
specified triangle.

\begin{definition}[Admissible coefficient families]\label{def:rg-admissible}
An admissible family consists of real planar coefficients $c_\beta(P)$
and real torus coefficients $c_{\beta,\Lambda}(P)$ on a collection
of honeycomb tori containing the diagonal tori under consideration,
$\Lambda=\Z^2/\sB_n\Z^2$ of bounded shape:
$\|n^{-1}\sB_n\|+\|n\sB_n^{-1}\|\leq C$.
The family is defined on all matching sectors of these tori, including
frozen sectors, and satisfies the following conditions. Constants are
uniform over the allowed shapes and over the stated compact
ranges of any further parameters.
\begin{enumerate}[label=(R\arabic*),ref=R\arabic*]
\item\label[condition]{it:rg-positive} \emph{Positive matching weights and fixed coefficients.}
The coefficients commute with translations and are independent of
auxiliary edge weights, boundary characters, and edge sources.
The exact matching weight is
\begin{equation}\label{eq:rg-pattern-completion}
 F_{\beta,\Lambda}(M)=
 \sum_{\substack{\{P_i\}:\ S(P_i)\text{ pairwise disjoint}\\P_i\subset M}}
                         \prod_i c_{\beta,\Lambda}(P_i)>0\,,
\end{equation}
on every matching. Individual coefficients
are allowed to have either sign.

\item\label[condition]{it:rg-envelopes} \emph{Absolute coefficient and entropy bounds.}
There are nonnegative envelopes $q^\mathrm{t},q^\mathrm{f}$, and constants
$c_0,C_0>0$, with
\[
 |c_{\beta,\Lambda}(P)|\leq q^\mathrm{t}_{\beta,\Lambda}(P)
                                  +q^\mathrm{f}_{\beta,\Lambda}(P),\qquad
 q^\mathrm{t}_{\beta,\Lambda}(P)\leq e^{-c_0\beta^{2/5}|P|},\qquad
 q^\mathrm{f}_{\beta,\Lambda}(P)\leq e^{\epsilon_\beta|P|}\,.
\]
The \texttt{thin} support family has at most $C_0^m$ rooted patterns of size
$m$. Every \texttt{fat} pattern has $|P|\geq B_\beta\geq c_0\sqrt\beta$, is
$(\epsilon_\beta,R_\beta/2)$-\texttt{fat}, and belongs to a family with at
most $e^{(\epsilon_\beta+C/n)m}$ rooted patterns of size $m$.
Here $\epsilon_\beta\to0$ and $R_\beta\to\infty$ as $\beta\to\infty$.
The planar coefficients obey the same bounds without $C/n$.

\item\label[condition]{it:rg-locality} \emph{Common bulk coefficients.}
There is a fixed $c_*>0$ such that every connected torus pattern
with fewer than $c_*n$ edges lifts to a planar pattern and has exactly
its planar coefficient.
The planar coefficients are the same for all volumes and period
sequences under comparison. All exceptional or noncontractible inputs
therefore have at least $c_*n$ edges. Finite planar windows are
obtained by restricting the complete pattern support to the window.

\item\label[condition]{it:rg-sources} \emph{Exact sources and sector algebra.}
For $K_y(e)=\omega_r y_r$ on an orientation-$r$ edge, $y_1=1$, use
paired edge order and the unnormalized polynomial
\begin{equation}\label{eq:rg-berezin}
 \fZ_{\beta,\Lambda}(y,A)=\int\cD\psi\,
 \exp\bigg(\sum_eK_y(e)e^{A(e)}\psi_b\psi_w+
 \sum_Pc_{\beta,\Lambda}(P)
                 \prod_{e=(b,w)\in P}K_y(e)e^{A(e)}\psi_b\psi_w\bigg)\,.
\end{equation}
It is the signed matching sum with weight \cref{eq:rg-pattern-completion}
and source $\exp(\sum_{e\in M}A(e))$; fixed counts are extracted
before taking a ratio. This identity includes singular twists.
Thus the source factor on an interaction pattern is exactly
$\prod_{e\in P}e^{A(e)}$, the free source is
$J_A=\sum_eK_y(e)(e^{A(e)}-1)\psi_b\psi_w$, and vertex gauge
identities hold for the full generating function.
In particular, the absolute Taylor factor at radius $\rho_0$ is
$e^{\rho_0|P|}$, with all marks on edges of $P$; a source bound
only for each fixed set of locations is not a substitute.
\end{enumerate}
The identities in \cref{it:rg-sources} follow algebraically from
\cref{it:rg-positive}; they specify the representation and source
convention used throughout the proof. Auxiliary weights are varied
in compact liquid sets and normalized by $\max_r\omega_r=1$
when estimating absolute monomials.

For the witness conclusion only, impose the following additional
condition on a translation-covariant family $\scW_\Lambda$ of
connected patterns of size at least $B_\beta$, counted without duplicates.
\begin{enumerate}[label=(W),ref=W]
\item\label[condition]{it:rg-marked} \emph{Absolute overlapping-collection bounds.}
Put $q=q^\mathrm{t}+q^\mathrm{f}$. For every root $x$, size threshold $s\geq1$
and connected partial matching $U$, form the nonnegative sum
\begin{equation}\label{eq:rg-absolute-insertion}
 \overline b_{x,s}(U)=
 \sum_{\substack{W\in\scW_\Lambda:\ x\in S(W),\ |W|\geq s\\W\subset U}}
 \ \sum_{\substack{\{P_i\}:\ S(P_i)\text{ pairwise disjoint}\\
      P_i\subset U,\ S(P_i)\cap S(W)\ne\varnothing\\
                       W\cup\bigcup_iP_i=U}}\prod_i q(P_i)\,.
\end{equation}
The empty inner collection is allowed. This sum has envelopes
$\overline b^\mathrm{t}+\overline b^\mathrm{f}$ with
\[
 \overline b^\mathrm{t}_{x,s}(U)\leq e^{-c_0\beta^{3/10}|U|},\qquad
 \overline b^\mathrm{f}_{x,s}(U)\leq e^{\epsilon'_\beta|U|}\,.
\]
Both vanish unless $x\in S(U)$ and $|U|\geq\max(s,B_\beta)$.
The \texttt{thin} family has at most $C_0^m$ choices of size $m$ through $x$.
The \texttt{fat} family is $(\epsilon'_\beta,R'_\beta)$-\texttt{fat} and has at most
$e^{(\epsilon'_\beta+C/n)m}$ such choices, where
$\epsilon'_\beta\to0$ and $R'_\beta\to\infty$.
All bounds are uniform in $x,s$ and include noncontractible supports.
\end{enumerate}
\end{definition}

\begin{remark*}
The conditions involve only the total coefficients
$c_{\beta,\Lambda}(P)$; the proof below never refers to their internal
structure. In the applications, each coefficient is a sum over
labelled terms (defect collections, histories, or contour-decorated
histories), and it is through these terms that \cref{it:rg-envelopes}
is verified: the envelope $q^{\mathrm t}$ is the sum of the absolute
values of the terms whose energy pays for the area of $P$, and
$q^{\mathrm f}$ the sum over the remaining terms, whose supports are
\texttt{fat}. No cancellation between terms is used; for the \HMSOS
and \SOS inputs the terms are nonnegative, and for the perturbed
interface input of \cref{sec:ising}, in which closed contours have
been integrated out, they are signed.
\end{remark*}

For diagonal tori we also write $F_{\beta,\Gamma}$ for
$F_{\beta,\Lambda}$, with $\Lambda=\Upsilon\Gamma$ in the fixed
diagonal lattice coordinates. The notation $\nu_{\beta,\Gamma}$
denotes its normalized restriction to the prescribed counts.
\Crefrange{it:rg-positive}{it:rg-sources} concern the original
coefficients, before any infrared integration; no critical flow,
covariance limit or Gaussianity is assumed. The extra condition
\cref{it:rg-marked} is needed only for the witness conclusion.

The following is their immediate Grassmann formulation.

\begin{proposition}[Pattern and fermion coefficients]\label{prop:fermion-input}
Write $\Psi_{P,A}=\prod_{e=(b,w)\in P}K_y(e)e^{A(e)}\psi_b\psi_w$
in paired edge order. Under \crefrange{it:rg-positive}{it:rg-sources},
the interaction potential $V_A=\sum_Pc(P)\Psi_{P,A}$ has balanced connected supports, uniformly
in the volume and unitary twists. At zero source, after converting
$|P|$ to $|S(P)|=2|P|$, it has the following properties:
\begin{enumerate}[(i)]
\item \texttt{Thin} labelled coefficients are at most $\epsilon_\beta^{|S|}$,
where $\epsilon_\beta\to0$.
\item \texttt{Fat} supports have diverging minimum size, are
$(\epsilon_\beta,R_\beta/2)$-\texttt{fat} with $\epsilon_\beta\to0$ and
$R_\beta\to\infty$, and have coefficient bound $e^{\epsilon_\beta|S|}$.
\item Their rooted labelled count is at most
$e^{(\epsilon_\beta+C/n)|S|}$. Labels distinguish patterns
with the same vertex support and their \texttt{thin}/\texttt{fat} parts;
no internal decomposition of a coefficient is required.
\item The source dependence is entire, with absolute Taylor factor
$e^{\rho|P|}$ at radius $\rho$, and the exact generating function
obeys vertex gauge covariance. The free-source nodes have size
$O(\rho)$ and are treated separately.
\end{enumerate}
The planar statements omit $C/n$. The fixed positive source radius
is chosen after reserving a negative exponent in the \texttt{fat}-node bound;
it is not asserted that the source contribution tends to zero with
$\beta$ at that radius.
\end{proposition}
\begin{proof}
\Cref{lem:pattern-algebra} gives the exact Grassmann expansion,
absolute source Taylor factor and gauge covariance, including
singular twists and the planar finite-window identity. To split the
total coefficients, put $a(P)=e^{-c_0\beta^{2/5}|P|}$ and
\[
 c^{\mathrm t}(P)=\sgn(c(P))\min\{|c(P)|,a(P)\},\qquad
 c^{\mathrm f}(P)=c(P)-c^{\mathrm t}(P)\,.
\]
Then $|c^{\mathrm t}(P)|\leq a(P)$ and
$|c^{\mathrm f}(P)|=(|c(P)|-a(P))_+\leq q^{\mathrm f}(P)$,
since $q^{\mathrm t}(P)\leq a(P)$. The fat part therefore has the
prescribed support family and bounds. The thin part may use patterns
from either original family; their union still has at most $C_1^m$
rooted patterns of size $m$, for a fixed $C_1$ and sufficiently large
$\beta,n$. This splitting depends only on $c(P)$ and $|P|$, so it
preserves the common bulk coefficients of \cref{it:rg-locality}.
The first three bounds now follow using
$|K_y(e)|\leq1$ (recall that the largest weight is normalized to one
and that the twists are unimodular) and $|S(P)|=2|P|$: a
\texttt{thin} coefficient $e^{-c_0\beta^{2/5}|P|}$ becomes
$(e^{-c_0\beta^{2/5}/2})^{|S|}$, which is of the form
$\epsilon_\beta^{|S|}$ with $\epsilon_\beta\to0$ (we use the same
symbol for the various quantities tending to zero with $\beta$), and
the \texttt{fat} bound and rooted count are rewritten in the same
way, using $|P|\leq|S|$.
Finally, $|e^{A(e)}-1|\leq e^\rho-1=O(\rho)$ on a fixed small source
polydisc, giving the separate free-source bound.
\end{proof}

\subsection{One step: HMSOS input}

\begin{theorem}[\HMSOS input to the \RG theorem]\label{thm:hmsos-rg-input}
Fix a compact subset of positive physical slopes. For sufficiently
large $\beta$, all sufficiently large $n$ depending on $\beta$,
and admissible periods
$\Gamma=\langle(-n,0,\Theta_1),(0,-n,\Theta_2)\rangle$ in this slope range,
the one-step \Tetris image of the \HMSOS ensemble, with the diagonal
rooting of \cref{sec:hmsos-rooting}, has law
\[
 \nu^{\textsc{hmsos}}_{\beta,\Gamma}(M)
 =\frac{Z_{\beta,M}}{\sum_{M'\in\cM_\Gamma}Z_{\beta,M'}},
 \qquad M\in\cM_\Gamma\,.
\]
It admits a translation-covariant canonical extension $F_{\beta,\Gamma}$
to all tilings of the diagonal torus of the form
\cref{eq:rg-pattern-completion}, with $F_{\beta,\Gamma}(M)=Z_{\beta,M}$
for $M\in\cM_\Gamma$. Here $P_i$ are connected lozenge patterns and $S(P_i)$ are their
triangle supports. The coefficients have a nonnegative splitting
$c=c^\mathrm{t}+c^\mathrm{f}$ with parameters
$\tau_\beta,B_\beta,R_\beta\to\infty$ and
$\epsilon_\beta\to0$ such that
\[
 c^\mathrm{t}_{\beta,\Gamma}(P)\leq e^{-\tau_\beta|P|},
 \qquad c^\mathrm{f}_{\beta,\Gamma}(P)\leq e^{\epsilon_\beta|P|}\,.
\]
There are at most $C_0^s$ rooted \texttt{thin} patterns of size $s$.
Every nonzero \texttt{fat} term has $|P|\geq B_\beta$, is
$(\epsilon_\beta,R_\beta/2)$-\texttt{fat}, and belongs to a family of at
most $e^{(\epsilon_\beta+C/n)s}$ rooted patterns of size $s$.
These coefficients sum all internal decompositions absolutely.
Small torus patterns agree exactly with their translation-invariant
planar coefficients, as specified in \cref{lem:periodic-locality}.

For liquid dimer weights normalized by $\max_r\omega_r=1$, the
resulting unnormalized Grassmann generating function satisfies
the \texttt{thin}/\texttt{fat} bounds and source conditions listed in
\cref{prop:fermion-input}, with the stated periodic entropy
correction. The pattern coefficients $c_{\beta,\Gamma}$ are
independent of the torus twists; twist and source factors enter
explicitly through the Kasteleyn edge monomials. The generating
function is analytic in the edge sources and obeys vertex gauge
covariance.
Its projection to the prescribed counts is exactly
\cref{eq:exact-sector-ratio}, including at singular twists.
These coefficients satisfy \crefrange{it:rg-positive}{it:rg-sources}
of \cref{def:rg-admissible}; the six-class witness family
of \cref{sec:hmsos-rg-verification} satisfies \cref{it:rg-marked}.
The relation to the physical rooted height law is
\cref{prop:physical-hmsos-law}.
\end{theorem}
The coefficient verification is in \cref{sec:hmsos-rg-verification};
the proof is completed after \cref{prop:physical-hmsos-law}.

The statement has three layers, which we comment on in turn, as
they recur for \SOS and for the Ising Interface model. The first
display is the exact image law: the weight $Z_{\beta,h}$ of a
tiling $h$ in the prescribed sector $\cM_\Gamma$ is the partition
function of its \HMSOS fiber, namely of all periodic \HMSOS surfaces
sorting to $h$, with the energy measured relative to the reference
surface; this is the content of \cref{prop:periodic-fiber}. The
second display is the \emph{canonical extension}: the Grassmann
polynomial of \cref{lem:pattern-algebra} sums over \emph{all}
perfect matchings of the honeycomb torus, including those whose
monotone lift has the wrong periods, and the \RG analysis of
\cref{sec:rg} integrates this polynomial without regard to sectors.
We therefore need a single family of pattern coefficients, defined
on every connected pattern of the torus, whose hard-core sum is
positive on every matching and reproduces $Z_{\beta,h}$ on the
physical sector. The same coefficients must moreover agree, for
patterns much smaller than the torus, with a fixed family of planar
coefficients independent of $n$: this is what allows the
finite-volume laws to be compared with one bulk limit in
\cref{thm:coefficient-rg}, and it is the consistency between volumes
referred to in \cref{thm:rg-typical}. The third layer consists of
the numerical inputs, the \texttt{thin}/\texttt{fat} splitting of the
coefficients and the rooted counts of their supports, which are
inherited from the mixed expansion of \cref{sec:hmsos-rg-input} and
its periodic version. The quantities $\tau_\beta,B_\beta,R_\beta/2$
and $\epsilon_\beta$ are the \texttt{thin} decay rate, the minimal
size of a \texttt{fat} pattern, the radius of the balls in the \texttt{fat} property and the common
small rate for the \texttt{fat} coefficient growth and entropy.

\begin{proposition}[Exact bulk \HMSOS fermion expansion]\label{prop:hmsos-fermions}
Normalize the liquid dimer weights so that $\max(a,b,c)=1$.
Let $\mu_0$ be either a planar finite-volume free dimer law with a
fixed boundary lift, or the infinite-volume liquid law. Restrict
defects to those whose enlarged patterns are contained in a fixed
finite geometric window $\Omega$; in the planar case this window
lies inside the graph. The restriction is determined by the enlarged
pattern and includes its entire reference footprint.
Let $G$ be the corresponding Kasteleyn inverse and
$\cE_G$ its normalized Grassmann expectation. For a finite edge-source set $F$, use $J_A$ and $\Psi_{P,A}$
from \cref{lem:pattern-algebra}.
With $c_\beta$ as in \cref{lem:mixed-locality}, define
\[
 V_{\beta,A}(\psi)=\sum_Pc_\beta(P)\Psi_{P,A}(\psi)\,.
\]
Then exactly
\begin{equation}\label{eq:direct-fermion-identity}
 \mu_0\bigg(Z_{\beta,h}
                e^{\sum_eA(e)\one_{\{e\in h\}}}\bigg)
       =\cE_G\big(e^{J_A+V_{\beta,A}}\big)\,.
\end{equation}
At $A=0$ the labelled monomials in $V_{\beta,0}$ satisfy the first
three bounds of \cref{prop:fermion-input}, after converting
lozenge counts to vertex counts. They are balanced and connected.
Their coefficients are entire in $A$. On $\|A\|_\infty\leq\rho$,
the bounds in \cref{prop:mixed-envelopes} acquire only the
factor
\begin{equation}\label{eq:direct-source-factor}
              e^{\rho|P\cap F|}\leq
                          e^{\rho\min\{|P|,|F|\}}\,.
\end{equation}
Every interaction monomial obeys the vertex gauge covariance in
\cref{lem:gauge-identity}.
\end{proposition}
\begin{proof}
Recall from \cref{lem:mixed-locality} that the window-restricted
fiber partition function above a tiling $h$ is the exact hard-core
gas \cref{eq:mixed-gas}: a sum over families of enlarged patterns
$P_i\subset h$ with pairwise disjoint supports, each weighted by the
coefficient $c_\beta(P_i)$ of \cref{eq:mixed-coefficient}, which
already sums all defect decompositions of $P_i$. This is precisely
the polynomial $F(h)$ of \cref{lem:pattern-algebra} with
$c=c_\beta$, and the restriction to the window $\Omega$ amounts to
setting $c_\beta(P)=0$ for patterns not contained in $\Omega$.
Multiplying by the source factor $e^{\sum_eA(e)\one_{\{e\in h\}}}$
and taking the $\mu_0$-expectation, that lemma gives
\cref{eq:direct-fermion-identity}, together with the exact source
Taylor factor and the vertex gauge covariance. The coefficients
$c_\beta(P)$ are sums over the defects; these
sums are finite in a finite window, and the removal of the
fiber truncation in the lemma, by letting $J\to\infty$, is justified by the
a priori bound of \cref{prop:positive-budget}, which bounds the
sum of all weights absolutely.

For the coefficient estimates, retain the \texttt{thin}/\texttt{fat}
labels of \cref{prop:mixed-envelopes}: each coefficient is split as
$c_\beta=c^\mathrm{t}_\beta+c^\mathrm{f}_\beta$ and the two parts are
carried as separate labelled monomials in $V_{\beta,0}$. The
Grassmann monomial of $P$ has $|S(P)|=2|P|$ fields and, since the
liquid weights are normalized by $\max(a,b,c)=1$, its Kasteleyn
prefactor satisfies $|\prod_{e\in P}K(e)|\leq1$. Hence the
\texttt{thin} decay \cref{eq:mixed-thin-bound} in $|P|$ becomes a
decay $\epsilon_\beta^{|S|}$ in the number of fields, with
$\epsilon_\beta=e^{-c\beta^{2/5}/2}$, and the \texttt{fat} bound
\cref{eq:mixed-fat-bound}, the minimum \texttt{fat} size $B_\beta$,
the \texttt{fat} parameters $(\epsilon_\beta,R_\beta/2)$ and the rooted entropy
$e^{\epsilon_\beta|P|}$ translate into the second and third
properties of \cref{prop:fermion-input} after changing the absolute
constants. The monomials are balanced, with one black and one white
field per edge, and connected, since $S(P)$ is connected.

It remains to discuss the sources. For a fixed finite source set
$F$, the factor \cref{eq:direct-source-factor} is at most
$e^{\rho|F|}$. Since every \texttt{fat} pattern has
$|P|\geq B_\beta$, this is at most $e^{(\rho|F|/B_\beta)|P|}$, so it
can be absorbed into the \texttt{fat} growth rate by replacing
$\epsilon_\beta$ with
$\epsilon_{\beta,F}=\epsilon_\beta+\rho|F|/B_\beta$, which still
tends to zero as $\beta\to\infty$ because $B_\beta\to\infty$. For
sources at arbitrary locations, however, no such absorption is
possible, and we retain the full Taylor factor $e^{\rho|P|}$ of
\cref{lem:pattern-algebra}, as required by \cref{it:rg-sources}.
\end{proof}

The source Taylor factor will be paid once, at the first integration
of \cref{sec:first-integration}, using a reserved negative exponent
in the \texttt{fat}-node bound.

The infinite-volume formulation in the last proposition is a
localized perturbation of the free liquid dimer law: only defects
inside the window $\Omega$ are allowed, and outside it the tiling is
free. Its coefficients $c_\beta(P)$ do not depend on the window
containing $P$, by \cref{lem:mixed-locality}, so letting the window
increase supplies the translation-invariant bulk coefficients
required by \cref{it:rg-locality}. The following construction
extends those coefficients to the torus. We emphasize that our
analytic estimates only ever use the bulk covariance, i.e.\ the
inverse Kasteleyn matrix of the full-plane liquid law, or its twisted
periodizations on the torus; no uniform claim about the inverse
Kasteleyn matrices of arbitrary planar regions with boundary is
needed.

\subsubsection{A canonical extension and exact sector extraction}
\label{sec:canonical-extension}

Keep the physical periods $\Gamma$ and the two tori of
\cref{sec:periodic-hmsos}. Call a torus tiling
\emph{$\Gamma$-compatible} if its monotone lift has these physical
periods. Such tilings exist: the boundary of the down-set
$\{x\in\Z^3:\Theta_1x_1+\Theta_2x_2+nx_3<t\}$ is one example, for any
nonintegral $t$, since the linear form $\Theta_1x_1+\Theta_2x_2+nx_3$ vanishes on
both generators $(-n,0,\Theta_1)$ and $(0,-n,\Theta_2)$ of $\Gamma$. Denote this
sector by $\cM_\Gamma$. Every other tiling of the diagonal torus
lifts to a monotone surface whose periods differ from those of
$\Gamma$ by integer multiples of $(1,1,1)$; these are the auxiliary
sectors on which the canonical extension must also be defined.

The coefficient of an overlap component in its full periodic fiber
will initially be denoted by $c^\mathrm{per}_{\beta,h}(P)$: here
$h\in\cM_\Gamma$ is a compatible tiling containing $P$, and the
coefficient is computed by the periodic version of
\cref{eq:mixed-coefficient}, summing over compatible families of
defects on the transverse torus whose enlarged patterns overlap to
form $P$. The first task is to show that this coefficient does not
depend on $h$, so that it is a function of the pattern alone, and
that it agrees with the planar coefficient when $P$ is small
compared with the torus. The dependence on $h$ is not obviously
absent: a defect on the transverse torus is read off from the lift
of $P$ to $\Z^3/\Gamma$, and the lift of a pattern on the diagonal
torus is a priori determined only up to the kernel of the diagonal
projection.

\begin{lemma}[Locality with prescribed monodromy]\label{lem:periodic-locality}
If a connected torus pattern $P$ occurs in two tilings
$h,h'\in\cM_\Gamma$, then
\[
 c^\mathrm{per}_{\beta,h}(P)=c^\mathrm{per}_{\beta,h'}(P)\,.
\]
The equality also holds separately for the \texttt{thin} and \texttt{fat} parts of
the mixed coefficient. There is a constant $\eta_*>0$, depending
only on the compact set of physical slopes, such that, whenever
$\diam_{\T_\mathrm{d}}P\leq\eta_*n$, this common coefficient
equals $c_\beta(\widetilde P)$ for the planar lift $\widetilde P$.
The same is true of the two parts.
\end{lemma}
\begin{proof}
We first show that the lift of $P$ is intrinsic. Lift one plaquette
of $P$ to $\Z^3/\Gamma$; there are countably many choices, differing
by multiples of $(1,1,1)$, the kernel of the diagonal projection
$\Upsilon$. Two adjacent lozenges of $P$ share an edge, and their
orientations determine the relative position in $\Z^3$ of the two
unit plaquettes they represent; hence, by connectedness, the lift of
one plaquette determines the lifts of all the others. Along a closed
path of adjacent lozenges in $P$ on the torus, the lifted path
closes up to a physical displacement, which is an element of
$\Gamma$ whose diagonal projection is the prescribed displacement of
the path on $\T_\mathrm{d}$. This element is unique because
$\Upsilon$ is injective on $\Gamma$: a relation
$k(-n,0,\Theta_1)+l(0,-n,\Theta_2)=m(1,1,1)$ forces $k=l=-m/n$ and then
$k(\Theta_1+\Theta_2)=-kn$, whence $k=l=m=0$ as $\Theta_1,\Theta_2,n>0$. Thus the two
surfaces $h,h'$ induce the same lift of $P$, up to a common
translation by an integer multiple of $(1,1,1)$. Such a translation
only permutes transverse columns and adds a common constant to their
charges; the energy and the monotonicity constraints are functions
of charge differences between neighboring columns, so no
defect count changes.

In this lift, the dots of $P$ have distinct transverse coordinates
modulo $\Gamma_\mathrm{s}$. Enumerate candidate supports $D$ among
these coordinates. Require every dot of $D^\ddagger$ and every complete
wall between adjacent columns of $D^\ddagger$ to belong to $P$, using
the prescribed period offsets on crossing edges. These tests
recover precisely the footprints contained in $P$.
The charges on the footprint determine the isolated activity,
including the neighbor monotonicities and all identifications
for a noncontractible support. For \texttt{bad} supports the radius-$R$ ball
cores, their tiling closures on the torus, and the removed pieces
are determined in the same way as in \cref{lem:mixed-locality}.
In particular, a tiling closure is contained in $P$ exactly when
its geometric cores are contained in the region covered by $P$.
The overlap, union and \texttt{thin}/\texttt{fat} tests now use only these recovered
data. This proves the first assertion, without cutting a noncontractible
component into independently lifted pieces.

For the last assertion choose $\eta_*$ small enough that a ball
of radius $4\eta_*n$ embeds in the diagonal torus. A small $P$
then has a unique connected planar lift up to a period translation.
The lift of any full-plane lozenge tiling to its stepped surface
is uniformly Lipschitz in diagonal coordinates: on every triangle
its affine lift has one of three fixed gradients, so integrating
along a straight segment gives a universal Lipschitz constant.
Consequently all lifted plaquettes of $P$, and their transverse
coordinates, have diameter at most $C\eta_*n$.
Decrease $\eta_*$ so that this is less than one tenth of the
shortest transverse period. Every candidate support whose
footprint is in $P$ then has an embedded planar lift; its neighborhood and
all the local constraints are the planar ones. In the \texttt{bad} case,
the ball cores of a witness contained in $P$ have the same lifts,
so the periodic and planar closures coincide as well. Conversely,
every planar candidate with enlarged pattern in $\widetilde P$
projects to such a periodic candidate. This is a weight-preserving
bijection of the families in the coefficient formula. It also
preserves $q(\cD)$ and hence the \texttt{thin}/\texttt{fat} split.
\end{proof}

Here and below $n$ is sufficiently large compared with the fixed
witness parameters at inverse temperature $\beta$. We call $P$
\emph{small} when it satisfies the diameter condition in the lemma.
For every connected partial matching $P$ of the torus define
\begin{equation}\label{eq:canonical-coefficients}
 c_{\beta,\Gamma}(P)=
 \begin{cases}
 c_\beta(\widetilde P),&P\text{ small and }\widetilde P
                      \text{ extends to a planar tiling},\\
 c^\mathrm{per}_{\beta,h}(P),&P\text{ not small and }
                      P\subset h\in\cM_\Gamma,\\
 0,&\text{otherwise}.
 \end{cases}\,.
\end{equation}
Use the same rule for $c^\mathrm{t}$ and $c^\mathrm{f}$. Extendability
in this definition is a property of the pattern, not a choice of
a tiling containing it, and the second case is well defined by
\cref{lem:periodic-locality}. The three cases should be read as
follows. A pattern occurring in some matching of the torus is either
small, in which case it receives the planar coefficient of its lift,
or large, in which case it receives the periodic coefficient when it
occurs in the physical sector. A large pattern occurring in no
physical-sector matching receives coefficient zero, as does a small
pattern whose lift has no planar extension. For small patterns there is no
requirement of extendability to the prescribed torus sector: a
small pattern that extends to some planar tiling receives its
planar coefficient even if it occurs in no tiling of $\cM_\Gamma$.
This convention is useful: it gives exact agreement with every bulk
coefficient without requiring a separate theorem on gluing a patch
into that sector, and it makes the small-pattern part of the torus
potential literally the periodization of the bulk potential, which
is the form in which \cref{it:rg-locality} is used. Coefficients of
patterns contained in no matching never enter a matching
polynomial, so this choice affects no partition function; it only
fixes the monomials of the Grassmann potential, and these are
controlled by the same bounds as their planar counterparts.

\begin{proposition}[Positive canonical completion]\label{prop:canonical-completion}
The coefficients \cref{eq:canonical-coefficients} are well-defined,
nonnegative and translation covariant. Their \texttt{thin} and \texttt{fat} bounds
are \cref{eq:mixed-thin-bound}--\cref{eq:mixed-fat-bound}, with
the same minimum \texttt{fat} size and \texttt{fat} parameters $(\epsilon_\beta,R_\beta/2)$. The rooted \texttt{fat}
count is at most $\exp[(\epsilon_\beta+C/n)s]$ at size $s$,
after enlarging $\epsilon_\beta\to0$.
For every torus tiling $M$ set
\begin{equation}\label{eq:canonical-completion}
 F_{\beta,\Gamma}(M)=
 \sum_{\substack{\{P_i\}:\ S(P_i)\text{ disjoint}\\P_i\subset M}}
                           \prod_i c_{\beta,\Gamma}(P_i)\,.
\end{equation}
This is a finite positive function, and
$F_{\beta,\Gamma}(h)=Z_{\beta,h}$ for all
$h\in\cM_\Gamma$. Coefficients differing from the planar
rule occur only on patterns with at least $c_*n$ lozenges.
\end{proposition}
\begin{proof}
We verify in turn that the coefficients are well defined, their
bounds and counts, the identity on the physical sector, positivity,
and the size of the exceptional patterns.

\Cref{lem:periodic-locality} shows that the second case of
\cref{eq:canonical-coefficients} does not depend on the tiling
$h\in\cM_\Gamma$ containing $P$. The two cases are also consistent
where both could apply: if a small pattern $P$ is contained in a
compatible tiling $h$, then its planar lift $\widetilde P$ is
contained in the lift of $h$, hence extends to a planar tiling, and
the lemma gives $c^\mathrm{per}_{\beta,h}(P)=c_\beta(\widetilde P)$.
Thus $c_{\beta,\Gamma}(P)$ equals the periodic coefficient on every
pattern occurring in a compatible tiling. Nonnegativity holds because
every periodic and planar coefficient is a sum of products of
positive activities. Translation of a pattern preserves smallness,
planar extendability, sector extendability and all the tests
entering the coefficient formula, proving covariance.

The coefficient bounds \cref{eq:mixed-thin-bound,eq:mixed-fat-bound}
follow from the planar and periodic versions of
\cref{prop:mixed-envelopes}, respectively; these have the same
minimum \texttt{fat} size and \texttt{fat} parameters. For the rooted
count of \texttt{fat} patterns, the
planar family has at most $e^{\epsilon_\beta s}$ rooted members of
size $s$ and the periodic family at most $e^{(\epsilon_\beta+C/n)s}$,
so their union has at most $2e^{(\epsilon_\beta+C/n)s}$; since
$s\geq B_\beta$, the factor two is absorbed by enlarging
$\epsilon_\beta$ by $\log2/B_\beta\to0$. Small \texttt{fat}
patterns and the balls in their \texttt{fat} decomposition embed in the
torus, since a ball of radius $4\eta_*n$ does, so a planar \texttt{fat} pattern is
also \texttt{fat} on the torus.

For the identity on $\cM_\Gamma$, fix $h\in\cM_\Gamma$. Its
complete periodic fiber is, by \cref{prop:periodic-fiber}, the
hard-core gas over compatible families of defect supports on the
transverse torus. Group the defects of a configuration according to
the overlap components of their enlarged patterns, exactly as in the
proof of \cref{eq:mixed-gas}: the unions $P_i$ of the components have
pairwise disjoint triangle supports, each $P_i$ is a connected
pattern contained in $h$, and the sum of the weights of all
configurations with prescribed components is the product of the
periodic coefficients $c^\mathrm{per}_{\beta,h}(P_i)$, which equal
$c_{\beta,\Gamma}(P_i)$ by the first paragraph. Conversely every
family of patterns in $h$ with disjoint supports and nonzero
coefficients arises this way. Summing over the families proves
$F_{\beta,\Gamma}(h)=Z_{\beta,h}$.

Positivity and finiteness on every sector are elementary: the empty
collection contributes one and all other terms are nonnegative, so
$F_{\beta,\Gamma}(M)\geq1$; there are finitely many patterns on the
torus, and each coefficient is finite by the a priori bound of
\cref{prop:periodic-fiber}. No assertion about an \HMSOS
interpretation of the other sectors is needed:
\cref{eq:canonical-completion} simply defines the auxiliary canonical
model there.

Finally, a connected union of $s$ unit tiles has diameter at most
$Cs$, so a pattern of diameter larger than $\eta_*n$ has at least
$c_*n$ lozenges with $c_*=\eta_*/C$. Since every pattern of diameter
at most $\eta_*n$ is small and carries its planar coefficient, all
exceptions to the planar rule have size at least $c_*n$.
\end{proof}

We now give the exact fixed-sector algebra, including singular
twists. The issue is the following. The Berezin integral of
\cref{lem:pattern-algebra} on the torus is a signed sum over all
perfect matchings, weighted by $F_{\beta,\Gamma}$, whereas the
\HMSOS law lives on the single sector $\cM_\Gamma$. To recover the
latter we exploit that the number $N_r(M)$ of lozenges of each
orientation $r$ is constant on a sector, and we introduce formal
variables $y_r$ marking these numbers, so that the sector becomes a
Fourier coefficient of the Berezin polynomial. Two facts are then
needed: that the orientation counts determine the sector, and that
the Kasteleyn sign is constant on the sector, so that the signed sum
restricted to it is a positive multiple of the physical partition
function. Write $N_0=n(n+\Theta_1+\Theta_2)$ for the total number of lozenges of
a torus tiling, let $\omega_1,\omega_2,\omega_3>0$ be arbitrary
dimer weights, and put $y_1=1$. In fixed black and white orderings
let $K_y$ be the honeycomb matrix with entry $\omega_r y_r$ on an
edge of orientation $r$. Define the unnormalized Berezin polynomial
$\fZ_{\beta,\Gamma}(y,A)$ by \cref{eq:rg-berezin}, with the
coefficients $c_{\beta,\Gamma}$ and with the orientation of
$\cD\psi$ chosen so that the quadratic integral is $\det K_y$. Sources may be present on any
subset of the edges.

\begin{proposition}[Exact projection to the \HMSOS sector]\label{prop:exact-sector}
Let $\nu_{\beta,\Gamma}$ be the law on $\cM_\Gamma$
with weights $Z_{\beta,h}$. Then
\begin{equation}\label{eq:exact-sector-ratio}
 \nu_{\beta,\Gamma}\!\bigg(e^{\sum_eA(e)\one_{\{e\in h\}}}\bigg)
 =\frac{[y_2^{n\Theta_2}y_3^{n^2}]\fZ_{\beta,\Gamma}(y,A)}
        {[y_2^{n\Theta_2}y_3^{n^2}]\fZ_{\beta,\Gamma}(y,0)}\,.
\end{equation}
The denominator is nonzero. Each coefficient extraction may
equivalently be performed by the integral over $|y_2|=|y_3|=1$,
or by the discrete Fourier sum on the $(N_0+1)^2$ roots-of-unity
grid. These identities hold even at twists where $K_y$ is singular.
The interaction coefficients are independent of $y$, and on the
unit torus their absolute bounds with sources acquire only the
factor $e^{\rho|P\cap F|}$ of \cref{eq:direct-source-factor}.
The vertex gauge identity holds before and after projection.
\end{proposition}
\begin{proof}
\Cref{lem:pattern-algebra} gives
\begin{equation}\label{eq:signed-canonical-sum}
 \fZ_{\beta,\Gamma}(y,A)=
 \sum_M\sgn(M)F_{\beta,\Gamma}(M)
        \prod_{r=1}^3\omega_r^{N_r(M)}y_r^{N_r(M)}
                           e^{\sum_{e\in M}A(e)}\,.
\end{equation}
Here $\sgn(M)$ is the permutation sign in the determinant, and the
monomial $\prod_ry_r^{N_r(M)}$ records the orientation counts. We
verify the two facts needed to select the sector: that the counts
$(N_1,N_2,N_3)=(n\Theta_1,n\Theta_2,n^2)$ are attained exactly on $\cM_\Gamma$,
and that $\sgn(M)$ is constant on $\cM_\Gamma$.

\emph{The orientation counts determine the physical monodromy.}
Let $\gamma_1=(-n,0,\Theta_1)$ and $\gamma_2=(0,-n,\Theta_2)$ be the prescribed
physical periods. The lift of any torus tiling to a monotone
surface in $\Z^3$ has the same diagonal periods $\Upsilon\gamma_i$,
and hence has periods $\gamma_i+k_i(1,1,1)$ for some integers
$k_1,k_2$; the sector $\cM_\Gamma$ is the case $k_1=k_2=0$. Cut the
lifted surface open along lifts of two cycles generating the
homology of the torus, to obtain a stepped fundamental surface whose
boundary is, up to the cutting curves, a parallelogram with sides
$\gamma_1+k_1(1,1,1)$ and $\gamma_2+k_2(1,1,1)$. The oriented area
vector of a surface with given boundary depends only on the
boundary, by Stokes' theorem applied to $\frac12\oint u\times\d u$,
and for this parallelogram it is the cross product of the two sides.
On the other hand the oriented area vector is the sum of the
oriented areas of the unit plaquettes, and a plaquette of a monotone
surface normal to the axis $r$ contributes $+\be_r$. Hence the area
vector is $(N_1,N_2,N_3)$, and in the physical sector it is
$\gamma_1\times\gamma_2=(n\Theta_1,n\Theta_2,n^2)$. For a tiling with periods
$\gamma_i+k_i(1,1,1)$, expanding the cross product gives
\[
 (\gamma_1+k_1(1,1,1))\times(\gamma_2+k_2(1,1,1))
 =\gamma_1\times\gamma_2
      +(1,1,1)\mathbin\times(k_1\gamma_2-k_2\gamma_1)\,.
\]
Thus equality of the counts with $(n\Theta_1,n\Theta_2,n^2)$ implies
\[
 (1,1,1)\mathbin\times(k_1\gamma_2-k_2\gamma_1)=0\,,
\]
i.e.\ $k_1\gamma_2-k_2\gamma_1$ is a multiple of $(1,1,1)$.
Projecting onto the diagonal plane kills $(1,1,1)$ and gives
$k_1\Upsilon\gamma_2=k_2\Upsilon\gamma_1$; since the two diagonal
periods are linearly independent, $k_1=k_2=0$. Thus the target
counts select exactly $\cM_\Gamma$. Two counts suffice in
\cref{eq:exact-sector-ratio} because $N_1+N_2+N_3=N_0$ for every
tiling, so prescribing $N_2$ and $N_3$ prescribes $N_1$ as well.

\emph{The determinant sign is constant at fixed counts.}
The toroidal Kasteleyn sign rule says that the relative sign of
two matchings depends only on the mod-$2$ homology class of their
symmetric difference; see \cite[\S2.1]{KSW16}. The all-positive
honeycomb signing is Kasteleyn, since each face has length six.
To apply the rule, orient the alternating cycles of $M\xor M'$
from black to white on the $M$-edges. If $d_r$ is the planar
black-to-white edge vector of orientation $r$, their total
displacement is
\[
                  \sum_r(N_r(M)-N_r(M'))d_r=0
\]
when the counts agree. The period vectors are linearly independent,
so the total integral homology class, and hence its mod-$2$
reduction, is zero. The sign rule therefore gives
$\sgn(M)=\sgn(M')$. Denote the common sign on $\cM_\Gamma$ by
$\sigma_\Gamma$.

\emph{Conclusion.} By the two facts just proved, the coefficient of
$y_2^{n\Theta_2}y_3^{n^2}$ in \cref{eq:signed-canonical-sum} is
\[
 \sigma_\Gamma\omega_1^{n\Theta_1}\omega_2^{n\Theta_2}\omega_3^{n^2}
       \sum_{h\in\cM_\Gamma}Z_{\beta,h}
                                e^{\sum_{e\in h}A(e)}\,,
\]
using $F_{\beta,\Gamma}(h)=Z_{\beta,h}$ on $\cM_\Gamma$ from
\cref{prop:canonical-completion}. At zero source the sum is a sum
of positive terms over the nonempty set $\cM_\Gamma$, so the
denominator in \cref{eq:exact-sector-ratio} is nonzero, and the
prefactor $\sigma_\Gamma\prod_r\omega_r^{N_r}$ cancels in the
ratio, which is therefore independent of the auxiliary positive
weights and equals the left-hand side by the definition of
$\nu_{\beta,\Gamma}$. Each variable $y_r$ has degree at most $N_0$
in the polynomial $\fZ_{\beta,\Gamma}(y,A)$, so the orthogonality
of characters on the roots of unity of order $N_0+1$ extracts the
coefficient with no aliasing; integration over the unit circles is
the same orthogonality identity. No inverse matrix or quotient of
determinants is used anywhere, which is why the identities hold at
twists where $K_y$ is singular. Finally, on the unit torus the
phases $y_r$ have modulus one, so $|K_y(e)|\leq\omega_r\leq1$ and
the absolute bounds on the interaction monomials with sources are
those of \cref{prop:hmsos-fermions}, with the factor
\cref{eq:direct-source-factor}; and since every matching uses
exactly one edge at each vertex, the vertex gauge identity of
\cref{lem:pattern-algebra} holds termwise in $M$, hence before and
after the coefficient extraction.
\end{proof}

This establishes the algebraic fixed-sector representation of
the fiber-weighted tiling law. We stress what it does and does not
give. It expresses the \HMSOS generating function as a ratio of two
Fourier coefficients of a Grassmann polynomial, each of which is a
signed sum over all sectors. Estimating such a coefficient after the
remaining integrations requires showing that the prescribed sector
carries a non-negligible fraction of the total weight, uniformly in
the twists; this is the analytic problem solved in
\cref{sec:fixed-sector-limit}, and the identity alone gives no lower
bound relative to the largest twisted partition function. The
relation to the physical height ensemble and its root is established
in \cref{sec:hmsos-rooting}.

\subsubsection{Verification of the marked input}
\label{sec:hmsos-rg-verification}

For witness estimates let $\scW_n$ consist of all connected
patterns of at least $B_\beta$ tiles admitting the six-class
representation of \cref{prop:periodic-witness}, with
fixed enlarged constants and duplicates removed. It contains the
\texttt{fat} coefficient patterns and the geometric witnesses used in
reconstruction. Recall the role of \cref{it:rg-marked}: in
\cref{sec:interacting-witnesses} the probability that a witness
$W\in\scW_n$ through $x$ is present in the interacting tiling is
bounded by inserting the indicator of $W$ into the pattern
expansion. The inserted witness may overlap ordinary expansion
patterns, and the insertion is made exact by merging $W$ with the
patterns it overlaps into a single connected pattern $U$ with the
merged coefficient $b_{x,s}(U)$ of
\cref{eq:marked-union-coefficient}; the condition asks for
\texttt{thin}/\texttt{fat} bounds and rooted counts for these
merged coefficients of the same type as for the original ones. The
next lemma verifies this absolute collection condition; its proof
is entirely geometric and combinatorial and uses no interacting
probability estimate.

\begin{lemma}[Absolute merged-collection bounds for six-class witnesses]
\label{lem:marked-union-envelopes}
For the \HMSOS coefficients and witnesses, the absolute merged
coefficients $b_{x,s}:=\overline b_{x,s}$ of
\cref{eq:rg-absolute-insertion} have a nonnegative split
$b_{x,s}=b^\mathrm{t}_{x,s}+b^\mathrm{f}_{x,s}$ such that
\begin{align}
 b^\mathrm{t}_{x,s}(U)&\leq e^{-c\beta^{3/10}|U|},
                                      \label{eq:marked-thin-envelope}\\
 b^\mathrm{f}_{x,s}(U)&\leq e^{\epsilon'_\beta|U|}\,,
                                      \notag
\end{align}
where $\epsilon'_\beta\to0$. The \texttt{thin} supports have at most
$C_0^m$ choices of size $m$ through $x$. The \texttt{fat} supports are
$(\epsilon'_\beta,R_\beta/2)$-\texttt{fat} and have at most
$e^{\epsilon'_\beta m}$ choices of size $m$ through $x$.
Both parts vanish unless $x\in S(U)$ and
$|U|\geq\max(s,B)$, as required in \cref{it:rg-marked}. These estimates hold
uniformly for noncontractible supports as well. More generally the same
conclusion holds for nonnegative envelopes with these coefficient
bounds and the same six-class geometry (with fixed changes of its
constants). Thus the statement also applies to absolute envelopes
of signed coefficients.
\end{lemma}
\begin{proof}
The proof has two parts. The numerical part bounds the total
weight of all merged decompositions of a given support $U$, with an
extra exponential enhancement on the \texttt{thin} patterns which is
then used to split the merged coefficient. The geometric part counts
the supports $U$ which can arise with a small \texttt{thin} content,
by showing that they inherit the six-class structure of their
constituents.

In this proof write $c^\mathrm{t},c^\mathrm{f}$ for the nonnegative
envelopes $q^\mathrm{t},q^\mathrm{f}$; for \HMSOS these are the actual
positive coefficient parts. Expand each
$c(P_i)=c^\mathrm{t}(P_i)+c^\mathrm{f}(P_i)$ in
\cref{eq:rg-absolute-insertion}, retaining these two labels, so
that a summand is a witness $W$ together with a family of labelled
patterns. Write $A_\beta=c\beta^{2/5}$ for a common \texttt{thin}
exponent in \cref{eq:mixed-thin-bound}. Enlarge $\epsilon_\beta\to0$
so that it bounds the \texttt{fat} coefficient exponent and the
rooted counts of both \texttt{fat} patterns and witnesses. For
large $n$ it also absorbs the torus entropy term $C/n$.
Set
\[
 \delta_\beta^\mathrm{m}=\beta^{-1/10},\qquad
 \zeta_\beta=4\epsilon_\beta+B^{-1/2}\,.
\]
The parameter $\zeta_\beta$ is a small damping inserted into every
activity to make the sum over patterns through a fixed triangle
convergent; it is chosen larger than twice the \texttt{fat} rate so
that the damping beats both the coefficient growth and the entropy,
yet it tends to zero. After reducing $A_\beta$ if necessary so that
$c^\mathrm{t}(P)\leq e^{-A_\beta|P|}$, the rooted counts $C_0^s$ for
\texttt{thin} and $e^{\epsilon_\beta s}$ for \texttt{fat} patterns
of size $s$ give, summing over the sizes,
\begin{align*}
 a_\beta
 &:={\sup_y}\sum_{P:y\in S(P)}
       \big(c^\mathrm{t}(P)e^{A_\beta|P|/2}+c^\mathrm{f}(P)\big)
                                          e^{-\zeta_\beta|P|}\\
 &\leq\sum_{s\geq1}\big(C_0e^{-A_\beta/2}\big)^s
   +\sum_{s\geq B}e^{(2\epsilon_\beta-\zeta_\beta)s}
 \leq Ce^{-A_\beta/3}
       +C\zeta_\beta^{-1}e^{-\zeta_\beta B/2}=o(1)\,,
\end{align*}
where the second sum used $2\epsilon_\beta-\zeta_\beta\leq-\zeta_\beta/2$
and $\sum_{s\geq B}e^{-\theta s}\leq C\theta^{-1}e^{-\theta B}$.
The constants here are independent of $n,x,s$. In particular the
small denominator $\zeta_\beta^{-1}\leq\sqrt B$ is paid by
$\zeta_\beta B\geq\sqrt B$, so the second term is at most
$C\sqrt Be^{-\sqrt B/2}\to0$.

Now fix a connected partial matching $U$ of size $m$ and consider a
summand of \cref{eq:rg-absolute-insertion}. The patterns $P_i$ have
disjoint supports inside $S(U)$, so $\sum_i|P_i|\leq m$; hence
inserting the damping $e^{-\zeta_\beta|P_i|}$ into every activity
costs at most a compensating factor $e^{\zeta_\beta m}$. Next drop
the requirements that the $P_i$ be disjoint, meet $S(W)$ and have
union $U$ with $W$; this only enlarges the sum, which becomes a sum
over arbitrary subsets of the patterns contained in $U$, i.e.\ the
product of $1+$activity over all patterns $P\subset U$. By the
elementary inequality $1+t\leq e^t$, the logarithm of this product
is at most the sum of the damped activities of all patterns in $U$,
and rooting each such pattern at one of its triangles, which lies in
$S(U)$, bounds this sum by $|S(U)|a_\beta=2ma_\beta$. The same
argument applies with the enhanced \texttt{thin} activities
$c^\mathrm{t}e^{A_\beta|P|/2}$, which were included in the
definition of $a_\beta$. We have therefore proved, for every $U$,
\begin{equation}\label{eq:disjoint-collection-budget}
 \sum_{\substack{\{P_i\}:\ P_i\subset U\\
                         S(P_i)\text{ pairwise disjoint}}}
  \prod_i\big(c^\mathrm{t}(P_i)e^{A_\beta|P_i|/2}
                       +c^\mathrm{f}(P_i)\big)
                 \leq e^{(\zeta_\beta+2a_\beta)|U|}\,.
\end{equation}
Finally the witness $W$ is a member of
$\scW_n$ with $x\in S(W)$, $W\subset U$ and $B\leq|W|\leq m$, so
there are at most $\sum_{B\leq s\leq m}e^{\epsilon_\beta s}
\leq me^{\epsilon_\beta m}$ choices for it. Since $m\geq B$ and
$\log m/m\leq\log B/B$ for large $B$, we obtain
\begin{equation}\label{eq:marked-collection-budget}
 \sum_W\sum_{\{P_i\}}
      e^{(A_\beta/2)\sum_{i:\,\rm t}|P_i|}
                        \prod_i c^\mathrm{label}(P_i)
 \leq e^{e_\beta m},\qquad
 e_\beta=\zeta_\beta+2a_\beta+\epsilon_\beta+\frac{\log B}{B}
                           \longrightarrow0\,.
\end{equation}
This argument uses only the coefficient bounds and rooted counts;
in particular it does not assume that the auxiliary canonical model
outside its prescribed sector is an \HMSOS fiber.

Assign a summand to $b^\mathrm{t}$ when
$\sum_{i:\,\rm t}|P_i|\geq\delta_\beta^\mathrm{m}m$, i.e.\ when a
fixed small fraction of the support is covered by \texttt{thin}
patterns, and to $b^\mathrm{f}$ otherwise. For a \texttt{thin}
summand the enhancement in \cref{eq:marked-collection-budget} is at
least $e^{A_\beta\delta_\beta^\mathrm{m}m/2}$, so extracting it gives
\[
 b^\mathrm{t}_{x,s}(U)\leq e^{(e_\beta-A_\beta\delta^\mathrm{m}_\beta/2)m}
                      \leq e^{-c\beta^{3/10}m}\,,
\]
since $A_\beta\delta_\beta^\mathrm{m}\asymp\beta^{2/5-1/10}=\beta^{3/10}$
while $e_\beta\to0$; this proves \cref{eq:marked-thin-envelope}.
Dropping the enhancement bounds $b^\mathrm{f}_{x,s}(U)$ by
$e^{e_\beta m}$, which is the asserted \texttt{fat} bound with
$\epsilon'_\beta\geq e_\beta$.

It remains to count the \texttt{fat} supports $U$, rather than
their decompositions. This is the only part using the six-class
geometry. Choose one contributing \texttt{fat} summand of $U$, i.e.\
one in which the \texttt{thin}-labelled patterns cover at most
$\delta_\beta^\mathrm{m}m$ tiles. Every \texttt{fat}-labelled $P_i$ has
the six-class representation constructed in
\cref{prop:mixed-envelopes}: up to a small remainder, it is a union
of thickened straight rows, grouped into six classes according to
the orientation of the row and its direction, each class being a
union of ball cores in which one lozenge orientation is
exceptional. Choose any such representation. The witness $W$ has
the same six classes by the definition of $\scW_n$. Unite the cores
separately in each class, over $W$ and all \texttt{fat}-labelled
$P_i$. The tiling closure of a union of cores is the union of their
closures, and the number of exceptional-orientation lozenges in a
class union is at most the sum of these numbers over the
constituents. Since the $P_i$ are disjoint subsets of $U$ and
$|W|\leq m$,
\[
                  |W|+\sum_{i:\,\rm f}|P_i|\leq2m\,.
\]
The class cores therefore have a total exceptional-orientation
count at most $CR_\beta^2m/L_\beta$, since each constituent of size
$s'$ contributes at most $CR_\beta^2s'/L_\beta$. The part of $U$ not
covered by the class unions consists of the removed pieces of the
constituents, at most $C(L_\beta/B_\beta)s'$ tiles for a constituent
of size $s'$, and of the \texttt{thin}-labelled patterns, at most
$\delta_\beta^\mathrm{m}m$ tiles in total; hence it has size at most
\[
 C\big(\beta^{-1/10}+L_\beta/B_\beta\big)m\,.
\]
Thus $U$ is determined by six class unions, a remainder of this size
and the exceptional orientations, which is the situation of
\cref{lem:periodic-unions}. Apply that lemma with
$l=\floor{\sqrt{R_\beta}}$. Its entropy rate
$C(l^{-2}+l/R_\beta+\delta l^2+\lambda^{-1}+\kappa\log(e/\kappa))$ tends
to zero: the remainder fraction $\delta$ is the last displayed
fraction, and $\delta l^2\leq C(\beta^{-1/10}+L_\beta/B_\beta)R_\beta\to0$
for the choices $B_\beta=\ceil{\sqrt\beta}$,
$R_\beta=\ceil{\beta^{1/100}}$, $L_\beta=\ceil{\beta^{1/10}}$ of
\cref{eq:history-scales}; the exceptional fraction $\kappa$ is at
most $CR_\beta^2/L_\beta\to0$; and $l^{-2}+l/R_\beta\to0$. The
ball-core boundary-layer estimate shows that the \texttt{fat} supports are
$(C/R_\beta+C\beta^{-1/10}+CL_\beta/B_\beta,R_\beta/2)$-\texttt{fat}, with
the first parameter tending to zero. All these estimates use only unions on the torus and
apply to noncontractible supports. For the \texttt{thin} part, the
unrestricted count of connected patterns of size $m$ through a fixed
triangle gives $C_0^m$. Enlarge $\epsilon'_\beta$ to absorb $e_\beta$
and the geometric count and the first \texttt{fat} parameter. Finally
the support restriction in \cref{it:rg-marked} holds because $U\supset W$ with
$x\in S(W)$ and $|W|\geq\max(s,B)$.
\end{proof}

The coefficient conditions are verified in this subsection. The proof of \cref{thm:hmsos-rg-input}
is completed after the physical identification in
\cref{prop:physical-hmsos-law}.

\subsection{Two steps: SOS input}

\begin{theorem}[\SOS input to the \RG theorem]
\label{thm:sos-rg-input}
For sufficiently large $\beta$ and all sufficiently large $n$
depending on $\beta$, uniformly for admissible periods
with slopes in a compact subset of $(0,\infty)^2$, the \Tetris
image $\varphi=\Tone\Ttwo\cS$ of the \SOS model, in the diagonal
translation quotient, has the fixed-sector law
\[
 \nu^{\textsc{sos}}_{\beta,\Gamma}(\varphi)
   =\frac{Z^{\textsc{sos}}_\beta(\varphi)}
          {\sum_{\psi\in\cM_\Gamma}Z^{\textsc{sos}}_\beta(\psi)},
 \qquad
 Z^{\textsc{sos}}_\beta(\varphi)
   =\sum_{\Tone\Ttwo\cS=\varphi}e^{-\beta(\sA(\cS)-N_0)},
 \quad N_0=n(n+\Theta_1+\Theta_2)\,.
\]
Here $\sA$ is the plaquette area. The sums are finite after
removing the height cutoff with a consistent
rooting. The full two-step weight admits a canonical extension and an
exact Grassmann expansion satisfying the four quantitative coefficient
and source requirements in \cref{prop:fermion-input}. The
bounds are uniform in volume and twists, with parameters tending to
zero or infinity as prescribed there. They include the periodic
\texttt{fat}-pattern entropy correction, agreement of small torus patterns
with bulk coefficients, and the fixed-count projection as in
\cref{thm:hmsos-rg-input}. Signed coefficients, if used, are bounded
after summing their internal decompositions absolutely.
In particular, \crefrange{it:rg-positive}{it:rg-sources} of
\cref{thm:coefficient-rg} hold. The six-class witness family
in \cref{lem:sos-marked-input} also satisfies \cref{it:rg-marked}.
\end{theorem}

The theorem is the two-step analogue of \cref{thm:hmsos-rg-input},
and its proof reduces to the one-step construction as follows.
Recall from \cref{def:two-step-history} that a two-step history consists of an
outer family of $x_1$-defects on $\varphi$, inverting the
last move $\Tone$ and producing the intermediate surface $\cH$, and
an inner family of $x_2$-defects on $\cH$, inverting $\Ttwo$ and
producing $\cS$, and that \cref{prop:sos-positive-expansion} carries
the footprints of the inner defects down to the final tiling and
groups overlapping footprints, giving an exact hard-core gas of
enlarged patterns on $\varphi$ whose coefficients
$\bar c^{\textsc{sos}}_\beta$ already sum over all intermediate
surfaces. From there the tasks are the same as for
\HMSOS: define coefficients on every torus pattern, compare small
torus coefficients with the bulk, and read off the Grassmann
representation from \cref{lem:pattern-algebra}. The only genuinely
new point is the comparison of small coefficients, treated first.

\subsubsection{Confinement, canonical completion and fermion sources}
\label{sec:sos-canonical}

The locality already proved in \cref{prop:sos-positive-expansion}
keeps the prescribed period identifications: the coefficient of a
pattern $P$ in a compatible torus tiling is determined by $P$.
To compare a small periodic coefficient with its bulk counterpart,
however, one must also exclude the following scenario, which has no
analogue in the one-step case: a small pattern $P$ on the final
tiling determines the outer supports through its dots, but the
inner defects live on the intermediate surface $\cH$, which is
reconstructed from $P$ and the outer states; one must show that this
reconstructed piece of $\cH$ is itself confined to a bounded
neighborhood of the lift of $P$, so that neither it nor the inner
supports read from it can wrap around the torus. The following
confinement observation supplies this missing step.

\begin{lemma}[Confinement and agreement of small \SOS coefficients]
\label{lem:sos-small-locality}
Let $P$ be a finite connected pattern in a planar reference tiling,
and insert any compatible family of outer states with
$P_\alpha\subset P$. Every changed cube and every new face of
the resulting $\cH$ lies in the coordinate bounding box of the
lift of $P$, enlarged by a fixed lattice distance. In particular,
the diameter of $P^{\cH}$ is at most $C(1+\diam P)$.

There is $\eta_*>0$, uniform over the specified compact slope
range, such that a periodic pattern with
$\diam_{\T_\mathrm{d}}P\leq\eta_*n$ has
\begin{equation}\label{eq:sos-small-locality}
 \bar c^{\textsc{sos},\mathrm{per}}_{\beta,\Gamma}(P)
                  =\bar c^{\textsc{sos}}_\beta(\widetilde P)\,,
\end{equation}
where $\widetilde P$ is its connected planar lift. The equality
holds separately for the \texttt{thin} and \texttt{fat} parts and is independent
of the compatible tiling containing $P$.
\end{lemma}
\begin{proof}
The first assertion is a sandwich argument: an outer defect changes
the spins of a column only between two unchanged columns on the
same sorting line, and monotonicity confines the changed cubes
between the reference heights of those two columns. Fix a changed
outer column $v=(x_2,x_3)\in D_\alpha$. On its $x_2$ coordinate
line, let $v_-,v_+$ be the first columns outside the consecutive run
of $D_\alpha$ containing $v$, in the two directions. These columns
belong to the neighborhood of $D_\alpha$ and to no other outer support:
otherwise two supports would be $*$-adjacent, contradicting
compatibility. They are therefore unchanged. Write $h_-,h_+$ for
their reference charges in direction $x_1$, so that the $\{0,1\}$-spin sequence
of these columns is $\one_{\{x_1<h_\pm\}}$ both before and after the
modification. Since $\cH$ is monotone in $x_2$, the $\{0,1\}$-spin sequence of
the column $v$, which lies between $v_-$ and $v_+$, is sandwiched
between those of $v_-$ and $v_+$, giving the pointwise inequalities
\[
 \one_{\{x_1<h_+\}}\leq I_{\cH}(x_1,v)
                              \leq\one_{\{x_1<h_-\}}\,.
\]
The same inequalities hold for the monotone reference surface.
Thus every changed cube on this column has $x_1$ between
$h_+$ and $h_-$. The endpoint dots, and the dot of $v$,
are in $P_\alpha\subset P$, so $h_\pm$ and the coordinates of $v$
are among the coordinates of lozenges of $P$. This bounds all three
coordinates of every changed cube by the coordinate ranges of the
lift of $P$. New faces are incident to changed cubes, proving the
first assertion, uniformly in the defects. The lift of a
full-plane tiling is uniformly Lipschitz in diagonal coordinates, as
in \cref{lem:periodic-locality}; hence the lift of $P$, and
so the bounding box, has diameter at most $C(1+\diam P)$.

Now take $P$ small on the torus. Choose $\eta_*$ so that
its connected lift embeds in a diagonal ball and its lifted
coordinate ranges are less than a sufficiently small fixed
multiple of both shortest transverse periods. Outer candidate
supports are read from the normal-$x_1$ dots of $P$. They,
their neighborhoods and the boundary columns just used therefore
have embedded planar lifts. The same sandwich proves confinement
of their intermediate faces on the cover. In particular all
normal-$x_2$ dots of $P^{\cH}$ lie in a box of diameter
$C\eta_*n+C$. Decrease $\eta_*$, and take $n$ large, so
that this box and its fixed neighborhood have diameter less than
one tenth of the shortest period of $\Gamma_2$.
Every inner support with footprint in $P^{\cH}$ then
has an embedded planar lift as well.

All inner reference charges, neighbor constraints and state
weights are consequently the planar ones. Arbitrarily distant
transitions on an individual sorting fiber do not introduce an
additional identification: projection of $\Gamma$ onto either
transverse lattice is injective, and distinct lifted support
copies have disjoint transverse neighborhoods. Equivalently, the
remaining vertical monotonicity sandwiches the inner states
between their unchanged boundary columns in direction $x_3$.
This bounds their moved faces by the coordinate ranges of
$P^{\cH}$ in the same way. Thus copying a planar history by the
periods, or lifting a periodic history, preserves every local
constraint and gives inverse operations on the allowed states.

The carrier tests and interval rows agree under this bijection.
Whenever a witness is contained in $P$, its ball cores lie in
the same embedded patch, so its torus and planar tiling closures
also agree. The enlarged overlap decomposition, exact union
and \texttt{thin}/\texttt{fat} tests are all preserved, as are $E_1,E_2$.
Summing this weight-preserving bijection proves
\cref{eq:sos-small-locality}. For larger patterns, use the
monodromy and state bijection in \cref{prop:sos-positive-expansion}.
The enlargement and \texttt{thin}/\texttt{fat} tests use only these states and
$P$, by \cref{prop:sos-mixed-envelopes}, so this bijection
also preserves the enlarged coefficients and their two parts.
\end{proof}

For an arbitrary connected partial matching $P$ on the
diagonal torus, define $c^{\textsc{sos}}_{\beta,\Gamma}(P)$
by the rule \cref{eq:canonical-coefficients}, using the
bulk and periodic enlarged \SOS coefficients
$\bar c^{\textsc{sos}}_\beta$ and
$\bar c^{\textsc{sos},\mathrm{per}}_{\beta,\Gamma}$ in place of
$c_\beta$ and $c^\mathrm{per}_{\beta,h}$, and the smallness
threshold of \cref{lem:sos-small-locality}. Apply the
same rule to their two parts. Thus small patterns use their
planar coefficients when planar-extendable, large patterns
use the periodic coefficients when extendable to
$\cM_\Gamma$, and all others are zero.
As in the one-step case, this is a definition on every torus
sector, whether or not the sector has an \SOS
interpretation, and the discussion following
\cref{eq:canonical-coefficients} applies verbatim.

\begin{proposition}[Canonical \SOS coefficients and exact fermion expansion]
\label{prop:sos-canonical-fermions}
The coefficients $c^{\textsc{sos}}_{\beta,\Gamma}(P)$ are
nonnegative, local and translation covariant. The completion
\begin{equation}\label{eq:sos-canonical-completion}
 F^{\textsc{sos}}_{\beta,\Gamma}(M)
   =\sum_{\substack{\{P_i\}:\ S(P_i)\text{ disjoint}\\P_i\subset M}}
                   \prod_i c^{\textsc{sos}}_{\beta,\Gamma}(P_i)\,,
\end{equation}
is finite, positive and equal to $Z^{\textsc{sos}}_\beta(M)$
on $\cM_\Gamma$. Its coefficients satisfy
\cref{eq:sos-mixed-envelopes} and the associated \texttt{fat} and
rooted entropy bounds, with torus entropy correction $C/n$.
Only patterns of size at least $c_*n$ can differ from the
bulk rule.

The potential
\[
 V^{\textsc{sos}}_{\beta,A,y}(\psi)
   =\sum_P c^{\textsc{sos}}_{\beta,\Gamma}(P)
                   \prod_{e=(b,w)\in P}K_y(e)e^{A(e)}\psi_b\psi_w\,,
\]
has the four coefficient and source properties in
\cref{prop:fermion-input}, in the source convention
specified there, uniformly in unit-modulus twists. With
\[
 \fZ^{\textsc{sos}}_{\beta,\Gamma}(y,A)
   =\int\!\cD\psi\,
       \exp\bigg(\sum_eK_y(e)e^{A(e)}\psi_b\psi_w
                         +V^{\textsc{sos}}_{\beta,A,y}(\psi)\bigg)\,,
\]
the exact fixed-sector generating function is
\begin{equation}\label{eq:sos-exact-sector}
 \nu^{\textsc{sos}}_{\beta,\Gamma}
     \!\bigg(e^{\sum_eA(e)\one_{\{e\in\varphi\}}}\bigg)
 =\frac{[y_2^{n\Theta_2}y_3^{n^2}]\fZ^{\textsc{sos}}_{\beta,\Gamma}(y,A)}
        {[y_2^{n\Theta_2}y_3^{n^2}]\fZ^{\textsc{sos}}_{\beta,\Gamma}(y,0)}\,.
\end{equation}
The denominator is nonzero; the polynomial formula includes
singular twists. The bulk analogue of
\cref{eq:direct-fermion-identity} holds with the enlarged
\SOS coefficients and the corresponding window-restricted
history partition function.
\end{proposition}
\begin{proof}
The proof follows that of
\cref{prop:canonical-completion,prop:exact-sector,prop:hmsos-fermions},
with the two-step ingredients substituted, and we indicate the
substitutions. \Cref{lem:sos-small-locality} makes the canonical
rule well defined and shows that it agrees with the full periodic
coefficient on every pattern in a compatible tiling: the periodic
coefficient is independent of the tiling containing the pattern, and a small
pattern in a compatible tiling has a planar-extendable lift whose
planar coefficient is the periodic one by
\cref{eq:sos-small-locality}. Nonnegativity holds because the
enlarged coefficients are sums of products of positive history
weights. Translation preserves all tests, including extendability,
giving covariance. For the identity on $\cM_\Gamma$, group the
histories of a tiling $M\in\cM_\Gamma$ by the overlap components of
their enlarged patterns, as in \cref{prop:sos-mixed-envelopes}: this
proves \cref{eq:sos-canonical-completion} with the right-hand side
equal to $Z^{\textsc{sos}}_\beta(M)$, and positivity and finiteness
on the other sectors follow from the empty collection and the
finiteness of the coefficients. The numerical bounds
\cref{eq:sos-mixed-envelopes}, the \texttt{fat} statement and the rooted
counts are those of \cref{prop:sos-mixed-envelopes}: for small
patterns use the planar entropy bound, for larger ones the periodic
bound with its correction $C/n$, and note that adding the two
families costs a factor two, absorbed into a vanishing entropy rate
since every \texttt{fat} pattern has size greater than
$2\sqrt\beta$. A connected pattern with $s$ unit tiles has diameter
at most $Cs$; hence a pattern outside the small class has
$s\geq c_*n$. These observations prove all completion and
comparison assertions.

\Cref{lem:pattern-algebra} applied to the coefficients
$c^{\textsc{sos}}_{\beta,\Gamma}$ gives the signed Berezin matching
sum
\[
 \fZ^{\textsc{sos}}_{\beta,\Gamma}(y,A)
 =\sum_M\sgn(M)F^{\textsc{sos}}_{\beta,\Gamma}(M)
   \prod_{r=1}^3\omega_r^{N_r(M)}y_r^{N_r(M)}e^{\sum_{e\in M}A(e)}\,,
\]
and its bulk counterpart, the analogue of
\cref{eq:direct-fermion-identity}, with the window-restricted
two-step history partition function in place of $Z_{\beta,h}$. The
orientation-count and sign argument of \cref{prop:exact-sector}
depends only on the honeycomb torus and on the positivity of the
weight on $\cM_\Gamma$, not on the origin of the weight, so it
applies unchanged to $F^{\textsc{sos}}_{\beta,\Gamma}$ and proves
\cref{eq:sos-exact-sector}, with a nonzero denominator even when
individual twisted matrices are singular.

For the Grassmann bounds, retain the \texttt{thin}/\texttt{fat}
labels of \cref{prop:sos-mixed-envelopes} on the monomials.
Normalize $\max_r\omega_r=1$ and use $|S(P)|=2|P|$ to convert the
lozenge-count bounds into the first three coefficient bounds of
\cref{prop:fermion-input}, exactly as in \cref{prop:hmsos-fermions}.
\Cref{lem:pattern-algebra} supplies its exact source and gauge
properties. For fixed source sets the \texttt{fat}-size argument in
\cref{prop:hmsos-fermions} absorbs the source factor into the
\texttt{fat} rate; the location-summed Taylor bound retains
$e^{\rho|P|}$. All estimates are on absolute internal weights, with
no unsummed history labels, since the enlarged coefficients already
sum over the intermediate surfaces.
\end{proof}

\begin{corollary}[Uniform pressure for the two-step fiber]
\label{cor:sos-fiber-pressure}
For every monotone periodic reference and sufficiently large
$t$,
\begin{equation}\label{eq:sos-fiber-pressure}
       1\leq Z^{\textsc{sos}}_t(\varphi)
                              \leq e^{Ce^{-t/64}N_0}\,.
\end{equation}
The complete two-step fiber is the total-variation limit
of its normalized fiber truncations at distance $J$ as $J\to\infty$.
More generally, for two sufficiently large inverse
temperatures,
\[
 \sum_{\substack{\Ttwo\cS=\cH\\\Tone\cH=\varphi}}
                   e^{-t_1E_1-t_2E_2}
       \leq e^{\epsilon_{t_2}N_0}
                         Z^{\textsc{hmsos}}_{t_1-\epsilon_{t_2}}(\varphi)\,,
\]
provided $t_1-\epsilon_{t_2}\geq t_0$.
\end{corollary}
\begin{proof}
Taking $W=\varphi$ in \cref{eq:sos-positive-collection-budget}
gives the upper bound in \cref{eq:sos-fiber-pressure}; the empty
history gives the lower bound. For two temperatures, sum over
the intermediate surface as in \cref{prop:sos-history-fiber}.
The inner bound \cref{eq:sos-first-fiber-pressure} and
$\sA(\cH)=N_0+E_1$ give directly
\[
 \sum_{\Tone\cH=\varphi}e^{-t_1E_1}
             \sum_{\Ttwo\cS=\cH}e^{-t_2E_2}
 \leq e^{\epsilon_{t_2}N_0}
       \sum_{\Tone\cH=\varphi}e^{-(t_1-\epsilon_{t_2})E_1},
\]
which is the stated inequality. Finally, let $Z_J$ be the partition
function with fiber truncation at distance $J$ at both inverse steps and let $Z$
be the unrestricted sum. The normalized restricted law is the full
law conditioned on both steps satisfying the truncation. Their total-variation distance
is $1-Z_J/Z\to0$, by finiteness of $Z$ and monotone convergence
of $Z_J$.
\end{proof}

The coefficient conditions are verified in this subsection. The proof of \cref{thm:sos-rg-input}
is completed after the physical identification in
\cref{prop:physical-sos-law}.

\subsubsection{SOS witnesses and the marked input}

Keep $B,R,L$ from \cref{sec:sos-mixed-coefficients} and put
$\delta^\mathrm{g}_\beta=\beta^{-1/10}$.
Let $\scW^{\textsc{sos}}_n$ consist of all connected torus patterns
of at least $B$ tiles having a six-class representation as in
\cref{lem:periodic-unions}, with remainder fraction at most
$C_*(\delta^\mathrm{g}_\beta+L/B)$ and total exceptional-orientation
fraction at most $C_*R^2/L$. Choose the fixed $C_*$ large enough
for unions of two families. Use the fixed orientation/direction
labels of \cref{prop:sos-witness} for the six classes, so their
exceptional orientations agree when class unions from different
patterns are combined. Distinct patterns are counted once. This
family contains all \texttt{fat} \SOS coefficient patterns. Write
$\cB_x^{\textsc{sos}}(s)$ for the event that one of its patterns
containing $x$, of size at least $s$, is present in the tiling. Its
rooted entropy rate tends to zero, and its patterns are
$(\epsilon_\beta,R/2)$-\texttt{fat} with $\epsilon_\beta\to0$.

\begin{lemma}[\SOS marked input]\label{lem:sos-marked-input}
The canonical \SOS coefficients and $\scW_n^{\textsc{sos}}$ satisfy
\cref{it:rg-marked}, uniformly in the root $x$ and the size
threshold $s$, including noncontractible supports.
\end{lemma}
\begin{proof}
The proof of \cref{lem:marked-union-envelopes} used two kinds of
input: coefficient bounds with rooted counts, for the numerical
estimate, and the six-class representation of \texttt{fat}
patterns and witnesses, for the count of \texttt{fat} supports. We
check that both are available here. Write $m=|U|$ and let $q$ be
the total size of the \texttt{thin}-labelled ordinary patterns in
a summand of \cref{eq:rg-absolute-insertion}. These patterns are
disjoint, so their total size, including the \texttt{fat} ones, is
at most $m$. The \SOS coefficients satisfy
\cref{eq:sos-mixed-envelopes} with the same rooted counts as in the
one-step case, so with $A_\beta=c\beta^{2/5}$ smaller than the \SOS
\texttt{thin} exponent, the product estimate
\cref{eq:marked-collection-budget} applies verbatim and gives
\[
 \sum_W\sum_{\{P_i\}}e^{A_\beta q/2}
                   \prod_i\bar c^{\mathrm{label}}_\beta(P_i)
                         \leq e^{e_\beta m},\qquad e_\beta\to0.
\]
The sum includes all internal histories through their already
summed coefficients. Declare the marked term \texttt{thin} when
$q\geq\delta^\mathrm{g}_\beta m$. Extracting its enhancement gives
$\exp[(e_\beta-A_\beta\delta^\mathrm{g}_\beta/2)m]
\leq e^{-c\beta^{3/10}m}$. Dropping the enhancement bounds the
complementary coefficient by $e^{e_\beta m}$.

For a complementary summand, unite the six classes of $W$ and
all its \texttt{fat}-labelled $P_i$. Their total size before union
is at most $|W|+\sum_i|P_i|\leq2m$. The fixed class labels ensure
that the exceptional orientation is the same in every union of a
given class. Their exceptional counts thus sum to at most
$CR^2m/L$. The old remainders and the \texttt{thin} patterns leave
at most $C(\delta^\mathrm{g}_\beta+L/B)m$ uncovered tiles.
With $l=\floor{\sqrt R}$, the hypotheses of
\cref{lem:periodic-unions} hold for large $\beta$, since
\[
 (\delta^\mathrm{g}_\beta+L/B)l^2\to0,
 \qquad R^2/L\to0,
 \qquad l^{-2}+l/R\to0.
\]
It gives vanishing rooted entropy; the boundary-layer estimate
shows that the patterns are
$(C/R+C(\delta^\mathrm{g}_\beta+L/B),R/2)$-\texttt{fat}, with the
first parameter tending to zero. The unrestricted connected-pattern count bounds
the marked \texttt{thin} family. Finally $U\supset W$ gives
$m\geq\max(s,B)$. These are precisely all the marked conditions,
and the quotient enumeration includes noncontractible unions.
\end{proof}

\subsection{Three steps: Ising Interface input}

\begin{theorem}[Ising Interface input to the \RG theorem]
\label{thm:interface-rg-input}
For sufficiently large $\beta$ and positive slopes in a fixed compact
set, the bubble-free Ising Interface model, with any admissible contact
perturbation $U$ of \cref{def:interface-perturbation}, admits a rooted
limit as the height cutoff $M\to\infty$, whose three-step image
$\varphi=\Tone\Ttwo\Tthree\cI$, in the diagonal translation quotient,
has law proportional to
\begin{equation}\label{eq:interface-image-weight}
 Z^U_{\beta,\Gamma}(\varphi)
   =\sum_{\Tone\Ttwo\Tthree\cI=\varphi}
             e^{-\beta(\sA(\cI)-N_0)+U_{\beta,\Gamma}(\cI)},
 \qquad \varphi\in\cM_\Gamma\,.
\end{equation}
All interface connectivity and no-bubble restrictions are retained.
This full three-step weight has a positive canonical extension and an
exact Grassmann expansion, allowing signed coefficients, with the
absolute \texttt{thin}/\texttt{fat} bounds, collection
control, source analyticity, gauge covariance, bulk consistency and
torus projection specified in \cref{thm:sos-rg-input}.
The corresponding constants have the required low-temperature
limits, uniformly in the prescribed slope range.
Thus the three-step image supplies the interface-model interaction
inputs \crefrange{it:rg-positive}{it:rg-sources} to
\cref{thm:coefficient-rg}. The witness family of
\cref{lem:interface-marked-input} satisfies \cref{it:rg-marked}.
\end{theorem}

Two features distinguish the three-step input from the previous
two. First, the perturbation $U$ of \cref{def:interface-perturbation},
which encodes the closed bubbles integrated out for full Ising in
\cref{sec:ising}, couples distinct defects of a history, so the
weight of a history is no longer a product of independent defect
weights; \cref{prop:interface-perturbed-absolute-history} restores a
hard-core structure by expanding the interaction into
inclusion--exclusion marks, at the price of signed terms. Second,
because the coefficients are signed, positivity of the matching
polynomial cannot be read off from the coefficients and must instead
come from the primitive sum of positive history weights
\cref{eq:interface-perturbed-positive-primitive-completion}, which is
why the following proposition constructs the positive weight on every
sector first and the coefficients afterwards.

\begin{proposition}[Positive canonical completion for a perturbed interface]
\label{prop:interface-perturbed-positive-completion}
For an admissible contact perturbation, the positive primitive sum
\cref{eq:interface-perturbed-positive-primitive-completion} is finite and obeys
\[
 e^{-\epsilon_tN_0}\leq F^+_{t,\Lambda}(M)
 \leq e^{(\epsilon_t+\epsilon^*_{t-\epsilon_t})N_0}\,,
\]
on every matching sector. On the prescribed sector it equals
$Z^U_{t,\Gamma}(\varphi)$, the sum \cref{eq:interface-image-weight}
with $t$ in place of $\beta$.
It has a signed connected-pattern expansion satisfying the absolute
\texttt{thin}/\texttt{fat} bounds, source identities, bulk consistency and torus
compatibility in \cref{thm:interface-rg-input}.
The coefficients are real and independent of the auxiliary dimer
weights. Their matching polynomial is $F^+$
on every sector, including at singular twists.
\end{proposition}
\begin{proof}
We proceed in five steps: the bounds on the primitive sum; its
rewriting as a signed hard-core gas of groups; the geometric
covering of a group and the resulting \texttt{thin}/\texttt{fat}
coefficient bounds; the identification of the coefficients as
functions of the pattern alone, with the canonical rule for small
patterns; and the Grassmann representation.

\begin{enumerate}[label=\textbf{Step~\arabic*.}, ref=\arabic*, wide=0pt, itemsep=1ex, listparindent=\parindent, parsep=0pt]
\item\label[step]{st:interface-completion-bounds}
(\emph{Bounds on the primitive sum.})
The empty history $\cX=\varnothing$ contributes
$e^{U_{t,\Gamma_M}(M)}$ to
\cref{eq:interface-perturbed-positive-primitive-completion}, and
$|U_{t,\Gamma_M}(M)|\leq\epsilon_tN_0$ by the contact variation
bound \cref{eq:interface-contact-variation}; this gives the lower
bound. For the upper bound, since $|\cI_{M,\cX}|=N_0+E(\cX)$, the
perturbation of a history is at most $\epsilon_t(N_0+E(\cX))$ in
absolute value by \cref{eq:interface-contact-variation}, so the sum
is at most $e^{\epsilon_tN_0}$ times the sum of
$e^{-(t-\epsilon_t)E(\cX)}$ over auxiliary histories, which
\cref{eq:interface-perturbed-auxiliary-budget} bounds by
$e^{\epsilon^*_{t-\epsilon_t}N_0}$. In particular the history sums
are absolutely convergent, which justifies the rearrangements below.
We construct the coefficients on every sector before any
coefficientwise extension.

\item\label[step]{st:interface-completion-gas}
(\emph{The signed hard-core gas of groups.})
Use the two types of marks in
\cref{prop:interface-perturbed-absolute-history}, now on the cylinder
$\Gamma_M$ of the matching itself: the reference marks, which
exponentiate the perturbation of the reference, and the interaction
marks, which carry the inclusion--exclusion coefficients
$D_X(\cV)$ of the change of $h_X$ produced by a subfamily $\cV$ of
components. The same argument with the auxiliary a priori bound gives, for every
finite set $W$ of reference plaquettes,
\begin{equation}\label{eq:interface-perturbed-all-sector-absolute-budget}
 \sum_{\text{objects in }W}|\mathrm{weight}|
 \leq\exp[(3\delta_t+\epsilon^*_{t-\delta_t})|W|],
 \qquad\delta_t\leq Ce^{-ct}\,.
\end{equation}
For $N_3=0$ there are no components and only reference marks. Group
components and mark occurrences by overlap of their reference
footprints, i.e.\ take the connected components of the graph joining
two objects when their footprints overlap. A nonzero interaction mark
$D_X(\cV)$ meets every component of $\cV$, by the vanishing assertion
of \cref{prop:interface-perturbed-absolute-history}, so all the
components it requires lie in its group, and its coefficient depends
only on these components. Hence the weight of a configuration factorizes
over its groups, and groups with disjoint footprints can be inserted
and deleted independently; a repeated mark and its factorial
$1/l!$ remain in one group. This proves that the primitive sum is an
exact hard-core gas of groups with signed activities.

\item\label[step]{st:interface-completion-geometry}
(\emph{Covering a group and the coefficient bounds.})
Let $H=E+M$ be a group's cost, with $E$ its history excess and $M$
the sum of its mark masses including multiplicity, and let $Q$ be
its footprint. The old interface footprints have the cover by rows
and dots of \cref{prop:interface-mixed-envelopes}. Expanded
footprints add at most $3E$ tiles, and the carriers of the marks at
most $CM$ tiles, since $|N_X|\leq C_0m_X$. Cover these additions by
singleton tiles, retaining each additional tile only once. The result is a cover of $Q$ by at most $CH$ rows
and singletons, of fixed multiplicity and total length at most
$C|Q|$. On an auxiliary sector a row of an active inverse step is
injective, because a pure period in its sorting direction is
excluded by \cref{lem:interface-perturbed-auxiliary-sorting}; an
inactive step contributes no rows. The incidence counts in that
lemma make the constants uniform even for short transverse periods.

Take $B=2\ceil{\sqrt\beta}$, $R=\ceil{\beta^{1/100}}$,
$L=\ceil{\beta^{1/10}}$ and $d_\beta=\beta^{-1/10}$, as in the
\SOS case. Call a group high-cost when $|Q|\leq BH$ and low-cost
otherwise. For a low-cost group, $|Q|>BH$, deleting the rows
shorter than $L$ from the cover leaves a remainder of at most
$CL|Q|/B$ tiles, since there are at most $CH<C|Q|/B$ rows. Apply
the six-class enlargement of \cref{prop:sos-witness} to the
remaining long rows and regroup overlaps of the enlarged patterns.
For an enlarged pattern $P$ of size $s$, let $q$ be the total size
of its high-cost base groups, and split its coefficient according
to whether $q\geq d_\beta s$. If $q\geq d_\beta s$, then the total
cost of the high-cost groups is at least $q/B\geq d_\beta s/B$, so
$H\geq d_\beta s/B$ for the whole collection. Every weight is
$B_0e^{-tH}$ with $B_0$ independent of $t$, so extracting half the
energy, i.e.\ writing $e^{-\beta H}=e^{-\beta H/2}e^{-\beta H/2}$ and
bounding the first factor by $e^{-\beta d_\beta s/(2B)}$, and then
applying \cref{eq:interface-perturbed-all-sector-absolute-budget}
at $t=\beta/2$ to the sum of the remaining weights, gives
\[
 |c_\beta^\mathrm{thin}(P)|
 \leq e^{-\beta d_\beta s/(2B)}
 e^{(3\delta_{\beta/2}+\epsilon^*_{\beta/2-\delta_{\beta/2}})s}
 \leq e^{-c\beta^{2/5}s}\,,
\]
since $\beta d_\beta/B\asymp\beta^{1-1/10-1/2}=\beta^{2/5}$. All
geometry and labels are held fixed at $\beta$ in this estimate.
Otherwise $q<d_\beta s$, so a low-cost group occurs, $s>B$, and the
six-class union has remainder fraction at most $d_\beta+CL/B$ and
exceptional count at most $CR^2s/L$. The union and net counts of
\cref{prop:sos-mixed-envelopes} show that the union is
$(\epsilon_\beta,R/2)$-\texttt{fat} with $\epsilon_\beta\to0$, and that
its rooted labelled entropy is $o_\beta(1)+C/n$. The
a priori bound \cref{eq:interface-perturbed-all-sector-absolute-budget}
at $t=\beta$ gives
\[
              |c_\beta^\mathrm{fat}(P)|\leq e^{\eta_\beta|P|},
                          \qquad\eta_\beta\longrightarrow0\,.
\]
These are bounds on the sums of absolute values of all labelled
terms, as required in \cref{it:rg-envelopes}.

\item\label[step]{st:interface-completion-canonical}
(\emph{Coefficients on patterns and the canonical rule.})
\Cref{lem:interface-perturbed-anchored-lifting} identifies the sums
defining $c_{\beta,\Lambda}(P)$ in any two matchings
containing $P$, preserving all marks, energies and geometric
labels. Thus they are coefficients on partial matchings, as
\cref{it:rg-positive} requires. Call a pattern small when its
diagonal diameter is at most $\eta_*n$, with $\eta_*$ as in
\cref{lem:sos-small-locality}, so that its lift embeds in a planar
ball and \cref{lem:interface-perturbed-anchored-lifting} applies
with trivial stabilizer $H_P=\{0\}$; the last assertion of that
lemma then identifies the torus coefficient of a small pattern with
its bulk coefficient. As in \cref{eq:canonical-coefficients}, for a
small pattern not extendable to a torus matching use its bulk
coefficient if it extends to a planar tiling, and zero otherwise;
put zero on other nonextendable patterns. These choices enter no
torus partition function. Large patterns extendable to auxiliary
sectors are retained, even if not extendable to the prescribed
physical sector, since the auxiliary sectors carry the primitive sum
as well. Small embedded patterns have exactly their bulk
coefficients, and every exception has size at least $c_*n$ by the
diameter bound used in \cref{prop:canonical-completion}; this is
\cref{it:rg-locality}. Absolute convergence
(\cref{st:interface-completion-bounds}) and the group factorization
(\cref{st:interface-completion-gas}) give
\cref{eq:rg-pattern-completion} for $F^+_{\beta,\Lambda}$ with the
coefficients $c_{\beta,\Lambda}$, simultaneously for every $M$.
Positivity of the left-hand side follows from the primitive sum, not
from the signed coefficients, which completes the verification of
\cref{it:rg-positive}.

\item\label[step]{st:interface-completion-grassmann}
(\emph{Grassmann representation.})
\Cref{lem:pattern-algebra}, which allows signed coefficients, gives
the Grassmann representation, the exact source factors and the gauge
covariance. The geometric sign and count argument of
\cref{prop:exact-sector} uses only the honeycomb torus and the
positivity of the weight on the prescribed sector, so it gives the
prescribed-sector ratio, including singular twists. The absolute
coefficient bounds of \cref{st:interface-completion-geometry} have
the same \texttt{thin}/\texttt{fat} parameters as in
\cref{prop:fermion-input}. In particular the Taylor factor
$e^{\rho|P|}$ is retained for the source-marked first integration;
no smallness at a fixed nonzero source radius is asserted here.
\qedhere
\end{enumerate}
\end{proof}

The coefficient conditions are verified in this subsection. The proof of \cref{thm:interface-rg-input}
is completed after the physical identification in
\cref{prop:physical-interface-law}.

\subsubsection{Ising Interface witnesses and the marked input}

For reconstruction we use positive geometric history components, not
the signed contact-marked expansion. Start from the expanded
components of \cref{prop:interface-perturbed-expanded-footprints},
with footprints $Q_r$,
declare one low-cost when $|Q_r|>BE_r$, enlarge it by the six-class
witness, and regroup overlaps. Denote a resulting group
by $\cG$, its enlarged pattern by $P(\cG)$, and its excess by
$E(\cG)$. Its data include all its original defects, but no
contact marks. This is a deterministic decomposition of each
primitive history. The groups are independently insertable and
deletable before the perturbation factor is applied. The cover by at
most $CE$ rows and dots, with bounded multiplicity, follows by adding
at most $3E$ singleton tiles to the old cover. The usual high-cost-mass
test gives \texttt{thin} and \texttt{fat} labels for these groups, using the same
$B,R,L,d_\beta$ as before.

Let $\scW_n^{U}$ contain all patterns of size at least $B$ with the
six-class representation, increasing its fixed constants to include
both these geometric \texttt{fat} groups and the \texttt{fat} coefficients of
\cref{prop:interface-perturbed-positive-completion}, and, as for
\SOS, to include unions of two such families, with remainder
fraction at most $C_*(d_\beta+L/B)$ and exceptional-orientation
fraction at most $C_*R^2/L$. Its rooted
entropy rate still vanishes and its patterns are
$(\epsilon_\beta,R/2)$-\texttt{fat} with $\epsilon_\beta\to0$. Let
$\cB_x^{U}(s)$ be the event that the tiling contains one of these
patterns through $x$ of size at least $s$.

\begin{lemma}[Interface marked input]\label{lem:interface-marked-input}
The interface coefficients and $\scW_n^{U}$ satisfy \cref{it:rg-marked},
uniformly in the root $x$ and the size threshold $s$, including noncontractible supports.
\end{lemma}
\begin{proof}
Use the absolute envelopes of
\cref{prop:interface-perturbed-positive-completion}, replacing every
internal labelled weight by its absolute value before forming
\cref{eq:rg-absolute-insertion}; this is permitted because
\cref{it:rg-marked} is formulated in terms of the envelopes
$q=q^\mathrm{t}+q^\mathrm{f}$ rather than the signed coefficients.
These envelopes satisfy the \texttt{thin} bound $e^{-c\beta^{2/5}|P|}$
and the \texttt{fat} bound $e^{\eta_\beta|P|}$ with the rooted counts
of \cref{st:interface-completion-geometry} of that proof, which is
all the numerical part of \cref{lem:marked-union-envelopes} uses.
That proof therefore applies to these envelopes: disjoint ordinary
patterns have total size at most the union size, and the witness
adds at most another copy of that size, giving the product estimate
\cref{eq:marked-collection-budget} and the split of the merged
coefficient according to the \texttt{thin} content. For the count
of \texttt{fat} supports, the \texttt{fat} envelope patterns and the
witnesses in $\scW_n^U$ have the fixed six-class representation, so
their class unions have the same vanishing geometric and entropy
errors as in the \SOS case, while the enhanced \texttt{thin}
coefficient gives the exponent $c\beta^{3/10}$. Every marked union
contains its witness and has size at least $\max(s,B)$, as required in
\cref{it:rg-marked}.
\end{proof}

\section{Renormalization analysis}\label{sec:rg}

\begin{remark*}[Reader's guide]
\normalfont
This section proves \cref{thm:coefficient-rg}, the precise version
of \cref{thm:rg-typical}, independently of the surface models.
The \RG analysis of interacting dimers in \cite{GMT17,GMT20}
supplies the multiscale expansion and its continuum reference model.
The refinement in \cite[\S8]{LaslierLubetzky24} treats exponentially
decaying interactions and prescribed slopes. Our coefficients are
different: a large coefficient may be supported on a rare, large
pattern, so the uniform interaction bound of that refinement does
not apply directly. The proof proceeds as follows.
\begin{enumerate}[(i)]
\item \emph{Turn rare patterns into small initial kernels.}
\Cref{lem:entropy-patch,lem:fat-coverage} give the free
\texttt{fat}-pattern probability bound of \cref{prop:fat_pattern}.
Combining this suppression with the coefficient bounds yields the
bulk and torus first integrations in
\cref{thm:bulk-first-step,thm:torus-first-step}.
\item \emph{Construct the critical flow for the actual interaction.}
\Cref{lem:infrared-localization} identifies the running parameters,
and \cref{thm:conditional-renormalized-trees} places the first-step
kernels in the induction of \cite[\S6]{GMT20}. The reference
cancellation and the two-history bounds of \cref{lem:two-history-beta}
give \cref{thm:critical-flow-fixed-dispersion}.
\Cref{thm:matched-bulk-flow} then matches the tuned quadratic data
to the original measure.
\item \emph{Obtain the local limit at the prescribed slope.}
The bulk source functional is constructed in
\cref{thm:bulk-source-functional}.
\Cref{lem:finite-size-source-kernels} compares torus and bulk
kernels, leading to \cref{thm:canonical-local-limit}.
The lower bound on the target sector in
\cref{prop:interacting-central-sector} and the uniqueness argument
in \cref{lem:pressure-equal-density} then give the fixed-sector
limit, \cref{thm:fixed-sector-limit}.
\item \emph{Control witnesses under the interacting measure.}
Under the additional marked-collection hypothesis,
\cref{lem:overlapping-witness-insertion} represents the number
of witnesses by an overlapping-pattern insertion. Its first integration
and \cref{lem:one-insertion-rg} give
\cref{thm:interacting-witness-tail}. This is the estimate needed
for the surface reconstruction in \cref{sec:transfer}.
\item \emph{Identify the Gaussian height limit.}
\Cref{lem:one-source-kernel-comparison} leads to the sharp source
flow of \cref{thm:sharp-source-flow}. The two-edge asymptotics and
\cref{lem:curl-covariance-algebra} identify the height covariance
in \cref{cor:height-covariance-limit}.
\Cref{lem:reference-density-gaussian,lem:integrated-source-blocks}
give the higher-cumulant estimates, and
\cref{lem:band-limited-height-current} converts them into the
height \GFF theorem, \cref{thm:dimer-height-gff}.
\end{enumerate}
The three conclusions of \cref{thm:coefficient-rg} are assembled
in \cref{sec:rg-proof-completion}.
\end{remark*}

\begin{theorem}[Interacting dimers with absolute \texttt{thin}/\texttt{fat} coefficients]
\label{thm:coefficient-rg}
Let $F_{\beta,\Lambda}$ be an admissible family of honeycomb-dimer
weights in \cref{def:rg-admissible}. Fix a compact
subset of the interior of the density simplex. For sufficiently large
$\beta$, uniformly on this subset, the following hold.
\begin{enumerate}[(i)]
\item On the diagonal tori with periods
$\Gamma=\langle(-n,0,\Theta_1),(0,-n,\Theta_2)\rangle$, condition the weight
$F_{\beta,\Lambda}$ on orientation counts $(n\Theta_1,n\Theta_2,n^2)$.
If $\Theta_1/n\to\theta_1>0$ and $\Theta_2/n\to\theta_2>0$, these laws converge
locally to a unique stationary mixing law $\nu_{\beta,\rho}$, where
$\rho=(\theta_1,\theta_2,1)/(1+\theta_1+\theta_2)$.
The limit is independent of the approximating integer windings.
It is a canonical infinite-volume law at suitable positive liquid
edge weights. Its edge cumulants are the bulk source coefficients
of \cref{thm:bulk-source-functional} and satisfy
\cref{eq:bulk-response-decay}.
\item For the centered diagonal height in $\nu_{\beta,\rho}$,
the rescaling of \cref{eq:test-field} converges, jointly on smooth
compactly supported integral-zero tests and in $H^{-s}_\mathrm{loc}$
modulo constants for every $s>1$, to
\[
 \upsigma_{\beta,\rho}\,\GFF\circ L_{\beta,\rho},\qquad
 \upsigma_{\beta,\rho}=\sqrt{a_*/|\det\sA_\mathrm{D}|}>0,
 \qquad L_{\beta,\rho}=\sA_\mathrm{D}^{-\trps}\,.
\]
Here $\sA_\mathrm{D}$ is the invertible real matrix of the matched
linear dispersion and $a_*>0$ is the covariance coefficient in
\cref{cor:height-covariance-limit}.
\item If a family $\scW_\Lambda$ additionally satisfies the
marked-collection condition \cref{it:rg-marked}, then, for all
sufficiently large tori, both the canonical laws at liquid weights
in a fixed compact set and the prescribed-sector laws satisfy
\[
 \P\big(\exists W\in\scW_\Lambda:\ x\in S(W),\
                   |W|\geq s,\ W\subset M\big)
                 \leq C e^{-c\max(s,B_\beta)}\,.
\]
The estimate is uniform in $x,s$ and includes noncontractible witnesses.
\end{enumerate}
\end{theorem}

This section proves \cref{thm:coefficient-rg} for the coefficient
class of \cref{def:rg-admissible}, independently of its surface
origin. The applications are verified in
\cref{sec:model-grassmann,sec:ising}; the proof of the abstract
theorem is assembled in \cref{sec:rg-proof-completion}.

\begin{remark}[The role of the conditions]\label{rem:rg-conditions-role}
We comment on what each requirement in \cref{def:rg-admissible}
is used for, and on which application forces it.
\Cref{it:rg-positive} prescribes the form of the weight: a hard-core
gas of disjoint patterns, rather than an exponential of a sum of local
functions as in \cite{LaslierLubetzky24}. This is the form in which
the fibers of the \Tetris histories arrive
(\cref{thm:hmsos-rg-input} and its two- and three-step counterparts):
the sum over the defects with a given footprint has already
been performed inside the coefficients, and, as explained in
\cref{rem:pattern-polynomial-vs-exponential}, no logarithm is taken.
Positivity of $F_{\beta,\Lambda}$ is needed because it is a
probability weight; the individual coefficients are nevertheless
allowed to be signed, since the perturbed interface expansion of
\cref{sec:ising}, in which closed bubbles are integrated out, produces
coefficients that need not be nonnegative. The coefficients are
produced by the surface fibers, which know nothing of the auxiliary
edge weights $\omega_r$, the twists $y$ and the edge sources $A$
introduced below to extract the prescribed sector and the
correlations; this independence is what makes the Grassmann
representation of \cref{it:rg-sources}, in which every source acts
edge by edge, an identity rather than an approximation.

\Cref{it:rg-envelopes} is the heart of the matter. The \texttt{thin}
part is a small-activity gas of the familiar kind: the coefficient
decays exponentially in the area at a rate diverging with $\beta$,
and connected patterns through a given triangle are counted by
$C_0^m$, so that their sum is small as soon as $\beta$ is large. The
\texttt{fat} part is new. Its coefficient may grow like
$e^{\epsilon_\beta|P|}$ and its rooted count like $e^{\epsilon_\beta m}$,
and the only quantity available to compensate both is the free dimer
law itself: by \cref{prop:fat_pattern} below, a prescribed pattern
whose support is almost entirely covered by large disks occurs under
the free law with probability at most $e^{-(P_*-\delta)|P|}$, at a rate
$P_*>0$ that does not depend on $\beta$. This is why the \texttt{fat} property, and not
merely size, is required. The minimum size $B_\beta\to\infty$ makes
the resulting geometric series $\sum_{m\geq B_\beta}e^{-pm}$ small,
and the term $C/n$ in the torus count accounts for noncontractible patterns,
cf.~\cref{prop:canonical-completion}. The bounds concern the total
coefficients only; see the remark following \cref{def:rg-admissible}
for how they are verified from signed labelled terms.

\Cref{it:rg-locality} is what allows the torus kernels to be compared
with the planar ones. A connected pattern with fewer than $c_*n$
edges cannot wrap around the torus, so it lifts to the plane and
carries its planar coefficient; all torus-specific inputs have size
at least $c_*n$ and fall into the tails of the convergent expansions.
This is used in \cref{thm:torus-first-step} and in the finite-volume
comparison leading to \cref{thm:canonical-local-limit}, and it is the
reason why the limit in \cref{thm:coefficient-rg} does not depend on
the approximating windings.

\Cref{it:rg-sources} fixes the generating function. The twists
$y_2,y_3$ allow the orientation counts to be read off by Fourier
projection, and the derivatives in the edge sources $A$ at $A=0$
produce the edge correlations, from which the height correlations
are assembled. The absolute Taylor factor $e^{\rho_0|P|}$ is the price
of a source on a pattern all of whose edges are marked; it enters
the node estimates of \cref{sec:first-integration} and is absorbed
by choosing $\rho_0$ small compared with the fixed negative margin
in the \texttt{fat}-node bound. The last sentence of
\cref{it:rg-sources} is a warning: a bound on the source dependence
for each fixed set of marked locations would not control the sum
over the locations, which is what the anchored norms below require.

The condition \cref{it:rg-marked} is used only for the witness
conclusion. It asks for the \texttt{thin}/\texttt{fat} bounds of the
coefficient of the pattern obtained by merging a witness with the
ordinary patterns it overlaps, which is precisely the object inserted
in \cref{lem:overlapping-witness-insertion}; the witness families
themselves are constructed in \cref{sec:preimage-witnesses}.
\end{remark}

\subsection{An exponential bound for \texorpdfstring{\texttt{fat}}{fat} dimer patterns}\label{sec:fat-determinant}

The purpose of this subsection is the estimate of
\cref{prop:fat_pattern}: under the free dimer law, a prescribed
\texttt{fat} pattern occurs with probability exponentially small in
its area, at a rate bounded below uniformly on compact sets of liquid
slopes. As explained in \cref{rem:rg-conditions-role}, this is the
mechanism that compensates both the possible growth of a
\texttt{fat} coefficient and the number of \texttt{fat} patterns.
The estimate is stated for the free law; the proof of
\cref{prop:conditional-first-step} below reproduces its mechanism
inside the Grassmann integral, where the free law is no longer
available, and this is why we prove it through determinants rather
than through a direct entropy argument.

\subsubsection*{Conventions}
We fix the conventions needed in the determinant argument. A pattern
$\cP$ is a finite set of disjoint dimers, $S$ is its vertex support,
and $|\cP|=|S|/2$ is the number of dimers. Work with positive honeycomb
weights $a,b,c$ in the liquid regime, that is, satisfying the strict
triangle inequalities $a<b+c$, $b<a+c$ and $c<a+b$, normalized by
$\max(a,b,c)=1$. The weighted Kasteleyn operator $K$ of the honeycomb
lattice is translation invariant, with symbol
\[
 \mu(k)=a+be^{ik_1}+ce^{ik_2}\,,
\]
so that the translation-invariant inverse $G$ is the Fourier inverse
of $1/\mu$. The equation $\mu(k)=0$ asks that three vectors of lengths
$a,b,c$ close into a triangle, which is possible precisely because
of the triangle inequalities; its two solutions $k_\pm$ correspond to
the two orientations of this triangle and are exchanged by
$k\mapsto-k$. Regarded as a map from $[-\pi,\pi]^2$ to $\R^2\cong\C$,
the symbol has an invertible real Jacobian at each zero (its
determinant is $bc\sin(k_2-k_1)$, which does not vanish since the
triangle is nondegenerate), which is what we mean by saying that the
two zeros are simple as real two-dimensional zeros. Hence
$|\mu(k)|\geq c|k-k_\pm|$ near the zeros, and consequently
$\int|\mu|^{-1}\,\d k/(2\pi)^2<\infty$. Let
\[
 P=P(a,b,c)=\int_{[-\pi,\pi]^2}\log|\mu(k)|\frac{\d k}{(2\pi)^2}\,.
\]
We use the standard free-dimer identities: the pattern probability is
$\prod_{e\in\cP}|K(e)|\,|\det G_S|$, and the entropy $h$ and edge
densities $\rho_i$ of this Gibbs measure satisfy
$P=h+\rho_a\log a+\rho_b\log b+\rho_c\log c$.
These are the Kasteleyn probability and pressure formulas (see
\cite{KOS06,Kenyon09}); the same Kasteleyn conventions are used in
\cite[\S8.1]{LaslierLubetzky24}. To evaluate $P$, relabel so that
$a=1$ is the largest weight, substitute $t=k_2-k_1$ and average first
over $k_1$: for fixed $t$ the integrand is $\log|1+e^{ik_1}w|$ with
$w=b+ce^{it}$, and Jensen's formula
$\frac1{2\pi}\int_{-\pi}^{\pi}\log|1+we^{i\theta}|\,\d\theta
=\log\max\{1,|w|\}$ gives
\[
 P=\frac1{2\pi}\int_{-\pi}^{\pi}
                  \log\max\{1,|b+ce^{it}|\}\,\d t>0\,.
\]
The integrand is nonnegative, and it is positive for $t$ near zero
because $b+c>1$ by the triangle inequality; hence $P>0$.
For equal weights, $P$ is the usual entropy per dimer. For arbitrary
slope the rate below is $P$ in this normalized gauge, not the unweighted
entropy at that slope: the probability of a pattern under the weighted
measure is a product of edge weights times a minor of $G$, and it is
the latter that the determinant argument below controls.

These facts are uniform on any compact set $\mathcal K$ of
normalized liquid weights. The two zeros stay separated and
their real Jacobians stay nonsingular. Uniform inverse-function
bounds give $|\mu(k)|\geq c\,\dist(k,\{k_+,k_-\})$ near
the zeros, while $|\mu|$ is uniformly bounded below away from
them. Integrating these bounds shows that
\[
 \sup_{w\in\mathcal K}\int|\mu_w|^{-1}\frac{\d k}{(2\pi)^2}
 <\infty,\qquad
 P_*:=\min_{w\in\mathcal K}P(w)>0.
\]
Here continuity of $P$ follows from the same bounds: the
integral of $|\log|\mu_w||$ in radius-$r$ neighborhoods of
the zeros is uniformly $O(r^2(1+|\log r|))$, and outside
these neighborhoods ordinary dominated convergence applies.
Compact sets of interior slopes give such compact weight
sets, since $w_i$ is proportional to $\sin(\pi\rho_i)$,
where $\rho_i$ are the three lozenge densities, and we
normalize the largest weight to one.

The first ingredient is a region whose partition function is as
large as the pressure predicts. An $L\times L$ square of fundamental
cells itself need not do: the boundary conditions imposed by its
shape can force a different slope, hence a smaller entropy, in a
macroscopic part of it, so that its partition function at the given
weights may be exponentially smaller than $e^{PL^2}$. The remedy is
to let the infinite-volume Gibbs measure itself select which
boundary vertices to remove.

\begin{lemma}[A large partition function with a selected boundary]\label{lem:entropy-patch}
For every $\delta>0$ and all sufficiently large $L$, an $L\times L$
array $Q_L$ of fundamental cells contains a balanced vertex set $D_L$
with at least $L^2-CL$ vertices of each color such that
\[
 |\det K_{D_L}|=Z(D_L)\geq e^{(P-\delta)L^2}\,.
\]
Only vertices in a fixed-width boundary layer of $Q_L$ are removed.
The threshold for $L$ can be chosen uniformly for weights in
a compact set $\mathcal K$ as above; the selected $D_L$ may
depend on the weights.
\end{lemma}
\begin{proof}
Sample the infinite-volume Gibbs matching. Let $X_L$ record the
matching edge at every black vertex of $Q_L$, and let $h$ be its
specific entropy. Entropy subadditivity gives
\[
                         H_e(X_L)\geq hL^2.
\]
Indeed, tile a much larger square by translates of $Q_L$, bound
its entropy by the sum of the block entropies and the entropy of
the remaining boundary strip, and divide by its area before taking
the limit. Thus no rate of convergence of block entropies is needed.

Let $\eta$ be the set of boundary-crossing edges. There are at most
$e^{CL}$ possible values of $\eta$. Removing their endpoints in $Q_L$ leaves
a balanced set $D(\eta)$, with only $O(L)$ vertices of each color
removed. Let $M$ be the matching inside this set. The pair $(M,\eta)$
determines $X_L$, so the entropy chain rule gives
\[
 H_e(M\mid\eta)=H_e(M,\eta)-H_e(\eta)
                                      \geq hL^2-CL.
\]
By the Gibbs property, conditionally on $\eta$ the matching $M$ has
partition function $Z(D(\eta))$. Averaging the Gibbs entropy identity
therefore yields
\[
 \E\log Z(D(\eta))
   =H_e(M\mid\eta)+\E\sum_{e\in M}\log|K(e)|
   \geq \Big(h+\sum_i\rho_i\log w_i\Big)L^2-CL
   =PL^2-CL.
\]
For the inequality, stationarity gives the mean log weight of $X_L$
as $L^2\sum_i\rho_i\log w_i$. Passing from $X_L$ to $M$ removes
only crossing edges, whose log weights are nonpositive because
$\max_i w_i=1$, so it can only increase this sum. Some value of
$\eta$ of positive probability thus satisfies
$\log Z(D(\eta))\geq PL^2-CL$. The constant $C$ counts the values
of $\eta$ and is independent of the weights; taking $L\geq C/\delta$
gives the claimed lower bound, uniformly on the stated compact set.

Finally, every deleted vertex is joined to the exterior of $Q_L$
by its recorded crossing edge. Deleting these vertices merges their
incident faces with the unbounded face and creates no new inner
face. The restricted planar Kasteleyn signing therefore remains
valid, and all perfect-matching terms have the same phase. Hence
$|\det K_{D(\eta)}|=Z(D(\eta))$, as required.
\end{proof}

\begin{lemma}[Interior coverage of a \texttt{fat} set]\label{lem:fat-coverage}
For every $\eta>0$ and fixed $L$, if $\epsilon$ is sufficiently small
and $R$ is sufficiently large compared with $L/\eta$, the support of
an $(\epsilon,R)$-\texttt{fat} pattern contains disjoint translates of $Q_L$
covering all but $\eta|S|$ vertices. Placing $D_L$ in these squares
leaves at most $C(\eta+L^{-1})|S|$ vertices outside their interiors,
where an interior vertex has all its lattice neighbors in its patch.
\end{lemma}
\begin{proof}
Let $U$ be the union of radius-$R$ balls in the \texttt{fat} support.
Retain the squares of a fixed $L\times L$ fundamental-cell grid
contained in $U$. The missed part of $U$ lies in its inner
$CL$-boundary layer. By \cref{eq:ball-union-boundary-layer},
including the exceptional part of the support, the number of
missed vertices is at most $C(\epsilon+L/R)|S|$.
Here fixed-width boundary layers convert planar area to vertex counts
and account for tiling closure, at an additional cost
$O(R^{-1})|S|$. Choosing $\epsilon\ll\eta$ and $R\gg L/\eta$
proves the first assertion. Each selected $D_L$ loses only
$O(L)$ vertices, including its boundary, from an $L\times L$
array of cells. Summing gives the additional $O(L^{-1})|S|$ loss.
\end{proof}

\begin{proposition}[\texttt{Fat}-pattern probability]\label{prop:fat_pattern}
For every $\delta>0$ there are $\epsilon>0$ and $R<\infty$ such that
every $(\epsilon,R)$-\texttt{fat} pattern $\cP$ satisfies
\[
             \mu(\cP\subset\varphi)\leq e^{-(P-\delta)|\cP|}\,.
\]
The constants can be chosen uniformly on compact sets of liquid
slopes. In particular, for $\delta<P_*$ the rate is uniformly
strictly positive on such a set.
\end{proposition}
\subsubsection*{Proof idea}
The pattern probability is $\prod_{e\in\cP}|K(e)|\,|\det G_S|$, and
the edge weights are at most one, so it suffices to show that the
$n\times n$ minor $G_S$ of the inverse Kasteleyn matrix is
exponentially small. A direct Gram--Hadamard bound only gives
$|\det G_S|\leq C_0^{n}$, with no decay. The idea is to multiply
$G_S$ by a matrix $K_I$ whose determinant is as large as
$e^{(P-\eta)n}$ and such that $K_IG_S$ is, up to a small minor, the
identity: on each entropy patch of \cref{lem:entropy-patch}, the
Kasteleyn matrix of the patch inverts $G$ exactly on interior rows,
by the identity $KG=1$ of the full lattice. Thus
$|\det G_S|=|\det(K_IG_S)|/|\det K_I|$, the numerator is at most
$C_1^{\eta n}$ and the denominator at least $e^{(P-\eta)n}$.

\begin{proof}
Let $n=|\cP|$ and select disjoint patches $D_L$ inside $S$ by
\cref{lem:fat-coverage}. Write $I$ for their union and define an
$n\times n$ matrix $K_I$, indexed by the black and white vertices of
$S$, to be block diagonal: the blocks on the patches are their
Kasteleyn matrices $K_{D_L}$; the remaining black and white vertices
are paired arbitrarily in an identity block. They have equal numbers
because every patch is balanced and $S$, being the support of a
partial matching, is balanced. By \cref{lem:fat-coverage}, the
patches cover $1-o(1)$ of $S$ as first $L\to\infty$ and then
$\epsilon\to0$, $R\to\infty$, and by \cref{lem:entropy-patch} each
patch contributes a factor at least $e^{(P-\delta')L^2}$ to
$|\det K_I|$. Thus, for any prescribed small $\eta>0$, the parameters
can be chosen so that
\[
 \log|\det K_I|\geq(P-\eta)n,
 \qquad\#\{\text{noninterior black rows}\}\leq\eta n\,.
\]
Now consider the product $K_IG_S$. For an interior black vertex $b$
of a patch, all three lattice neighbors of $b$ lie in the patch, so
the row of $K_I$ at $b$ coincides with the row of the full Kasteleyn
operator $K$ at $b$ restricted to $S$; since $KG=1$ on the full
lattice, the corresponding row of $K_IG_S$ is the identity row
$\delta_{b,\cdot}$. Expanding the determinant successively along
these identity rows leaves a square minor of $K_IG_S$ of size at most
$\eta n$, indexed by the noninterior black rows and the columns not
eliminated.

To bound this minor use the Gram--Hadamard inequality. The inverse
$G$ has white rows and black columns. In
$L^2(\T^2,\d k/(2\pi)^2)$, with the inner product conjugate-linear
in its first slot, put
\[
 a_w(k)=\frac{e^{-ik\cdot w}}{\sqrt{|\mu(k)|}},\qquad
 b_b(k)=\frac{e^{-ik\cdot b}\sqrt{|\mu(k)|}}{\mu(k)}\,.
\]
Then $\langle a_w,b_b\rangle=G(w,b)$ by Fourier inversion, and
both squared norms equal
$C_0=\int|\mu|^{-1}\,\d k/(2\pi)^2<\infty$. A
remaining row of $K_IG_S$ is either a row of $G_S$ (for a vertex in
the identity block) or a linear combination of at most three such
rows with coefficients of modulus at most one (for a noninterior
vertex of a patch, whose row of $K_{D_L}$ has at most three nonzero
entries). The remaining minor is therefore a Gram matrix with row
norms at most $3C_0^{1/2}$ and column norms at most $C_0^{1/2}$, and
the Gram--Hadamard inequality bounds its determinant by
$C_1^{\eta n}$ with $C_1=\max(1,3C_0)$. Since
$\prod_{e\in\cP}|K(e)|\leq1$, we conclude that
\[
 \mu(\cP\subset\varphi)
 \leq\frac{|\det(K_I G_S)|}{|\det K_I|}
 \leq \exp(-[P-\eta-\eta\log C_1]n)\,,
\]
and it remains to choose $\eta(1+\log C_1)<\delta$.

For the uniformity, on a compact weight set choose $C_1$ using the
uniform bound on $\int|\mu_w|^{-1}$ established above, and choose
$\eta$ once. The uniform version of \cref{lem:entropy-patch} then
supplies one sufficiently large $L$, and \cref{lem:fat-coverage},
which is purely geometric, supplies common $\epsilon,R$. This proves
the asserted uniformity, including the positive rate $P_*-\delta$.
\end{proof}

\begin{corollary}[Annealed pressure and excess energy]\label{cor:mixed-pressure}
Under the infinite-volume free liquid tiling law $\mu_0$, restrict
defects by $Q_h(D)\subset\Omega$, where $\Omega$ is a finite set of
triangles. There are $c,C>0$ such that
\begin{equation}\label{eq:mixed-pressure}
 0\leq\log\mu_0(Z_{\beta,h})
          \leq C|\Omega|e^{-c\beta^{2/5}}\,.
\end{equation}
The normalized joint law of the reference tiling and its defect
family, with density proportional to $e^{-\beta E}$ relative to
$\mu_0$ and counting measure on defects, satisfies
\begin{equation}\label{eq:annealed-energy}
       \E_{\beta,\Omega}E
           \leq C|\Omega|\beta^{-1}e^{-c\beta^{2/5}}\,.
\end{equation}
In the second assertion the allowed supports are fixed at the
stated value of $\beta$ when comparing temperatures.
\end{corollary}
\begin{proof}
Expand the positive gas \cref{eq:mixed-gas} into disjoint patterns
labelled \texttt{thin} or \texttt{fat}. A union of disjoint
\texttt{fat} patterns is fat with the same parameters: take the
union of their ball cores and add their exceptional areas.
By \cref{prop:fat_pattern}, the free probability that the entire
collection occurs is therefore at most
$e^{-p\sum_{i\text{ \texttt{fat}}}|P_i|}$ for some fixed $p>0$,
uniformly on compact liquid slope sets. Requiring the thin patterns
as well only decreases this probability. Assign a deterministic root
to each pattern and drop disjointness. The collection sum is then
bounded by a product of $1+$activity, and $\log(1+x)\leq x$ and
\cref{prop:mixed-envelopes} give
\[
 0\leq\log\mu_0(Z_{\beta,h})
 \leq C|\Omega|\left(
       \sum_{m\geq1}C_0^m e^{-c\beta^{2/5}m}
       +\sum_{m\geq B_\beta}e^{-(p-o(1))m}\right)
 \leq C|\Omega|e^{-c'\beta^{2/5}}.
\]
Here each root lies in $\Omega$, the empty collection gives the lower
bound, and $B_\beta=\lceil\sqrt\beta\rceil$ makes the second series
smaller than the stated error.

For the energy bound keep the supports, enlarged patterns and split
fixed at their values for $\beta$, and put
$f(t)=\log\mu_0(Z_{t,h})$. The coefficient proof uses the extracted
factor $t/(4B_\beta)$ and the a priori bound at $t/2$, so the same
bound holds uniformly for $\beta/2\leq t\leq\beta$. Since
$\E_{\beta,\Omega}e^{\beta E/2}
 =e^{f(\beta/2)-f(\beta)}<\infty$, Jensen's inequality gives
\[
 \E_{\beta,\Omega}E
 \leq\frac2\beta\log\E_{\beta,\Omega}e^{\beta E/2}
 =\frac2\beta\big(f(\beta/2)-f(\beta)\big)
 \leq C|\Omega|\beta^{-1}e^{-c\beta^{2/5}}.
\]
This is a total-energy bound, not a local insertion estimate.
\end{proof}

\subsection{Grassmann expansion and its required bounds}\label{sec:grassmann-input}

We follow the Grassmann conventions and sector-extraction strategy of
\cite[\S8.2--8.4]{LaslierLubetzky24}, with the exact pattern
coefficients of \cref{def:rg-admissible}. We recall the setting.
The partition function of a weighted dimer model on a finite planar
bipartite graph is, by Kasteleyn's theorem, the absolute value of
the determinant of the Kasteleyn matrix $K$, and a determinant is a
Gaussian Grassmann integral: in a planar finite domain with an
invertible Kasteleyn matrix, choose the order of Berezin integration
so that $\int e^{\psi_bK\psi_w}\,\d\psi=\det K$, where
$\psi_bK\psi_w$ abbreviates $\sum_{b,w}\psi_bK(b,w)\psi_w$. The
interaction, which multiplies the weight of a matching $M$ by
$F(M)$, becomes an even element $V(\psi)$ of the Grassmann algebra
by \cref{lem:pattern-algebra}, and the normalized Gaussian
expectation $P_{K^{-1}}$ with covariance $K^{-1}$, characterized
by Wick's rule, then satisfies
\begin{equation}\label{eq:normalized-grassmann}
 Z_\mathrm{int}=\int e^{\psi_bK\psi_w+V(\psi)}\,\d\psi
       =\det K\int e^{V(\psi)}\,\d P_{K^{-1}}(\psi)\,.
\end{equation}
There is no additional factor $Z_\mathrm{dimers}$ in the first
integral: the free partition function is exactly $\det K$. On a
torus the situation is more involved. Kasteleyn's theorem there
expresses the partition function as a signed combination of four
determinants, each of which is itself a signed sum over the winding
sectors of the matchings, and the orientation counts that our
surfaces prescribe correspond to one such sector. Recovering fixed
lozenge counts therefore requires a Fourier projection in the
twists $y_2,y_3$ of \cref{it:rg-sources}, together with a tuning of
the canonical weights $\omega_r$ to the interacting slope, so that
the prescribed sector carries a non-negligible fraction of the
canonical measure; this is carried out in
\cref{sec:fixed-sector-limit} following
\cite[\S8.6]{LaslierLubetzky24}. This step requires estimates that
are uniform in the twists, and it cannot be replaced by one planar
determinant.

The coefficient bounds consumed by the first integration are recorded
in \cref{prop:fermion-input}. They apply to signed coefficients through
the absolute envelopes of \cref{it:rg-envelopes}.

In the source convention here $w_A$ denotes interaction-derived
coefficients. The free edge source is kept explicitly as
$\sum_e K(e)(e^{A(e)}-1)\psi_b\psi_w$, or equivalently in the quadratic
matrix $K_A$. These extra two-field activities are $O(\rho)$ on a
polydisc of radius $\rho$; they are not $O(\epsilon_\beta)$, since
they do not come from the interaction, and this is why they are
treated separately from the interaction nodes in all the estimates
below. A small fixed $\rho$ is chosen separately. Gauge covariance of
$w_A$ below is understood with the quadratic part $K_A$.

The correlations we eventually need are derivatives in the sources
at $A=0$, so we must convert bounds on the source-dependent
coefficients into bounds on their derivatives. The following lemma
is our only tool for this purpose, and the remark after its
proof explains why a more convenient relative estimate is not
available.

\begin{lemma}[Cauchy bound for source derivatives]\label{lem:source-envelope}
If $w_A(S)$ is analytic for $|A(e)|<\rho$ on a finite source set $F$
and $|w_A(S)|\leq M(S)$ there, then for every multi-index $m$,
\[
 |\partial_A^m w_0(S)|\leq m!\rho^{-|m|}M(S)\,.
\]
A relative estimate $|\partial_Aw_A(S)|\leq C|w_A(S)|$ does not follow
and in general is false for signed coefficients.
\end{lemma}
\begin{proof}
Apply Cauchy's formula successively in the variables $A(e)$ on circles
of radius less than $\rho$, then let their radii increase to $\rho$.
For the last assertion the analytic function $w(A)=A$ has $w(0)=0$
and $w'(0)=1$, so no bound of the form $|w'(0)|\leq C|w(0)|$ can
hold. Cancellation in a resummed signed coefficient can produce
precisely this issue: a coefficient may vanish at $A=0$ because of a
cancellation between its internal terms, while its source derivative
does not vanish. This is why all source estimates below are
formulated through envelopes $M(S)$ valid on a whole polydisc, and
never relative to the value at $A=0$.
\end{proof}

The next lemma records the exact gauge identity of the dimer
partition function and its Grassmann counterpart. Its role is to
provide the lattice Ward identities used later in the identification
of the two-edge asymptotics and of the height covariance; we state it
here because its hypothesis, the gauge covariance of the resummed
coefficients, is a property that must be verified for the actual
expansion at every step rather than deduced from size bounds.

\begin{lemma}[Gauge identity]\label{lem:gauge-identity}
The exact dimer partition function with arbitrary fixed interaction
and edge sources satisfies
\[
 Z(A+t\one_{\{e\ni v\}})=e^t Z(A),\qquad
 \sum_{e\ni v}\partial_{A(e)}\log Z(A)=1\,.
\]
If a Grassmann representation with auxiliary sources $J_v$ has
$w_{A+t\one_{\{e\ni v\}}}(S)=e^{t\one_{\{v\in S\}}}w_A(S)$,
then
\[
 Z(A+t\one_{\{e\ni v\}},J)=e^tZ(A,J^{(v,t)}),
 \quad J^{(v,t)}_v=e^{-t}J_v\,,
\]
and other source coordinates are unchanged. Differentiating gives
the lattice Ward identities, including the contact term $-1$ when a
fermion insertion is at $v$.
\end{lemma}
\begin{proof}
Every perfect matching has exactly one edge incident to $v$.
Increasing $A$ by $t$ on these edges therefore multiplies each
matching weight by $e^t$, proving the first identity term by
term. The second follows by differentiating its logarithm at
$t=0$. For the Grassmann identity set $\psi_v=e^{-t}\xi_v$ and leave
all other variables unchanged. In every monomial $\psi_S$ with
$v\in S$ this produces a factor $e^{-t}$, which cancels the assumed
gauge factor $e^{t\one_{\{v\in S\}}}$ of its coefficient, while the
Berezin differential transforms as $\d\psi_v=e^t\,\d\xi_v$,
producing the overall factor $e^t$. The source term at $v$ becomes
$e^{-t}\xi_vJ_v$, which is the stated rescaling of $J_v$. This proves
the identity; differentiating its logarithm in $t$ at $t=0$, and then
in the source variables, proves the stated contact term: the
derivative of $J^{(v,t)}_v=e^{-t}J_v$ contributes $-1$ exactly when
a fermion insertion sits at $v$. We stress that gauge covariance
must be checked for the actual resummation; it is not implied by a
bound on $|w_A(S)|$.
\end{proof}

\subsection{The first integration}\label{sec:first-integration}

The multiscale analysis of \cite{GMT17,GMT20} integrates the
Grassmann field scale by scale, from the ultraviolet to the infrared,
and the infrared scales are those of the momenta near the two zeros
$p^\pm$ of the symbol $\mu$, where the free covariance $1/\mu$ is
singular. The first integration treats all the remaining momenta at
once. Its purpose, in our setting, is to turn the interaction of
\cref{prop:fermion-input}, whose \texttt{fat} coefficients may be
large, into an effective interaction on the slow modes which is
small in a norm suitable for the infrared induction. To this end,
choose a smooth cutoff $0\leq\chi_d\leq1$ supported in balls of radius
$2d$ around the two zeros of $\mu$, equal to one in radius-$d$ balls,
and of Gevrey order $2$ (that is, with derivatives of order $m$
bounded by $C^{m+1}(m!)^2$ times the appropriate power of $d$).
Put
\[
 \widehat g_0=(1-\chi_d)/\mu,\qquad
 \widehat g_1=\chi_d/\mu\,,
\]
so that $g_0+g_1$ is the free covariance $G$, and write the field
as $\psi+\xi$ with $\xi$ the fast field, of covariance $g_0$, and
$\psi$ the slow field, of covariance $g_1$. The fast multiplier
$\widehat g_0$ is smooth, since $1-\chi_d$ vanishes identically
where $\mu$ does, and its Gevrey regularity translates into
stretched exponential decay of the kernel $g_0$ on scale $d^{-1}$;
this is made quantitative in the proof of
\cref{prop:conditional-first-step} below. The exact identity for an
interior transformed black field, in the sense explained in the
next paragraph, is
\begin{equation}\label{eq:correct-fast-covariance}
       Kg_0=1-\chi_d,\qquad
       (Kg_0)(x,y)=\delta_{xy}-\check\chi_d(x-y)\,.
\end{equation}
Indeed, the symbol of $Kg_0$ is $\mu\cdot(1-\chi_d)/\mu=1-\chi_d$,
whose kernel is the identity minus the inverse Fourier transform
$\check\chi_d$ of the cutoff. We emphasize that this is not
$(1-\epsilon_d)\delta_{xy}$ for a small number $\epsilon_d$: the
second term is nonlocal: it is a smooth kernel spread over the scale
$d^{-1}$, whose $\ell^1$ norm is of order one, and it must be
carried along in the estimates.

We now explain how the determinant gain of
\cref{prop:fat_pattern} is realized inside the Grassmann integral.
For each \texttt{fat} monomial $\psi_S$ separately use the patches of
\cref{sec:fat-determinant}. Anchor the patch grid at the
lexicographically first vertex of its support and use fixed local
tie-breaking rules; these choices commute with lattice translations,
which preserves translation covariance of the resulting kernels.
Replace its black column vector by
$K_I^{\trps}\psi_b$, with the identity block on the complement.
(The transpose fixes the convention $K=K(b,w)$; the resulting interior
form at $w$ is $\sum_bK(b,w)\psi_b$.) Since the exterior product of
the $n$ linear forms $(A\psi)_1,\ldots,(A\psi)_n$ equals $\det A$
times the exterior product of $\psi_1,\ldots,\psi_n$, this gives
the algebraic identity
\[
 w(S)\psi_S=\frac{w(S)}{\det K_I}\,\widetilde\psi_S\,,
\]
where $\widetilde\psi_S$ is the monomial in the transformed forms.
The coefficient has thus been divided by $\det K_I$, which is at
least $e^{(P-\eta)|S|/2}$, and the price is that the fields of the
monomial are no longer the original ones; the point of the identity
\cref{eq:correct-fast-covariance} is that the covariance of an
interior transformed black form with a white field has multiplier
$1-\chi_d$, of modulus at most one, so that this price is bounded.
The interior contraction can equivalently be written $g_0K$ instead
of $Kg_0$ in \cref{eq:correct-fast-covariance}; their scalar Fourier
multipliers coincide. This is an identity of exterior products,
valid even when different monomials overlap. It is not one
simultaneous global change of variables, which would not be
available when the patches of different monomials overlap. Interior
black forms have multiplier $\mu$; the exceptional black forms, at
the noninterior vertices of the patches, are sums of at most three
neighboring fields, or unchanged fields in the identity block. For
any fixed small $\eta>0$, patches can be chosen so that
$|\det K_I|\geq e^{(P-\eta)|S|/2}$ and at most $\eta|S|$ black
forms are exceptional; this is \cref{lem:fat-coverage} together
with \cref{lem:entropy-patch}, exactly as in the proof of
\cref{prop:fat_pattern}.

The following lemma records the Gram representation of the fast
covariance between these linear forms. Its purpose is the
Gram--Hadamard inequality: a determinant whose entries are inner
products $\langle u_i,v_j\rangle$ is bounded by
$\prod_i\|u_i\|\prod_j\|v_j\|$, with no factorial, and this is the
standard way in which fermionic expansions avoid the combinatorial
divergence of their bosonic counterparts; see, e.g.,
\cite{GentileMP,GMT17}. The parameter $u$ below
trades the norm of the white vectors against that of the ordinary
black vectors; it will be chosen small, so that the ordinary white
and interior black occurrences, which make up all but a fraction
$\eta$ of a \texttt{fat} monomial, cost essentially nothing.

\begin{lemma}[Gram bound after the cutoff]\label{lem:correct-gram}
For every $u>0$, the covariance $g_0$ between the above linear forms
has a Gram representation in which white vectors have norm at most
$\sqrt{1+u^2}$, interior black vectors have norm at most $1$, and all
other black vectors have norm at most $C/u$. These bounds also apply
after any row or column deletion and after tensoring with unit vectors
attached to different monomials.
\end{lemma}
\begin{proof}
Write the ordinary covariance as $g_0(x,y)=\langle a_y,b_x\rangle$
with $\|a_y\|,\|b_x\|\leq\sqrt{C_0}$, by splitting its Fourier
multiplier $\widehat g_0$ into two factors of equal modulus
$|\widehat g_0|^{1/2}$ and multiplying them by the characters of
$y$ and $x$, as in the proof of \cref{prop:fat_pattern}; here
$C_0$ bounds $\int|\widehat g_0|\,\d k/(2\pi)^2$, which is finite
since $\widehat g_0$ is bounded. Work in
$\ell^2(\Z^2)\oplus L^2(\T^2)$. Represent a white field at $y$ by
$(\delta_y,u a_y/\sqrt{C_0})$ and an ordinary black field at $x$ by
$(0,\sqrt{C_0}b_x/u)$. The inner product of these two vectors is
$\langle a_y,b_x\rangle=g_0(x,y)$, as required, and their norms are
at most $\sqrt{1+u^2}$ and $C_0/u$ respectively. For an interior
black form use $((1-\chi_d)\delta_x,0)$ in the first summand,
interpreting $\chi_d$ as a Fourier multiplier on $\ell^2(\Z^2)$. Its
norm is at most one, because the multiplier $1-\chi_d$ has modulus
at most one, and its inner product with the white vector at $y$ is
$\langle(1-\chi_d)\delta_x,\delta_y\rangle=(1-\chi_d)(x,y)$, which
is exactly the covariance \cref{eq:correct-fast-covariance} of the
interior form with the white field. Exceptional black forms are
linear combinations of at most three ordinary black vectors with
coefficients of modulus at most one, so have norm at most
$3C_0/u=C/u$. This proves the representation. Deleting vectors
preserves it, and tensoring with unit vectors does not change their
norms; both operations therefore preserve the Gram bound.
\end{proof}

The first integration is organized by the cumulant expansion: the
logarithm of the fast Gaussian expectation of $e^{V}$ is the sum
over $n\geq1$ of the connected expectations
$\cE^T_{g_0}(V;\ldots;V)/n!$. We use the Brydges--Kennedy Taylor forest formula in the form of
\cite{AbdesselamRivasseau95}; the fermionic version is also recalled
in \cite[Lem.~3]{GMT17}. We record the formula and explain its
application to vertices with external fields. Its remaining
determinant admits the Gram bound of \cref{lem:correct-gram}.

\begin{theorem}[Connected Gaussian expectations]\label{thm:truncated_expectation_bound}
Let $X_i$ be Grassmann monomials, allowing external Grassmann
coefficients so that each complete vertex is even. Suppose the
integrated fields have a total of $q$ black and $q$ white occurrences.
Each individual integrated submonomial may be unbalanced or odd. If
external fields are present, both sides below denote the coefficient
of their prescribed product. For $n\geq2$,
\[
 \cE_g^T(X_1;\ldots;X_n)
 =\sum_T\pm(T)\prod_{e\in T}g_e
          \int_{[0,1]^{n-1}}\det G_T(s)\,\d s\,.
\]
The sum is over trees on the $n$ monomials with distinct black/white
occurrences marking their edges. The remaining determinant has size
$q-n+1$, after deleting the occurrences used by the tree. Its block
$(i,j)$ is the remaining covariance multiplied by
$t_{ij}=\min_{e\in[i,j]_T}s_e$, and $t_{ii}=1$.
An empty input makes the cumulant zero when $n\geq2$; for $n=1$,
$\cE_g^T(1)=1$.
\end{theorem}
\begin{proof}
Multiply every covariance between vertices $i,j$ by a variable
$t_{ij}$ and apply the Taylor forest formula
\cite{AbdesselamRivasseau95} to the resulting Gaussian moment,
viewed as a polynomial with coefficients in the external Grassmann
algebra. At a forest interpolation point, covariances between
components vanish; since the complete vertices are even, the moment
factorizes over those components. The moment--cumulant relation
therefore selects exactly the connected forests, namely trees.
Each tree derivative selects a black and a white occurrence and
contracts them. Repeated occurrences give zero, and Wick's rule
for the remaining fields gives the determinant of size $q-n+1$.
Extract the prescribed external coefficient only after this
operation; this permits odd or unbalanced integrated submonomials.

The forest formula gives the stated minimum-path interpolation.
Its positivity property \cite{AbdesselamRivasseau95} makes
$(t_{ij})$ a Gram matrix of unit vectors. Tensoring the covariance
vectors with these unit vectors preserves their norms, so
\cref{lem:correct-gram} applies. An input without integrated fields
has no possible incident tree edge, giving the last assertions.
\end{proof}

We next specify the norm in which the effective interaction after
the first integration will be measured. An effective kernel
$w^{(1)}(S)$ is a function of a finite set $S$ of field occurrences,
and we need to express both that it is small and that it decays when
$S$ is spread out. For kernels after integration, define the tree
length
\[
 \ell(S)=\min_{T\text{ spanning }S}\sum_{xy\in T}\sqrt{1+|x-y|}\,,
\]
the minimum being over trees on the points of $S$; the square root
reflects the stretched exponential decay of the fast covariance.
We require an \emph{anchored} norm, in which $S$ is summed subject to
containing a fixed point $x$, with a \emph{positive} weight
$e^{a\ell(S)}$, together with exponential field-degree weights
$c^{|S|}$. Each of these features is forced. An unanchored sum
diverges simply by translation invariance, since every translate of
a support contributes the same amount. A negative spatial weight
makes summability weaker than at $a=0$ and cannot express the needed
decay. The field-degree weight is what allows a monomial of degree
$k$ to be estimated by $\epsilon^{k}$, as the infrared induction
requires; see \cref{prop:infrared-input}. Replacing $\ell(S)$ by
$\sqrt{\operatorname{Animal}(S)}$, the square root of the size of
the smallest lattice animal containing $S$, is also insufficient, as
explained in the remark preceding \cite[\S8.6]{LaslierLubetzky24}:
the latter weight can be much smaller than the sum of square roots
of the connecting lengths, and therefore does not supply the
tree-length control needed in the subsequent integration.

The convergence of the cumulant expansion is then a consequence of
the following elementary lemma on rooted trees, which packages the
standard Kotecký--Preiss type argument~\cite{KoteckyPreiss86} in the
form in which it will be applied: the nodes of the tree are the
inputs, with activities $b(S)$, and its edges are the contraction
lines, with the propagator bound $k(x,y)$.

\begin{lemma}[Rooted tree summation]\label{lem:rooted-tree}
Let $b(S)\geq0$ be activities of finite supports, and let $k(x,y)\geq0$
with $\sup_x\sum_yk(x,y)=G<\infty$. Suppose
\[
 \sup_x\sum_{S\ni x}b(S)e^{u|S|}\leq\delta,\qquad G\delta\leq u\,.
\]
The sum of rooted marked trees with root support $S$, node activities
$b$, and edge factors $k$, with the exponential-series symmetry
factors for unordered children, is at most $b(S)e^{u|S|}$.
\end{lemma}
\begin{proof}
Let $T_m(S)$ be the sum restricted to trees of height at most $m$.
A tree of height at most $m+1$ with root $S$ consists of the root
together with an unordered collection of subtrees of height at most
$m$, each attached to the root by an edge from a point $x\in S$ to a
point $y$ of the child's support. The sum over unordered collections
of children, with the symmetry factors $1/j!$ for $j$ children, is
the exponential of the sum over a single child, so
\[
 T_{m+1}(S)\leq b(S)\exp\bigg(
      \sum_{x\in S}\sum_{S'}\sum_{y\in S'}k(x,y)T_m(S')\bigg)\,.
\]
Assume inductively that $T_m(S')\leq b(S')e^{u|S'|}$ for all $S'$.
Then the inner sum at each $x\in S$ is at most
$\sum_y k(x,y)\sum_{S'\ni y}b(S')e^{u|S'|}\leq G\delta\leq u$, and
the exponent is at most $u|S|$, which proves
$T_{m+1}(S)\leq b(S)e^{u|S|}$; the induction starts with $T_0=b$.
Monotone convergence in $m$ completes the proof. We note where the
factorials have gone. Rooting a cumulant series, that is, fixing one
of its $n$ inputs to contain the anchor, changes $1/n!$ into at most
$1/(n-1)!$, and $1/(n-1)!$ is precisely the symmetry factor of a
rooted labelled tree on $n$ nodes that the above recursion encodes
through its exponential generating function. Thus no unaccounted
factor $n^n$ remains.
\end{proof}

\begin{proposition}[Conditional first-step summability]\label{prop:conditional-first-step}
Fix a liquid slope and assume the first three bounds in
\cref{prop:fermion-input}. Then for sufficiently small $\epsilon_\beta$
and sufficiently large $R_\beta$ the first integration is absolutely
convergent, uniformly per unit volume. Its nonconstant effective
kernels on the remaining Fourier modes have arbitrarily small
anchored norm
\begin{equation}\label{eq:correct-rg-norm}
 \sup_x\sum_{S\ni x}|w^{(1)}(S)|c^{|S|}e^{a\ell(S)}\,,
\end{equation}
for any prescribed fixed $c$, for some $a>0$ depending on the initial
cutoff. The constant coefficient is bounded per unit volume. If the
same bounds for interaction-derived coefficients hold on a source
polydisc, source analyticity holds on every smaller sufficiently small
polydisc; pure edge-source terms contribute an additional $O(\rho)$
to the norm bound.
\end{proposition}
\begin{proof}
We apply the cumulant expansion to the fast field, with the slow
field as an external field, and estimate its terms by
\cref{thm:truncated_expectation_bound,lem:correct-gram}, after the
exterior-product transformation of the \texttt{fat} monomials. The
proof is organized as follows. We first determine what an external
field occurrence costs (\cref{st:first-external}), then use the
balance between black and white occurrences to distribute these
costs evenly (\cref{st:first-balance}), and derive the envelope of
a single node (\cref{st:first-node}). The order in which the
parameters $u,d,\eta$ and $\beta$ are chosen is then explained
(\cref{st:first-parameters}), after which the trees are summed with
their spatial weight by \cref{lem:rooted-tree}
(\cref{st:first-trees}). Analyticity in the sources and the torus
case are addressed last (\cref{st:first-sources}).
\begin{enumerate}[label=\textbf{Step~\arabic*.}, ref=\arabic*, wide=0pt, itemsep=1ex, listparindent=\parindent, parsep=0pt]
\item\label[step]{st:first-external}
\emph{(The cost of an external field.)}
We specify how to pay for the external fields, which are essential
in this estimate. Let $\zeta_d$ be a Gevrey cutoff equal to one on the
support of $\chi_d$, supported in slightly larger balls. On the slow
field space an ordinary external field can be written as
$\zeta_d\psi$, since the slow field has all its Fourier support
where $\chi_d\ne0$. Its convolution kernel has $\ell^1$ norm at most
$C$; an interior black external form is $K\zeta_d\psi$, with kernel
norm at most $Cd$, which is the gain from the transformation of the
\texttt{fat} monomials on the external legs. Both assertions follow
by writing $k=p^\pm+dq$: the smooth rescaled symbols and their
Gevrey bounds are uniform, and $\mu(p^\pm+dq)=O(d)$ because $\mu$
vanishes at $p^\pm$ and is smooth. More explicitly, for either symbol $f_d$ with amplitude $d^j$
($j=0$ for $\zeta_d$ and $j=1$ for $\mu\zeta_d$), its support has
area $O(d^2)$ and its derivatives satisfy
$\|\partial^m f_d\|_\infty\leq C^{m+1}(m!)^2d^{j-m}$.
Integrating by parts $m$ times in a coordinate with
$|x_i|\geq|x|/\sqrt2$ yields
\[
 |\check f_d(x)|\leq
 C d^{j+2}\inf_{m\geq0}
          \frac{C^m(m!)^2}{(d|x|)^m}
 \leq C d^{j+2}e^{-c_0\sqrt{d|x|}}\,.
\]
For the last inequality take $m$ to be the integer part of
$c_1\sqrt{d|x|}$ with $c_1$ small; the estimate for bounded $d|x|$
uses $m=0$. Summation over lattice shells, or the substitution
$t=\sqrt{dr}$ in the radial integral, gives weighted $\ell^1$ norm
at most $Cd^j$ with weight $e^{a\sqrt{|x|}}$ for
$a<c_0\sqrt d/2$. The constants below include these weighted norms.

\item\label[step]{st:first-balance}
\emph{(Balancing the black and white costs.)}
Every surviving monomial has equal total black and white degree:
the inputs are balanced, since every pattern monomial is a product
of factors $\psi_b\psi_w$, and the fast contractions pair a black
occurrence with a white one, so the external occurrences left over
are balanced as well. In estimates we
may therefore multiply the cost of every external white occurrence
by $\sqrt d$ and of every external black occurrence by $1/\sqrt d$
without changing a term, since the two factors cancel over the
whole monomial. Both ordinary white and interior black external
occurrences now cost at most $C\sqrt d$: the former had cost $C$ and
the latter $Cd$. Exceptional black occurrences cost at most
$C/\sqrt d$, and they are at most a fraction $\eta$ of all the
occurrences of a \texttt{fat} node.

\item\label[step]{st:first-node}
\emph{(The envelope of a node.)}
Apply the connected expectation formula to the fast fields. Its
remaining determinant, for every choice of external fields and tree
endpoints, has the bound of \cref{lem:correct-gram}. Summing over
whether a non-tree occurrence is external or in the determinant adds
the two corresponding bounds. Tree endpoints may be overcounted by
the same positive bounds. Thus a \texttt{thin} node of size $s$ has a bound
$(C_{u,d,c}\epsilon_\beta)^s$, while a \texttt{fat} node has a bound
\begin{equation}\label{eq:node-envelope}
 \exp\big(\big[-(P-\eta)/2+\epsilon_\beta
       +O(u^2+\sqrt d)+\eta\log C_{u,d,c}\big]s\big)\,.
\end{equation}
The factor $c$ for the external degree has been included in the
external costs. Here at most $\eta s$ black occurrences are
exceptional. The determinant denominator supplies the term
$-(P-\eta)/2$.

To make the small loss in this estimate explicit, fix a contraction
tree and its field endpoints. A node has $s/2$ white and $s/2$
black occurrences before choosing the external fields. For each
occurrence not used by the tree, sum its external-field cost and
its Gram-vector norm. For an ordinary white occurrence this costs
at most $\sqrt{1+u^2}+C_c\sqrt d$, and for an interior black
occurrence at most $1+C_c\sqrt d$. An exceptional black occurrence
instead costs at most $C_{u,d,c}$. Increasing this last constant
above one, all the displayed factors are at least one. We may
thus insert the same factors for the tree endpoints, whose
contractions are estimated separately by the tree kernel below.
The sum over external-field choices is thus bounded by
\[
 (\sqrt{1+u^2}+C_c\sqrt d)^{s/2}
 (1+C_c\sqrt d)^{s/2}C_{u,d,c}^{\eta s}.
\]
Its logarithm is $O_c(u^2+\sqrt d)s+\eta s\log C_{u,d,c}$.
Thus the large cost of an exceptional occurrence is paid only
on a small fraction of the support; the determinant gain is paid
on the whole support. This also explains the order of parameter
choices in the next step.

\item\label[step]{st:first-parameters}
\emph{(Choice of the parameters.)}
Choose $u>0$ small, next $d>0$ small, and then $\eta>0$ small; finally
choose $\beta$ so large that the \texttt{fat} sets allow these patches and
$\epsilon_\beta$ is small. The order matters: the constant
$C_{u,d,c}$ in \cref{eq:node-envelope} depends on $u$ and $d$, so
$\eta$ must be chosen after them for the term $\eta\log C_{u,d,c}$
to be small, and the \texttt{fat} parameters $\epsilon_\beta,R_\beta/2$,
which must be small and large enough for \cref{lem:fat-coverage} at
the chosen $\eta$ and patch size, are only available for large
$\beta$. The exponent in \cref{eq:node-envelope} is then at most
$-Ps/4$. \texttt{Thin} supports have at most $C^s$ rooted choices by
connected-set enumeration; \texttt{fat} supports have at most
$e^{\epsilon_\beta s}$ rooted labelled choices and size at least
$2B_\beta$. Consequently the sum of these node envelopes, anchored
at a site and multiplied by $e^{v s}$ for a fixed sufficiently small
$v>0$, tends to zero with $\beta$: the \texttt{thin} sum is a
geometric series with ratio $C\epsilon_\beta e^v$, and the
\texttt{fat} sum is at most
$\sum_{s\geq2B_\beta}e^{-(P/4-\epsilon_\beta-v)s}\leq Ce^{-cB_\beta}$.

\item\label[step]{st:first-trees}
\emph{(Summation of the trees and the spatial weight.)}
All transformed fast propagators are finite linear combinations of
$g_0$ and therefore are bounded by a nonnegative translation-invariant
kernel $k$ with finite weighted $\ell^1$ norm $G$. Fix the parameters
above before taking $\beta$ large; $G$ need not be small. The anchored
node sum can be made smaller than $v/G$, so
\cref{lem:rooted-tree} sums the cumulant trees absolutely.
The same proof pays the spatial weight in \cref{eq:correct-rg-norm}:
connect every input support by nearest-neighbor edges, add the
contraction tree edges and the external convolution legs, and double
this tree to traverse it. The triangle inequality for
$\sqrt{1+|x-y|}$ bounds the minimal tree length of its external
vertices by twice the total length of this traversal tree. Choosing
a smaller $a$ absorbs the contraction and external-leg lengths in
the weighted kernel norms, and the internal nearest-neighbor lengths
in the factor $e^{v s}$. Root at an external occurrence (or at an
input site for the constant coefficient); its convolution is summed
with the same weighted norm. This proves the claimed anchored bound.

\item\label[step]{st:first-sources}
\emph{(Analyticity, sources and the torus.)}
The cumulant expansion is initially an identity of formal series in
the activities. The absolutely convergent tree bounds on a complex
polydisc extend the identity to the stated activities and allow
termwise differentiation. The empty external set gives a constant
in the logarithm, with one free translation, hence its stated
per-volume bound. Pure edge-source terms have bounded support and
activity $O(\rho)$, so choosing $\rho$ small preserves the
rooted-tree criterion. Source estimates then follow by the same
argument and \cref{lem:source-envelope}. For finite tori the fast
covariance is a periodization of a smooth Fourier multiplier and has
the same bounds uniformly in the twists; modes at a zero of $K$ are
kept in the slow polynomial integral and are never inverted in this
step. This completes the proof.
\end{enumerate}
\end{proof}

\begin{theorem}[First integration for bulk coefficient families]\label{thm:bulk-first-step}
Fix a liquid slope and an external field-degree weight $c<\infty$.
Use the infinite-volume free covariance and restrict the pattern
supports of \cref{def:rg-admissible} to a finite window $\Omega$.
For all sufficiently large $\beta$ there is an absolutely convergent
first integration
\begin{equation}\label{eq:bulk-first-step}
 \cE_{g_0}
 e^{J_A(\psi+\xi)+V_{\beta,A}(\psi+\xi)}
       =e^{C_{\beta,\Omega}(A)+V^{(1)}_{\beta,\Omega,A}(\psi)}\,,
\end{equation}
where $\xi$ is the fast field, $\psi$ is the remaining field, and
$V^{(1)}$ has no constant term. There are $a,\rho_0,c_1,C>0$,
independent of $\Omega$ and of the finite source set $F$, such that,
with $\delta_\beta=Ce^{-c_1\beta^{2/5}}$,
\begin{align}
 \sup_x\sum_{S\ni x}|w^{(1)}_{\beta,\Omega,0}(S)|
                 c^{|S|}e^{a\ell(S)}&\leq\delta_\beta,
                                      \label{eq:bulk-small-kernel}\\
 |C_{\beta,\Omega}(0)|&\leq\delta_\beta|\Omega|\,.
\end{align}
On $\|A\|_\infty\leq\rho\leq\rho_0$, the first bound is replaced by
$C(\delta_\beta+\rho)$ and the second by
$C\delta_\beta|\Omega|+C\rho|F|$. All coefficients are analytic
on the open source polydisc. As $\Omega$ increases to the plane,
the nonconstant kernels converge in the anchored norm at each fixed
root; at $A=0$ their limit is translation covariant.
The exact generating function keeps the vertex gauge identity.
A bound keeping the source locations in the spatial weight also
holds. Write
\[
 C_{\beta,\Omega}(A)+V^{(1)}_{\beta,\Omega,A}(\psi)
       =\sum_{S,\fm}w^{(1)}(S;\fm)
                              \psi_S A^{\fm},
 \qquad
 \scS(S,\fm)
       =S\cup\bigcup_{e:\fm_e>0}\{b(e),w(e)\}\,,
\]
where $\fm$ is a nonnegative edge multi-index. With possibly
smaller $a,\rho_0$, for $0\leq\rho\leq\rho_0$,
\begin{equation}\label{eq:source-marked-norm}
 \sup_x\sum_{\substack{S,\fm:\,x\in\scS(S,\fm)}}
 |w^{(1)}(S;\fm)|c^{|S|}\rho^{|\fm|}
                e^{a\ell(\scS(S,\fm))}
       \leq C(\delta_\beta+\rho)\,.
\end{equation}
Thus the derivatives are controlled with their source locations,
including terms with no remaining fermion field. The coefficient
of $A^{\fm}$ includes $1/\fm!$.
\end{theorem}
\begin{proof}
The theorem is \cref{prop:conditional-first-step} with the rate of
smallness made explicit, together with the source bounds and the
window limit. We prove \cref{eq:bulk-small-kernel} and the constant
bound first, then the bounds on the source polydisc, then
\cref{eq:source-marked-norm}, and finally the convergence as the
window grows and the identity \cref{eq:bulk-first-step} itself.
\Cref{it:rg-envelopes,prop:fermion-input}
provide the three coefficient hypotheses of
\cref{prop:conditional-first-step}. In its proof fix
$u,d,\eta,a$ before increasing $\beta$. After these choices, \texttt{thin}
node weights, including their rooted pattern count and the fixed
degree and spatial factors, sum to at most
$Ce^{-c\beta^{2/5}}$. \texttt{Fat} nodes have an exponential envelope
$e^{-p|S|}$ with fixed $p>0$, rooted count $e^{o(1)|S|}$, and
$|S|\geq2B_\beta$. Their rooted sum is at most $Ce^{-cB_\beta}$.
Both estimates remain true after multiplication by $e^{v|S|}$ for
a fixed sufficiently small $v>0$. The tree lemma therefore bounds
the rooted cumulant sum by $\delta_\beta$, after changing constants.
This proves \cref{eq:bulk-small-kernel}. Rooting a constant term
at an input vertex gives its bound by $\delta_\beta|\Omega|$.

Uniformity in the source set does not follow merely by absorbing
$e^{\rho|F|}$ into a coefficient: that factor grows with the number
of source edges, whereas the asserted bounds do not. Instead use
\cref{it:rg-sources} in the form $e^{\rho|S|/2}$: the source factor
of an interaction pattern is at most $e^{\rho|P|}=e^{\rho|S|/2}$,
whatever the source set, and it is paid node by node.
Choose $\rho_0>0$ so small that this changes the
\texttt{fat}-node exponent by less than half of its fixed negative margin.
It also preserves the \texttt{thin} estimate for large $\beta$. The same
interaction-node sums are then bounded by $C\delta_\beta$, uniformly
in $F$. The additional nodes from $J_A$ have degree two, bounded
diameter and rooted sum at most $C\rho$. Decrease $\rho_0$ once
more so that the rooted-tree criterion holds for their combined
activity. This proves the asserted source bound. For the constant
coefficient, root each tree at any one of its inputs: roots from
the interaction contribute at most $C\delta_\beta|\Omega|$, and
roots from $J_A$ at most $C\rho|F|$. This overcounts and hence
is valid for the absolute sum. Locally uniform convergence of
analytic cumulants proves analyticity, and Cauchy's formula gives
the derivative bounds without a factor depending on the volume.

For \cref{eq:source-marked-norm}, perform the same estimates in the
algebra of source power series with the absolute Taylor norm at
radius $\rho$. At an interaction input this replaces
$\prod_{e\in P}e^{A(e)}$ by
\[
 \prod_{e\in P}\sum_{m\geq0}\frac{\rho^m}{m!}
                         =e^{\rho|P|}\,,
\]
and at a free-source input it replaces $e^{A(e)}-1$ by
$e^\rho-1$. These are exactly the envelopes just used. Every source
mark lies on an input edge. Include its endpoints in the connected
nearest-neighbor tree of that input before adding contraction edges
and external convolution legs. The spatial-weight argument in
\cref{prop:conditional-first-step} then bounds the tree
length of $\scS$, with the same constants. A root $x$ is
either an external occurrence or a source endpoint in an input;
both cases have the preceding rooted bounds. The absolute Taylor
norm is submultiplicative, so all products and sums of source
coefficients are bounded before any cancellation. The tree lemma
now proves \cref{eq:source-marked-norm}, including the terms
$S=\varnothing$, $\fm\ne0$. Its excluded term
$S=\varnothing$, $\fm=0$ is precisely the constant already
bounded per unit volume.

For increasing windows, every fixed pattern eventually occurs with its
unrestricted coefficient, by \cref{it:rg-locality}. The bounds just proved give a summable envelope for
the rooted labelled trees, including their external convolution
legs. Dominated convergence therefore proves convergence in the
weighted sum anchored at each fixed root. Translation covariance
at zero source follows from that of the unrestricted coefficients
and of the two covariances. This is a local anchored limit; no
uniform convergence over roots escaping with the window is claimed.

Finally the finite-window fermion identity and Gaussian convolution
give \cref{eq:bulk-first-step} on the slow field space, with the
constant chosen by analytic continuation from zero activities.
The absolutely convergent connected series fixes this branch.
Integrating the remaining fields recovers exactly
the generating function in \cref{it:rg-sources}; hence the vertex gauge identity
of the original generating function is unchanged. A momentum
cutoff need not make individual effective kernels covariant under
a local change of variables; the assertion concerns the exact
generating function and its resulting identities.
\end{proof}

Thus the bulk first integration has a small interaction in
the positive spatial norm \cref{eq:correct-rg-norm}, including
external fields and a source polydisc independent of the volume.
We next establish the corresponding result for the canonical
completion on the torus. Neither result supplies an interacting
insertion bound after the remaining modes have been integrated.

\begin{lemma}[Entropy patches on the torus]\label{lem:torus-entropy-patches}
Fix liquid weights normalized by $\max_r\omega_r=1$ and $\eta>0$.
If $\epsilon$ is sufficiently small, $r$ sufficiently large and
$n$ sufficiently large, the support $S$ of any
$(\epsilon,r)$-\texttt{fat} torus pattern admits the exterior-product
transformation of \cref{sec:first-integration} with
\[
 |\det K_I|\geq e^{(P-\eta)|S|/2},\qquad
 \#\{\text{exceptional black forms}\}\leq\eta|S|\,.
\]
The assertion includes noncontractible patterns and the entire torus.
It holds with the same bounds for arbitrary unitary boundary twists,
or unit phases on the two edge orientations.
All parameter thresholds are uniform on compact sets of liquid
slopes.
\end{lemma}
\begin{proof}
Choose the planar patch size $L$ and its selected boundary as in
\cref{lem:entropy-patch}, with errors much smaller than $\eta$.
The proof of \cref{lem:fat-coverage} applies to the periodic
family of ball centers: for $n$ large enough the individual balls
are local, although their union may wrap. There is just one
additional loss. A uniformly translated fundamental parallelogram
has a $CL$-neighborhood of its boundary containing any given vertex with probability
at most $CL/n$, since its perimeter is $O(n)$ and the torus area
is at least $cn^2$. Some cut therefore loses at most
$(CL/n)|S|$ vertices. Retain grid squares contained in the ball
union and the cut-open domain, and place $D_L$ in each. The
fraction of vertices not interior to a selected patch is at most
\[
                 C(\epsilon+L/R+L/n+L^{-1})\,.
\]
Choose first $L$, then $\epsilon,R$, and finally $n$, to make
this loss small enough. The patch determinant bound and the
identity block on unused balanced vertices give the assertion.
The uniform threshold in \cref{lem:entropy-patch} and the
uniform torus shape bounds make these choices uniform on compact
slope sets. Finally, every flat unitary twist is a vertex gauge
on each selected square: it changes only the determinant phase
and the phases of the interior forms. Thus the same bounds hold
without a global gauge or a planar lift of the entire support.
\end{proof}

\begin{theorem}[First integration uniformly in the torus twist]
\label{thm:torus-first-step}
Use the coefficients of \cref{def:rg-admissible}, and fix positive liquid
weights with $\max_r\omega_r=1$ and a field-degree weight $c<\infty$.
For all sufficiently
large $\beta$ and all $n\geq n_0(\beta)$, the first integration
of \cref{eq:rg-berezin} is absolutely convergent,
uniformly for $|y_2|=|y_3|=1$.
More precisely, write the perturbation of its quadratic part as
$J_A+V_{\beta,\Gamma,A}$, with $J_A$ defined as in
\cref{it:rg-sources}, with the torus coefficients in $V$. On the
slow field space, \cref{eq:bulk-first-step} holds with $g_{0,n}$ and
$\Gamma$ in place of $g_0$ and $\Omega$.
The scalar $C$ here is generated by the interaction and sources;
the free Gaussian normalization is retained separately.
The effective nonconstant kernels satisfy
\cref{eq:bulk-small-kernel} with torus distance in $\ell$ and
$\delta_\beta=Ce^{-c_1\beta^{2/5}}$. The constant is bounded by
$C\delta_\beta N_0$ at zero source and by
$C\delta_\beta N_0+C\rho|F|$ on a fixed source polydisc.
The nonconstant and source-marked norms on that polydisc are
bounded by $C(\delta_\beta+\rho)$, with constants independent of
$n$, the twist and the number of source edges.

For fixed $\beta$, in a local gauge removing the twist, these
kernels converge at each fixed root to the bulk kernels of
\cref{thm:bulk-first-step}, uniformly in the twists.
Convergence holds in every smaller positive spatial norm and,
for a fixed finite source set, on every smaller source polydisc.
Modes at which the quadratic matrix vanishes remain in an
unnormalized slow integral; they are not divided out.
\end{theorem}
\begin{proof}
The proof repeats that of \cref{thm:bulk-first-step} on the torus,
and we concentrate on the three points which differ: the twisted
fast covariance and its kernels, the torus version of the
Gram bound and of the entropy patches, and the algebraic meaning of
the Gaussian operation at a singular twist, where the quadratic
form has a zero mode. We then prove the convergence to the bulk
kernels.

Multiplying two edge orientations by phases shifts the Fourier
argument of $\mu$ and is equivalent, after a vertex gauge, to a
unitary character on the period lattice. Use this shifted cutoff,
or equivalently the usual cutoff on the twisted momentum grid.
The fast multiplier $(1-\chi_d)/\mu$ is smooth everywhere,
including the two zeros, where it is defined to be zero.
For a boundary character $\vartheta$, its kernel is the twisted
periodization
\[
 g_{0,n}^{\vartheta}(x,y)=
     \sum_{\lambda\in\Upsilon\Gamma}
          e^{i\vartheta(\lambda)}g_0(x-y+\lambda)\,,
\]
in lattice coordinates and a fixed local gauge. The same formula
holds for the external-leg kernels $\zeta_d$ and $K\zeta_d$.
The triangle inequality and the inequality
$\dist_{\T_\mathrm{d}}(x,y)\leq|x-y+\lambda|$ bound their
weighted torus $\ell^1$ norms by the planar ones, uniformly in
$\vartheta$. On bounded sets they converge to their planar kernels
uniformly in $\vartheta$, since all nonzero image terms have
distance tending to infinity.

\Cref{lem:correct-gram} also holds on the torus. Replace its
Fourier integrals by averages over the twisted momentum grid.
The ordinary fast covariance has a uniformly bounded Gram norm,
since its smooth multiplier is bounded for the fixed cutoff.
The interior vector has norm at most one because its multiplier
is $1-\chi_d$ with modulus at most one. Use
\cref{lem:torus-entropy-patches} for every \texttt{fat} input.
All estimates of individual transformed nodes, external legs,
and contraction edges now hold with the same fixed constants.

\Cref{it:rg-envelopes} supplies the two coefficient
envelopes. Its extra \texttt{fat}-pattern entropy $C/n$ is
absorbed in the fixed negative margin of
\cref{eq:node-envelope} by increasing $n_0(\beta)$.
The anchored \texttt{thin}-node sum is at most
$Ce^{-c\beta^{2/5}}$, and the \texttt{fat}-node sum at most
$Ce^{-cB_\beta}$. The torus version of the rooted-tree lemma
therefore proves the asserted convergence and the zero-source
norm bounds. The proof of \cref{eq:source-marked-norm} applies
unchanged: all marks lie on input edges, the absolute Taylor
factor is $e^{\rho|P|}$, and the extra free-source nodes cost
$O(\rho)$. Choosing the same fixed small $\rho_0$ proves
analyticity and the bounds including source locations. Rooting
constant terms at interaction or source inputs gives the stated
separate volume and source-set factors.

For clarity, the Gaussian operation at a singular twist is
algebraic. First integrate modes with $\chi_d=0$, whose quadratic
eigenvalues are nonzero. On a mode with $0<\chi_d<1$, both
covariances in the split are nonsingular, so ordinary Gaussian
convolution applies. A mode with $\chi_d=1$ is wholly retained.
In particular, every zero mode is retained with its original
quadratic coefficient zero. Equivalently, the remaining normalized
expectation on a mode with $\chi_d>0$ can be combined with the
original Gaussian normalization as
\[
 \mu\,\cE_{\chi_d/\mu}(f)
      =\chi_d\int\d\psi\,
                    e^{(\mu/\chi_d)\psi_b\psi_w}f(\psi)\,.
\]
The right-hand side is defined also at $\mu=0$. Together with
the scalar factors from the other modes, it gives exactly the
original unnormalized polynomial. The first integration estimates
use only the fast expectation and never divide by this slow
integral or by a twisted partition function.

Finally consider a fixed root, with any fixed finite source set.
For every bounded collection of planar input patterns, the
canonical coefficients equal their bulk values once $n$ is
large enough, by \cref{it:rg-locality}. The fast
and external-leg kernels converge uniformly in the twist on
every bounded collection of positions. Thus each fixed finite
piece of the connected expansion converges to the corresponding
bulk piece. To justify passing to the whole series, keep a
strictly positive reserve in the node-size exponents and in the
weighted propagator bounds. Terms with an input of size greater
than $M$, a contraction or external-leg displacement greater
than $M$, or more than $M$ nodes then have a rooted absolute
tail tending to zero with $M$, uniformly in $n$ and the twist.
This follows directly from the same tree recursion, using the
reserved size, spatial and activity factors, respectively.
Nodes which cannot be lifted into a bounded neighborhood fall
into these tails; in particular a coefficient using the non-small
rule has size at least $c_*n$. The remaining terms have finitely
many positions and degrees after a common local lift. Dominated
convergence proves the anchored convergence, using a smaller
spatial exponent and a smaller Taylor radius to retain the
reserves. A monomial's exterior-product transformation is an
identity, so differing choices of auxiliary entropy patches
do not affect these coefficient limits.
\end{proof}

\Cref{thm:torus-first-step} closes the finite-volume
coefficient and first-integration part of the abstract argument. It does not estimate the remaining infrared integral,
its fixed-count Fourier coefficient, or the interacting witness
probabilities. Those issues require the subsequent flow and
the tuning of the auxiliary weights to the interacting slope.

The infrared induction of \cite[\S6]{GMT20} starts from an
interaction which is small in every field degree, with a quadratic
part whose symbol is a small smooth perturbation of $\mu$, and with
two Fermi points near which the slow modes live. The next
proposition extracts these data from \cref{thm:bulk-first-step}. We
stress that its zeros are the zeros of the dispersion after one
integration; they will be corrected at every subsequent scale.

\begin{proposition}[Small initial infrared data]\label{prop:infrared-input}
Fix the liquid weights. For every $\epsilon>0$, all sufficiently
large $\beta$ give a bulk first-step interaction satisfying, for
some $a_\epsilon>0$ and every even $k\geq2$,
\begin{equation}\label{eq:degree-small-input}
 \sup_x\sum_{\substack{S\ni x\\|S|=k}}
        |w^{(1)}_{\beta,0}(S)|e^{a_\epsilon\ell(S)}
                         \leq\epsilon^{k+1}\,.
\end{equation}
The cutoff and $a_\epsilon$ are fixed before increasing $\beta$.
At any such fixed cutoff, let $\Sigma_\beta(k)$ be the Fourier
symbol of the two-field part. It is smooth and satisfies
\begin{equation}\label{eq:self-energy-derivatives}
 \|\partial^m\Sigma_\beta\|_\infty
       \leq C_a^{|m|+1}(|m|!)^2\delta_\beta,
 \qquad
 \Sigma_\beta(-k)=\overline{\Sigma_\beta(k)}\,.
\end{equation}
Here $\delta_\beta\to0$ is the small norm at that cutoff.
The one-step dispersion $\mu+\Sigma_\beta$ has exactly two simple
zeros $p_\beta^{(1),\pm}$, with
\[
 p_\beta^{(1),-}=-p_\beta^{(1),+},\qquad
 |p_\beta^{(1),\pm}-p^\pm|\leq C\delta_\beta\,.
\]
They lie in the region where $\chi_d=1$. The quadratic Taylor
remainder at either zero is bounded by $C\delta_\beta|q|^2$
for the self-energy, and the quartic kernel differs from its
localization at $p_\beta^{(1),\pm}$ by at most
$C\delta_\beta\sum_i|q_i|$.
\end{proposition}
\begin{proof}
Choose the degree weight $c=\epsilon^{-1}$ in
\cref{thm:bulk-first-step}, decreasing $\epsilon$ below one
if needed. Choose $\beta$ after its cutoff and constants so that
its norm bound is at most $\epsilon$. Keeping only degree $k$
in the anchored sum, the factor $c^{|S|}=\epsilon^{-k}$ is
available, and the whole sum is at most $\delta_\beta\leq\epsilon$,
which gives \cref{eq:degree-small-input}. Thus smallness in arbitrary
field degree is available; it is not necessary to infer a power of
$\delta_\beta$ for every field at one previously fixed cutoff.
Splitting each remaining field into its two Fermi neighborhoods
costs at most $2^k$ in this absolute sum. Choosing $c=2/\epsilon$
at the outset absorbs this cost as well. Passing to ordered
antisymmetric kernels uses the usual $1/k!$ normalization, which
cancels the number of orderings and adds no factorial cost.

Choose paired black--white order for quadratic monomials and write
their translation-invariant kernel as $\sigma_\beta(y-x)$.
The anchored norm implies
\[
 \sum_r|\sigma_\beta(r)|e^{a\sqrt{1+|r|}}\leq C\delta_\beta\,.
\]
Differentiate its absolutely convergent Fourier series. Since
\[
 \sup_{t\geq0}t^j e^{-a\sqrt t}
       =\Big(\frac{2j}{ae}\Big)^{2j}
       \leq C_a^{j+1}(j!)^2\,,
\]
for $j\geq1$, all derivatives have the asserted bounds; order zero
uses the same summability. The initial matching coefficients and
the real-space Kasteleyn operator are real in the honeycomb gauge
used here. Choose even real cutoffs about the conjugate pair of
zeros. The fast covariance and the external convolution kernels
are then real in real space, and all finite determinant
transformations can be chosen real. The convergent cumulant
coefficients, including $\sigma_\beta$, are real. This proves the
conjugation symmetry in \cref{eq:self-energy-derivatives}.

Regard the complex symbol as the real map
$k\mapsto(\Re\mu(k),\Im\mu(k))$. At each liquid zero its Jacobian
is invertible. Choose two small disjoint convex neighborhoods on
which its derivative is uniformly close to the respective
invertible Jacobian. The $C^1$ bound on $\Sigma_\beta$ makes the
same true for the perturbed symbol. Applying the contraction
mapping theorem to
\[
 k\longmapsto k-
   [D(\Re\mu,\Im\mu)(p^\pm)]^{-1}
      (\Re(\mu+\Sigma_\beta)(k),\Im(\mu+\Sigma_\beta)(k))\,,
\]
in a ball of radius $C\delta_\beta$ gives one zero near each $p^\pm$.
The derivative bound makes the map injective on each chosen
neighborhood. On their complement $|\mu|$ has a positive minimum,
so the uniform perturbation bound excludes other zeros.
Conjugation exchanges the two zeros. Shrinking $\delta_\beta$
places them where $\chi_d=1$, so these are also the zeros of the
remaining quadratic action at this stage.

The quadratic remainder follows from the second derivative bound
and Taylor's integral formula. For the quartic kernel, anchor one
of its four vertices and take the Fourier transform in the other
three relative coordinates. The positive tree weight controls all
their first moments: each relative distance is bounded by the
square of the tree length. Its Fourier gradient is therefore
$O(\delta_\beta)$, and the mean-value formula proves the stated
remainder. Antisymmetry makes a localized quartic term with two
black fields at the same Fermi point vanish, and likewise for two
white fields. Thus the only potentially nonzero localized quartic
term has one field of each color at each Fermi point. These
arguments identify the initial quadratic and quartic data without
assuming boundedness of their subsequent \RG flow.
\end{proof}

The degree estimate \cref{eq:degree-small-input} verifies the
field-degree requirement at the start of the infrared expansion
(compare \cite[Eq.~(6.30)]{GMT20}); its spatial norm is stronger
than the one displayed there. The two zeros just constructed are
zeros after one integration, not the final dressed Fermi points.
The remaining relevant and marginal flow is constructed below,
culminating in the physical matching of
\cref{thm:matched-bulk-flow}. Uniform finite-volume
projection and the continuum comparison for correlations still
require proofs. We first establish the scale estimates and
cancellations without assuming that this flow stays bounded.

\subsection{Infrared expansion and the critical flow}
\label{sec:infrared-flow}

We construct a bounded critical history and then match its quadratic
data to the original dimer measure.

\subsubsection{Infrared scales and localization}
\label{sec:infrared-scales}

The first integration leaves a small interaction and a covariance
supported near the two simple Fermi zeros. We use the shell
decomposition of \cite[\S2.2.1]{GMT17} and \cite[\S6.1]{GMT20},
recording its uniform form for perturbed dispersions and arbitrary
torus characters. The latter uniformity is needed when extracting
the prescribed slope.

We use a decreasing momentum scale $s_j=d_*2^{-j}$, $j\geq0$;
thus increasing $j$ means going farther into the infrared. Fix a
Gevrey-$2$ periodic dispersion $D$ with two distinct simple zeros
$p^\omega$, $\omega\in\{+,-\}$. Its differential at either zero,
viewed as a real map to $\C$, is invertible. All constants below
are uniform when the Gevrey bounds, inverse differentials, separation
of the zeros, and lower bound for $|D|$ outside their neighborhoods
are uniform. Choose $d_*>0$ small enough that the radius-$4d_*$
neighborhoods are disjoint and admit Euclidean coordinates. Set
\[
 D_\omega(q)=D(p^\omega+q),\qquad
 L_\omega(q)=\nabla D(p^\omega)\cdot q\,.
\]
In these neighborhoods $|D_\omega(q)|,|L_\omega(q)|\geq c|q|$.
Choose an even, radially nonincreasing Gevrey-$2$ function $\chi$,
equal to one on the unit ball and zero outside the radius-two ball,
and define
\[
 f_j(q)=\chi(q/s_j)-\chi(q/s_{j+1}),\qquad
 \widehat g_{j,\omega}(q)=\frac{f_j(q)}{D_\omega(q)},\qquad
 \widehat g^\mathrm{lin}_{j,\omega}(q)=\frac{f_j(q)}{L_\omega(q)}\,.
\]
These are smooth compactly supported multipliers in a coordinate
chart, extended by zero; $g_{j,\omega}$ denotes their inverse
Fourier transform on $\Z^2$, with measure $\d q/(2\pi)^2$.
These are quasi-particle covariances: the phase
$e^{ip^\omega\cdot x}$ is removed from the lattice field and retained
in the kernels. The shell is supported on
$s_j/2\leq|q|\leq2s_j$, where $|D_\omega|\asymp s_j$.
For $q\ne0$, telescoping gives
$\sum_j f_j(q)=\chi(q/d_*)$; zero modes are retained separately.

\begin{lemma}[Scale bounds and linear comparison]
\label{lem:infrared-covariance}
Let $r_{j,\omega}=g_{j,\omega}-g^\mathrm{lin}_{j,\omega}$.
For each fixed nonnegative integer $m$, there are $C_m,c>0$ such that
\begin{align*}
 |\partial_x^m g_{j,\omega}(x)|
     &\leq C_m s_j^{1+m}e^{-c\sqrt{s_j|x|}},\\
 |\partial_x^m r_{j,\omega}(x)|
     &\leq C_m s_j^{2+m}e^{-c\sqrt{s_j|x|}}\,.
\end{align*}
Here the derivatives refer to the continuous Fourier interpolation;
any multi-index of order $m$ is allowed. Consequently the weighted
$\ell^1$ norms, with weight $e^{(c/2)\sqrt{s_j|x|}}$, are bounded
by $C_ms_j^{m-1}$ and $C_ms_j^m$, respectively.
The covariance $g_{j,\omega}$ has a Gram representation with both
vector norms at most $C\sqrt{s_j}$; for $r_{j,\omega}$ the norms
are at most $Cs_j$.
Moreover,
\[
 g^\mathrm{lin}_{j,\omega}(-x)=-g^\mathrm{lin}_{j,\omega}(x),\qquad
 \sum_{x\in\Z^2}g_{j,\omega}(x)=0\,.
\]
\end{lemma}
\begin{proof}
On the common rescaled annulus, Taylor's formula and the simple-zero
bounds give
\[
 \frac{f_j(q)}{D_\omega(q)}=s_j^{-1}a_j(q/s_j),\qquad
 f_j(q)\left(\frac1{D_\omega(q)}-\frac1{L_\omega(q)}\right)
                                      =b_j(q/s_j),
\]
where $a_j,b_j$ have uniform Gevrey-$2$ bounds. The second identity
uses $D_\omega-L_\omega=O(|q|^2)$ and
$|D_\omega L_\omega|\asymp s_j^2$ on the shell.
The Gevrey Fourier estimate of \cite[\S2.2.1, Lemma~2]{GMT17}
therefore gives the stated bounds: the Fourier measure contributes
$s_j^2$, and $m$ spatial derivatives contribute $s_j^m$.
The weighted lattice sum of the remaining stretched exponential
is $O(s_j^{-2})$, proving the $\ell^1$ estimates.

For a Gram representation split each multiplier into two factors
of equal modulus. Their squared $L^2$ norms are its $L^1$ norm,
bounded by $Cs_j$ and $Cs_j^2$, respectively. Finally, the linear
multiplier is odd and the full shell multiplier vanishes at zero,
which give the two identities. All constants are uniform under
the dispersion and cutoff bounds stated above.
\end{proof}

On a torus the reciprocal grid is shifted by its boundary character.
We integrate shells down to scale $n^{-1}$ and retain the remaining
modes in an unnormalized polynomial integral, so no estimate divides
by a possibly vanishing eigenvalue.

\begin{corollary}[Scale bounds for arbitrary twists]
\label{cor:infrared-torus}
Let the period lattice have a basis whose matrix divided by $n$
and its inverse have uniformly bounded norms. On a torus with
any unitary boundary character, the periodized bounds of
\cref{lem:infrared-covariance} hold whenever $ns_j\geq c_*>0$,
with torus distance in the exponential and constants depending on
$c_*$. In particular the supremum, $\ell^1$ and squared Gram
norms of $g_{j,\omega}$ are $Cs_j$, $Cs_j^{-1}$ and $Cs_j$,
uniformly in the character. The corresponding bounds for the
remainder are $Cs_j^2$, $C$ and $Cs_j^2$.
After stopping at $s_j\asymp n^{-1}$, only a bounded number of
Fourier modes remain in the two infrared balls. They can be retained
in an unnormalized Grassmann integral, including any zero modes.
\end{corollary}
\begin{proof}
Restoring the phase $e^{ip^\omega\cdot x}$ merely shifts the boundary character. In a
local gauge the covariance is a sum of image terms
$e^{i\vartheta(\gamma)}g_{j,\omega}(x+\gamma)$ over the period
lattice. The triangle inequality gives the weighted $\ell^1$
bound directly from its planar counterpart. For the supremum bound,
the images have separation at least $cn$; summing the stretched
exponential over these images, using $ns_j\geq c_*$, gives
$Ce^{-c'\sqrt{s_j\dist(x,0)}}$. The same proof treats the remainder.

For the Gram estimate use the discrete Fourier representation.
The shifted reciprocal grid has at most $C(1+n^2s_j^2)$ points
in a radius-$2s_j$ ball. Its normalized counting measure has
mass $O(s_j^2+n^{-2})\leq Cs_j^2$ there. The multiplier is bounded
by $C/s_j$, and its remainder by $C$, proving the claims by the
equal-modulus factorization. The same counting estimate gives
the bounded number of final modes. No division by their eigenvalues
is needed to leave them in the final polynomial integral.
\end{proof}

These estimates also apply after absorbing the first quadratic
kernel $\Sigma_\beta$. Indeed, the slow covariance and its dressed
dispersion are
\[
 g_{\mathrm{slow}}=\chi_d/D,\qquad D=\mu+\chi_d\Sigma_\beta,
\]
since the inverse covariance on the slow field space is
$\mu/\chi_d+\Sigma_\beta$. Near the zeros $\chi_d=1$, so
\cref{prop:infrared-input} gives two simple zeros; elsewhere the
small perturbation creates none, and
\cref{eq:self-energy-derivatives} gives uniform Gevrey bounds.
Choose the new infrared balls inside $\{\chi_d=1\}$.
The part of $\chi_d/D$ outside their cutoff is smooth and
integrated in one fixed-scale step using
\cref{sec:first-integration}. This is an exact reorganization of
the unnormalized integral, also at singular twists. Further
quadratic corrections still require the running localization.

A local monomial with $2m$ fields and $l$ derivatives has engineering
dimension $2-m-l$: each field costs $s_j^{1/2}$, each derivative
$s_j$, and the spatial sum $s_j^{-2}$. Thus the constant and linear
quadratic jets and the local quartic term must be retained. The
following is the localization of \cite[\S6.2.3]{GMT20}, with an
explicit check that no arithmetic nonresonance hypothesis is needed.

\begin{lemma}[The local relevant and marginal monomials]
\label{lem:infrared-localization}
For a translation-invariant, balanced effective interaction with
summable spatial moments, write its low-momentum fields as
$\psi_{b,\omega},\psi_{w,\omega}$, using the opposite phases
$e^{\pm ip^\omega\cdot x}$ for black and white. Apart from the constant term, its local
relevant and marginal part has precisely the form
\begin{equation}\label{eq:infrared-local-part}
 \sum_{\omega=\pm}\sum_x
    \psi_{b,\omega}(x)
       [\nu_\omega+v_\omega\cdot(-i\nabla)]\psi_{w,\omega}(x)
 +\lambda\sum_x
       \psi_{b,+}(x)\psi_{w,+}(x)
       \psi_{b,-}(x)\psi_{w,-}(x)\,.
\end{equation}
Derivatives here are Fourier multipliers on the quasi-particle fields.
No arithmetic nonresonance assumption on $p^\pm$ is required.
For a fixed ordered $2m$-field kernel, let
\[
 M_{m,r}=\sum_{x_2,\ldots,x_{2m}}
 |W_{2m}(0,x_2,\ldots,x_{2m})|
       \bigg(1+\sum_{i=2}^{2m}|x_i|\bigg)^r\,,
\]
including the finite sum over flavor labels. At external quasi-particle
momenta of modulus at most $s$, subtracting the displayed local
part leaves quadratic and quartic amplitudes bounded, respectively,
by $CM_{1,2}s^2$ and $CM_{2,1}s$. These are single-kernel estimates;
they do not assume uniform bounds on the moments at later \RG scales.
\end{lemma}
\begin{proof}
This is the bulk localization of \cite[\S6.2.3]{GMT20}, with
black and white fields in place of its $+$ and $-$ fields. A
quadratic monomial with flavors $\omega,\omega'$ carries the
phase $e^{i(p^\omega-p^{\omega'})\cdot x}$. Translation invariance requires
its total momentum to be a reciprocal lattice vector. By choosing
the momentum neighborhoods smaller than the separation of the two
Fermi points on the momentum torus, this is impossible for
$\omega\ne\omega'$. Thus only the two same-flavor constant and
linear quadratic jets remain. For a local quartic term, antisymmetry
requires two distinct black and two distinct white flavors, giving
exactly the quartic monomial in \cref{eq:infrared-local-part}.
This argument also excludes any additional resonant quartic:
a repeated flavor makes the local exterior product zero, regardless
of arithmetic relations between the Fermi points.

The remainder estimates are the Taylor gains used in
\cite[\S6.1.4]{GMT17}. Anchor the first coordinate at zero and
absorb the fixed phases $e^{\pm ip^\omega\cdot x}$ into the kernel; they have modulus
one. In the quadratic case the remaining Fourier factor is
$e^{iq\cdot x_2}$. Subtracting its constant and linear terms
leaves modulus at most $s^2|x_2|^2/2$. In the quartic case
the factor is $e^{i\sum_{i=2}^4q_i\cdot x_i}$, with signs
absorbed into $q_i$. Subtracting its value at zero leaves
modulus at most $s\sum_{i=2}^4|x_i|$. Summing these bounds
against the absolute kernels gives $CM_{1,2}s^2$ and
$CM_{2,1}s$. If a quartic flavor choice has zero local exterior
product, subtracting that product changes nothing, so the same
first-order gain still applies. The engineering dimension is $2-m$ for
$2m$ fields, lowered by one per derivative. Thus the two remainders
and all terms with at least six fields have dimension at most $-1$.
The required initial moments follow from
\cref{eq:bulk-small-kernel}.
\end{proof}

\subsubsection{Reference model and bulk comparison}
\label{sec:infrared-cancellation}
\label{sec:conditional-bulk-comparison}

The marginal quartic flow is controlled by comparison with the
linear continuum reference model of \cite{BM04,GMT20}. Its symmetries
eliminate the constant quadratic jet, and its Ward identities give
the quartic cancellation used below. We first record the symmetries
for the bulk localization; they hold coefficientwise before any
uniform convergence assertion.

\begin{lemma}[Symmetries of the linear reference interaction]
\label{lem:dirac-parity}
Let $g_\omega$ be any finite sum of the linear shell covariances
in \cref{lem:infrared-covariance}, with independent Gaussian
fields for the two flavors. For
\[
 Q_x(\psi)=\psi_{b,+}(x)\psi_{w,+}(x)
                 \psi_{b,-}(x)\psi_{w,-}(x),\qquad
 V(\psi)=\lambda\sum_x Q_x(\psi)\,,
\]
the formal connected effective interaction
$\log\cE_g e^{V(\psi+\xi)}$ is invariant under
$\psi_{c,\omega}(x)\mapsto i\psi_{c,\omega}(-x)$.
Its $2m$-field kernel changes by $(-1)^m$ under simultaneous
reversal of all relative positions. In particular the quadratic
Fourier kernel is odd, so its constant term vanishes at every
order in $\lambda$; the quartic kernel is even, so its first
derivatives at zero vanish. Moreover, the local quartic coefficient
is an odd formal function of $\lambda$: all its even-order
coefficients vanish. These identities also hold for any successive
integrations with such covariances.
\end{lemma}
\begin{proof}
First restrict the vertices to a finite set invariant under
$x\mapsto-x$. The transformed Gaussian field has the same
covariance: the covariance of the pair
$\psi_{b,\omega}(x),\psi_{w,\omega}(y)$ is $g_\omega(x-y)$, and
after the transformation it becomes
\[
 i^2g_\omega(-x+y)=-g_\omega(y-x)=g_\omega(x-y)\,,
\]
by the oddness of the linear shell covariances recorded in
\cref{lem:infrared-covariance}, which is inherited by their finite
sums. Wick's rule therefore gives the same expectation for every
polynomial. Each quartic vertex acquires the factor $i^4=1$,
and reflection permutes its spatial sum. Apply this change to
both the external and integrated fields in the moment generating
polynomial and then take its formal logarithm. This proves the
identity coefficient by coefficient. At any fixed order the
anchored connected coefficients have an infinite-volume limit:
their graphs are connected, their covariances are bounded and
summable, and a spanning tree sums all relative positions.
Reflection invariance thus passes to the bulk kernels.
The factor on $2m$ external fields is $i^{2m}=(-1)^m$,
which proves the assertions about their Fourier transforms.

There is a second Gaussian symmetry: exchange the black and white
fields of the plus flavor, leaving the minus flavor unchanged.
If $T$ denotes this exchange, its plus-flavor contraction is
$-g_+(y-x)=g_+(x-y)$, so Wick's rule is again unchanged.
But $Q_x(T\psi)=-Q_x(\psi)$. Writing $W_\lambda$ for the
connected effective interaction, the change of variables gives
$W_\lambda(T\psi)=W_{-\lambda}(\psi)$. Extract the local quartic
coefficient; by \cref{lem:infrared-localization} its only
monomial is $Q_x$, on which $T$ acts by minus one. Its coefficient
at $-\lambda$ is therefore the negative of that at $\lambda$.
This proves the assertion about all even orders.
Successive integrations preserve both identities by induction.
This argument is a formal symmetry proof and makes no claim of
convergence of the full series uniformly as the infrared cutoff
is removed.
\end{proof}

The reference results of \cite{GMT20} are proved for a specific
class of infrared cutoffs, radial in the coordinates in which the
linear dispersion is isotropic. In order to import them verbatim we
use the same cutoff, which we now specify. For the real dispersions
used below, $D(-k)=\overline{D(k)}$ and
$p^-=-p^+$, so $L_-(q)=-\overline{L_+(q)}$. Put
$Aq=(\Re L_+(q),\Im L_+(q))$, so that $|L_+(q)|=|Aq|$; the real
matrix $A$ is invertible since the differential of $D$ at its zero
is. The invertibility of $A$ allows us
to use the adapted radial cutoff $\chi(Aq/s_j)$ throughout the
infrared analysis. A fixed invertible linear change preserves the
Gevrey shell bounds of \cref{lem:infrared-covariance}; the cutoff
remains even, so the symmetry above also holds. Changing the outer
cutoff adds a smooth fixed-scale covariance.
For definiteness we choose this cutoff within the class specified
in \cite[\S6.1, Eq.~(6.6)]{GMT20}. After exchanging the two
flavor labels if necessary, put $\Delta=\det A>0$.
The matrix $\mathcal M$ of their Eq.~(4.1) satisfies
$|\mathcal M^{-1}q|=|Aq|/\sqrt\Delta$. Taking their scalar
profile to be $\bar\chi(r)=\chi(\sqrt\Delta\,r/s_0)$
therefore gives
\[
 \bar\chi_h(q)=\chi(Aq/s_j),\qquad h=-j,
\]
with the harmless initial scale shift used below. Its transition
radii are $c_0/2,c_0$, where $c_0=2s_0/\sqrt\Delta$,
exactly as required there. Thus the reference analysis uses
one and the same admissible cutoff on the two sides; no
invariance under a change of infrared cutoff profile is invoked.

\begin{remark}[Cutoff corrections in the reference comparison]
\label{rem:cutoff-ward}
The continuum Ward identities of \cite[\S4]{GMT20} apply to the
cutoff-removed reference model. A local phase change in the lattice
integral also differentiates the momentum cutoff and the nonlocal
interaction, so those identities cannot be applied directly to its
effective kernels. We instead transfer the reference cancellations
by the kernel comparisons below. In particular, the quartic
cancellation is the external input
\cref{eq:reference-constant-history}, and the source exponents
are transferred in \cref{sec:sharp-source-flow}.
\end{remark}

One more difference between our reference interaction and that of
\cite[\S4]{GMT20} must be addressed: there, the quartic interaction
is a smooth density--density potential between the two flavors,
of the form $\lambda\int v(x-y)\rho_{+,x}\rho_{-,y}$, rather than
the strictly local term $\lambda\sum_xQ_x$, since the Ward identity
analysis of the continuum model requires a regular potential. At
the local quartic jet the two are compatible: for even $v$ with
$\int v=1$ and finite second moment,
\[
 |\widehat v(q)-1|
   \leq\tfrac12|q|^2\int |x|^2|v(x)|\,\d x\,.
\]
Thus replacing the local quartic vertex changes its amplitude by
$O(|\lambda|s^2)$ at quasi-particle momenta of size at most $s$, which is
an irrelevant correction of the kind treated in
\cref{lem:infrared-localization}, and it preserves the independent
density phase symmetries, since the potential still depends on the
fields only through the densities.

\paragraph{Conditional all-orders bulk comparison}

A prescribed history assigns local coefficients at every scale,
without requiring them to solve the flow equations. The renormalized
trees of \cite[\S6]{GMT20} define analytic beta functions of such
histories. We verify that our first-step kernels are admissible
endpoints, compare these functions with the linear reference ones,
and obtain the estimates needed to construct a critical history.

We use the renormalized determinant-tree estimate of
\cite[Prop.~3 and the paragraph following Eq.~(6.68)]{GMT20} as a
published input. It gives convergence for small running
coefficients and a factor $2^{-\vartheta(j-k)}$ when a tree
reaching infrared scale $j$ is marked at scale $k\leq j$;
any fixed $\vartheta<1$ is allowed after reducing the smallness
threshold. We use $\vartheta=3/4$ below, and reserve the stronger
exponent $\alpha=7/8$ when differentiating a recursion.
Its endpoint estimates concern
absolute spatial sums, not pointwise bounds. We verify explicitly
that the initial kernels here are admissible endpoints and that
the resulting comparison is uniform in their range and degree.
No conclusion about a solution of the running recursion or about
source correlations is included in this input verification.

Write the initial, source-free interaction in ordered form as
\[
 V_0(\psi)=\sum_{\substack{\fn\geq2\\\fn\text{ even}}}\frac1{\fn!}
   \sum_{\boldsymbol\alpha,\bx}
      W_{0,\fn}(\boldsymbol\alpha;\bx)
                       \prod_{i=1}^{\fn}\psi_{\alpha_i}(x_i)\,,
\]
where $\alpha_i$ includes color and Fermi label and the kernels
are antisymmetric. The phases $e^{\pm ip^\omega\cdot x}$ are included in the kernels.
For $s_0>0$, let
$\ell_{s_0}(\bx)=\min_T\sum_{uv\in T}\sqrt{s_0|u-v|}$.
Fix $\kappa>0$ small and a degree parameter $0<\epsilon<1$.
The endpoint norm we use is
\begin{equation}\label{eq:infrared-endpoint-norm}
 \fE(V_0)=
 \sum_{\substack{\fn\geq2\\\fn\text{ even}}}\epsilon^{-\fn}s_0^{\fn/2-2}
 \sup_{x_1}\frac1{(\fn-1)!}
 \sum_{\boldsymbol\alpha,x_2,\ldots,x_{\fn}}
 |W_{0,\fn}|\,(1+s_0\diam\bx)^2
                    e^{2\kappa\ell_{s_0}(\bx)}\,.
\end{equation}
The factorial removes nonroot orderings, $\epsilon^{-\fn}$ pays
field choices, and $s_0^{\fn/2-2}$ makes the kernel dimensionless.
The spatial weights supply the two localization moments and a
remaining tree-length weight. The finite root-label sum is included.
Let $\cL$ extract the local quadratic and quartic part of
\cref{lem:infrared-localization}, and put $R_0=(1-\cL)V_0$.
The norm is imposed before this subtraction,
avoiding spatial absolute values of derivatives of delta kernels.

Write $Z_j$ for the running field normalization obtained by absorbing
the scalar parts of the reference linear quadratic jets into the
Gaussian part, with $Z_0=1$. Its recursion is given in
\cref{eq:reference-normalization-recursion} below.

\begin{theorem}[Bridge to the GMT induction]
\label{thm:conditional-renormalized-trees}
For any fixed sufficiently small $\epsilon$, the bulk dimer
initial data can be chosen with
$\fE(V_0)\leq\delta$, where $\delta\to0$ as
$\beta\to\infty$. In particular, after splitting into Fermi
flavors and separating the local jets, the first-step kernels
are admissible initial endpoints for the induction and tree
bound of \cite[Eqs.~(6.30)--(6.31), Prop.~3]{GMT20}; the positive tree-length
weight in \cref{eq:bulk-small-kernel} also supplies the spatial
moments needed for localization. The source-marked norm
\cref{eq:source-marked-norm} supplies the corresponding spatial
moments for source endpoints, including coincident and repeated
source marks. The free edge source is an $O(\rho)$ input, not
an interaction endpoint tending to zero with $\beta$.

Assume $\delta\leq\epsilon$ and that a
prescribed history, through scale $j$, has local dimensionless
coefficients of modulus at most $\epsilon$, with normalization
ratios satisfying $|Z_{k+1}/Z_k-1|\leq\epsilon$.
Assume also the uniform shell and differentiated Gram bounds of
\cref{lem:infrared-covariance} for the normalized propagators.
Then the source-free renormalized expansion and its localized
beta functions are analytic in this history and in $V_0$ on a
smaller polydisc, uniformly in $j$.

After the engineering dimension of a root coefficient is removed,
the sum of trees containing an initial $R_0$ endpoint is bounded
by $C\delta 2^{-\vartheta j}$. A marked local endpoint or
propagator on scale $k$ has the additional bound
$C2^{-\vartheta(j-k)}$ times its dimensionless marked norm and
the activities of the other endpoints. These estimates also hold
after one parameter derivative, with the differentiated activity
replaced by its derivative. Constants are independent of any
cutoff on the degree or spatial range of $V_0$.
\end{theorem}
\begin{proof}
We apply \cite[Prop.~3 and the ensuing short-memory bound]{GMT20};
see also \cite[Prop.~8, Cor.~1 and Rem.~16]{GMT17} for the tree
summation and analyticity in a prescribed history. Our scale $j$
corresponds to their infrared scale $h=-j$, up to the fixed initial
cutoff. It remains to verify the endpoint norms and their uniformity.

First, the quadratic convention is the one in the cited induction.
In \cite[Eqs.~(6.41)--(6.49)]{GMT20}, only the scalar linear jet
of $W^{(h),R}_{2,0}$ is absorbed into $Z$; the remaining two
linear coefficients are $a_{h,\omega},b_{h,\omega}$ and remain
endpoints in Eq.~(6.51). Under $h=-j$, our $\nu_{j,\omega}$
and the two components of $v_{j,\omega}$ are these
constant and linear coefficients, with the fixed factor $s_0$
included in the units of $\nu$. The scalar absorption is the
Gaussian identity of their Eqs.~(6.43)--(6.45), with the same
cutoff multiplier in the full and reference covariances.
We retain the resulting normalization ratios at each vertex.
Thus keeping $v$ as an endpoint and using the reference scalar
does not require a different induction from Proposition~3.

Since $\ell_{s_0}\leq\sqrt{s_0}\ell$ and
$\sqrt{\diam\bx}\leq\ell(\bx)$, choose $s_0,\kappa$ so that
$2\kappa\sqrt{s_0}$ is strictly smaller than the spatial exponent
in \cref{eq:bulk-small-kernel}. The remaining exponent pays
for the second moment. The factors
$\epsilon^{-\fn}s_0^{\fn/2-2}$, the flavor choices and polynomial
factors in $\fn$ are absorbed by the available fixed degree weight.
\Cref{thm:bulk-first-step} then gives
$\fE(V_0)\leq C\delta_\beta$; all choices precede increasing
$\beta$. Any intervening smooth covariance is integrated at a
fixed scale by the same first-step estimate.

At an initial endpoint, sum its relative coordinates absolutely.
The quadratic and quartic Taylor subtractions cost respectively
two and one moments, supplied by the norm; higher degrees need no
subtraction. Here is the scale gain, with the degree dependence
retained. Let $I_{\fn}$ be the $\fn$th spatial sum in
\cref{eq:infrared-endpoint-norm}, including its factorial,
moments and exponential weight, and put
$b_{\fn}=\epsilon^{-\fn}s_0^{\fn/2-2}I_{\fn}$; then
$\sum_{\fn} b_{\fn}\leq\delta$. If its fields have momenta of order
$s_k$, its dimensionless remainder costs at most
\[
 C^{\fn} b_{\fn}\epsilon^{\fn}
       (s_k/s_0)^{\fn/2-2+r(\fn)},\qquad
 r(2)=2,\quad r(4)=1,\quad r(\fn)=0\ (\fn\geq6).
\]
For $\fn=2,4$ this is the Taylor remainder formula, with the
rescaled moments in $I_{\fn}$; for $\fn\geq6$ it is the
engineering rescaling. In every case the exponent is at least
one. The ordered-kernel factorial accounts for permutations of
the nonroot fields; choices of labels, contracted fields and
anchors cost $C^{\fn}$ times polynomial factors in $\fn$, absorbed
by choosing the fixed degree reserve before $\beta$.
Thus summing initial endpoints costs $C\sum_{\fn} b_{\fn}\leq C\delta$,
uniformly in degree and range. Joining each endpoint's
internal spatial tree to the contraction tree pays the output
weight $e^{\kappa\sqrt{s_j\diam\bx}}$ by
$\sqrt{a+b}\leq\sqrt a+\sqrt b$. This is a bound on fully summed
kernels, so it introduces no further entropy sum over supports.

The remaining hypotheses of the cited estimate are exactly our
shell and Gram bounds and the small prescribed local coefficients.
Its localization gains are $r(2)=2$, $r(4)=1$ and $r(\fn)=0$ for
$\fn\geq6$, giving $2-\fn/2-r(\fn)\leq-1$ at every nonconstant
internal vertex. The normalization loss in the exponent at a
vertex with $\fn$ external fields is at most $C\epsilon \fn$.
Choose $\epsilon$ so that
$\fn/2+r(\fn)-2-C\epsilon \fn\geq15/16$ for all even $\fn\geq2$;
the worst cases have $\fn=2,4,6$, so one choice suffices for all
degrees. A one-child chain of length $d$ thus contributes
at most $2^{-15d/16}$. Reserve $2^{-\alpha d}$ with
$\alpha=7/8$ along a marked path; the remaining geometric
factor sums its intermediate scales. The branching-vertex
factorials and field-degree reserve are those of
\cite[Eq.~(6.67)]{GMT20} and sum the remaining labels and trees.
This leaves short memory $2^{-\vartheta(j-k)}$ for a mark
on scale $k$. An initial remainder is on scale zero. Parameter
derivatives mark endpoints or propagators, as in
\cite[\S6.4.1]{GMT20}; differentiated determinants are estimated
by their Gram minors. The degree reserve pays the choices of marks.
A slightly stronger memory exponent pays normalization ratios and
leaves the stated exponent $\vartheta$.

For completeness, first truncate endpoint degree and range.
The estimates above are uniform in these truncations. The difference
of two expansions is bounded by the same marked estimate times the
endpoint-norm difference. Thus they converge in the claimed root
norms, uniformly on smaller complex polydiscs; analyticity and
parameter differentiation pass to the limit. The source norm
\cref{eq:source-marked-norm} has the same spatial-tree property,
with source positions included as vertices. Its absolute Taylor
sum pays for coincident and repeated marks before localization.
The free source is kept as its separate $O(\rho)$ endpoint.
\end{proof}

The bridge verifies the initial data and the conditional induction;
it does not assert that an arbitrary small initial quadratic jet
produces a bounded flow. The counterterm choice and its matching to
the physical model are made in \cref{sec:critical-flow,sec:physical-matching}.
Likewise, uniform twisted-sector normalization and the sharp source
comparison are separate from the bulk induction.

\begin{lemma}[Sampling the reference expansion at every order]
\label{lem:all-order-reference-sampling}
Choose the adapted initial cutoff so that
$8s_0\|A^{-1}\|<2\pi$. For the pure linear reference expansion,
the bulk localized lattice beta functions agree with the
continuum beta functions, coefficient by coefficient at every
order. This includes field normalization and any finite Taylor
jet at zero external quasi-particle momenta. In the conditional
analytic domain of \cref{thm:conditional-renormalized-trees},
the equality holds for the summed functions.
\end{lemma}
\begin{proof}
The point is that the reference model has only strictly local
vertices and propagators with compact Fourier support in a small
ball, so that the lattice and the continuum model have the same
amplitudes: by the Poisson summation formula, the lattice sum over
the position of a vertex of a function whose Fourier support
contains no nonzero reciprocal-lattice vector
equals the corresponding continuum integral. Every reference
endpoint is a local quartic vertex or a local
quadratic counterterm, possibly with one derivative. All
propagators have Fourier support in a fixed ball of radius
$R=2s_0\|A^{-1}\|$, by the choice of the adapted cutoff
$\chi(Aq/s_0)$. Expand a fixed-order coefficient into its
finitely many contraction graphs before making estimates.
At any original vertex at most four propagators occur. Therefore
the Fourier support in that vertex position lies in the ball
of radius $4R$, independently of the total number of vertices.
Connectedness and a positive infrared cutoff make the integrand
Schwartz in all relative vertex positions. Poisson summation
in each such position has no nonzero reciprocal-lattice term.
It replaces all lattice sums by continuum integrals exactly.
Taking an external Taylor jet multiplies by position polynomials
or differentiates Fourier multipliers, without enlarging support.

Apply this argument recursively to the Taylor subtractions and
the quadratic counterterm insertions defining the renormalized
expansion. A localized subgraph becomes a vertex of degree at
most four with at most one derivative, so the same support
bound continues to apply. The scalar Gaussian normalization
terms are treated as connected vacuum graphs. Formal division
by the field-normalization series, whose constant coefficient
is one, preserves equality. Thus equality holds at every
coefficient of the prescribed running history. The uniform
analytic tree bounds make these coefficientwise equalities
equalities of functions. This argument only concerns the local
linear reference endpoints; nonlocal admissible endpoints are handled
by the marked estimates instead.
\end{proof}

The reference propagators use the linear dispersions $L_\omega$,
whereas the actual ones use $D_\omega=L_\omega+E_\omega$ with a
curvature term $E_\omega(q)=O(|q|^2)$. On the $k$th shell the
curvature is of relative size $s_k$, so the difference between the
two propagators is an irrelevant correction. The running field
normalization, however, divides the curvature by $Z_k$ while the
linear part is rescaled by the reference multiplier, and the
normalization ratios may drift by a factor $1\pm\epsilon$ per
scale. The next lemma shows that this drift only lowers the
geometric decay rate of the curvature correction from $2^{-k}$ to
$2^{-(1-\eta)k}$, which remains larger than the short-memory
exponent $\vartheta$ used below.

\begin{lemma}[Curvature after running field normalization]
\label{lem:running-curvature}
Let $D_\omega=L_\omega+E_\omega$, with
$E_\omega(q)=O(|q|^2)$ and uniform Gevrey bounds. Suppose
$Z_0=1$ and $|Z_{k+1}/Z_k-1|\leq\epsilon$. In a normalized
shell, allow a common multiplier $a_{k,\omega}(q)$ of
$L_\omega(q)$, bounded away from zero with uniform rescaled
Gevrey bounds, in both the full and reference propagators.
Their multipliers are
\[
 \frac{f_k(q)}{a_{k,\omega}(q)L_\omega(q)+E_\omega(q)/Z_k},
 \qquad
 \frac{f_k(q)}{a_{k,\omega}(q)L_\omega(q)}\,.
\]
For sufficiently small fixed $s_0$, the difference has the shell
and Gram bounds of an ordinary normalized propagator multiplied
by
\begin{equation}\label{eq:running-curvature-gain}
 d_k=Cs_0\,2^{-(1-\eta)k},\qquad
 \eta=-\log_2(1-\epsilon)\,.
\end{equation}
The constants are uniform in the history. The multipliers used
when absorbing a small quadratic jet satisfy the assumptions on
$a_{k,\omega}$, by the rescaled cutoff bounds.
\end{lemma}
\begin{proof}
The product of the normalization ratios gives
$|Z_k|\geq(1-\epsilon)^k=2^{-\eta k}$. On the $k$th shell,
the ratio of $E_\omega/Z_k$ to $a_{k,\omega}L_\omega$ is at
most $Cs_k/|Z_k|\leq Cs_0 2^{-(1-\eta)k}$.
For $s_0$ small, both denominators have modulus at least $cs_k$.
Subtracting their reciprocals gives
\[
 -\frac{f_k E_\omega/Z_k}
  {(a_{k,\omega}L_\omega)
   (a_{k,\omega}L_\omega+E_\omega/Z_k)}\,.
\]
After rescaling $q=s_ku$, its Gevrey bounds are those of a
shell propagator times $d_k$. Fourier integration and the
equal-modulus Gram factorization, as in
\cref{lem:infrared-covariance}, give the claimed spatial
and Gram bounds. The same proof applies to convex interpolations
between the two propagators, used in a telescoping difference.
\end{proof}

\begin{proposition}[All-orders comparison at a prescribed history]
\label{prop:conditional-beta-comparison}
Use the hypotheses of \cref{thm:conditional-renormalized-trees,lem:running-curvature}, with $1-\eta>\vartheta$.
Let $\lambda_k$ denote the local quartic history. Normalize
fields by the scalar quadratic jet of the pure linear reference
trees; let $u_k$ be the maximum modulus of the remaining
dimensionless quadratic coefficients (constant and linear).
Write $B_j^\lambda$ for the increment of the normalized quartic
coefficient, and $B_j^R$ for its pure linear reference counterpart
at the same quartic and normalization history. Then
\begin{equation}\label{eq:conditional-beta-comparison}
 |B_j^\lambda-B_j^R|
 \leq C\delta 2^{-\vartheta j}
       +C\epsilon\sum_{k=0}^j
          2^{-\vartheta(j-k)}(u_k+d_k)\,.
\end{equation}
The same right side bounds the nonreference quadratic beta
functions, after removing their engineering scaling and the
propagated local terms. All bounds sum every perturbative order.
They are conditional on the small prescribed history, not on
the smallness of any fixed-order truncation.
\end{proposition}
\begin{proof}
The reference quadratic constant vanishes by
\cref{lem:dirac-parity}. The adapted reference symmetries of
\cite[App.~A]{GMT20} make its linear jet a common scalar multiple
of $L_\omega$. Their hypotheses hold here: the two dispersions
are conjugate, the cutoff is radial in $Aq$, and the reference
contains only local quartic vertices. Absorbing the scalar jet
preserves these symmetries, since its multiplier is a common radial
function; every other quadratic jet is retained among the $u_k$.
\Cref{lem:all-order-reference-sampling} transfers the continuum
identities to the reference lattice expansion. Thus the remaining
quadratic contributions are all nonreference.

First collect the trees containing an initial $R_0$ endpoint.
\Cref{thm:conditional-renormalized-trees} bounds their sum
by $C\delta 2^{-\vartheta j}$. Among the other trees, a remaining
quadratic endpoint at scale $k$ contributes at most
$C\epsilon u_k2^{-\vartheta(j-k)}$ after summation.
Indeed a lone local quadratic endpoint propagates unchanged
(with the relevant engineering factor) and is excluded from
the beta increment. Every contributing such tree has at least
one other endpoint, supplying the factor $\epsilon$.

The trees left have only local quartic endpoints, and differ
between the two expansions only through their propagators. Compare
their full and reference propagators by interpolating the
covariances: on every scale $k$ replace the propagator by the
convex combination of the reference propagator and the full one,
with a common interpolation parameter, and write the difference of
the two tree values as the integral of the derivative with respect
to that parameter. Each derivative marks one
propagator or one entry of a Gram determinant on some scale $k$,
replacing it by the difference of the two propagators on that
scale, whose shell and Gram norms are at most $d_k$ times those of
an ordinary propagator by \cref{lem:running-curvature}; a
differentiated determinant is expanded along the marked entry,
leaving a cofactor to which the Gram bound still applies.
The cofactor argument in
\cref{thm:conditional-renormalized-trees,eq:running-curvature-gain} give
$C\epsilon d_k2^{-\vartheta(j-k)}$; there is at least one
quartic endpoint. Summing over the position and scale of the
mark proves the second term of
\cref{eq:conditional-beta-comparison}. The same argument
applies to the root quadratic Taylor jets, with their engineering
factors removed. Field normalization uses the same reference
scalar in both expansions, so its factors and their Taylor
series are already included in these tree estimates.
\end{proof}

We next eliminate the reference normalization history. Recall that
the reference expansion absorbs, at each scale, the scalar part of
the generated linear quadratic jet into the Gaussian part, which
rescales the fields; the product of these rescalings is the running
field normalization $Z_j$ of \cite[\S6.2]{GMT20}. In
\cref{prop:conditional-beta-comparison} the normalization ratios
were part of the prescribed history. They are, however, determined
by the quartic history, and the fixed-point argument below is
cleaner if the quartic history is the only reference unknown. At
step $j$, let $\cZ_j(\lambda_{\leq j},z_{\leq j})$ be the
scalar reference quadratic jet, with its sign chosen so that
absorbing it multiplies the kinetic coefficient by
$1+\cZ_j$. Set
\begin{equation}\label{eq:reference-normalization-recursion}
 z_0=0,\quad Z_0=1,\qquad
 z_{j+1}=\cZ_j(\lambda_{\leq j},z_{\leq j}),\qquad
 Z_{j+1}=Z_j(1+z_{j+1})\,.
\end{equation}
All fields in $\cZ_j$ are normalized by $Z_j$ before
this step. Thus the recursion is causal: its right side does
not contain $z_{j+1}$. If $b_j^4$ is the quartic contribution
generated at the step in these same units, the normalized
reference beta function before eliminating $z$ is
\begin{equation}\label{eq:normalized-reference-beta}
 \cB_j(\lambda,z)=
 \frac{\lambda_j+b_j^4(\lambda,z)}
      {(1+\cZ_j(\lambda,z))^2}-\lambda_j\,.
\end{equation}
This also specifies the normalization convention in the preceding
proposition. Denote the function after substitution of
\cref{eq:reference-normalization-recursion} by
$\widehat B_j^R(\lambda)=\cB_j(\lambda,z(\lambda))$.

\begin{lemma}[Derivative bounds after reference normalization]
\label{lem:reference-history-derivatives}
If $\sup_k|\lambda_k|\leq\epsilon$ is sufficiently small, the
recursion \cref{eq:reference-normalization-recursion} is well defined at
every scale, with $|z_{j+1}|\leq C\epsilon^2$. Each of its
coefficients is analytic in the finite history on which it
depends. Uniformly in $0\leq k\leq j$,
\begin{equation}\label{eq:reference-history-derivatives}
 \Big|\frac{\partial z_{j+1}}{\partial\lambda_k}\Big|
       \leq C\epsilon2^{-\vartheta(j-k)},\qquad
 \Big|\frac{\partial\widehat B_j^R}{\partial\lambda_k}\Big|
       \leq C\epsilon2^{-\vartheta(j-k)}\,.
\end{equation}
The estimates hold on a smaller complex polydisc, independently
of the number of scales.
\end{lemma}
\begin{proof}
These are the reference-normalization estimates of
\cite[\S6.4.1]{GMT20}. The proof there differentiates a marked
renormalized tree; the derivative estimate uses smallness of the
history, not decay of its increments. The reference functional is
identified with ours by \cref{lem:all-order-reference-sampling}.
We record the causal estimate needed for the uniform exponent here.

The one-quartic-vertex quadratic contraction has no linear momentum
jet. Hence every tree for $\cZ_j$ has at least two quartic endpoints,
and $|\cZ_j|\leq C\epsilon^2$. This closes the smallness induction
for $z$. With the short-memory exponent $\alpha=7/8$, the marked
estimate of \cite[\S6.4.1]{GMT20} gives
\begin{align*}
 |\partial_{\lambda_k}\cZ_j|
       &\leq C\epsilon2^{-\alpha(j-k)},\\
 |\partial_{z_l}\cZ_j|
       &\leq C\epsilon^2 2^{-\alpha(j-l)},\qquad 1\leq l\leq j\,.
\end{align*}
Here a $\lambda$ derivative removes one endpoint activity, whereas
a $z$ derivative marks a propagator or a crossing of normalization
ratios and removes none. Thus no absolute bound on $Z_l$ is used.
For $D_{j,k}=|\partial_{\lambda_k}z_{j+1}|$, causality gives
\[
 D_{j,k}\leq C\epsilon2^{-\alpha(j-k)}
  +C\epsilon^2\sum_{l=k+1}^j2^{-\alpha(j-l)}D_{l-1,k}\,.
\]
In the norm $\sup_{k\leq j\leq J}2^{\vartheta(j-k)}D_{j,k}$,
the convolution has norm at most
$C\epsilon^2 2^{\vartheta}/(1-2^{-(\alpha-\vartheta)})<1/2$.
This proves the first bound uniformly in the terminal scale $J$.

The propagated one-vertex quartic term is subtracted in
\cref{eq:normalized-reference-beta}. The same marked estimate
therefore gives
$|\partial_{\lambda_k}\cB_j|\leq C\epsilon2^{-\alpha(j-k)}$
and $|\partial_{z_l}\cB_j|\leq C\epsilon2^{-\alpha(j-l)}$.
Composing with the bound just proved gives the second inequality
of \cref{eq:reference-history-derivatives} by the same geometric
convolution. The tree sums are analytic on the smaller complex
polydisc and their normalization denominators stay away from zero,
so all these estimates hold there with constants independent of $J$.
\end{proof}

The one further published input is the continuum cancellation
\cite[Eq.~(6.96)]{GMT20}. It states that, along a constant quartic
history, the quartic beta function of the reference model vanishes
asymptotically in the scale: the marginal coupling neither grows
nor decays, and the reference model has a line of fixed points
parametrized by the quartic coupling. This is the property that
distinguishes the present problem from a generic marginal flow, and
its proof, via the Ward identities of the continuum model, is the
deepest external input of this section; see \cite[Thm.~3.1]{BM04}
and \cite[\S4]{GMT20}. In the normalization just specified it is
\begin{equation}\label{eq:reference-constant-history}
 |\widehat B_j^R(t,t,\ldots,t)|
       \leq C|t|^2 2^{-\theta j},\qquad \theta=1/2\,,
\end{equation}
for small real $t$. The dominant reference functional in that
statement has only local quartic endpoints and linear
propagators; its smooth nonlocal potential is part of the initial
irrelevant interaction and is excluded from this functional.
Use the cutoff fixed above, the zero-momentum Taylor
localization of \cite[Eq.~(6.37)]{GMT20}, and its scalar
Gaussian absorption at each step. To specify the latter
without suppressing a cutoff factor, write
$\chi_j(q)=\chi(Aq/s_j)$ and $D_\omega=L_\omega+E_\omega$.
Absorbing the scalar jet $z=z_{j+1}$ replaces the slow covariance by
\[
 \frac{\chi_j(q)}{Z_j(1+z\chi_j(q))L_\omega(q)+E_\omega(q)}.
\]
Since $\chi_j=1$ on the support of $\chi_{j+1}$, subtracting
the next slow covariance gives the shell
\[
 \frac{\chi_j(q)-\chi_{j+1}(q)}
      {Z_j(1+z\chi_j(q))L_\omega(q)+E_\omega(q)},
 \qquad Z_{j+1}=Z_j(1+z).
\]
The reference step uses the same formula with $E_\omega=0$.
In particular, the common cutoff-dependent scalar multiplier
allowed in \cref{lem:running-curvature} is retained, rather
than treated as a decaying curvature error. These are the
scalar-absorption and shell prescriptions of
\cite[Eqs.~(6.36), (6.43)--(6.45)]{GMT20}.
The reference construction in \cite{BM04} also retains this
cutoff-dependent scalar factor; its Theorem~3.1 is the
constant-history cancellation input cited by GMT20.
The change of variables $p=Aq$, with constant field rescaling,
identifies the continuum integrations and their zero-momentum
jets. Induction over the integrations identifies their scalar
normalization recursion as well.
\Cref{lem:all-order-reference-sampling}
then identifies the lattice functional with it. All fixed
Jacobian factors are absorbed in the units of $t$ and the
constant $C$. Thus \cref{eq:reference-constant-history} is a
precisely specified external cancellation theorem, whereas
\cref{eq:reference-history-derivatives} has been proved above.
For the physical real coefficients, complex conjugation exchanges
the two Fermi labels and preserves the ordered quartic monomial
(it exchanges two pairs of fields). The quartic coefficient is
therefore real; the following application uses real histories.

\begin{corollary}[Summable beta functions on slowly varying histories]
\label{cor:conditional-beta-summability}
Under the preceding hypotheses, use the reference normalization
recursion for a real quartic history. Suppose, in addition,
\[
 u_k\leq U2^{-\theta k},\qquad
 |\lambda_{k+1}-\lambda_k|\leq M2^{-\theta k},\qquad
 |\lambda_k|\leq\epsilon,
 \quad U,M\leq\epsilon\,.
\]
Then the full quartic beta increment obeys
\begin{equation}\label{eq:conditional-beta-summability}
 |B_j^\lambda|\leq
 C[\delta+\epsilon(U+M+s_0)+\epsilon^2]2^{-\theta j}\,.
\end{equation}
The nonreference quadratic beta functions obey the same bound
without the terms $\epsilon M+\epsilon^2$. In particular these
all-orders increments are absolutely summable in $j$ on the
stated class of histories.
\end{corollary}
\begin{proof}
The idea is to compare the given history with the constant history
equal to its value $\lambda_j$ at the current scale: the external
cancellation \cref{eq:reference-constant-history} controls the
reference beta function on the constant history, the derivative
bounds of \cref{lem:reference-history-derivatives} control the
change caused by the variation of the history, which is small
because the history varies slowly and the derivatives have short
memory, and \cref{prop:conditional-beta-comparison} controls the
difference between the full and the reference beta functions.
For $k\leq j$, telescoping gives
$|\lambda_k-\lambda_j|\leq\sum_{l=k}^{j-1}M2^{-\theta l}
\leq C_\theta M2^{-\theta k}$.
Interpolate from the given history to the constant history
$\lambda_j$. \Cref{lem:reference-history-derivatives}, with its
normalization recursion included, bounds the difference of beta
functions by
\[
 C\epsilon M\sum_{k=0}^j
       2^{-\vartheta(j-k)}2^{-\theta k}
       \leq C\epsilon M2^{-\theta j}\,.
\]
The constant-history bound gives the additional
$C\epsilon^2 2^{-\theta j}$. The equality of the reference
functionals permits using these continuum estimates here.
For the full-reference difference apply
\cref{eq:conditional-beta-comparison,eq:running-curvature-gain}. The geometric sums satisfy
\[
 \sum_{k=0}^j2^{-\vartheta(j-k)}2^{-\theta k}
       \leq C2^{-\theta j},\qquad
 \sum_{k=0}^j2^{-\vartheta(j-k)}2^{-(1-\eta)k}
       \leq C2^{-\vartheta j}\,,
\]
because $\theta<\vartheta<1-\eta$. This proves
\cref{eq:conditional-beta-summability} and its quadratic
version. Summing the final geometric bound proves absolute
summability, with constants independent of the last scale.
\end{proof}

The next subsubsection constructs a history satisfying the recursion,
when the initial quadratic coefficients may be tuned independently.
Matching those coefficients to the physical model, by choosing its
dressed Fermi points and velocities, is a separate step. Source
normalization and comparison require further estimates.

\subsubsection{Construction of the critical flow at fixed dispersion}
\label{sec:critical-flow}

The relevant recursion $\nu_{j+1}=2\nu_j+F_j^\nu$ selects its
initial value through decay at infinite scale. We construct that
value, and the marginal quadratic data, by the backward tail map
of \cite[\S6.4.2--6.4.4]{GMT20}. The required uniform bounds on
its forcing and differences are proved below for our endpoints;
matching the resulting data to the original model is deferred to
\cref{sec:physical-matching}.

Fix the dispersion, adapted cutoff and initial kernel space of
\cref{sec:conditional-bulk-comparison}. Restrict the
kernel space to translation-covariant balanced kernels with the
physical reality symmetry; complex estimates use its complexification.
The data are an
initial quartic coefficient $t$ and an initial remainder
$R_0=(1-\cL)V_0$, where
$\fE(V_0)\leq\delta$. Its representative $V_0$ is measured
in \cref{eq:infrared-endpoint-norm}; its local coefficients are
not used a second time. In this subsubsection the initial local
quadratic coefficients are free parameters to be determined by
the infrared boundary condition. Thus this is a construction of
critical quadratic data for a fixed Gaussian dispersion, not yet
an identification of the quadratic data of the original model.

Use $r=1/8$ for the history norm. The strict inequalities
$2r<\theta=1/2<\vartheta=3/4$ leave room for the following
real-variable derivative estimate: the reference cancellation is
only asserted for real couplings, and differentiating it loses
half its decay exponent.

\begin{lemma}[A differentiated constant-history cancellation]
\label{lem:constant-history-derivative}
Let $f_j(t)=\widehat B_j^R(t,\ldots,t)$. For sufficiently small
$\epsilon$, uniformly for real $|t|\leq\epsilon$,
\begin{equation}\label{eq:constant-history-derivative}
 |f_j'(t)|\leq C\epsilon 2^{-j/4}\,.
\end{equation}
For complex histories in a fixed larger polydisc, and any direction
$h=(h_0,\ldots,h_j)$,
\begin{equation}\label{eq:reference-mixed-derivative}
 \big|D_h\partial_{\lambda_k}\widehat B_j^R(\lambda)\big|
   \leq C\|h\|_\infty2^{-\vartheta(j-k)},\qquad k\leq j\,.
\end{equation}
All constants are independent of $j$.
\end{lemma}
\begin{proof}
Apply \cref{lem:reference-history-derivatives} with a
fixed enlargement of the smallness parameter, leaving a complex
neighborhood of radius $\epsilon$. Summing its derivative bound
over $k$ gives $|f_j'|\leq C\epsilon$ on that neighborhood.
Cauchy's estimate in the common scalar parameter gives
$|f_j''|\leq C$, uniformly in $j$. For real $|t|\leq\epsilon$
put $h_j=\epsilon2^{-j/4}$. Both $t$ and $t+h_j$ belong to
the real interval $[-2\epsilon,2\epsilon]$. Taylor's formula
and \cref{eq:reference-constant-history} give
\[
 |f_j'(t)|\leq
 \frac{|f_j(t+h_j)|+|f_j(t)|}{h_j}
       +\tfrac12 h_j\sup_{|s|\leq2\epsilon}|f_j''(s)|
 \leq C\epsilon2^{-j/4}\,.
\]
For \cref{eq:reference-mixed-derivative}, apply Cauchy's estimate
to the analytic function
$\zeta\mapsto\partial_{\lambda_k}\widehat B_j^R(\lambda+\zeta h)$
on a circle of radius $c\epsilon/\|h\|_\infty$.
The first derivative estimate costs $C\epsilon
2^{-\vartheta(j-k)}$, so division by this radius proves the
claim. This argument has no factor equal to the number of scales.
\end{proof}

Let $q_j=(\nu_j,v_j)$ collect the nonreference quadratic
coefficients: $\nu_j$ is the constant jet divided by $s_j$,
and $v_j$ is the linear jet. Both include their Fermi labels and
are in units of the normalized fields. Use the maximum norm on
this finite-dimensional space. Set
\[
 \|q\|_r=\sup_{j\geq0}2^{rj}|q_j|,\qquad
 \|\lambda\|_{\mathrm{inc},r}
   =\sup_{j\geq0}2^{rj}|\lambda_{j+1}-\lambda_j|,
 \qquad c_r=(1-2^{-r})^{-1}\,.
\]
The norm $\|q\|_r$ requires the nonreference quadratic
coefficients to decay geometrically in the scale, which is the
boundary condition at infinite scale discussed above; the increment
norm requires the quartic coefficient to converge geometrically,
without prescribing its limit. When $\lambda_0=t$, these norms imply
$\sup_j|\lambda_j|\leq |t|+c_r\|\lambda\|_{\mathrm{inc},r}$,
since $\sum_{k\geq0}2^{-rk}=c_r$; thus a small ball in these norms
is contained in the domain of small histories on which the
estimates of \cref{sec:conditional-bulk-comparison} hold.
All histories below obey the reality relations of the original
kernels, and the quartic history is real.

Write the exact flow equations, after eliminating the reference
normalization by \cref{eq:reference-normalization-recursion}, as
\begin{equation}\label{eq:critical-flow-recursion}
 \nu_{j+1}=2\nu_j+F_j^\nu(q,\lambda),\qquad
 v_{j+1}=v_j+F_j^v(q,\lambda),\qquad
 \lambda_{j+1}=\lambda_j+F_j^\lambda(q,\lambda)\,.
\end{equation}
Here $F_j^\lambda=B_j^\lambda$. The quadratic functions include
the generated jets and the normalization drift of the propagated
jet. In particular the exact propagation factors before extracting
this drift are $2/(1+z_{j+1})$ and $1/(1+z_{j+1})$.
Their differences from $2$ and $1$ are included in $F_j^\nu$
and $F_j^v$. They are bounded by $C\epsilon^2|q_j|$ on the
histories under consideration. This convention gives precisely
\cref{eq:critical-flow-recursion}, without omitting a wave-function
factor in the relevant recursion.

\begin{lemma}[Two-history estimates for the full beta functions]
\label{lem:two-history-beta}
On histories satisfying
\[
 \sup_j|\lambda_j|\leq\epsilon,\qquad
 \|q\|_r,\|\lambda\|_{\mathrm{inc},r}\leq\rho\leq\epsilon,
\]
the beta functions satisfy
\begin{equation}\label{eq:flow-beta-size}
 \sup_j2^{rj}\max_{\#=\nu,v,\lambda}|F_j^\#|
       \leq C[\delta+\epsilon(\rho+s_0)+\epsilon^2]\,.
\end{equation}
For two such histories with the same initial $t$ and $R_0$, put
\[
 d_r=\max\{\|q-q'\|_r,
                   \|\lambda-\lambda'\|_{\mathrm{inc},r}\}\,.
\]
Then
\begin{equation}\label{eq:flow-beta-lipschitz}
 \sup_j2^{rj}\max_\#
 |F_j^\#(q,\lambda)-F_j^\#(q',\lambda')|
       \leq C(\epsilon+s_0+\delta/\epsilon)d_r\,.
\end{equation}
The constants are uniform for sufficiently small $s_0$, after
measuring the initial kernels in the rescaled endpoint norm.
\end{lemma}
\begin{proof}
The proof of \cref{cor:conditional-beta-summability}
works with $r$ in place of its history exponent $\theta$:
the constant-history estimate still decays as $2^{-j/2}$,
and the convolutions use $r<\vartheta<1-\eta$.
The added quadratic normalization drift has norm at most
$C\epsilon^2\rho$. This proves \cref{eq:flow-beta-size}.

First estimate the reference contribution to the difference.
For $k\leq j$ put $a_k=\lambda_k-\lambda_j$ and
$w_k(s)=\lambda_j+sa_k$. The exact identity
\begin{equation}\label{eq:reference-history-telescope}
 \widehat B_j^R(\lambda)=f_j(\lambda_j)
  +\int_0^1\sum_{k=0}^j
       \partial_{\lambda_k}\widehat B_j^R(w(s))a_k\,\d s\,,
\end{equation}
and its primed version separate the constant part from the
variation of the history. We have
\[
 |a_k|,|a'_k|\leq c_r\rho2^{-rk},\qquad
 |a_k-a'_k|\leq c_r d_r2^{-rk},\qquad
 \|w(s)-w'(s)\|_\infty\leq c_r d_r\,.
\]
The last inequality uses the common value at scale zero.
The difference of constant terms is at most
$C\epsilon d_r2^{-j/4}$ by
\cref{lem:constant-history-derivative}.
In the integral, the difference of the $a_k$ factors costs
$C\epsilon d_r2^{-rj}$ by
\cref{eq:reference-history-derivatives}. The difference of
the derivative factors costs $C\rho d_r2^{-rj}$ by
\cref{eq:reference-mixed-derivative}. In both cases the sum is
bounded using
\[
 \sum_{k=0}^j2^{-\vartheta(j-k)}2^{-rk}
       \leq (1-2^{-(\vartheta-r)})^{-1}2^{-rj}\,.
\]
Thus the reference part of \cref{eq:flow-beta-lipschitz}
has coefficient $C\epsilon$.

For the remaining functions put
$H_j^\lambda=F_j^\lambda-\widehat B_j^R$ and
$H_j^{\nu,v}=F_j^{\nu,v}$. They are analytic also for complex
histories. \Cref{prop:conditional-beta-comparison},
the curvature estimate and the normalization drift give
\begin{equation}\label{eq:nonreference-complex-bound}
 |H_j^\#|\leq C(\delta+\epsilon^2+\epsilon s_0)2^{-rj}\,,
\end{equation}
on a fixed enlargement of the domain
$\|q\|_r\leq\epsilon$, $\|\lambda\|_\infty\leq\epsilon$.
Here no increment condition on $\lambda$ is needed. The
reference recursion is analytic for arbitrary bounded small
quartic histories by \cref{lem:reference-history-derivatives}.

Apply Cauchy's estimate along a joint direction $(h^q,h^\lambda)$,
with circle radius
$c\epsilon/(\|h^q\|_r+\|h^\lambda\|_\infty)$.
\Cref{eq:nonreference-complex-bound} gives
\[
 |D_hH_j^\#|\leq
 C(\delta/\epsilon+\epsilon+s_0)
  (\|h^q\|_r+\|h^\lambda\|_\infty)2^{-rj}\,.
\]
Integrate this estimate along the segment between the histories
and use $\|\lambda-\lambda'\|_\infty\leq c_r d_r$.
This proves \cref{eq:flow-beta-lipschitz}. The complex
enlargements used here have fixed multiplicative size; all
smallness thresholds in the preceding lemmas are reduced once
to allow them. Uniformity in $s_0$ follows from their rescaled
shell estimates and the factors in
\cref{eq:infrared-endpoint-norm}.
\end{proof}

\begin{theorem}[Critical flow with freely tuned quadratic data]
\label{thm:critical-flow-fixed-dispersion}
Fix a dispersion satisfying the uniform hypotheses of
\cref{sec:infrared-scales} and the conjugation identity
$D(-k)=\overline{D(k)}$, with the paired Fermi points and
adapted reference model of \cref{sec:conditional-bulk-comparison}.
Choose a sufficiently small
adapted initial scale $s_0$, then a sufficiently small $\epsilon$.
There is $c_*>0$ such that, for real $|t|\leq c_*\epsilon$
and $\fE(V_0)\leq\delta\leq\epsilon^2$, the recursion
\cref{eq:critical-flow-recursion} has a unique solution in
\[
 X=\{(q,\lambda):\lambda_0=t,\quad
       \|q\|_r\leq\rho,\quad
       \|\lambda\|_{\mathrm{inc},r}\leq\rho\},
       \qquad \rho=c_*\epsilon\,.
\]
In particular the solution determines the initial quadratic
coefficients $q_0$; all its quadratic coefficients tend to zero
and its quartic coefficient has a limit $\lambda_\infty$.
More quantitatively, with
$a=\delta+s_0|t|+t^2$,
\begin{equation}\label{eq:critical-flow-estimate}
 |q_j|+|\lambda_{j+1}-\lambda_j|
       +|\lambda_j-\lambda_\infty|
       \leq Ca2^{-rj},\qquad
 |\lambda_\infty-t|\leq Ca\,.
\end{equation}
The reference normalizations are well defined and satisfy
$|Z_{j+1}/Z_j-1|\leq C(|t|+a)^2$. The effective interactions
are the convergent renormalized tree sums along this solution.
\end{theorem}
\begin{proof}
Use coordinates $m_j=\lambda_{j+1}-\lambda_j$ and
$\lambda_j=t+\sum_{k<j}m_k$, so that $X$ is a closed Banach ball.
Choose $(1+c_r)c_*<1/4$. Every history in $X$ then lies in the
small-history polydisc, with $|z_{j+1}|\leq C\epsilon^2$.
The absorption denominators stay away from zero, and
\cref{lem:running-curvature} verifies the shell hypotheses on
this entire ball.

Solving the quadratic recursions backwards from a terminal scale
and using their decay in $X$ gives the tail map
\begin{align*}
 (T\nu)_j&=-\sum_{k=j}^\infty2^{j-k-1}F_k^\nu,
 & (Tv)_j&=-\sum_{k=j}^\infty F_k^v,\qquad
 (Tm)_j=F_j^\lambda.
\end{align*}
Conversely every fixed point solves the recursions with
$\lambda_0=t$. From the forcing norm
$\sup_j2^{rj}|F_j|$ to the $(q,m)$ norm, the displayed linear
operators have norm at most
$K_r=\max\{(2-2^{-r})^{-1},c_r,1\}$, by geometric summation.
Thus \cref{lem:two-history-beta} gives
\[
 \|Tx\|\leq K_rC[\delta+\epsilon(\rho+s_0)+\epsilon^2],
 \qquad
 \|Tx-Tx'\|\leq K_rC(\epsilon+s_0+\delta/\epsilon)\|x-x'\|.
\]
Choose $s_0$ small compared with $c_*$, then $\epsilon$ small,
and use $\delta\leq\epsilon^2$. The first bound is at most
$\rho/2$ and the second has coefficient at most $1/2$.
The map preserves the reality relations. Banach's theorem gives
its unique fixed point in $X$; every decaying solution satisfies
the same tail identities and hence is this fixed point.

For the sharper bound, apply $T$ to the flat history
$(q,m)=(0,0)$, with quartic coefficient $t$. Initial remainders,
curvature marks and the reference cancellation cost respectively
$C\delta2^{-\vartheta j}$, $Cs_0|t|2^{-\vartheta j}$ and
$Ct^2 2^{-j/2}$. Its distance from its image is therefore at most
$Ca$, where $a=\delta+s_0|t|+t^2$. The contraction bounds the
fixed point's distance by $2Ca$. Summing the quartic increments
gives \cref{eq:critical-flow-estimate}. The two-endpoint bound
for the reference scalar jet gives
$|Z_{j+1}/Z_j-1|\leq C(|t|+a)^2$. Finally, the convergent tree
expansion along this history recovers the successive Gaussian
integrations, as in \cite[\S6.2--6.4]{GMT20}.
\end{proof}

\Cref{thm:critical-flow-fixed-dispersion} constructs the
bulk flow for tuned quadratic data without assuming the decay
of an unknown solution. At fixed cutoff and degree weight,
\cref{thm:bulk-first-step} supplies the required smallness
of $t$ and $R_0$ for all sufficiently large $\beta$. We next match
the initial quadratic data to the physical model. Continuity in
the dispersion, rather than differentiability of the entire
infinite history in moving Fermi coordinates, suffices for this
existence problem.

\subsubsection{Matching the physical quadratic data}
\label{sec:physical-matching}

The preceding theorem tunes the quadratic data for a fixed
dispersion. To obtain the actual dimer model, we move its Fermi
points and differentials and compensate this change in the
interaction, following \cite[\S6.1]{GMT20}. The physical jets
must then equal the tuned ones. We solve this finite-dimensional
matching problem by Brouwer's theorem. The additional work for our
coefficient class is to verify that the first integration remains
small and continuous under these moving cutoffs.

\paragraph{The dressed dispersions}
Let $p^0$ be the plus zero of the bare dispersion $\mu$, and let
$\ell^0=\nabla\mu(p^0)\in\C^2$. Identify $\C$ with $\R^2$.
The inverse function theorem gives a smooth
local inverse $\pi$ of $\mu$, with $\pi(0)=p^0$: for $w$ in a
neighborhood of zero, $\pi(w)$ is the unique point near $p^0$
at which $\mu$ takes the value $w$. We parametrize the position
of the plus Fermi point and the differential there by the
following six real parameters:
\begin{equation}\label{eq:physical-parameter-chart}
 y=(y_0,y_1)\in\C\times\C^2,\qquad
 p(y)=\pi(s_0y_0),\qquad
 \ell(y)=\nabla\mu(p(y))-y_1\,.
\end{equation}
The constant jet is measured in units of $s_0$. Restrict to
$|y|\leq b_*s_0$, for fixed small $b_*$, using the Euclidean norm.
Set $p^+=p(y)$, $p^-=-p(y)$ and
\[
 L_{y,+}(q)=\ell(y)\cdot q,\qquad
 L_{y,-}(q)=-\overline{\ell(y)}\cdot q,\qquad
 A_yq=(\Re L_{y,+}(q),\Im L_{y,+}(q))\,.
\]
Use the radial Gevrey function of \cref{sec:infrared-scales}
to define $\chi_y(q)=\chi(A_yq/s_0)$; this is the adapted cutoff
of \cref{sec:infrared-scales} for the pair $(p(y),\ell(y))$, so
that its level sets are the ellipses on which the linearized
dispersion has constant modulus. The two translated
supports are disjoint. On these supports put
\begin{align}
 J_{y,+}(q)&=s_0y_0+y_1\cdot q,
 &J_{y,-}(q)&=\overline{J_{y,+}(-q)},\notag\\
 X_y(k)&=\sum_{\omega=\pm}\chi_y(k-p^\omega),
 &D_y(k)&=\mu(k)-\sum_{\omega=\pm}
                   \chi_y(k-p^\omega)J_{y,\omega}(k-p^\omega)\,.
                                      \label{eq:physical-dispersion}
\end{align}
Extend each summand periodically by zero outside its chart.
Since $\chi_y=1$ near zero,
\[
 D_y(p(y))=0,\qquad \nabla D_y(p(y))=\ell(y),
\]
and conjugation gives the minus zero and differential. Away from
the two charts $D_y=\mu$. The compensating slow-field interaction
is $J_y$ without a cutoff factor; the exact Gaussian identity is
verified in \cref{lem:moving-first-integration}.

\begin{lemma}[A uniform family of dressed Gaussian dispersions]
\label{lem:dressed-dispersion-family}
For sufficiently small $s_0,b_*$, the dispersions $D_y$ have exactly
the two simple zeros $p^\pm(y)$, have differentials $L_{y,\pm}$
there, and satisfy $D_y(-k)=\overline{D_y(k)}$.
Their adapted shell covariances obey the hypotheses and estimates
of \crefrange{sec:infrared-scales}{sec:critical-flow}
uniformly for $|y|\leq b_*s_0$. In particular the curvature gain
in \cref{eq:running-curvature-gain} has a common constant.
The bare zeros lie in the interior of $\{X_y=1\}$, and
\begin{equation}\label{eq:moving-fast-covariance}
 \widehat g_{\mathrm f,y}(k)=\frac{1-X_y(k)}{\mu(k)}\,,
\end{equation}
is a smooth Gevrey multiplier, continuous in $y$, with uniform
spatial and Gram bounds at the initial scale.
\end{lemma}
\begin{proof}
The point of the lemma is that the whole multiscale analysis of
\crefrange{sec:infrared-scales}{sec:critical-flow} depends on
the dispersion only through a finite list of quantitative
hypotheses (Gevrey bounds, invertible differentials at two
separated simple zeros, a lower bound away from them, and the
curvature bounds of \cref{lem:running-curvature}), and that
these hold with common constants over the family $D_y$.

The chart gives $|p(y)-p^0|\leq Cs_0|y|$ and
$|\ell(y)-\ell^0|\leq C|y|$, since $\pi$ is smooth and
$\nabla\mu$ is Lipschitz. Hence $A_y$ and its inverse
are uniformly bounded. All supports stay in two fixed disjoint
charts, and the bare zeros are at adapted distance $O(s_0|y|)$
from their centers. For small $b_*,s_0$ they lie in a common
inner region where the cutoff is identically one.

The main point is the curvature bound. In the plus chart write
\[
 R_{\mu,p}(q)=\mu(p+q)-\mu(p)-\nabla\mu(p)\cdot q\,,
\]
for the Taylor remainder of the bare dispersion at $p=p(y)$.
Then, exactly,
\[
 D_y(p+q)-L_{y,+}(q)
   =R_{\mu,p}(q)+(1-\chi_y(q))J_{y,+}(q)\,.
\]
Indeed, $\mu(p+q)=\mu(p)+\nabla\mu(p)\cdot q+R_{\mu,p}(q)
=s_0y_0+(\ell(y)+y_1)\cdot q+R_{\mu,p}(q)$, while
$J_{y,+}(q)=s_0y_0+y_1\cdot q$, so subtracting
$\chi_y(q)J_{y,+}(q)$ and then $L_{y,+}(q)=\ell(y)\cdot q$
leaves the displayed expression.
The first term is $O(|q|^2)$ uniformly, with Gevrey bounds,
as required of a curvature term. The second vanishes
when $|A_yq|\leq s_0$, where the cutoff equals one, so it is
invisible on all shells but the outermost. On the transition
annulus, $|q|\asymp s_0$ and
$|J_{y,+}(q)|\leq Cs_0|y|\leq Cb_*s_0^2$; that is, there the
polynomial is as small as a genuine quadratic remainder.
After rescaling $q=s_0u$, all its derivatives satisfy the
Gevrey bounds of a function of size $Cs_0^2$.
On smaller shells it vanishes identically. Thus the curvature
has the uniform rescaled bounds used in
\cref{lem:running-curvature}, including its differentiated
bounds. The minus chart follows by conjugation. For small
$s_0$ this also gives $|D_y(p^\omega+q)|\geq c|q|$
throughout each cutoff support, because the linear term
dominates the curvature there. Outside both supports $D_y=\mu$;
the only zeros of $\mu$ have already been included in the
inner regions, where we have just seen that $D_y$ vanishes only
at $p^\pm$. This proves the zero and differential assertions.
Conjugation follows directly from the paired definitions.

The denominator in \cref{eq:moving-fast-covariance} is used
only away from the common neighborhoods of the bare zeros.
There it is bounded away from zero at the fixed scale. Extend
the multiplier as zero through those neighborhoods. Differentiating
the rescaled cutoffs gives uniform Gevrey bounds, and Fourier
integration by parts gives a common stretched exponential
envelope. The equal-modulus Fourier factorization gives the
Gram bounds. These multipliers and every fixed spatial matrix
element depend continuously on $y$ by dominated convergence.
The same arguments apply to every fixed infrared shell.
All constants needed for shell power counting are uniform as
$s_0$ decreases: in the exceptional outer shell the estimate
above is $O(s_0^2)$, and on the inner shells only the bare
Taylor remainder remains. Fixed-scale spatial constants may
depend on the chosen $s_0$, which is held fixed before $\beta$
is increased.
\end{proof}

We next redo the first integration of \cref{sec:first-integration}
with the slow modes defined by the moving cutoff $X_y$ instead
of the fixed cutoff used there. This produces an interaction
$U_{\beta,y}$ on the slow fields, and the three quantities that
the matching must control are its \emph{physical quadratic data},
namely the value and the gradient at the Fermi point of the
Fourier symbol $\sigma_y$ of its two-field part (the first moves
the zero of the effective dispersion, the second its
differential), and its localized quartic coefficient $t_y$.

\begin{lemma}[The first integration with moving slow fields]
\label{lem:moving-first-integration}
At fixed sufficiently small $s_0,\epsilon,b_*$, the first
integration using $g_{\mathrm f,y}$ gives a bulk interaction
$U_{\beta,y}$, uniformly for $|y|\leq b_*s_0$. Its kernels
satisfy the endpoint estimates with a fixed reserve in spatial
weight and field-degree weight. There is
$d_\beta=Ce^{-c\beta^{2/5}}\to0$, with constants allowed to
depend on these fixed choices, such that
\begin{equation}\label{eq:moving-input-smallness}
 \fE_y(U_{\beta,y})+|t_y|
  +s_0^{-1}|\sigma_y(p(y))|
  +|\nabla\sigma_y(p(y))|\leq d_\beta\,.
\end{equation}
Here $\sigma_y$ is its physical quadratic Fourier symbol,
$t_y$ its localized quartic coefficient, and
$\fE_y$ means \cref{eq:infrared-endpoint-norm} after
splitting into the Fermi fields centered at $p^\pm(y)$.
Every fixed-position coefficient and each of the jets in
\cref{eq:moving-input-smallness} is continuous in $y$.

Replacing the remaining bare covariance $X_y/\mu$ by
$X_y/D_y$ adds precisely
\begin{equation}\label{eq:physical-compensating-term}
 \sum_{\omega,x}\psi_{b,\omega}(x)
             J_{y,\omega}(-i\nabla)\psi_{w,\omega}(x)\,,
\end{equation}
to $U_{\beta,y}$, and changes only an explicitly retained
Gaussian determinant factor. The resulting representation is
an exact reorganization of the original model.
\end{lemma}
\begin{proof}
The free covariance splits as
$\mu^{-1}=\widehat g_{\mathrm f,y}+X_y/\mu$.
The proof of \cref{thm:bulk-first-step} applies uniformly
to the fast part. To check its two cutoff-dependent ingredients,
the transformed covariance obeys the exact multiplier identity
$\mu\widehat g_{\mathrm f,y}=1-X_y$, and an external cutoff
$\zeta_y$ equal to one on the support of $X_y$ can be chosen
in slightly larger adapted balls. Its convolution norm is $O(1)$,
while that of $\mu\zeta_y$ is $O(s_0)$, since these balls stay
within $O(s_0)$ of the bare zeros. These are exactly the
interior Gram and external-leg gains in
\cref{prop:conditional-first-step}.
The \texttt{fat}-pattern determinant denominator still uses the original
Kasteleyn operator; no determinant estimate for $D_y$ is needed.
All remaining fast contractions have a common weighted summable
envelope by \cref{lem:dressed-dispersion-family}.
The same strictly negative \texttt{fat}-node exponent and small \texttt{thin}-node
sum therefore apply uniformly in $y$.

Choose the first-step degree weight larger than the one required
by $\fE_y$, and choose its spatial exponent with a
strict reserve over $2\kappa\sqrt{s_0}$. The proof of
\cref{thm:conditional-renormalized-trees} then gives the
endpoint bound. The first two spatial moments also bound the
displayed quadratic jets and the quartic localization, after
enlarging the fixed constant by $s_0^{-1}$.
The coefficients before inserting the phases $e^{\pm ip^\omega\cdot x}$ are continuous
in $y$: truncate the absolutely convergent first-step trees,
their spatial sums and their input patterns, use continuity of
each remaining covariance, and remove the truncations using
the uniform envelopes. These phases are continuous at each
fixed set of coordinates. The moment reserves give the same
conclusion for the Fourier jets. We do not assert continuity
in a supremum norm over all translates of these phases, which
would be false for a varying nonzero $p^\omega$.
The paired cutoff $X_y$ is real and even. Hence the fast
covariance is real in position space, and the same reality
argument as in \cref{prop:infrared-input} gives
$\sigma_y(-k)=\overline{\sigma_y(k)}$ and real $t_y$.
Differentiating this identity supplies the required opposite
sign in the conjugation relation for the two linear jets.

For the algebraic assertion first take a finite Fourier grid
on which the relevant bare and dressed eigenvalues are nonzero.
Recall that in the convention \cref{eq:normalized-grassmann}
a Gaussian Grassmann measure with covariance $C$ has density
$e^{\psi_bC^{-1}\psi_w}$ divided by $\det C^{-1}$, so that the
quadratic form in the exponent is the \emph{inverse} covariance.
The slow covariance $X_y/\mu$ thus has quadratic form
$\mu/X_y$, and the dressed one $X_y/D_y$ has quadratic form
$D_y/X_y$. On a slow mode in the $\omega$ chart, where
$X_y=\chi_y(\cdot-p^\omega)$, the definition
\cref{eq:physical-dispersion} gives the exact identity of
quadratic forms
\begin{equation}\label{eq:exact-slow-quadratic-identity}
 \frac{\mu(p^\omega+q)}{\chi_y(q)}
   =\frac{D_y(p^\omega+q)}{\chi_y(q)}+J_{y,\omega}(q),
                         \qquad \chi_y(q)>0\,,
\end{equation}
simply because $\mu-D_y=\chi_yJ_{y,\omega}$ on that chart.
Multiplying by $\psi_{b,\omega}\psi_{w,\omega}$ at that mode and
summing over the slow modes, the term $J_{y,\omega}(q)$ is the
Fourier representation of \cref{eq:physical-compensating-term}.
The ratio of the normalized Gaussian densities is consequently
the exponential of \cref{eq:physical-compensating-term},
times $\prod_{X_y(k)>0}D_y(k)/\mu(k)$, the latter being the
ratio of the two normalizing determinants. Multiplying by the
original Gaussian determinants retains this factor exactly.
Modes with $X_y=0$ have $D_y=\mu$ and contribute the factor one.
At a singular mode use the unnormalized Berezin integral:
\cref{eq:exact-slow-quadratic-identity} is an equality of its
exponents and never divides by an eigenvalue. Equivalently,
the determinant-weighted identities extend from nonsingular
parameters by continuity. A cutoff-zero mode is omitted from
the slow field space. This proves the identity also at those
degeneracies, with zero modes left as polynomial integrations.
The bulk kernels follow by the already dominated first-step
limit. This algebraic statement does not estimate the final
finite-volume integral or any projected partition function.
\end{proof}

For each $y$ we now have a dispersion $D_y$ and an interaction
$U_{\beta,y}$ with quartic coefficient $t_y$, to which
\cref{thm:critical-flow-fixed-dispersion} applies. That theorem
returns the initial quadratic jets $q_0^*(y)$ that the critical
flow for $D_y$ requires; these are to be compared with the
physical jets of $U_{\beta,y}$. Since the comparison will be
closed by a fixed-point argument, what we need from the tuned
data is a uniform bound and continuity in $y$. The latter is
delicate, because the tuned data are defined through an
infinite recursion whose differentiability in $y$ we do not
control; the proof below obtains continuity from uniqueness of
the fixed point and compactness instead.

\begin{lemma}[Continuity of the tuned bulk data in this family]
\label{lem:critical-data-moving-dispersion}
For all sufficiently large $\beta$, apply
\cref{thm:critical-flow-fixed-dispersion} to $D_y$, with
initial quartic coefficient $t_y$ and remainder
$(1-\cL_y)U_{\beta,y}$. Denote its solution by
$(q^*(y),m^*(y))$, with
$q_j^*=(\nu_j^*,v_j^*)$ and $m_j^*=\lambda_{j+1}^*-\lambda_j^*$.
It exists uniformly for $|y|\leq b_*s_0$ and obeys
\begin{equation}\label{eq:moving-critical-data-size}
 |q_j^*(y)|+|m_j^*(y)|\leq Cd_\beta2^{-rj}\,.
\end{equation}
Each coordinate, in particular $q_0^*(y)$, and
$\lambda_\infty^*(y)$ are continuous in $y$.
\end{lemma}
\begin{proof}
The preceding shell and endpoint bounds give a common history ball
and contraction constant. The reference cancellation is uniform too:
put $\Delta_y=|\det A_y|$ and change momentum variables to $A_yq$.
A quartic graph with $N$ vertices has $N-1$ independent loop
momenta, giving Jacobian $\Delta_y^{1-N}$; the scalar quadratic
derivative has factor $\Delta_y^{-N}$. Including the recursive
normalization, this gives
\[
 \widehat B^R_{j,y}(\lambda_0,\ldots,\lambda_j)
   =\Delta_y\widehat B^R_{j,I}
          (\lambda_0/\Delta_y,\ldots,\lambda_j/\Delta_y).
\]
Here $I$ is the canonical linear pair with the same radial cutoff.
\Cref{lem:all-order-reference-sampling} transfers the identity
to the lattice. The uniform bounds on $A_y^{\pm1}$ give common
cancellation constants and a common no-aliasing margin.
For large $\beta$, \cref{thm:critical-flow-fixed-dispersion}
therefore applies uniformly and proves
\cref{eq:moving-critical-data-size}.

For fixed scale, the beta functions are jointly continuous in $y$
and the finite history: truncate tree number, endpoint degrees and
relative coordinates, use continuity of the resulting Fourier
integrals, and remove the truncations with the uniform envelopes.
One may do this in real space, where the phases $e^{\pm ip^\omega\cdot x}$ are pointwise
continuous and have modulus one; no nonresonance assumption enters.

If $y_n\to y$, the uniform bound on $(q_j^*,m_j^*)$ gives a
coordinatewise convergent subsequence. Finite-scale continuity
passes every recursion to its limit, which lies in the same history
ball and hence equals the unique solution for $y$. This identifies
every subsequential limit and proves continuity of each coordinate.
The uniform tail $\sum_{j\geq J}|m_j^*|\leq Cd_\beta2^{-rJ}$
also proves continuity of $\lambda_\infty^*$.
\end{proof}

\begin{theorem}[Matching the quadratic data and the bulk flow]
\label{thm:matched-bulk-flow}
Fix liquid canonical weights. For all sufficiently large
$\beta$, there exist paired dressed Fermi points $p_\beta^\pm$
and dressed differentials $\ell_\beta,-\overline{\ell_\beta}$
such that the zero-source bulk dimer interaction has the
critical flow of \cref{thm:critical-flow-fixed-dispersion}
with its actual initial quadratic data. At the fixed cutoff,
\begin{align*}
 |p_\beta^+-p^0|+|\ell_\beta-\ell^0|&\leq Cd_\beta,
                                                   \\
 |q_j|+|\lambda_{j+1}-\lambda_j|
       +|\lambda_j-\lambda_\infty|&\leq Cd_\beta2^{-rj},
 &|\lambda_\infty|&\leq Cd_\beta,\\
 |Z_{j+1}/Z_j-1|&\leq Cd_\beta^2\,.
\end{align*}
The effective interactions are the convergent source-free bulk
tree expansions of the original model, after the exact Gaussian
reorganization above. Existence of a matching is asserted;
uniqueness or differentiability of its choice is not needed here.
\end{theorem}
\begin{proof}
The moving first integration gives the slow interaction
$U_{\beta,y}+J_y$. Its dimensionless quadratic jets in the plus
chart are
\[
 q_0^{\mathrm{phys}}(y)=
 \big(y_0+s_0^{-1}\sigma_y(p(y)),\,
                    y_1+\nabla\sigma_y(p(y))\big).
\]
There are no cutoff derivatives: the compensating term is exactly
$J_y$ by \cref{eq:exact-slow-quadratic-identity}, and the cutoff
is one near the localization point. Thus the matching equation is
$\Phi(y)=y$, where
\[
 \Phi(y)=\big(\nu_{0,+}^*(y)-s_0^{-1}\sigma_y(p(y)),\,
                      v_{0,+}^*(y)-\nabla\sigma_y(p(y))\big).
\]
The preceding two lemmas make $\Phi$ continuous with
$|\Phi(y)|\leq C_0d_\beta$ on $|y|\leq b_*s_0$.
Choose $M>C_0$ and then $\beta$ so large that $Md_\beta<b_*s_0$.
Brouwer's theorem gives a fixed point in $|y|\leq Md_\beta$.

At this point the physical and tuned plus jets agree, and reality
gives the minus jets. The compensating polynomial changes neither
the quartic coefficient nor the irrelevant remainder, so all
initial data equal the critical ones. The exact Gaussian identity
shows that the convergent trees represent the original model.
\Cref{eq:critical-flow-estimate,eq:moving-input-smallness} gives
the flow and normalization bounds; the parameter chart gives the
Fermi-point and differential bounds.
\end{proof}

\begin{remark}[Why a fixed-point theorem, and why uniqueness is not needed]
A contraction-mapping or implicit-function argument would
require additional quantitative regularity of the tuned data
$y\mapsto q_0^*(y)$ with respect to the dressed dispersion. Since our first integration involves the
\texttt{fat}-pattern determinant estimates, we have only proved
continuity of these data in $y$
(\cref{lem:critical-data-moving-dispersion}), which is exactly
what the Brouwer theorem needs. The price is that the matched
dispersion need not be unique. This causes no difficulty
downstream: every later statement about the dimer law is proved
for an arbitrary fixed matched dispersion, and the finite-volume
law being compared to it does not depend on that choice.
\Cref{thm:canonical-local-limit} below shows that, as a
consequence, the limiting source functional is the same for every
matched dispersion.
\end{remark}

This closes the matching problem for the bulk zero-source flow.
\Cref{sec:bulk-sources,sec:finite-volume-sources} construct the bulk source series and
identify it with normalized canonical correlations. The sharper
Ward-identity cancellations needed for the \GFF are additional
steps; they do not follow from the quadratic matching.

\subsection{Canonical and prescribed-slope local limits}
\label{sec:dimer-local-limits}

We construct the bulk source functional, identify the canonical local
limit, and then pass to prescribed orientation counts.

\subsubsection{Bulk source expansion and spatial decay}
\label{sec:bulk-sources}

We carry the logarithmic edge sources through the matched bulk
flow, following \cite[\S6.2.3]{GMT20}. This constructs an analytic
source functional and spatially decaying candidate edge cumulants.
The torus comparison in \cref{sec:finite-volume-sources} identifies
them with the limiting law. A coarse source-growth estimate suffices
here; the sharp density and mass flows are established in
\cref{sec:sharp-source-flow}.

Fix one of the matched dispersions in
\cref{thm:matched-bulk-flow}. The dispersion and
all zero-source counterterms are held fixed when sources are
varied. Write $e=(x,r)$ for the edge from the black vertex at
$x$ to the white vertex at $x+d_r$, with weight $K_r$.
Here $r$ ranges over the three honeycomb edge types; it is
unrelated to the decay exponent $r=1/8$ used above.
The source is the original logarithmic edge source $A(e)$ of
\cref{it:rg-sources}, which multiplies the weight of a matching
by $\exp(\sum_{e\in M}A(e))$; thus differentiating in $A(e)$
inserts the edge indicator $\one_{\{e\in M\}}$ in the
exact finite-volume identity, and mixed derivatives at $A=0$ of
the logarithm of the generating function are the joint cumulants
of edge indicators.

We use normalized fermion fields at scale $j$. Recall that the
lattice fields are decomposed into their two Fermi components,
with phases $e^{ip^\omega\cdot x}$ for the black fields and
$e^{-ip^\omega\cdot x}$ for the white fields. A source bilinear at the
edge $(x,r)$ therefore acquires, after extracting these phases,
the phase $e^{i(p^\omega-p^{\omega'})\cdot x}$ when it couples
the flavors $\omega,\omega'$. Denote by
$\Gamma_{j,r;\omega\omega'}$ the coefficient of
\begin{equation}\label{eq:local-edge-source}
 \sum_x A(x,r)e^{i(p^\omega-p^{\omega'})\cdot x}
                  \psi_{b,\omega}(x)\psi_{w,\omega'}(x)\,.
\end{equation}
Thus $\Gamma_j$ includes the inverse field-normalization
factor: in the notation $Y_j$ for coefficients in unnormalized
fields, $\Gamma_j=Y_j/Z_j$. No boundedness assumption on
$\Gamma_j$ is made. All four flavor pairs in
\cref{eq:local-edge-source} must be retained, since the
source can carry the momentum $p^\omega-p^{\omega'}$: the
diagonal pairs describe the non-oscillating part of the edge
correlations, and the off-diagonal pairs their part oscillating
at momentum $p^+-p^-$.

\begin{lemma}[Source endpoints and localization]
\label{lem:source-endpoints-localization}
At the fixed matched cutoff the first integration, including
sources, has a common positive spatial norm of the form
\cref{eq:source-marked-norm}, with a reserve for two spatial
moments and any prescribed fermion-degree weight. Its part
containing at least one original interaction endpoint
has norm $O(d_\beta)$ on a smaller fixed source polydisc.
In particular,
\begin{equation}\label{eq:initial-source-coefficient}
 \Gamma_{0,r;\omega\omega'}
       =K_r e^{-ip^{\omega'}\cdot d_r}+O(d_\beta)\,.
\end{equation}
Apart from source-only terms, which are extracted into the
logarithm of the generating function at each step, the only
additional marginal local terms are
\cref{eq:local-edge-source}. The remainder of a two-fermion,
one-source kernel gains one power of the external fermion
momenta. All terms with at least one source and at least two
fermion fields are irrelevant after this subtraction.
\end{lemma}
\begin{proof}
The moving-cutoff integration of \cref{lem:moving-first-integration}
preserves the source-marked norm of \cref{thm:bulk-first-step}.
Every tree in the interacting-minus-free difference contains an
original interaction endpoint; marking it gives $O(d_\beta)$ on a
smaller source polydisc, at all source orders.

The free first integration is a Gaussian integral with quadratic
source vertices. A connected contraction is a cycle with no external
fermion or a path with two. At first order in $A$, the latter is
$K_rA(x,r)\psi_b(x)\psi_w(x+d_r)$; removing the phases $e^{\pm ip^\omega\cdot x}$ and
localizing at $x$ gives \cref{eq:initial-source-coefficient}.
The external cutoffs equal one at the localization momenta. Free
quadratic kernels at higher source orders are retained as initial
endpoints.

The source localization is precisely that of
\cite[\S6.2.3]{GMT20}: all four bilinears are kept, with source
momentum $p^\omega-p^{\omega'}$. Localization replaces a kernel
with one source and two fermion fields by its value at zero
quasi-particle momenta, which is the local term
\cref{eq:local-edge-source}; the remainder, in Fourier space, is
the kernel multiplied by the difference between the phase at the
actual momenta and the phase at zero. For our nonlocal kernels,
whose fermion positions $u,v$ relative to the source are summed
rather than fixed, this remainder is bounded using
\[
 |e^{iq_1\cdot u-iq_2\cdot v}-1|
       \leq |q_1||u|+|q_2||v|\,,
\]
which is paid by the absolute first moment of the source norm:
this is why a reserve for spatial moments was required. The power
counting in \cite[Prop.~3]{GMT20} assigns dimension $2-\fn/2-m$
to a kernel with $\fn$ fermions and $m$ source occurrences
(each fermion field carries dimension $1/2$ and each source
occurrence dimension $1$). For $\fn\geq2,m\geq1$ this dimension
is at most $2-1-1=0$, with equality only for $(\fn,m)=(2,1)$;
so only $(2,1)$ is marginal, and the subtraction of its local
part, which gains one power of the external momenta, makes it
irrelevant. Source-only terms, $\fn=0$, have no fields and
factor out of subsequent integrations; they are extracted into
the logarithm of the generating function.

Finally, sum the relative source and fermion coordinates in the
positive endpoint norm before applying that tree bound. Its degree
reserve pays the field labels, its spatial reserve pays the Taylor
moments, and its source Taylor factorials pay multiple or repeated
marks. These are the endpoint checks required beyond the local
sources of the cited construction.
\end{proof}

\begin{proposition}[A bound on the running source coefficients]
\label{prop:source-flow-bound}
Along the physical bulk flow, the finite-dimensional vector
$\Gamma_j$ satisfies an exact causal affine linear recursion
\begin{equation}\label{eq:source-flow-recursion}
 \Gamma_{j+1}=(1+z_{j+1})^{-1}\Gamma_j
                 +\sum_{k=0}^j M_{j,k}\Gamma_k+h_j\,,
\end{equation}
where the matrices and forcing depend only on the zero-source
flow and the initial source remainders. For sufficiently large
$\beta$,
\[
 \|M_{j,k}\|\leq Cd_\beta2^{-\vartheta(j-k)},\qquad
 |h_j|\leq Cd_\beta2^{-\vartheta j},\qquad \vartheta=3/4\,.
\]
The recursion has a unique solution for its prescribed initial
value, and
\begin{equation}\label{eq:source-growth-bound}
                  |\Gamma_j|\leq C2^{\eta_\beta j},
                  \qquad 0\leq\eta_\beta\leq Cd_\beta\,.
\end{equation}
These estimates do not identify the source critical exponents.
\end{proposition}
\begin{proof}
Differentiate once in the source at zero. Each connected tree has
one source endpoint, either a local $\Gamma_k$ or an initial
remainder, giving the affine recursion
\cref{eq:source-flow-recursion}. The lone local bilinear is
propagated unchanged before field normalization, producing
$(1+z_{j+1})^{-1}$. Source-only constants have already been extracted.

Apply the source-tree estimate and short-memory bound of
\cite[Prop.~3 and \S6.4.9]{GMT20}. The source localization and
initial endpoint hypotheses are verified in
\cref{lem:source-endpoints-localization}; ordinary endpoints
are covered by \cref{thm:conditional-renormalized-trees}.
Every contribution to $M_{j,k}$ has an ordinary endpoint of size
$Cd_\beta$. Marking the local source at scale $k$ gives
$Cd_\beta2^{-\vartheta(j-k)}$. Normalization ratios cost at most
$2^{Cd_\beta^2(j-k)}$, paid by using memory exponent $7/8$ before
weakening it to $\vartheta=3/4$.

An interaction-derived initial source remainder costs $Cd_\beta$.
A free one-source remainder is quadratic with zero local jet:
alone it can only keep that zero jet or contract to a constant.
Thus any contribution from it to $h_j$ also contains an ordinary
endpoint. Its scale-zero mark gives
$|h_j|\leq Cd_\beta2^{-\vartheta j}$. Degree and range truncations
are removed using the uniform endpoint bounds just cited.

The recursion is valid at every finite scale without assuming that
$\Gamma$ is small: its trees have just one source occurrence,
so they are linear in $\Gamma$, and the bounds on the matrices
$M_{j,k}$ come from the ordinary endpoints alone.
Since $|z_{j+1}|\leq Cd_\beta^2$, we have
$|(1+z_{j+1})^{-1}|\leq1+Cd_\beta^2$, and the summable memory
kernel gives, for $B_j=\max\{1,\max_{k\leq j}|\Gamma_k|\}$,
\[
 |\Gamma_{j+1}|
 \leq(1+Cd_\beta^2)B_j+Cd_\beta\sum_{k=0}^j2^{-\vartheta(j-k)}B_j
                +Cd_\beta
 \leq(1+C_1d_\beta)B_j\,,
\]
whence $B_{j+1}\leq(1+C_1d_\beta)B_j$, with $B_0\leq C$ by
\cref{eq:initial-source-coefficient}. Iterating gives
$|\Gamma_j|\leq C(1+C_1d_\beta)^j$, which
proves \cref{eq:source-growth-bound} with
$\eta_\beta=\log_2(1+C_1d_\beta)$. Causality proves uniqueness:
$\Gamma_{j+1}$ is determined by $\Gamma_0,\ldots,\Gamma_j$, and
no division by a source coefficient is involved.
\end{proof}

\begin{remark*}
The exponent $\eta_\beta$ is an upper bound for the anomalous
dimension of the edge observable; its exact value is a nontrivial
function of $\beta$, which in the setting of \cite{GMT17,GMT20}
is identified through Ward identities. We do not need its value
here: the only property of $\Gamma_j$ used below is that it grows
slower than $2^{\vartheta j}$ for $\vartheta=3/4$, which leaves
room in the dimensional bounds of the next lemma. The sharp
information needed for the height \GFF is obtained in
\cref{sec:sharp-source-flow}.
\end{remark*}

\begin{lemma}[Marked comparisons of source trees]
\label{lem:marked-source-comparison}
Consider renormalized expansions with the endpoint and shell bounds
of \cref{thm:conditional-renormalized-trees,lem:source-endpoints-localization}.
Write $\delta$ for a common small bound on the ordinary running
coefficients, and assume the normalization increments are $O(\delta^2)$.
Fix a source order $m$, keeping its Taylor factorials, and write
$\|\cdot\|_{1,j}$ for the positive spatial norm with one coordinate
fixed, including localization moments and a fermion-degree reserve.
A root with $\fn$ fermions and $m$ sources has dimensional factor
$s_j^{2-\fn/2-m}$. Suppose local source coefficients are bounded by
$C2^{\eta j}$, with $\eta$ sufficiently small.

Telescope a difference of two expansions at an initial remainder,
a covariance, or a normalized running coefficient. If the initial
mark has norm $a_{\mathrm{in}}$ and a mark on scale $t$ has
dimensionless norm $a_t$, then, after removing the root dimensional
factor, its summed-kernel bound is at most
\begin{equation}\label{eq:marked-source-convolution}
 C_m2^{m\eta j}\left[
 a_{\mathrm{in}}2^{-\gamma j}
       +\sum_{t=0}^j2^{-\gamma(j-t)}a_t\right],
 \qquad \gamma=1/4.
\end{equation}
This applies to single-scale scalar terms, irrelevant kernels, and
localized beta coefficients; retained local coefficients are compared
separately. For a one-source recursion kernel, factor out its local
source coefficient on scale $k$. A difference in the remaining
kernel has the additional factor $2^{-\gamma(j-k)}$ in the bound
above, with the source-growth factor omitted.
A factor from another required interaction endpoint can be retained
in both estimates. The same marking argument preserves reserved
spatial connecting weights in a fixed-position tree estimate.
The constants are uniform in endpoint degree and range, and the
smallness threshold for ordinary couplings is independent of $m$.
\end{lemma}
\begin{proof}
Use \cite[Prop.~3]{GMT20} with the endpoint extension already
verified in \cref{thm:conditional-renormalized-trees}.
At a nonconstant vertex with $\fn\geq2$ fermions and $m\geq1$
sources, the gain after localization is
\begin{equation}\label{eq:source-gain}
 d(\fn,m)=\fn/2+m-2+\one_{\{(\fn,m)=(2,1)\}}\geq1\,.
\end{equation}
The source-free gains were verified in the bridge theorem. After
field-normalization losses, a common gain $7/8$ remains; the minimal
source cases are $(2,1),(2,2),(4,1)$, and the gain increases with
degree and source order. Choose the couplings small once for these
cases. Local source factors on scales at most $j$ cost at most
$C_m2^{m\eta j}$; their scale growth can be paid per source leaf.

A mark at scale $t$ must be propagated to the root at scale $j$.
Reserve an exponent $\gamma$ along this path, which gives
$2^{-\gamma(j-t)}$. For the recursion-kernel assertion, the
source endpoint at scale $k$ must also be propagated to the root.
Let $s$ be the scale of the first common vertex of these two
paths. Before they meet we reserve $\gamma$ on each path;
after they meet we must reserve $2\gamma$ on the common path.
The resulting factor is exactly
\[
 2^{-\gamma(s-t)}2^{-\gamma(s-k)}
 2^{-2\gamma(j-s)}
       =2^{-\gamma(j-t)}2^{-\gamma(j-k)}\,.
\]
This also covers the case where one marked vertex lies on the
other path, with one of the first two factors equal to one.
The available exponent is at least $3/4$, whereas the common
path uses only $2\gamma=1/2$. Thus at least $1/4$ remains there
for summing the intermediate scales and the locations of the
marks. For example, on a stretch of $d$ scales two distinguished
levels have at most $(d+1)^2$ choices, and
$\sum_{d\geq0}(d+1)^2 2^{-d/8}<\infty$ pays for these choices
while leaving a further exponent $1/8$ in reserve. The separate
paths have still more decay available. Multiplying by the size
$a_t$ of the marked difference and summing over $t$ gives
\cref{eq:marked-source-convolution}; an initial remainder is a
scale-zero mark. In the recursion-kernel bound the local source
coefficient has been factored out, so its scale growth is omitted.

For clarity, telescoping normalization ratios uses
\[
 \frac{Z_a}{Z_b}=\prod_{l=a}^{b-1}(1+z_{l+1})^{-1},\qquad
 (1+z)^{-1}-(1+z')^{-1}
       =\frac{z'-z}{(1+z)(1+z')}.
\]
The term marked at $l$ is bounded by
$C|z_{l+1}-z'_{l+1}|2^{C\delta^2(b-a)}$; the same chain reserve
pays this loss and the choice of $l$. For determinants interpolate
$G_u=(1-u)G+uG'$ in the direct sum of the two Gram spaces and use
\begin{equation}\label{eq:marked-gram-difference}
 \det G'-\det G
 =\int_0^1\sum_{a,b}(G'-G)_{ab}
                       \operatorname{cof}_{ab}(G_u)\,\d u.
\end{equation}
A cofactor is a Gram minor, so the marked entry replaces one ordinary
pair without division by a determinant. The degree reserve pays its
polynomial number of choices.

Sum the relative coordinates of each initial endpoint in its
positive norm before contracting it. The first or second moments
pay localization; its internal spatial tree joins the contraction
tree and pays any reserved connecting weight. Multiple sources in
one endpoint form a block at the initial scale, with their Taylor
factorials retained. At fixed $m$ there are finitely many source
partitions, and all further field choices are paid by the degree
reserve. Marking one of $N$ ordinary endpoints differentiates its
convergent envelope, giving $\sum_N N^2(C\delta)^{N-1}<\infty$.
These estimates permit removal of range and degree truncations.
They also show that a separately required ordinary endpoint keeps
its small factor. This proves the assertions for the stated norms;
fixed-position estimates still require their source-merger factors.
\end{proof}

Let $S_\mathrm{UV}(A)$ be the source-dependent constant from the
first integration, with its value at zero subtracted. At
infrared step $j$, write the extracted source-only term as
\[
 S_j(A)=\sum_{m\geq1}\sum_{e_1,\ldots,e_m}
                 W_{j,m}(e_1,\ldots,e_m)\prod_{i=1}^m A(e_i)\,.
\]
We use symmetric ordered kernels: $W_{j,m}$ includes $1/m!$.
Thus an $m$th derivative, including repeated edges, is
$m!W_{j,m}$. The sums are restricted to the support of $A$.
For edge anchors $x_1,\ldots,x_m$ put
$D(\be)=\max_{a,b}|x_a-x_b|$, with $D=0$ for
$m=1$. Replacing anchors by edge endpoints changes this
distance by a bounded amount.

\begin{lemma}[Summable source-only shells]
\label{lem:source-shell-summability}
Set $\vartheta=3/4$. For all sufficiently large $\beta$, there
are fixed $C,c>0$ such that, for every $m\geq1$ and $j\geq0$,
\begin{equation}\label{eq:source-shell-pointwise}
 |W_{j,m}(\be)|
   \leq C^m2^{-\vartheta j}
               e^{-c\sqrt{s_jD(\be)}}\,.
\end{equation}
For two sources the sharper bound is
\begin{equation}\label{eq:two-source-shell-pointwise}
 |W_{j,2}(e,e')|
   \leq C2^{-\vartheta j}(1+s_0D(e,e'))^{-\vartheta}
                         e^{-c\sqrt{s_jD(e,e')}}\,.
\end{equation}
The initial kernels in $S_\mathrm{UV}$ obey the same estimates
with $j=0$, after changing constants.
\end{lemma}
\begin{proof}
We apply the spatially resolved source-tree estimate
\cite[Eq.~(6.132), and its summation in Eq.~(6.133)]{GMT20}.
The source positions must remain fixed, so we retain the merger
factors in that estimate rather than applying only its summed norm.

First consider local source endpoints. Prune the renormalized
tree to the smallest subtree containing its source endpoints,
and suppress vertices with one surviving child. Let $v_*$
be its first branching vertex, let $k_*$ be its scale, and
write $s_v\geq2$ for the number of surviving children at a
branching vertex $v$. Scales decrease toward the ultraviolet:
$0\leq k_v\leq k_*$ for descendants of $v_*$, and
$k_*\leq j$. Fixing the source positions, rather than summing
them, saves an integration for every merger of source branches.
The source-tree bound retains the corresponding factors
$s_{k_v}^{2(s_v-1)}$ and spatial connecting weights. After
the ordinary vertices and the one-child chains are summed,
a weaker form of that bound is
\begin{equation}\label{eq:pruned-source-bound}
 C^m2^{-\vartheta(j-k_*)}2^{-\vartheta s_{v_*}k_*}
       \prod_{v\ne v_*}2^{-\vartheta(s_v-1)k_v}\,.
\end{equation}
Here is the direct reduction from the cited formula. Denote its
small dimensional reserve by $\epsilon_{\rm G}$, to distinguish
it from the endpoint parameter. Equation~(6.133) has the stem
factor $2^{-(1-\epsilon_{\rm G})(j-k_*)}$ and vertex factors
$2^{-\alpha_vk_v}$, where
\[
 \alpha_{v_*}=(1-\epsilon_{\rm G})s_{v_*}
                       +\epsilon'_{\rm G}l_{v_*},\qquad
 \alpha_v=(1-\epsilon_{\rm G})(s_v-1)
                       +\epsilon'_{\rm G}l_v\quad(v\ne v_*).
\]
The integer $l_v\geq0$ counts the immediate source leaves
in that formula. Choose $\epsilon_{\rm G}=1/8$ and reduce the
ordinary coupling so that the source-growth and normalization
losses are at most $1/16$ of each leaf's scale exponent.
The derivation of Eq.~(6.133)
then permits $\epsilon'_{\rm G}=1/16$: a source-growth factor
subtracts its exponent from the reserved leaf exponent.
Since $k_v\geq0$, discarding the $l_v$ terms and replacing
$1-\epsilon_{\rm G}=7/8$ by $\vartheta=3/4$ gives
\cref{eq:pruned-source-bound} factor by factor. The passage
from Eq.~(6.132) to Eq.~(6.133) already sums the ordinary
branches and the one-child chains; no additional scale sum
is hidden in this replacement. Powers of the fixed $s_0$
are absorbed in $C^m$. In particular, bounded source
coefficients or cancellations between source channels are
not needed for this bound.

Here is the summation explicitly. A tree with $m$ source
leaves has at most $m-1$ branching vertices. The root factor
in \cref{eq:pruned-source-bound} equals
\[
 2^{-\vartheta j}2^{-\vartheta(s_{v_*}-1)k_*}\,.
\]
Drop the ordering constraints on all the branching scales
and sum each from zero to infinity. Each geometric sum
is at most $(1-2^{-\vartheta})^{-1}$. The plane tree shapes
with $m$ leaves and no one-child vertices number at most
$C^m$. Our symmetric Taylor convention averages over the
source labelings, so no further $m!$ is incurred. This proves
the bound $C^m2^{-\vartheta j}$. For $m=1$ there is no branching
subtree: the source-only root has dimension one and yields
$s_j$ times a source coefficient. The same bound follows
from $1-\eta_\beta>\vartheta$ and the initial-endpoint estimates.

Initial endpoints require two additional checks. Their
fermion and source coordinates are summed with the positive
weights of \cref{lem:source-endpoints-localization},
before counting the ordinary trees. Join these coordinates
by a minimizing spatial tree, choosing one deterministically
at ties. Its weight is already paid by the
initial norm. If an initial endpoint contains several
sources, retain them as separate leaves connected inside
this endpoint at the ultraviolet scale zero. These extra
connections have no infrared scale cost; any missing
fixed-scale integration factor is bounded by a constant
per added leaf. Thus they contribute $C^m$ and fit
\cref{eq:pruned-source-bound} with the added vertices at
scale zero. A whole source subtree contained in one initial
endpoint is treated likewise with $k_*=0$.
Sources coinciding at one edge are allowed and have zero
connecting length. The analytic source norm controls the
sum over source multiplicities and their partitions among
endpoints, with the Taylor factorials retained. Equivalently,
in the plane-tree summation one sums compositions of the
$m$ ordered source occurrences into endpoint blocks; their
number and the coefficients of the convergent endpoint
envelope are bounded by $C^m$. The degree reserve similarly
absorbs all fermion labels. These are the only changes to
the finite local endpoint estimate. In particular no entropy
sum over supports is applied to a pointwise kernel bound.

To recover spatial decay, reserve part of the stretched
exponential weight on every contraction edge and every
initial connecting tree. All their scales are at most $j$.
The connected graph contains a path joining a pair of sources
at distance $D$. If its segment lengths are $t_i$, then
\[
       \sum_i\sqrt{s_{k_i}t_i}
           \geq\sqrt{s_j\sum_i t_i}\geq\sqrt{s_jD}\,.
\]
The reserved weights therefore supply
$e^{-c\sqrt{s_jD}}$. The remaining weights give the preceding
scale summation. This proves
\cref{eq:source-shell-pointwise}. Localization differences
use points on these connecting segments and satisfy the
same path bound with a changed constant.

For two sources the pruned tree has one merger, at scale
$k\leq j$. Reserving its connecting weight at that scale
gives the stronger sum
\[
 C\sum_{k=0}^j
   2^{-\vartheta(j-k)}2^{-2\vartheta k}e^{-c\sqrt{s_kD}}\,.
\]
This includes a merger inside an initial endpoint at $k=0$.
Reserve half the exponential at scale $j$ and use
\begin{equation}\label{eq:dyadic-source-sum}
 \sum_{k\geq0}2^{-\vartheta k}e^{-c\sqrt{2^{-k}R}}
                  \leq C(1+R)^{-\vartheta},\qquad R\geq0\,.
\end{equation}
For $R>1$, split at the first $k$ with $2^{-k}R\leq1$.
The tail is geometric of order $R^{-\vartheta}$; writing the
earlier indices as that index minus $l$ bounds their sum
by $CR^{-\vartheta}\sum_{l\geq0}2^{\vartheta l}e^{-c'2^{l/2}}$.
For $R\leq1$ the geometric sum suffices. This proves
\cref{eq:two-source-shell-pointwise}. Finally the initial
source-only kernels have a direct stretched exponential
connecting weight by \cref{eq:source-marked-norm}, which
implies both asserted bounds for them.
\end{proof}

\begin{theorem}[An analytic bulk source functional]
\label{thm:bulk-source-functional}
For every finite set of edges $F$, there is $\rho_F>0$ such
that the bulk source series
\[
 S_\infty(A)=S_\mathrm{UV}(A)+\sum_{j=0}^{\infty}S_j(A),
                \qquad \supp A\subset F\,,
\]
converges absolutely and locally uniformly on the complex
polydisc $\|A\|_\infty<\rho_F$. It is analytic, vanishes
at zero, and its truncation
$S_{<J}=S_\mathrm{UV}+\sum_{j<J}S_j$ satisfies, on a smaller
closed polydisc,
\begin{equation}\label{eq:source-infrared-tail}
                  |S_\infty(A)-S_{<J}(A)|\leq C_F2^{-3J/4}\,.
\end{equation}
Its connected response coefficients
\[
 \cK_m(e_1,\ldots,e_m)
       =\big.\partial_{A(e_1)}\cdots\partial_{A(e_m)}
                                             S_\infty(A)\big|_{A=0}\,,
\]
are translation covariant and satisfy, for $m\geq2$,
\begin{align}
 |\cK_m(\be)|
       &\leq m!C^m(1+s_0D(\be))^{-3/4},\notag\\
 |\cK_2(e,e')|
       &\leq C(1+s_0D(e,e'))^{-3/2}\,.
                                      \label{eq:bulk-response-decay}
\end{align}
The one-source coefficient is uniformly bounded. These are
coefficients of the convergent bulk connected expansion;
identification with cumulants of a probability law requires the
finite-volume comparison proved for the canonical completion in
\cref{sec:finite-volume-sources}.
\end{theorem}
\begin{proof}
Choose $\rho_F$ below the first-step source radius and so that
$C|F|\rho_F<1/2$. \Cref{lem:source-shell-summability}
gives
\[
 \sum_{m\geq1}\sum_{e_1,\ldots,e_m\in F}
 |W_{j,m}(\be)|\rho_F^m
       \leq 2^{-3j/4}\sum_{m\geq1}(C|F|\rho_F)^m
       \leq C_F2^{-3j/4}\,.
\]
The initial source term obeys the analogous convergent
Taylor bound. The Weierstrass criterion gives the claimed
analyticity and the geometric tail estimate. Derivatives
can be taken termwise on every smaller polydisc, either
by the same absolute Taylor sums or by Cauchy's formula.

For the spatial estimates, sum
\cref{eq:source-shell-pointwise} over $j$ and use
\cref{eq:dyadic-source-sum} with $R=s_0D$.
For two sources, first factor out
$(1+s_0D)^{-3/4}$ from
\cref{eq:two-source-shell-pointwise}, and apply that same
dyadic estimate to the remaining sum. The source Taylor
normalization restores the factor $m!$. The ultraviolet
term has faster spatial decay and is covered as well.
Translation covariance follows at every finite scale from
the physical source convention and the translation-covariant
zero-source kernels. Although the localized source vertices
carry phases, a source-only coefficient has no net external
phase. Translation covariance therefore passes
to its absolutely convergent limit.
\end{proof}

We will use the following local-limit criterion. Suppose the source
generating functions of probability laws on matchings satisfy,
uniformly on a fixed smaller polydisc,
\[
 \frac{\cZ_N(A)}{\cZ_N(0)}
       =e^{S_{<J_N}(A)}(1+r_N(A)),\qquad
 J_N\to\infty,\qquad \sup_A|r_N(A)|\to0.
\]
For large $N$, $|r_N|<1/2$, so Cauchy's formula applied to
$\log(1+r_N)$ and \cref{eq:source-infrared-tail} gives convergence
of all edge cumulants to $\cK_m$. Joint moments, and then every
atom probability on a finite set of edge indicators, are
polynomials in these cumulants. This identifies a unique local
limit. It is mixing: expand two local functions as multilinear
polynomials and use the moment--cumulant formula. The partitions
with no block joining the two supports cancel in the covariance;
each remaining partition contains a cumulant tending to zero by
\cref{eq:bulk-response-decay}.

The next subsubsection proves this finite-volume comparison. The
sharper source and height estimates needed for Gaussianity are
proved in \crefrange{sec:sharp-source-flow}{sec:height-gaussianity}.

\subsubsection{Finite-volume source comparison and the canonical local limit}
\label{sec:finite-volume-sources}

We identify the bulk source functional with the positive canonical
law, whose orientation counts fluctuate. Its torus partition function
is a signed combination of four Grassmann integrals, while extraction
of prescribed counts requires a continuous family of characters.
We therefore need bounds uniform even when a reciprocal-grid point
is a Fermi zero. Integrations stop at scale $N^{-1}$; the nearest
mode to each zero is retained in an unnormalized polynomial integral.
The finite-size comparison, last-mode formula and deleted-determinant
estimate below then give the canonical local limit. The prescribed
sector is treated in \cref{sec:fixed-sector-limit}.

Fix the positive liquid edge weights and one matched
physical dispersion $D$ from \cref{thm:matched-bulk-flow}.
Let $\Lambda_N=\Z^2/\sB_N\Z^2$, where both
$N^{-1}\sB_N$ and its inverse are uniformly bounded, and put
$V_N=|\det\sB_N|\asymp N^2$. This includes the diagonal tori in a compact positive slope range. A boundary character
$\theta\in\R^2/(2\pi\Z)^2$ has momentum grid
\[
 \cP_N(\theta)
 =\sB_N^{-\trps}(2\pi\Z^2+\theta)
                          \pmod {2\pi\Z^2}\,;
\]
that is, the Fourier modes of fields which are
$e^{i\theta}$-periodic along the two periods of the torus.
Write $\fZ_N^\theta(A)$ for the unnormalized polynomial
with this character, in a gauge which is trivial near the fixed
finite source set $F$. A source counts edge indicators and is
unchanged by this gauge. All constants below are uniform in
$\theta$, in the allowed torus shapes, and in the choice of the
matched dispersion. We retain the exponent $\kappa=1/2$, leaving
room in the exponent $3/4$ of the preceding subsubsection: the
finite-size marks cost additional dimensional factors, which
are paid by the difference $3/4-\kappa$.
All assertions in this subsubsection are for sufficiently large $N$,
depending on the fixed $\beta$ and source set $F$.

\begin{lemma}[Finite-size remainders and final source kernels]
\label{lem:finite-size-source-kernels}
The finite-volume integration can use the bulk running local
coefficients and bulk field normalizations. Choose $J=J_N$ with
$s_J\asymp N^{-1}$, with a sufficiently large fixed proportionality
constant. Up to that scale, differences from the bulk expansion
carry a finite-size mark $e^{-c\sqrt{Ns_j}}$ at root scale $j$.
More precisely, source-only shell coefficients at fixed source
positions satisfy
\begin{equation}\label{eq:finite-source-shell-difference}
 |W^{N,\theta}_{j,m}(\be)-W_{j,m}(\be)|
       \leq C_F^m2^{-\kappa j}e^{-c\sqrt{Ns_j}},
                  \qquad 0\leq j<J\,,
\end{equation}
with the analogous ultraviolet estimate. For kernels with
$\fn=2,4$ external normalized fermions and $m\geq1$ fixed sources,
the mixed norm, summed in all fermion positions and pointwise in
the sources, is bounded by
\begin{equation}\label{eq:finite-mixed-source-kernel}
 \sum_{\bx,\,\boldsymbol\alpha}
       |W^{N,\theta}_{j;\fn,m}(\bx;\be)|
                       \leq C_F^{m+\fn}s_j^{\kappa-\fn/2}\,.
\end{equation}
Derivative fields in a localized representation are included
with their scale factors. At zero source, the anchored
$\fn$-fermion norm is $Cd_\beta s_j^{2-\fn/2}$ for $\fn=2,4$.
These estimates continue through a final integration of all
but one nearest mode to each Fermi point.
\end{lemma}
\begin{proof}
Use bulk localization as in \cite[App.~B]{GMT20}. Finite-size
errors are marked propagator or endpoint differences; their extra
spatial decay will also pay for local jets that do not vanish.

Here is the explicit difference estimate. Let
$\lambda_N\geq cN$ be the shortest nonzero period length.
Lift the contraction tree, together with the
internal connecting trees of its endpoints, starting at a fixed
root. If their total length is at least $\lambda_N/4$, their
reserved spatial weights already give
$e^{-c\sqrt{Ns_j}}$, by the path inequality recalled below.
Otherwise every pair of lifted positions
has distance less than $\lambda_N/4$ and the lifts are consistent
in a single planar neighborhood. On this latter class, for any
derivative orders $a,b$ occurring in localization, periodization gives
\[
 |\partial_x^a\partial_y^b
       (g^{N,\theta}_k-g_k)(x,y)|
 \leq C_{a,b}s_k^{1+a+b}e^{-c\sqrt{Ns_k}}.
\]
Indeed, the periodized propagator with character $\theta$ is
the sum over the period lattice of the planar propagator at the
translated positions, multiplied by unit phases, so that the
difference is the sum over the nonzero images; every nonzero
image has distance at least $3\lambda_N/4$, where the stretched
exponential decay of \cref{lem:infrared-covariance} applies,
and summing images uses their uniform separation and
$Ns_k\geq c_*$. This also treats the tails
omitted when planar coordinate sums are restricted to the
chosen neighborhood.

Apply the marked Gram identity
\cref{eq:marked-gram-difference} to each contraction determinant.
For an $r\times r$ block on scale $k$, its cofactors and the
preceding difference bound give
$Cr^2(Cs_k)^r e^{-c\sqrt{Ns_k}}$, with the usual derivative powers.
Thus the finite-size factor replaces one Gram pair without an
inverse determinant. Telescope tree lines likewise.
\Cref{lem:marked-source-comparison} pays the marked choices and
endpoint coordinates; characters have modulus one.

By \cref{it:rg-locality}, a nonplanar initial coefficient has
size at least $cN$, and its reserved size exponent gives
$Cd_\beta e^{-cN}$. For all other endpoints, join their internal
spatial trees to the lifted contraction tree. Source positions are
included by \cref{lem:source-endpoints-localization}. A path
of length $cN$, whether in these trees or in a nonprincipal image,
pays the finite-size factor because
\[
 \sum_i\sqrt{s_{k_i}t_i}\geq\sqrt{s_j\sum_i t_i}
                                      \geq c\sqrt{Ns_j},
 \qquad k_i\leq j\,.
\]
All character factors have modulus one, and
\cref{cor:infrared-torus} gives uniform Gram bounds for the
allowed torus shapes. The degree reserve absorbs the choices of a
marked entry. Reserving part of the spatial exponent therefore
extends the cited comparison to our summed endpoints, with factor
$e^{-c\sqrt{Ns_j}}$ (and $e^{-c\sqrt N}$ at the first integration).

Localization is defined using the bulk jets. Thus a finite-size
difference need not have the vanishing jets of a bulk remainder.
It must not be assigned a Taylor gain it does not possess. Instead,
mark its earliest finite-size difference. A missing Taylor gain
along a path from $k$ to $j$ costs at most $2^{b(j-k)}$, with a
fixed $b$ (one may take $b=2$ for the quadratic and quartic
localizations here). For $Ns_j$ bounded below,
\begin{equation}\label{eq:winding-pays-relevance}
 \sum_{k=0}^j2^{b(j-k)}e^{-c\sqrt{Ns_k}}
       \leq C_b e^{-c'\sqrt{Ns_j}},\qquad 0<c'<c\,.
\end{equation}
Indeed, put $l=j-k$ and use $Ns_k=2^lNs_j$; the resulting
series is bounded by a constant times its reserved exponential
at $l=0$. Each additional finite-size insertion has its own mark
and its ordinary small endpoint factor. This proves, inductively
in the integration scale and then by the convergent tree envelope,
that the finite-size local terms remain $O(d_\beta)$ in scaled
units. It also proves the marked bounds just described. No
decaying finite-volume trajectory, or finite-volume tuning of the
bulk counterterms, has been assumed.

For a fixed set of sources, retain the source subtree from
\cref{lem:source-shell-summability}. The finite-size mark can
be reserved before summing its scales. Weakening $3/4$ to
$\kappa$ pays the additional marked choices and localization
moments; the same geometric sums give
\cref{eq:finite-source-shell-difference}. With $\fn$ external
fermions, summing their positions leaves the extra root factor
$s_j^{-\fn/2}$. All internal dimensions and source-merger factors
are the same. This gives \cref{eq:finite-mixed-source-kernel},
including a single source, for which
$|\Gamma_j|\leq Cs_j^{-\eta_\beta}$ and
$1-\eta_\beta>\kappa$. The source-free root has dimension
$2-\fn/2$ and at least one ordinary endpoint, giving the last
asserted norm bound. Source Taylor factorials, degree reserves,
and the positive sums over nonlocal initial endpoints are kept
throughout, exactly as in \cref{lem:source-shell-summability}.

Finally choose a nearest grid point $k_\theta^\omega$ to each
$p^\omega$, resolving ties by a fixed rule. Grid separation and
covering radius are respectively bounded below by $c/N$ and above
by $C/N$. Every other grid point has distance at least $c/(2N)$
from that zero: if a second point were closer, so would the chosen
nearest point, contradicting separation. Take the proportionality
constant defining $J$ large enough that both selected modes lie
where the last cutoff equals one. There are only a bounded number
of other remaining modes. Their normalized quadratic eigenvalues
have modulus at least $c/N$, by the simple-zero bounds and the
curvature estimate. Integrating these modes together has covariance
supremum and squared Gram norm $C/N$ and spatial $\ell^1$ norm
$CN$. Torus diameter is $O(N)$, so these are the same scale bounds
as before, with changed constants. The connected expansion of this
one final step obeys the same estimates. Neither selected mode is
inverted, even when its eigenvalue is zero. The bounds also justify
the analytic source series on a fixed smaller polydisc, initially
computed as Taylor identities.
\end{proof}

\begin{proposition}[The unnormalized last-mode formula]
\label{prop:last-mode-source-formula}
Let
\[
 \widetilde D_N^\theta
   =\prod_{k\in\cP_N(\theta)\setminus
                   \{k_\theta^+,k_\theta^-\}}D(k)\,.
\]
There is a bulk scalar $\Delta$, with $|\Delta|\leq Cd_\beta$,
and nonzero source-independent factors $P_N^\theta$ such that
\begin{equation}\label{eq:last-mode-prefactor}
 P_N^\theta
 =\xi_N^\theta e^{V_N\Delta}N^{-2}
       \widetilde D_N^\theta e^{r_N^\theta},
       \qquad |\xi_N^\theta|=1,\quad |r_N^\theta|\leq Cd_\beta\,.
\end{equation}
In four Grassmann variables $\chi_{b,\omega},\chi_{w,\omega}$,
put $B_\omega=\chi_{b,\omega}\chi_{w,\omega}$ and orient the
integral by $\int\d\chi\,B_+B_-=1$. Then, for sources on a fixed
finite $F$ and a fixed sufficiently small complex polydisc,
\[
 \fZ_N^\theta(A)
 =P_N^\theta e^{S_\infty(A)+\epsilon_N^\theta(A)}
       \int\d\chi\,\exp(H_N^\theta(\chi)+R_N^\theta(\chi,A))\,,
\]
where $\epsilon_N^\theta(0)=R_N^\theta(\chi,0)=0$ and
\begin{align}
 H_N^\theta&=(a_++u_+)B_++(a_-+u_-)B_-+vB_+B_-,\notag\\
 a_\omega&=N\big[L_\omega(q_\omega)
                        +R_\omega(q_\omega)/Z_J\big],
       \qquad q_\omega=k_\theta^\omega-p^\omega,\notag\\
 |a_\omega|&\leq C,\qquad |u_+|+|u_-|+|v|\leq Cd_\beta,\notag\\
 \sup_A\big(|\epsilon_N^\theta(A)|
            +\|R_N^\theta(\cdot,A)\|_\mathrm{alg}\big)
                         &\leq C_F N^{-\kappa}\,.
                                           \label{eq:last-mode-bounds}
\end{align}
Here $R_\omega(q)=D(p^\omega+q)-L_\omega(q)$, and
$\|\cdot\|_\mathrm{alg}$ is the sum of absolute monomial
coefficients. In particular, writing
$b_N^\theta=(a_++u_+)(a_-+u_-)+v$, one has
\begin{equation}\label{eq:uniform-twist-factorization}
 |b_N^\theta|\leq C,\qquad
 \sup_A|\fZ_N^\theta(A)
             -e^{S_\infty(A)}\fZ_N^\theta(0)|
                     \leq C_F N^{-\kappa}|P_N^\theta|\,.
\end{equation}
No lower bound on $b_N^\theta$ is asserted for every twist.
\end{proposition}

Before the proof we explain the shape of the formula. The
prefactor $P_N^\theta$ collects everything that does not involve
the sources or the two retained modes: the Gaussian determinants
of all integrated modes, which produce the deleted determinant
$\widetilde D_N^\theta$, and the vacuum energy accumulated by the
integration, which is extensive and produces $e^{V_N\Delta}$ up to
a bounded error $e^{r_N^\theta}$. The factor $e^{S_\infty(A)}$ is
the bulk source-only functional, with the finite-volume error
$\epsilon_N^\theta(A)$. The integral over the four Grassmann
variables $\chi_{b,\pm},\chi_{w,\pm}$ represents the mode nearest
to each Fermi point, rescaled by $\sqrt N$; its effective action
$H_N^\theta$ is quadratic, with coefficients $a_\omega+u_\omega$,
plus the one possible quartic monomial $vB_+B_-$, and $R_N^\theta$
is the source-dependent part. Since $B_\omega^2=0$, expanding the
exponential gives at zero source
\[
 \int\d\chi\,e^{H_N^\theta}
 =\int\d\chi\,\big(1+(a_++u_+)B_+\big)\big(1+(a_-+u_-)B_-\big)
                                   \big(1+vB_+B_-\big)
 =b_N^\theta\,,
\]
which explains the definition of $b_N^\theta$. The quantities
$a_\omega$ are the normalized bare eigenvalues at the retained
modes, multiplied by $N$ to compensate the rescaling; since the
retained mode is within $C/N$ of the Fermi point, $a_\omega$ is
bounded, but it may be arbitrarily small, and in the free case
($u=v=0$) the integral $b_N^\theta=a_+a_-$ vanishes exactly at a
singular twist. This is why no lower bound on $b_N^\theta$ can be
asserted for every twist, and why the estimates are formulated
relative to $|P_N^\theta|$ rather than to $|\fZ_N^\theta(0)|$.

\begin{proof}
Use \cref{lem:finite-size-source-kernels}, retaining the two
selected momenta. In orthonormal Fourier coordinates, rescale a
remaining normalized field by $\sqrt N\chi$; its value at a
lattice vertex is therefore a sum of two terms of modulus
$\sqrt{N/V_N}$ times a Grassmann variable and a phase of modulus one.
A source kernel with $\fn$ external fermions consequently costs
\[
 C^{\fn}N^{-\fn/2}\,C_F^{m+\fn}s_J^{\kappa-\fn/2}
                                      \leq C_F^{m+\fn}N^{-\kappa}\,.
\]
Only $\fn=2,4$ survive. Sum in $m$ on a smaller source polydisc
to obtain the asserted algebra norm of $R_N^\theta$. This
includes off-diagonal flavor pairs, higher source orders, and
repeated source edges. With zero sources, translation invariance
forbids a quadratic term between the distinct selected momenta.
The source-free bounds in the lemma give $u_\omega,v=O(d_\beta)$:
the free anchor contributes $V_N$, so that the bound at degree
$\fn$ is
$CV_N(N/V_N)^{\fn/2}d_\beta s_J^{2-\fn/2}\leq Cd_\beta$.
The quadratic Gaussian part gives
$a_\omega$ as displayed. Since $|Z_J|^{-1}\leq CN^{Cd_\beta^2}$
and $R_\omega(q)=O(|q|^2)$, these coefficients are bounded.

The source-only scalar differs from $S_\infty$ by the omitted
bulk tail, the finite-size shell differences, and the final
bounded-mode integration. The first and last are $O_F(N^{-\kappa})$.
For the middle term, use \cref{eq:finite-source-shell-difference}
and $s_j=s_02^{-j}$:
\[
 \sum_{j<J}2^{-\kappa j}e^{-c\sqrt{Ns_j}}
 \leq CN^{-\kappa}
       \sum_{l\geq0}2^{\kappa l}e^{-c'2^{l/2}}
 \leq CN^{-\kappa}\,.
\]
The initial scalar has the same finite-size estimate. This proves
the bound on $\epsilon_N^\theta$, with its value at zero
subtracted.

We give the normalization accounting, since discarding it would
invalidate a comparison of partition functions. Factor out the
Gaussian determinant on all unselected modes, using the exact
quadratic reorganization of \cref{lem:moving-first-integration}.
At every further rescaling retain the Gaussian normalization.
The bulk vacuum density on scale $j$ is bounded by
$Cd_\beta s_j^2$. Its finite-size difference is at most
$Cd_\beta s_j^2e^{-c\sqrt{Ns_j}}$, by the marked-tree proof
above, now with a vacuum root and one ordinary endpoint.
Consequently its sum contributes a bulk density to $\Delta$
and a bounded remainder, because
\[
 V_N\sum_{j<J}d_\beta s_j^2e^{-c\sqrt{Ns_j}}
       +V_N\sum_{j\geq J}d_\beta s_j^2\leq Cd_\beta\,.
\]
The ultraviolet density and its remainder have the same bounds.

For completeness, the quadratic-normalization terms require a
slightly different estimate near a Fermi point. Their logarithm
has the form
\[
 f_j(q)=\log\Big(1+
       \frac{z_{j+1}\chi_j(q)L_\omega(q)}
                    {L_\omega(q)+R_\omega(q)/Z_j}\Big)\,.
\]
Here $\chi_j$ is the remaining adapted cutoff at that rescaling,
of scale $s_j$, equal to one at the two selected modes.
Split $f_j=\log(1+z_{j+1}\chi_j)+h_j$. The first term is a
smooth Gevrey cutoff. Its grid-sum error, after multiplying the
normalized average by $V_N$, is
$C|z_{j+1}|(Ns_j)^2e^{-c\sqrt{Ns_j}}$.
For the second, away from zero,
\[
 |\partial^l h_j(q)|
    \leq C_l|z_{j+1}|\,|Z_j|^{-1}|q|^{1-l},
          \qquad l=0,1,2,\quad |q|\leq Cs_j\,,
\]
with the corresponding cutoff bounds at the boundary of its
support. Use centered reciprocal cells of diameter $C/N$.
The cells at distance $O(N^{-1})$ from zero contribute
$C|z_{j+1}|\,|Z_j|^{-1}N^{-1}$ to the unnormalized sum error.
On the other cells the linear Taylor term averages to zero;
the integral of the Hessian bound is
$C|z_{j+1}|\,|Z_j|^{-1}s_j$. These quantities bound the total
sum error, since $Ns_j\geq c$. Thus no smoothness at the zero
of the quotient $R_\omega/L_\omega$ is being assumed.
Their sum in $j$ is $O(d_\beta^2)$, using
$|z_{j+1}|\leq Cd_\beta^2$ and
$|Z_j|^{-1}\leq C2^{Cd_\beta^2j}$.

Deleting each selected mode subtracts
$\log(1+z_{j+1})$ from the smooth part at every scale before
the stopping scale. The discrepancy from $h_j$ is bounded
by $C|z_{j+1}|\,|Z_j|^{-1}/N$ and is summable with bound
$Cd_\beta^2$. The two smooth subtractions telescope to
$-2\log Z_J$. At the last integration, the change from the
original Fourier variables to $\chi$ has Berezin Jacobian
$(Z_J/N)^2$. Its $Z_J^2$ cancels that scalar exactly,
leaving $N^{-2}$. The bounded final integrations add only
$O(d_\beta)$ to the logarithmic remainder. Together these
computations prove \cref{eq:last-mode-prefactor}; any Fourier
ordering or gauge phase is denoted by $\xi_N^\theta$.
All logarithms used here are the branches obtained from
$\log(1+w)$ with $|w|<1/2$. No selected eigenvalue is in a
denominator. Singular twists follow from the same unnormalized
Berezin identities, not from a normalized Gaussian expectation.

The algebra norm is submultiplicative and Berezin integration
has norm one. Since $\|H_N^\theta\|_\mathrm{alg}\leq C$,
\[
 \Big|\int\d\chi\,e^{H_N^\theta}
                          (e^{R_N^\theta}-1)\Big|
 \leq e^C\big(e^{\|R_N^\theta\|_\mathrm{alg}}-1\big)
 \leq C_F N^{-\kappa}\,.
\]
The zero-source integral is exactly $b_N^\theta$.
The scalar estimate for $\epsilon_N^\theta$ now gives
\cref{eq:uniform-twist-factorization}.
\end{proof}

The prefactors $P_N^\theta$ of different twists must be compared,
since the canonical partition function is a combination of four
of them. The only $\theta$-dependent factor of unbounded size is
the deleted determinant $\widetilde D_N^\theta$, whose logarithm
is a Riemann sum of $\log|D|$ over the shifted grid with the two
nearest modes removed. A naive quadrature estimate near the
logarithmic singularities would lose a power of $\log N$, which
would be fatal to the comparison; the following lemma shows that
the deletion of the two nearest modes compensates the
singularities exactly, up to a bounded factor.

\begin{lemma}[Deleted determinants on shifted grids]
\label{lem:deleted-determinant-comparison}
Set
\[
 p_D=\int_{[-\pi,\pi]^2}\log|D(k)|\,\frac{\d k}{(2\pi)^2}\,.
\]
For the two simple zeros and the bounded-shape tori above,
\begin{equation}\label{eq:deleted-determinant-comparison}
 C^{-1}N^2e^{V_Np_D}\leq|\widetilde D_N^\theta|
                              \leq CN^2e^{V_Np_D}\,.
\end{equation}
In particular, with $\cM_N=e^{V_N(p_D+\Re\Delta)}$,
\[
       C^{-1}\cM_N\leq|P_N^\theta|
                                      \leq C\cM_N\,.
\]
\end{lemma}
\begin{proof}
Put $f=\log|D|$. At distance $r$ from either zero,
$|\partial^l f|\leq C_l r^{-l}$ for $1\leq l\leq3$, and
$f(k)=\log|k-p^\omega|+O(1)$ locally, uniformly in angle.
Tile the Fourier torus by centered reciprocal parallelograms,
one around each point of $\cP_N(\theta)$, each of area
$(2\pi)^2/V_N$ and diameter at most $C/N$.
Let $U$ be the union of all cells whose centers are within
$C_0/N$ of either zero, for a large fixed $C_0$. It consists
of a bounded number, say $t$, of cells; its boundary stays
at distance between $c/N$ and $C/N$ from the zeros.
All unselected centers in $U$ have distance at least $c/N$
from the zeros, by the nearest-point argument in the preceding
lemma. Thus
\[
 \sum_{k\in U,\ k\ne k_\theta^\pm}f(k)
       =-(t-2)\log N+O(1),\qquad
 \frac{V_N}{(2\pi)^2}\int_U f(k)\,\d k=-t\log N+O(1)\,.
\]
The integral estimate follows by scaling by $N$; the logarithmic
singularities are integrable on a uniformly bounded region.

On each other cell $C_k$, Taylor's formula and central symmetry
give
\[
 f(k)-\frac1{|C_k|}\int_{C_k}f
   =-\frac12 Q_N:\nabla^2f(k)
       +O\big(N^{-3}\sup_{C_k}|\nabla^3f|\big)\,,
\]
where $Q_N$ is the common second-moment matrix, with
$\|Q_N\|\leq CN^{-2}$. The sum of the remainders is $O(1)$:
annular summation bounds it by
$CN^{-1}\int_{c/N}^1r^{-3}r\,\d r+C/N$.
Replace the sum of Hessians by its integral. The error after
multiplication by $Q_N$ is again $O(1)$, using the same third
derivative bound. Finally,
\[
 \bigg|\int_{\T^2\setminus U}\partial_i\partial_jf\bigg|
   =\bigg|\int_{\partial U}\partial_jf\,n_i\bigg|\leq C\,,
\]
up to the immaterial sign of the outward normal: periodic
boundary terms cancel, $|\partial U|=O(N^{-1})$, and
$|\nabla f|\leq CN$ there. This cancellation avoids the
spurious $\log N$ loss from summing $|\nabla^2f|$.
It follows that the total contribution of the cells outside
$U$ to the quadrature error is $O(1)$. Adding the cells in
$U$ gives
$\log|\widetilde D_N^\theta|-V_Np_D=2\log N+O(1)$,
which proves the lemma and the assertion about $P_N^\theta$.
\end{proof}

\begin{theorem}[The canonical local limit]
\label{thm:canonical-local-limit}
Let $\mu_N^\mathrm{can}$ be the probability law on all torus
matchings with positive weights
$F_{\beta,\Gamma}(M)\prod_r\omega_r^{N_r(M)}$.
For fixed liquid weights and sufficiently large $\beta$, it
converges locally to a translation-invariant mixing law
$\mu_\infty^\mathrm{can}$. For every fixed finite edge set $F$,
on a sufficiently small fixed complex source polydisc,
\begin{equation}\label{eq:canonical-source-identification}
 \sup_A\bigg|\mu_N^\mathrm{can}
       \bigg(e^{\sum_e A(e)\one_{\{e\in M\}}}\bigg)
                     -e^{S_\infty(A)}\bigg|
                             \leq C_F N^{-\kappa}\,.
\end{equation}
Its edge cumulants are exactly $\cK_m$ from
\cref{thm:bulk-source-functional}, and therefore obey
\cref{eq:bulk-response-decay}. The functional $S_\infty$ is
independent of the choice of matched dispersion. This theorem
by itself does not identify the law conditioned on the prescribed counts;
that identification is proved in \cref{sec:fixed-sector-limit}.
It does not prove a height \GFF.
\end{theorem}
\begin{proof}
Write $\cZ_N^\mathrm{can}(A)$ for the positive canonical
sum with sources. We first record the algebra which converts
between it and the boundary twists. Fix a reference matching $M_0$
and classify each matching by the mod-$2$ homology class
$h\in(\Z/2\Z)^2$ of $M\xor M_0$. By the toroidal Kasteleyn
sign rule \cite[\S2.1]{KSW16}, the determinant sign is constant in
each class; write it as $\epsilon_h$, with an arbitrary sign for
an empty class. The pattern weight $F_{\beta,\Gamma}(M)$ multiplies
each matching term without changing this sign, by
\cref{lem:pattern-algebra}.
A parity boundary character $\tau\in\{0,\pi\}^2$ multiplies
the contribution of $M$ by $(-1)^{(\tau/\pi)\cdot h}$, up to a
common unit phase from the reference matching, since the
character changes the sign of a cycle exactly when the cycle
crosses the corresponding period an odd number of times. Fourier
inversion on $(\Z/2\Z)^2$ therefore yields
\begin{equation}\label{eq:four-character-canonical}
 \cZ_N^\mathrm{can}(A)
      =\sum_{\tau\in\{0,\pi\}^2}c_\tau
                       \fZ_N^\tau(A),\qquad |c_\tau|\leq1\,.
\end{equation}
For example, after removing the common phase,
$c_\tau=\frac14\sum_h\epsilon_h
                      (-1)^{(\tau/\pi)\cdot h}$.
The pattern weights cause no change to this identity: each
complete matching still has its original determinant sign
in the polynomial \cref{eq:rg-berezin}.
For every unitary character, positivity also gives
\begin{equation}\label{eq:twist-dominated-by-positive}
               |\fZ_N^\theta(0)|
                              \leq\cZ_N^\mathrm{can}(0)\,.
\end{equation}

The key step is a lower bound on the positive partition function
$\cZ_N^\mathrm{can}(0)$ of the same order $\cM_N$ as the
prefactors. It is obtained from a single well-chosen twist,
at which the last-mode integral is bounded below, together with
the domination \cref{eq:twist-dominated-by-positive}; no
cancellation between the four parity characters needs to be
controlled. Choose a character $\theta_*$ for which both nearest
grid points have distance at least $\delta/N$ from their Fermi
points. Such a character exists for a fixed small $\delta>0$:
as the grid is translated uniformly over a reciprocal cell,
the probability of a point in either radius-$\delta/N$ ball
is at most $C\delta^2$ (the cell has area of order $N^{-2}$ and
each ball has area of order $\delta^2N^{-2}$), which is less
than one. The upper nearest-point distance remains $C/N$. The
simple-zero bound $|L_\omega(q)|\geq c|q|$ and
$|R_\omega(q)/Z_J|\leq CN^{-2+Cd_\beta^2}$ give
$|a_\omega|\geq c\delta$ for all large $N$. In
\cref{prop:last-mode-source-formula},
$b_N^{\theta_*}=a_+a_-+O(d_\beta)$, so, after making
$d_\beta$ sufficiently small, $|b_N^{\theta_*}|\geq c_*>0$.
Since $\fZ_N^{\theta_*}(0)=P_N^{\theta_*}b_N^{\theta_*}$,
\cref{eq:twist-dominated-by-positive,eq:deleted-determinant-comparison} now imply
\begin{equation}\label{eq:canonical-normalization-bound}
 c\cM_N\leq\cZ_N^\mathrm{can}(0)
                                      \leq C\cM_N\,.
\end{equation}
Here the lower bound is
$\cZ_N^\mathrm{can}(0)\geq|\fZ_N^{\theta_*}(0)|
=|P_N^{\theta_*}|\,|b_N^{\theta_*}|\geq c\cM_N$, and the upper
bound follows from \cref{eq:four-character-canonical}, which
gives $\cZ_N^\mathrm{can}(0)\leq\sum_\tau|P_N^\tau|\,|b_N^\tau|$,
together with $|b_N^\tau|\leq C$. This argument uses no sign or
noncancellation assertion about dressed Gaussian determinants.

Now sum \cref{eq:uniform-twist-factorization} over the four
parity characters with the coefficients in
\cref{eq:four-character-canonical}: since the bulk factor
$e^{S_\infty(A)}$ is the same for every twist,
\[
 \sup_A\big|\cZ_N^\mathrm{can}(A)
            -e^{S_\infty(A)}\cZ_N^\mathrm{can}(0)\big|
 \leq\sum_\tau|c_\tau|\,C_FN^{-\kappa}|P_N^\tau|
 \leq C_FN^{-\kappa}\cM_N\,.
\]
Division by \cref{eq:canonical-normalization-bound} proves
\cref{eq:canonical-source-identification}. On a smaller polydisc
$e^{S_\infty}$ is bounded away from zero, so Cauchy's formula
also gives convergence of every derivative of the logarithm,
as explained after \cref{thm:bulk-source-functional}.
The criterion following \cref{thm:bulk-source-functional} now
gives the local limit and mixing. Matching constraints and
translation invariance pass to the limit.
Finally the finite-volume law is independent of the choice of
matched dispersion. Applying the proved convergence to each
choice identifies $e^{S_\infty}$, and then its analytic logarithm
vanishing at zero, proving the last assertion.
\end{proof}

\begin{corollary}[Continuity in the auxiliary weights]
\label{cor:canonical-weight-continuity}
Use logarithmic weight ratios
\[
 u=(\log(\omega_2/\omega_1),\log(\omega_3/\omega_1))\,,
\]
in a compact subset of the liquid region, keeping the canonical
coefficients fixed. On a sufficiently small fixed source polydisc,
\[
 \sup_A|S_\infty(A;u)-S_\infty(A;u')|
                         \leq C_F|u-u'|^{1/5}\,.
\]
All fixed edge cumulants are continuous in $u$. In particular,
if $u_N\to u$, the canonical laws with weights $u_N$ converge
locally to $\mu_\infty^\mathrm{can}(u)$, with source error bounded
by $C_F(N^{-1/2}+|u_N-u|^{1/5})$.
\end{corollary}
\begin{proof}
The idea is to interpolate between the two weights at a finite
volume $N$, where the dependence on the weights is explicit and
Lipschitz with a constant of order $N^2$, and to use the uniform
finite-volume rate $N^{-\kappa}$ at both endpoints; the choice of
$N$ balances the two errors. The constants in the preceding
theorem are uniform on the compact weight set: the
first-integration reserves, simple-zero bounds, bulk contraction
and all finite-size estimates can be chosen uniformly there. No
continuous choice of the matched dispersion was used. Fix one
admissible sequence of bounded-shape tori to compare different
weights. The canonical coefficients $c_{\beta,\Gamma}$ do not
depend on the auxiliary weights, so the weight $u$ enters the
canonical law only through the factor
$\prod_r\omega_r^{N_r(M)}=\omega_1^{V_N}e^{u_1N_2(M)+u_2N_3(M)}$,
and, for $G_A(M)=\exp(\sum_eA(e)\one_{\{e\in M\}})$,
\[
 \partial_{u_r}\mu_N^\mathrm{can}(G_A)
       =\cov_{\mu_N^\mathrm{can}}(G_A,N_{r+1}),
 \qquad
 |\partial_{u_r}\mu_N^\mathrm{can}(G_A)|
                  \leq 2V_N e^{\rho_F|F|}\leq C_FN^2\,.
\]
This estimate holds for complex $A$, since the probability
weights themselves are positive and $0\leq N_{r+1}\leq V_N$.
For nearby $u,u'$ their joining segment stays in a slightly
enlarged compact liquid set. Integrate the bound along that
segment and apply \cref{eq:canonical-source-identification}
at both endpoints to obtain
\[
 \sup_A|e^{S_\infty(A;u)}-e^{S_\infty(A;u')}|
                   \leq C_F(N^{-\kappa}+N^2|u-u'|)\,.
\]
Choose $N$ of order $|u-u'|^{-1/(2+\kappa)}$.
Since $\kappa=1/2$, the exponent is $1/5$; larger parameter
differences are covered by boundedness. Decrease $\rho_F$ so
that $e^{\rho_F|F|}-1<1/2$. The finite-volume generating
functions, and their limits, then lie in the disk of radius
$1/2$ about one. The analytic logarithm on this disk is
Lipschitz, giving the claimed estimate for $S_\infty$.
Cauchy's formula gives continuity of its derivatives.
The last assertion follows by a triangle inequality using
the uniform finite-volume rate.
\end{proof}

\begin{corollary}[Local limit of the fixed-sector law]
\label{cor:sector-source-reduction}
Let $\cC_N$ be the prescribed orientation-count sector,
and put $p_N=\mu_N^\mathrm{can}(\cC_N)$. At fixed
auxiliary weights, if $N^{-\kappa}/p_N\longrightarrow0$, then
the law conditioned on $\cC_N$ has the same local limit
$\mu_\infty^\mathrm{can}$.
More quantitatively, its source generating function differs from
$e^{S_\infty(A)}$ by at most $C_FN^{-\kappa}/p_N$.
For the counts in \cref{thm:coefficient-rg}, this is exactly
the fixed-sector law appearing there.
The same conclusion holds for weights $u_N\to u$ in a compact
liquid set, with limiting law $\mu_\infty^\mathrm{can}(u)$ and
error $C_F(N^{-1/2}/p_N+|u_N-u|^{1/5})$.
\end{corollary}
\begin{proof}
Recall from \cref{it:rg-sources} that the polynomial
$\fZ_{\beta,\Lambda}(y,A)$ carries a phase $y_r$ on every
edge of type $r$, so that the matchings with prescribed counts
$b_2,b_3$ of types two and three are exactly the coefficient of
$y_2^{b_2}y_3^{b_3}$. Use the normalized circle integral for this
coefficient extraction. More explicitly, this operator is
\[
 \mathcal E_{\cC_N}H=
 \int_{[0,2\pi]^2}e^{-i(b_2t_2+b_3t_3)}
             H(e^{it_2},e^{it_3})\,
                         \frac{\d t_2\,\d t_3}{(2\pi)^2}.
\]
Its norm on bounded functions of the phases is at most one.
Consequently a uniform absolute error before extraction remains
an absolute error of the same size afterwards; extracting a
coefficient does not require dividing by any twisted determinant.
This is where the twist uniformity of the preceding
estimates is used.

We now express the polynomial at unit phases through the twisted
polynomials of \cref{prop:last-mode-source-formula}. Unit edge
phases are locally a vertex gauge and globally a boundary
character. Indeed their alternating product around a hexagon is
one, since each hexagon has exactly one edge of each type in
each of its two alternating classes, so a vertex gauge removes
them on the universal cover and leaves only their holonomies
along the two periods, which constitute a character $\theta(y)$.
Every pattern monomial contains each of its incident vertices
once. The same Grassmann change of variables therefore transforms
both the quadratic part and all interaction terms, leaving the
edge sources unchanged. Its Berezin Jacobian is a common unit
phase, independent of $A$:
\[
 \fZ_N(y,A)=\zeta_N(y)\fZ_N^{\theta(y)}(A),
                  \qquad |\zeta_N(y)|=1.
\]
Consequently \cref{eq:uniform-twist-factorization,lem:deleted-determinant-comparison}
give a uniform absolute error $C_FN^{-\kappa}\cM_N$ before
extraction, with the same bulk factor $e^{S_\infty(A)}$ for
every $y$. Taking this factor outside the integral gives
\[
 \big|\mathcal E_{\cC_N}\fZ_N(\cdot,A)
       -e^{S_\infty(A)}\mathcal E_{\cC_N}\fZ_N(\cdot,0)\big|
                       \leq C_FN^{-\kappa}\cM_N.
\]
The extracted zero-source coefficient is the signed sum of the
weights of the matchings in the sector. The determinant sign is
constant on the sector, since the prescribed counts fix the
winding of the matchings, hence their relative homology class,
and the positive weights are those of the canonical law; so the
coefficient has modulus $\cZ_N^\mathrm{can}(0)p_N$, and likewise
the extracted coefficient at source $A$ is, up to the same sign,
the sector generating function. Divide by
$\cZ_N^\mathrm{can}(0)p_N$ and use
\cref{eq:canonical-normalization-bound} to obtain the error
$C_FN^{-\kappa}/p_N$. Derivatives and local
laws follow by the preceding theorem's argument. For varying
weights apply \cref{cor:canonical-weight-continuity}.
\end{proof}

The corollary is a reduction, not a lower bound on $p_N$.
The next subsubsection chooses auxiliary weights with prescribed
Fermi points and proves the stronger bound $p_N\geq c>0$.
\Cref{cor:canonical-weight-continuity} supplies the
required source continuity for convergent weights. A separate
convexity argument there removes the need for a unique or
continuous selection of the tuned weights or Brouwer fixed point.
The canonical local-limit theorem resolves the finite-volume
identification at fixed weights. The sharp Ward/source comparison
for the height \GFF is proved in
\crefrange{sec:sharp-source-flow}{sec:height-gaussianity}.

\subsubsection{Prescribed slopes and the fixed-sector local limit}
\label{sec:fixed-sector-limit}

The results of the preceding subsubsections concern the canonical
laws $\mu_N^\mathrm{can}$ of \cref{thm:canonical-local-limit}, in
which the numbers $N_r(M)$ of lozenges of each orientation are
allowed to fluctuate, and are only biased by the auxiliary weights
$\prod_r\omega_r^{N_r(M)}$. The surfaces we started from, however,
live on a torus with prescribed periods, so the counts of all
three lozenge orientations are prescribed exactly: in the notation
of \cref{thm:coefficient-rg}, they are $(n\Theta_1,n\Theta_2,n^2)$. The law of
interest is therefore the microcanonical one, the canonical law
conditioned on the sector $\cC_N$ of \cref{cor:sector-source-reduction}.
That corollary reduces the identification of its local limit to
a single estimate, namely a lower bound on the canonical probability
$p_N=\mu_N^\mathrm{can}(\cC_N)$ of the target sector which beats
the source error $N^{-\kappa}$. The same issue arose in
\cite[\S8.6]{LaslierLubetzky24}, where it was resolved by tuning the
auxiliary weights so that the target counts become typical.
We follow the same strategy, but the auxiliary weights can now be
tuned without a differentiable choice of dressed dispersion, which
is not available in the present setting.

The plan is as follows. For free dimers, the counts are typical
exactly when the zeros of the dispersion sit at explicit points
$\pm p(\rho)$ determined by the density vector $\rho$. We first
prescribe the zeros of the matched dispersion $D$ at these points,
solving for the weights and the quadratic data in one
fixed-point problem (\cref{prop:fixed-zero-tuning}). We then compare
the resulting signed partition function, whose leading term is a
product of $D$ over the Fourier points of a boundary character,
with a free dimer determinant having the same zeros
(\cref{lem:same-zero-determinants}); the free torus winding theorem
of Kenyon, Sun and Wilson~\cite{KSW16} then gives the required
lower bound on the sector mass (\cref{lem:free-central-sector,prop:interacting-central-sector}).
Finally, the tuned weights need not be unique, and may vary with
$N$; convexity of the pressure and mixing show that any two
canonical laws with the same density coincide
(\cref{lem:pressure-equal-density}), proving the fixed-sector
limit (\cref{thm:fixed-sector-limit}). Throughout this subsubsection
$s_0$ is fixed, $\beta$ is sufficiently large, and all density
vectors lie in a fixed compact subset of the open simplex.
Constants are uniform on that set and on the bounded-shape tori of
the preceding subsubsection.

Normalize weights by $\omega_1=1$. For a density vector
$\rho=(\rho_1,\rho_2,\rho_3)$, put
\[
 \bar\omega(\rho)=\Big(1,
       \frac{\sin(\pi\rho_2)}{\sin(\pi\rho_1)},
       \frac{\sin(\pi\rho_3)}{\sin(\pi\rho_1)}\Big),\qquad
 p(\rho)=\big(\pi(1-\rho_3),-\pi(1-\rho_2)\big)\,.
\]
Recall that the free dispersion at weights $\omega$ is
$\mu_\omega(k)=1+\omega_2e^{ik_1}+\omega_3e^{ik_2}$. The triangle
construction of Kenyon (the three weights are the sines of the
angles of a triangle, and the two zeros of $\mu_{\bar\omega}$ are
read off its angles), or direct substitution, gives
$\mu_{\bar\omega}(p)=0$. Indeed, since $\rho_1+\rho_2+\rho_3=1$,
\[
 \sin(\pi\rho_2)e^{ip_1}+\sin(\pi\rho_3)e^{ip_2}
 =-\big(\sin(\pi\rho_2)e^{-i\pi\rho_3}
          +\sin(\pi\rho_3)e^{i\pi\rho_2}\big)
 =-\sin(\pi(\rho_2+\rho_3))=-\sin(\pi\rho_1)\,,
\]
as the imaginary parts cancel. These are the free weights of
density $\rho$: under the free dimer law at weights
$\bar\omega(\rho)$, the probability that a given black vertex is
matched by an edge of type $r$ is exactly $\rho_r$, by the local
statistics of Kenyon, Okounkov and Sheffield~\cite{KOS06}. This also
follows by differentiating the Jensen formula for the free
pressure in $\log\omega_r$, which produces the angle $\pi\rho_r$
divided by $\pi$. The two real vectors $e^{ip_1},e^{ip_2}$ are
linearly independent, uniformly on the compact density set:
they are parallel exactly when $p_1-p_2=\pi(1+\rho_1)$ is an
integer multiple of $\pi$, and $\rho_1$ is bounded away from
$0$ and $1$ on that set.

\begin{proposition}[Simultaneous tuning of weights and quadratic data]
\label{prop:fixed-zero-tuning}
For every such $\rho$ there are positive liquid weights
$\omega=\omega(\rho,\beta)$ and a physically matched dispersion
$D$ of \cref{thm:matched-bulk-flow} such that
\begin{equation}\label{eq:fixed-zero-tuning}
 D^{-1}(0)=\{p(\rho),-p(\rho)\},\qquad
 |\omega-\bar\omega(\rho)|+
       \|D-\mu_{\bar\omega(\rho)}\|_{C^3}\leq Cd_\beta\,.
\end{equation}
The constants are uniform. No unique or continuous selection of
$\omega$ or $D$ is asserted; \cref{lem:pressure-equal-density}
below shows that none is needed.
\end{proposition}

The proof uses the same matching map, with the weights now moving
so that the dressed zeros stay at their prescribed positions.

\begin{proof}
Fix $p=p(\rho)$ and the uniformly invertible real-linear map
$T_p(a,b)=ae^{ip_1}+be^{ip_2}$. Set
\[
 (\omega_2(y),\omega_3(y))
       =(\bar\omega_2,\bar\omega_3)+s_0T_p^{-1}y_0.
\]
Then $\mu_{\omega(y)}(p)=s_0y_0$ exactly, so the chart
\cref{eq:physical-parameter-chart} has $p(y)=p$ and
$\ell(y)=\nabla\mu_{\omega(y)}(p)-y_1$. For small $y$ the weights
stay in a common compact liquid set.

The proofs of \cref{lem:moving-first-integration,lem:critical-data-moving-dispersion}
apply jointly to these weights and $y$: fast denominators,
differentials and endpoint reserves have common bounds; truncated
coefficients are continuous; the same summable envelopes remove
truncations. Thus the matching map becomes the continuous map
\[
 \Psi(y)=\big(\nu^*_{0,+}(\omega(y),y)
                  -s_0^{-1}\sigma_{\omega(y),y}(p),\,
             v^*_{0,+}(\omega(y),y)
                  -\nabla\sigma_{\omega(y),y}(p)\big),
 \qquad |\Psi(y)|\leq C_0d_\beta.
\]
Apply the Brouwer argument of \cref{thm:matched-bulk-flow} on
$|y|\leq Md_\beta$, $M>C_0$. Its verification of the physical
jets and exact Gaussian compensation applies unchanged.
The resulting zeros are $\pm p$ by construction. Finally
$|\omega-\bar\omega|\leq Cs_0d_\beta$ and
$\|D-\mu_\omega\|_{C^3}\leq C_{s_0}d_\beta$ by
\cref{eq:physical-dispersion}; since $s_0$ is fixed, these imply
\cref{eq:fixed-zero-tuning}.
\end{proof}

\begin{lemma}[Comparison of determinants with the same zeros]
\label{lem:same-zero-determinants}
Suppose $D$ satisfies \cref{eq:fixed-zero-tuning}, and write
$\bar\mu=\mu_{\bar\omega(\rho)}$. There is a scalar $c_D$,
$|c_D|\leq Cd_\beta$, such that, uniformly over all boundary
characters,
\[
 \prod_{k\in\cP_N(\theta)}D(k)
 =e^{V_Nc_D+\eta_N^\theta}
                   \prod_{k\in\cP_N(\theta)}\bar\mu(k),
                 \qquad |\eta_N^\theta|\leq Cd_\beta\,.
\]
For the products with the nearest point to each zero deleted,
the same assertion holds with another remainder of size
$Cd_\beta$. In particular $p_D=p_{\bar\mu}+\Re c_D$.
At a singular twist the full-product equality means that both
sides vanish; no quotient at a zero is required.
\end{lemma}

The point of the lemma is this. The leading term of the
interacting partition function at a boundary character $\theta$ is,
up to the bulk scalar $e^{V_N\Delta}$ of
\cref{prop:last-mode-source-formula}, the product of the matched
dispersion $D$ over the Fourier points $\cP_N(\theta)$, whereas the
object whose sector mass we understand is the free determinant,
the product of $\bar\mu$ over the same points. The two products
are comparable only because $D$ and $\bar\mu$ vanish at the same
two points, so that their ratio stays bounded near the zeros; this
is what \cref{prop:fixed-zero-tuning} arranged.

\begin{proof}
Let $h=D-\bar\mu$ and $d(k)=\dist(k,\{p,-p\})$.
Both dispersions vanish at these points. Hence $h$ vanishes at
$\pm p$ and has $C^1$ norm at most $Cd_\beta$ by
\cref{eq:fixed-zero-tuning}, so that $|h(k)|\leq Cd_\beta d(k)$
by the mean value theorem; and $\bar\mu$ has exactly two zeros,
both simple, so that $|\bar\mu(k)|\geq c d(k)$ on the torus, with
$c>0$ uniform on the compact density set.
Consequently, on the punctured torus,
\[
 u(k)=\log\Big(1+\frac{h(k)}{\bar\mu(k)}\Big),\qquad
 |u(k)|\leq Cd_\beta,\qquad
 |\partial^l u(k)|\leq C_l d_\beta d(k)^{-l}\quad(1\leq l\leq3)\,,
\]
using the branch defined by $\log(1+z)$ for $|z|<1/2$; the
derivative bounds follow from the quotient rule, each derivative
of $h/\bar\mu$ costing one power of $d(k)^{-1}$.
Set $c_D=\int u(k)\,\d k/(2\pi)^2$, so that $|c_D|\leq Cd_\beta$,
and note that $D=\bar\mu e^{u}$ on the punctured torus, so that
the asserted product identity is the statement that the Riemann
sum of $u$ over $\cP_N(\theta)$ equals $V_Nc_D$ up to an error
$O(d_\beta)$.

Apply the centered-cell proof of
\cref{lem:deleted-determinant-comparison} to $u$.
The boundedness of $u$ makes the near-zero cells contribute
$O(d_\beta)$, with no logarithmic correction. The derivative
bounds multiply all the other errors in that proof by $d_\beta$:
the cubic remainder and the Hessian sampling error are summable,
while the signed Hessian integral is a boundary integral with
length $O(N^{-1})$ and integrand $O(d_\beta N)$. Hence
\[
 \sum_{k\in\cP_N(\theta)}u(k)=V_Nc_D+O(d_\beta)\,,
\]
when no point is a zero. If a point is a zero, set $u=0$
there; the estimate is unchanged
and the asserted product equality is still valid. Deleting two
points changes this sum by $O(d_\beta)$, proving the deleted
version without any zero assignments. Taking the real part of
the integral identity proves the pressure relation.
\end{proof}

\begin{lemma}[Uniform free central-sector mass]
\label{lem:free-central-sector}
Let $\bn=(n_1,n_2,n_3)$ be attainable orientation counts
on $\Lambda_N$, with $\sum n_r=V_N$, and let
$\rho=\bn/V_N$ lie in the prescribed compact simplex.
For free dimers at weights $\bar\omega(\rho)$,
\begin{equation}\label{eq:free-central-sector}
 \mu_{N,\bar\omega}^\mathrm{free}
                    (N_r=n_r\text{ for all }r)\geq c_0>0\,,
\end{equation}
for all sufficiently large $N$, uniformly in the counts and shape.
\end{lemma}
\begin{proof}
We use the free torus winding theorem of
\cite[Thm.~4 and Cor.~4.4]{KSW16}, which describes the law of the
homology class of a free torus matching as an asymptotically
discrete Gaussian. The proof has four parts: we identify the
counts with the homology class and locate the center of that
discrete Gaussian at the target counts; we check that its
precision matrix is uniformly bounded above and below; we verify
that the asymptotic statement of \cite{KSW16} holds uniformly over
the weights, torus shapes and counts allowed here, since
fixed-weight convergence alone would not suffice (the weights
$\bar\omega(\rho)$ vary with the counts); and we conclude by a
compactness argument.

Use the all-type-1 matching as reference and identify its white
and black cell coordinates. Summing the displacements of the
oriented difference cycles gives $(N_2,N_3)$ in these coordinates.
It is a period vector, and their total homology is
\[
 H_N=\sB_N^{-1}(N_2,N_3)^{\trps}\in\Z^2\,.
\]
Thus fixing $H_N$ fixes the counts. In this convention the center
parameter in the free winding theorem is
\[
 m_N=\sB_N^{-1}V_N(\rho_2,\rho_3)^{\trps}\,.
\]
Here is an explicit check of the convention. Choose
$\det\sB_N>0$ by ordering the period generators. The
row-period matrix of \cite{KSW16} is
$E=\sB_N^{\trps}=\bigg(\begin{smallmatrix}u&v\\x&y\end{smallmatrix}\bigg)$.
Its reference-oriented polynomial is
$1-\bar\omega_2z-\bar\omega_3w$; the root angles in
\cite[Eq.~(57)]{KSW16} are $(r_0,s_0)=(\rho_3,-\rho_2)$.
Here $s_0$ is only the root-angle notation of that formula,
not our infrared cutoff. By \cite[Cor.~4.4]{KSW16}, the center is
\[
 (xr_0+ys_0,-ur_0-vs_0)^{\trps}
       =-\adj(\sB_N)(\rho_2,\rho_3)^{\trps}\,.
\]
Reversing the orientation of the difference cycles gives our
$H_N$ convention and hence the displayed $m_N$, since
$\adj(\sB_N)=V_N\sB_N^{-1}$.
The common all-type-1 reference has zero homology, so there is
no further integer offset. Attainability gives $m_N\in\Z^2$.

The theorem says that $H_N-m_N$ has asymptotically the discrete
Gaussian with masses proportional to
\begin{equation}\label{eq:free-centered-gaussian}
          \exp(-\tfrac12 j^{\trps} Q_Nj),\qquad j\in\Z^2,
                         \qquad cI\leq Q_N\leq CI\,.
\end{equation}
The covariance parameter in \cite[Eq.~(44)]{KSW16} and its inverse
are uniformly bounded: the honeycomb Hessian in
\cite[Eq.~(58)]{KSW16} is uniformly positive definite on the
compact liquid set, and $\sB_N/N$ and its inverse are bounded.
This gives the asserted bounds on the Gaussian precision matrix $Q_N$.

Here is why the asymptotic statement is uniform in the parameters
used here. For the honeycomb polynomial and its conjugate product,
the two zeros are uniformly separated, their linear maps have
uniformly bounded inverses, and all fixed-order derivatives are
uniformly bounded. The quadrature errors in
\cite[Lem.~6.3--6.5]{KSW16} therefore have common constants:
the local Taylor comparison and its fourth-derivative bound use
only these bounds, and the complementary region has a common
positive lower bound. More quantitatively, let $K$ be a slightly
enlarged compact set of weights and normalized period matrices.
In the coordinates of their \S6, the mesh is comparable to
$N^{-1}$, uniformly on $K$. Near each node the characteristic
polynomial $P$ is at least $c_K|x|^2$, its quadratic part $P_2$
has eigenvalues in $[c_K,C_K]$, and
\[
 |\log P(x)-\log P_2(x)|\leq C_K|x|,\qquad
 |\partial^\alpha\log P(x)|\leq C_K|x|^{-4}
                       \quad(|\alpha|=4).
\]
The same estimates hold for the quadratic comparison polynomial
of their Eq.~(61). With $\mu=N^{-4/5}$, the inner comparison
and the outer Simpson-rule error in Lemma~6.3 are bounded by
\[
 C_K\big(N^2\mu^3+N^{-2}\mu^{-2}\big)
                                      \leq C_KN^{-2/5}.
\]
The first term sums the $O(|x|)$ error in the node boxes;
the second is $N^{-2}\int_{|x|>\mu}|x|^{-4}\,\d x$.
All other regions have uniform smooth bounds. In Lemmas~6.4
and~6.5, the transformed Hessians stay in a compact positive
set, so $\Im\tau$ is bounded above and away from zero.
The only nonuniform denominator is the distance $r$ from a
singular character in Lemma~6.5. Thus their error is uniformly
$C_KN^{-2/5}/r$. Away from a boundary character within
$N^{-1/5}$ of a singular one, the error
$O(N^{-2/5}/r)$ in \cite[Thm.~1]{KSW16} is $O(N^{-1/5})$.
Near such a character, \cref{lem:deleted-determinant-comparison}
with $D=\bar\mu$ gives instead
\[
 e^{-V_Np_{\bar\mu}}
       \bigg|\prod_{k\in\cP_N(\theta)}\bar\mu(k)\bigg|
       \leq C\,r_+(\theta)r_-(\theta)\,,
\]
where $r_\pm$ is the distance of the boundary character to that
of the corresponding zero. Indeed the selected eigenvalues
are at most $Cr_\pm/N$ and the deleted product is at most
$CN^2e^{V_Np_{\bar\mu}}$. The limiting theta-function factors
have the same bound: their holomorphic theta factors vanish
simply and have uniformly bounded derivatives, and the remaining
prefactors are bounded on the compact set of conformal shapes.
Writing $r=\min(r_+,r_-)$, both normalized terms are at most
$C_KN^{-1/5}$ when $r\leq N^{-1/5}$; when $r>N^{-1/5}$,
the relative error just proved and boundedness of the limiting
factor give the same absolute error. Thus the determinant
asymptotic has absolute error at most $C_KN^{-1/5}$ uniformly
also at singular characters. Summing the four signed
terms as in \cite[\S1.3.2 and Prop.~2.2]{KSW16} gives a uniform error for
the positive partition function. Its limiting finite-size factor
is the positive Gaussian sum \cite[Thm.~2(b), Eq.~(11)]{KSW16}, uniformly bounded above
and below on this compact shape and phase set.

This also applies after changing logarithmic weights by
$\sB_N^{-\trps}t$, with $t$ in a fixed real neighborhood of zero.
The free pressure has uniformly bounded third derivatives on a
slightly larger compact liquid set, so its cubic Taylor remainder,
multiplied by $V_N$, is $O(N^{-1})$. The variations of the zero
locations and Hessians in \cite[Lem.~4.2--4.3]{KSW16} likewise
have common constants. In the Gaussian sums used there, the
tails outside $|j|\leq C\sqrt{\log N}$ are uniformly $o(1)$,
by \cref{eq:free-centered-gaussian}; on the retained terms the
Taylor errors are at most $C_KN^{-1}(1+|j|^3)$, hence at most
$C_KN^{-1}(1+(\log N)^{3/2})$. Their computation
of the Laplace transform consequently yields, uniformly for
such real $t$,
\[
 \E_\mathrm{free}e^{t\cdot(H_N-m_N)}
  =\frac{\sum_{j\in\Z^2}e^{-j^TQ_Nj/2+t\cdot j}}
          {\sum_{j\in\Z^2}e^{-j^TQ_Nj/2}}+o(1)\,.
\]
These estimates explain the needed uniformity directly from the
free theorem's proof, including its possible singular twists.

To conclude without presuming a fixed limiting shape, suppose
\cref{eq:free-central-sector} failed along a sequence. Pass to
a subsequence on which $Q_N$ converges to a positive definite
$Q$. The displayed Laplace transforms converge on a neighborhood
of zero to those of the discrete Gaussian with matrix $Q$;
hence the laws converge weakly to it. The open square
$(-1/2,1/2)^2$ has boundary of zero limiting mass and contains
only the lattice point zero. Its probability therefore converges
to $(\sum_j e^{-j^TQj/2})^{-1}$, which is bounded below by
$(\sum_j e^{-c|j|^2/2})^{-1}>0$. This contradiction proves a
uniform positive $c_0$.
\end{proof}

\begin{proposition}[Positive mass of the interacting target sector]
\label{prop:interacting-central-sector}
Write $\bN=(N_1,N_2,N_3)$ for the random vector of orientation counts.
For attainable counts $\bn$ as above, choose any weights
and matched dispersion supplied by
\cref{prop:fixed-zero-tuning} for
$\rho=\bn/V_N$. Then
\begin{equation}\label{eq:interacting-central-sector}
             \mu_{N,\omega}^\mathrm{can}(N_r=n_r\text{ for all }r)
                                      \geq c>0\,.
\end{equation}
In particular, for each fixed finite source set,
\begin{equation}\label{eq:tuned-sector-source}
 \sup_A\bigg|\mu_{N,\omega}^\mathrm{can}
       \bigg(e^{\sum_e A(e)\one_{\{e\in M\}}}\given\bN=\bn\bigg)
                       -e^{S_\infty(A;\omega)}\bigg|
                                 \leq C_FN^{-1/2}\,.
\end{equation}
Both estimates are uniform over all choices of the tuned data.
\end{proposition}

\paragraph{Proof idea}
Recall from \cref{cor:sector-source-reduction} that the sector
partition function is extracted from the signed character sum by
the normalized two-circle integral $\mathcal E_{\cC_N}$. By the
last-mode formula, the interacting partition function at each
character is, up to a bulk scalar and a small error, the product of
$D$ over the Fourier points; by \cref{lem:same-zero-determinants}
this product is, in turn, a constant multiple of the free
determinant at the weights $\bar\omega(\rho)$. Applying the
extraction operator to both sides, the leading term is a constant
multiple of the free sector partition function, which
\cref{lem:free-central-sector} bounds below by a constant fraction
of the free partition function, and the error is a small fraction
of the canonical normalization. The lower bound follows. The one
subtlety is that the comparison must be made in a common gauge, so
that the extracted coefficients, and not just their moduli, are
compared; we address this first.

\begin{proof}
We keep track of the determinant phase before making a comparison.
In a boundary-character Fourier basis, use the same cell ordering
for black and white vertices and pair their Fourier transforms.
Their determinants cancel. Permuting the two retained black-white
pairs has even Grassmann parity. Thus the phase
$\xi_N^\theta$ in \cref{eq:last-mode-prefactor} is the same
geometric gauge and ordering phase as for the free determinant;
it can be set to one in this basis. Reversing to the edge-phase
gauge multiplies both expressions by the same unit phase. There
is no interaction-dependent phase to discard: the Gaussian
rescalings cancel as in that proposition, and their remaining
complex scalar logarithms are retained in $\Delta,r_N^\theta$.

At each selected mode, \cref{eq:last-mode-bounds} implies
\[
 a_\omega=ND(k_\theta^\omega)+O(N^{-1+Cd_\beta^2})\,.
\]
The bound follows by writing their difference as
$NR_\omega(q_\omega)(Z_J^{-1}-1)$; no eigenvalue is divided
out. For large $\beta$ the error is $O(N^{-1/2})$.
Expanding $b_N^\theta$ and $e^{r_N^\theta}$ therefore gives,
in the just specified common gauge,
\begin{align*}
 \fZ_N^\theta(0)
   &=e^{V_N\Delta}\bigg\{
        \prod_{k\in\cP_N(\theta)}D(k)+E_N^\theta\bigg\},
       \\
 |E_N^\theta|&\leq
        C(d_\beta+N^{-1/2})e^{V_Np_D}\,.
\end{align*}
Here we used the deleted-product bound and the factor $N^{-2}$.
In particular the estimate includes twists where the leading
product vanishes.

Write $\bar\omega=\bar\omega(\rho)$. The preceding determinant
lemma, and $|\det K_{\bar\omega}^\theta|\leq Ce^{V_Np_{\bar\mu}}$,
now imply
\[
 \fZ_N^\theta(0)
 =e^{V_N(\Delta+c_D)}\det K_{\bar\omega}^\theta
                       +\cE_N^\theta,\qquad
 |\cE_N^\theta|
                 \leq C(d_\beta+N^{-1/2})\cM_N\,,
\]
where $\cM_N=e^{V_N(p_D+\Re\Delta)}$ as before.
For circle edge phases the same formula holds with their common
geometric phase restored; this is exactly the gauge appearing
in the coefficient extraction of \cref{it:rg-sources}.

Extract the target coefficient by the normalized two-circle
integral $\mathcal E_{\cC_N}$ of \cref{cor:sector-source-reduction},
applied to the four-character combination. On the interacting
side, the extracted coefficient is, up to the common unit phase,
the positive weighted sector partition function
$Z_N^\mathrm{sec}(\omega)$, the sum of
$F_{\beta,\Gamma}(M)\prod_r\omega_r^{N_r(M)}$ over the matchings
with the target counts. On the free side, the extracted
coefficient of the determinant combination is
$e^{V_N(\Delta+c_D)}$ times the free sector partition function
$Z_N^\mathrm{sec,free}(\bar\omega)$, again up to the same unit
phase: the sign with which a matching enters a torus Kasteleyn
determinant depends only on its homology class, which is fixed by
the counts, so the determinant sign is constant on that sector
and the four-character combination reproduces the positive sum.
Since the circle integral has total mass one, applying
$\mathcal E_{\cC_N}$ to the previous display and using
$|\mathcal E_{\cC_N}\cE^\theta_N|\leq\sup_\theta|\cE_N^\theta|$,
the triangle inequality yields
\[
 Z_N^\mathrm{sec}(\omega)
 \geq e^{V_N\Re(\Delta+c_D)}Z_N^\mathrm{sec,free}(\bar\omega)
                    -C(d_\beta+N^{-1/2})\cM_N\,.
\]
By \cref{lem:free-central-sector}, the free sector partition
function is at least $c_0$ times the free partition function at
$\bar\omega$, and the free instance of
\cref{eq:canonical-normalization-bound} bounds the latter below by
$ce^{V_Np_{\bar\mu}}$; since $\cM_N=e^{V_N(p_D+\Re\Delta)}$ and
$p_D=p_{\bar\mu}+\Re c_D$, the first term is at least $c_*\cM_N$,
with a uniform $c_*>0$.
Choose $\beta$ so that $Cd_\beta<c_*/4$, then $N$ large so
that $CN^{-1/2}<c_*/4$. The sector partition function is then at
least $c_*\cM_N/2$. Dividing by the canonical upper bound
$\cZ_N^\mathrm{can}(0)\leq C\cM_N$ in
\cref{eq:canonical-normalization-bound} proves
\cref{eq:interacting-central-sector}. In particular
$N^{-1/2}/p_N\to0$, so \cref{cor:sector-source-reduction}
applies and proves \cref{eq:tuned-sector-source}.
\end{proof}

\begin{lemma}[Pressure and uniqueness of the law at a given density]
\label{lem:pressure-equal-density}
On a convex compact neighborhood of any liquid weight vector, use
$u=(\log\omega_2,\log\omega_3)$ and keep the canonical coefficients
fixed. The bulk pressure $P_\beta(u)$ is convex and continuously
differentiable, and
\[
 \log\cZ_N^\mathrm{can}(u)=V_NP_\beta(u)+O(1),\qquad
 \nabla P_\beta(u)=
   \big(\mu_\infty^\mathrm{can}(u)(\one_{e_2}),
        \mu_\infty^\mathrm{can}(u)(\one_{e_3})\big)\,.
\]
Here $e_r$ is the type-$r$ edge at a fixed black vertex, and the
$O(1)$ is uniform. If two weights in this neighborhood have the
same density vector, their canonical infinite-volume laws coincide.
This conclusion does not require strict convexity of $P_\beta$.
\end{lemma}

The lemma addresses the nonuniqueness left open by
\cref{prop:fixed-zero-tuning}: the tuned weights for the target
counts are not canonical, and along a sequence of tori they may
have several limit points. Every such limit point has the same
density, and the lemma shows that this already determines the
infinite-volume law. The natural route, strict convexity of the
pressure, would identify the weights themselves, but it is not
available without a differentiable fixed point; the argument below
replaces it by an affinity bound between the two laws together
with ergodicity.

\begin{proof}
The canonical normalization bound gives
$\log\cZ_N^\mathrm{can}(u)=V_NP_\beta(u)+O(1)$ uniformly on
the compact weight neighborhood, with $P_\beta=p_D+\Re\Delta$.
As a limit of finite-volume pressures, $P_\beta$ is convex
and independent of the matched dispersion. Translation invariance
and weight independence of coefficients give
\[
 \nabla\big(V_N^{-1}\log\cZ_N^\mathrm{can}(u)\big)
 =\big(\mu_N^\mathrm{can}(u)(\one_{e_2}),
       \mu_N^\mathrm{can}(u)(\one_{e_3})\big).
\]
By \cref{thm:canonical-local-limit,cor:canonical-weight-continuity},
these gradients converge uniformly to the stated continuous
vector field. Integrating along line segments proves that it
is $\nabla P_\beta$.

If $\nabla P_\beta(u)=\nabla P_\beta(v)$, the supporting-plane
inequalities at the two endpoints are equalities, so convexity
makes $P_\beta$ affine on $[u,v]$. The same positive pattern
weight occurs in the two canonical laws, giving the Hellinger affinity
\begin{equation}\label{eq:canonical-affinity}
 \sum_M\sqrt{\mu_N^\mathrm{can}(u)(M)\mu_N^\mathrm{can}(v)(M)}
 =\frac{\cZ_N^\mathrm{can}((u+v)/2)}
        {\sqrt{\cZ_N^\mathrm{can}(u)\cZ_N^\mathrm{can}(v)}}
 \geq c>0.
\end{equation}
The volume terms in its logarithm cancel by affinity of the
pressure; the uniform $O(1)$ errors give $c$. Cauchy--Schwarz
on each fiber shows that marginalization increases affinity.
Local convergence therefore gives the same lower bound for
every pair of finite-window marginals of the limiting laws.

If these ergodic laws differed, a bounded cylinder function
would have different means under them. By the ergodic theorem,
thresholding its spatial averages between those means gives
cylinder events $E_L$ with $\mu_\infty(u)(E_L)\to1$ and
$\mu_\infty(v)(E_L)\to0$. Cauchy--Schwarz on $E_L$ and its
complement bounds the affinity of a window supporting $E_L$ by
\[
 \sqrt{\mu_\infty(u)(E_L)\mu_\infty(v)(E_L)}
 +\sqrt{\mu_\infty(u)(E_L^c)\mu_\infty(v)(E_L^c)}\longrightarrow0,
\]
contradicting \cref{eq:canonical-affinity}.
\end{proof}

\begin{theorem}[The dimer limit at every prescribed density]
\label{thm:fixed-sector-limit}
Let $\Theta_i/n\to\theta_i>0$ for $i=1,2$. For sufficiently large $\beta$, the
laws with weights $F_{\beta,\Gamma}$ restricted to counts
$(n\Theta_1,n\Theta_2,n^2)$ converge locally to a unique
translation-invariant mixing law $\nu_{\beta,\rho}$, where
\[
 \rho=\frac{(\theta_1,\theta_2,1)}{1+\theta_1+\theta_2}\,.
\]
The limit is independent of the approximating integer windings.
It is a canonical law $\mu_\infty^\mathrm{can}(u)$ with density
$\rho$, at some $u$ within $Cd_\beta$ of the logarithmic free weight
ratios for that density. Its edge cumulants are the bulk source coefficients
and satisfy \cref{eq:bulk-response-decay}. This is a local limit
for the dimer measure; the \GFF and the physical height transfer
are not asserted here.
\end{theorem}
\begin{proof}
Apply the results for general tori with $N=n$. The diagonal tori have
bounded shape along the stated sequences, and their target counts
and volume are
\[
 \bn_N=(n\Theta_1,n\Theta_2,n^2),\qquad V_N=n(n+\Theta_1+\Theta_2)\,.
\]
Their density vectors $\rho_N$ tend to $\rho$. Choose tuned
weights $u_N$ by \cref{prop:fixed-zero-tuning}.
For large $N$ they all lie in one closed convex liquid neighborhood
of the free weights for $\rho$, since $|\omega-\bar\omega(\rho_N)|\leq Cd_\beta$
and $\bar\omega(\rho_N)\to\bar\omega(\rho)$. Their conditioning on
the target counts is exactly $\nu_{\beta,\Gamma}$, independent of
$u_N$: on the sector the auxiliary factor
$\prod_r\omega_r^{N_r(M)}$ is constant and cancels between the
numerator and the normalization, leaving the weights
$F_{\beta,\Gamma}$ restricted to the counts $(n\Theta_1,n\Theta_2,n^2)$.
\Cref{prop:interacting-central-sector} gives the uniform
source comparison \cref{eq:tuned-sector-source}.

Any subsequence of $u_N$ has a further subsequence converging to
some $u$. The uniform source estimate and
\cref{cor:canonical-weight-continuity} identify the local
limit on that subsequence as $\mu_\infty^\mathrm{can}(u)$.
The conditional finite-volume law is translation invariant, so
its probability of the type-$r$ edge at any black vertex is
exactly $n_r/V_N$. Thus the limiting canonical law has density
$\rho$, or $\nabla P_\beta(u)=(\rho_2,\rho_3)$.
\Cref{lem:pressure-equal-density} shows that all possible
subsequential limits are the same law. Compactness of the weight
neighborhood and the local source comparison then prove convergence
of the full sequence. The bulk canonical functional and pressure
are independent of the torus sequence: its non-small completion
coefficients disappear in the bulk limit, as proved in
\cref{lem:finite-size-source-kernels}. The same equal-density
argument therefore compares any two sequences of integer windings.
Mixing and the cumulant estimates follow from
\cref{thm:canonical-local-limit}.
\end{proof}

This proves the prescribed-sector identification, including varying
slopes, without a uniqueness or differentiability
claim for the quadratic fixed point. No conclusion about logarithmic
height covariance follows from the coarse edge-cumulant bounds.
The higher-cumulant comparison in \cref{sec:height-gaussianity} is
required for the dimer \GFF. We next prove the interacting
witness estimate in \cref{sec:interacting-witnesses}, then the
physical transfer in \cref{sec:transfer}.

\subsection{Interacting witness tails}
\label{sec:interacting-witnesses}

The reconstruction of the original surfaces in \cref{sec:transfer}
is driven by the geometric witnesses of
\cref{sec:preimage-witnesses}: a large deviation of the surface
from its tiling forces a large witness pattern in the tiling.
The first integration of \cref{sec:first-integration} already
exploits the rarity of such patterns, but only under the free
dimer law, through the determinant bound of
\cref{prop:fat_pattern}. For the reconstruction we need the
corresponding exponential tails under the interacting tiling law,
that is, under the canonical completion and under its sector
conditioning; this is the third conclusion of
\cref{thm:coefficient-rg}, which we prove here.

The natural way to compute the probability of a pattern under the
interacting law is to insert its Grassmann monomial into the
generating function and carry the insertion through the \RG as
a source. This is not quite correct for the present purpose. A
witness $W$ need not be one of the patterns selected in the pattern
expansion of $F_{\beta,\Gamma}$, and an expansion pattern $P$ may
share edges with it. In the Grassmann representation, the monomials
$\Psi_P$ and $\Psi_W$ then contain a common pair of variables, so
their product vanishes by nilpotency, although the event that both
$P$ and $W$ occur in the matching may have positive weight. Thus
inserting just the monomial $\Psi_W$ would discard exactly those
terms of the expansion overlapping the witness, and would compute
a different quantity. We first give the exact
insertion which retains those patterns, by merging the witness with
the expansion patterns it overlaps into a single monomial.
All constants in this subsection are uniform on a compact set of
positive liquid weights and bounded torus shapes. The torus size is
sufficiently large depending on $\beta$, but no bound on the size or
diameter of the inserted pattern is imposed; this is essential, as
the witnesses of \cref{sec:preimage-witnesses} can be as large as
the torus.
We normalize the largest edge weight to one. A common rescaling
multiplies every matching weight by the same factor and leaves
both the canonical law and its sector conditioning unchanged.

For this subsection, assume in addition \cref{it:rg-marked}, and let
$\scW_N$ be its witness family, a translation-covariant family of
connected patterns of size at least $B_\beta$. No positivity of
the individual pattern coefficients is required: the matching
weight is positive by \cref{it:rg-positive}, and the merged
coefficients defined below are bounded in absolute value by the
nonnegative sums \cref{eq:rg-absolute-insertion}.
For a fixed triangle $x$ and $s\geq1$, put
\[
 \cB_x(s)=\{M:\exists W\in\scW_N,\ x\in S(W),\
                                  |W|\geq s,\ W\subset M\},
\]
and define its witness count by
\[
 T_{x,s}(M)=\sum_{\substack{W\in\scW_N:\ x\in S(W)\\
                            |W|\geq s}}\one_{\{W\subset M\}}\,.
\]
Here $S(W)$ is its triangle support; in particular
$\one_{\cB_x(s)}\leq T_{x,s}$, since a matching in $\cB_x(s)$
contains at least one witness counted by $T_{x,s}$. We shall
bound the expectation of the count $T_{x,s}$, which is linear in
the witnesses and hence accessible to a single Grassmann
insertion, rather than the probability of the event $\cB_x(s)$
directly. This remains true for the witness events allowing
patterns of zero coefficient, since membership in $\scW_N$ does
not depend on the coefficients.

\begin{lemma}[Exact overlapping insertion]
\label{lem:overlapping-witness-insertion}
Write $c(P)=c_{\beta,\Gamma}(P)$. For a connected partial matching
$U$ define
\begin{equation}\label{eq:marked-union-coefficient}
 b_{x,s}(U)=
 \sum_{\substack{W\in\scW_N,\ x\in S(W),\ |W|\geq s\\W\subset U}}
 \ \sum_{\substack{\{P_i\}:\ S(P_i)\text{ pairwise disjoint}\\
       P_i\subset U,\ S(P_i)\cap S(W)\ne\varnothing\\
                        W\cup\bigcup_iP_i=U}}
                      \prod_i c(P_i)\,.
\end{equation}
The empty inner collection is permitted. Set
\[
 B_{x,s}(\psi)=\sum_U b_{x,s}(U)\Psi_U(\psi),\qquad
 \Psi_U=\prod_{e=(b,w)\in U}K(e)\psi_b\psi_w\,,
\]
in the paired convention.
In the exact unnormalized torus Grassmann identity, replacement of
$e^{V_\beta}$ by $e^{V_\beta}(1+tB_{x,s})$ replaces the weight
$F_{\beta,\Gamma}(M)$ by
$F_{\beta,\Gamma}(M)(1+tT_{x,s}(M))$.
This is a polynomial identity, also at singular twists.
Every nonzero $b_{x,s}(U)$ has
\begin{equation}\label{eq:marked-union-size}
                x\in S(U),\qquad |U|\geq\max(s,B)\,.
\end{equation}
\end{lemma}
\begin{proof}
Fix a matching $M$ and expand both sides. By the definition of
$F$ in \cref{eq:rg-pattern-completion} and of $T_{x,s}$, a term of
$F(M)T_{x,s}(M)$ is specified by a witness $W\in\scW_N$ with
$x\in S(W)$, $|W|\geq s$ and $W\subset M$, together with a family
$\{P_i\}$ of expansion patterns contained in $M$ with pairwise
disjoint supports, and it carries the weight $\prod_ic(P_i)$.
Split the family into the members meeting $S(W)$ and the others,
and merge $W$ with exactly the former, obtaining
\[
 U=W\cup\bigcup_{i:\ S(P_i)\cap S(W)\ne\varnothing}P_i\,,
\]
a partial matching contained in $M$ which is connected, since
each merged member meets the connected support $S(W)$. Every
non-merged member misses $S(U)$: it misses $S(W)$ by definition,
and it misses each merged $S(P_i)$ because the original family had
pairwise disjoint supports. Thus the data $(W,\{P_i\})$ determine
a connected $U\subset M$, a summand of
\cref{eq:marked-union-coefficient} for this $U$, and a family of
expansion patterns with pairwise disjoint supports contained in
$M\setminus U$ and missing $S(U)$. Conversely, these three data
recover $(W,\{P_i\})$ uniquely, since the merged members are
exactly those of the inner family in
\cref{eq:marked-union-coefficient}. Summing the weights, for every
matching $M$,
\[
 F(M)T_{x,s}(M)
 =\sum_{U\subset M\text{ connected}}b_{x,s}(U)
   \sum_{\substack{\{P_i\}:\ S(P_i)\text{ pairwise disjoint}\\
                   P_i\subset M,\ S(P_i)\cap S(U)=\varnothing}}
                                 \prod_ic(P_i)\,.
\]
This is precisely the expansion of $e^{V_\beta}B_{x,s}$ under the
unnormalized Berezin identity of \cref{lem:pattern-algebra}:
expanding $e^{V_\beta}=\prod_P(1+c(P)\Psi_P)$ and multiplying by
the monomial $b_{x,s}(U)\Psi_U$, the surviving terms are exactly
those whose patterns miss $S(U)$ and each other, by nilpotency,
and the Gaussian integral then completes the forced edges to a
matching $M$ exactly as in that lemma. The monomial
$\Psi_U$ assigns each forced edge its Kasteleyn weight once,
including edges shared by the witness and the merged members, so
each edge of $M$ carries one weight and the complete-matching
sign is unchanged. This proves the identity for the coefficient
of $t$; the coefficient of $t^0$ is the zero-insertion identity.

The union is essential even for a single overlapping pattern:
if $P$ and $W$ share an edge, $\Psi_P\Psi_W=0$ because its
Grassmann variables are repeated, although the event
$\{P,W\subset M\}$ may have positive weight. The correct term
is $c(P)\Psi_{P\cup W}$, which retains the coefficient of $P$
and forces the common edge only once.
No division by a determinant is used, so singular twists are
included. Finally, $U\supset W$ together with $x\in S(W)$,
$|W|\geq s$ and $|W|\geq B$ gives \cref{eq:marked-union-size}.
\end{proof}

After the signed character sum is converted back to the positive
canonical partition function, denote the result with this insertion
by $\cZ_N^{\mathrm{can}}[t]$. The identity just proved says
\[
 \frac{\cZ_N^{\mathrm{can}}[t]}{\cZ_N^{\mathrm{can}}[0]}
       =\frac{\sum_MF_{\beta,\Gamma}(M)\prod_r\omega_r^{N_r(M)}
                       (1+tT_{x,s}(M))}
             {\sum_MF_{\beta,\Gamma}(M)\prod_r\omega_r^{N_r(M)}}
       =1+t\,\mu_N^{\mathrm{can}}(T_{x,s})
                                      \pmod {t^2}\,.
\]
Thus the quantity to be bounded is the coefficient of $t$ after
normalization. The following two estimates keep the numerator
unnormalized through every Grassmann integration, and use the
positive canonical denominator only at the end, through its lower
bound \cref{eq:canonical-normalization-bound}; in this way no
sign or cancellation property of the intermediate Gaussian
determinants is needed.

The absolute value of \cref{eq:marked-union-coefficient} is bounded
by \cref{eq:rg-absolute-insertion}, by the triangle inequality and
$|c(P)|\leq q(P)$. Thus \cref{it:rg-marked} supplies its
\texttt{thin}/\texttt{fat} envelopes, including the rooted counts,
and the merged coefficient can be treated by the first integration
exactly like an interaction coefficient, with the one difference
that its support is forced to contain $x$ and to have size at least
$\max(s,B)$.

For the following bounds one mark is enough: since only the
coefficient of $t$ is needed, we work in the coefficient ring
$\C[t]/(t^2)$, so that every term contains the insertion at most
once, and no summation over multiple insertions or choice of a
source radius is involved.
For a polynomial in slow fields with a distinguished triangle $x$,
write
\begin{equation}\label{eq:single-insertion-norm}
 \|Q_x\|_{a,c;x}
   =|Q_x(\varnothing)|+
      \sum_{S\ne\varnothing}|Q_x(S)|c^{|S|}
                     e^{a\ell(S\cup\{x\})}\,.
\end{equation}
The coefficients here are in the same antisymmetric kernel
convention as the first-step norms; field flavours are also summed.
We use torus distance on a torus. The exponent and degree weight
are chosen with the reserves for two spatial moments and field
localization used in \cref{lem:source-endpoints-localization}.
For localization, identify a triangle with its honeycomb cell;
this changes spatial distances by at most a fixed constant.

\begin{proposition}[The first integration of a witness insertion]
\label{prop:witness-first-integration}
Perform the moving first integration of
\cref{lem:moving-first-integration}, with the matched
parameters held fixed. With the insertion of
\cref{lem:overlapping-witness-insertion}, its torus output is
\[
        e^{C_N^\theta+U_N^\theta}(1+tQ_{x,s})\,,
\]
where $C_N^\theta+U_N^\theta$ is the zero-source first-step output.
There are fixed $c_1,C>0$ and a norm as in
\cref{eq:single-insertion-norm} such that
\begin{equation}\label{eq:witness-insertion-smallness}
          \|Q_{x,s}\|_{a,c;x}
                   \leq C e^{-c_1\max(s,B_\beta)}\,.
\end{equation}
The estimate is uniform in the boundary character and includes
the scalar term. Noncontractible witnesses are included.
\end{proposition}
\begin{proof}
Work modulo $t^2$. The connected fast expansion has exactly
one marked node, with coefficient $b_{x,s}(U)$; every other node
is an ordinary interaction node. By \cref{it:rg-marked}, the same
patch and Gram estimate as in \cref{eq:node-envelope} bounds
the sum of marked \texttt{thin} and \texttt{fat} nodes of size $m$ by $Ce^{-pm}$,
with a fixed reserve and $m\geq\max(s,B_\beta)$. This uses the
original Kasteleyn operator, also for the moving cutoff.

Root the tree at the marked node. The ordinary branches cost at
most $e^{v|S(U)|}$ by \cref{lem:rooted-tree}; since
$|S(U)|=2|U|$, take $2v<p$. Since $x\in S(U)$, a tree inside
$U$ also pays the weight in \cref{eq:single-insertion-norm}.
Writing $m_0=\max(s,B_\beta)$, the resulting numerical estimate is
\[
 \|Q_{x,s}\|_{a,c;x}
       \leq C\sum_{m\geq m_0}e^{-(p-2v)m}
       =\frac{Ce^{-(p-2v)m_0}}{1-e^{-(p-2v)}}.
\]
The marked-node estimate already includes the rooted count of
patterns of size $m$. Thus there is no remaining translation or
decomposition factor in this sum. The same rooted expansion
includes the case in which all marked fields are contracted, so
it controls the scalar term as well. The torus patch and periodization estimates
of \cref{thm:torus-first-step} are uniform in the character
and apply to noncontractible supports as well. The compensating
quadratic term is independent of $t$.
\end{proof}

\begin{lemma}[Uniform bound for one localized insertion]
\label{lem:one-insertion-rg}
Keep the matched zero-source flow fixed. Suppose a torus
insertion after its first integration has equal black and white
degree in every monomial and has
$\|Q_x\|_{a,c;x}\leq q$ in
\cref{eq:single-insertion-norm}. The corresponding unnormalized
Grassmann partition function, denoted by $\fZ_N^\theta[t]$,
satisfies
\begin{equation}\label{eq:one-insertion-unnormalized}
       |[t]\fZ_N^\theta[t]|
                                \leq Cq\cM_N\,,
\end{equation}
uniformly in the character, including singular ones. Here
$\cM_N$ denotes the normalization in
\cref{lem:deleted-determinant-comparison}.
If the original insertion is a matching observable $T$, this implies
\begin{equation}\label{eq:one-insertion-positive}
                         |\mu_N^\mathrm{can}(T)|\leq Cq\,.
\end{equation}
The constant is independent of the support of the insertion.
\end{lemma}
\begin{proof}
Treat $x$ as one source position, with distinguished endpoint
$Q_x$. Its first moment pays the quadratic localization and its
higher degrees are irrelevant. The source-flow proof applies
linearly in this endpoint: $|\Gamma_0^Q|\leq Cq$ and the forcing
is $O(q2^{-3j/4})$, while the recursion matrices are unchanged.
Thus \cref{prop:source-flow-bound,lem:source-shell-summability}
give $|\Gamma_j^Q|\leq Cq2^{\eta_\beta j}$ and summable scalar
shells $Cq2^{-3j/4}$.

For the torus comparison, lift a minimizing spatial tree of each
monomial from $x$, with shortest representatives and its character
phases. Its periodization is $Q_x$ and its planar norm is at most
$q$, although it may depend on $N$ and the character.
\Cref{lem:finite-size-source-kernels}, including
\cref{eq:winding-pays-relevance}, applies uniformly to this
single endpoint. At fermion degrees $\fn=2,4$ it gives the bound
$Cq s_j^{\kappa-\fn/2}$ on its summed kernel. Each selected
last mode contributes $N^{-1/2}$, so the last-mode formula is
\[
 \fZ_N^\theta[t]
  =P_N^\theta\int\d\chi\,e^{H_N^\theta(\chi)}
              \big(1+t(C_N^Q+R_N^Q(\chi))\big),\qquad
 |C_N^Q|\leq Cq,\quad \|R_N^Q\|_\mathrm{alg}\leq CqN^{-\kappa}.
\]
Only degrees two and four survive. The bounded algebra norm of
$H_N^\theta$ and \cref{lem:deleted-determinant-comparison} give
the unnormalized bound, including singular characters. For a
matching observable, use the four-character identity and divide
by \cref{eq:canonical-normalization-bound}.
\end{proof}

\begin{theorem}[Interacting witness tails]
\label{thm:interacting-witness-tail}
For sufficiently large $\beta$, the canonical completion at
liquid weights in a fixed compact set satisfies
\[
                  \mu_N^\mathrm{can}(T_{x,s})
                            \leq Ce^{-c\max(s,B_\beta)}\,.
\]
For the prescribed-sector laws $\nu_n$, with physical
slopes in a fixed compact positive set, uniformly for all triangles
$x$, all $s\geq1$ and all sufficiently large $n$,
\begin{equation}\label{eq:interacting-witness-exponential}
                 \nu_n(\cB_x(s))
                            \leq Ce^{-c\max(s,B_\beta)}\,.
\end{equation}
In particular the witness tail is bounded by $C_q(1+s)^{-q}$
for every $q>0$. The estimate includes noncontractible witnesses and witnesses
whose generating support has zero activity.
\end{theorem}
\begin{proof}
For the canonical assertion, combine the three preceding steps:
the exact insertion of \cref{lem:overlapping-witness-insertion}
represents $\mu_N^\mathrm{can}(T_{x,s})$ as the normalized
coefficient of $t$; its first integration
(\cref{prop:witness-first-integration}) produces a localized
insertion of norm $q\leq Ce^{-c_1\max(s,B_\beta)}$; and
\cref{lem:one-insertion-rg}, applied with this $q$, gives
$|\mu_N^\mathrm{can}(T_{x,s})|\leq Cq$, which is the claim with
$c=c_1$.

For the prescribed-sector laws, recall that $\nu_n$ is the law
with weights $F_{\beta,\Gamma}$ restricted to the target counts.
Choose the tuned weights of \cref{prop:fixed-zero-tuning} for
these counts. They lie in a compact liquid set, so the canonical
assertion applies to them, and
\cref{prop:interacting-central-sector} gives target sector
probability $\mu_N^\mathrm{can}(\cC_N)\geq c_*>0$. The canonical
law conditioned on $\cC_N$ is exactly $\nu_n$, and is independent
of the auxiliary weights, since the factor
$\prod_r\omega_r^{N_r(M)}$ is constant on the sector. Positivity
of all the weights and $\one_{\cB_x(s)}\leq T_{x,s}$ therefore
give
\[
 \nu_n(\cB_x(s))
   =\frac{\mu_N^\mathrm{can}(\cB_x(s)\cap\cC_N)}
         {\mu_N^\mathrm{can}(\cC_N)}
   \leq\frac{\mu_N^\mathrm{can}(T_{x,s})}
                    {\mu_N^\mathrm{can}(\cC_N)}
   \leq\frac{C}{c_*}e^{-c\max(s,B_\beta)}\,,
\]
which is \cref{eq:interacting-witness-exponential}.
All support sizes were retained in the positive marked-node sums,
so this proof also controls witnesses growing with the torus.
In particular there is no exchange of a local limit with an
unbounded union of witness events, which would not be justified
for witnesses comparable to the torus. Averaging over $x$ and
using $e^{-cs}\leq C_q(1+s)^{-q}$ proves the last assertion.
\end{proof}

\subsection{Gaussian height fluctuations}
\label{sec:gaussian-height-fluctuations}

We sharpen the source estimates, identify the limiting height covariance,
and prove that the higher smeared cumulants vanish.

\subsubsection{Matching the running density and mass sources}
\label{sec:sharp-source-flow}

We now sharpen the coarse growth bound on the four source bilinears.
The normalized density coefficients will converge, whereas the mass
coefficients may have an anomalous exponent. Throughout this
subsubsection $\delta=d_\beta$, the matched dispersion is fixed, and
$r_*=1/8$ is the critical-history decay exponent.
The source labels are
\[
 a=(\omega,\omega'),\qquad
 a\in\cD\ \Longleftrightarrow\ \omega=\omega',
 \qquad a\in\cM\ \Longleftrightarrow\ \omega=-\omega'\,.
\]
The letters $\cD,\cM$ distinguish density and mass
bilinears, which pair a black and a white field at the same or
at opposite Fermi points respectively; they do not refer to
diagonal and physical coordinates. All phases $e^{\pm ip^\omega\cdot x}$ remain as
in \cref{eq:local-edge-source}.

We compare the lattice source recursion with the scalar reference
recursions and apply a stability estimate to their summable
difference. This transfers the continuum Ward-identity exponents
without dividing by a possibly vanishing lattice source coefficient.

We use the following established facts about the continuum reference
model of \cite[\S4]{GMT20}, with the same linear dispersions
$L_\omega$ and an adapted radial, smooth density--density potential.
Its infrared quartic coupling is an analytic function
$\lambda_*^R(g)=g+O(g^2)$ of its bare coupling $g$. There are
exponents $\eta$, $\eta_\mathrm{d}$ and $\eta_\mathrm{m}$ such that
its running field normalization $Z_j^R$, its running unnormalized
density-source coefficients ($a\in\cD$) and its running unnormalized
mass-source coefficients ($a\in\cM$), divided by $2^{\eta j}$,
$2^{\eta_\mathrm{d}j}$ and $2^{\eta_\mathrm{m}j}$ respectively, have
nonzero limits $1+O(g)$ with an error $O(|g|2^{-j/2})$. Moreover,
by \cite[Thm.~3]{BM02} (see also \cite[\S6.4.9]{GMT20}),
\begin{equation}\label{eq:continuum-source-input}
 \eta_\mathrm{d}=\eta,
 \qquad \eta=O(g^2),\qquad \eta_\mathrm{m}=O(g)\,.
\end{equation}
These are the reference-model statements used in
\cite[\S6.4.6--6.4.9, Eqs.~(6.103)--(6.109) and~(6.159)]{GMT20}.
In particular, the equality in \cref{eq:continuum-source-input}
is an input from the cutoff-removed continuum Ward identities,
with the cutoff qualification of \cref{rem:cutoff-ward}. We only use these continuum facts;
the extension to the present lattice endpoints is proved next.
Changing the fixed initial scale changes the limiting constants
and not these exponents or the stated geometric convergence.

\begin{lemma}[A reference model with the same infrared coupling]
\label{lem:matched-source-reference}
There is a real continuum bare coupling $g=O(\delta)$ for which
$\lambda_*^R(g)=\lambda_\infty$, the limiting physical coupling.
Choose the reference normalization initially to be one. With
scales and localization conventions matched, its histories satisfy
\begin{equation}\label{eq:source-reference-history}
 |\lambda_j-\lambda_j^R|+|z_{j+1}-z_{j+1}^R|
                         \leq C\delta2^{-r_*j}\,.
\end{equation}
There is a nonzero constant $a_Z=1+O(\delta)$ such that
\begin{equation}\label{eq:source-reference-Z}
 Z_j/Z_j^R=a_Z\big(1+O(\delta2^{-r_*j})\big)\,.
\end{equation}
After rescaling each bare reference bilinear so that its normalized
coefficient at scale zero equals one, let $q_{a,j}$ be its running
normalized coefficient. Then $q_{a,0}=1$, and
\begin{equation}\label{eq:reference-source-asymptotic}
 q_{a,j}=b_a2^{\xi_a j}
               \big(1+O(\delta2^{-j/2})\big),\qquad
 b_a=1+O(\delta),\qquad
 \xi_a=\begin{cases}0,&a\in\cD,\\
                     \eta_\mathrm{m}-\eta,&a\in\cM.
       \end{cases}\,.
\end{equation}
In particular $q_{a,j}$ never vanishes, $|\xi_a|\leq C\delta$,
and
$|q_{a,k}/q_{a,j}|\leq C2^{C\delta|j-k|}$.
The reference one-source flow is diagonal in the four labels $a$.
\end{lemma}
\begin{proof}
The inverse function theorem applied to $g\mapsto\lambda_*^R(g)$
gives the unique small real $g$ with the required limit. The physical
history approaches that limit at rate $2^{-r_*j}$ by
\cref{thm:critical-flow-fixed-dispersion}; the reference
history does so at rate $2^{-j/2}$. This gives the first term in
\cref{eq:source-reference-history}.

For the normalization increments compare their reference trees.
Replacing a quartic coefficient on scale $k$ costs at most
$C\delta2^{-3(j-k)/4}|\lambda_k-\lambda_k^R|$, by
\cref{lem:reference-history-derivatives}.
The initial smooth reference interaction differs from its local
quartic jet by an irrelevant endpoint with a finite second moment.
The resulting extra contribution to a normalization increment is
$O(\delta^2 2^{-3j/4})$: a linear quadratic jet needs two quartic
vertices, since the one-vertex contraction has zero momentum
derivative. Integrating the finitely many intervening ultraviolet
scales has the same marked bound. The convolution with
$2^{-r_*k}$ is at most $C2^{-r_*j}$. This proves the second term
in \cref{eq:source-reference-history}, with the weaker factor
$\delta$ sufficient here. No estimate on the absolute size of
$Z_j$ was used.

The ratios $1+z_{j+1}$ and $1+z_{j+1}^R$ lie in a fixed
neighborhood of one. Therefore
\[
 \log(Z_j/Z_j^R)=
 \sum_{k<j}\big[\log(1+z_{k+1})-\log(1+z_{k+1}^R)\big]\,,
\]
converges absolutely, with a tail $O(\delta2^{-r_*j})$ and a
sum $O(\delta)$. Exponentiation proves
\cref{eq:source-reference-Z}.

For completeness, diagonality here is a property of the reference
model, not an assumption on the lattice source recursion.
Independent constant changes
$b_\omega\mapsto e^{t_\omega}b_\omega$,
$w_\omega\mapsto e^{-t_\omega}w_\omega$ distinguish the two
mass bilinears from each other and from the densities. In the
adapted continuum coordinates the rotations used in
\cref{prop:conditional-beta-comparison} give
$b_\omega w_\omega$ rotation character $e^{i\omega t}$.
Thus they distinguish the two density bilinears as well.
The radial cutoff and potential preserve these symmetries;
localization at zero quasi-particle momentum commutes with them.
Consequently no reference local source coefficient mixes two
of the four labels. The same argument applies to the effective
source kernels obtained by integrating its ultraviolet scales.
Reflection exchanges the two labels within each class and gives
identical real exponents.

Divide the reference unnormalized source asymptotics by those
of $Z_j^R$. \Cref{eq:continuum-source-input} gives
\cref{eq:reference-source-asymptotic}. The finite normalization
at scale zero is $1+O(\delta)$ and is nonzero, so it can be
removed by rescaling the bare source. Uniformly small errors in
\cref{eq:reference-source-asymptotic}, including the finitely
many initial scales, prove nonvanishing and the ratio bound.
\end{proof}

\begin{lemma}[Comparison of the one-source recursion kernels]
\label{lem:one-source-kernel-comparison}
For each edge type $r$, the exact physical recursion can be written
\begin{equation}\label{eq:source-comparison-recursion}
 \Gamma_{j+1,r;a}=a_j^R\Gamma_{j,r;a}
       +\sum_{k=0}^j m^R_{a;j,k}\Gamma_{k,r;a}
       +e_{j,r;a},\qquad a_j^R=(1+z_{j+1}^R)^{-1}\,,
\end{equation}
where the same scalar reference kernels satisfy
\begin{equation}\label{eq:scalar-reference-recursion}
 q_{a,j+1}=a_j^Rq_{a,j}+\sum_{k=0}^j m^R_{a;j,k}q_{a,k},
 \qquad |m^R_{a;j,k}|\leq C\delta2^{-(j-k)/2}\,.
\end{equation}
For a sequence with the coarse bound
$|\Gamma_{j,r}|\leq C2^{\eta_\beta j}$ of
\cref{prop:source-flow-bound},
\begin{equation}\label{eq:source-comparison-error}
 |e_{j,r}|\leq C\delta2^{-(r_*-\eta_\beta)j}\,.
\end{equation}
These estimates include all nonlocal initial lattice source kernels;
the constants are uniform in their degree and spatial range.
\end{lemma}
\begin{proof}
Define $m^R_{a;j,k}$ by differentiating the reference integrations
in one bare source and retaining its local bilinear on scale $k$.
An initial source remainder is proportional to the same bare source;
since $q_{a,0}=1$, include it in $m^R_{a;j,0}$. A free quadratic
remainder has zero local jet and cannot generate one alone. Thus
every increment has an interaction endpoint, and the one-source
bound gives $|m^R_{a;j,k}|\leq C\delta2^{-(j-k)/2}$, with room
for normalization ratios. The preceding lemma gives diagonality.

Compare physical and reference kernels before substituting source
coefficients. Pure local reference graphs agree with their continuum
counterparts by \cref{lem:all-order-reference-sampling}: adding a
bilinear source keeps every vertex degree at most four and preserves
the no-aliasing margin. Local source jets do not enlarge momentum
support. The edge type $r$ is carried only by the source coefficient.

All other graphs have one of the marks covered by
\cref{lem:marked-source-comparison}. Initial ordinary and
interaction-derived source remainders cost $O(\delta)$. A free
source remainder, or a curvature mark, needs an additional ordinary
endpoint to contribute to a local source jet, supplying the same
factor. The remaining quadratic coefficients and the differences
of quartic and normalization histories cost
$O(\delta2^{-r_*t})$ on scale $t$, by
\cref{thm:matched-bulk-flow,eq:source-reference-history}; curvature
has the faster decay of \cref{lem:running-curvature}.
Holding the source at $k$, the two-path estimate therefore gives
\[
 C\delta2^{-(j-k)/4}
 \left[2^{-j/4}+\sum_{t=0}^j2^{-(j-t)/4}2^{-r_*t}\right]
 \leq C\delta2^{-(j-k)/4}2^{-r_*j},
 \qquad r_*=1/8.
\]
The convolution is bounded by $2^{-j/8}\sum_{d\geq0}2^{-d/8}$.
This includes off-diagonal physical kernels. The difference of the
propagated scalar coefficients is a $k=j$ term. Initial source
remainders not represented by a local source coefficient give
forcing $O(\delta2^{-j/4})$.

Substituting the physical source coefficients now yields
\[
 |e_{j,r}|\leq C\delta2^{-j/4}
       +C\delta2^{-r_*j}
                 \sum_{k=0}^j2^{-(j-k)/4}|\Gamma_{k,r}|.
\]
The coarse bound \cref{eq:source-growth-bound} and a geometric
sum prove \cref{eq:source-comparison-error}. The common marked
lemma supplies uniformity in degree and spatial range, including
normalization ratios and nonlocal source endpoints.
\end{proof}

\begin{lemma}[Stability of a scalar reference solution]
\label{lem:source-volterra-stability}
Suppose $q_j\ne0$ solves
$q_{j+1}=a_jq_j+\sum_{k\leq j}m_{j,k}q_k$, where
$|m_{j,k}q_k/q_{j+1}|\leq\epsilon2^{-c(j-k)}$.
Suppose $u_{j+1}=a_ju_j+\sum_{k\leq j}m_{j,k}u_k+e_j$
and $|e_j/q_{j+1}|\leq E2^{-b_0j}$.
For any $0<b<\min(b_0,c)$ and sufficiently small $\epsilon$
(depending only on $b,c$), there exists $u_*$ such that
\begin{equation}\label{eq:source-volterra-limit}
 |u_j/q_j-u_*|\leq C E2^{-bj},\qquad
 |u_*-u_0/q_0|\leq CE\,.
\end{equation}
The coefficients and sequences may be complex.
\end{lemma}
The lemma is the elementary heart of the matching: a perturbed
solution of a Volterra recursion with short memory, divided by a
nonvanishing reference solution, converges geometrically. Note
that only the reference solution is used as a divisor; the
perturbed solution may vanish.

\begin{proof}
Put $x_j=u_j/q_j$. Subtract $x_j$ times the reference recursion
from the perturbed one to obtain
\[
 x_{j+1}-x_j=\sum_{k<j}\frac{m_{j,k}q_k}{q_{j+1}}(x_k-x_j)
                                      +\frac{e_j}{q_{j+1}}.
\]
For $d_j=|x_{j+1}-x_j|$, telescope $x_k-x_j$ and interchange
the finite sums:
\[
 d_j\leq C_c\epsilon\sum_{l<j}2^{-c(j-l)}d_l+E2^{-b_0j},
 \qquad C_c=(1-2^{-c})^{-1}.
\]
On $j\leq J$, taking the norm $\max 2^{bj}d_j$ gives an operator
norm at most $C_c\epsilon\sum_{n\geq1}2^{-(c-b)n}<1/2$ for
small $\epsilon$. Hence $d_j\leq2E2^{-bj}$, uniformly in $J$.
Summing its tail proves both conclusions. Only the nonvanishing
reference sequence $q$ was divided out.
\end{proof}

\begin{theorem}[Bounded density sources and matched mass exponent]
\label{thm:sharp-source-flow}
For sufficiently large $\beta$, there are coefficients
$G_{r;a}$ such that, with $b=1/32$,
\[
 \begin{aligned}
 \Gamma_{j,r;\omega\omega}
       &=G_{r;\omega\omega}+O(\delta2^{-bj}),\\
 \Gamma_{j,r;\omega,-\omega}
       &=2^{(\eta_\mathrm{m}-\eta)j}
                \big(G_{r;\omega,-\omega}+O(\delta2^{-bj})\big),
 \end{aligned}\,,
\]
where the exponents belong to the matched reference model.
Moreover
\begin{equation}\label{eq:sharp-source-initial-match}
 G_{r;\omega\omega'}
                =K_re^{-ip^{\omega'}\cdot d_r}+O(\delta)\,.
\end{equation}
The estimates hold for all three edge types and both Fermi labels,
with constants uniform on compact positive slope sets.
In particular, the normalized smooth density sources are bounded
and converge; their anomalous exponent is zero.
\end{theorem}
\begin{proof}
Fix $r,a$ and apply \cref{lem:source-volterra-stability}
to \cref{eq:source-comparison-recursion} with
$u_j=\Gamma_{j,r;a}$ and $q_j=q_{a,j}$. The hypotheses of the
lemma are supplied by the two preceding lemmas. Indeed, by
\cref{eq:scalar-reference-recursion} and the ratio bound
$|q_{a,k}/q_{a,j+1}|\leq C2^{C\delta(j+1-k)}$ of
\cref{eq:reference-source-asymptotic}, the reference convolution
satisfies
\[
 \Big|\frac{m^R_{a;j,k}q_{a,k}}{q_{a,j+1}}\Big|
   \leq C\delta2^{-(1/2-C\delta)(j-k)}\,,
\]
which is the required bound with $\epsilon=C\delta$ and, for
large $\beta$, $c=1/4$. Also, by
\cref{eq:source-comparison-error} and that ratio bound,
\[
 |e_{j,r;a}/q_{a,j+1}|
       \leq C\delta2^{-(r_*-\eta_\beta-C\delta)j}\,.
\]
Increase $\beta$ so that the exponent is at least $r_*/2=1/16$.
We may therefore take $E=C\delta$, $b_0=1/16$ and $b=1/32$.
The lemma gives constants $c_{r;a}$ with
\[
 \Gamma_{j,r;a}/q_{a,j}=c_{r;a}+O(\delta2^{-bj}),\qquad
 c_{r;a}=\Gamma_{0,r;a}+O(\delta)\,.
\]
Multiply by \cref{eq:reference-source-asymptotic} and set
$G_{r;a}=b_ac_{r;a}$. Since $b_a=1+O(\delta)$ and the
initial coefficients are bounded, this proves both assertions,
using \cref{eq:initial-source-coefficient} for the last one.
This argument allows the physical recursion to mix all four
labels and allows its coefficients to vanish at intermediate
scales. Only the nonvanishing scalar reference solutions were
used as divisors. The comparison is at the already matched
zero-source dispersion; no source-dependent retuning occurs.
\end{proof}

The theorem closes the source-exponent matching step. Boundedness
of the density source alone does not imply Gaussianity of the
height. The next subsubsection proves the two-point comparison and
uses the exact curl identity of height gradients to determine their
homogeneous covariance, retaining the short-distance contribution.
\Cref{sec:height-gaussianity} then proves the cancellation
of the higher smeared height cumulants.

\subsubsection{Two-edge asymptotics and the height covariance}
\label{sec:height-covariance}

We now prove the two-point comparison and determine the height
covariance. The plan of this subsubsection is as follows. The
multiscale comparison of the preceding subsubsections yields, in
\cref{prop:sharp-two-edge}, the large-distance asymptotics of the
covariance of two edge indicators: a term decaying like the
inverse square of a complex linear form $L^\perp$ of the separation,
its complex conjugate, two oscillating terms carrying the momenta
$\pm2p^+$ and an anomalous exponent, and a remainder which is
summable in the plane. The increments of the diagonal height are,
up to a sign and a centering, edge indicators, by the exact
identities \cref{eq:height-gradient-anchor-shift} below; the
asymptotics therefore determine the low-frequency behavior of
the spectral density of the increments, which is computed in
\cref{lem:covariance-fourier-form}. The leading angular terms
obtained in this way have the form
$uu^{\trps}\,\overline{L_+(p)}/L_+(p)$ and its conjugate, which is
not of the form $pp^{\trps}/Q(p)$ expected for the gradient of a
\GFF. The missing information is supplied by the exact curl
identity satisfied by the increments of any function on $\Z^2$:
it forces, in \cref{lem:curl-covariance-algebra}, the two angular
terms and the constant Fourier term to combine into
$a_*pp^{\trps}/Q(p)$ for a single scalar $a_*>0$. This is why the
constant Fourier term contributed by short separations must be
kept track of throughout and cannot simply be discarded: it is
the only term at our disposal to complete the chiral leading terms
into the gradient of a logarithmically correlated field.
\Cref{thm:height-spectral-density,cor:height-covariance-limit}
then convert the spectral density into the logarithmic variance
of the height, the covariance limit of its smooth averages and
tightness in local negative Sobolev spaces. Gaussianity is a
separate matter, taken up in \cref{sec:height-gaussianity}.

Fix the diagonal lattice coordinates, and write
\[
 L_+(p)=\ell_1p_1+\ell_2p_2,\qquad
 \sA_\mathrm{D}p=(\Re L_+(p),\Im L_+(p)),\qquad
 Q(p)=|L_+(p)|^2,
\]
and define the rotated linear form
\[
 L^\perp(x):=\ell_2x_1-\ell_1x_2=L_+(-x_2,x_1)\,.
\]
Here $\ell=(\ell_1,\ell_2)\in\C^2$ is the complex gradient of the
matched dispersion of \cref{sec:physical-matching} at its plus
zero $p^+$, expressed in the diagonal height basis, so that
$L_+$ is the linearization of that dispersion at $p^+$ and
$L_-(p)=-\overline{L_+(p)}$ is the linearization at $p^-=-p^+$.
The real matrix $\sA_\mathrm{D}$ is invertible, because the
differential of the dispersion at its zero is invertible as a
real map to $\C$; it is the matrix appearing in the statement of
\cref{thm:coefficient-rg}, and $Q(p)=|\sA_\mathrm{D}p|^2$ is the
positive definite quadratic form which replaces $|p|^2$ in the
limiting \GFF covariance. The linear form $L^\perp$ is the
position-space counterpart of $L_+$: writing $z$ for the complex
coordinate of $\sA_\mathrm{D}^{-\trps}x$, one has
$L^\perp(x)=i(\det\sA_\mathrm{D})z$, as verified in the proof of
\cref{lem:covariance-fourier-form}, and $1/L^\perp(x)^2$ is, up to a
constant factor, the square of the free chiral propagator, that
is, the position-space form of the density--density bubble.
Below all constants are uniform on compact liquid parameter sets.
Let $n_r(x)$ be the indicator of the type-$r$ edge anchored at
$x$, and let $C_{rs}(x)=\cov(n_r(x),n_s(0))$ in
the infinite-volume tiling law. By the canonical and sector
identifications this is the coefficient $\cK_2$ already
constructed in \cref{sec:bulk-sources}, whose coarse decay
\cref{eq:bulk-response-decay} is now to be sharpened into an
asymptotic expansion.

\begin{proposition}[Two-edge comparison with a summable remainder]
\label{prop:sharp-two-edge}
For sufficiently large inverse temperature, there are nonzero
complex numbers $v_r$, bounded coefficients $B_{rs}^{\pm}$, and
$\nu=1+\eta-\eta_\mathrm{m}=1+O(d_\beta)>0$ such that, for
$x\in\Z^2\setminus\{0\}$,
\begin{equation}\label{eq:sharp-two-edge}
 C_{rs}(x)=\frac{v_rv_s}{L^\perp(x)^2}
       +\frac{\overline{v_rv_s}}{\overline{L^\perp(x)}^2}
       +\sum_{\omega=\pm}
          B_{rs}^{\omega}
               \frac{e^{2i\omega p^+\cdot x}}{|L^\perp(x)|^{2\nu}}
       +\cR_{rs}(x),
 \qquad
 |\cR_{rs}(x)|\leq C(1+|x|)^{-2-\kappa}\,,
\end{equation}
where one may take $\kappa=1/128$. The vector $(v_r)_r$ is a
common nonzero scalar multiple of $(G_{r;++})_r$. The scalar
and its reciprocal are uniformly bounded. The two oscillation
momenta $\pm2p^+$ are nonzero modulo $2\pi\Z^2$.
No assertion at coincident edges is implicit in the displayed
leading terms; their values at zero are set separately below.
\end{proposition}
\begin{proof}
Use the density/mass/remainder decomposition of
\cite[\S6.6, Eq.~(6.152)]{GMT20}. Its density pairs carry no
phase, its mass pairs carry the phases $e^{\pm2ip^+\cdot x}$, and
the remaining trees are compared by
\cref{lem:marked-source-comparison}.

For the first step, give reference bilinear $a$ the multiplier
$c_{r;a}=G_{r;a}/b_a$, where $G_{r;a}$ are the limiting source
coefficients of \cref{thm:sharp-source-flow} and $b_a$ the
reference normalizations of \cref{eq:reference-source-asymptotic}.
\Cref{thm:sharp-source-flow,eq:reference-source-asymptotic} give
\[
 |\Gamma_{j,r;a}-c_{r;a}q_{a,j}|
       \leq Cd_\beta2^{-bj}2^{Cd_\beta j},\qquad b=1/32\,.
\]
The ordinary history and curvature comparisons are
\cref{eq:source-reference-history,lem:running-curvature}. These provide exactly the
marked differences used in the cited decomposition; only ratios
of field normalizations occur, so the overall normalization
constant $a_Z$ of \cref{eq:source-reference-Z} cancels.

The two-source tree bound of \cite[\S6.6, Eq.~(6.153)]{GMT20} retains
$s_h^2e^{-c\sqrt{s_h|x|}}$ at the first merger scale $h$ of
the source branches. The two normalized sources cost at most
$C2^{2Cd_\beta h}$. Every difference tree has one of the marks
in \cref{lem:one-source-kernel-comparison}, or the source
coefficient difference above. A mark at scale $k\leq h$ gives
\[
 \sum_{k=0}^h2^{-(h-k)/4}2^{-bk}\leq C2^{-bh}\,,
\]
and a mark at $k\geq h$ has $2^{-bk}\leq2^{-bh}$ and a geometric
remaining branch sum. Reserving half this gain yields
\begin{equation}\label{eq:two-edge-comparison-shell}
 C2^{-(2+b/2-2Cd_\beta)h}e^{-c\sqrt{s_h|x|}}\,.
\end{equation}
Its constant need not vanish in the free case, since free curvature
and free source remainders are included.

The positive endpoint norms and reserved connecting weights are
handled by \cref{lem:marked-source-comparison}. If both sources
belong to one initial endpoint, retain its fixed-scale weight
between $0$ and $x$; \cref{lem:source-shell-summability} sums the
remaining branches and gives $Ce^{-c\sqrt{s_0|x|}}$.

Choose $d_\beta$ with $b/2-2Cd_\beta>\kappa$ and sum
\cref{eq:two-edge-comparison-shell} using
\cref{eq:dyadic-source-sum}. The resulting difference is
$O((1+|x|)^{-2-\kappa})$. Pure local reference graphs are exactly
continuum graphs by \cref{lem:all-order-reference-sampling}:
the two bilinear source vertices keep the maximum vertex degree
at four, so Poisson summation still excludes nonzero reciprocal
modes. Smooth reference-potential remainders belong to the marked
class just estimated.

The continuum asymptotics are
\cite[Eqs.~(4.17)--(4.20)]{GMT20}: the equal-chirality densities
have leading term $a_\mathrm d/L^\perp(x)^2$ and its conjugate,
with $a_\mathrm d\ne0$, and the opposite mass correlations have
leading term $a_\mathrm m/|L^\perp(x)|^{2\nu}$, where
$\nu=1+\eta-\eta_\mathrm m$ is the anomalous exponent produced
by the different running of the mass and density sources in
\cref{eq:reference-source-asymptotic}. Their remainders are
$O(|x|^{-3})$ and $O(|x|^{-2-1/2})$, respectively, at sufficiently
small coupling. To account
also for the mixed density term, that is, the covariance between a
density of chirality $+$ and a density of chirality $-$, which
vanishes in the free model but not in the interacting reference
model, choose a normalized nonnegative
Gaussian reference potential in adapted coordinates. The mixed
covariance follows by solving the density Ward identities of
\cite[Eq.~(5.11)]{GMTH17}. Explicitly, writing $D_\omega$ for the chiral
symbol, $\widehat S_{\alpha,\omega}$ for the reference density
covariances, and $c_{\rm ref}$ for their normalization constant,
these identities read
\[
\begin{split}
 D_\omega\widehat S_{\omega,\omega}
   -\tau\widehat v D_{-\omega}\widehat S_{-\omega,\omega}
     &=-c_{\rm ref}D_{-\omega},\\
 D_{-\omega}\widehat S_{-\omega,\omega}
   -\tau\widehat v D_\omega\widehat S_{\omega,\omega}&=0.
\end{split}
\]
Eliminating $D_\omega\widehat S_{\omega,\omega}$ gives, for $p\ne0$,
\[
 \widehat S_{-\omega,\omega}(p)
   =-c_{\rm ref}\frac{\tau\widehat v(p)}
                         {1-\tau^2\widehat v(p)^2}.
\]
Since $|\tau|<1$ and $0\leq\widehat v\leq1$ is Gaussian,
this extends smoothly through zero and is a Schwartz function;
its inverse Fourier transform is therefore rapidly decreasing.
Thus the mixed density term contributes only to the remainder
$\cR_{rs}$, and in fact decays faster than any power.
Independent flavor charges
exclude density--mass and equal-charge mass pairs: a density
bilinear carries no charge, whereas a mass bilinear carries
charge $\pm2$, and the reference model conserves this charge.

Multiplying these correlations by the two source multipliers
$c_{r;a}$ and $c_{s;a}$ and
restoring the phases $e^{\pm2ip^+\cdot x}$ of the mass bilinears, which were
factored out of the mass bilinears in the reference model, gives
\cref{eq:sharp-two-edge},
with $v_r=\sqrt{a_\mathrm d}\,G_{r;++}/b_{++}$. Reality gives
the conjugate density term: the covariance is real, and the
chirality $-$ density is the complex conjugate of the chirality
$+$ density. In particular the vector $(v_r)_r$ is the scalar
multiple $\sqrt{a_\mathrm d}/b_{++}$ of $(G_{r;++})_r$, as
asserted. The reference amplitudes stay nonzero
and bounded; \cref{eq:sharp-source-initial-match} and positive
liquid weights give $v_r\ne0$, since the leading term
$K_re^{-ip^+\cdot d_r}$ of $G_{r;++}$ has modulus bounded below
by the smallest edge weight and the $O(\delta)$ correction is
small. Finally $p^-=-p^+$ and the Fermi
points are distinct, so $\pm2p^+$ are nonzero modulo the reciprocal
lattice, uniformly on compact liquid parameter sets.
\end{proof}

We use the Fourier convention $\widehat a(p)=\sum_xa(x)e^{-ip\cdot x}$
on the lattice, with inverse measure $\d p/(2\pi)^2$.
A Fourier series which is not absolutely convergent is initially
understood as a periodic distribution. The chiral leading terms in
\cref{eq:sharp-two-edge} are not absolutely summable, since they
decay like $|x|^{-2}$ in the plane; the mass terms decay like
$|x|^{-2\nu}$ and need not be absolutely summable either. The next lemma
computes their Fourier series near zero and proves the regularity
needed to use the full series as a spectral density. The outcome
is that the chiral term $1/L^\perp(x)^2$ produces a bounded but
discontinuous angular multiplier at $p=0$, whereas the oscillating
mass terms are smooth near zero because their singularity sits at
the nonzero momenta $\pm2p^+$.

\begin{lemma}[The low-frequency form of the covariance]
\label{lem:covariance-fourier-form}
Set $L^\perp(x)^{-2}$ and $|L^\perp(x)|^{-2\nu}$ equal to zero at
$x=0$. The Fourier series of the first sequence has, near
$p=0$, the form
\begin{equation}\label{eq:chiral-fourier-form}
 \sum_{x\ne0}\frac{e^{-ip\cdot x}}{L^\perp(x)^2}
     =\frac{\pi}{|\det\sA_\mathrm{D}|}
                       \frac{\overline{L_+(p)}}{L_+(p)}+H(p)\,,
\end{equation}
where $H$ is smooth. It is a bounded function on the momentum
torus. For any $q\notin2\pi\Z^2$ and $\nu>0$, the Fourier
series of $e^{iq\cdot x}|L^\perp(x)|^{-2\nu}$ is integrable on
the torus and smooth in a neighborhood of zero not containing
$q$ modulo the reciprocal lattice.
The remainder in \cref{eq:sharp-two-edge} has a Fourier
transform of class $C^{0,\epsilon}$ for every
$0<\epsilon<\kappa$.
\end{lemma}
\begin{proof}
Put $y=\sA_\mathrm{D}^{-\trps}x$ and $z=y_1+iy_2$.
Then $L^\perp(x)=i(\det\sA_\mathrm{D})z$. The standard distributional
identities
$\partial_{\bar z}(1/z)=\pi\delta_0$ and
$\partial_z(1/z)=-\pv(1/z^2)$ give
$\cF[\pv(1/z^2)](q)=-\pi\bar q/q$ in our Fourier convention.
Changing variables $x=\sA_\mathrm{D}^{\trps}y$ therefore gives
\[
 \cF[\pv(1/(L^\perp)^2)](p)
       =\frac{\pi}{|\det\sA_\mathrm{D}|}
                          \frac{\overline{L_+(p)}}{L_+(p)}.
\]
This multiplier is bounded, with its only possible discontinuity
at zero.

It remains to pass from the continuum Fourier transform to the
lattice Fourier series. Choose a smooth cutoff which vanishes
near zero and equals one
outside a disk containing no other lattice point, so that the
lattice sequence $1/L^\perp(x)^2$, with the value zero at the origin,
is the restriction to $\Z^2$ of the continuum function
$a(x)$ obtained by multiplying $1/(L^\perp)^2$ by this cutoff. Multiplying
$1/(L^\perp)^2$ by this cutoff changes its continuum principal value
by a compactly supported distribution. Its Fourier transform is
smooth. The resulting function $a(x)$ is smooth everywhere;
each of its derivatives of positive order is integrable, since
the $k$th derivative of $1/(L^\perp)^2$ is $O(|x|^{-2-k})$.
Repeated integration by parts, also after multiplying by any
fixed polynomial in $x$, shows that $\widehat a$ and its
derivatives decrease faster than every power away from $p=0$.
Poisson summation, first with a Schwartz regularization and then
in distributions, gives
\[
                  \sum_{x\in\Z^2}a(x)e^{-ip\cdot x}
                      =\sum_{\ell\in2\pi\Z^2}\widehat a(p+\ell)\,.
\]
The nonzero translates on the right and all their derivatives
converge uniformly near zero, by the rapid decrease just
established, and their sum is the smooth function $H$ in
\cref{eq:chiral-fourier-form}, together with the smooth
contribution of the compactly supported correction. This proves
the formula and boundedness, since zero is its only possible
singular point on the torus and the multiplier is bounded there.

For the mass sequence use the same cutoff and decompose its
continuum extension into smooth dyadic annuli of radius $R=2^j$.
A piece has the form $R^{-2\nu}a_0(x/R)$ with $a_0$ smooth
on a fixed annulus. Its periodized Fourier transform has
$L^1$ norm at most $CR^{-2\nu}$, by scaling and the Schwartz
bounds on $\widehat a_0$. These norms sum when $\nu>0$.
Away from zero in momentum, the differentiated series converges
uniformly: integration by parts gives any power of
$(R|p+\ell|)^{-1}$, paying the powers of $R$ from the
momentum derivatives. The finite initial piece is smooth.
Frequency translation by $q$ proves the assertions
for the oscillating sequence.

Finally, absolute summability of $\cR$ follows from
\cref{eq:sharp-two-edge}. Since
$\sum_x|x|^\epsilon|\cR_{rs}(x)|<\infty$ for
$\epsilon<\kappa$, use
$|e^{-ip\cdot x}-e^{-ip'\cdot x}|
 \leq C_\epsilon|p-p'|^\epsilon|x|^\epsilon$
to prove the last assertion. Changing the value at zero in
any of the three sequences only changes a constant Fourier
term, which has been retained throughout.
\end{proof}

\begin{lemma}[Curl determines the homogeneous covariance]
\label{lem:curl-covariance-algebra}
Let $L$ be a real-invertible complex linear form on $\R^2$,
let $Q=|L|^2$, and suppose
\[
 S_0(p)=C+uu^{\trps}\frac{\overline{L(p)}}{L(p)}
              +\bar u\bar u^{\trps}
                                      \frac{L(p)}{\overline{L(p)}}
                  \qquad(p\ne0)\,,
\]
where $C$ is a real symmetric two-by-two matrix and
$u\in\C^2\setminus\{0\}$. Suppose $S_0(p)$ is positive
semidefinite and $p_2(S_0)_{1j}=p_1(S_0)_{2j}$ for $j=1,2$.
Then there is $a_*>0$ such that
\[
                           S_0(p)=a_*\frac{pp^{\trps}}{Q(p)}\,.
\]
\end{lemma}
\begin{proof}
The entries of $N(p)=Q(p)S_0(p)$ are real homogeneous quadratic
polynomials. The curl identities give
$p_2N_{11}=p_1N_{12}$ and $p_2N_{12}=p_1N_{22}$.
Since the coordinate polynomials $p_1,p_2$ have no common
nonconstant factor in $\R[p_1,p_2]$, these identities imply successively
$N_{11}=p_1\ell$, $N_{12}=p_2\ell$ with $\ell$ linear,
and then $\ell=a_*p_1$. Hence $N_{11}=a_*p_1^2$,
$N_{12}=a_*p_1p_2$, $N_{22}=a_*p_2^2$ for a real constant $a_*$.
Positivity gives $a_*\geq0$. If $a_*=0$, choose $i$ with $u_i\ne0$.
The $ii$ entry would assert that
$C_{ii}+u_i^2 e^{-2it}+\bar u_i^2 e^{2it}=0$ for every
$t$, since the argument of $L(p)$ ranges over the full circle.
The independence of these three Fourier modes gives $u_i^2=0$,
a contradiction. Thus $a_*>0$.
\end{proof}

\begin{theorem}[The diagonal height spectral density]
\label{thm:height-spectral-density}
Let $f(0)=0$ be the diagonal height and set
$g(x)=f(x)+\rho_1x_1+\rho_2x_2$ on $\Z^2$. The centered
stationary increments $\zeta_i(x)=g(x+e_i)-g(x)$ have an
integrable matrix spectral density $S(p)$. There are $a_*>0$
and $\epsilon>0$ such that
\begin{equation}\label{eq:height-gradient-spectrum}
               S(p)=a_*\frac{pp^{\trps}}{Q(p)}+O(|p|^\epsilon)
                           \qquad(p\to0)\,.
\end{equation}
For $d_i(p)=e^{ip_i}-1$ and $p\ne0$ in the momentum torus,
\begin{equation}\label{eq:height-scalar-spectrum}
 S(p)=d(p)d(p)^*s(p),\qquad s(p)\geq0,\qquad
 s(p)=\frac{a_*}{Q(p)}+O(|p|^{-2+\epsilon})\quad(p\to0)\,.
\end{equation}
Away from zero $s$ is integrable. For finitely supported real
weights $(w_x)$ with $\sum_xw_x=0$,
\begin{equation}\label{eq:height-spectral-test}
 \E\bigg|\sum_xw_xg(x)\bigg|^2
        =\int_{[-\pi,\pi]^2}
                  \bigg|\sum_xw_xe^{ip\cdot x}\bigg|^2
                                 s(p)\frac{\d p}{(2\pi)^2}\,.
\end{equation}
\end{theorem}
\begin{proof}
In the diagonal coordinates of \cref{sec:tetris}, the projections of
the three coordinate directions of $\Z^3$ are $e_1,e_2,-e_1-e_2$.
Choose the black elementary triangles
$B_x=\conv\{x,x+e_1,x+e_1+e_2\}$ and the
white triangles
$W_x=\conv\{x,x+e_2,x+e_1+e_2\}$.
This fixes the lattice basis used for $L_+$ above; the covariance
bounds are invariant under the fixed change from the initial
Kasteleyn basis. The matching edges of orientations $1,2,3$
from $B_x$ end at $W_{x-e_2},W_{x+e_1},W_x$, respectively.
Indeed their shared sides have directions $e_1,e_2,e_1+e_2$;
the corresponding lozenges are normal to the $x_1$-, $x_2$- and
$x_3$-axes, respectively.

On a lozenge normal to the $i$-th axis, $i=1,2$, the height has gradient
$-e_i$, and on the third orientation it has gradient zero.
The segment from $x$ to $x+e_1$ is a side of $B_x$, whereas
that from $x$ to $x+e_2$ is a side of $B_{x-e_1}$.
Consequently the exact identities are
\begin{equation}\label{eq:height-gradient-anchor-shift}
 \zeta_1(x)=-\big(n_{r_1}(x)-\rho_1\big),\qquad
 \zeta_2(x)=-\big(n_{r_2}(x-e_1)-\rho_2\big)\,.
\end{equation}
Thus these are bounded centered stationary fields. The anchor
shift must be retained in the Fourier covariance: its $ij$
entry is multiplied by $e^{ip\cdot(s_i-s_j)}$, with
$s_1=0,s_2=-e_1$.

Apply \cref{prop:sharp-two-edge,lem:covariance-fourier-form} to these two edge types.
All mass terms are smooth near momentum zero; the remainder
is H\"older there. Thus, for $0<\epsilon<\kappa$,
\[
 S(p)=C+uu^{\trps}\frac{\overline{L_+(p)}}{L_+(p)}
       +\bar u\bar u^{\trps}
                              \frac{L_+(p)}{\overline{L_+(p)}}
       +O(|p|^\epsilon),
 \qquad
 u_i=\sqrt{\pi/|\det\sA_\mathrm{D}|}\,v_{r_i}\,.
\]
The common minus signs cancel in the covariance, and the anchor
multipliers equal $1+O(|p|)$ near zero. Thus they only affect the
H\"older remainder in this display. Here $u\ne0$, and $C$ includes the short-distance coefficients,
the smooth mass contribution, and all nonzero Poisson images.
Reality and $S_{ij}(p)=\overline{S_{ji}(p)}$ make $C$ real
symmetric. More explicitly the leading angular terms are real,
even and symmetric, while $S(-p)=\overline{S(p)}$, so the
constant limit of their difference has these properties too.

The Fourier formulas give an $L^1$ matrix density on the entire
torus. Its Fourier coefficients are the actual covariances,
so uniqueness of Fourier coefficients identifies it with their
positive spectral measure. One can also see positivity directly:
for every vector trigonometric polynomial $P$, the integral
$\int P^*SP$ is the variance of the corresponding finite linear
combination of increments, hence is nonnegative. Approximation
then gives positive semidefiniteness almost everywhere.

We now bring in the curl identity. Since $\zeta_1,\zeta_2$ are the
increments of the single function $g$, both sides of
\[
 \zeta_1(x+e_2)-\zeta_1(x)
 =\zeta_2(x+e_1)-\zeta_2(x)
\]
equal $g(x+e_1+e_2)-g(x+e_1)-g(x+e_2)+g(x)$, so the identity
holds pointwise. Taking the covariance of both sides with
$\zeta_j(0)$ and transforming in $x$, a shift by $e_i$
multiplies the Fourier series by $e^{ip_i}$, whence
\[
                         d_2(p)S_{1j}(p)=d_1(p)S_{2j}(p)\,.
\]
For $p=tq$ with $q\ne0$ fixed, divide by $it$ and let
$t\downarrow0$: since $d_i(tq)=itq_i+O(t^2)$ and the H\"older
remainder in the display above vanishes in the limit, the
leading matrix $S_0(q)$, of the form assumed in
\cref{lem:curl-covariance-algebra} with $L=L_+$,
satisfies $q_2(S_0)_{1j}(q)=q_1(S_0)_{2j}(q)$. It is continuous
on $q\ne0$ and positive semidefinite there, being the limit of
the positive semidefinite matrices $S(tq)$, so it satisfies all
the hypotheses of \cref{lem:curl-covariance-algebra}.
That lemma proves \cref{eq:height-gradient-spectrum} with
$a_*>0$. This step is why the constant $C$ could not be omitted:
without it the two angular terms alone cannot satisfy the curl
identities, and it is the specific value of $C$, which collects
the short-distance covariances, that completes them into
$a_*pp^{\trps}/Q(p)$.

At almost every nonzero momentum the exact curl identities
make every column of $S$ proportional to $d$. Hermitian
symmetry and positivity give $S=dd^*s$ with
$s=\operatorname{tr}S/|d|^2\geq0$. Since
$|d(p)|^2=|p|^2+O(|p|^4)$, taking the trace in
\cref{eq:height-gradient-spectrum} proves
\cref{eq:height-scalar-spectrum}, reducing $\epsilon$ below
one if needed. Integrability away from zero follows from
that of $S$ and the nonvanishing of $|d|$ there.

Finally express $g(x)$ as the sum of its increments along a
fixed finite lattice path from zero to $x$. The Fourier
multiplier of that path contracted with $d(p)$ is
$e^{ip\cdot x}-1$, by telescoping, independently of the path.
The covariance formula for finite linear combinations of
increments therefore gives
\[
 \E[g(x)g(y)]=\int
 (e^{ip\cdot x}-1)(e^{-ip\cdot y}-1)s(p)\frac{\d p}{(2\pi)^2}\,.
\]
The product removes the singularity at zero. Multiply by
$w_xw_y$ and sum to prove \cref{eq:height-spectral-test}.
\end{proof}

\begin{corollary}[Logarithmic variance and height covariance]
\label{cor:height-covariance-limit}
In the unconditioned diagonal tiling law,
\begin{equation}\label{eq:logarithmic-height-variance}
 \E|g(x)|^2
     =\frac{a_*}{\pi|\det\sA_\mathrm{D}|}
           \log\big(2+|\sA_\mathrm{D}^{-\trps}x|\big)+O(1)
                      \qquad(x\in\Z^2)\,,
\end{equation}
with a uniform $O(1)$ error. For real $F,G\in C_c^\infty(\R^2)$
of integral zero, set
\[
 X_k(F)=\sum_xg(x)\int_{x/k+[0,1/k)^2}F(u)\,\d u\,.
\]
Then
\begin{equation}\label{eq:height-covariance-limit}
 \lim_{k\to\infty}\E[X_k(F)X_k(G)]
       =a_*\int_{\R^2}\frac{\widehat F(p)\overline{\widehat G(p)}}
                                  {Q(p)}\frac{\d p}{(2\pi)^2}\,.
\end{equation}
This is the covariance of
$\sqrt{a_*/|\det\sA_\mathrm{D}|}\,\GFF
                \circ\sA_\mathrm{D}^{-\trps}$,
when the standard \GFF has Fourier covariance $|p|^{-2}$.
The centered height fields are tight, modulo constants, in
$H^{-s}_\mathrm{loc}$ for every $s>1$. These assertions concern
second moments and tightness; they do not assert Gaussianity.
\end{corollary}
\begin{proof}
Use \cref{eq:height-scalar-spectrum} in the last covariance
formula of the theorem with $y=x$. The error near zero is
integrable, since $\epsilon>0$, and the part away from zero
is integrable as well. Both contribute $O(1)$ uniformly in
$x$. In the leading integral change variables
$q=\sA_\mathrm{D}p$. Replacing the transformed small
neighborhood by a fixed disk costs $O(1)$. For
$y=\sA_\mathrm{D}^{-\trps}x$, the remaining integral is
\[
 \frac{2a_*}{|\det\sA_\mathrm{D}|}
    \int_{|q|<c}\frac{1-\cos(q\cdot y)}{|q|^2}
                                             \frac{\d q}{(2\pi)^2}\,.
\]
In polar coordinates its value is
$a_*(\pi|\det\sA_\mathrm{D}|)^{-1}\log(2+|y|)+O(1)$.
For clarity, the part $|q||y|<1$ is bounded using
$1-\cos t\leq t^2/2$. In the other part the angular mean of
$\cos(r|y|\cos t)$ is $O((r|y|)^{-1/2})$: split off angular
intervals of width $(r|y|)^{-1/2}$ at its stationary points
and integrate by parts on the complement. Its radial integral
against $\d r/r$ is bounded. The constant angular term
integrates to the stated logarithm. This proves
\cref{eq:logarithmic-height-variance}.

Let $w_{k,F}(x)$ be the cell integral defining $X_k(F)$ and set
\[
                 W_{k,F}(p)=\sum_x w_{k,F}(x)e^{ip\cdot x}\,.
\]
Since $F$ has integral zero, $\sum_xw_{k,F}(x)=0$. For every fixed $M$,
\begin{equation}\label{eq:height-test-fourier-bound}
 |W_{k,F}(p)|\leq C_{F,M}
             \min\{k|p|,(1+k|p|)^{-M}\}
                          \qquad(p\in[-\pi,\pi]^2)\,.
\end{equation}
The first bound follows by subtracting one from the exponential
and using $\sum_x|x||w_{k,F}(x)|\leq C_Fk$.
For the second, finite differences of order $M$ in a coordinate
with $|p_i|\geq|p|/\sqrt2$ have summed modulus at most
$C_{F,M}k^{-M}$. Discrete summation by parts divides by
$|e^{ip_i}-1|^M\geq c_M|p_i|^M$. Combine this with
$\sum_x|w_{k,F}(x)|\leq C_F$ to obtain the displayed bound.
Also $W_{k,F}(q/k)\to\widehat F(-q)$, by Riemann sums.

In \cref{eq:height-spectral-test}, restrict first to a fixed
small momentum neighborhood and put $p=q/k$.
Then $k^{-2}s(q/k)\to a_*/Q(q)$ and is bounded by
$C|q|^{-2}$. \Cref{eq:height-test-fourier-bound}
provides an integrable envelope, both at zero and infinity.
Dominated convergence proves the variance version of
\cref{eq:height-covariance-limit}. Outside that neighborhood,
\cref{eq:height-test-fourier-bound} and integrability of $s$
make the integral tend to zero. Polarization gives the asserted
mixed covariance. For $M=\sA_\mathrm{D}^{-\trps}$, the Fourier covariance of
$\GFF\circ M$ is $|\det\sA_\mathrm{D}|/Q(p)$ by the change of
variables $q=\sA_\mathrm{D}p$. This gives the stated normalization.

We finish with tightness. Fix a nonnegative
$\psi\in C_c^\infty$ of integral one and use the pinned field
\[
 Y_k(u)=g(\floor{ ku})
                    -\int\psi(v)g(\floor{ kv})\,\d v\,.
\]
Put $V(x)=\E g(x)^2$. Stationary increments give
$\E[g(x)g(y)]=(V(x)+V(y)-V(x-y))/2$. Indeed, $g(x)-g(y)$ has
the law of $g(x-y)$, so $\E(g(x)-g(y))^2=V(x-y)$, and expanding
the square gives the formula.
In the covariance of $Y_k$, the first two terms cancel, because
$Y_k$ is invariant under adding a constant to $g$, and
only the doubly pinned kernel $-V(x-y)/2$ remains.
By \cref{eq:logarithmic-height-variance}, subtracting the
common multiple of $\log k$ from this kernel leaves, on any
fixed bounded set, an absolute bound
\[
 C\big(1+|\log(|u-v|+k^{-1})|\big)\,.
\]
The same bound with a changed constant holds after pinning:
integrating either variable against $\psi$ is uniformly
bounded. Its integral in either variable over a bounded
rectangle $U$ is uniformly bounded. The Schur estimate
therefore gives
$\E|\langle Y_k,e_j\rangle|^2\leq C_U$ for every
$L^2(U)$-normalized Dirichlet eigenfunction $e_j$.
For $1<t<s$, summing with weights $(1+\lambda_j)^{-t}$
gives a uniform bound on
$\E\|Y_k\|_{H^{-t}(U)}^2$. The embedding into
$H^{-s}(U)$ is compact: in this eigenbasis its tail operator
norm is at most $(1+\lambda_{J+1})^{-(s-t)/2}\to0$.
Markov's inequality proves tightness on $U$. Exhaustion by
rectangles proves the local assertion, and changing the
pinning changes only an additive constant.
\end{proof}

The affine interpolation has the same covariance limit: the
$L^2$ comparison is included in the proof of
\cref{thm:dimer-height-gff}, using mixing and bounded
height increments. The tested cellwise fields have bounded variances
by the corollary.

\subsubsection{Smooth sources and Gaussianity of the height limit}
\label{sec:height-gaussianity}

It remains to prove Gaussianity. Following \cite[\S7]{GMT17},
we compare higher source cumulants with the continuum reference
model, whose density functional is quadratic. We retain the
cancellation at coincident sources, since smooth height averages
pair the kernels with tests that need not vanish on collision
diagonals. Restricting first to annular Fourier support removes
local mass sources exactly and gives an exact current representation
of the height. The covariance estimate then extends the limit to
all smooth integral-zero tests.

\begin{lemma}[Density cancellation in the reference model]
\label{lem:reference-density-gaussian}
In the cutoff-removed continuum reference model, the connected
density-source coefficients of every order $m\geq3$ vanish.
In particular, their pairing with Schwartz sources whose Fourier
transforms have compact support away from zero is zero. The
assertion includes contact contributions at coincident source
positions. Fixed linear combinations and translations of the
two density sources have the same property.
\end{lemma}
\begin{proof}
We give a distributional version of the density cancellation
\cite[Prop.~12, Eq.~(6.94)]{GMT17}. A cancellation only at
pairwise distinct positions would not exclude contact terms.
The smooth, nonlocal potential in the reference model is
important for the following argument.

\begin{enumerate}[label=\textbf{Step~\arabic*.}, ref=\arabic*, wide=0pt, itemsep=1ex, listparindent=\parindent, parsep=0pt]
\item\label[step]{st:free-currents} \emph{(Free currents, including collisions.)}
After the adapted linear change of coordinates, the propagator
of chirality $\omega$ is a constant times $1/z_\omega(x)$,
where $z_\omega$ is a real-linear isomorphism from $\R^2$ to
$\C$. At a fixed source order $m\geq3$, the free connected
density coefficient is, up to an overall constant and sign,
\begin{equation}\label{eq:free-density-cycle-cancellation}
 \sum_{\sigma\text{ an }m\text{-cycle}}
        \prod_{i=1}^m\frac{1}{z_\omega(x_i)-z_\omega(x_{\sigma(i)})}
                                                        =0\,.
\end{equation}
For distinct points this is the cyclic cancellation of
\cite[\S3.3]{GMT17}. Mixed chiralities have no free connected cycle.
We still need its extension through collision
diagonals, i.e., we must show that the identity holds as an
identity between distributions on $(\R^2)^m$ and not only
pointwise off the diagonals.

Each cycle product is locally integrable in the relative
coordinates when $m\geq3$. To see this directly, anchor one
vertex and restrict the others to a fixed ball. In adapted
coordinates its absolute integral is bounded, up to a constant, by
$K^{*m}(0)$, where $K(x)=\one_{\{|x|\leq R\}}/|x|$ and $R$
is large enough to contain all differences. Young's convolution
inequality followed by H\"older's inequality gives
\[
 K^{*m}(0)\leq\|K\|_{m/(m-1)}^m<\infty,
\]
since $m/(m-1)<2$. This controls all simultaneous and partial
collisions at once.

For the common adapted smooth radial ultraviolet cutoff, its
propagator $g_{\omega,\Lambda}$ satisfies
$|g_{\omega,\Lambda}(x)|\leq C/|x|$ for $x\ne0$, uniformly
in $\Lambda$, and converges to the uncut propagator there.
For example write it as the convolution of $1/z_\omega$
with the rescaled Schwartz inverse transform of the cutoff.
Scaling and splitting at $|y-x|=|x|/2$ give the bound; the
same argument applies to a difference of ultraviolet and
infrared cutoffs. The preceding integrable envelope therefore
proves convergence of every cycle as a distribution, by
dominated convergence. Thus the sum is zero as a distribution,
also on the collision diagonals. The one-current coefficient
is zero with the symmetric cutoff. The two-current coefficient
has its cutoff-defined principal-value extension; its Fourier
multiplier $C_\omega(p)$ is a constant times
$D_{-\omega}(p)/D_\omega(p)$ and is bounded. This is the same
free bubble appearing in the continuum reference analysis of
\cite[\S4]{GMT20}.
Here $D_\omega$ denotes the adapted chiral dispersion.
Consequently the free normalized source functional is, as
a Taylor series on Schwartz tests,
\[
 \cZ_0(J)=\exp(\tfrac12\langle J,CJ\rangle),
                 \qquad C=\operatorname{diag}(C_+,C_-)\,.
\]
The pairing in this formula is bilinear, not a positive
Gaussian probability covariance.

\item\label[step]{st:smooth-density-interaction} \emph{(Adding the smooth density interaction.)}
Absorb its sign and fixed normalization into the symmetric
operator $V$, whose only nonzero entries are
$V_{+-}(x,y)=V_{-+}(y,x)=\gamma v(x-y)$, with $\gamma=O(g)$.
Initially multiply the potential by
$\chi_R(x)\chi_R(y)$, where $0\leq\chi_R\leq1$ is a smooth
compact cutoff. At every fixed perturbative order the exact
Grassmann identity for inserting this interaction is
\begin{equation}\label{eq:reference-current-reweighting}
 \cZ_V(J)=
 \frac{\displaystyle
   \exp(\tfrac12\langle\partial_J,V\partial_J\rangle)
                         \cZ_0(J)}
      {\displaystyle
   \big.\exp(\tfrac12\langle\partial_J,V\partial_J\rangle)
                         \cZ_0(J)\big|_{J=0}}\,.
\end{equation}
This identity is used coefficientwise before any infinite
sum. Ultraviolet removal commutes with its coefficient
integrals at fixed $R$ (initially truncate the sources
spatially as well): the free density moment is a finite
sum of tensor products of connected free cycles on disjoint
sets of variables, and these distributions are paired with
smooth compactly supported functions. A tensor factor of
order at least three tends to zero by the preceding argument.
The two-point factors converge in the principal-value sense.
Crucially, $v$ is smooth: it does not restrict a tensor product
to a collision diagonal. This would not justify replacing
$v$ by a delta interaction before removing the cutoff.

The surviving contractions in
\cref{eq:reference-current-reweighting} are pairings of free
currents. An interaction vertex has degree two, and an external
source has degree one. A connected graph is therefore either
a cycle with no external sources, which cancels between
numerator and denominator, or a chain with exactly two external
sources. Summing the chains gives
\[
 \log\cZ_V(J)=\tfrac12\langle J,C_VJ\rangle,
 \qquad C_V=C+CVC+CVCVC+\cdots=(1-CV)^{-1}C\,.
\]
This series converges in operator norm on $L^2(\R^2;\C^2)$
whenever $\|C\|\,|\gamma|\|v\|_1<1$. Its bound is uniform
in $R$, since multiplication by $\chi_R$ has norm at most one.
As $R\to\infty$, $V_R\to V$ strongly, so the convergent
series and each of its pairings with Schwartz sources converge.
The $L^2$ operator bound also removes the spatial truncations
of those sources. At every fixed higher source order the
vanishing distribution extends from compact tests to
Schwartz tests by the continuity of the reference kernels.
The reference construction is analytic in its small bare
coupling, with locally uniform cutoff limits of each fixed
source coefficient; this is the reference expansion of
\cite[\S6.3.2]{GMT17}, whose kernels are analytic in the running
sequences by \cite[Cor.~1]{GMT17} and hence in the small bare
coupling, used already in
\cref{lem:matched-source-reference}. Hence the
coefficientwise identity identifies the constructed reference
functional, rather than merely a formal perturbation series.
The smallness condition above is independent of the source
order. All its connected coefficients of order at least
three are zero, including contact contributions.
\end{enumerate}
Fixed source rescalings, translations, linear combinations and the
adapted coordinate change preserve this quadratic functional.
\end{proof}

For a family of real Schwartz profiles $a_{k,r}$, define
\begin{equation}\label{eq:smooth-edge-source}
 A_k(x,r)=k^{-1}a_{k,r}(x/k),\qquad
 T_k=\sum_{x,r}A_k(x,r)(n_r(x)-\rho_r)\,.
\end{equation}
A uniformly bounded Schwartz family means that every Schwartz
seminorm is bounded independently of $k$. The sum defining
$T_k$ converges absolutely, since the edge indicators are bounded.

\begin{lemma}[Integrated source blocks at large scales]
\label{lem:integrated-source-blocks}
Start two of the source expansions above at scale $J$ and
use the coordinates $u=s_Jx$, with cell volume $s_J^2$.
Suppose their ordinary histories obey the renormalized-tree
bounds, their normalized density sources are bounded, and
their source-dependent initial kernels have bounded
dimensionless positive norms, including the field-degree and
localization reserves of
\cref{lem:source-endpoints-localization}, at each fixed
source order.
Suppose also that the norms of differences of these initial
kernels, ordinary kernels and normalized histories are bounded by $\mu$ times
the corresponding fixed-order bounds, for $0<\mu\leq1$.
The propagator differences are understood in the shell and
Gram norms used in the preceding comparisons.
Let the rescaled source profiles $a_r$ satisfy
\[
                 \max_r(\|a_r\|_1+\|a_r\|_\infty)\leq M\,,
\]
where the first norm includes the cell volume. For every
fixed $m\geq1$, the difference of the scalar Taylor
coefficients extracted at relative scale $\ell\geq0$ obeys
\begin{equation}\label{eq:integrated-source-shell}
             |\Delta S_{J+\ell}^{(m)}(a)|
                         \leq C_{m,M}\mu\,2^{-\ell/4}\,.
\end{equation}
The constants are uniform in the mesh $s_J$.
The ordinary small-coupling threshold is independent of $m$.
\end{lemma}
\begin{proof}
Work at fixed order $m$, grouping the bounded local-source flow
into its running coefficient. Partition source occurrences into
$p\leq m$ initial blocks of orders $q_i\geq1$, with
$\sum_iq_i=m$. A block consists of the sources belonging to
one initial endpoint, before that endpoint is contracted with
others. A block with fermions has at least two; scalar blocks
are extracted at relative scale zero. We will keep one spatial
integral for each block, rather than one for each of its sources.

For example, consider a block with two sources. Write $\mu_J$
for counting measure on the rescaled lattice multiplied by its
cell volume $s_J^2$, and put $a(u)=\max_r|a_r(u)|$. After
summing the other relative coordinates and the field labels in
the positive kernel norm, its absolute kernel is bounded by a
nonnegative kernel $b_2$ satisfying
$\sup_{u_1}\int b_2(u_1,u_2)\,\d\mu_J(u_2)\leq B_2$.
Consequently
\[
 \iint a(u_1)a(u_2)b_2(u_1,u_2)
                    \,\d\mu_J(u_1)\,\d\mu_J(u_2)
       \leq B_2\|a\|_1\|a\|_\infty\,.
\]
In particular, two sources at the same point cause no loss in
the mesh: the delta kernel is
$s_J^{-2}\one_{\{u_1=u_2\}}$, whose integral in $u_2$ is one,
and the left-hand side is then $\int a(u)^2\,\d\mu_J(u)$.
For a block of order $q_i$ the same argument bounds $q_i-1$
profiles in supremum norm and keeps the integral of the remaining
profile. Its positive block measure has mass at most
$B_{q_i}\|a\|_1\|a\|_\infty^{q_i-1}$, with the finite source
labels absorbed in $B_{q_i}$. This argument uses the integrated
kernel norms only, and applies equally to nonlocal initial blocks.

Prune the tree to these blocks and apply
\cref{eq:pruned-source-bound} with initial scale $J$. The gain
$d(\fn,q)$ of \cref{eq:source-gain} is nondecreasing in $q$, so
counting each block as one source only weakens its dimensional gain.
The borderline cases $(\fn,q)=(2,1)$ and $(2,2)$ both have
gain one; additional sources increase the gain. Thus merging
the sources within an initial block does not consume any of
the decay needed to sum the tree above it.
The connecting weights are at most one and can be integrated
against the block measures. For $p\geq2$ the remaining bound is
\[
 C_{m,M}2^{-3\ell/4}\prod_v2^{-3(d_v-1)k_v/4},
\]
where $d_v\geq2$ are the branching degrees and $k_v$ their relative
scales. Each scale sum is geometric and there are at most $p-1$
branching vertices. For $p=1$, the one-source scalar bound of
\cref{lem:source-shell-summability} gives the same stem decay.
Cell volume was included in the block measures, so the bounds
are uniform in the mesh.

Telescope the two expansions using
\cref{lem:marked-source-comparison}. Each difference costs $\mu$;
its endpoint choices and normalization-ratio crossings are paid
by the degree and chain reserves. Reserving part of the stem and
branching exponents leaves $2^{-\ell/4}$, since
$(1+\ell)^2 2^{-3\ell/4}\leq C2^{-\ell/4}$.
This proves \cref{eq:integrated-source-shell}. The finite source
partitions change $C_{m,M}$ but not the ordinary-coupling threshold.
\end{proof}

\begin{proposition}[Vanishing higher cumulants for smooth edge sources]
\label{prop:smooth-source-cumulants}
Suppose the profiles in \cref{eq:smooth-edge-source} form a
uniformly bounded Schwartz family and their Fourier supports
lie in a fixed compact annulus
$\{q:\epsilon_0\leq|q|\leq R_0\}$, where
$0<\epsilon_0<R_0<\infty$. For sufficiently large $\beta$,
for every fixed $m\geq3$,
\begin{equation}\label{eq:smooth-source-cumulants}
                         |\operatorname{cum}_m(T_k)|
                                      \leq C_m k^{-b},\qquad b=1/128\,.
\end{equation}
The constants may depend on $m$, the annulus and the profile
bounds, but the smallness requirement on $d_\beta$ does not
depend on $m$. The conclusion holds for joint cumulants of
finitely many such families as well.
\end{proposition}
\begin{proof}
Work at fixed source order $m$, retaining Taylor factorials;
the final factor $m!$ is absorbed into $C_m$. Since
$\|A_k\|_1<\infty$, the coarse cumulant estimate of
\cref{thm:bulk-source-functional} makes all source sums
absolutely convergent and identifies them with the cumulants of
$T_k$ by finite-support approximation. The estimates below use
the full band-limited sources.

\begin{enumerate}[label=\textbf{Step~\arabic*.}, ref=\arabic*, wide=0pt, itemsep=1ex, listparindent=\parindent, parsep=0pt]
\item\label[step]{st:remove-mass-sources} \emph{(Removing local mass sources.)}
Integrate finitely many initial scales so that the sum of two
remaining fermion momenta has norm less than
$\dist(2p^+,2\pi\Z^2)/4$. All endpoint and matching bounds
persist, with changed constants. Recall from
\cref{eq:local-edge-source} that the Fermi momenta are carried by
the explicit phases $e^{i(p^\omega-p^{\omega'})\cdot x}$, so that
the fields $\psi_{b,\omega},\psi_{w,\omega'}$ themselves have
small momenta, and that for a mass bilinear, $\omega'=-\omega$,
the carried momentum is $p^\omega-p^{-\omega}=2\omega p^+$.
Poisson summation gives
\[
 \widehat A_k(p,r)=k\sum_{\ell\in2\pi\Z^2}
                         \widehat a_{k,r}(k(p+\ell))\,,
\]
whose torus support lies in $|p|\leq R_0/k$. For large $k$,
momentum conservation therefore gives the polynomial identity
\[
 \sum_x A_k(x,r)e^{2i\omega p^+\cdot x}
             \psi_{b,\omega}(x)\psi_{w,-\omega}(x)=0\,.
\]
The momentum $2\omega p^+$ of the phase and the two remaining fermion momenta stay at distance
at least $3\dist(2p^+,2\pi\Z^2)/4$ from the reciprocal lattice,
while $\widehat A_k$ is supported in $|p|\leq R_0/k$.
Thus local mass-source endpoints vanish before contraction on
every remaining scale. Nonlinear and nonlocal endpoints remain
and are estimated instead.

\item\label[step]{st:compare-kernels} \emph{(Comparing effective kernels.)}
The reference density sources are
\[
 J_{k,\omega}=\sum_r c_{r;\omega\omega}A_k(\cdot,r),
 \qquad c_{r;a}=G_{r;a}/b_a\,,
\]
as in \cref{prop:sharp-two-edge}.
The density coefficients differ by $O(d_\beta2^{-j/32})$ by
\cref{thm:sharp-source-flow}; ordinary histories,
propagators and initial remainders are compared by
\cref{lem:one-source-kernel-comparison}.

Write $\|\cdot\|_{1,j}$ for the positive spatial kernel norm
with one coordinate fixed, including label sums, Taylor
factorials, field-degree reserves and, for fermionic kernels,
the localization moments of
\cref{lem:source-endpoints-localization}. The dimension of
a kernel with $\fn$ fermions and $m$ sources is $2-\fn/2-m$.
After separating relevant and marginal terms, the differences
obey
\begin{equation}\label{eq:smooth-source-kernel-comparison}
 s_j^{\fn/2+m-2}\|\Delta W_{j;\fn,m}\|_{1,j}
                      \leq C_m2^{-bj},\qquad \fn\geq2,\ m\geq1\,,
\end{equation}
including the separately compared density coefficient at
$(\fn,m)=(2,1)$. Extracted scalar kernels satisfy
\begin{equation}\label{eq:smooth-source-scalar-comparison}
 \|\Delta W_{j;0,m}\|_{1,j}
                      \leq C_m s_j^{2-m}2^{-bj},\qquad m\geq3\,.
\end{equation}
Here and below $b=1/128$. These are kernels of the source
functional after the identically zero local mass terms in
\cref{st:remove-mass-sources} have been removed. The same
momentum-support test removes any regenerated local mass
term at a later scale. Thus every retained local source
coefficient is a bounded density coefficient; the estimates
do not assert boundedness of the mass-source flow for
arbitrary, unsmeared sources.

Apply \cref{lem:marked-source-comparison} with bounded density
coefficients, hence $\eta=0$, and the root dimension
$2-\fn/2-m$. The marks have rates at least $1/4$ for initial
irrelevant endpoints, $1-C\delta^2$ for curvature, $1/8$ for
ordinary histories and $1/32$ for density coefficients. Consequently
\[
 2^{-j/4}+\sum_{l=0}^j2^{-(j-l)/4}2^{-l/32}
                         \leq C2^{-j/32}.
\]
The common lemma has already paid normalization-ratio crossings
and marked choices from the unused path reserve. Weakening to
$b=1/128$ gives
\cref{eq:smooth-source-kernel-comparison,eq:smooth-source-scalar-comparison}.
Free multi-source endpoints need not be small: those with
$\fn\geq2,m\geq2$ have dimension at most $-1$, so they are
initial irrelevant marks of bounded norm. Initial scalar terms
are extracted only at the initial scale. Repeated sources are
included by the positive Taylor norm.

Pure reference graphs agree with continuum graphs by
\cref{lem:all-order-reference-sampling}; after decreasing the
fixed cutoff, the source momentum $O(k^{-1})$ preserves the
no-aliasing margin. Smooth-potential remainders and the integrated
reference ultraviolet kernels are initial marks. This completes
the comparison, with a smallness threshold independent of $m$.

\item\label[step]{st:observation-scale} \emph{(Splitting at the observation scale.)}
Choose $J$ with $1\leq ks_J<2$, up to a fixed multiplicative
constant. Since $\|A_k\|_1\leq Ck$ and
$\|A_k\|_\infty\leq C/k$, testing a scalar kernel against
one $\ell^1$ source and $m-1$ supremum norms gives
\[
 \sum_{j<J}|\Delta S_j^{(m)}(A_k)|
 \leq C_m k^{2-m}\sum_{j<J}s_j^{2-m}2^{-bj}
 \leq C_m k^{-b}\,.
\]
The sum is geometric because $m-2-b>0$ and includes all source
collisions. The reference scales above the common ultraviolet
cutoff contribute $O(C_mk^{2-m})$, by scalar dimension $2-m<0$
and the bounded normalized reference sources.

For $j\geq J$, rescale space by $s_J$ and sources by $s_J^{-1}$.
The profiles have bounded supremum and integrated norms:
\[
 s_J^{-1}\|A_k\|_\infty\leq C,\qquad
 s_J^2\sum_x|s_J^{-1}A_k(x)|\leq C\,.
\]
\Cref{eq:smooth-source-kernel-comparison} bounds the
initial fermionic source-kernel differences by $C_m2^{-bJ}$.
The bulk comparison gives the same bound for ordinary kernels
and propagators; \cref{thm:sharp-source-flow} and
\crefrange{eq:source-reference-history}{eq:source-reference-Z}
give it for later normalized histories. These are precisely the
hypotheses of \cref{lem:integrated-source-blocks}, with
$\mu=C2^{-bJ}$ (finitely many small $J$ are absorbed in $C_m$).
Using cell volume $s_J^2$, it gives
$|\Delta S_{J+\ell}^{(m)}(A_k)|\leq C_m2^{-bJ}2^{-\ell/4}$,
and hence
\[
 |\Delta S_{\geq J}^{(m)}(A_k)|\leq C_m2^{-bJ}\leq C_mk^{-b}\,.
\]
\end{enumerate}
The full coefficient thus differs from the cutoff-removed
reference coefficient by $O(C_mk^{-b})$. The latter vanishes
by \cref{lem:reference-density-gaussian}, including contact
terms. Centering changes only the first cumulant. This proves
\cref{eq:smooth-source-cumulants}; linear combinations and
polarization give the joint assertion.
\end{proof}

\begin{lemma}[An exact smooth-current representation]
\label{lem:band-limited-height-current}
Let $F$ be real Schwartz with smooth Fourier transform supported
in a compact annulus disjoint from zero. Define $X_k(F)$ by
the cell integrals in \cref{cor:height-covariance-limit}.
For all sufficiently large $k$, it has the form $T_k$ in
\cref{eq:smooth-edge-source}, with a uniformly bounded
Schwartz family of profiles whose Fourier supports lie in
that same annulus. In particular its higher cumulants satisfy
\cref{eq:smooth-source-cumulants}.
\end{lemma}
\begin{proof}
Write $w_k(x)=\int_{x/k+[0,1/k)^2}F$ and use the lattice
Fourier convention of the preceding subsubsection. If
$B(p)=\int_{[0,1)^2}e^{ip\cdot v}\,\d v$, Poisson
summation gives, on the momentum torus and for large $k$,
\[
                         \widehat w_k(p)=B(p)\widehat F(kp)\,.
\]
There are no nonzero reciprocal translates because the
Fourier support is fixed. In particular $\widehat w_k(0)=0$,
also following directly from $\int F=0$.
For $d_i(p)=e^{ip_i}-1$, put
\[
 \widehat j_{k,i}(p)=\frac{d_i(p)}{|d(p)|^2}\widehat w_k(p)\,,
\]
with value zero near zero and outside the indicated support.
This is a smooth periodic function. Its inverse transform is
real, and
\begin{equation}\label{eq:height-current-divergence}
                 w_k(x)=\sum_{i=1}^2
                             \big(j_{k,i}(x-e_i)-j_{k,i}(x)\big)\,.
\end{equation}
Indeed the Fourier multiplier of the right side is
$\sum_i\bar d_i d_i/|d|^2=1$ on the support.
After $p=q/k$, Fourier inversion shows
$j_{k,i}(x)=k^{-1}a_{k,i}(x/k)$, where
\[
 \widehat a_{k,i}(q)=
 \frac{k^{-1}d_i(q/k)}{|d(q/k)|^2}B(q/k)\widehat F(q)\,.
\]
On the fixed annulus the multiplier and all its derivatives
are uniformly bounded and converge to $iq_i/|q|^2$.
Thus the profiles are uniformly Schwartz and retain that
Fourier support.

The height has deterministically bounded increments, hence
$|g(x)|\leq C|x|$. All sums below are absolutely convergent.
Summation by parts in \cref{eq:height-current-divergence}
gives
\[
                       X_k(F)=\sum_{x,i}j_{k,i}(x)\zeta_i(x)\,.
\]
Using the exact anchor shift in
\cref{eq:height-gradient-anchor-shift}, set
$A_k(x,r_1)=-j_{k,1}(x)$,
$A_k(x,r_2)=-j_{k,2}(x+e_1)$ and set the third source to zero.
The fixed shift preserves the uniform Schwartz bounds and
Fourier support. This proves the representation and the
last assertion.
\end{proof}

\begin{theorem}[The diagonal dimer height has a \GFF limit]
\label{thm:dimer-height-gff}
Fix a positive slope and take $\beta$ sufficiently large,
uniformly on compact sets of positive slopes. The
centered diagonal height field in the infinite-volume tiling
law of \cref{thm:fixed-sector-limit} converges,
modulo constants, to
\begin{equation}\label{eq:dimer-height-gff}
            \sqrt{a_*/|\det\sA_\mathrm{D}|}\,
                              \GFF\circ\sA_\mathrm{D}^{-\trps}\,.
\end{equation}
Convergence holds jointly on smooth compactly supported
integral-zero tests and in $H^{-s}_\mathrm{loc}$ modulo constants
for every $s>1$. The covariance and normalization are those
of \cref{cor:height-covariance-limit}.
\end{theorem}
\begin{proof}
First use the annular Fourier class of the preceding lemma.
Its first cumulant is zero. Its second converges to
\cref{eq:height-covariance-limit}, whose proof applies also
to Schwartz tests: the summed finite-difference and moment
bounds in \cref{eq:height-test-fourier-bound} follow from
Schwartz seminorms in exactly the same way. All cumulants of
fixed order at least three tend to zero by
\cref{prop:smooth-source-cumulants,lem:band-limited-height-current}. The finite
moment--cumulant formula, which expresses the $n$th moment as a
polynomial in the cumulants of order at most $n$, therefore
proves convergence of every
moment to its Gaussian value. This supplies the moment bounds
as well; no bound on pointwise higher height moments has been
assumed. The Gaussian distribution is moment-determinate by
Carleman's criterion:
its even moments are $(2m-1)!!\,V^m$, whose inverse
$2m$th roots are of order $m^{-1/2}$ and have divergent sum.
The method of moments thus
gives convergence in law. Linear combinations remain in the
same test class (enlarge the annulus if necessary), so the
Cram\'er--Wold argument proves joint convergence of finitely
many such smoothed heights.

To pass to an arbitrary real $F\in C_c^\infty$ with
$\int F=0$, choose real even smooth frequency cutoffs $\chi_n$
supported on compact annuli, bounded between zero and one,
and tending to one at every nonzero momentum. Let
$\widehat F_n=\chi_n\widehat F$. Then $F_n$ belongs to the
annular Schwartz class, and
\[
       \int\frac{|\widehat F-\widehat F_n|^2}{Q(p)}\,\d p
                                                    \longrightarrow0\,.
\]
Near zero this follows from $\widehat F(p)=O(|p|)$;
at infinity it follows from rapid decrease. The extended
covariance limit, applied to the fixed Schwartz test $F-F_n$,
gives
\[
 \lim_{n\to\infty}\limsup_{k\to\infty}
                       \E|X_k(F)-X_k(F_n)|^2=0\,.
\]
The inequality $|e^{itx}-e^{ity}|\leq|t||x-y|$ transfers
the characteristic-function limit from $F_n$ to $F$.
Apply the same approximation to finitely many tests to get
joint convergence with covariance
\cref{eq:height-covariance-limit}.

Tightness in the stated local spaces is already proved in
\cref{cor:height-covariance-limit}. Pinning by a fixed
smooth probability density turns every test into an
integral-zero test. The just-proved finite-dimensional laws
determine the unique possible subsequential distributional
law, which is \cref{eq:dimer-height-gff}. Tightness and
uniqueness prove convergence in each of the stated spaces.

The affine interpolation has the same limit. The difference of
the two centered interpolations is bounded, since tiling-height
increments are bounded. Its integral against a fixed smooth test
is, up to an $O(k^{-1})$ error, a finite sum of averages
$k^{-2}\sum_x F(x/k)\eta_x$, where each $\eta_x$ is a bounded,
centered local linear combination of height increments. This follows
by integrating the fixed affine basis functions over each mesh cell;
replacing $F((x+u)/k)$ by $F(x/k)$ costs $O(k^{-1})$ after summation.
Mixing makes $\cov(\eta_x,\eta_y)$ tend to zero with
$|x-y|$. For a fixed $R$, the pairs with $|x-y|\leq R$ contribute
$O(R^2/k^2)$ to the variance of this average, and the remaining pairs
contribute at most $C_F\sup_{|z|>R}|\cov(\eta_0,\eta_z)|$.
Let $k\to\infty$ and then $R\to\infty$ to get convergence to zero
in $L^2$. On each bounded domain the interpolation difference lies
in a fixed $L^2$ ball, relatively compact in $H^{-t}$ for $t>0$.
The vanishing smooth pairings identify every subsequential limit as
zero. The same applies after pinning by a smooth probability density,
so the difference vanishes in probability in the asserted local
negative Sobolev spaces modulo constants.
\end{proof}

\subsection{Completion of the abstract theorem}\label{sec:rg-proof-completion}

\begin{proof}[Proof of \cref{thm:coefficient-rg}]
We assemble the preceding results, recording the hypotheses used.
\begin{enumerate}[label=\textbf{Step~\arabic*.}, ref=\arabic*, wide=0pt, itemsep=1ex, listparindent=\parindent, parsep=0pt]
\item\label[step]{st:completion-first-integration} \emph{(First integration and the bulk flow.)}
The absolute envelopes and source bounds in
\cref{it:rg-envelopes,it:rg-sources} give the bulk and torus
first integrations, \cref{thm:bulk-first-step,thm:torus-first-step}.
Their small kernels satisfy the hypotheses of
\cref{thm:conditional-renormalized-trees,thm:critical-flow-fixed-dispersion}.
\Cref{thm:matched-bulk-flow} matches this flow to the original
measure and supplies its dressed dispersion matrix $\sA_\mathrm{D}$.
All estimates use absolute envelopes, so signed coefficients are allowed.

\item\label[step]{st:completion-local-limit} \emph{(The prescribed-sector local limit.)}
Bulk consistency, positivity and the exact source algebra in
\cref{it:rg-locality,it:rg-positive,it:rg-sources} yield
\cref{thm:canonical-local-limit}. Independence from auxiliary
weights permits their tuning in \cref{prop:fixed-zero-tuning}.
The sector bound \cref{prop:interacting-central-sector} and the
equal-density uniqueness lemma then give
\cref{thm:fixed-sector-limit}, which is conclusion (i), including
mixing, the cumulant bounds and independence of the winding sequence.
No unique or continuous choice of matching is required.

\item\label[step]{st:completion-height-gff} \emph{(The height limit.)}
\Cref{thm:dimer-height-gff} is conclusion (ii), with covariance
and positive normalization identified in
\cref{cor:height-covariance-limit}. It includes joint convergence
on integral-zero tests and convergence in local $H^{-s}$ modulo
constants for every $s>1$.

\item\label[step]{st:completion-witness-tails} \emph{(Witness tails.)}
Under the additional hypothesis \cref{it:rg-marked},
\cref{thm:interacting-witness-tail} gives conclusion (iii), for
both canonical and prescribed-sector laws, including noncontractible witnesses.
\end{enumerate}
All estimates are uniform on the stated compact density sets and
bounded torus shapes. Conclusions (i)--(ii) also give the informal
\cref{thm:rg-typical}.
\end{proof}

\section{From the tiling to the original surface}\label{sec:transfer}

\begin{remark*}[Reader's guide]
\normalfont
This section proves the four model limits in \cref{thm:gff-limits}
from the dimer theorem. The order of limits is to construct the
infinite-volume surface law first and then rescale its height.
The main steps are the following.
\begin{enumerate}[(i)]
\item \emph{Separate local reconstruction from height transfer.}
\Cref{lem:joint-limit} gives a unique joint local limit from
annealed local approximations to the conditional surface law.
Once this coupling is stationary and ergodic,
\cref{lem:ergodic-correction,cor:gff-transfer} show that an
$L^2$ displacement between the surface and its reference height
disappears upon rescaling. \Cref{lem:deletion,prop:pi_decay}
supply the deletion and switching estimates used to localize the
conditional reconstruction.
\item \emph{Carry out the reconstruction for \HMSOS.}
\Cref{lem:rooting-volume} identifies the relation between diagonal
and physical rooting. \Cref{lem:displacement-support} locates a
height displacement in one defect, while
\cref{prop:hmsos-conditional-locality} makes the reconstruction
local using witness tails. With
\cref{thm:interacting-witness-tail}, these yield the joint diagonal
limit and displacement moments of
\cref{cor:hmsos-joint-exponential}.
\Cref{lem:positive-density-proper-lift,lem:physical-reconstruction}
then give local decoding, ergodicity and mean slopes in physical
coordinates, proving \cref{thm:hmsos-physical-local-limit}.
\item \emph{Change coordinates for the fluctuating field.}
\Cref{lem:diagonal-height-second-moment} controls the error in
replacing the random coordinate map by its limiting linear map.
\Cref{lem:height-palm-conditioning} handles the conditioning used
for physical rooting, and \cref{prop:reference-coordinate-change}
performs the change of coordinates. \Cref{cor:physical-height-transfer}
combines this with displacement transfer, proving
\cref{thm:hmsos-gff-complete} and applying to the remaining models.
\item \emph{Extend reconstruction to longer histories.}
For \SOS, \cref{prop:sos-conditional-locality} controls the complete
history and \cref{lem:sos-column-carrier} bounds both physical
displacements by one component. For the perturbed Ising Interface
model, \cref{lem:interface-perturbed-positive-marks} gives the
positive representation used in
\cref{prop:interface-perturbed-conditional-locality}.
The resulting physical limits and \GFF transfers are
\cref{thm:sos-gff-complete,thm:interface-gff-complete}.
\item \emph{Reduce full Ising to the perturbed interface.}
\Cref{prop:ising-one-sheet-contours} integrates out closed
contours, and \cref{prop:ising-admissible-perturbation} verifies
the resulting contact interaction. At $M_n=\ceil{e^{n^2}}$,
\cref{prop:ising-joint-cutoff} identifies the same local interface
law. \Cref{prop:ising-contour-kernel} reconstructs the bulk spins,
completing the full Ising limit in \cref{thm:ising-gff-complete}.
\end{enumerate}
\end{remark*}

We first specify the meaning of the target limits, since two
different notions of convergence are involved. Local convergence
concerns height gradients: for every finite window, their joint law
in that window must converge. For interfaces with overhangs, we also
prove convergence of the full interface rooted modulo a common
vertical translation, which is a stronger assertion. For irrational slopes take admissible integer
winding numbers whose slopes converge to the prescribed slope; exact
equality $\Theta_i=n\theta_i$ for every $n$ is not required, and indeed
cannot be required when $\theta_i$ is irrational. The scaling limit
is a statement about the rescaled height viewed as a random
distribution. Recall that the full-plane \GFF is a random
distribution modulo constants, defined only through its pairings with
smooth compactly supported test functions of integral zero. For
$f\in C_c^\infty(\R^2)$ with $\int f=0$, set
\begin{equation}\label{eq:test-field}
 \langle\fh_k,f\rangle
 =\sum_{x\in\Z^2}(\fh(x)-\E\fh(x))
                 \int_{x/k+[0,1/k)^2}f(u)\,\d u\,.
\end{equation}
Thus the lattice height is viewed as the piecewise constant
function $u\mapsto\fh(\floor{ku})$ on $\R^2$, and the rescaled field
is this function tested against $f$. The integral convention cancels
additive constants exactly: adding a constant to $\fh$ adds that
constant times $\sum_x\int_{x/k+[0,1/k)^2}f=\int f=0$ to the
right-hand side, which is why a limit modulo constants is the
natural object. Pointwise values $\fh(kx)$ do not converge to
pointwise values of a two-dimensional \GFF, whose pointwise variance
is infinite; only the pairings above converge. An amplitude $\upsigma$
must also be specified in the statement of the limit, since
precomposition of the \GFF by a linear map alone does not change the
coefficient of the logarithmic singularity of its covariance and so
cannot in general absorb the stiffness of the model.

We write $\phi$ for the vertical height of the monotone reference
surface $\varphi$, and $\fh$ for the height observable of the model
under study. For \SOS and \HMSOS the height in \cref{eq:test-field} is
the usual column height of $\cS$ or $\cH$. For the Ising Interface case of
\cref{thm:gff-limits}, the height is the vertical
column magnetization, namely the height of $\Tthree(\cI)$.
For full Ising, $\cI$ is the unique end-separating interface with
probability tending to one. Its two-phase filling defines the same
observable; \cref{def:ising-erasure} fixes a convention
on the exceptional configurations as well.
\Cref{prop:ising-joint-cutoff} compares its law at
$M_n=\ceil{e^{n^2}}$ with the perturbed interface law, jointly with
bulk spins in fixed rooted windows. The height is defined from the
erased surface rather than by summing raw spins in an infinite column.

The coordinate-change argument of \cite[\S7]{LaslierLubetzky24}
uses quantitative height moments as well as a limiting plane.
\Cref{sec:quantitative-coordinate-change} obtains the needed
moment directly from \cref{cor:height-covariance-limit}.

\subsection{Main results of this section}
\label{sec:transfer-main-results}

Throughout, fix $\theta_1,\theta_2>0$, let
$\rho=(\rho_1,\rho_2,\rho_3)=(\theta_1,\theta_2,1)/(1+\theta_1+\theta_2)$
be the corresponding face densities, as in \cref{thm:coefficient-rg},
and put
\begin{equation}\label{eq:physical-coordinate-matrix}
 A=I-e\alpha^{\trps},\qquad \alpha=(\rho_1,\rho_2)^{\trps},\qquad
 e=(1,1)^{\trps},\qquad \det A=\rho_3>0\,,
\end{equation}
the matrix of the change from diagonal to physical coordinates
constructed in \cref{sec:quantitative-coordinate-change}. The
following four theorems are proved in the course of this section,
each at the point where its ingredients are complete. All four are
stated for fixed positive slopes $\theta_1,\theta_2$, all sufficiently
large $\beta$, and along every admissible winding sequence approaching
these slopes. In each, the \emph{scaling limit} of a height $\fh$ in
the limiting law means that its centered rescaled fields
\cref{eq:test-field} converge, jointly on smooth compactly supported
integral-zero tests and in local $H^{-s}$ modulo constants for every
$s>2$, to
\begin{equation}\label{eq:physical-gff-form}
 \frac{\upsigma}{\rho_3}\,\GFF\circ(LA^{-1})\,,
\end{equation}
where $\upsigma$ and $L$ are the amplitude and linear map that
\cref{thm:coefficient-rg}(ii) assigns to the coefficient family there.

\begin{theorem}[The \HMSOS local limit and \GFF scaling limit]
\label{thm:hmsos-gff-complete}
The physically rooted \HMSOS heights $\cH-\cH(0)$ converge locally to
a unique full-plane law, which is stationary and ergodic under
physical translations followed by rerooting, with mean gradients
$\E\nabla_i\cH=-\theta_i$, $i=1,2$. The height has the scaling limit
\cref{eq:physical-gff-form} for the \HMSOS family of
\cref{thm:hmsos-rg-input}, with $\upsigma=\sqrt{a_*/|\det\sA_\mathrm{D}|}$
and $L=\sA_\mathrm{D}^{-\trps}$.

This implies the statement of \cref{thm:gff-limits} on the \HMSOS model.
\end{theorem}

\begin{theorem}[The \SOS physical local limit and \GFF]
\label{thm:sos-gff-complete}
The physical \SOS histories $(\cS,\cH,\phi)$, rooted by $\phi(0)=0$,
converge locally to a unique law, which is stationary and ergodic
under physical translations followed by the common rerooting,
satisfies $\Ttwo\cS=\cH$ and $\Tone\cH=\phi$ almost surely, and has
mean gradients $\E\nabla_i\cS=\E\nabla_i\cH=\E\nabla_i\phi=-\theta_i$,
$i=1,2$. Its marginal $\cS-\cS(0)$ is the unique local limit of the
physically rooted \SOS heights, and $\cS$ has the scaling limit
\cref{eq:physical-gff-form} for the \SOS family of
\cref{thm:sos-rg-input}, with the parameters
$\upsigma_{\mathrm d}^{\textsc{sos}},L_{\mathrm d}^{\textsc{sos}}$
of \cref{eq:sos-diagonal-gff}.

This implies the statement of \cref{thm:gff-limits} on the \SOS model.
\end{theorem}

\begin{theorem}[Perturbed interface local limit and \GFF]
\label{thm:interface-gff-complete}
Let the slopes lie in a compact subset of $(0,\infty)^2$ and let $U$
be an admissible contact perturbation as in
\cref{def:interface-perturbation}. The physical histories
$(\cI,\cS,\cH,\phi)$ rooted by $\phi(0)=0$ converge locally to a
unique stationary ergodic law, with $\Tthree\cI=\cS$, $\Ttwo\cS=\cH$
and $\Tone\cH=\phi$ almost surely; both spin phases of $\cI$ have
paths to their respective infinite vertical ends, and the mean
gradients of $\cS,\cH,\phi$ are $-\theta_i$, $i=1,2$. The height
$\cS=\Tthree\cI$ has the scaling limit \cref{eq:physical-gff-form}
for the perturbed interface family of \cref{thm:interface-rg-input},
with the parameters $\upsigma_{\mathrm d}^{U},L_{\mathrm d}^{U}$ of
\cref{eq:interface-diagonal-gff}. The same local limit is obtained
from the finite-height restrictions with $M_n/n^2\to\infty$, by
\cref{eq:interface-height-cutoff-tv}.

The case $U=0$ implies the statement of \cref{thm:gff-limits} on the Ising
Interface model.
\end{theorem}

\begin{theorem}[The full Ising local limit and \GFF]
\label{thm:ising-gff-complete}
Consider the 3\Dim Ising model with the surface erasure of
\cref{def:ising-erasure} and the joint cutoff $M_n=\ceil{e^{n^2}}$.
The rooted erased histories $(\cI,\cS,\cH,\phi)$ converge locally to
a unique stationary ergodic law, with $\Tthree\cI=\cS$,
$\Ttwo\cS=\cH$ and $\Tone\cH=\phi$ almost surely and mean gradients
$-\theta_i$ for $\cS,\cH,\phi$. The erased interface and the original
bulk spins, with the same vertical rooting, converge jointly, and the
limiting conditional spin law given $\cI$ is the kernel $\cK(\cI)$ of
\cref{prop:ising-contour-kernel}, with orientation fixed by the
boundary conditions. The height function $\fh$ of
\cref{def:height-function-in-models}, which equals $\cS=\Tthree\cI$
by \cref{def:tetris-charge}, has the scaling limit
\cref{eq:physical-gff-form} with the parameters of
\cref{thm:interface-gff-complete} for the contact perturbation of
\cref{prop:ising-admissible-perturbation}.

This implies the statement of \cref{thm:gff-limits} on the full 3\Dim Ising model.
\end{theorem}

\subsection{Local limits, field transfer and interaction costs}
\label{sec:transfer-criteria}

\subsubsection{Local limits and transfer criteria}

The reconstruction of the original surface from its tiling is carried
out in two stages, and this subsection isolates the criteria
for each stage. The first criterion, \cref{lem:joint-limit}, produces
a joint local limit of the surface and its tiling from a local limit
of the tiling alone, provided the conditional law of the surface given
the tiling is \emph{annealed-local}: the window of the surface around
a point is approximately determined, in expectation over the tiling,
by the tiling in a large ball around that point. The approximation is
only required on average over the tiling law, which is what the
witness estimates of \cref{sec:interacting-witnesses} supply; a
uniform approximation in the tiling is neither available nor needed.
The second criterion, \cref{lem:ergodic-correction} and
\cref{cor:gff-transfer}, transfers the \GFF limit from the reference
height to the surface height once the joint limit is stationary and
ergodic and the displacement between the two heights has a finite
second moment. No uniform pointwise bound on this displacement is
required.

\begin{lemma}[Annealed local approximation implies a joint limit]\label{lem:joint-limit}
Suppose $\varphi_N$ converges locally to $\varphi$, and let $\cX_N$
denote the reconstructed surface data (including the height field
$\fh_N$ when appropriate). All finite coordinate spaces are assumed
countable and discrete. For each finite
window $A$ suppose there is a kernel $K_{A,r}(\varphi|_{B_r})$ on the
surface data in $A$, depending on finitely many base coordinates, such
that
\[
 \limsup_{N\to\infty}\E\big[
 \|\cL(\cX_N|_A\given\varphi_N)
           -K_{A,r}(\varphi_N|_{B_r})\|_\mathrm{TV}\big]
       \leq\delta_A(r),\qquad \delta_A(r)\longrightarrow0\,.
\]
Then $(\cX_N,\varphi_N)$ has a unique local limit. If the finite-volume
coupling is stationary, so is the limiting coupling. If, in addition
to this stationarity, the base limit is mixing and, for bounded surface cylinder
observables $F,G$,
\[
 \lim_{|x|\to\infty}\limsup_{N\to\infty}
 \E\big|\cov(F(\cX_N),G(\tau_x\cX_N)\given\varphi_N)\big|=0\,,
\]
then the joint limit is mixing and hence ergodic.
\end{lemma}
\begin{proof}
Fix a surface window $A$ and a base window $C$, and let $\mu_N$
be their joint law. For fixed $r$, replace the conditional surface
law by $K_{A,r}$, obtaining $\widehat\mu_{N,r}$. Then
\[
 \limsup_N\|\mu_N-\widehat\mu_{N,r}\|_{\mathrm{TV}}
                                      \leq\delta_A(r).
\]
The base law on $B_r\cup C$ converges in total variation, since
these coordinates take values in a countable discrete space.
Contraction under the fixed kernel therefore makes
$(\widehat\mu_{N,r})_N$ Cauchy. The triangle inequality gives
\[
 \limsup_{N,M\to\infty}\|\mu_N-\mu_M\|_{\mathrm{TV}}
                                      \leq2\delta_A(r).
\]
Let $r\to\infty$. Completeness in total variation gives a unique
probability limit on each window, including the full countable
defect data. Marginalization and translation commute with this
limit, giving consistency and stationarity. Kolmogorov's extension
theorem supplies the joint law.

For mixing it suffices, by approximation, to use products
$a(\varphi)F(\cX)$ of bounded base and surface cylinder functions.
The law of total covariance splits their covariance into the
expected conditional covariance and the covariance of conditional
means. The former tends to zero by hypothesis, after extracting
the bounded base factors. For the latter, replace
$\E[F(\cX_N)\mid\varphi_N]$ by
$K_{A,r}(\varphi_N|_{B_r})(F)$, and likewise for the translated
observable $G$. The error is at most
$C(\delta_A(r)+\delta_{A'}(r))+o_N(1)$, uniformly in the translation
by stationarity. For fixed $r$, the replacements are base cylinder
functions; their limiting covariance tends to zero by base mixing.
Send $N\to\infty$, then the separation to infinity, and finally
$r\to\infty$. This proves mixing of the joint law, hence ergodicity.
\end{proof}

The next lemma is the reason that the height transfer requires no
uniform bound on the displacement between the surface and its
reference tiling. The rescaled field of \cref{eq:test-field} is an
average of the lattice height over the order of $k^2$ sites in the
support of $f(\cdot/k)$, each with weight of order $k^{-2}$. A
displacement which is a stationary ergodic field with a finite second
moment is averaged in the same way, and the ergodic theorem makes
such an average converge to the constant $\E R_0$, which the
centering removes. The rescaled reference height, in contrast, has
fluctuations of order one in this normalization (its variance grows
logarithmically in the window), so it survives the averaging while
the displacement does not.

\begin{lemma}[An ergodic correction vanishes upon rescaling]\label{lem:ergodic-correction}
Let $(R_x)_{x\in\Z^2}$ be a stationary ergodic field with
$\E R_0^2<\infty$. Define
$R^{(k)}(u)=R_{\floor{ ku}}-\E R_0$.
Then for every compactly supported smooth $f$,
\[
                  \langle R^{(k)},f\rangle\longrightarrow0
                           \quad\mbox{in }L^2\,.
\]
On a bounded rectangle $U$ the same convergence holds in
$L^2(\Omega;H^{-s}(U))$ for each $s>0$.
\end{lemma}
\begin{proof}
Work on the stationary field's probability space, with commuting
unitary shifts $U_i$ on $L^2$, and put $r=R_0-\E R_0$.
The closed span of the differences $g-U_i g$, $g\in L^2$,
$i=1,2$, is the mean-zero subspace: its orthogonal complement
consists exactly of functions fixed by both shifts, hence of
constants by ergodicity. We may therefore approximate $r$ in
$L^2$ by finite sums of such differences.

For a smooth compactly supported test $f$, set
$w_k(x)=\int_{x/k+[0,1/k)^2}f(u)\,\d u$ and
$T_k(f)g=\sum_xw_k(x)U_xg$. Stationarity and Minkowski give
$\|T_k(f)g\|_2\leq\|f\|_1\|g\|_2$, uniformly in $k$.
For a single difference, summation by parts gives
\[
 T_k(f)(g-U_i g)
   =\sum_x\big(w_k(x)-w_k(x-e_i)\big)U_xg,
 \qquad
 \|T_k(f)(g-U_i g)\|_2
       \leq k^{-1}\|\partial_i f\|_1\|g\|_2.
\]
Indeed the sum of the absolute weight differences is at most
$\|f-f(\cdot-e_i/k)\|_1\leq k^{-1}\|\partial_i f\|_1$.
The approximation of $r$, together with the uniform operator bound,
now proves $T_k(f)r=\langle R^{(k)},f\rangle\to0$ in $L^2$.
The same bound extends this convergence to every $L^1$ test
by density of smooth compactly supported functions.

For the Sobolev assertion use the orthonormal Dirichlet eigenbasis
$(e_j,\lambda_j)$ on $U$. Stationarity gives
$\E\|R^{(k)}\|_{L^2(U)}^2=|U|\var(R_0)$, so Parseval yields
\[
 \E\|R^{(k)}\|_{H^{-s}(U)}^2
 \leq\sum_{j\leq J}(1+\lambda_j)^{-s}
                      \E|\langle R^{(k)},e_j\rangle|^2
       +(1+\lambda_{J+1})^{-s}|U|\var(R_0).
\]
Each fixed coefficient tends to zero by the first part, applied to
$e_j$ extended by zero, an $L^1$ function. Let $k\to\infty$ and
then $J\to\infty$; since $\lambda_J\to\infty$, this proves the
claim for every $s>0$.
\end{proof}

\begin{corollary}[Transfer of the \GFF limit]\label{cor:gff-transfer}
Suppose a stationary ergodic coupling of the original height $\fh$ and its
approximating tiling $\varphi$, in the same coordinate convention and with a
common choice of vertical origin, has
$R_x=\fh(x)-\phi(x)\in L^2$. If the centered reference height field converges
to $\upsigma\,\GFF\circ L$, then the centered original height has the
same limit. This holds jointly for test functions, and as random
distributions whenever the tiling convergence is in that topology
(or in a common local negative Sobolev topology).
\end{corollary}
\begin{proof}
Centering the identity $\fh=\phi+R$ at every site gives
$\fh(x)-\E\fh(x)=(\phi(x)-\E\phi(x))+(R_x-\E R_0)$, where we used
that $\E R_x=\E R_0$ by stationarity. Testing against $f$ as in
\cref{eq:test-field} therefore gives, for every $k$,
\[
 \langle\fh_k,f\rangle
     =\langle\phi_k,f\rangle+\langle R^{(k)},f\rangle\,,
\]
with $R^{(k)}$ as in \cref{lem:ergodic-correction}. The first term
converges in law to $\langle\upsigma\,\GFF\circ L,f\rangle$ by
hypothesis and the second tends to zero in $L^2$, hence in
probability, by \cref{lem:ergodic-correction}; Slutsky's theorem
gives the convergence of the sum, jointly over finitely many tests,
and the Sobolev version of the lemma gives the distributional
statement. No above/below symmetry of the displacement and no
assertion $\E R=0$ is needed: the centering removes whatever
constant mean the displacement has. One may root $\fh$ afterwards,
for instance by subtracting $\fh(0)$; zero-mean tests remove the
resulting common additive constant. In particular, separate rooting
of both fields is not needed to obtain a stationary displacement
field, and this is why the coupling is required to use a common
vertical origin.
\end{proof}

\begin{lemma}[A sufficient annealed moment condition]\label{lem:annealed-moment}
Suppose for every $r\geq1$ there is a base event $G_r$ such that
$\P(G_r^c)\leq C e^{-cr}$ and
$\P(|R_0|>r\given\varphi)\one_{G_r}\leq C e^{-cr}$.
Then $\P(|R_0|>r)\leq2Ce^{-cr}$ and $R_0$ has all finite moments.
\end{lemma}
\begin{proof}
Split the expectation of the conditional tail over $G_r$ and its
complement. For $p>0$, integrate
$\E|R_0|^p=p\int_0^\infty t^{p-1}\P(|R_0|>t)\,\d t$.
\end{proof}

Two remarks on the hypotheses are in order. First, the \texttt{good}
events $G_r$ must have small complements under the
\emph{interacting} tiling law, which is the law of the reference
tiling in the coupling; a bound under the free dimer law, such as
\cref{prop:fat_pattern}, does not supply this hypothesis, and this
is why the witness tails are carried through the \RG in
\cref{sec:interacting-witnesses}. Second, a small conditional
variance of $\fh$ given $\varphi$ would not suffice by itself: the
displacement $\fh-\phi$ also contains the random conditional mean
$\E[\fh\given\varphi]-\phi$, which is a function of the tiling and
could in principle carry fluctuations of the same order as the
tiling height. The ergodic-correction criterion handles both
contributions at once, because it only uses that the total
displacement is a stationary ergodic field with an annealed second
moment, which is what \cref{lem:joint-limit,lem:annealed-moment}
provide.

\subsubsection{Interaction costs and decay in the conditional law}\label{sec:interaction-cost}

In finite volume, write $\pi_{\beta,h}$ for the conditional law on the complete
one-step fiber over the reference surface $h$: each admissible
preimage $\cI$ has probability $Z_{\beta,h}^{-1}e^{-\beta E(\cI;h)}$.
This law includes the internal configuration of each defect; its
support marginal is the hard-core gas of \cref{lem:polymer-factorization}.
Reconstruction requires two kinds of estimates on this law: that a
prescribed defect is unlikely (a Peierls-type bound), and that the
law is insensitive to what happens far away (a decay of
correlations, used to localize the reconstruction of the surface
from its tiling). The natural way to
quantify the second property is through an ``interaction cost''
between two sets $A,B$ which is large when the two sets are far
apart. Such a cost must take into account the entropy of chains of
defects connecting $A$ to $B$, not merely the weight of one such
chain, and it must be formulated for normalized probabilities. An unnormalized partition function with a defect prescribed at $A$
still sums over arbitrary distant defects and can grow with the
volume. Dividing by the full partition function removes this growth
and gives an occurrence probability. To control the influence of
$B$ on $A$, we instead need the normalized probability of a chain
joining the two sets. The definition below includes its entropy.
It is a cost rather than a
metric, and in particular we never use a triangle inequality for it:
the two-copy argument of \cref{lem:chain-connection} below bounds
the probability of a connection directly by a sum over chains.

Let $\cD_\Lambda$ be the nonempty connected supports contained in a
finite region $\Lambda$. Defects are incompatible when their supports
intersect or are $*$-adjacent. The following statements hold for the exact gas
in \cref{lem:polymer-factorization}, and for any hard-core gas with
nonnegative activities; we write $\pi$ for its law, so that the
probability of a family $\{D_i\}$ of compatible supports is
proportional to $\prod_i z_{\beta,h}(D_i)$.

\begin{lemma}[Deletion estimate]\label{lem:deletion}
For compatible prescribed supports $D_1,\ldots,D_k$,
\[
 \pi(D_1,\ldots,D_k\text{ present})
 \leq\prod_{j=1}^k z_{\beta,h}(D_j)\,.
\]
The bound also holds after fixing any compatible collection of other
polymers.
\end{lemma}
\begin{proof}
Deleting $D_1,\ldots,D_k$ from a compatible family containing them
is injective, and the weight of the original family is
$\prod_j z_{\beta,h}(D_j)$ times the weight of the remaining family.
Summing over the images and bounding their total weight by the full
partition function proves the bound after normalization.
Conditional on the presence of other compatible supports, deleting
them and their incompatible neighbors leaves a hard-core gas with
the same activities, so the same argument applies. If they exclude
a required support, its conditional probability is zero.

For later use, if $\xi$ is a specified $D$-defect with weight
$z_{\beta,h}(D,\xi)$, so that
$z_{\beta,h}(D)=\sum_\xi z_{\beta,h}(D,\xi)$, the same deletion
argument gives
\[
 \pi((D,\xi)\text{ present})
 \leq z_{\beta,h}(D,\xi)\,.
\]
This is the version used before summing over the $D$-defects.
\end{proof}

\begin{definition}[Chain interaction cost]\label{def:chain-cost}
For sets $A,B\subset\Lambda$, a colored chain is a finite sequence of
distinct pairs $(D_j,c_j)$, $c_j\in\{1,2\}$, such that $D_1$ meets $A$,
$D_k$ meets $B$, consecutive supports are incompatible, colors
alternate, and supports of the same color are pairwise compatible.
Let
\[
 Q_{\beta,h}^{\Lambda}(A,B)
   =\sum_{\text{colored chains from }A\text{ to }B}\prod_j z_{\beta,h}(D_j),
 \qquad
 \fd_{\beta,h}^{\Lambda}(A,B)
   =-\log\min\{1,Q_{\beta,h}^{\Lambda}(A,B)\}\,.
\]
The value is $+\infty$ when the sum vanishes. In a finite region the
sum is finite; the infinite-volume sum is the increasing limit over
finite regions. This is a cost, not a metric. It includes both the
activity and the number of possible chains, is determined by $h$ on
the two-square neighborhood of $\Lambda$, and is well defined even in bad
regions, where it is simply small and the estimates below remain
true but uninformative.
\end{definition}

The reason for the two colors is the following. The chains which
will matter are not chains of defects in one sample of the gas
(whose supports are pairwise compatible, so that no two of them are
incompatible), but chains alternating between two independent samples
of the gas, as in the next lemma. Within one sample the supports must
be compatible, while consecutive supports in the chain come from
different samples and must be incompatible, which is exactly the
constraint in the definition. The product of the $z_{\beta,h}(D_j)$ is the
bound of \cref{lem:deletion} for the simultaneous presence of the
supports of each color.

\begin{lemma}[Two-copy connection bound]\label{lem:chain-connection}
For two independent copies of the gas (possibly restricted to different
subregions of $\Lambda$),
\[
 (\pi_1\otimes\pi_2)(A\longleftrightarrow B
          \text{ in the union of defect supports})
       \leq e^{-\fd_{\beta,h}^{\Lambda}(A,B)}\,.
\]
\end{lemma}
\begin{proof}
Suppose that $A$ and $B$ are connected in the union of the supports
present in the two copies. Form the graph whose vertices are these
supports, with an edge between two of them when they are
incompatible (intersecting or $*$-adjacent). A connection from $A$
to $B$ in the union of supports yields a path in this graph from a
support meeting $A$ to a support meeting $B$, and we may take it to
be a simple path. No edge of that path joins polymers
of the same copy, since the supports present in one copy are
pairwise compatible. The path therefore alternates between the two
copies, and coloring its supports by their copy makes it
one of the colored chains in the definition. For a fixed colored
chain, the supports of color $1$ are compatible and present in the
first copy, those of color $2$ are compatible and present in the
second copy; by independence between the
copies and \cref{lem:deletion} applied to each copy, the probability
of this event is at most the product
of $z_{\beta,h}(D_j)$ along the chain. Take the union bound over the colored
chains and then the minimum with one.
\end{proof}

\begin{proposition}[Switching, with the boundary specified]\label{prop:pi_decay}
Let $\Lambda'\subset\Lambda$ and use the same activities in both
regions, with empty exterior polymer boundary conditions. Let $F$ be
a bounded observable supported on $A\subset\Lambda'$. Define $B$ to
be the set of vertices of $\Lambda'$ $*$-adjacent to its complement.
Then
\[
 |\pi_\Lambda(F)-\pi_{\Lambda'}(F)|
 \leq2\|F\|_\infty
       e^{-\fd_{\beta,h}^{\Lambda}(A,B)}\,.
\]
The same switching identity applies to an integrable observable
attached to the defect component at a point, with the bound
\begin{equation}\label{eq:weighted-switch}
 |\pi_\Lambda(F)-\pi_{\Lambda'}(F)|
 \leq\E[(|F_1|+|F_2|)\one_{\{A\leftrightarrow B\}}]\,.
\end{equation}
More generally, the same estimate holds for two allowed polymer
families with common activities. For an observable determined by
polymer types in $\mathcal A$, replace the connection event by the
event that a union component meets both $\mathcal A$ and the
symmetric difference of the allowed families.
For two independent copies of one gas and bounded observables $F,G$
determined by polymer types in $\mathcal A,\mathcal B$, respectively,
let $\mathcal E_{\mathcal A,\mathcal B}$ be the event that a union
component meets both type sets. Then
\begin{equation}\label{eq:switching-covariance}
 |\cov_\pi(F,G)|\leq
 2\|F\|_\infty\|G\|_\infty
       (\pi\otimes\pi)(\mathcal E_{\mathcal A,\mathcal B})\,.
\end{equation}
\end{proposition}
\begin{proof}
We prove the assertion for two allowed families with common activities;
the spatial version is the special case in which the second family
consists of the polymers contained in $\Lambda'$, so that the
\emph{forbidden} types are the polymers leaving $\Lambda'$. Sample the
two gases independently, keeping the full description of each polymer and not only its support, and
color the polymers of the first copy red and those of the second blue.
Form the graph whose vertices are the polymer occurrences of both
copies and whose edges join incompatible pairs; it is bipartite, since
each copy is a compatible family. Let $\mathcal E$ be the event that
some connected component meeting $\mathcal A$ contains a polymer whose
type is allowed in only one of the two families. Off $\mathcal E$,
exchange the two colors on every component meeting $\mathcal A$; on
$\mathcal E$, do nothing. The result is again a pair of allowed
families: every swapped type is allowed in both gases, and a swapped
polymer is compatible with every unswapped one, since an
incompatibility between them would be an edge placing them in the same
component. The pair weight is the product of the activities of all
polymers present, which the swap permutes, and the uncolored components
are unchanged, so the swap is a measure-preserving involution of the
product gas which maps $\mathcal E$ to itself. Off $\mathcal E$ every
polymer of a type in $\mathcal A$ lies in a swapped component, so the
swap exchanges $F_1$ and $F_2$. Hence
\[
 \pi_1(F)-\pi_2(F)
 =\E[(F_1-F_2)\one_{\mathcal E}],\qquad
 |\pi_1(F)-\pi_2(F)|
 \leq\E[(|F_1|+|F_2|)\one_{\mathcal E}]\,,
\]
the contribution of the complement of $\mathcal E$ vanishing by the
involution. For integrable $F$, truncate first and use dominated
convergence. In the spatial setting, a component meeting $\mathcal A$
and containing a forbidden polymer is a union component joining $A$ to
a polymer leaving $\Lambda'$, and such a component crosses the inner
$*$-boundary $B$ of $\Lambda'$. Hence $\mathcal E\subset\{A\leftrightarrow B\}$,
and \cref{lem:chain-connection} gives the stated bounds, with the factor
$2$ from $|F_1|+|F_2|\leq2\|F\|_\infty$.

For \cref{eq:switching-covariance}, take two independent copies of the
same gas, so that every type is allowed in both, and write
$\cov_\pi(F,G)=\frac12\E[(F_1-F_2)(G_1-G_2)]$. Exchange the colors on
every component meeting $\mathcal A$, unless some component meets
both $\mathcal A$ and $\mathcal B$, i.e., unless
$\mathcal E_{\mathcal A,\mathcal B}$ occurs, in which case do nothing.
As above, this is a measure-preserving involution, and it preserves the
uncolored components, hence the event $\mathcal E_{\mathcal A,\mathcal B}$.
Off $\mathcal E_{\mathcal A,\mathcal B}$, every polymer of a type in
$\mathcal A$ lies in a swapped component and no polymer of a type in
$\mathcal B$ does, so the swap exchanges $F_1$ with $F_2$ and fixes $G_1$
and $G_2$; it therefore reverses the sign of $(F_1-F_2)(G_1-G_2)$, and the
contribution of the complement of $\mathcal E_{\mathcal A,\mathcal B}$
to the expectation vanishes. On $\mathcal E_{\mathcal A,\mathcal B}$ we
bound the integrand by
$|F_1-F_2|\,|G_1-G_2|\leq4\|F\|_\infty\|G\|_\infty$, which gives
$|\cov_\pi(F,G)|\leq2\|F\|_\infty\|G\|_\infty
(\pi\otimes\pi)(\mathcal E_{\mathcal A,\mathcal B})$, i.e.,
\cref{eq:switching-covariance}.
\end{proof}

\subsection{One step: HMSOS reconstruction}

\subsubsection{Physical rooting and finite-volume compactness}
\label{sec:hmsos-rooting}

Continue to use the integers
$n,\Theta_1,\Theta_2$ of \cref{sec:periodic-hmsos}. Let
$\fH_\Gamma$ be the set of integer heights on $\Z^2$ such that
\begin{equation}\label{eq:physical-hmsos-periods}
 \cH(x+ne_1)=\cH(x)-\Theta_1,\qquad
 \cH(x+ne_2)=\cH(x)-\Theta_2,\qquad \cH(x+e_2)\leq \cH(x)\,.
\end{equation}
We use the occupied cubes $x_3<\cH(x_1,x_2)$, so that every such
height is an \HMSOS configuration. Conversely these are all the
periodic \HMSOS heights with the prescribed windings. Write
$Q_n=\{0,\ldots,n-1\}^2$ and
$\nabla_i\cH(z)=\cH(z+e_i)-\cH(z)$. Its excess energy is
\[
 E_\Gamma(\cH)=\sum_{z\in Q_n}
       (|\nabla_1\cH(z)|+|\nabla_2\cH(z)|)-n(\Theta_1+\Theta_2)
          =2\sum_{z\in Q_n}(\nabla_1\cH(z))_+\,.
\]
By \cref{eq:physical-hmsos-periods}, the signed increments telescope:
$\sum_{z\in Q_n}\nabla_i\cH(z)=-n\Theta_i$ for $i=1,2$.
Since the second gradient is
nonpositive, $\sum_z|\nabla_2\cH(z)|=n\Theta_2$, whereas
$|\nabla_1\cH|=2(\nabla_1\cH)_+-\nabla_1\cH$ gives
$\sum_z|\nabla_1\cH(z)|=2\sum_z(\nabla_1\cH(z))_++n\Theta_1$. Adding the
two sums and subtracting $n(\Theta_1+\Theta_2)$ yields the displayed identity.
Thus the excess energy counts twice the vertical plaquettes in the
``wrong'' direction along $x_1$, those not present in any monotone
surface with the same windings. The horizontal plaquettes
contribute the additional constant $n^2$ to the surface area, and
the minimal number $n(\Theta_1+\Theta_2)$ of vertical plaquettes is attained
exactly by the monotone surfaces.

The physical law is the probability measure on
$\fH_\Gamma^0=\{\cH\in\fH_\Gamma:\cH(0)=0\}$ with
weights $e^{-\beta E_\Gamma(\cH)}$. Thus it counts each height
modulo a vertical translation once. This specifies the rooting
in \cref{def:HMSOS}; finiteness is proved below. The point of this
subsection is to relate this physical rooting, by a
vertical translation, to the rooting that the \Tetris map
produces, which is modulo the diagonal translation.
For $v\in\Z^3$ let
\[
 (\tau_v\cH)(x)=\cH(x-(v_1,v_2))+v_3,
 \qquad d=(1,1,1)\,.
\]
Both the energy and the map $\Tone$ commute with these translations.
The quotient by $\tau_d$, in contrast to the physical quotient
by $\tau_{e_3}$, is the one naturally parametrized by diagonal
tilings and their fibers.

\begin{lemma}[Volume charge and the two translation quotients]
\label{lem:rooting-volume}
Set $V(\cH)=\sum_{z\in Q_n}\cH(z)$ and $N_0=n(n+\Theta_1+\Theta_2)$. Then
\begin{equation}\label{eq:volume-translation}
 V(\tau_v\cH)-V(\cH)=n\Theta_1 v_1+n\Theta_2 v_2+n^2v_3\,.
\end{equation}
If $\phi=\Tone\cH$ is written as a vertical height, then
$V(\phi)=V(\cH)$. There is a bijection
\[
 \fH_\Gamma/\langle\tau_d\rangle
 \longleftrightarrow
 \coprod_{M\in\cM_\Gamma}\{\cH:\Tone\cH=\phi_M\}\,,
\]
where $\phi_M$ is the unique lift of $M$ for which
$0\leq V(\phi_M)<N_0$. The bijection preserves excess energy.
In particular, no division by a configuration-dependent translation
stabilizer occurs in this parametrization.
\end{lemma}
\begin{proof}
Shifting the first coordinate in the sum defining $V$ replaces
one boundary column by its translate, increasing the sum by $n\Theta_1$;
the other horizontal shift increases it by $n\Theta_2$. A vertical unit
shift increases it by $n^2$. Iteration, also for negative shifts,
proves \cref{eq:volume-translation}.

The difference between the $\{0,1\}$-spin configurations of $\cH$ and of $\phi$
has finite support on $\Z^3/\Gamma$: both heights have the same
periods and are bounded relative to the prescribed affine plane.
Summing this signed difference by vertical fibers gives
$V(\cH)-V(\phi)$. Summing by the finitely many sorting-fiber orbits
gives zero, since $\Tone$ preserves the charge of each such fiber.
Both sums are of finitely many nonzero terms, so they may be
reordered. This proves volume preservation.

A $\Gamma$-compatible tiling lifts to a monotone surface uniquely
up to translation by an integer multiple of $d$: relative lifts
are fixed by adjacent tiles, and the kernel of diagonal projection
on $\Z^3$ is $\Z d$. \Cref{eq:volume-translation} gives
$V(\tau_d\phi)=V(\phi)+N_0$, so precisely one lift satisfies the
specified inequality. The image of any class of $\cH$ under $\Tone$
selects its tiling $M$, and translating the class along $d$ selects
precisely one preimage of $\phi_M$. Conversely such a preimage
defines a class. This proves the bijection and, by translation
invariance, preservation of its weight. A nonzero power of
$\tau_d$ cannot fix a height, as it changes $V$.
\end{proof}

\begin{proposition}[Exact physical \HMSOS law]\label{prop:physical-hmsos-law}
For all $\beta>0$, the physical partition function
is finite, with
\begin{equation}\label{eq:physical-partition-ratio}
 \cZ^\mathrm{phys}_{\beta,\Gamma}
 :=\sum_{\cH\in\fH_\Gamma^0}e^{-\beta E_\Gamma(\cH)}
   =\frac{n^2}{N_0}
                   \sum_{M\in\cM_\Gamma}Z_{\beta,\phi_M}\,.
\end{equation}
The following sampling procedure gives this physical law exactly:
sample $M$ with law $\nu_{\beta,\Gamma}$ of
\cref{prop:exact-sector}, sample its complete fiber
with law $\pi_{\beta,\phi_M}$, choose an independent uniform
$Z\in Q_n$, and output
\[
                    \cH^{[Z]}(x)=\cH(x+Z)-\cH(Z)\,.
\]
The result is independent of the choice of lift $\phi_M$.
The physical gradient law is stationary under all base translations.
Every increasing finite height cutoff exhausting
$\fH_\Gamma^0$ converges to it in total variation.
These identities require no rationality condition beyond the
integer windings $\Theta_1,\Theta_2$. In particular they apply along sequences
approaching arbitrary positive real slopes.
\end{proposition}
\begin{proof}
We first prove finiteness, then count both translation quotients
inside a common strip of volumes to obtain
\cref{eq:physical-partition-ratio}, and finally extend the same
count to observables to identify the law.

First, on the set $E_\Gamma(\cH)=k$ the sum of absolute gradients
is $n(\Theta_1+\Theta_2)+k$. An edge path inside $Q_n$ from $0$ to $z$ therefore
gives $|\cH(z)|\leq n(\Theta_1+\Theta_2)+k$ for every physically rooted height.
There are at most $(2n(\Theta_1+\Theta_2)+2k+1)^{n^2-1}$ assignments of these
values. Summing this polynomial bound against $e^{-\beta k}$
proves finiteness for every $\beta>0$.

Let $L$ be a positive common multiple of $N_0$ and $n^2$, and
consider all heights satisfying $0\leq V(\cH)<L$. Every orbit of
$\tau_d$ has exactly $L/N_0$ representatives in this strip;
every orbit of $\tau_{e_3}$ has exactly $L/n^2$. These facts follow
from \cref{eq:volume-translation} and do not require horizontal
translation orbits to have a common size. Summing the nonnegative
weights over the strip and using \cref{lem:rooting-volume}
gives
\[
 \frac{L}{N_0}\sum_M Z_{\beta,\phi_M}
     =\frac{L}{n^2}\sum_{\cH\in\fH_\Gamma^0}
                                             e^{-\beta E_\Gamma(\cH)}\,.
\]
The right side is finite by the preceding bound. Thus both sums
are finite and the identity proves
\cref{eq:physical-partition-ratio}.

For a bounded function $f$ of a physically rooted height set
\[
 \overline f(\cH)=\frac1{n^2}\sum_{z\in Q_n}f(\cH^{[z]})\,.
\]
The winding identities show that $\cH^{[z]}$ depends on $z$ modulo
$n\Z^2$. Consequently $\overline f$ is invariant under every
$\tau_v$, including both translations used in the strip count.
That count applies with the extra factor $\overline f$, and
division by the partition functions gives equality of its
expectations in the two quotient measures. On the physical
quotient each map $\cH\mapsto \cH^{[z]}$ is an energy-preserving
bijection, with inverse the shift by $-z$. Hence its Gibbs law
is stationary and has expectation of $f$ equal to that of
$\overline f$. On the diagonal quotient the expectation of
$\overline f$ is exactly the stated sampling procedure. Changing
the lift translates by $kd$, which just permutes the average
over $z$. This proves all the law and stationarity assertions.
Finally, the cutoff partition functions increase to the finite
unrestricted sum. The total variation distance of their normalized
restrictions from the full law is the probability of the omitted
states, which tends to zero.
\end{proof}

\begin{proof}[Proof of \cref{thm:hmsos-rg-input}]
The exact image law is \cref{prop:periodic-fiber}: the weight of
$M\in\cM_\Gamma$ is its complete periodic fiber partition function
$Z_{\beta,M}$. The mixed expansion of this fiber, with its planar
and periodic bounds, is \cref{prop:mixed-envelopes}, and
one may take $\tau_\beta=c\beta^{2/5}$, $B_\beta=\ceil{\sqrt\beta}$
and $R_\beta$ a constant multiple of $\beta^{1/100}$, enlarging
$\epsilon_\beta\to0$ as needed. We check the conditions of
\cref{def:rg-admissible} one by one.

\Cref{it:rg-positive} asks for translation-covariant coefficients,
independent of the auxiliary edge weights and sources, whose
hard-core sum is positive on every matching. These are the canonical
coefficients \cref{eq:canonical-coefficients}:
\cref{prop:canonical-completion} gives positivity on all sectors and
translation covariance, and the construction uses no auxiliary edge
weight, boundary character or source, since it only involves the
\HMSOS fibers. Since the coefficients are nonnegative, the absolute
summation of internal labels required there is the summation
already performed in \cref{eq:mixed-coefficient}.

\Cref{it:rg-envelopes} asks for the \texttt{thin}/\texttt{fat}
envelopes and rooted counts; here the envelopes can be taken to be
the coefficient parts themselves. \Cref{prop:mixed-envelopes} and its
periodic completion in \cref{prop:canonical-completion} verify it
with $B_\beta=\ceil{\sqrt\beta}$, including the torus entropy
correction $C/n$.

\Cref{it:rg-locality} asks that small torus patterns carry fixed
planar coefficients; this is the first case of
\cref{eq:canonical-coefficients}, justified by
\cref{lem:periodic-locality}, together with the last assertion of
\cref{prop:canonical-completion} that all exceptions have at least
$c_*n$ lozenges.

\Cref{it:rg-sources} is the Grassmann representation with exact
sources and sector algebra; \cref{prop:hmsos-fermions,prop:exact-sector}
give the precise source factors, gauge covariance and the polynomial
sector extraction \cref{eq:exact-sector-ratio}, including at
singular twists.

\Cref{it:rg-marked} is verified by \cref{lem:marked-union-envelopes}
for the six-class family $\scW_n$. The coefficient estimates used here were proved in
\cref{sec:model-grassmann}, independently of the \RG conclusions. Finally \cref{prop:physical-hmsos-law}
identifies the physical \HMSOS ensemble with the rooted sampling of
$\nu^{\textsc{hmsos}}_{\beta,\Gamma}$, with its rooting multiplicity
and its height-cutoff limit.
\end{proof}

There is also a local description of the root change that will
be useful for the height transfer. The sampling procedure of
\cref{prop:physical-hmsos-law} re-roots the pair at a uniformly
chosen base square, which is a global operation. The following
corollary expresses the same physical law instead as a conditioning
of the diagonal law on a cylinder event at the origin, in the
manner of a Palm distribution: conditioning on the presence of a
horizontal reference face at a fixed location. Let $\Q_{\beta,\Gamma}$
denote the joint diagonal law of $(\cH,\phi)$ modulo $\tau_d$,
obtained by sampling the tiling and then its fiber as above.
Choose as $o$ the black triangle in the projection of the horizontal
face over $0$ at height zero. Let
$A_3(o)$ be the event that its matching lozenge has normal $x_3$,
i.e., that the reference surface has a horizontal face projecting
to $o$. On this event choose the common lift in which that lozenge
is the horizontal face over the base square $0$ at height zero.
Thus $\phi(0)=0$, while $\cH$ uses the same vertical origin, so
that the displacement $\cH-\phi$ is well defined.

\begin{corollary}[The horizontal-face root]
\label{cor:physical-palm}
The law of this lifted pair under
$\Q_{\beta,\Gamma}(\,\cdot\given A_3(o))$ is exactly the
law obtained from the physical \HMSOS height $\cH^0$ by setting
\[
 \phi^0=\Tone\cH^0,\qquad
          (\cH,\phi)=(\cH^0-\phi^0(0),\phi^0-\phi^0(0))\,.
\]
Moreover
\begin{equation}\label{eq:physical-palm-probability}
 \Q_{\beta,\Gamma}(A_3(o))=\frac{n^2}{N_0},\qquad
 \E_\mathrm{phys}F(\cH,\phi)
   =\frac{N_0}{n^2}\,
       \Q_{\beta,\Gamma}
          \big(\one_{A_3(o)}F(\cH,\phi)\big)\,,
\end{equation}
for every bounded observable of the pair with this common lift.
The displacement field $\cH-\phi$ is stationary in physical base
coordinates. Conditioning on $A_3(o)$ preserves any local
convergence of the diagonal decorated law, since its probability
is bounded away from zero in the prescribed slope range.
Here local convergence after conditioning still refers to diagonal
coordinates; converting to physical windows is a separate step.
\end{corollary}
\begin{proof}
The joint diagonal law is invariant under the translations of
$\Z^3/(\Gamma+\Z d)$. These act transitively on its $N_0$ black
triangles. Every reference tiling has precisely $n^2$ lozenges
of normal $x_3$, proving the first equality in
\cref{eq:physical-palm-probability}.

In any common lift, mark a uniformly chosen normal-$x_3$ face
$q$ of $\phi$. Translate the pair so that $q$ is the horizontal
face over $0$ at height zero, and denote the resulting observable
by $F_q$. There is exactly one such face above each base square
modulo $n\Z^2$. Averaging $F_q$ is therefore the sampling rule
of \cref{prop:physical-hmsos-law}, followed by the
common change of vertical origin displayed in the statement.
On the other hand, translation invariance on black triangles gives
\[
 \Q_{\beta,\Gamma}
    \bigg(\frac1{n^2}\sum_{q:\,\mathrm{normal}(q)=x_3}F_q\bigg)
 =\frac{N_0}{n^2}\Q_{\beta,\Gamma}
                             (\one_{A_3(o)}F_o)\,.
\]
After the projected root has been fixed, the lift of the marked
face is unique, so the last observable is precisely the one in
the statement. This proves the identity. Physical translation
of the pair re-roots both heights by the same value of $\phi$;
their difference is merely shifted, proving stationarity.
Finally, $A_3(o)$ is a cylinder event. Expectations conditioned
on it are ratios of local expectations, with denominator tending
to $1/(1+\theta_1+\theta_2)>0$ along convergent slopes.
\end{proof}

\subsubsection{A direct displacement estimate for HMSOS}\label{sec:displacement-hmsos}

For this subsection write $\cH$ for the height of the \HMSOS input
in direction $x_3$, and $\phi$ for the height of its monotone image
in that same direction. The sorting direction is $x_1$. Thus $\cH$
and $\phi$ are the physical column heights, not the column magnetizations
in the sorting direction. Let $q_z$ be the reference plaquette of
normal $x_3$ over the base square $z=(x_1,x_2)$.

\begin{lemma}[A displacement crosses a single defect]\label{lem:displacement-support}
For a finite planar fiber, if $|\cH(z)-\phi(z)|\geq r$ with integer
$r\geq1$, some defect component contains $r$ consecutive transverse
squares, and its footprint contains $q_z$. Consequently this
component has $|D|\geq r$ and $E(D)\geq2r$.
For the periodic fiber of \cref{sec:periodic-hmsos}, the same
conclusion holds for $1\leq r\leq \Theta_1$. If the displacement is at
least $\Theta_1$, the component contains a noncontractible transverse loop.
\end{lemma}
\begin{proof}
Apply \cref{lem:sos-column-component} with the inner stage empty,
$\cS=\cH$. Its proof shows that the differing vertical interval
crosses consecutive transverse squares in one outer defect $D$,
with $q_z\in P_h(D)$. The interval has
$|\cH(z)-\phi(z)|$ cubes. In the plane the consecutive transverse
squares are distinct. On the periodic transverse torus their first
repetition is after exactly $\Theta_1$ steps, since
$(0,k)\in\Gamma_\mathrm{s}$ holds precisely for $k\in \Theta_1\Z$.
The $\Theta_1$ distinct squares with their closing edge form a
noncontractible loop. The energy conclusion follows from
$E(D)\geq2|D|$.
\end{proof}

For the rest of the \HMSOS reconstruction, use the witness event
$\cB_x(s)$ of \cref{sec:interacting-witnesses} with the family
restricted to the \texttt{bad} witnesses $W_h(D)$. Thus the event requires
one such witness containing triangle $x$ and at least $s$ lozenges.
This is a subfamily of \cref{sec:hmsos-rg-verification}, so all
the interacting witness bounds remain valid. Here $D$ ranges over all nonempty connected
transverse supports with $|P_h(D)|\geq B_\beta|D|$, without requiring
$z_{\beta,h}(D)>0$. The construction is geometric and applies to
all these supports; this convention also includes unions of supports
from two copies of the gas. The corresponding event
$\cB_q(s)$ for a reference lozenge requires that it contain
that lozenge. These events concern possible witnesses in the
reference tiling, whether or not their defects are present in the
sampled fiber.

\begin{proposition}[Displacement tail reduced to reference witnesses]\label{prop:displacement-witness}
For all sufficiently large $\beta$, uniformly in the reference
surface and the allowed planar supports,
\begin{equation}\label{eq:quenched-displacement-witness}
 \pi_{\beta,h}(|\cH(z)-\phi(z)|\geq r)
       \leq Ce^{-c\beta r}
                     +\one_{\cB_{q_z}(B_\beta r)}\,.
\end{equation}
For the full periodic fiber the same assertion holds for
$1\leq r\leq \Theta_1$, provided $n$ is large enough for the witnesses
of \cref{prop:periodic-witness}.

Let $\nu_n$ be any probability law on the periodic reference
tilings in the specified sector, and define the averaged witness tail
\[
 \alpha_n(s)=\frac1{2N_0}\sum_{x\in\T_\mathrm{d}}
                                      \nu_n(\cB_x(s)),
 \qquad N_0=n(n+\Theta_1+\Theta_2)\,,
\]
where $x$ runs over elementary triangles. In the joint law obtained
by sampling $h$ from $\nu_n$ and then its fiber, choose a uniform
physical base square $Z$ independently. Then, for $1\leq r\leq \Theta_1$,
\begin{equation}\label{eq:annealed-displacement-witness}
 \P(|\cH(Z)-\phi(Z)|\geq r)
       \leq Ce^{-c\beta r}+
                      \frac{N_0}{n^2}\alpha_n(B_\beta r)\,.
\end{equation}
For a translation-invariant joint law the same bound holds at every
fixed base square. In particular it applies with the interacting
image law for $\nu_n$.
\end{proposition}
\begin{proof}
By \cref{lem:displacement-support}, the displacement event
requires a component of size at least $r$ containing $v_+$ or
$v_-$. If this component is \texttt{bad}, its witness contains $q_z$ and
has size at least $|P_h(D)|\geq B_\beta|D|\geq B_\beta r$.
Otherwise $\Arr_h(D^\dagger)\leq |P_h(D)|<B_\beta|D|$, and the activity bound
gives, for all sufficiently large $\beta$,
\[
 z_{\beta,h}(D)\leq
       2e^{-(\beta-2B_\beta e^{-\beta/64})|D|}
                  \leq2e^{-\beta|D|/2}\,.
\]
The deletion estimate of \cref{lem:deletion} bounds the probability
that a specified component is present by its activity. There are at
most $2\cdot64^m$ connected supports of size $m$ through the two
specified squares $v_\pm$, in the plane or on the quotient. Summing
$2\cdot64^m\cdot2e^{-\beta m/2}$ over $m\geq r$ gives at most
$Ce^{-c\beta r}$ for large $\beta$, proving
\cref{eq:quenched-displacement-witness}: either the component
through $v_\pm$ is \texttt{bad}, and then the witness event holds,
or it is \texttt{good}, and then its activity is exponentially
small in its size.

There are exactly $n^2$ reference plaquettes of normal $x_3$, one
for each physical base square, and $2N_0$ reference triangles.
For each fixed reference tiling,
\[
 \frac1{n^2}\sum_z\one_{\cB_{q_z}(s)}
   =\frac1{2n^2}\sum_x
        \one_{\{x\text{ belongs to a normal-}x_3\text{ lozenge}\}}
                         \one_{\cB_x(s)}
   \leq\frac1{2n^2}\sum_x\one_{\cB_x(s)}\,.
\]
A witness is a tiling pattern, so containing one triangle of a
reference lozenge is equivalent to containing the whole lozenge;
this justifies the equality. Average this identity over $\nu_n$
and use the quenched bound. This proves
\cref{eq:annealed-displacement-witness} without assuming stationarity
or replacing a randomly located reference face by a deterministic
one. Stationarity makes the spatial average equal to each single-site
probability.
\end{proof}

\subsubsection{Conditional reconstruction from witness tails}
\label{sec:hmsos-reconstruction}

We next derive conditional locality from the interacting witness
estimate, that is, the annealed local approximation hypothesis of
\cref{lem:joint-limit} for the \HMSOS fiber. Recall from
\cref{lem:polymer-factorization} that, given the reference surface
$h$, the fiber is an exact hard-core gas of defects; the surface
data to be reconstructed in a window are the defects whose
footprints meet that window. In the diagonal coordinates encode a
present defect by its reference footprint and its
spin changes relative to the reference surface. These data
are invariant under a common lift translation, so they are
well-defined data of the pair modulo $\tau_d$. For a finite set $A$
of reference triangles let $\scD_A$ be the list of the present
defects whose footprints meet $A$, each with its full spin change.
There are at most $25|A|$ such defects by
\cref{lem:footprint-overlap}. The list takes values in a countable
set; an empty list is allowed.

The approximating kernel is the law of the same list in the gas
truncated to a ball. For $A\subset B(o,r/2)$ let $\pi_h^{(r)}$ be
the gas restricted to defects with footprint in $B(o,r)$,
and let $K_{A,r}$ be its law of $\scD_A$. Its activities, possible
defects and normalization are determined by the reference tiling in
$B(o,r+C)$: the footprint of a defect contains the charges and all
neighbor constraints needed to compute its isolated activity, and
compatibility of two defects is read off their footprints. This is
the same locality argument as in \cref{lem:periodic-locality}. For
fixed $r$ and sufficiently large $n$, this ball has a planar lift,
so $K_{A,r}$ is the same kernel in the plane and on the torus; this
is what allows the finite-volume kernels to converge to the kernel
of the limiting planar tiling. Noncontractible defects retain their period
identifications as part of the countable data; they are absent from
every fixed-radius kernel once $n$ is large, and the proposition
below shows that they are also unlikely in the complete fiber.

The mechanism of the proof is the one of \cref{prop:pi_decay}: two
independent copies of the gas are compared by exchanging the
defects of one copy with those of the other on the components, in
the union of the two defect configurations, which meet $A$. The
exchange is a weight-preserving involution unless such a component
reaches the boundary of $B(o,r)$, and a component reaching that far
is a chain of overlapping defects, a ``disagreement chain'', which
is either \texttt{thin}, and then costs its energy by the deletion
estimate of \cref{lem:deletion}, or \texttt{fat}, and then produces
a witness in the reference tiling. This is where the witness events
enter the reconstruction.

\begin{proposition}[Conditional locality including noncontractible defects]
\label{prop:hmsos-conditional-locality}
For sufficiently large $\beta$, the complete periodic fiber and
the restricted gas satisfy, for $r\leq c_*n$,
\begin{equation}\label{eq:conditional-locality-witness}
 \big\|\cL_{\pi_{\beta,h}}(\scD_A)
                       -K_{A,r}(h)\big\|_\mathrm{TV}
 \leq C|A|e^{-c\beta r/B_\beta}
               +C\sum_{x\in A}\one_{\cB_x(cR_\beta r)}\,.
\end{equation}
For bounded observables $F,G$ of $\scD_A,\scD_{A'}$,
respectively, whose supports have torus distance at least $r$,
\begin{equation}\label{eq:conditional-covariance-witness}
 |\cov_{\pi_{\beta,h}}(F,G)|
 \leq C\|F\|_\infty\|G\|_\infty
 \bigg(|A|e^{-c\beta r/B_\beta}
                +\sum_{x\in A}\one_{\cB_x(cR_\beta r)}\bigg)\,.
\end{equation}
In particular, for any translation-invariant law $\nu_n$ of the
reference tiling, the expectation of the right side of
\cref{eq:conditional-locality-witness} is at most
\begin{equation}\label{eq:annealed-conditional-locality}
 C|A|\big(e^{-c\beta r/B_\beta}
                              +\alpha_n(cR_\beta r)\big)\,.
\end{equation}
The same averaged bound, multiplied by
$C\|F\|_\infty\|G\|_\infty$, applies to the covariance.
\end{proposition}
\begin{proof}
Apply \cref{prop:pi_decay} to the full and restricted allowed
families, with $\mathcal A$ the defect types whose footprints meet
$A$. The forbidden types are those whose footprints leave $B(o,r)$.
Thus it suffices to bound the probability, in two independent
copies, of a union component joining these two type sets. The
defects keep their weights under the swap.

Any exceptional component has a footprint containing some $x\in A$
and another footprint leaving $B(o,r)$. Its footprint union is
connected and has diameter at least $r/2$, since
$A\subset B(o,r/2)$. Let $U$ be the union of its transverse
supports, $m=|U|$, and $s=|P_h(U)|$. The set $U$ is connected,
since consecutive defects of the chain are incompatible and
therefore have intersecting or $*$-adjacent supports, and $P_h(U)$
contains all the component's footprints, so
$\diam P_h(U)\geq r/2$. We now distinguish whether the footprint of
$U$ is large relative to its support.

If $s\geq B_\beta m$, the union $U$ is \texttt{bad} in the sense of
\cref{prop:periodic-witness}, which constructs a witness containing
$x$ with at least $cR_\beta r$ lozenges. This is precisely one of
the events on the right side of
\cref{eq:conditional-locality-witness}, including when $U$ is noncontractible.
Note that no admissibility of $U$ as a single defect is used: the
witness construction is geometric and applies to every connected
transverse support, as stipulated in
\cref{sec:displacement-hmsos}.

Otherwise $s<B_\beta m$. For a specified $U$ there are at most
$3^m$ color-occupancy assignments (each square of $U$ belongs to a
defect of color $1$, of color $2$, or of both), which determine the
individual defect supports in each copy. Their total sizes lie
between $m$ and $2m$. Applying the neighborhood bound separately in the
two copies gives $\sum_D \Arr_h(D^\dagger)\leq Cs<CB_\beta m$ for the sum of
the neighborhood gradient sums of the defects in the chain. The deletion estimate
of \cref{lem:deletion}, applied in each copy, bounds the probability
of a specified compatible collection by the product of its
activities, and the activity bound gives, for each assignment,
\[
 \prod_D z_{\beta,h}(D)
 \leq\exp([-\beta+2\log2+CB_\beta e^{-\beta/64}]m)
 \leq e^{-\beta m/2}\,,
\]
for all sufficiently large $\beta$. It remains to count the
possible $U$. Since a footprint of one of these defects contains
$x$, $U$ meets one of the $25$ transverse squares specified in
\cref{lem:footprint-overlap}. A connected set of $m$ squares
containing a specified square can be chosen in at most $64^m$ ways,
so there are at most $25\cdot64^m$ choices for $U$. A connected
pattern of $s$ unit tiles has diameter at most $Cs$, also on the
quotient, so $r/2\leq\diam P_h(U)\leq Cs<CB_\beta m$ and hence
$m\geq cr/B_\beta$. Summing $3^m\cdot25\cdot64^m\cdot e^{-\beta m/2}$
over $m\geq cr/B_\beta$ gives $Ce^{-c\beta r/B_\beta}$ for large
$\beta$. A union bound over $x\in A$ proves
\cref{eq:conditional-locality-witness}. The argument uses the
complete periodic gas throughout; a noncontractible component is included
in exactly the same \texttt{good}-or-witness alternative.

For covariance apply \cref{eq:switching-covariance} to the types
whose footprints meet $A$ and $A'$. A connecting component has
diameter at least $r$, so the same good-or-witness estimate applies.
Finally translations act transitively on black triangles and on
white triangles separately. Each of their two witness probabilities
is at most $2\alpha_n(s)$, which proves the averaged assertions.
\end{proof}

We now apply the dimer theorem and the preceding locality estimate
to the interacting \HMSOS image law. This gives the diagonal height
limit and the joint limit with its defects. The conversion to
physical-coordinate windows is carried out in
\cref{sec:physical-local-limit}, and the field's coordinate change
in \cref{sec:quantitative-coordinate-change}.

\begin{proposition}[The diagonal \HMSOS tiling has a \GFF limit]
\label{prop:hmsos-diagonal-gff}
\Cref{thm:dimer-height-gff} applies to the \HMSOS \Tetris-image
law, which is the fixed-sector law of the coefficient family just
verified by \cref{prop:exact-sector}: for every positive slope and
all sufficiently large $\beta$, its centered diagonal height field
converges, modulo constants, to the field \cref{eq:dimer-height-gff},
in the topology stated there.
\end{proposition}
\begin{proof}
Apply \cref{thm:coefficient-rg}(i)--(ii) to the coefficient
family just verified.
\end{proof}

\begin{corollary}[Joint diagonal limit and displacement moments]
\label{cor:hmsos-joint-exponential}
For every positive limiting physical slope and all sufficiently
large $\beta$, the \HMSOS reference tiling together with its local
defects has a unique translation-invariant mixing limit in
diagonal coordinates, independent of the approximating windings.
Its conditioning at a horizontal reference face also converges.
The physical-root displacements $\cH(0)-\phi(0)$ are tight. Every
subsequential limiting displacement $R$ satisfies, for integers
$r\geq1$,
\begin{equation}\label{eq:limiting-displacement-exponential}
 \P(|R|\geq r)
       \leq Ce^{-c\beta r}+Ce^{-cB_\beta r},\qquad
 \E|R|^p\leq C_p e^{-c'B_\beta}\quad(p>0)\,.
\end{equation}
No uniform moment assertion about finite-volume displacements
beyond a transverse period is needed here.
\end{corollary}
\begin{proof}
\Cref{thm:fixed-sector-limit} supplies the unique mixing tiling
limit. By \cref{thm:interacting-witness-tail}, the annealed local
approximation error and conditional covariance in
\cref{prop:hmsos-conditional-locality} are bounded, with the
observable norms included for the covariance, by
\cref{eq:annealed-conditional-locality} with
$\alpha_n(cR_\beta r)\leq Ce^{-cR_\beta r}$.
The approximating kernels depend on finitely many reference tiles,
and this bound tends to zero uniformly in sufficiently large volumes
as $r\to\infty$. The finite-volume coupling is translation invariant,
since the reference law is translation invariant and the complete-fiber
kernel is translation covariant. Thus \cref{lem:joint-limit} gives the unique
stationary mixing joint limit, including tightness of the countable
defect lists. The horizontal-face event has probability tending to
$\rho_3>0$ by \cref{cor:physical-palm}, so conditioning preserves
local convergence.

For each fixed integer $r$, eventually $r\leq \Theta_1$, so
\cref{eq:annealed-displacement-witness,eq:interacting-witness-exponential} give
\[
 \P(|\cH(0)-\phi(0)|\geq r)
              \leq Ce^{-c\beta r}+Ce^{-cB_\beta r}\,.
\]
We used the physical stationarity in
\cref{prop:physical-hmsos-law} and the bounded ratio
$N_0/n^2$. First choose $r$ large and then take the volume large;
this proves tightness. Since the variables are integer-valued,
the bound passes to any subsequential limit at each integer $r$.
For $p>0$, summation of
$\big(r^p-(r-1)^p\big)\P(|R|\geq r)$ gives
\cref{eq:limiting-displacement-exponential}. This summation is
performed for the limit law; it makes no assertion about rare
finite-volume displacements exceeding the period.
\end{proof}

\begin{remark}[Polynomial witness bounds]
For stationary reference laws in a compact positive-slope range,
the same reconstruction argument works if $\nu_n$ converge locally
to a mixing law and, for some $q>0$,
\[
 \limsup_{n\to\infty}\alpha_n(s)\leq C(1+s)^{-q},\qquad\mbox{for all $s\geq1$}\,.
\]
Indeed the error in \cref{eq:annealed-conditional-locality} still
tends to zero with $r$, so \cref{lem:joint-limit} gives the mixing
joint diagonal limit. For the displacement at a uniform physical
root, \cref{eq:annealed-displacement-witness} gives tightness and,
in any limiting law,
\[
 \P(|R|\geq r)\leq Ce^{-c\beta r}+C(1+B_\beta r)^{-q},\qquad
 \E|R|^p\leq C_p e^{-c\beta}+C_{p,q}B_\beta^{-q}\quad(0<p<q)\,.
\]
The tail passes to the limit at each fixed integer $r$, since
eventually $r\leq\Theta_1$; summing it against
$r^p-(r-1)^p$ proves the moment bound. The same holds at a fixed
physical root when the joint physical laws are stationary. These
are moments of the limiting displacement, with no assertion about
finite-volume moments beyond a transverse period.
\end{remark}

\subsubsection{The physical HMSOS local limit and displacement transfer}
\label{sec:physical-local-limit}

The results of \cref{sec:interacting-witnesses} verify the witness
requirements for \HMSOS. We first check that physical rooting is
an almost surely local operation on the limiting decorated tiling.
This is a qualitative coordinate statement; a quantitative change
of coordinates for a distribution-valued height field is a separate
issue.

\begin{lemma}[The limiting plane and finite coordinate fibers]
\label{lem:positive-density-proper-lift}
Let a stationary ergodic lozenge law have orientation densities
$\rho_1,\rho_2,\rho_3>0$. Its stepped-surface lift, after fixing
one lifted vertex, is almost surely asymptotic to the plane with
normal $\rho=(\rho_1,\rho_2,\rho_3)$. Every coordinate fiber has
occupied and vacant ends. In particular the lift is a vertical
height $\phi:\Z^2\to\Z$, and, up to its additive constant,
\begin{equation}\label{eq:physical-reference-plane}
 \sup_{|z|\leq R}
    \Big|\phi(z)+\frac{\rho_1}{\rho_3}z_1
                         +\frac{\rho_2}{\rho_3}z_2\Big|=o(R)
                      \quad\text{almost surely}\,.
\end{equation}
After marking one horizontal face, the heights of $\phi$ on any
fixed finite set of physical base squares are determined by a
finite, configuration-dependent diagonal patch.
\end{lemma}
\begin{proof}
Use coordinates $t=(x_1-x_3,x_2-x_3)$ on the diagonal plane.
The lifted surface is
\[
        L(t)=(t_1+f(t),t_2+f(t),f(t))\,,
\]
where $f$ is continuous, affine on each lozenge, and universally
Lipschitz. On faces with normals $x_1,x_2,x_3$ its gradients are,
respectively, $(-1,0),(0,-1),(0,0)$. Fix $f(0)=0$.
The increments of $f$ on the integer lattice form a bounded
stationary additive cocycle. Write
$m_i=\E(f(e_i)-f(0))$. Stationarity gives
$\E f(t+k)=\E f(t)+m\cdot k$ for $k\in\Z^2$.
Consequently $\E f(t)-m\cdot t$ is periodic and its
gradient has zero average over a fundamental cell. Averaging the
three possible gradients above over that cell gives
\[
                         m=(-\rho_1,-\rho_2)\,.
\]
Indeed all three projected lozenge types have area one, so their
mean area fractions are their matching densities: assign the
unit area of each lozenge to its black-triangle anchor and use
translation invariance. Bounded gradients justify the averages
and differentiation of expectations.

Set $f_R(t)=f(Rt)/R$. The cell averages of $\nabla f$ form a
bounded stationary ergodic field with mean $m$. The pointwise
ergodic theorem for the two lattice shifts gives, almost surely,
\[
           \int_Q\nabla f_R(t)\,\d t\longrightarrow |Q|m
\]
simultaneously for all rectangles $Q$ with rational endpoints.
Indeed complete cells in $RQ$ give the ergodic averages, and
its $O(R)$ boundary cells contribute $O(R^{-1})$ after rescaling.
The uniform gradient bound then gives distributional convergence
$\nabla f_R\to m$ by approximation with rectangular step functions.
The functions $f_R$ are equi-Lipschitz and vanish at zero, so
Arzel\`a--Ascoli makes them precompact locally uniformly.
Every subsequential limit has gradient $m$ and value zero at the
origin, hence equals $m\cdot t$. Uniqueness of this limit gives
$\sup_{|t|\leq R}|f(t)-m\cdot t|=o(R)$ almost surely.

It follows that $\rho\cdot L(t)=o(|t|)$, uniformly in direction.
More directly, at a point $x$ the signed height above the lifted
surface along its diagonal fiber is
\[
 x_3-f(x_1-x_3,x_2-x_3)
            =\rho\cdot x+o(|x|)\,.
\]
Along any coordinate line this tends to opposite infinities at
its two ends, since the corresponding $\rho_i$ is positive.
The occupied set is decreasing in every coordinate, so every
coordinate fiber has precisely the required ends and one
transition. The use of cube centers avoids a boundary convention
at vertices; the error from doing so is bounded. Thus the vertical
height $\phi$ is finite on every base square.

For points on its horizontal faces we have
$\rho_1z_1+\rho_2z_2+\rho_3\phi(z)
=o(|z|+|\phi(z)|)+O(1)$. Since $\rho_3>0$, first absorb the
error proportional to $|\phi(z)|$ to obtain
$|\phi(z)|\leq C(1+|z|)$. Substitution then proves
\cref{eq:physical-reference-plane}, uniformly on balls.
Finally, lift connected diagonal patches starting from the marked
horizontal face. Each desired physical square has a horizontal
face at a finite diagonal location. A sufficiently large connected
patch contains these finitely many faces and paths to the root.
Their lifted positions, hence their heights, are fixed by that
patch. Any other monotone tiling with the same lifted patch has
the same horizontal faces over those base squares; uniqueness
of a coordinate-fiber transition prevents a competing value.
\end{proof}

\begin{lemma}[Physical reconstruction from finitely many defects]
\label{lem:physical-reconstruction}
Suppose stationary decorated diagonal laws $\Q_n$ converge locally
to an ergodic law $\Q$ with positive orientation densities $\rho_i$.
Assume that their physical laws are obtained by conditioning on
$A_3(o)$ and choosing the common lift $\phi(0)=0$, and are stationary
under physical translations followed by this common rerooting.
Suppose that under $\Q(\cdot\given A_3(o))$, almost surely the data
in every finite physical window are determined by a finite diagonal
patch and finitely many finite planar defects, in the
countable spaces of \cref{sec:hmsos-reconstruction}, through a map
we call the \emph{decoding}, which must be translation
covariant and agree in every sufficiently large volume with the
same patch and defect values.

Then the physical laws converge locally to the decoded law under
$\Q(\cdot\given A_3(o))$, which is stationary and ergodic. A further
vertical rerooting by a decoded integer preserves local convergence
to the correspondingly rerooted limit. The reference increments are
integrable and
\begin{equation}\label{eq:palm-reference-slopes}
       \E\nabla_i\phi=-\rho_i/\rho_3,\qquad i=1,2.
\end{equation}
Every reconstructed height $\fh=\phi+R$ with $R_0\in L^1$ has these
same mean gradients. If $\P(|R_0|\geq r)\leq Ce^{-cr}$, then $R_0$
has all finite moments, $|R_z|=O_\omega(\log(2+|z|))$, and $\fh$ has
the same limiting plane as $\phi$.
\end{lemma}
\begin{proof}
Conditioning preserves diagonal local convergence because $A_3(o)$
is a cylinder event with probability tending to $\rho_3>0$.
By \cref{lem:positive-density-proper-lift}, finitely many physical
reference faces are found in an almost surely finite diagonal
patch. Restrict this patch to a fixed ball and its defects
to finitely many values. The resulting cylinder events exhaust the
conditioned limit law, and on each the physical decoding is the
same finite computation in every sufficiently large volume.
Local convergence on these events, followed by exhaustion, proves
physical local convergence. Specified finite planar defects embed
in every sufficiently large torus; no noncontractible defect has
been discarded in advance.
For a further integer rerooting, first restrict that integer to
$[-L,L]$. A rerooted cylinder test then uses finitely many translated
windows, so the same argument applies; let $L\to\infty$ afterwards.
This also permits passage of the finite-volume physical translation
action to the limit, proving stationarity.

For ergodicity, pull a physically invariant event back to $A_3(o)$
and saturate it under diagonal translations. On the proper lift,
any two horizontal faces differ by a physical translation followed
by common rerooting. The saturation therefore agrees with the
original event on $A_3(o)$, and is trivial by diagonal ergodicity.
For invariance modulo null sets, first remove their countably many
translates; stationarity and the conditioning identity keep this
exceptional set null.

The mean slopes follow by counting the wall plaquettes. For $i=1,2$,
send one unit from each normal-$x_i$ plaquette to the horizontal
reference face over its higher adjacent base square. A face over
$z$ receives exactly $\phi(z)-\phi(z+e_i)\geq0$ units. Let $m_i(x,y)$
be this transport between the black-triangle anchors of the two
lozenges. It is translation covariant, so stationarity and Tonelli give
\[
 \rho_i=\E_\Q\sum_y m_i(o,y)
       =\E_\Q\sum_x m_i(x,o)
       =\rho_3\E_{\Q(\cdot\mid A_3(o))}
                              [\phi(0)-\phi(e_i)].
\]
This proves integrability as well as \cref{eq:palm-reference-slopes}.
The displacement $R$ is stationary under the common rerooting, so
its integrable increments have mean zero. Finally, an exponential
tail gives all moments and, for $K>3/c$,
\[
 \sum_{z\in\Z^2}\P(|R_z|>K\log(2+|z|))<\infty.
\]
Borel--Cantelli and \cref{eq:physical-reference-plane} give the
last assertions. The argument applies to any fixed finite list
of reconstructed heights.
\end{proof}

\begin{theorem}[Unique ergodic physical \HMSOS local limit]
\label{thm:hmsos-physical-local-limit}
Fix positive slopes $\theta_1,\theta_2$ and sufficiently large
$\beta$. Along every admissible sequence of windings converging
to these slopes, the physical joint laws of $(\cH,\phi)$, rooted by
$\phi(0)=0$ in one common lift, converge locally to the same law.
The limit is stationary and ergodic under physical base
translations followed by this common rerooting. Its marginal
$\cH-\cH(0)$ is the unique local limit of the physical \HMSOS laws.
Moreover,
\begin{equation}\label{eq:physical-limit-slopes}
 \E\nabla_i\phi=\E\nabla_i\cH=-\theta_i,
                   \qquad i=1,2\,.
\end{equation}
The stationary displacement $R_z=\cH(z)-\phi(z)$ satisfies
\cref{eq:limiting-displacement-exponential}, is in every
$L^p$, and almost surely
\begin{equation}\label{eq:physical-displacement-logarithmic}
                |R_z|\leq C(\omega)+C_\beta\log(2+|z|)\,.
\end{equation}
Both physical heights have the limiting plane
$-\theta_1z_1-\theta_2z_2$, and $\Tone\cH=\phi$ almost surely.
\end{theorem}
\begin{proof}
Let $\Q$ be the mixing decorated diagonal limit of
\cref{cor:hmsos-joint-exponential}. Its densities are
$\rho=(\theta_1,\theta_2,1)/(1+\theta_1+\theta_2)$, by the
prescribed orientation counts. Condition on $A_3(o)$ and use the
common lift of \cref{cor:physical-palm}. By
\cref{lem:positive-density-proper-lift}, the reference horizontal
faces $q_z$ over any finite set of physical squares are determined
in an almost surely finite diagonal patch.

Read the defects at these faces. If $\Delta_D$ is
the $\{0,1\}$-spin change specified by a defect, the exact decoding is
\[
 \cH(z)-\phi(z)
       =\sum_{D:\,P_\phi(D)\ni q_z}
                          \sum_{k\in\Z}\Delta_D(z_1,z_2,k)\,.
\]
Indeed the difference of the two graph columns is one interval
of cubes ending at the reference transition, and its signed
cardinality is $\cH(z)-\phi(z)$. Its consecutive transverse
squares lie in one defect component whose footprint contains
$q_z$, by \cref{lem:displacement-support}. Thus no remote defect
contributes, and the displayed sum is finite.

This decoding agrees with the finite-volume one on each specified
patch and list of defects. The rooting and stationarity
hypotheses of \cref{lem:physical-reconstruction} are
\cref{cor:physical-palm}, so that lemma gives the unique stationary
ergodic physical limit and its marginal rooted by $\cH(0)=0$.
Monotonicity in $x_2$ and the height constraints are closed local
conditions and pass to the limit. The displacement tail is
\cref{cor:hmsos-joint-exponential}; the same reconstruction lemma
therefore gives its moments, \cref{eq:physical-displacement-logarithmic},
the common limiting plane and the mean slopes
\cref{eq:physical-limit-slopes}. In particular, since $\theta_1>0$,
every sorting column of $\cH$ has occupied and vacant ends.

It remains to check $\Tone\cH=\phi$. At each cube the spin change
is the recorded change, by the local decoding. Take the countable
intersection of these probability-one decoding events over all
integer cubes and reference dots, and of the events that each
defect is finite and admissible. On this single event the
identities hold simultaneously on entire columns. For a fixed
sorting column all its changes lie in the single finite defect
at its reference dot, whose signed charge is zero. Thus the whole
column has the reference charge; the ends and \cref{eq:tetris-map}
give $\Tone\cH=\phi$. No infinite charge sum is passed through a
weak limit.
\end{proof}

The physical coupling just constructed has a stationary ergodic
$L^2$ displacement. By \cref{lem:ergodic-correction}, centering
makes its rescaled field vanish both on smooth tests and in local
$H^{-s}$ for $s>0$. It remains to change the coordinates of the
reference field; \cref{cor:physical-height-transfer} below combines
these two steps for all three models.

\subsection{Quantitative change of coordinates}
\label{sec:quantitative-coordinate-change}

The \RG theorem describes the height of the tiling in diagonal
coordinates, as a function on the plane orthogonal to $(1,1,1)$,
whereas the physical height of \cref{thm:gff-limits} is the
vertical height over the base plane. The two are related by a
random change of coordinates: the base square below a point of the
surface at diagonal position $t$ is $t+e f(t)$, where $f(t)$ is the
height of the surface there, so the map from diagonal to physical
coordinates depends on the surface itself. The strategy, following
\cite[\S7]{LaslierLubetzky24}, is to replace this random map by its
deterministic linear approximation, obtained by replacing $f$ by the
limiting plane, and to show that the error of this replacement is
negligible on the scale of the test functions. The error is
controlled by a sublinear second moment of the height at distance
$|t|$ from the root. The logarithmic variance already proved in
\cref{cor:height-covariance-limit} supplies this bound directly.
Keep the diagonal height $f$ and coordinates
of \cref{lem:positive-density-proper-lift}, with $f(0)=0$, and set
\[
       g(t)=f(t)+\rho_1t_1+\rho_2t_2\,,
\]
the diagonal height with its limiting plane subtracted. The point
$0$ is a vertex of the horizontal face used for Palm rooting, so on
$A_3(o)$ this is exactly the lift with $\phi(0)=0$.

\begin{lemma}[A height moment for the coordinate change]
\label{lem:diagonal-height-second-moment}
In the fixed-sector tiling limits of \cref{thm:coefficient-rg},
\begin{equation}\label{eq:diagonal-height-second-moment}
 \E|g(t)|^2\leq C\log(2+|t|)\leq C(1+|t|)^{1/2}
                                      \qquad(t\in\R^2)\,.
\end{equation}
The same bound, with a changed constant, holds after conditioning
on $A_3(o)$.
\end{lemma}
\begin{proof}
For lattice vertices this is \cref{eq:logarithmic-height-variance},
with the fixed linear coordinate and height normalizations absorbed
in $C$. The affine height has a deterministic Lipschitz constant,
so for a nearest lattice vertex $x$,
$|g(t)-g(x)|\leq C$ and
$\E|g(t)|^2\leq2\E|g(x)|^2+2C^2$.
This extends the logarithmic bound to every $t$.
Finally $\P(A_3(o))=\rho_3>0$, and conditioning any nonnegative
variable on this event increases its expectation by at most
$\rho_3^{-1}$.
\end{proof}

For a smooth compactly supported $F$ of integral zero write
\[
                   X_k(F)=k^{-2}\int g(t)F(t/k)\,\d t\,.
\]
This is the diagonal height field in the normalization $f=x_3$
of the monotone surface over the tiling. Other fixed normalizations of dimer height
require the corresponding constant conversion.
The affine-interpolation conclusion of \cref{thm:dimer-height-gff}
identifies this field with the cellwise height field in the limit.
Its comparison in probability also survives conditioning on $A_3(o)$,
since that conditioning has density at most $\rho_3^{-1}$.

\begin{lemma}[Microscopic conditioning does not change a height limit]
\label{lem:height-palm-conditioning}
If the diagonal fields $X_k$ converge in law to a field modulo
constants, their laws conditioned on $A_3(o)$ have the same
limit. This holds jointly for test functions and in any
distribution topology in which the unconditioned fields are tight.
\end{lemma}
\begin{proof}
The point is that a fixed lattice translation has a vanishing
effect on the rescaled field: it moves the test function by
$O(k^{-1})$, and the additive change in the height disappears
against a test function of integral zero. Bounded Lipschitz
functions of the field are therefore asymptotically invariant
under these translations. Ergodicity then makes their correlation
with the fixed centered local observable
$\one_{A_3(o)}-\rho_3$ tend to zero. We verify these two steps
below; no rate of mixing is needed.

Let $T_z$ translate the diagonal tiling by a fixed lattice vector.
The cocycle identity and $\int F=0$ give
\[
 X_k(F)(T_zM)-X_k(F)(M)
   =k^{-2}\int g(u)
            \big[F((u-z)/k)-F(u/k)\big]\,\d u\,.
\]
By Minkowski's inequality and the logarithmic bound in
\cref{eq:diagonal-height-second-moment}, for each fixed $z$,
\begin{equation}\label{eq:microscopic-height-shift}
 \|X_k(F)\circ T_z-X_k(F)\|_2
       \leq C_F|z|\frac{\sqrt{\log(2+k)}}{k}\longrightarrow0\,.
\end{equation}
The test difference is bounded by $C_F|z|/k$ and supported on an
area of order $k^2$. For $k\geq1+|z|$, the height has $L^2$ norm
$O_F(\sqrt{\log(2+k)})$ on this support.

For a finite vector of tests put $Y_k=\Psi(X_k)$, where $\Psi$
is bounded Lipschitz. Then $Y_k$ is bounded in $L^2$ and
$\|Y_k\circ T_z-Y_k\|_2\to0$ for each fixed $z$.
The density argument in \cref{lem:ergodic-correction} says that,
under the ergodic diagonal law, the mean-zero subspace of $L^2$
is the closed span of $h-h\circ T_{e_i}$. Stationarity gives
\[
 \big|\E[(h-h\circ T_{e_i})Y_k]\big|
   =\big|\E[h(Y_k-Y_k\circ T_{-e_i})]\big|
   \leq\|h\|_2\|Y_k-Y_k\circ T_{-e_i}\|_2\longrightarrow0.
\]
Approximate $\one_{A_3(o)}-\rho_3$ by finite sums of these
differences; the uniform $L^2$ bound on $Y_k$ controls the
approximation error. Hence
$\E[(\one_{A_3(o)}-\rho_3)Y_k]\to0$. Division by $\rho_3>0$
proves equality of the limiting finite-test laws. Finally
$\P(E\mid A_3(o))\leq\rho_3^{-1}\P(E)$ transfers tightness,
proving the distributional assertion.
\end{proof}

\begin{proposition}[The reference-height coordinate change]
\label{prop:reference-coordinate-change}
Put $\alpha=(\rho_1,\rho_2)^{\trps}$, $e=(1,1)^{\trps}$ and
\[
                  A=I-e\alpha^{\trps},\qquad \det A=\rho_3>0\,.
\]
Let the physical limiting coupling be obtained, as in
\cref{cor:physical-palm}, by horizontal-face conditioning of a
stationary mixing diagonal limit with positive densities, whose
reference height has mean gradients $\E\nabla_i\phi=-\theta_i$ and
satisfies \cref{eq:diagonal-height-second-moment}. Then, for every
smooth compactly supported $F$ of integral zero,
\begin{equation}\label{eq:reference-coordinate-change}
          \langle\fh_k^\phi,F\rangle-X_k(F\circ A)
                                  \longrightarrow0
                               \quad\text{in probability}\,.
\end{equation}
The error is uniform for tests of bounded $C^1$ norm in a fixed
compact set. It also tends to zero locally in negative Sobolev
norms of order $s>2$, modulo constants.
\end{proposition}
\begin{proof}
Write $\phi^\mathrm{pc}(u)=\phi(\floor{ u})$ and
$\theta=\alpha/\rho_3$, so that $-\theta\cdot u$ is the limiting
plane of the physical height and $A$, $\alpha$ and $e$ are as in
the statement. The proof expresses both the physical field and the
diagonal field as integrals of the test function over the region
between the surface and its limiting plane, in the two coordinate
systems; the two regions coincide, and the two integrals differ
only because the test function is evaluated at points which differ
by a displacement of order the height of the surface, small on
the scale $k$ of the test.
\begin{enumerate}[label=\textbf{Step~\arabic*.}, ref=\arabic*, wide=0pt, itemsep=1ex, listparindent=\parindent, parsep=0pt]
\item\label[step]{st:cc-volume} \emph{(The signed-volume identity.)}
We compare the signed volume between the occupied set and its
limiting half-space in two coordinate systems. For
$(t,r)\in\R^2\times\R$ put
\[
             z=r-\alpha\cdot t,\qquad
             u=t+ez=At+er\,.
\]
Here $(u,z)\in\R^2\times\R$ are the physical coordinates (base
point $u=(x_1,x_2)$ and vertical height $z=x_3$), and
$t=(u_1-z,u_2-z)$ is the diagonal
coordinate of \cref{lem:positive-density-proper-lift}, so that a
point of the surface at diagonal position $t$ has $z=f(t)$ and
$r=f(t)+\alpha\cdot t=g(t)$. The map $(t,r)\mapsto(u,z)$ is the
composition of the two shears $(t,r)\mapsto(t,z)$ and
$(t,z)\mapsto(u,z)$, so its Jacobian is one. Up to boundaries of
cells, the occupied set $z<\phi^\mathrm{pc}(u)$ becomes $r<g(t)$,
by the preceding description of the surface. The reference
half-space $z<-\theta\cdot u$ becomes $r<0$: indeed
$z+\theta\cdot u=z+\theta\cdot t+(\theta\cdot e)z
=z(1+\theta\cdot e)+\alpha\cdot t/\rho_3=(z+\alpha\cdot t)/\rho_3
=r/\rho_3$, using $\rho_1+\rho_2+\rho_3=1$. The $z$-integral of the
difference of the indicators of the two sets, at fixed $u$, equals
$\phi^\mathrm{pc}(u)+\theta\cdot u$. Integrating this difference
against $k^{-2}F(u/k)$, first in $z$ and then after the change of
variables in $r$, gives the exact signed-volume identity
\[
 I_k(F):=k^{-2}\int
       (\phi^\mathrm{pc}(u)+\theta\cdot u)F(u/k)\,\d u
 =k^{-2}\int\!\int_0^{g(t)}F((At+er)/k)\,\d r\,\d t\,.
\]
The inner integral is oriented when $g(t)<0$. Fubini is
legitimate: the support of $F(u/k)$ meets finitely many physical
base squares, their heights are finite, and the volume between
the two graphs over this compact set is finite. Recall that
\[
 X_k(F\circ A)=k^{-2}\int g(t)F(At/k)\,\d t
             =k^{-2}\int\!\int_0^{g(t)}F(At/k)\,\d r\,\d t\,.
\]
Subtracting this from
$I_k(F)$ thus replaces the integrand by $F((At+er)/k)-F(At/k)$,
whose absolute value is at most
\[
                    |e|\,\|\nabla F\|_\infty |r|/k\,,
\]
by the mean value theorem.

\item\label[step]{st:cc-support} \emph{(Localizing the integral and bounding the error.)}
Fix a compact support set for $F$. Since $g(t)=o(|t|)$ uniformly
at infinity, by \cref{lem:positive-density-proper-lift}, there is
a deterministic $C$ such that, on events $E_k$ of probability
tending to one, all the preimages of this set under the maps
$t\mapsto(At+\lambda eg(t))/k$, $0\leq\lambda\leq1$, lie in
$B(0,Ck)$. Explicitly, choose $\epsilon$ smaller than one half the
smallest singular value of $A$ divided by $|e|$. Almost surely
$|g(t)|\leq\epsilon|t|$ for all $|t|$ beyond a finite random
threshold. Once $Ck$ exceeds that threshold,
$|At+\lambda eg(t)|\geq|At|-|e||g(t)|\geq c|t|$ outside
$B(0,Ck)$, so these points are mapped outside the support of
$F(\cdot/k)$ for large enough fixed $C$; the event $E_k$ is that
the random threshold is at most $Ck$. On $E_k$ the integrand of
\cref{st:cc-volume} vanishes for $|t|>Ck$, and integrating the
bound $|e|\|\nabla F\|_\infty|r|/k$ over $0\leq|r|\leq|g(t)|$
gives $|e|\|\nabla F\|_\infty g(t)^2/(2k)$, whence
\begin{equation}\label{eq:coordinate-error-second-moment}
 |I_k(F)-X_k(F\circ A)|
       \leq\frac{C\|F\|_{C^1}}{k^3}
                          \int_{|t|\leq Ck}g(t)^2\,\d t\,.
\end{equation}
By \cref{lem:diagonal-height-second-moment}, including its
conditioned version, the expectation of the right-hand side under
the horizontal-face conditioning is at most
$C\|F\|_{C^1}k^{-3}\int_{|t|\leq Ck}(1+|t|)^{1/2}\,\d t
\leq C\|F\|_{C^1}k^{-1/2}$. Markov's inequality and
$\P(E_k^c)\to0$ prove convergence in probability of
$I_k(F)-X_k(F\circ A)$ to zero, uniformly over the indicated
tests. No moment bound outside $E_k$ and no rate for the
probability of its complement is being assumed; this is why the
conclusion is convergence in probability rather than in $L^2$.

\item\label[step]{st:cc-center} \emph{(Centering.)}
It remains to compare $I_k(F)$ with the centered physical field.
By \cref{eq:physical-limit-slopes} and $\phi(0)=0$,
$\E\phi^\mathrm{pc}(u)=-\theta\cdot\floor{ u}$. Thus
$I_k(F)-\langle\fh_k^\phi,F\rangle$ is the integral of the bounded
periodic function $\theta\cdot(u-\floor{ u})$ against
$k^{-2}F(u/k)$. Its cell average, a constant, multiplies
$\int F=0$ and contributes nothing. Subtract that average on every
unit square; the remaining function has zero mean on each cell, so
its integral against $F(u/k)$ over a cell with corner $x$ equals
its integral against $F(u/k)-F(x/k)$, which is at most
$C\|\nabla F\|_\infty/k$ in absolute value. Summing over the
$O(k^2)$ cells meeting the support, with the weight $k^{-2}$, bounds
the total by $O(k^{-1}\|F\|_{C^1})$. This proves
\cref{eq:reference-coordinate-change}. In particular no
convergence of the expectation of the random error in
\cref{eq:coordinate-error-second-moment} is required.
The uniform $C^1$ bound also gives convergence in local
$H^{-s}$ modulo constants for $s>2$: on a fixed bounded domain in
two dimensions, $H^s$ embeds continuously into $C^1$ when $s>2$,
and one tests on the closed subspace of integral-zero functions.
Equivalently, subtract a fixed smooth test of integral one to
choose a representative modulo constants.
\qedhere
\end{enumerate}
\end{proof}

\begin{corollary}[\GFF transfer to physical coordinates]
\label{cor:physical-height-transfer}
Under the hypotheses of \cref{prop:reference-coordinate-change},
suppose the centered diagonal height has limit
$\upsigma_{\mathrm d}\GFF\circ L_{\mathrm d}$, jointly on smooth
integral-zero tests and in local $H^{-s}$ modulo constants for
all $s>1$. Then the centered physical reference height has limit
\[
 \frac{\upsigma_{\mathrm d}}{\rho_3}\,
                  \GFF\circ(L_{\mathrm d}A^{-1})\,.
\]
The same holds for any height $\fh=\phi+R$ in the common lift if
$R$ is a stationary ergodic field with $\E R_0^2<\infty$.
Both convergences hold jointly on tests and in local $H^{-s}$
modulo constants for every $s>2$.
\end{corollary}
\begin{proof}
By \cref{lem:height-palm-conditioning}, horizontal-face
conditioning preserves the diagonal limit
$Y=\upsigma_{\mathrm d}\GFF\circ L_{\mathrm d}$.
\Cref{prop:reference-coordinate-change} identifies the physical
reference limit on a test $F$ with $\langle Y,F\circ A\rangle$.
The test still has integral zero, and the distributional
change-of-variables convention is
\[
 \langle Y\circ A^{-1},F\rangle
             =|\det A|\langle Y,F\circ A\rangle.
\]
Since $\det A=\rho_3>0$, the resulting field is
$\rho_3^{-1}Y\circ A^{-1}$. The coordinate error tends to zero
in local $H^{-s}$ for $s>2$, while conditioning preserves tightness
and a fixed invertible linear map preserves the local Sobolev
topology. This proves the reference assertion. For $\fh=\phi+R$,
centering gives the correction $R-\E R_0$, whose rescaled field
vanishes in local $H^{-s}$ for every $s>0$ by
\cref{lem:ergodic-correction}. Slutsky's theorem proves the second
assertion.
\end{proof}

\begin{proof}[Proof of \cref{thm:hmsos-gff-complete}]
\Cref{thm:hmsos-physical-local-limit} supplies the unique physical
coupling, its mean slopes and its stationary ergodic $L^2$
displacement. Its reference tiling has the diagonal \GFF of
\cref{prop:hmsos-diagonal-gff} and the height moment of
\cref{lem:diagonal-height-second-moment}. Thus
\cref{cor:physical-height-transfer} applies with
$\upsigma_{\mathrm d}=\sqrt{a_*/|\det\sA_{\mathrm D}|}$ and
$L_{\mathrm d}=\sA_{\mathrm D}^{-\trps}$, giving the stated field
and topology.
\end{proof}

\subsection{Two steps: SOS reconstruction}

The \SOS reconstruction repeats the one-step program of the
previous subsection, with the single \HMSOS fiber replaced by the
complete two-step history $(\cS,\cH,\phi)$ of
\cref{sec:sos}. The three stages are the same: a physical
rooting that identifies the \SOS law with a tiling law decorated
by its history (\cref{prop:physical-sos-law}); a conditional
locality statement for that history, obtained by the same
switching argument but now applied to the gas of enlarged history
components
(\cref{prop:sos-conditional-locality}); and a local decoding of
the physical heights from the limiting history components
(\cref{prop:sos-physical-joint}). The new feature is that the
displacement $\cS-\phi$ between the \SOS height and its reference
height now splits into two terms,
$(\cS-\cH)+(\cH-\phi)$, carried by the inner and the outer
stage of the history, respectively, and the central geometric
input, \cref{lem:sos-column-carrier}, is that both terms at a
given column are read off a single history component.
In this subsection, unless explicitly called a base component,
a history component means an enlarged history component of
\cref{sec:sos-mixed-coefficients}, obtained by regrouping base
components whose enlarged patterns overlap. These regrouped
components have disjoint patterns and retain all the defects of
their constituent base components.
Throughout, the interacting tiling law is the \SOS image law
$\nu^{\textsc{sos}}_{\beta,\Gamma}$, whose witness tails are
supplied by \cref{thm:coefficient-rg}(iii) via the marked input
of \cref{lem:sos-marked-input}.

\subsubsection{Physical rooting and finite-volume SOS compactness}
\label{sec:sos-rooting}

Let $\fS_\Gamma$ be the integer heights with the two
period identities in \cref{eq:physical-hmsos-periods}, without
the monotonicity requirement. For $\cS\in\fS_\Gamma$,
\begin{align}\label{eq:physical-sos-excess}
 E_\Gamma(\cS)=\sA(\cS)-N_0
       &=2\sum_{z\in Q_n}\big((\nabla_1\cS(z))_+
                                      +(\nabla_2\cS(z))_+\big)\notag\\
       &=E_2(\cS;\cH)+E_1(\cH;\phi),
 \qquad \cH=\Ttwo\cS,\quad\phi=\Tone\cH\,.
\end{align}
To see the first equality, recall that the area of the surface
$\cS$ over one period cell consists of the $n^2$ horizontal
plaquettes together with
$\sum_{z\in Q_n}(|\nabla_1\cS(z)|+|\nabla_2\cS(z)|)$ vertical
ones, and that the period identities force the signed sums
$\sum_{z\in Q_n}\nabla_1\cS(z)=-n\Theta_1$ and
$\sum_{z\in Q_n}\nabla_2\cS(z)=-n\Theta_2$. Writing $|u|=-u+2u_+$ for
each gradient and recalling $N_0=n^2+n(\Theta_1+\Theta_2)$ gives the displayed
sum of positive parts. The last equality is the telescoping of
the area along the two sorting steps,
$\sA(\cS)-\sA(\phi)=(\sA(\cS)-\sA(\cH))+(\sA(\cH)-\sA(\phi))$,
together with $\sA(\phi)=N_0$ for the monotone surface $\phi$,
exactly as in the proof of \cref{prop:sos-history-fiber}.
In particular the excess vanishes precisely when both gradients
are nonpositive, that is, when $\cS$ is itself monotone.

\begin{proposition}[The physical \SOS law and its diagonal root]
\label{prop:physical-sos-law}
For every $\beta>0$, the \SOS partition function on
$\fS_\Gamma^0=\{\cS\in\fS_\Gamma:\cS(0)=0\}$ is finite,
and
\begin{equation}\label{eq:physical-sos-partition-ratio}
 \cZ^{\textsc{sos},\mathrm{phys}}_{\beta,\Gamma}
   =\frac{n^2}{N_0}
             \sum_{M\in\cM_\Gamma}
                          Z^{\textsc{sos}}_\beta(\phi_M)\,.
\end{equation}
Here $\phi_M$ is the lift with $0\leq V(\phi_M)<N_0$.
Sampling a tiling with these fiber weights, sampling its
complete two-step history, and then choosing an independent
uniform physical square $Z\in Q_n$ gives the physical law
by outputting $\cS(x+Z)-\cS(Z)$. The two intermediate surfaces
are rooted using that same translation and additive constant.
All three gradient fields are jointly stationary.

For the joint diagonal law of $(\cS,\cH,\phi)$, the
horizontal-face event $A_3(o)$ has probability $n^2/N_0$.
Conditioning on this event and choosing its face to have
height zero gives the physical triple rooted by the common
constant $\phi(0)$. Increasing exhausting height cutoffs
converge in total variation to the physical law.
\end{proposition}
\begin{proof}
The rooting arguments of \cref{lem:rooting-volume},
\cref{prop:physical-hmsos-law,cor:physical-palm} require only the
period identities, a finite rooted partition function,
translation covariance, and preservation of volume by sorting.
They do not use the monotonicity of the input height.
For \SOS, the first and third properties are immediate. For
finiteness, the same path bound used in
\cref{prop:physical-hmsos-law} applies: at excess $k$, the total
absolute gradient is $n(\Theta_1+\Theta_2)+k$ by
\cref{eq:physical-sos-excess}, giving at most
$(2n(\Theta_1+\Theta_2)+2k+1)^{n^2-1}$ rooted assignments.
Their sum against $e^{-\beta k}$ is finite for every $\beta>0$.
Finally, each $\Ti$ preserves the signed $\{0,1\}$-spin sum on the
quotient, by the sorting-fiber argument in
\cref{lem:rooting-volume}; hence $V(\cS)=V(\cH)=V(\phi)$.

Apply the volume-strip count and averaged-observable argument of
\cref{prop:physical-hmsos-law} to the whole triple and the map
$\Tone\Ttwo$. This proves \cref{eq:physical-sos-partition-ratio},
the sampling rule and joint stationarity. The horizontal-face
count of \cref{cor:physical-palm} involves only $\phi$, so it gives
the same probability $n^2/N_0$ and the common rooting by $\phi(0)$.
Finiteness also gives convergence of exhausting cutoffs: their
normalized laws are restrictions of the full law, at total-variation
distance equal to the probability of the omitted states.
\end{proof}

\begin{proof}[Proof of \cref{thm:sos-rg-input}]
The exact image weight $Z^{\textsc{sos}}_\beta(\varphi)$ is
\cref{prop:sos-history-fiber}, its finiteness is
\cref{cor:sos-fiber-pressure}, and \cref{prop:physical-sos-law}
gives the diagonal quotient, the physical-rooting correspondence and
the removal of the height cutoff. The positive history
expansion, the interval-stack witnesses and the mixed bounds are
proved in
\cref{prop:sos-positive-expansion,prop:sos-witness,prop:sos-mixed-envelopes};
they play the roles of
\cref{prop:periodic-fiber,prop:periodic-witness,prop:mixed-envelopes}
in the one-step case. Finally
\cref{lem:sos-small-locality,prop:sos-canonical-fermions} give the
canonical completion, exact Grassmann expansion, quantitative source
bounds, gauge covariance and fixed-count projection, including the
periodic entropy correction and agreement of small coefficients with
the bulk; these verify \crefrange{it:rg-positive}{it:rg-sources} in
the same way as in the proof of \cref{thm:hmsos-rg-input}. The
marked assertion \cref{it:rg-marked} follows from
\cref{lem:sos-marked-input}.
\end{proof}

\begin{corollary}[Small \SOS excess and joint gradient tightness]
\label{cor:sos-gradient-tightness}
For sufficiently large $\beta$,
\[
 \frac1{n^2}\E_{\textsc{sos},\mathrm{phys}}E_\Gamma(\cS)
       \leq\frac{C}{\beta}e^{-\beta/128}
                            \Big(1+\frac{\Theta_1}n+\frac{\Theta_2}n\Big)\,.
\]
For $i=1,2$ one has
\[
 \E|\nabla_i\cS(0)|
       =\frac{\Theta_i}{n}+2\E(\nabla_i\cS(0))_+,
 \qquad
 \sum_{i=1}^2\E|\nabla_i\cS(0)|
       =\frac{\Theta_1+\Theta_2}{n}+\frac1{n^2}\E E_\Gamma(\cS)\,.
\]
The joint gradient laws of $(\cS,\cH,\phi)$ are tight along
sequences of the prescribed slopes. Their subsequential
limits are stationary curl-free integer gradient fields;
$\cH$ is monotone in the second direction and $\phi$ in both.
No uniqueness, displacement moment bound or preservation
of the mean slopes in the limit is asserted by this
compactness statement.
\end{corollary}
\begin{proof}
The ground states are the monotone heights, whose number is
$G_\Gamma=(n^2/N_0)|\cM_\Gamma|$ by the volume-strip count.
Writing $\cZ_t=\cZ^{\textsc{sos},\mathrm{phys}}_{t,\Gamma}$,
\cref{eq:sos-fiber-pressure,eq:physical-sos-partition-ratio} give
\[
 G_\Gamma\leq\cZ_t\leq G_\Gamma e^{Ce^{-t/64}N_0}.
\]
As in \cref{cor:mixed-pressure}, Jensen's inequality gives directly
\[
 \E_\beta E_\Gamma
 \leq\frac2\beta\log\E_\beta e^{\beta E_\Gamma/2}
 =\frac2\beta\log\frac{\cZ_{\beta/2}}{\cZ_\beta}
 \leq\frac{C}{\beta}e^{-\beta/128}N_0.
\]
Divide by $n^2$ to obtain the first assertion. Stationarity and
the winding identities give $\E\nabla_i\cS(0)=-\Theta_i/n$.
The two gradient identities now follow from $|u|=-u+2u_+$
and \cref{eq:physical-sos-excess}.

Sorting decreases area, so joint stationarity gives uniform
single-edge first-moment bounds for all three gradient fields.
These imply tightness on their countable product space: enumerate
the coordinates as $Y_j$, choose $L_j$ with
$\sup_n\P(|Y_j|>L_j)\leq\varepsilon2^{-j}$, and use the compact
set $\prod_j([-L_j,L_j]\cap\Z)$ and a union bound.
Stationarity, curl-freeness and the stated monotonicities pass to
limits, the latter two being closed local conditions. These
first-moment bounds alone assert neither uniform integrability
nor tightness of the relative heights at the root.
\end{proof}

The identification of the limit, the displacement moments and the
mean slopes are proved in
\cref{prop:sos-physical-joint,thm:sos-gff-complete}.

The remaining \SOS argument uses \cref{thm:sos-rg-input} to
obtain the interacting tiling limit, then reconstructs its complete
two-step history. Both terms in
$\cS-\phi=(\cS-\cH)+(\cH-\phi)$ must be controlled under this joint law.
The compactness just proved guarantees subsequential limits; the
uniqueness of the limit and the quantitative control of the
displacement come from the interacting analysis.

\subsubsection{The SOS tiling limit and interacting reconstruction}
\label{sec:sos-rg-application}

We apply \cref{thm:coefficient-rg} to the verified input. Its
conclusions (i)--(ii) give the limit of the tiling alone; the
reconstruction estimates below use its additional witness
conclusion (iii), which is the only information about the
\emph{interacting} law that enters the conditional locality
argument. As explained after \cref{lem:annealed-moment}, a bound
under the free dimer law would not suffice here.

\begin{proposition}[The \SOS image: local limit and diagonal \GFF]
\label{prop:sos-diagonal-gff}
For positive limiting physical slopes and sufficiently large
$\beta$, the laws $\nu^{\textsc{sos}}_{\beta,\Gamma}$ converge
locally to a unique stationary mixing tiling law
$\nu^{\textsc{sos}}_{\beta,\rho}$, independently of the
approximating windings, with
$\rho=(\theta_1,\theta_2,1)/(1+\theta_1+\theta_2)$.
Its edge cumulants are the bulk \SOS source coefficients and
satisfy \cref{eq:bulk-response-decay}. Its centered diagonal
height has the limit
\begin{equation}\label{eq:sos-diagonal-gff}
 \upsigma_{\mathrm d}^{\textsc{sos}}\,
                  \GFF\circ L_{\mathrm d}^{\textsc{sos}},
 \qquad
 \upsigma_{\mathrm d}^{\textsc{sos}}
       =\sqrt{a_{\textsc{sos}}/|\det\sA_{\mathrm D}^{\textsc{sos}}|},
 \qquad
 L_{\mathrm d}^{\textsc{sos}}
       =(\sA_{\mathrm D}^{\textsc{sos}})^{-\trps}\,.
\end{equation}
Here $a_{\textsc{sos}}>0$ and
$\sA_{\mathrm D}^{\textsc{sos}}$ is the invertible real
matrix of its matched linear dispersion, in the conventions of
\cref{cor:height-covariance-limit}. Convergence holds
jointly on smooth compactly supported integral-zero tests and
in local $H^{-s}$ modulo constants for every $s>1$.
\end{proposition}
\begin{proof}
\Cref{thm:sos-rg-input} identifies the fixed-sector image law and
verifies the coefficient conditions. Apply
\cref{thm:coefficient-rg}(i)--(ii): its density is the limit of the
prescribed orientation proportions, and its covariance parameters
are those of the \SOS coefficient family.
\end{proof}

For reconstruction use the witness family
$\scW_n^{\textsc{sos}}$ and events $\cB_x^{\textsc{sos}}(s)$
of \cref{lem:sos-marked-input}. Recall that
$\cB_x^{\textsc{sos}}(s)$ is the event that the tiling contains a
member of $\scW_n^{\textsc{sos}}$ through the triangle $x$ with at
least $s$ lozenges; it is a property of the reference tiling alone,
whether or not the corresponding history component is present.

\begin{proposition}[Interacting \SOS witness and component tails]
\label{prop:sos-interacting-tails}
Uniformly for sufficiently large tori in the compact slope range,
all triangles $x$, and $s\geq1$,
\begin{equation}\label{eq:sos-interacting-witness-tail}
 \nu^{\textsc{sos}}_{\beta,\Gamma}
          (\cB_x^{\textsc{sos}}(s))
                              \leq Ce^{-c\max(s,B)}\,.
\end{equation}
In the joint law of the tiling and its history
components, the probability that a component whose enlarged
pattern contains $x$ has at least $s$ tiles is at most
\begin{equation}\label{eq:sos-component-tail}
             Ce^{-c\beta^{2/5}s}+Ce^{-c\max(s,B)}\,.
\end{equation}
The estimates include noncontractible patterns and sum all histories.
\end{proposition}
\begin{proof}
\Cref{lem:sos-marked-input} verifies \cref{it:rg-marked}, so
\cref{thm:coefficient-rg}(iii)
proves \cref{eq:sos-interacting-witness-tail}.

For the component tail, condition on the tiling and use the exact
gas of enlarged patterns of \cref{prop:sos-mixed-envelopes}, with each
component labelled by its full history and its
\texttt{thin}/\texttt{fat} class. This is a hard-core gas with
nonnegative activities on patterns with disjoint triangle supports,
so the deletion estimate of \cref{lem:deletion} applies: the
conditional probability that a specified labelled component is
present is at most its activity. We treat the two classes
separately. For the \texttt{thin} class, the activity of a pattern of
size $m$ is at most $e^{-c\beta^{2/5}m}$ and there are at most
$C_0^m$ rooted \texttt{thin} patterns of size $m$ through $x$, by
\cref{prop:sos-mixed-envelopes}; hence the probability that some
\texttt{thin} component of size at least $s$ has its pattern through
$x$ is at most
\[
 \sum_{m\geq s}C_0^me^{-c\beta^{2/5}m}\leq Ce^{-c'\beta^{2/5}s}\,,
\]
for $\beta$ large enough that $C_0e^{-c\beta^{2/5}}<1/2$.
For the \texttt{fat} class we do not use the activity at all, which
may be large: a present \texttt{fat} component of size at least $s$
whose pattern contains $x$ entails the event
$\cB_x^{\textsc{sos}}(s)$, since by construction its pattern
belongs to the chosen witness family. The first bound
\cref{eq:sos-interacting-witness-tail} then controls this event
under the interacting law. Adding the two contributions gives
\cref{eq:sos-component-tail}.
\end{proof}

For a finite set $A$ of reference triangles, let
$\scR_A$ be the list of the enlarged history components whose
pattern meets $A$, each with all its outer and inner defects.
Thus such a component contains the outer
defects on $\phi$, the inner defects on the intermediate surface
$\cH$, and the complete spin changes they specify, so that
the history above the pattern can be reconstructed from the
component alone, by \cref{lem:sos-carried-components}. There are at
most $|A|$ such components, since enlarged patterns have disjoint
triangle supports.
If $A\subset B(o,r/2)$, let $K^{\textsc{sos}}_{A,r}(\varphi)$
be their law in the gas of \cref{prop:sos-mixed-envelopes} restricted to patterns contained
in $B(o,r)$, the analogue of the restricted kernel
$K_{A,r}$ of \cref{sec:hmsos-reconstruction}. For fixed $r$ it
depends only on the tiling in $B(o,r+C)$ and agrees with the
planar kernel for large periods, by
\cref{lem:sos-small-locality}: the coefficients of small patterns
are determined by the pattern and its neighborhood, and a small periodic
pattern has an embedded planar lift.

\begin{lemma}[An enhanced a priori bound for connecting unions]
\label{lem:connecting-union-budget}
Fix the reference tiling. Suppose an exceptional connection from
$A$ to distance $r$ entails a connected pattern $U$ containing some
$x\in A$, with $\diam U\geq r$. Let $\xi$ range over its possible
two-copy collections, with nonnegative probability envelopes
$w(\xi)$ and costs $Q(\xi)$. Suppose that, for some $a>0$, $b\geq0$
and $0<\delta<1$, every fixed $U$ satisfies
\[
       \sum_{\xi:\,\mathrm{union}(\xi)=U}
               e^{aQ(\xi)}w(\xi)\leq e^{b|U|},
\]
and $Q(\xi)<\delta|U|$ forces a reference event
$\cB_x(cRr)$. If connected patterns of size $m$ through $x$
number at most $C_0^m$ and
$\lambda=a\delta-b-\log C_0\geq1$, then the connection probability
is at most
\[
       C|A|e^{-c\lambda r}
                    +\sum_{x\in A}\one_{\cB_x(cRr)}.
\]
The assertion holds on finite quotients with the same constants.
\end{lemma}
\begin{proof}
Off the stated witness events, every contributing collection has
$Q\geq\delta|U|$. Extracting its enhancement bounds the sum for
a fixed $U$ by $e^{-(a\delta-b)|U|}$. A connected union of $m$
unit tiles has diameter at most $Cm$, so $m\geq cr$. The union
bound over $x\in A$ and the rooted patterns gives
$|A|\sum_{m\geq cr}e^{-\lambda m}\leq C|A|e^{-c\lambda r}$.
All internal decompositions are already summed in the hypothesis.
\end{proof}

\begin{proposition}[Conditional locality for the complete \SOS history]
\label{prop:sos-conditional-locality}
For $r\leq c_*n$, the conditional law in the full \SOS fiber obeys
\begin{align}\label{eq:sos-conditional-locality}
 \big\|\cL(\scR_A\given\varphi)
                -K^{\textsc{sos}}_{A,r}(\varphi)\big\|_\mathrm{TV}
 &\leq C|A|e^{-c\beta^{3/10}r}
       +C\sum_{x\in A}
               \one_{\cB_x^{\textsc{sos}}(cRr)}\,.
\end{align}
For bounded functions $F,G$ of $\scR_A,\scR_{A'}$
with $\dist(A,A')\geq r$, the absolute conditional covariance
has the same bound multiplied by $C\|F\|_\infty\|G\|_\infty$.
After averaging over the \SOS image law both bounds are at most
\begin{equation}\label{eq:sos-annealed-locality}
                  C|A|\big(e^{-c\beta^{3/10}r}+e^{-cRr}\big)\,,
\end{equation}
with the same observable factor for the covariance.
\end{proposition}
\begin{proof}
Apply \cref{prop:pi_decay} to the exact gas of enlarged patterns of
\cref{prop:sos-mixed-envelopes}, keeping the complete two-step
histories and their \texttt{thin}/\texttt{fat} labels. Incompatibility
is overlap of triangle supports. The observable types meet $A$,
and the forbidden types leave $B(o,r)$. Thus total variation is
bounded by twice the probability of a two-copy union component
joining these type sets. Its pattern $U$ contains some $x\in A$
and has diameter at least $r/2$. The same connection estimate gives
the covariance bound by \cref{eq:switching-covariance}.

Put $m=|U|$, and let $q$ be the sum of the sizes of its
\texttt{thin} patterns in both copies. We verify
\cref{lem:connecting-union-budget}. Choose $A_\beta=c\beta^{2/5}$
smaller than the \texttt{thin} exponent, and let $\eta_\beta\to0$
bound the \texttt{fat} coefficient and rooted entropy rates,
including the periodic correction $C/n$. Apply
\cref{eq:disjoint-collection-budget} with
$\zeta_\beta=4\eta_\beta+B^{-1/2}$ and its $a_\beta=o(1)$.
For each copy, the patterns are disjoint inside $U$, so that
estimate bounds the enhanced sum by
$e^{(\zeta_\beta+2a_\beta)m}$. Squaring gives
\begin{equation}\label{eq:sos-two-copy-budget}
 \sum_{\cP_1,\cP_2\subset U}
       e^{A_\beta q/2}\prod_{i=1}^2
             \prod_{P\in\cP_i}\bar c^\mathrm{label}_\beta(P)
                   \leq e^{(2\zeta_\beta+4a_\beta)m}.
\end{equation}
Dropping the union and connection restrictions only increases
this sum. Deletion in each copy, by \cref{lem:deletion}, bounds
the probability of any specified pair by its product of activities.
Thus the displayed bound includes all the required probability
envelopes and sums over histories.

It remains to check the witness alternative. Set
$\delta=\delta^\mathrm{g}_\beta=\beta^{-1/10}$, as fixed before
\cref{lem:sos-marked-input}. If $q<\delta m$, unite the six
classes of all \texttt{fat} patterns in both copies. Their total
size is at most $2m$; the uncovered remainder, including the
\texttt{thin} patterns, has size at most $C(\delta+L/B)m$, and
the exceptional-orientation count is at most $CR^2m/L$.
There is a \texttt{fat} pattern, since otherwise $q\geq m$.
Hence $m>B$, and $U$ belongs to the witness family
$\scW_n^{\textsc{sos}}$, whose constants admit these unions of
two six-class families. The net and boundary-layer estimate of
\cref{prop:periodic-witness} gives $m\geq cR\diam U\geq cRr$:
the remainder fraction times $R^2$ tends to zero. Thus this
alternative forces $\cB_x^{\textsc{sos}}(cRr)$.

\Cref{lem:connecting-union-budget}, at distance $r/2$, now applies with
$a=A_\beta/2$, $b=2\zeta_\beta+4a_\beta$ and
$\lambda\asymp\beta^{3/10}$. It proves
\cref{eq:sos-conditional-locality} and the covariance assertion.
Averaging and using \cref{eq:sos-interacting-witness-tail} gives
\cref{eq:sos-annealed-locality}. The a priori bounds, the witness
construction and the connected-pattern count all hold on the
full quotient, so noncontractible unions are included.
\end{proof}

\begin{corollary}[The joint \SOS history in diagonal coordinates]
\label{cor:sos-joint-diagonal}
The \SOS reference tiling and its local history components have
a unique stationary mixing limit, independently of the
approximating windings. All limiting components are finite and
planar. Conditioning on a horizontal reference face also converges.
\end{corollary}
\begin{proof}
We verify the hypotheses of \cref{lem:joint-limit} with
$\varphi_N$ the \SOS image tiling and $\cX_N$ its local
history components. \Cref{prop:sos-diagonal-gff} supplies the mixing
tiling limit. The fixed-radius kernels
$K^{\textsc{sos}}_{A,r}$ depend on finitely many tiling
coordinates and are probability laws on finite planar components:
\cref{lem:sos-small-locality} confines first $\cH$ and then
its inner changes to a finite box for each finite set of reference plaquettes, so a
component with pattern inside $B(o,r)$ takes values in a countable
set of finite planar components. The annealed approximation bound
\cref{eq:sos-annealed-locality} is the hypothesis
$\delta_A(r)\to0$ of \cref{lem:joint-limit}, and its covariance
counterpart is the conditional-covariance hypothesis for mixing;
the finite-volume coupling is stationary by
\cref{prop:physical-sos-law}. \Cref{lem:joint-limit} therefore
proves the unique stationary mixing joint limit, independent of
the windings because the tiling limit is. Exhausting these
kernels, that is, letting $r\to\infty$ in the limiting law, shows
that every limiting component is finite and planar: for each $r$ the
limiting law of $\scR_A$ is within $\delta_A(r)$ of a law
supported on such components. Finally the horizontal-face event is a
cylinder event with limiting probability $\rho_3>0$ by
\cref{prop:physical-sos-law}, so conditioning on it preserves
local convergence, exactly as in \cref{cor:physical-palm}.
\end{proof}

\subsubsection{The physical SOS local limit and GFF transfer}
\label{sec:sos-physical-limit}

All three surfaces in this subsection use the same lift and
vertical origin. Write
$R^{(1)}_z=\cH(z)-\phi(z)$ and $R^{(2)}_z=\cS(z)-\cH(z)$
for the two displacements, so that $\cS-\phi=R^{(1)}+R^{(2)}$.
The plan is the one of \cref{sec:physical-local-limit}: locate
each displacement in a single component of the limiting history,
so that the physical heights can be decoded locally from the
diagonal components (\cref{prop:sos-physical-joint}); deduce
stationarity, ergodicity and the sorting identities for the
limiting triple; and then remove the stationary $L^2$ field
$R^{(1)}+R^{(2)}$ by \cref{lem:ergodic-correction} after
changing coordinates by \cref{prop:reference-coordinate-change}
(\cref{thm:sos-gff-complete}). We first locate both displacements
in the complete history. The point of the next lemma is that,
although $R^{(1)}_z$ is created by the outer stage and
$R^{(2)}_z$ by the inner stage, the inner defect responsible for
$R^{(2)}_z$ sits on a face that is either a reference face or a
face created by the very outer defect responsible for
$R^{(1)}_z$, so the two are never carried by different
components.

\begin{lemma}[Both physical displacements are carried by one component]
\label{lem:sos-column-carrier}
Let $q_z$ be the horizontal face of $\phi$ over $z$.
If either $R^{(1)}_z$ or $R^{(2)}_z$ is nonzero, all
spin changes on that vertical column, in both stages,
belong to one base history component, hence to one component
of the gas of \cref{prop:sos-mixed-envelopes}, whose enlarged
pattern $P$ contains $q_z$.
In a finite planar history, all changes of this
component in $\cH$ and $\cS$ lie in the coordinate bounding box
of the lift of $P$, enlarged by a fixed distance. In particular,
\begin{equation}\label{eq:sos-column-confinement}
                 |R^{(1)}_z|+|R^{(2)}_z|\leq C|P|\,.
\end{equation}
The periodic bound holds for $\diam P\leq\eta_* n$.
Here $\eta_*>0$ can be smaller than the constant in
\cref{lem:sos-small-locality}.
\end{lemma}
\begin{proof}
The common-carrier assertion, including periodic histories, is
\cref{lem:sos-column-component}.

Delete all other components. The retained history
is admissible and has exactly the same changes on the column,
by \cref{lem:sos-carried-components}. The outer confinement
of \cref{lem:sos-small-locality} bounds its intermediate
replacement patch $P^{\cH}$ by the coordinate box of $P$, with
a fixed enlargement. For an inner defect $D_\gamma$, fix
a changed transverse square $(x_1,x_3)$ and move in both
$x_3$ directions to the first squares outside its consecutive
run. These boundary squares belong to its neighborhood and to no
other inner support. Their unchanged normal-$x_2$ dots are
in $Q_\gamma\subset P^{\cH}$. Vertical monotonicity sandwiches
the entire changed binary $x_2$-fiber between the two
unchanged step functions. Hence all its changed cubes have
$x_2$ between these two dot heights. Their other two
coordinates already belong to the transverse projection of
$P^{\cH}$. This proves the corresponding confinement for all
inner changes and new faces of $\cS$.

The lift of the monotone reference is uniformly Lipschitz
in diagonal coordinates, so this box has diameter
$C(1+\diam P)\leq C'|P|$. Each nontrivial vertical
displacement has its endpoint faces in this box, including
$q_z$, which proves \cref{eq:sos-column-confinement}.
For small periodic patterns the supports of both stages and
their neighborhoods have embedded planar lifts by
\cref{lem:sos-small-locality}, so this same proof applies.
\end{proof}

\begin{proposition}[Joint physical \SOS convergence and displacement tails]
\label{prop:sos-physical-joint}
For each positive limiting slope and sufficiently large $\beta$,
the physical triples $(\cS,\cH,\phi)$, rooted by $\phi(0)=0$,
have a unique local limit, independent of the winding sequence.
For $1\leq r\leq c_*n$, their finite-volume displacements obey
\begin{equation}\label{eq:sos-displacement-tail}
 \P\big(|R^{(1)}_0|+|R^{(2)}_0|\geq r\big)
       \leq Ce^{-c\beta^{2/5}r}+Ce^{-c\max(B,r)}\,.
\end{equation}
The limit satisfies this bound for every integer $r\geq1$,
and for each $p>0$,
\begin{equation}\label{eq:sos-displacement-moments}
 \E\big(|R^{(1)}_0|+|R^{(2)}_0|\big)^p
                                    \leq C_p e^{-c'\beta^{2/5}}\,.
\end{equation}
\end{proposition}
\begin{proof}
We first prove the finite-volume tail bound. On the event in
\cref{eq:sos-displacement-tail}, let $P$ be the enlarged
pattern of the component at $q_0$ given by \cref{lem:sos-column-carrier}. If
$P$ is small, in the sense $\diam P\leq\eta_*n$ of that lemma,
then \cref{eq:sos-column-confinement} gives
$r\leq|R^{(1)}_0|+|R^{(2)}_0|\leq C|P|$, that is, $|P|\geq cr$.
If $P$ is not small, then, since a connected pattern of $|P|$
tiles has diameter at most $C|P|$, we get $|P|\geq c n$, which
also implies $|P|\geq c'r$ in the range $r\leq c_*n$. In either
case the event requires a component through $q_0$ of size at
least $c'r$. By \cref{prop:physical-sos-law}, the common physical
root is obtained by conditioning the diagonal joint law on a
horizontal face at the root, an event of probability $n^2/N_0$,
which is bounded below on the compact slope range; hence the
conditional probability of any event is at most $N_0/n^2$ times
its diagonal probability. The component-size bound
\cref{eq:sos-component-tail}, applied at either triangle of the
reference face $q_0$ with $s=c'r$, therefore proves
\cref{eq:sos-displacement-tail} after changing the fixed
constants. Physical stationarity, from
\cref{prop:physical-sos-law}, gives the same bound at any fixed
base square. This conditioning step accounts for the random
location of the reference face in the diagonal law.

For local convergence, use \cref{lem:physical-reconstruction} with
\cref{cor:sos-joint-diagonal,prop:physical-sos-law}. Its only new
hypothesis to check is the decoding. At each physical square $z$,
read the history component whose enlarged pattern contains the
reference face $q_z$.
The signed $\{0,1\}$-spin changes in its two stages give
\[
 \begin{split}
 \cH(z)-\phi(z)&=\sum_{\alpha\text{ in the component at }q_z}
                       \sum_{k\in\Z}\Delta_\alpha(z,k),\\
 \cS(z)-\cH(z)&=\sum_{\gamma\text{ in the component at }q_z}
                       \sum_{k\in\Z}\Delta_\gamma(z,k).
 \end{split}
\]
Here $\Delta_\alpha,\Delta_\gamma\in\{-1,0,1\}$ are the $\{0,1\}$-spin
changes relative to $\phi$ and the specified intermediate $\cH$;
an absent component gives zero. By \cref{lem:sos-column-carrier}, no
remote component contributes. Every limiting component is finite by
\cref{cor:sos-joint-diagonal}, and on each specified patch and
list of components these finite sums agree with the finite-volume
construction. The reconstruction lemma gives the unique physical
limit, its stationarity and ergodicity, including the marginal
rooted by $\cS(0)=0$.

For the moment bound, fix an integer $r\geq1$; eventually
$r\leq c_*n$ along the sequence, so
\cref{eq:sos-displacement-tail} holds at this $r$ for all large
$n$. The displacements are integer-valued, so the event
$\{|R^{(1)}_0|+|R^{(2)}_0|\geq r\}$ is a cylinder event and the
finite-volume tail bound passes to the limiting triple at that
$r$. Using the discrete tail identity
$\E X^p=\sum_{r\geq1}(r^p-(r-1)^p)\P(X\geq r)$, the first term of
the tail bound contributes $C_pe^{-c\beta^{2/5}}$, while the
second term equals $Ce^{-cB}$ for $r\leq B$, contributing at most
$CB^pe^{-cB}$, and is $Ce^{-cr}$ for $r>B$, contributing at most
$C_pe^{-cB}$. Hence
\[
 \E\big(|R^{(1)}_0|+|R^{(2)}_0|\big)^p
   \leq C_p e^{-c\beta^{2/5}}
       +C_p(1+B^p)e^{-cB}\leq C_p'e^{-c'\beta^{2/5}}\,,
\]
since $B=\sqrt\beta\geq\beta^{2/5}$. The summation is for the
limit law; no uniform moment estimate above the finite-volume
threshold $c_*n$ is needed.
\end{proof}

\begin{proof}[Proof of \cref{thm:sos-gff-complete}]
The limiting triple is the one constructed in
\cref{prop:sos-physical-joint}, whose decoding verifies
\cref{lem:physical-reconstruction}. Thus the limiting triple is
stationary and ergodic. Its displacement tails and moments are
\cref{eq:sos-displacement-tail,eq:sos-displacement-moments}, so the
same lemma gives logarithmic displacement growth, the common
limiting plane and all three mean slopes. It remains to verify
the sorting identities and apply the field-transfer criterion.
\begin{enumerate}[label=\textbf{Step~\arabic*.}, ref=\arabic*, wide=0pt, itemsep=1ex, listparindent=\parindent, parsep=0pt]
\item\label[step]{st:sos-sorting-identities} \emph{(The sorting identities.)}
The common limiting plane and $\theta_1,\theta_2>0$ give the
occupied and vacant ends in both sorting directions. On an $x_1$
sorting column, every outer change lies in the finite zero-charge
defect at its reference dot, as in the proof of
\cref{thm:hmsos-physical-local-limit}. On an $x_2$ sorting column,
the dot of $\cH$ is either common with $\phi$ or carried by the
footprint $P_\alpha$ of the outer defect $\alpha$ that created it, by
\cref{lem:sos-column-carrier}. Any inner defect changing that
column contains this dot in its footprint, so its complete
zero-charge change lies in the same history component. Work on the
countable intersection of the events that every component is finite
and admissible and the decoding holds at every cube. Both sorting
columns then retain their charges simultaneously, proving
$\Tone\cH=\phi$ and $\Ttwo\cS=\cH$.

\item\label[step]{st:sos-gff-transfer} \emph{(The \GFF transfer.)}
The diagonal \GFF is \cref{prop:sos-diagonal-gff}, and
\cref{lem:diagonal-height-second-moment} supplies the reference
height moment. The physical coupling is stationary and ergodic,
with the stated mean slopes, and $\cS-\phi=R^{(1)}+R^{(2)}$ is in
$L^2$ by \cref{eq:sos-displacement-moments}. Thus
\cref{cor:physical-height-transfer} gives the scaling limit
asserted in \cref{thm:sos-gff-complete}. Uniqueness
and independence of the windings are those of
\cref{prop:sos-physical-joint}.
\end{enumerate}
\end{proof}

\subsection{Three steps: Ising Interface reconstruction}

\subsubsection{Physical rooting for the Ising Interface model}
\label{sec:interface-rooting}

The three-step reconstruction follows the two-step one, with two
additions. First, the history now has a third stage, the
bubble-free interface $\cI$ above its \SOS surface $\cS=\Tthree\cI$,
whose defects are entire binary columns rather than heights; the
height observable of \cref{thm:gff-limits} is $\cS$ itself, so
the third stage enters the displacement estimates only through the
positivity and locality of its conditional law. Second, we carry
throughout an admissible contact perturbation
$U=U_{\beta,\Gamma}$ in the sense of
\cref{def:interface-perturbation}. This is what makes the result
applicable to the full Ising model in \cref{sec:ising}, where $U$
is the effective interaction produced by integrating out the
closed contours; the bubble-free Ising Interface model itself is the
case $U=0$. The perturbation has two consequences for the
argument. Its weight is not a product over history components, so
the deletion estimates must be rederived with the variation bound
\cref{eq:interface-contact-variation}; and the pattern expansion
used as \RG input has signed coefficients, so it cannot serve as a
probability space for the switching argument. The latter issue is
resolved in \cref{sec:interface-perturbed-conditional-reconstruction}
by a positive marking of the contact terms. Finally, the
interface model is defined with a vertical cutoff $M$, and we
show that the cutoff can be removed before the limit $n\to\infty$
provided $M_n/n^2\to\infty$.

\begin{proposition}[Rooted interface law and cutoff removal]
\label{prop:physical-interface-law}
For sufficiently large $\beta$, the bubble-free interface
law with weights $e^{-\beta(\sA(\cI)-N_0)+U_{\beta,\Gamma}(\cI)}$,
for any admissible contact perturbation, has a finite partition function
when rooted by $(\Tthree\cI)(0)=0$.
Writing $\cZ^{\mathrm{I},\mathrm{phys}}_{\beta,\Gamma}$ for its
partition function with excess energy $\sA(\cI)-N_0$,
\begin{equation}\label{eq:physical-interface-partition}
 \cZ^{\mathrm{I},\mathrm{phys}}_{\beta,\Gamma}
       =\frac{n^2}{N_0}\sum_{M\in\cM_\Gamma}
                                  Z^U_{\beta,\Gamma}(\phi_M)\,.
\end{equation}
The diagonal history law conditioned on a horizontal
reference face gives the physical joint law of $(\cI,\cS,\cH,\phi)$
rooted by $\phi(0)=0$; that event has probability $n^2/N_0$.
The joint law is stationary under physical translations with
this common rerooting. Increasing rooted cutoffs converge in total
variation to the corresponding rooted laws at fixed $n$.
More quantitatively, let $\mu_{n,M}^U$ be the interface law in the
tilted band of half-height $M$, with the same translation-invariant
weight $U_{\beta,\Gamma}$, followed by rooting, and let $\mu_n^U$
be its infinite-height rooted law. Uniformly on compact positive
slope ranges,
\begin{equation}\label{eq:interface-height-cutoff-tv}
 \|\mu_{n,M}^U-\mu_n^U\|_{\mathrm{TV}}
           \leq C\min\{1,N_0/M\},\qquad N_0=n^2+n(\Theta_1+\Theta_2)\,.
\end{equation}
Thus $M_n/n^2\to\infty$ gives the same local limit as removing
$M$ first. This is a sufficient condition, with no upper bound
on $M_n$.
\end{proposition}
\begin{proof}
Initially root by $\cS(0)=0$, where $\cS=\Tthree\cI$.
We first show that the unweighted rooted sum is finite. Group
the interfaces according to $\cS=\Tthree\cI$ and write
$\sA(\cI)-N_0=(\sA(\cI)-\sA(\cS))+(\sA(\cS)-N_0)$. For fixed
$\cS$, the first-fiber bound \cref{eq:interface-first-pressure}
gives $\sum_{\cI:\Tthree\cI=\cS}e^{-t(\sA(\cI)-\sA(\cS))}
\leq e^{\epsilon_t\sA(\cS)}=e^{\epsilon_tN_0}e^{\epsilon_t(\sA(\cS)-N_0)}$,
whence
\[
 \sum_{\cI} e^{-t(\sA(\cI)-N_0)}
 \leq e^{\epsilon_tN_0}
       \sum_{\cS:\cS(0)=0}e^{-(t-\epsilon_t)(\sA(\cS)-N_0)}<\infty\,,
\]
the last sum being the rooted \SOS partition function at inverse
temperature $t-\epsilon_t>0$, which is finite by
\cref{prop:physical-sos-law}. For the perturbed weight use the
first bound in \cref{eq:interface-contact-variation},
$|U(\cI)|\leq\epsilon_\beta\sA(\cI)=\epsilon_\beta(N_0+\sA(\cI)-N_0)$,
to bound it by
$e^{\epsilon_\beta N_0}
 e^{-(\beta-\epsilon_\beta)(\sA(\cI)-N_0)}$, and apply the
unweighted bound at $t=\beta-\epsilon_\beta$.

Define $V(\cI)=V(\Tthree\cI)$. All three sorting steps preserve
this volume, by the fiber-orbit argument in the proof of
\cref{prop:physical-sos-law}. The translation formula
\cref{eq:volume-translation} therefore holds for $\cI$ as well,
and the physical and diagonal orbit representatives have the
same multiplicities $n^2$ and $N_0$ in a common volume strip as
in \cref{prop:physical-sos-law}; the weights are translation
invariant since both the area and $U_{\beta,\Gamma}$ are. The
strip count and the averaged-observable argument prove
\cref{eq:physical-interface-partition}, the sampling
correspondence and the horizontal-face conditioning, and
translation invariance of the weights proves stationarity of the
joint law under translation followed by common rerooting.

Finiteness of the rooted partition function proves
total-variation convergence of exhausting rooted restrictions,
as before. For the quantitative assertion put
$E(\cI)=\sA(\cI)-N_0$ and keep $U=U_{\beta,\Gamma}$ fixed while
varying the coefficient of $E$. Put
$F(t)=\log\sum_{\cI} e^{-tE(\cI)+U(\cI)}$, and let $G_\Gamma$ count the physically
rooted monotone interfaces, which are the interfaces with
$E(\cI)=0$. For the lower bound at $t=\beta$, keep only these
ground states and use $|U(\cI)|\leq\epsilon_\beta N_0$ on them.
For the upper bound at $t=\beta/2$, use
$|U(\cI)|\leq\epsilon_\beta(N_0+E(\cI))$ to compare with the
unperturbed sum at inverse temperature $\beta/2-\epsilon_\beta$,
and bound that sum through the partition identity
\cref{eq:physical-interface-partition} with $U=0$ and the fiber
bound $Z^{\textsc{interface}}_t(\varphi)\leq e^{\epsilon^{\mathrm I}_tN_0}$
of \cref{eq:interface-positive-budget}, noting
$G_\Gamma=(n^2/N_0)|\cM_\Gamma|$. This gives
\[
 G_\Gamma e^{-\epsilon_\beta N_0}\leq e^{F(\beta)},\qquad
 e^{F(\beta/2)}\leq G_\Gamma
       e^{(\epsilon_\beta+\epsilon^{\mathrm I}_{\beta/2-\epsilon_\beta})N_0}\,.
\]
Jensen's inequality, as in \cref{cor:sos-gradient-tightness}, gives
\[
 \frac\beta2\E_{\mu_n^U}E
 \leq\log\E_{\mu_n^U}e^{\beta E/2}
 =F(\beta/2)-F(\beta)
 \leq(2\epsilon_\beta+
              \epsilon^{\mathrm I}_{\beta/2-\epsilon_\beta})N_0.
\]
Since $\sA(\cI)=N_0+E(\cI)$, it follows that
$\E_{\mu_n^U}\sA(\cI)\leq CN_0$.

We now compare the band law with the infinite-height law. For a
rooted interface $\cI$, let $b_M(\cI)$ count its vertical
translates fitting the tilted band of half-height $M$. The band
law weighs each rooted interface by its weight times $b_M(\cI)$,
since every translate fitting the band is a distinct band
configuration with the same weight; thus $\mu_{n,M}^U$ is the
reweighting of $\mu_n^U$ by $r_M(\cI)=b_M(\cI)/(2M)$, normalized.
To estimate $b_M$, note that the periodic tilted height
$\ell(x)=x_3+(\Theta_1/n)x_1+(\Theta_2/n)x_2$ changes by at most a fixed
constant between adjacent plaquettes. Since $\cI$ is edge-connected,
by \cref{lem:phase-tree}, its $\ell$-range is at most $C\sA(\cI)$.
Counting the integer translates between the two cutoff levels,
including the case where none fits, gives
\begin{equation}\label{eq:interface-translate-count}
                  |b_M(\cI)-2M|\leq C(\sA(\cI)+1)\,.
\end{equation}
Together with $\E_{\mu_n^U}\sA(\cI)\leq CN_0$ this implies
$\delta:=\E_{\mu_n^U}|r_M-1|\leq CN_0/M$. Writing $\bar r=\E r_M$,
the total variation distance between $\mu_n^U$ and its
normalized reweighting is
$\frac12\E|r_M/\bar r-1|\leq\frac1{2\bar r}(\E|r_M-1|+|1-\bar r|)
\leq\delta/(1-\delta)$ when $\delta\leq1/2$, hence at most
$2\delta$; otherwise use the trivial bound one. This proves
\cref{eq:interface-height-cutoff-tv}
and the stated cutoff limits. The same estimate holds for the
joint sorting history, since it is a deterministic function
of the rooted interface.
\end{proof}

\begin{proof}[Proof of \cref{thm:interface-rg-input}]
\Cref{prop:physical-interface-law} identifies the rooted
vertical-cutoff limit of the perturbed interface law and its
three-step \Tetris image with the law proportional to
$Z^U_{\beta,\Gamma}(\varphi)$ on $\cM_\Gamma$. The geometric
decomposition of the histories keeps all interface connectivity and
no-bubble restrictions, since these are imposed on the reconstructed
surface $\cI_{M,\cX}$ itself.
\Cref{prop:interface-perturbed-absolute-history,prop:interface-perturbed-positive-completion}
give the a priori collection bound, the positive canonical
completion, the signed \texttt{thin}/\texttt{fat} coefficients with
their absolute bounds, bulk consistency, and the exact source
identities and torus projection; this verifies
\crefrange{it:rg-positive}{it:rg-sources} as in the proofs of
\cref{thm:hmsos-rg-input,thm:sos-rg-input}, the only difference
being that the envelopes $q^\mathrm{t},q^\mathrm{f}$ of
\cref{it:rg-envelopes} are now the sums of absolute values of the
labelled terms rather than the coefficients themselves.
\Cref{lem:interface-marked-input} verifies the marked assertion.
Taking the empty contact family, $U\equiv0$, gives the original
bubble-free interface model.
\end{proof}

\subsubsection{The perturbed interface image and its reconstruction}
\label{sec:interface-rg-application}

\begin{proposition}[Perturbed interface image \GFF]
\label{prop:interface-diagonal-gff}
For every admissible contact perturbation $U$, positive limiting
slopes and sufficiently large $\beta$, the image with weights
$Z^U_{\beta,\Gamma}$ has a unique stationary mixing local limit
of prescribed density $\rho$, independent of the approximating
windings. Its centered diagonal height converges jointly on smooth
integral-zero tests and in local $H^{-s}$ modulo constants, $s>1$, to
\begin{equation}\label{eq:interface-diagonal-gff}
 \upsigma_{\mathrm d}^{U}\,\GFF\circ L_{\mathrm d}^{U},\qquad
 \upsigma_{\mathrm d}^{U}
   =\sqrt{a_U/|\det\sA_{\mathrm D}^{U}|},\qquad
 L_{\mathrm d}^{U}=(\sA_{\mathrm D}^{U})^{-\trps}\,,
\end{equation}
where $a_U>0$ and the matched real dispersion matrix is invertible.
The constants are uniform under the bounds of
\cref{def:interface-perturbation}.
\end{proposition}
\begin{proof}
Apply \cref{thm:coefficient-rg}(i)--(ii) to the perturbed image
identified in \cref{thm:interface-rg-input}. That input theorem
verifies all coefficient conditions uniformly under
\cref{def:interface-perturbation}, including the absolute internal
weights. The resulting dispersion and covariance coefficient give
\cref{eq:interface-diagonal-gff}.
\end{proof}

For reconstruction use the positive geometric groups, witness family
$\scW_n^{U}$ and events $\cB_x^{U}(s)$ of
\cref{lem:interface-marked-input}. Recall from the discussion
preceding that lemma that a geometric history group $\cG$ is
obtained from the expanded components of a primitive history by
enlarging each low-cost component by its six-class witness and
regrouping overlaps; it carries an enlarged pattern $P(\cG)$ on
the final tiling, an excess $E(\cG)$ and its
defects, but no contact marks, and it is labelled \texttt{thin}
or \texttt{fat} by the high-cost-mass test, from this geometry
alone. Recall also that $\cB_x^{U}(s)$ is the event that the
reference tiling contains a member of $\scW_n^{U}$ through $x$
with at least $s$ lozenges.

\begin{proposition}[Interacting witnesses and primitive history tails]
\label{prop:interface-perturbed-interacting-tails}
Let $\nu^{U}_{\beta,\Gamma}$ be the perturbed interface
image law. Uniformly in sufficiently large tori in the prescribed
compact physical slope range,
\begin{equation}\label{eq:interface-perturbed-interacting-witness-tail}
 \nu^{U}_{\beta,\Gamma}(\cB_x^{U}(s))
                              \leq Ce^{-c\max(s,B)}\,.
\end{equation}
In the positive joint law with conditional history weights
$e^{-\beta E(\cX)+U_{\beta,\Gamma}(\cI_{\cX})}$,
the geometric groups above satisfy
\begin{equation}\label{eq:interface-perturbed-primitive-component-tail}
 \P\big(\exists\cG:\ x\in S(P(\cG)),\
                          |P(\cG)|\geq s\big)
       \leq Ce^{-c\beta^{2/5}s}+Ce^{-c\max(s,B)}\,.
\end{equation}
These estimates retain all histories and noncontractible groups.
\end{proposition}
\begin{proof}
\Cref{lem:interface-marked-input} verifies \cref{it:rg-marked}.
Together with the base conditions verified above,
\cref{thm:coefficient-rg}(iii) proves
\cref{eq:interface-perturbed-interacting-witness-tail}.
Positivity is required for the matching measure, as supplied by its
primitive completion; individual expansion coefficients may be signed.

For the group tail we work directly in the positive primitive
history law, whose weights are
$e^{-\beta E(\cX)+U(\cI_{\cX})}$ on compatible families $\cX$ of
geometric groups. The contact term couples different groups, so
the law is not a product gas and \cref{lem:deletion} does not
apply as stated; we therefore repeat the deletion argument with
the variation bound \cref{eq:interface-contact-variation} in place
of exact factorization.
Fix the tiling and a specified geometric group $\cG$
with footprint $P$ and excess $E$. If $\cX$ is a
compatible remaining family, deleting $\cG$ preserves
all its states, and the interface $\cI_{\cX}$ differs from
$\cI_{\cX\cup\cG}$ exactly on the faces changed by $\cG$, because
distinct expanded components change disjoint sets of faces. Every
reference face removed by $\cG$ lies in its footprint, and the
faces it adds number at most the removed ones plus its excess, so
\[
 |\cI_{\cX\cup\cG}\xor \cI_{\cX}|
       =|\cI_{\cG}\xor\varphi|\leq 2|P|+E\,.
\]
Using the second bound in \cref{eq:interface-contact-variation},
with $\epsilon_\beta\leq Ce^{-c\beta}$, the contact terms of the
two histories differ by at most $\epsilon_\beta(2|P|+E)$, which
gives the uniform ratio
\begin{equation}\label{eq:interface-perturbed-positive-deletion-ratio}
 \frac{e^{-\beta(E(\cX)+E)+U(\cI_{\cX\cup\cG})}}
      {e^{-\beta E(\cX)+U(\cI_{\cX})}}
       \leq e^{2\epsilon_\beta|P|}
                    e^{-(\beta-\epsilon_\beta)E}\,.
\end{equation}
The map deleting this specified group is injective from the
histories containing $\cG$ to the histories not containing it.
Summing \cref{eq:interface-perturbed-positive-deletion-ratio}
over the image, which is a subset of all remaining histories, and
dividing by the full positive conditional partition function
therefore bounds the conditional probability that $\cG$ is
present by the right side of
\cref{eq:interface-perturbed-positive-deletion-ratio}. Both sides
retain the perturbation factor; deletion is performed in the
positive history law, not in the signed expansion.

For a \texttt{thin} group of size $m=|P|$, its high-cost
components have total footprint size at least $d_\beta m$, by the
definition of the \texttt{thin} class, and a high-cost component has
excess at least its footprint size divided by $B$. Thus
$E\geq d_\beta m/B$. Now sum the preceding bound over all
groups $\cG$ with the same pattern $P$: extract half of
the energy factor using this lower bound, and bound the remaining
sum $\sum_{\cG:P(\cG)=P}e^{-(\beta-\epsilon_\beta)E(\cG)/2}$ by
the a priori bound for unweighted histories,
\cref{eq:interface-positive-budget} at inverse temperature
$(\beta-\epsilon_\beta)/2$ and $W=P$, keeping the geometric
classification fixed at its value for $\beta$. For large $\beta$
this gives
\[
 e^{2\epsilon_\beta m}
 \sum_{\substack{\cG\text{ \texttt{thin}}\\P(\cG)=P}}
             e^{-(\beta-\epsilon_\beta)E(\cG)}
 \leq
 e^{2\epsilon_\beta m
     -(\beta-\epsilon_\beta)d_\beta m/(2B)
     +\epsilon^{\mathrm I}_{(\beta-\epsilon_\beta)/2}m}
 \leq e^{-c\beta^{2/5}m}\,,
\]
since $\beta d_\beta/B\asymp\beta^{2/5}$ while the other two
exponents are $o(1)$.
There are at most $C_0^m$ connected labelled lozenge patterns
of size $m$ through $x$, so summing over $m\geq s$ shows that the
\texttt{thin} probabilities are bounded by the first term
of~\cref{eq:interface-perturbed-primitive-component-tail}.
A \texttt{fat} geometric group has size at least $B$ and its pattern
belongs to $\scW_n^{U}$. Its presence therefore
implies the corresponding witness event, regardless of its
conditional activity, which may be large. Applying
\cref{eq:interface-perturbed-interacting-witness-tail}
proves the second term and completes the proof.
\end{proof}

\subsubsection{Positive contact marks and conditional reconstruction}
\label{sec:interface-perturbed-conditional-reconstruction}

The switching argument of \cref{prop:sos-conditional-locality}
requires a positive product-form weight on the objects being
swapped: it compares two independent samples and exchanges
components between them, and the exchange must preserve the
weight exactly. The signed expansion used for the \RG input is not
a history probability space, since its contact marks carry
inclusion--exclusion coefficients of either sign; and the positive
history law $e^{-\beta E(\cX)+U(\cI_{\cX})}$, while positive, is
not of product form, because $U$ couples all groups whose changed
faces meet a common carrier. We therefore give a different
expansion, solely for coupling the positive conditional
histories: each contact term $e^{a_\beta(X)h_X}$ is written as a
constant times $1+w_X\chi_X$ with $w_X\geq0$ and $\chi_X$ an
indicator, and the resulting product is expanded into a positive
law on histories decorated by independent marks. The marks are
then treated as additional polymers in the switching, joined to
the groups whose changes they see. Throughout this subsection the
geometric groups are those preceding
\cref{prop:interface-perturbed-interacting-tails}; in particular they
are formed before adding any contact marks.

\begin{lemma}[Positive marking of the contact interaction]
\label{lem:interface-perturbed-positive-marks}
For finitely many terms in
\cref{eq:interface-contact-functional}, set
\[
 w_X=e^{|a_\beta(X)|}-1,\qquad
 \chi_X(\cI)=\begin{cases}h_X(\cI),&a_\beta(X)\geq0,\\
                       1-h_X(\cI),&a_\beta(X)<0.
          \end{cases}\,.
\]
The positive history law is the marginal of the measure with
weight
\begin{equation}\label{eq:interface-perturbed-positive-mark-weight}
 e^{-\beta E(\cX)}
                  \prod_{X\in\cM}w_X\chi_X(\cI_{\cX})\,,
\end{equation}
where each contact label can occur at most once in
$\cM$. Conditional on the history, its marks are
independent, with probabilities
$(1-e^{-|a_\beta(X)|})\chi_X(\cI_{\cX})$.
The construction extends to all contact labels. Only finitely many
marks meet the changed faces of any fixed finite history,
almost surely. For some $c_0>0$,
\begin{equation}\label{eq:interface-perturbed-positive-mark-budget}
 \sup_p\sum_{X:p\in N_X}w_Xe^{c_0\beta m_X}
                                  \leq Ce^{-c\beta}\,.
\end{equation}
Repeated carriers with distinct labels remain distinct marks.
\end{lemma}
\begin{proof}
For $h\in\{0,1\}$ and real $a$, with $\chi(h)=h$ when $a\geq0$
and $\chi(h)=1-h$ when $a<0$, we have the exact identity
\[
 e^{a h}=e^{\min(a,0)}
       \big[1+(e^{|a|}-1)\chi(h)\big]\,.
\]
Indeed, for $a\geq0$ both sides equal $1$ at $h=0$ and $e^a$ at
$h=1$; for $a<0$ both sides equal $e^a$ at $h=1$ and
$e^a\cdot e^{-a}=1$ at $h=0$. Taking the product over a finite
list of contact terms and expanding each bracket, the weight
$e^{-\beta E(\cX)+U(\cI_{\cX})}$ becomes a constant times the sum,
over subsets $\cM$ of the list, of
\cref{eq:interface-perturbed-positive-mark-weight}; this proves
the finite-list assertion. The prefactor
$\prod_Xe^{\min(a_\beta(X),0)}$ is independent of the interface
and cancels from the probability law. Conditional on the history,
the weight of $\cM$ is a product over labels of
$w_X\chi_X(\cI_{\cX})$ or $1$, so the marks are independent, and
label $X$ is marked with probability
$w_X\chi_X/(1+w_X\chi_X)$, which equals $1-e^{-|a_\beta(X)|}$ when
$\chi_X(\cI_\cX)=1$ and zero otherwise, as claimed. The
anchored bound \cref{eq:interface-contact-norm} implies
$|a_\beta(X)|\leq Ce^{-c\beta}$ for each nonzero term;
thus $w_X\leq2|a_\beta(X)|$, and the same bound, with the factor
$e^{-3\beta m_X/4}$ in \cref{eq:interface-contact-norm} compared
to the factor $e^{-\beta m_X}$ in $a_\beta(X)$, gives
\cref{eq:interface-perturbed-positive-mark-budget} with the
reserved mass exponent $c_0\beta m_X$ for any $c_0<1/4$.
Infinitely many marks far from the interface may be present. They
are not projected into a planar Poisson process: only marks
meeting a changed face will enter the coupling below.

For the exhaustion, let $U_J$ retain the first $J$ contact labels
and subtract $U_J(\varphi)$ before normalizing the history weights.
For histories whose group patterns lie in a fixed finite set $W$ of
reference plaquettes, the contact bound and the area identity give
\cref{eq:contact-variation-history} at $t=\beta$ for
$|U_J(\cI_{\cX})-U_J(\varphi)|$, uniformly in $J$.
The resulting weights are bounded by
$e^{2\epsilon_\beta|W|}e^{-(\beta-\epsilon_\beta)E(\cX)}$.
This envelope is summable over histories by
\cref{eq:interface-positive-budget}, applied at
$t=\beta-\epsilon_\beta$. Restricting the expanded footprints to
$W$ also restricts the old footprints to $W$, so that bound applies.
Each relative potential converges absolutely, so dominated convergence
gives convergence of both the partition functions and the normalized
history laws in total variation. On a finite quotient one takes
$W=\varphi$; the same argument applies in each finite planar window.
Given a limiting history, use the countable product of the stated
Bernoulli mark laws. The sum of the probabilities of marks meeting
its changed faces is finite by
\cref{eq:interface-perturbed-positive-mark-budget}, proving their
almost sure finiteness.

The same construction is the limit of the finite-list markings
for all marks used in the coupling below. To see this, write
$D(\cX)=\cI_{\cX}\xor\varphi$. Given $\cX$, the probability
of a mark outside the first $J$ labels meeting $D(\cX)$ is at most
\[
 r_J(\cX):=\min\bigg\{1,
       \sum_{\substack{X\text{ outside the first }J\text{ labels}\\
                         N_X\cap D(\cX)\ne\varnothing}}w_X\bigg\}.
\]
The anchored bound makes this sum finite and its tail tend to
zero for each finite history. Couple the history marginals using
their total-variation convergence and, when they agree, use the
same Bernoulli variables for the first $J$ labels. The remaining
disagreement probability is at most the expectation of $r_J$
under the limiting history law, which tends to zero by bounded
convergence. Thus the history together with all marks touching its
changed faces converges in total variation. For the two-copy
coupling, use $D(\cX_1)\cup D(\cX_2)$ in the same bound for
each copy; this also retains marks touching only the other copy's
changed faces. No convergence in total variation is claimed for
the marks arbitrarily far from these faces.
This constructs the required marking without taking an infinite
product of the interface-independent prefactors.
\end{proof}

Let $W=B(o,r)$ in the diagonal tiling and $A\subset B(o,r/2)$.
Restrict geometric history groups to patterns contained in $W$.
Keep only contact terms with $A_X^\varphi\subset W$, and write
their sum as $U^W$. Define $K^{U}_{A,r}(\varphi)$ by
the weights
\begin{equation}\label{eq:interface-perturbed-restricted-kernel}
 e^{-\beta E(\cX)+U^W(\cI_{\cX})-U^W(\varphi)}\,,
\end{equation}
recording the complete geometric groups meeting $A$. This
normalization of the exponent does not change the kernel.
The bound
$|U^W(\cI_{\cX})-U^W(\varphi)|
 \leq\epsilon_\beta(2|W|+E)$
proves summability. Its groups form a countable space; their
transition heights are not truncated. Anchored lifting gives
the same kernel for every reference tiling agreeing on
$B(o,r+C)$, and exact agreement with its planar version
when this ball embeds. Contact terms may extend arbitrarily
far along the diagonal line, but their projection stays in
$W$ and their lift is anchored at a reference or changed face.

The constants fixed in the definition of $\scW_n^{U}$, before
\cref{lem:interface-marked-input}, already admit unions of two
six-class families with remainder fraction
$C\delta^\mathrm{g}_\beta+CL/B$, as in
\cref{prop:sos-conditional-locality}; its vanishing entropy rate
and its interacting tail were proved for these constants.

\begin{proposition}[Conditional locality for perturbed interface histories]
\label{prop:interface-perturbed-conditional-locality}
For $r\leq c_*n$, the complete perturbed interface conditional
history satisfies
\begin{align}\label{eq:interface-perturbed-conditional-locality}
 \|\cL(\scR_A\given\varphi)
                     -K^{U}_{A,r}(\varphi)\|_\mathrm{TV}
 &\leq C|A|e^{-c\beta^{3/10}r}
       +C\sum_{x\in A}\one_{\cB_x^{U}(cRr)}\,.
\end{align}
Conditional covariances of bounded observables of the groups on sets
$A,A'$ at distance at least $r$ have the same bound, multiplied
by $C\|F\|_\infty\|G\|_\infty$. Both bounds have annealed
error at most
\begin{equation}\label{eq:interface-perturbed-annealed-locality}
                C|A|\big(e^{-c\beta^{3/10}r}+e^{-cRr}\big)\,.
\end{equation}
\end{proposition}
\begin{proof}
The proof is that of \cref{prop:sos-conditional-locality}, with
the marked histories of
\cref{lem:interface-perturbed-positive-marks} in place of the
gas of enlarged patterns. We first explain the switching in the marked
representation, then bound the exceptional event by the same two
alternatives as before, with the mark masses added to the
\texttt{thin} mass, and finally remove the truncations.
\begin{enumerate}[label=\textbf{Step~\arabic*.}, ref=\arabic*, wide=0pt, itemsep=1ex, listparindent=\parindent, parsep=0pt]
\item\label[step]{st:interface-marked-switching} \emph{(Switching in the marked representation.)}
First use finite label lists and a finite common list of
allowed geometric group states; restrict compatible families
from that list, rather than imposing a cutoff on their total
energy. Such a truncation is preserved by color switching, since
the switching only exchanges allowed objects between the copies.
Take two independent positively marked histories, one full
and one restricted. Form a graph on the groups and mark
occurrences of both copies: join incompatible geometric groups,
and join a mark occurrence $X$ to every group in either copy
whose changed plaquette set meets $N_X$. There are no
mark--mark edges. Marks having no such neighbor can be
ignored, since their predicate is evaluated on the reference and
is the same in both copies. Each connected component has a
color-swapping involution: interchange all its groups and marks
between the copies. Group compatibility is preserved because all
incompatibilities were joined into one component. For a swapped
mark, every history change affecting $h_X$ is swapped with it,
because every group whose changes meet $N_X$ is its neighbor. Its
predicate $\chi_X$ therefore has exactly its former value in the
other copy. An unswapped mark has no changed predicate, since
none of its neighbors was swapped. The reference is the same in
both copies. Thus the product weight
\cref{eq:interface-perturbed-positive-mark-weight} is preserved
exactly by the swap: the energies add over groups, and each mark
contributes the same factor $w_X\chi_X$ as before. In other words,
the switching is performed in a positive extension of the
perturbed law, with all predicates preserved, which is what the
signed expansion could not provide.

\item\label[step]{st:interface-exceptional-event} \emph{(The exceptional event and its enumeration.)}
Swap components meeting $A$, unless one contains a group or
mark forbidden in the restricted copy. Off this exceptional event
the swap exchanges $\scR_A$ between the copies, so the total
variation distance is at most twice the probability of the
exceptional event, as before. On this exceptional
event a connected union $U$ of group patterns and mark
footprints contains $x\in A$ and has diameter at least $r/2$.
A mark footprint is connected, has at most $C_1m_X$ tiles,
and meets every neighboring group footprint, by expanded
projection control. Write
\[
 m=|U|,\quad q=\sum_{\text{\texttt{thin} groups in both copies}}|P|,
 \quad M=\sum_{\text{mark occurrences}}m_X,\quad
 Q=q+C_1M\,.
\]
Within each copy the group patterns are disjoint, so their
total size is at most $m$ and the two-copy total at most $2m$.
Marks need not have disjoint footprints.

We now give the probability and enumeration bounds, including
their sums over histories; the quantity $Q$ plays the role of the
\texttt{thin} mass $q$ in \cref{prop:sos-conditional-locality},
with the marks counted as \texttt{thin} material because their
weights are exponentially small in their masses. For specified
families of groups, simultaneous deletion in each positive history
law, by the argument of
\cref{eq:interface-perturbed-positive-deletion-ratio} applied
successively to each group, gives
\[
 \P(\text{specified groups})
 \leq\prod_{i=1}^2
       \exp(2\epsilon_\beta\sum_{\cG\in\cX_i}
                           |P(\cG)|)
       e^{-(\beta-\epsilon_\beta)E(\cX_i)}\,.
\]
The same estimate holds for the restricted potential $U^W$.
Given the histories, specified marks have probability at
most the product of their $w_X$. For
$\boldsymbol{\cX}=(\cX_1,\cX_2)$ and
$\boldsymbol{\cM}=(\cM_1,\cM_2)$ define
\[
 \fa_\beta(\boldsymbol{\cX},\boldsymbol{\cM})
 =\exp\bigg(2\epsilon_\beta\sum_{i,\cG}|P(\cG)|
                 -(\beta-\epsilon_\beta)\sum_iE(\cX_i)\bigg)
                       \prod_{i=1}^2\prod_{X\in\cM_i}w_X\,.
\]
This bounds the probability of all specified objects. In
sums below each mark has at least one neighbor among the
specified groups. Summing all possible marks
attached to the selected groups, in either copy, with a
mass enhancement $e^{c_0\beta M/2}$ gives at most
\[
 \exp(C\epsilon_\beta
            (\sum_{i,\cG}|P(\cG)|+
                                      \sum_iE(\cX_i)))\,.
\]
Indeed every such mark has a changed face as anchor;
\cref{eq:interface-perturbed-positive-mark-budget} and
$\sum|\cI_{\cG}\xor\varphi|
 \leq2\sum|P|+\sum E$ bound the product of $1+w_X$.
Dropping the connection and union constraints only increases
this sum. Marks are summed in three dimensions at these
anchors, not over all their possible heights at a planar tile.

Absorb the last energy factor by replacing $\beta$ with
$\beta-C\epsilon_\beta$, keeping the geometric labels
fixed at their original value. The proof of the \texttt{thin}
envelope, with half the energy reserved, still gives
$e^{-A_\beta|P|}$, $A_\beta=c\beta^{2/5}$; the \texttt{fat}
envelope remains $e^{o(1)|P|}$ with rooted entropy $o(1)$.
The factors $e^{C\epsilon_\beta|P|}$ are absorbed there.
Repeating the explicit product-of-$1+$activity estimate
\cref{eq:sos-two-copy-budget} therefore gives, for fixed $U$,
\begin{equation}\label{eq:interface-perturbed-marked-two-copy-budget}
 \sum_{\substack{\boldsymbol{\cX},\boldsymbol{\cM}:\\
                         \text{footprint union }U}}
 e^{A_\beta Q/4}
 \fa_\beta(\boldsymbol{\cX},\boldsymbol{\cM})
                                      \leq e^{o_\beta(1)m}\,.
\end{equation}
To check the mark enhancement, $A_\beta C_1/4<c_0\beta/2$
for large $\beta$. For the group enhancement, insert
$e^{-\zeta_\beta|P|}$ as in that estimate; the total
compensating exponent is at most $2\zeta_\beta m$.
This proves the displayed inequality without bounding the
number of components or mark occurrences by fiat.

\item\label[step]{st:interface-two-alternatives} \emph{(The connection estimate.)}
Apply \cref{lem:connecting-union-budget} with the envelopes
$\fa_\beta$, cost $Q=q+C_1M$, $a=A_\beta/4$, $b=o_\beta(1)$ and
$\delta=\delta^\mathrm{g}_\beta$. The enhanced a priori bound is
\cref{eq:interface-perturbed-marked-two-copy-budget} and
$\lambda\asymp\beta^{3/10}$. To check the witness alternative,
if $Q<\delta m$ then the tiles outside the \texttt{fat} group
patterns are at most $Q$: the \texttt{thin} groups contribute
at most $q$ and the marks at most $C_1M$. There must be a
\texttt{fat} group, since otherwise $m\leq Q$. The six-class
verification in \cref{prop:sos-conditional-locality} therefore
applies, with the marks included in the remainder. It makes $U$
a member of $\scW_n^U$ with $m\geq cR\diam U\geq cRr$, so
$\cB_x^U(cRr)$ holds. The lemma proves the required exceptional-event
bound and hence \cref{eq:interface-perturbed-conditional-locality}.

\item\label[step]{st:interface-covariance-exhaustion} \emph{(Covariances, removal of the truncations and the annealed bound.)}
For the covariance use two full copies and swap all components
meeting $A$ unless one also meets $A'$. Off that connection
event $(F_1-F_2)(G_1-G_2)$ changes sign with unchanged weight,
so only the connection event contributes to
$\cov(F,G)=\frac12\E[(F_1-F_2)(G_1-G_2)]$, and the same
enumeration bounds this event, whose union has diameter
at least $r$. To remove the truncations, exhaust the finite label
lists and the defect truncations: by
\cref{lem:interface-perturbed-positive-marks}, the marked
histories converge in total variation together with all marks
touching their changed faces, using the a priori bounds and the
anchored mark bound, so every local coupling event and its bound
survive this exhaustion. Finally, averaging over the image law and
using \cref{eq:interface-perturbed-interacting-witness-tail}
with $s=cRr$ gives the annealed estimate
\cref{eq:interface-perturbed-annealed-locality}. Unions were
counted on the quotient, so no noncontractible component was omitted.
\end{enumerate}
\end{proof}

\begin{corollary}[The complete perturbed interface history limit]
\label{cor:interface-perturbed-joint-diagonal}
The perturbed interface tiling and its geometric history
groups converge jointly to a unique stationary mixing
law, independently of the approximating windings. Every limiting
group is finite and planar. Horizontal-face conditioning
has the same convergence.
\end{corollary}
\begin{proof}
This is the analogue of \cref{cor:sos-joint-diagonal}. For fixed
$r$, the kernel defined by
\cref{eq:interface-perturbed-restricted-kernel} is a probability
law on countable finite planar groups, determined by a finite
tiling patch, and agrees with its planar version for large
periods, by the anchored lifting discussed after
\cref{lem:interface-perturbed-positive-marks}.
\Cref{prop:interface-diagonal-gff} supplies the mixing tiling
limit, and the approximation
\cref{eq:interface-perturbed-annealed-locality} together with its
covariance counterpart are the two hypotheses of
\cref{lem:joint-limit}, with the finite-volume coupling
stationary by \cref{prop:physical-interface-law}. This proves the
unique stationary mixing limit. Since the groups are countable
data including the third-stage transition heights, the tightness
supplied by \cref{lem:joint-limit} includes tightness of these
heights, which is not assumed elsewhere. Exhausting the restricted
kernels, as $r\to\infty$, gives finiteness and planarity of all
limiting groups. Conditioning preserves convergence since the
horizontal-face event is a cylinder event with probability tending
to $\rho_3>0$ by \cref{prop:physical-interface-law}.
\end{proof}

\begin{proof}[Proof of \cref{thm:interface-gff-complete}]
Apply \cref{lem:physical-reconstruction} to the joint diagonal
limit of \cref{cor:interface-perturbed-joint-diagonal}, using the
physical rooting and stationarity of
\cref{prop:physical-interface-law}. We check the decoding and
displacement hypotheses, then verify the sorting and no-bubble
properties and transfer the field limit.
\begin{enumerate}[label=\textbf{Step~\arabic*.}, ref=\arabic*, wide=0pt, itemsep=1ex, listparindent=\parindent, parsep=0pt]
\item\label[step]{st:interface-decoding} \emph{(Local decoding of the spin interface.)}
We specify a local decoding of the entire spin interface.
For finitely many physical columns, first locate their
reference faces in an almost surely finite diagonal patch.
The geometric group whose pattern contains $q_z$, if any, contains
every lower-history change at $z$, by
\cref{lem:sos-column-carrier}. It also contains every
third-stage change there, because $z\in D_\delta$ implies
$q_z\in F_\delta$. There is at most one third-stage defect
on a fixed physical column. Read $\cH(z)$ and $\cS(z)$ from
the lower defects and then read the full binary column of
$\cI$ from this third-stage defect, or use the graph column
of $\cS$ if that defect is absent. All limiting groups are
finite and their column sequences have finite transition
sets. Intersect the probability-one finiteness and admissibility
events over the countably many reference dots, and the local
decoding events over all finite sets of integer cubes and
columns, including compatibility of each finite family of
groups. On the resulting probability-one event this decoding
is well defined for every column simultaneously, and all
spin identities hold at every cube. The decoding agrees with the
finite-volume construction on every specified patch and finite list
of group values. This verifies
\cref{lem:physical-reconstruction}, which gives the unique
stationary ergodic physical limit. In particular, the tightness
of the countable group data in \cref{cor:interface-perturbed-joint-diagonal}
supplies the needed finite transition heights; no confinement of
these heights was assumed. The same lemma permits rerooting by
the decoded integer $\cS(0)$ instead of $\phi(0)$, giving the
original interface rooting. Finally
\cref{eq:interface-height-cutoff-tv} extends the local limit to
$M_n/n^2\to\infty$.

\item\label[step]{st:interface-displacements} \emph{(Displacement tails and mean slopes.)}
Put $R_1=\cH-\phi$ and $R_2=\cS-\cH$. Whenever either is nonzero
at $z$, its base \SOS component lies in one geometric group,
whose enlarged pattern $P$ contains $q_z$.
The confinement proof of \cref{lem:sos-column-carrier}
uses only the lower history and gives
$|R_1(z)|+|R_2(z)|\leq C|P|$ for planar groups, or for
periodic groups of diameter at most $\eta_* n$.
Larger connected periodic patterns have size at least $c_*n$.
After horizontal-face conditioning, which costs at most the
bounded factor $N_0/n^2$, the group tail
\cref{eq:interface-perturbed-primitive-component-tail} therefore
gives \cref{eq:sos-displacement-tail} for $R_1(0),R_2(0)$ and
$1\leq r\leq c_*n$.
Pass to the physical limit at fixed $r$, using that the
displacements are integer-valued, and then sum the tails as in
\cref{prop:sos-physical-joint}. Both displacements have every
finite moment, with
$\E(|R_1(0)|+|R_2(0)|)^p\leq C_pe^{-c'\beta^{2/5}}$.
\Cref{lem:physical-reconstruction} now gives logarithmic growth
of both displacements, the common limiting plane and all three
mean slopes.

\item\label[step]{st:interface-sorting-identities} \emph{(Sorting identities and absence of bubbles.)}
The limiting planes just established give the required ends
along the sorting directions. Every third-stage column change is read in
its entirety from one third-stage defect and has zero signed charge,
so $\Tthree\cI=\cS$ on every column, by \cref{def:tetris-charge}.
The lower identities $\Ttwo\cS=\cH$ and $\Tone\cH=\phi$ follow from
the finite zero-charge changes by
\cref{st:sos-sorting-identities} in the proof of
\cref{thm:sos-gff-complete}. For the no-bubble condition, each
third-stage defect is an admissible isolated configuration on a finite
planar support. Apply the first-exit argument of
\cref{lem:polymer-factorization}: a same-sign path starting at a
cube of a changed column either stays in this state until it
reaches the end of that state, or exits into an unchanged graph
column and then continues vertically along that column to its
infinite end. The argument remains valid after inserting all
other defects, since distinct defects have pairwise
non-$*$-adjacent supports, so that the columns adjacent to a
support are unchanged in the full configuration. Unchanged
columns have their vertical graph paths. This proves that both
spin phases of $\cI$ have paths to their respective infinite
vertical ends.
These assertions are made on the common event just constructed:
the defects at every dot are admissible and the supports
are pairwise non-$*$-adjacent. There are only countably many
columns and starting cubes, so the zero-charge identities and
these path constructions hold simultaneously for all of them.
Thus neither column-charge conservation nor absence of bubbles
is being inferred from weak convergence of spin cylinders alone.

\item\label[step]{st:interface-gff-transfer} \emph{(The \GFF transfer.)}
\Cref{prop:interface-diagonal-gff,lem:diagonal-height-second-moment}
give the diagonal \GFF and reference height moment.
\Cref{lem:physical-reconstruction,st:interface-displacements}
give the physical stationary ergodic coupling, the reference mean
slopes and the $L^2$ correction $\cS-\phi=R_1+R_2$. Apply
\cref{cor:physical-height-transfer} to obtain the scaling limit
asserted in \cref{thm:interface-gff-complete}, with amplitude
$\upsigma_{\mathrm d}^U/\rho_3$, map $L_{\mathrm d}^UA^{-1}$ and
local $H^{-s}$ convergence for $s>2$. Taking $U=0$ gives the
Ising Interface case of \cref{thm:gff-limits}.
\end{enumerate}
\end{proof}

\subsection{The full 3D Ising model}\label{sec:ising}
\label{sec:ising-physical-limit}

The full 3\Dim Ising model of \cref{def:3dIsing} differs from the
Ising Interface model of \cref{def:interface-model} in two respects.
First, a spin configuration may contain ``bubbles'': finite
islands of one phase surrounded by the other, whose boundaries are
closed plaquette surfaces that do not separate the two ends of the
cylinder. Second, it may contain several end-separating surfaces,
stacked one above the other, with the phases alternating between
them. Neither feature is present in the three-step analysis of
\cref{sec:interface,sec:model-grassmann}, which was carried out for
bubble-free interfaces. The purpose of this subsection is to reduce
the full model to that analysis. The reduction has two parts, and
the two parts dictate the choice of the height cutoff
\cref{eq:ising-height-cutoff}.

The first part concerns the bubbles. Conditional on the
end-separating surface $\cI$, the bubbles form a gas of
closed contours with the Peierls weights $e^{-\beta|\gamma|}$,
interacting with one another and with $\cI$ only through a
hard-core constraint. At low temperature such a gas has a
convergent cluster expansion, and summing it out produces a
correction $U_{\beta,\Gamma}(\cI)$ to the energy of $\cI$ which
is a sum of exponentially decaying local terms, each supported on
the plaquettes near a cluster of contours touching $\cI$
(\cref{prop:ising-one-sheet-contours}). This correction is
precisely an admissible contact perturbation in the sense of
\cref{def:interface-perturbation}
(\cref{prop:ising-admissible-perturbation}), so that the
perturbed interface analysis of \cref{sec:interface,sec:model-grassmann}
applies verbatim to the single-surface law with weights
$e^{-\beta(\sA(\cI)-N_0)+U_{\beta,\Gamma}(\cI)}$.

The second part concerns the comparison between the finite
cylinder and this single-surface law. Two errors must be
controlled. A configuration with three or more end-separating
surfaces costs at least $e^{-2\beta N_0}$ relative to one with
a single surface, since each additional pair of surfaces has
area at least $2N_0$, while the number of places where such a
pair may sit is polynomial in the cylinder volume; this error is
small as long as the half-height $M$ is much smaller than
$e^{(\beta-C)N_0}$ (\cref{prop:ising-multiple-surfaces}).
Conversely, the single surface is free to translate vertically
inside the cylinder, and rooting it by its column
magnetization averages its weight over roughly $2M$ translates.
The contour error decays exponentially with the distance to the
cutoff, so its sum over these translates is bounded by
$C\sA(\cI)$, as is the error from translates that do not fit.
After averaging over interfaces and normalizing, the error is at
most $C(N_0+1)e^{\delta_\beta N_0}/M$, with $\delta_\beta\to0$
(\cref{prop:ising-joint-cutoff}).
Since $N_0=n^2+n\Theta_1+n\Theta_2$ is of order $n^2$ on a compact slope range,
the choice $M_n=\ceil{e^{n^2}}$ lies in the admissible window
$e^{\delta_\beta N_0}\ll M\ll e^{(\beta-C)N_0}$ once $\beta$ is
large, and with this choice the full Ising cylinder has a single
end-separating surface with probability $1-O(e^{-cn^2})$.

The height observable of the full model is the vertical column magnetization of the end-separating surface,
with its closed bubbles filled; this is the surface erasure of
\cref{def:ising-erasure}. The erased surface is compared, in
total variation and jointly with the bulk spins in any fixed
window, with the rooted perturbed interface law
(\cref{thm:ising-rg-input}). The \GFF for the full model is then
the perturbed-interface \GFF of \cref{thm:interface-gff-complete},
and the bulk spins are recovered from the limiting interface by
the closed-contour kernel of \cref{prop:ising-contour-kernel}
(\cref{thm:ising-gff-complete}).

\begin{theorem}[3\Dim Ising reduction to the \RG input]
\label{thm:ising-rg-input}
Fix a compact positive slope range. For sufficiently large $\beta$,
let $n\to\infty$ in the full Ising cylinders with
$M_n=\ceil{e^{n^2}}$ and admissible windings in this range.
The erased surface, rooted by $(\Tthree\cI)(0)=0$, is at total-variation
distance tending to zero from the normalized rooted surface law
with weights $e^{-\beta(\sA(\cI)-N_0)+U_{\beta,\Gamma}(\cI)}$.
With the diagonal rooting of \cref{prop:physical-interface-law},
the latter has the \Tetris image weights $Z^U_{\beta,\Gamma}(\varphi)$
of \cref{eq:interface-image-weight}.
They admit a positive canonical extension and exact Grassmann
coefficients with the absolute \texttt{thin}/\texttt{fat} bounds, source identities
and torus compatibility of \cref{thm:interface-rg-input}.
The contour factor is included in all absolute bounds. Thus
\crefrange{it:rg-positive}{it:rg-sources} and~\cref{it:rg-marked}
of \cref{thm:coefficient-rg} hold for these comparison weights.
The comparison also holds jointly with the original bulk spins in
every fixed rooted window, using the closed-contour kernel
$\cK_\Gamma(\cI)$ and the orientation prescribed by the boundary conditions.
\end{theorem}

The proof is completed after \cref{prop:ising-joint-cutoff}.

\subsubsection{Closed contours in the full Ising model}
\label{sec:ising-bubbles}

Recall from \cref{lem:phase-tree} that a binary configuration on
the periodic cylinder, with occupied lower end and vacant upper
end, determines a finite tree: its vertices are the $*$-connected
spin components and its edges are the edge-connected components of
the disagreement plaquettes, each such component being incident to
exactly one occupied and one vacant spin component. A bubble-free
configuration is one whose tree has a single edge. For a general
configuration we need a canonical way to select one end-separating
surface and to fill in everything else, so that the result is
bubble-free and the construction commutes with the symmetries of
the model.

\begin{definition}[Surface erasure]\label{def:ising-erasure}
In a finite tilted cylinder, take the unique path in the
phase tree of \cref{lem:phase-tree} from the occupied
lower-end vertex to the vacant upper-end vertex. Its length
is odd, because the spin components along the path alternate in
sign, starting from an occupied and ending at a vacant one.
When the path has one edge,
let $F_*$ be that edge; otherwise select its middle edge to fix a
convention on the exceptional configurations. Define $\cE(\sigma)$
to be the binary
configuration with disagreement set $F_*$ and occupied lower
end. Call a plaquette component \emph{end-separating} when
its edge lies on this tree path, and write $K(\sigma)$ for
the number of these components.
\end{definition}

Thus $K(\sigma)$ is always odd, $K(\sigma)=1$ exactly when the
configuration has a single end-separating surface, and on that
event $\cE(\sigma)$ retains this surface and fills every bubble
with the phase that surrounds it. The following proposition
records that $\cE$ is well defined and computes the exact weight
of its image.

\begin{proposition}[A canonical finite-volume erasure]
\label{prop:ising-erasure}
The map $\cE$ is well defined, takes values in the
bubble-free interface model, is the identity on that model,
and commutes with translations and reversal of the two ends
together with spin reversal. Extending the fixed ends does
not change it for a fixed configuration. Its exact marginal is
\begin{equation}\label{eq:ising-erasure-weight}
 W_\beta(\cI)=e^{-\beta\sA(\cI)}\Xi_\beta(\cI),\qquad
 \Xi_\beta(\cI)=\sum_{\sigma:\cE(\sigma)=\cI}
                     e^{-\beta(\scH(\sigma)-\sA(\cI))}\,.
\end{equation}
Here the energy is the disagreeing-bond count of
\cref{def:3dIsing}.
\end{proposition}
\begin{proof}
The cochain construction in the proof of \cref{lem:phase-tree}
associates with the plaquette component $F_*$ a periodic binary
function $\tau_{F_*}$ whose disagreement set is exactly $F_*$ and
which is constant on each of the two ends. Since $F_*$ is a
bridge of the phase tree whose deletion separates the two end
vertices, the function $\tau_{F_*}$ takes different values on the
two ends; replacing it by its complement if necessary, we may
take the lower end to be occupied. This is the configuration
$\cE(\sigma)$, which is therefore well defined.

To see that $\cE(\sigma)$ is bubble-free, observe that its two
phases are the unions of the original spin components on the two
sides of the tree cut at $F_*$. Each side is connected: adjacent
tree vertices are spin components joined by a plaquette component,
hence they contain cubes sharing a face, and once both are filled
with a common sign these cubes are $*$-adjacent cubes of the same
phase. Thus every cube of the new occupied phase is connected to
the lower end and every cube of the new vacant phase to the upper
end, which is the no-bubble condition of
\cref{def:interface-model}; its disagreement set is exactly $F_*$,
whence $\cE(\sigma)$ belongs to the bubble-free interface model.
If $\sigma$ is itself bubble-free, its phase tree has a single
edge by \cref{lem:phase-tree}, so $F_*$ is its full disagreement
set and $\cE(\sigma)=\sigma$.

The constructions are intrinsic to the configuration, so $\cE$
commutes with translations. Reversing the two ends together with
the spins reverses the tree path, and the middle edge of an odd
path is unchanged under reversal; this proves the remaining
symmetry. Extending the fixed ends does not change the phase
tree of a fixed configuration, hence does not change $\cE(\sigma)$.
Finally, the marginal of $\cE$ under the Ising law is obtained by
grouping the finite Boltzmann sum $\sum_\sigma e^{-\beta\scH(\sigma)}$
according to the value of $\cE(\sigma)$ and factoring out
$e^{-\beta\sA(\cI)}$; this gives \cref{eq:ising-erasure-weight}
exactly, since each $\sigma$ is counted once in the fiber of its
image.
\end{proof}

On $\{K=1\}$, erasure leaves the unique end-separating interface
unchanged and fills precisely the closed contours. The middle-edge
selection is needed only on $\{K\geq3\}$, whose probability is
$O(e^{-cn^2})$ under the cutoff \cref{eq:ising-height-cutoff},
by \cref{prop:ising-multiple-surfaces}. We first compute
the contribution with $K=1$ to the weight $\Xi_\beta(\cI)$ of
\cref{eq:ising-erasure-weight}, which will turn out to be the
partition function of a hard-core gas of closed contours.

We define the contours as follows. For a finite cylinder
$\Lambda$, let $\scC_\Lambda$ be the possible nonempty
edge-connected disagreement sets $\gamma$ of binary configurations
equal to zero outside $\Lambda$ and with equal values in the two
ends. Equivalently, by the cochain argument of
\cref{lem:phase-tree}, they are the connected closed plaquette
cochains with zero end parity; each such $\gamma$ has a unique
binary \emph{filling} $F_\gamma$, namely the configuration with
disagreement set $\gamma$ which vanishes outside $\Lambda$, and
we fix this filling to be zero outside. These are the classical
Peierls contours of the low-temperature expansion (see, e.g.,
\cite[\S3.7]{FriedliVelenik18}), with the difference that
noncontractible contours are allowed: a contour may wind around the
periodic directions of the cylinder, so its filling need not be a ball. Two contours are
incompatible if they share a plaquette or if some of their
plaquettes share an edge. The same rule defines incompatibility
with a fixed interface $\cI$. Put
\[
 z_\beta(\gamma)=e^{-\beta|\gamma|},\qquad
 Z^0_\Lambda=\sum_{\{\gamma_i\}\text{ compatible}}
                                      \prod_i z_\beta(\gamma_i)\,,
\]
where the sum is over families of pairwise compatible contours,
including the empty family with weight one. Thus $Z^0_\Lambda$ is
the partition function of the unrestricted closed-contour gas. Equivalently, it is the Ising partition function with both
ends occupied, restricted to configurations with no end-separating
disagreement component. Without this restriction a pair of separating
surfaces could occur and would not be counted by this gas.
The normalization is independent of the interface and cancels in
all comparisons below.

\begin{proposition}[The single end-separating surface contour factor]
\label{prop:ising-one-sheet-contours}
For every bubble-free $\cI$ fitting in $\Lambda$,
\begin{equation}\label{eq:ising-one-sheet-gas}
 \Xi^{(1)}_{\beta,\Lambda}(\cI)
 :=\sum_{\substack{\cE(\sigma)=\cI\\K(\sigma)=1}}
             e^{-\beta(\scH(\sigma)-\sA(\cI))}
 =\sum_{\substack{\{\gamma_i\}\text{ compatible}\\
                         \gamma_i\text{ compatible with }\cI}}
                                      \prod_i z_\beta(\gamma_i)\,.
\end{equation}
For all sufficiently large $\beta$ the normalized logarithm
\[
 U_{\beta,\Lambda}(\cI)
           =\log\Xi^{(1)}_{\beta,\Lambda}(\cI)-\log Z^0_\Lambda\,,
\]
is nonpositive and has an absolutely convergent cluster expansion,
uniformly in $\cI$, the cutoff and the periods. For some
$\epsilon_\beta\leq Ce^{-c\beta}$,
\begin{align}
 |U_{\beta,\Lambda}(\cI)|&\leq\epsilon_\beta\sA(\cI),
                                                \label{eq:ising-U-size}\\
 |U_{\beta,\Lambda}(\cI)-U_{\beta,\Lambda}(\cI')|
                   &\leq\epsilon_\beta|\cI\xor \cI'|\,.
                                                \label{eq:ising-U-difference}
\end{align}
It has a translation-covariant infinite-height limit
$U_{\beta,\Gamma}(\cI)$ for every fixed periodic finite-area
interface. The error from a boundary at distance $r$ is at
most $\epsilon_\beta\sA(\cI)e^{-c\beta r}$.
\end{proposition}
\begin{proof}
We identify the contour gas and apply the Koteck\'y--Preiss theorem
with an exponential reserve for the distance estimates.
\begin{enumerate}[label=\textbf{Step~\arabic*.}, ref=\arabic*, wide=0pt, itemsep=1ex, listparindent=\parindent, parsep=0pt]
\item\label[step]{st:ising-gas-representation} \emph{(The gas representation.)}
In the phase tree, $K=1$ means that $\cI$ is the unique
end-separating boundary component. Every other boundary
component is an edge of the tree not on the path between the
ends; deleting it separates the tree into two parts both of
whose end vertices lie on the same side, so the component has
zero end parity and is an element of $\scC_\Lambda$.
Distinct components of one configuration share no plaquette and
no edge, hence are compatible with one another and with $\cI$.
Conversely, take such a compatible family and add its cochains
modulo two to the binary interface configuration. Their plaquette
sets are disjoint, so the disagreement set of the sum is precisely
the union of the contours with $\cI$, and its edge-connected
components are exactly the individual contours and $\cI$; the phase
tree of the sum therefore has exactly the one end-separating edge
$\cI$. These two maps are mutually inverse. Note that there is no
independent choice of signs for nested contours: in the
low-temperature expansion of the Ising model the sign of the phase
inside each contour is determined by the parity of the number of
contours surrounding it, and here the occupied lower end fixes all
the signs. Since the energy $\scH(\sigma)$ is the total number of
disagreement plaquettes, it equals $\sA(\cI)+\sum_i|\gamma_i|$,
proving~\cref{eq:ising-one-sheet-gas}. Its sum is restricted to a
subset of the contour families counted by $Z^0_\Lambda$, so
$\Xi^{(1)}_{\beta,\Lambda}(\cI)\leq Z^0_\Lambda$ and $U\leq0$.

\item\label[step]{st:ising-cluster-convergence} \emph{(Convergence of the cluster expansion.)}
Apply the Koteck\'y--Preiss theorem~\cite{KoteckyPreiss86}, in the
rooted form of \cite[Thm.~5.4]{FriedliVelenik18}, with size function
$a(\gamma)=|\gamma|$. Put
$w(\gamma)=e^{-3\beta|\gamma|/4}$, so that
$z_\beta(\gamma)=w(\gamma)e^{-\beta|\gamma|/4}$, and let $k(p,q)$
indicate that the plaquettes $p,q$ coincide or share an edge.
The plaquette graph has uniformly bounded degree, also on a
quotient, so $G=\sup_p\sum_q k(p,q)\leq C$ and the number of
connected $m$-plaquette sets containing $p$ is at most $C_0^m$,
by the standard lattice-animal bound.
Discarding the cochain and filling constraints, which can only
decrease the number of contours, therefore gives
\[
 \delta_\beta:=\sup_p\sum_{\gamma\ni p}w(\gamma)e^{|\gamma|}
 \leq\sum_{m\geq1}C_0^m e^{-(3\beta/4-1)m}
 \leq Ce^{-c\beta}\,.
\]
For large $\beta$, $G\delta_\beta\leq1$. Since a contour
incompatible with $\gamma$ contains a plaquette $q$ with
$k(p,q)=1$ for some $p\in\gamma$, summing over $p\in\gamma$ and
then over $q$ gives
\begin{equation}\label{eq:ising-polymer-tree-budget}
 \sum_{\gamma':\gamma'\not\sim\gamma}
              e^{-3\beta|\gamma'|/4}e^{|\gamma'|}
                    \leq G\delta_\beta|\gamma|\leq|\gamma|\,.
\end{equation}
Here $\not\sim$ denotes incompatibility, including equality.

Write
\[
 \phi^T(\gamma_1,\ldots,\gamma_k)
 =\sum_{G\text{ connected on }[k]}
        \prod_{ij\in E(G)}(-\one_{\{\gamma_i\not\sim\gamma_j\}}),
 \qquad \phi^T(\gamma)=1.
\]
The cited theorem gives the absolutely convergent logarithmic
expansion with coefficients $\phi^T/k!$, and its rooted estimate is
\[
 \sum_{k\geq1}\frac1{(k-1)!}
 \sum_{\gamma_2,\ldots,\gamma_k}
 |\phi^T(\gamma,\gamma_2,\ldots,\gamma_k)|
                     \prod_{i=2}^k w(\gamma_i)\leq e^{|\gamma|}.
\]
The hypothesis is precisely \cref{eq:ising-polymer-tree-budget}.
For every plaquette $p$, distinguish one contour touching $p$
and apply this rooted estimate. Summing that root gives at most
$G\delta_\beta$, hence
\begin{equation}\label{eq:ising-anchored-clusters}
 \sum_{k\geq1}\frac1{k!}
 \sum_{\substack{\boldsymbol\gamma:\
          \exists i\;\gamma_i\not\sim\{p\}}}
 |\phi^T(\boldsymbol\gamma)|
        e^{-\beta\sum_i|\gamma_i|}
        e^{(\beta/4)\sum_i|\gamma_i|}
                         \leq Ce^{-c\beta}\,.
\end{equation}
Indeed, choosing the root changes $1/k!$ to $1/(k-1)!$, and the
activities on the left are exactly $w$. Thus a quarter of the
original exponent remains available. Apply the theorem first to
finite contour families; the rooted bounds persist under
exhaustion. In a finite cylinder the animal count also gives
$\sum_\gamma w(\gamma)e^{|\gamma|}<\infty$, so the full logarithmic
series converges there. Removing allowed contour types only
improves the criterion.

\item\label[step]{st:ising-restricted-gas} \emph{(The restricted gas and the bounds on $U$.)}
The restricted gas in \cref{eq:ising-one-sheet-gas} is the same
hard-core gas with the contours incompatible with $\cI$ removed,
and its logarithm has the same expansion with the sum restricted
to clusters all of whose contours are compatible with $\cI$. The
clusters that avoid $\cI$ contribute identically to both
logarithms. Subtracting the unrestricted logarithm therefore
leaves the formula
\begin{equation}\label{eq:ising-U-clusters}
 U_{\beta,\Lambda}(\cI)
 =-\sum_{k\geq1}\frac1{k!}
   \sum_{\substack{\boldsymbol\gamma:\
                  \exists i\;\gamma_i\not\sim \cI}}
          \phi^T(\boldsymbol\gamma)e^{-\beta\sum_i|\gamma_i|}\,,
\end{equation}
in which every surviving cluster touches $\cI$. Root each
contributing cluster at a face of $\cI$ that it touches; the sum
over the faces of $\cI$ of the anchored bound
\cref{eq:ising-anchored-clusters} (with the reserved factor
dropped) yields \cref{eq:ising-U-size}. For
\cref{eq:ising-U-difference}, observe that a cluster contributes
to $U(\cI)$ and not to $U(\cI')$ only if it touches
$\cI\setminus\cI'$, and symmetrically, so the difference is a sum
over clusters touching $\cI\xor\cI'$, and anchoring at a face
of the symmetric difference gives the bound.

For the distance tail, note that a cluster with nonzero Ursell
function has a connected incompatibility graph, and that
incompatible contours touch, so a cluster of total
plaquette number $m$ has diameter at most $Cm$. A cluster
touching both $\cI$ and a set at distance $r$ from $\cI$ thus
has $m\geq cr$, and the reserved factor $e^{-\beta m/4}$ in
\cref{eq:ising-anchored-clusters} gives it the weight
$e^{-c\beta r}$. Increasing the cylinder admits every finite
cluster on the quotient, and the absolute summability just
proved allows dominated convergence, which gives the
infinite-height limit $U_{\beta,\Gamma}$; the terms missing at a
finite cutoff are the clusters connecting $\cI$ to the boundary,
which pay the stated distance tail. All counts include noncontractible
contours, since the animal bound and the bounded degree of the
plaquette graph hold on the quotient.
\end{enumerate}
\end{proof}

\begin{remark}[Distant variations of the contour interaction]
A cluster can touch the interface in several places, but its
contribution to a mixed variation supported on disjoint sets $A,B$
vanishes unless it touches both. More precisely, if
$\cI,\cI_A,\cI_B,\cI_{AB}$ are admissible interfaces with the same
respective modifications on $A$ and $B$, then
\[
 |U(\cI_{AB})-U(\cI_A)-U(\cI_B)+U(\cI)|
 \leq C\epsilon_\beta\min\{|A|,|B|\}e^{-c\beta\dist(A,B)}\,.
\]
Indeed, a cluster's compatibility indicators factor as $c_Xa_i b_j$,
where $c_X$ concerns the common part outside $A\cup B$ and $a_i,b_j$
the corresponding parts inside $A,B$. Their mixed
difference is $c_X(a_1-a_0)(b_1-b_0)$, of modulus at most one.
A contributing cluster has mass at least $c(1+\dist(A,B))$, since
its unit neighborhood is connected and has diameter at most a
constant times its mass. Anchoring at a touched face of the smaller
set and using the reserve in \cref{eq:ising-anchored-clusters}
gives the bound; absolute convergence justifies the cancellation.
\end{remark}

Since $U\leq0$, the normalized single-surface weights
$e^{-\beta(\sA(\cI)-N_0)+U_{\beta,\Gamma}(\cI)}$ have a finite
physically rooted partition function, bounded by the unperturbed
one of \cref{prop:physical-interface-law}.

The erased surface discards the bubbles. To recover the original
spins from the limiting interface, we need the conditional law of
the bubbles given $\cI$, which by \cref{eq:ising-one-sheet-gas} is
the compatible closed-contour gas, and we need this law to be local
in $\cI$ and to have an infinite-height limit. Both properties are
proved by the switching argument of \cref{prop:pi_decay}, applied
to two independent copies of the contour gas.

\begin{proposition}[Local reconstruction of the closed-contour spins]
\label{prop:ising-contour-kernel}
For sufficiently large $\beta$, the closed-contour gas compatible
with a fixed interface $\cI$ has a unique local limit under exhaustion
of the height cutoff. Denote its spin reconstruction kernel by
$\cK_\Gamma(\cI)$: sample the limiting gas and add the contour
fillings modulo two to the binary configuration of $\cI$.
The same construction gives a kernel $\cK(\cI)$ on the full
lattice. These kernels are translation covariant. If the binary
reference configurations associated with $\cI,\cI'$ by the lower-end
convention agree on the radius-$r$ neighborhood of a finite set of
cubes $A$, then
\begin{equation}\label{eq:ising-contour-kernel-locality}
 \|\cK(\cI)|_A-\cK(\cI')|_A\|_\mathrm{TV}
                    \leq C|A|e^{-c\beta r}\,.
\end{equation}
The same bound compares a cutoff kernel with its infinite-height
kernel when $A$ is at distance at least $r$ from the boundary.
It also compares a periodic kernel with the full-lattice kernel
when the radius-$r$ neighborhood embeds and the reference spins
agree there. For bounded spin observables $f,g$ supported on $A,B$,
\begin{equation}\label{eq:ising-contour-kernel-mixing}
 |\cov_{\cK(\cI)}(f,g)|
 \leq C\|f\|_\infty\|g\|_\infty(|A|+|B|)
                         e^{-c\beta\dist(A,B)}\,.
\end{equation}
All constants are independent of $\cI$, the cutoff and the periods.
In particular, a unique mixing local limit of erased interfaces
would give a unique mixing joint limit with these bulk spins.
\end{proposition}
\begin{proof}
We apply the allowed-family version of \cref{prop:pi_decay}.
The geometric input is a uniform diameter bound for contour
fillings, needed because spins depend on fillings, not just on the
contour plaquettes.

Write $m(\gamma)=|\gamma|$ and let $F_\gamma$ be the unique
zero-end filling. A connected contour of size $m$ has diameter at
most $Cm$ in the quotient metric. To control its filling uniformly,
including auxiliary sectors with a zero coordinate face density,
use diagonal rays in direction $d=(1,1,1)$. On a matching cylinder put
$\ell(x)=\rho\cdot x$, where
$\rho_i=N_i(M)/N_0\geq0$ and $\sum_i\rho_i=1$.
This function annihilates the period lattice, hence descends to
the quotient, satisfies
$\ell(d)=1$, and has uniformly bounded Lipschitz constant.
Thus the contour's $\ell$-range has length at most $Cm$.
The zero-end filling vanishes outside that range: from a cube
outside the range a diagonal ray to the nearer end crosses no
plaquette of the contour, and the filling is zero at that end.
From any filled cube a
generic diagonal ray must therefore hit the contour within
distance $Cm$. Combining these paths with connectedness of the
contour shows that its entire filling has diameter at most $Cm$.
The planar assertion follows using $\ell(x)=(x_1+x_2+x_3)/3$. In
particular, if $F_\gamma(v)=1$, some contour face is within $Cm$
of $v$. The bounded-degree animal count, with at most $Cm^3$
possible root faces within this distance, gives the following
two bounds. For $a=\beta/8$ and $\beta$ sufficiently large,
\begin{align}
 \sum_{\gamma'\not\sim\gamma}
       z_\beta(\gamma')e^{a m(\gamma')}m(\gamma')
       &\leq\kappa_\beta m(\gamma),\qquad \kappa_\beta<\tfrac12,
                                           \label{eq:ising-ancestor-budget}\\
 \sum_{\gamma:F_\gamma(v)=1}
       z_\beta(\gamma)e^{a m(\gamma)}m(\gamma)
       &\leq Ce^{-c\beta}\,.                 \label{eq:ising-filling-budget}
\end{align}
Dropping interface restrictions only increases these sums.

Compare finite gases with allowed contour families
$\mathcal C'\subset\mathcal C$ and the same reference spins.
The spins on $A$ depend only on the contour types whose fillings
meet $A$. By \cref{prop:pi_decay}, their total-variation distance
is bounded by the probability of a union component joining one
of these types to $\mathcal C\setminus\mathcal C'$. Such a
component contains a simple alternating path between the two
classes. By \cref{lem:deletion} and independence, a specified
simple colored path has probability at most the product of its
activities. We sum over such paths by a transfer-operator bound.
Set
\[
 T_a(\gamma,\gamma')=
 \one_{\{\gamma'\not\sim\gamma\}}
                  z_\beta(\gamma')e^{a m(\gamma')}.
\]
By \cref{eq:ising-ancestor-budget}, $T_a m\leq\kappa_\beta m$.
Dropping the path restrictions and summing over its two starting
colors, the total path weight, with the extra factor
$e^{a\sum_i m(\gamma_i)}$, is at most
\[
 2\sum_{v\in A}\sum_{\gamma:F_\gamma(v)=1}
 z_\beta(\gamma)e^{a m(\gamma)}
             \sum_{j\geq0}(T_a^j1)(\gamma)
 \leq\frac{C|A|e^{-c\beta}}{1-\kappa_\beta}.
\]
Here $1\leq m$ and \cref{eq:ising-filling-budget} give the last
inequality. Repeated contour types in this larger sum are harmless.

Let $\widehat\gamma$ consist of the contour, its filling and their
unit neighborhood; it has diameter at most $Cm(\gamma)$. Restrict the
smaller gas to $\widehat\gamma\subset D$, with $D$ the
radius-$r$ neighborhood of $A$. An exceptional path starts in $A$
and leaves $D$. Consecutive contours touch, so its total
size is at least $cr$. The weighted bound above yields
\[
 \|\cL_\pi(\sigma|_A)-\cL_{\pi'}(\sigma|_A)\|_{\mathrm{TV}}
                              \leq C|A|e^{-c\beta r}.
\]
This estimate is uniform in the larger contour family. The
restricted gas is finite and depends on the reference spins only in
$D$: their disagreements determine every compatibility test because the unit neighborhoods were retained, and their values on $A$ specify the spins
to which the contour fillings are added. Two large cutoffs therefore
give spin laws within
$C|A|e^{-c\beta r}$ of it. Letting $r\to\infty$
proves existence and uniqueness of the local exhaustion limit.
The same argument proves convergence of the local contour configurations:
start paths at contours meeting $A$, for which the root sum has
the same bound by the animal count. Deletion and
\cref{eq:ising-filling-budget} give a uniform exponentially small
probability of any contour of size greater than $m$ whose filling
contains a fixed cube. In the limit only finitely many fillings
contain each cube, and the spin is exactly their parity added to
the reference spin. Countability makes this reconstruction valid
simultaneously at all cubes.

The restricted laws also agree when the reference spins agree in $D$,
when two cutoffs contain $D$, or when $D$ embeds in a periodic
quotient and the full lattice. The same estimate and the triangle
inequality prove \cref{eq:ising-contour-kernel-locality} and its
cutoff and periodic versions. This treats noncontractible contours as
well: a contour for which $\widehat\gamma$ is not contained in the
common neighborhood can contribute only through a path paying
the same size cost. The local approximation makes the limiting
kernels measurable and translation covariant.

For mixing, restrict to contours whose enlarged sets
$\widehat\gamma$ lie in either of the two neighborhoods of
$A,B$ of radius one quarter of their separation. The two
restricted gases have no mutual incompatibilities and hence
factorize. The switching estimate, applied to both starting
sets, bounds the total-variation error for the joint spins by
$C(|A|+|B|)e^{-c\beta\dist(A,B)}$. Comparing the joint expectation
and the two marginals with this product law proves
\cref{eq:ising-contour-kernel-mixing}. For bounded separations,
compare instead with the deterministic reference spins: the filling
estimate bounds this error by $C(|A|+|B|)e^{-c\beta}$.
Finally, approximate conditional means by these finite local
kernels. Local convergence of the reference interfaces gives joint
convergence with the spins. Splitting a joint covariance into its
conditional covariance and the covariance of conditional means,
as in \cref{lem:joint-limit}, proves mixing whenever the rooted
reference law is mixing. The argument commutes with common
rerooting under physical translations.
\end{proof}

\subsubsection{Verification of the interface perturbation hypotheses}
\label{sec:ising-perturbation-application}

We now verify that the bubble correction belongs to the
contact-perturbation class used in the interface analysis.
Recall from \cref{def:interface-perturbation} that an admissible
contact perturbation is specified by a countable family of labels
$X$, each with an edge-connected plaquette carrier $N_X$, a mass
$m_X\geq1$ with $|N_X|\leq C_0m_X$ and a real coefficient $b_X$;
the label contributes $b_Xe^{-tm_X}$ to the energy of an interface
exactly when the interface meets its carrier. The requirements are
the anchored norm \cref{eq:interface-contact-norm} and the
compatibility of the labels under lifting a pattern to the plane.
In the present application the labels are the contour clusters,
and the parameter $t$ is the inverse temperature $\beta$ itself.

\begin{proposition}[The bubble correction is an admissible perturbation]
\label{prop:ising-admissible-perturbation}
The infinite-height closed-contour correction $U_{t,\Gamma}$ satisfies
\cref{def:interface-perturbation}, uniformly on all matching
cylinders. Hence the normalized single-surface law and its \Tetris
image satisfy \cref{thm:interface-rg-input}.
\end{proposition}
\begin{proof}
For an ordered contour cluster $X=(\gamma_1,\ldots,\gamma_k)$ set
\[
 m_X=\sum_i|\gamma_i|,\qquad b_X=-\phi^T(X)/k!,
 \qquad
 N_X=\{p:p\text{ equals or is edge-adjacent to a face of }X\}\,.
\]
Retain clusters with connected incompatibility graph and nonzero
$\phi^T$. Their carriers are connected, because incompatible
contours touch, and $|N_X|\leq Cm_X$ by the bounded degree of the
plaquette graph. An interface meets $N_X$ exactly when some contour
of the cluster is incompatible with it, so the indicator
$h_X(\cI)$ of \cref{eq:interface-contact-functional} is the
indicator $\one_{\{\exists i\;\gamma_i\not\sim\cI\}}$ in
\cref{eq:ising-U-clusters}, and the activity $a_t(X)=b_Xe^{-tm_X}$
is the cluster weight $-\phi^T(X)e^{-\beta m_X}/k!$ at $t=\beta$.
Thus \cref{eq:ising-U-clusters} is precisely
\cref{eq:interface-contact-functional}; the factorial is already
in $b_X$ and no extra permutation multiplicity is introduced,
since ordered clusters are the labels.
The uniform norm \cref{eq:interface-contact-norm}
is~\cref{eq:ising-anchored-clusters}, with $t$
in place of $\beta$: the anchoring plaquette $p\in N_X$ is either
a face of the cluster or edge-adjacent to one, which costs at most
a bounded factor in the animal count. Its animal and tree counts
use only the bounded plaquette degree and the two diagonal ends,
and hence also hold on auxiliary cylinders with zero coordinate
face densities.

For the lifting condition, fix a component $\widetilde P$ of the
planar lift of the partial matching $P$, with stabilizer $H_P$.
As in \cref{lem:interface-perturbed-anchored-lifting}, lift this
component to a monotone surface after fixing one face height.
Its face labels determine every height increment along a path in
$\widetilde P$, and in particular the physical period map $g_M$
on $H_P$. Any other matching containing $P$ has the same
period map there. A cluster whose entire carrier projection lies
in this component can therefore be copied with its relative integer
heights and these identifications. Periods outside $H_P$ move
$\widetilde P$ to other components of the lifted pattern; they
create no additional identifications within the copied carrier.
The edge-neighborhoods included in $N_X$ also retain every contact
and incompatibility test. A new contact between different lifted
copies would force their carrier projections to share a triangle,
placing them in the same component and hence under an identification
already recorded in $H_P$.

The contours must still have zero-end fillings after transport.
For each $\gamma_i$, every generic diagonal line through a filled
cube crosses $\gamma_i$ before reaching its zero upper end. Thus
the projection of its filling lies inside its boundary projection,
and hence inside $\widetilde P$. Transport this filling with the
same period map, taking it to be zero outside the component.
Every boundary incidence is retained in the carrier, and there are
no new incidences between distinct lifted components, as just
checked. Its boundary is therefore the copied $\gamma_i$.
Along each diagonal line its relative heights are unchanged, so
its values at both ends remain zero. Uniqueness of the zero-end
filling identifies it with the filling of the transported contour.
This proves that no additional period constraint outside $H_P$
is needed.

Copying back gives the inverse correspondence. It preserves the
ordered list of contours, masses, incompatibility graph and hence
$\phi^T/k!$. For an embedded pattern $H_P=\{0\}$, this is exactly
lifting to the plane and copying back. Relative heights may be
arbitrarily large but finite; the anchored norm justifies summing
the bijection. This verifies the remaining requirement of
\cref{def:interface-perturbation}.
\end{proof}

\subsubsection{Comparison at the joint height cutoff}
\label{sec:ising-multiple-surfaces}

The weighted a priori bound~\cref{eq:interface-perturbed-history-budget} applies
to the contour correction by \cref{prop:ising-admissible-perturbation}.
The remaining issue is the possible presence of other separating
surfaces. The following Peierls-type estimate shows that, since any
end-separating surface has area at least $N_0$ (the area of a
monotone surface of the prescribed windings, which is the minimal
area), removing two of them gains a factor $e^{-2\beta N_0}$ which
beats the polynomial volume factor as soon as the cylinder volume
is at most exponential in $N_0$. This is the upper constraint on
the cutoff $M_n$ discussed at the beginning of \cref{sec:ising}.

\begin{proposition}[Additional end-separating surfaces]
\label{prop:ising-multiple-surfaces}
Let $V_\Lambda$ be the number of bonds meeting the finite
tilted cylinder. For sufficiently large $\beta$,
\begin{equation}\label{eq:ising-multiple-surfaces}
 \P_{\beta,\Lambda}(K\geq3)
              \leq C V_\Lambda^2 e^{-2(\beta-C)N_0}\,.
\end{equation}
For $M=M_n=\ceil{e^{n^2}}$ and positive slopes in a fixed compact
set, this probability is at most $Ce^{-cn^2}$ after increasing the
lower threshold on $\beta$. Consequently conditioning on $K=1$
changes the full finite-volume law by a vanishing total-variation error.
\end{proposition}
\begin{proof}
By the cochain construction of \cref{lem:phase-tree}, every
end-separating component $F$ is the disagreement set of a periodic
binary configuration with occupied lower and vacant upper ends.
Its three-step \Tetris image is a monotone surface with the same
windings, of area $N_0$. Since sorting cannot increase area
(\cref{prop:sorting}), $|F|\geq N_0$.

Next, for two specified end-separating components $F_1,F_2$ of a
configuration $\sigma$, add their cochains
modulo two to the spin configuration. Each of $\tau_{F_1},\tau_{F_2}$
has opposite ends, so their sum has zero
end parity, and the modified configuration satisfies the same
boundary condition. The modification removes exactly these two
disagreement components, leaving the others unchanged, and thus
lowers the energy by $|F_1|+|F_2|$. For the specified pair the
map $\sigma\mapsto\sigma+\tau_{F_1}+\tau_{F_2}$ is injective, being
an involution on binary configurations. Comparing the Boltzmann
weights of $\sigma$ and its image, the probability that both
$F_1,F_2$ occur as disagreement components is at most
$e^{-\beta(|F_1|+|F_2|)}$; this is the Peierls argument in its
classical form~\cite{Dobrushin72a}.
Connected plaquette sets of size $m$ have at most
$CV_\Lambda C_0^m$ choices, by the animal bound with the root among
the $O(V_\Lambda)$ plaquettes of the cylinder. If $K\geq3$, the
tree path contains at least three end-separating components, and
any two of them can be removed in this way. Summing the Peierls
bound over pairs of connected plaquette sets with both
sizes at least $N_0$ proves~\cref{eq:ising-multiple-surfaces},
with the entropy $C_0^m$ absorbed into the constant $C$ in the
exponent.

For the stated cutoff, $V_\Lambda\leq Cn^2(M_n+n)$ and $N_0\geq n^2$.
Thus, for $\beta$ large enough that $\beta-C\geq2$,
\[
 V_\Lambda^2 e^{-2(\beta-C)N_0}
 \leq Cn^4\exp(2n^2-2(\beta-C)N_0)
 \leq Cn^4e^{-2n^2}\,,
\]
where we used $M_n\leq e^{n^2}+1$ in the first inequality.
Finally, the total-variation distance between the law and its
conditional law given $K=1$ equals
$\P(K\ne1)=\P(K\geq3)$, since $K$ is positive and odd.
\end{proof}

\begin{proposition}[Identification in an explicit joint cutoff regime]
\label{prop:ising-joint-cutoff}
Fix a compact positive slope range and sufficiently large $\beta$.
Set $M_n=\ceil{e^{n^2}}$ as in \cref{eq:ising-height-cutoff}. Erase the full Ising
configuration in this finite cylinder by
\cref{def:ising-erasure}, then root its interface
by $(\Tthree\cI)(0)=0$. As $n\to\infty$, its law is at
total-variation distance tending to zero from the probability
law on rooted bubble-free interfaces proportional to
\[
               e^{-\beta(\sA(\cI)-N_0)+U_{\beta,\Gamma}(\cI)}\,.
\]
For every fixed finite set of cubes $A$, the same total-variation
comparison holds jointly with the rooted spins in $A$, when the
comparison interface is decorated by the kernel
$\cK_\Gamma(\cI)$ of \cref{prop:ising-contour-kernel}.
This reduces joint local convergence to local convergence of
the contour-weighted interface law.
\end{proposition}
\begin{proof}
Condition first on $K=1$. Write $A(\cI)=\sA(\cI)$,
$E(\cI)=A(\cI)-N_0$ and $U(\cI)=U_{\beta,\Gamma}(\cI)$.
Let $\cA_M(\cI)$ be the integer vertical translates of the rooted
interface fitting the band. By \cref{eq:ising-one-sheet-gas}, after
dividing by the common factor $2MZ^0_\Lambda e^{-\beta N_0}$, its
unnormalized rooted weight is
\[
 w_M(\cI)=\frac{e^{-\beta E(\cI)}}{2M}
             \sum_{k\in\cA_M(\cI)}e^{U_{\beta,\Lambda}(\cI+k)}.
\]
We compare this with $w(\cI)=e^{-\beta E(\cI)+U(\cI)}$ by summing
boundary errors over the translates, without choosing a boundary
layer.

The allowed translates form an integer interval $[k_-,k_+]$,
possibly empty. Indeed, translation by $k$ adds $k$ to the
minimum and maximum of the tilted height
$\ell(x)=x_3+(\Theta_1/n)x_1+(\Theta_2/n)x_2$ on $\cI$.
For a nonempty interval put $j_k=\min(k-k_-,k_+-k)$.
The distance $d_k$ of $\cI+k$ from the cutoff boundaries is
at least $\max(0,cj_k-C)$, with uniform constants on the compact
slope range: $\ell$ has bounded Lipschitz constant, and each unit
vertical translation increases the lower clearance and decreases
the upper clearance by one. There are at most two translates at
each value of $j_k$, so, for sufficiently large $\beta$,
\[
             \sum_{k\in\cA_M(\cI)}e^{-c\beta d_k}\leq C.
\]
The contour tail of \cref{prop:ising-one-sheet-contours} bounds
$|U_{\beta,\Lambda}(\cI+k)-U(\cI)|$ by
$\epsilon_\beta A(\cI)e^{-c\beta d_k}$. Both potentials are
nonpositive, so the exponential is 1-Lipschitz on their range.
Together with
$\big||\cA_M(\cI)|-2M\big|\leq C(A(\cI)+1)$ from
\cref{eq:interface-translate-count}, this gives
\begin{equation}\label{eq:ising-cutoff-weight-error}
 |w_M(\cI)-w(\cI)|
       \leq \frac{C}{M}e^{-\beta E(\cI)}(A(\cI)+1).
\end{equation}
The same inequality holds for an empty interval, by the translate
count alone. Thus it also covers interfaces too large to fit.

Let $\cZ^{(0)}_{\beta,\Gamma}=\sum_\cI e^{-\beta E(\cI)}$ and
$Z_\Gamma=\sum_\cI w(\cI)$. The area estimate in the proof of
\cref{prop:physical-interface-law}, applied with $U=0$, gives
\[
 \sum_\cI(A(\cI)+1)e^{-\beta E(\cI)}
                   \leq C(N_0+1)\cZ^{(0)}_{\beta,\Gamma}.
\]
If $G_\Gamma$ counts the physically rooted monotone interfaces,
\cref{eq:physical-interface-partition,eq:interface-positive-budget}
with $U=0$ and \cref{eq:ising-U-size} give, respectively,
\[
 \cZ^{(0)}_{\beta,\Gamma}\leq G_\Gamma e^{\epsilon_\beta^{\mathrm I}N_0},
 \qquad Z_\Gamma\geq G_\Gamma e^{-\epsilon_\beta N_0}.
\]
Consequently the relative total error is
\[
 \eta_n:=Z_\Gamma^{-1}\sum_\cI|w_M(\cI)-w(\cI)|
 \leq \frac{C(N_0+1)}{M}e^{\delta_\beta N_0},
 \qquad \delta_\beta=\epsilon_\beta+\epsilon_\beta^{\mathrm I}\to0.
\]
Write $N_0/n^2\leq\kappa_+$ on the prescribed slope range. Choose
$\beta$ with $\delta_\beta\kappa_+\leq1/2$. At
$M=M_n\geq e^{n^2}$, the bound is $Cn^2e^{-n^2/2}\to0$.
It also controls normalization: if $Z_M=\sum_\cI w_M(\cI)$, then
$|Z_M-Z_\Gamma|\leq\eta_n Z_\Gamma$, and adding and subtracting
$w_M/Z_\Gamma$ bounds the total-variation distance of the normalized
laws by $\eta_n$.

For the joint spin assertion, let $A$ now denote the fixed rooted
window. Conditional on $\cI+k$, its spins have the compatible
closed-contour law by \cref{eq:ising-one-sheet-gas}. The root
charge lies between the extreme transitions in its column.
Hence the distance of $A+k$ from the cutoff boundaries is at
least $\max(0,cj_k-C_A)$, apart from the bounded number of translates
for which the window crosses a boundary. Use the bound one for
these translates. For all others,
\cref{prop:ising-contour-kernel} bounds the error relative to
$\cK_\Gamma(\cI)|_A$ exponentially in this distance. Summing over
$k$ as above gives a total conditional error at most $C_A$.
Since $e^{U_{\beta,\Lambda}}\leq1$, the extra unnormalized error
for a fixed $\cI$ is at most $C_Ae^{-\beta E(\cI)}/M$.
The preceding summation and normalization therefore prove the
joint comparison as well.

Finally \cref{prop:ising-multiple-surfaces} gives
$\P(K\ne1)\leq Ce^{-cn^2}$. Erasure, rooting and restriction to
$A$ cannot increase total variation, so removing the conditioning
on $K=1$ proves both assertions.
\end{proof}

\begin{proof}[Proof of \cref{thm:ising-rg-input}]
By \cref{prop:ising-admissible-perturbation}, the infinite-height
contour correction $U_{\beta,\Gamma}$ is an admissible contact
perturbation. Apply \cref{thm:interface-rg-input} to it: this gives
the exact \Tetris image, its positive canonical extension and all
absolute coefficient, source and marked-collection bounds. The
physical rooting is that of \cref{prop:physical-interface-law}.
\Cref{prop:ising-joint-cutoff} gives the total-variation comparison
with this rooted surface law, jointly with the bulk spins sampled
from the kernel of \cref{prop:ising-contour-kernel}. Sorting and
rooting cannot increase total variation, so the comparison also
holds for the image laws. These are precisely the assertions of
the theorem.
\end{proof}

\subsubsection{Application of the abstract dimer theorem}
\label{sec:ising-rg-application}

\Cref{prop:ising-admissible-perturbation} places the bubble
correction in the contact-perturbation class of \cref{sec:interface}.
Consequently all conditions of \cref{def:rg-admissible},
including the marked-collection condition, hold for the normalized
single-surface image by \cref{thm:interface-rg-input,lem:interface-marked-input}. \Cref{thm:coefficient-rg}
gives its unique mixing diagonal limit and its \GFF, with the parameters
of this bubble-decorated interaction. We emphasize that the bubble
correction is not treated as a small perturbation of the
coefficients after the fact: it enters the absolute
\texttt{thin}/\texttt{fat} bounds, the source identities and the
witness estimates of the three-step analysis from the outset, which
is why those estimates were proved for a general admissible contact
perturbation.

For the full Ising law with $M_n=\ceil{e^{n^2}}$, the conclusion about
the local image limit follows from the total-variation comparison in
\cref{thm:ising-rg-input}: pushforward by sorting and restriction
to a finite window cannot increase total variation, so every local
limit of the full-model \Tetris images coincides with the local
limit of the single-surface images. This identifies
the same infinite-volume tiling law. Its \GFF is then the conclusion
of \cref{thm:coefficient-rg}(ii) for that law. This follows the order of limits in \cref{thm:gff-limits}: first
construct the infinite-volume law, then rescale its height. Total
variation also controls the law of any measurable finite-volume
field, even one depending on a growing number of sites; however, a
simultaneous finite-torus scaling limit is not asserted.
The proof of \cref{thm:ising-gff-complete} below completes the physical transfer and
reconstructs the bulk spins.

\begin{proof}[Proof of \cref{thm:ising-gff-complete}]
\Cref{prop:ising-admissible-perturbation} verifies the
admissible contact hypotheses for the closed-contour correction.
Apply \cref{thm:interface-gff-complete} with this $U$.
It gives the unique physical single-surface history law, its sorting
identities, ergodicity, displacement moments, mean slopes and the
\GFF limit of the height function. Here $\upsigma_\mathrm{d}^\mathrm{Is}=\upsigma_\mathrm{d}^{U}$
and $L_\mathrm{d}^\mathrm{Is}=L_\mathrm{d}^{U}$. The height is the
column magnetization $\cS=\Tthree\cI$ of \cref{def:height-function-in-models}; no moment bound for the highest transition or
for an infinite sum of un-erased spins is required.

\Cref{prop:ising-contour-kernel} supplies the conditional spin
law $\cK(\cI)$ and local kernels with error
$C|A|e^{-c\beta r}$ on a fixed window $A$, uniformly once the
radius-$r$ neighborhood embeds in the quotient. Apply
\cref{lem:joint-limit} with the whole interface history as the
base and these spin kernels. The base limit is given by
\cref{thm:interface-gff-complete}; hence the histories and spins
converge jointly, with conditional spin law $\cK(\cI)$.

For the column-magnetization rooting of \cref{thm:gff-limits},
shift vertically by the decoded integer $\cS(0)-\phi(0)$.
The rerooting assertion of \cref{lem:physical-reconstruction}
applies also to the joint history and spin law: restricting this
integer to a finite interval leaves only finitely many translated
spin windows. Its limiting displacement tail makes the omitted
probability tend to zero.

\Cref{prop:ising-joint-cutoff} identifies the local limit
of the full Ising cylinders with this joint interface and spin law.
Its total-variation comparison is used only for local limits; the
\GFF is proved in the resulting infinite-volume height law by the
perturbed interface theorem.
\end{proof}

\subsection*{Acknowledgments}
The research of BL was supported by the grant \textsc{anr-18-ce40-0033 dimers}. The research of EL was supported by the NSF grant \textsc{dms}-2451083.

The project began in Sep 2024, soon after posting the first version of \cite{LaslierLubetzky24} on the arXiv. By Aug 2026, we were in a position where we had outlined the full argument that should lead to the \GFF limit of \HMSOS, then \SOS, then 3\Dim Ising; we had 
verified all the main calculations supporting the \GFF limit of the first step, the \HMSOS model, which we viewed as the main one (the \SOS and 3\Dim Ising models follow from similar arguments---once one controls fibers of one iteration of the \Tetris, this can be repeated relatively~painlessly). We had intended to carefully fill in the details over the next year(s), but the AI-related events of Oct 2026 forced us to change course. We prepared a 20-page draft of the paper stating all the main lemmas/propositions and sketching their proofs, and asked \textsc{ChatGPT-6 Astra} to check them and fill in all the missing details. It~did so while staying within our proof structure, even after we asked it to find simpler, more direct arguments. We then used \textsc{Claude Fable 5.1} for audits, and have been revising the paper since then. We are hastening to post this version, albeit still in a less-than-ideal form.

\bibliographystyle{abbrv}
\bibliography{tilted_interfaces}

\end{document}